\documentclass[fullpage]{book}

\usepackage{amsthm,amsfonts,amssymb,amsmath}
\usepackage{enumerate}
\usepackage[colorlinks=true]{hyperref}
\newcommand{\BA}{{\mathbb {A}}}

\newcommand{\BC}{{\mathbb {C}}}

\newcommand{\BF}{{\mathbb {F}}}

\newcommand{\BP}{{\mathbb {P}}}
\newcommand{\BQ}{{\mathbb {Q}}}
\newcommand{\BR}{{\mathbb {R}}}

\newcommand{\BZ}{{\mathbb {Z}}}

\newcommand{\CA}{{\mathcal {A}}}
\newcommand{\CB}{{\mathcal {B}}}
\newcommand{\CD}{{\mathcal {D}}}
\newcommand{\CE}{{\mathcal {E}}}
\newcommand{\CF}{{\mathcal {F}}}
\newcommand{\CG}{{\mathcal {G}}}
\newcommand{\CH}{{\mathcal {H}}}

\newcommand{\CK}{{\mathcal {K}}}
\newcommand{\CL}{{\mathcal {L}}}
\newcommand{\CM}{{\mathcal {M}}}
\newcommand{\CN}{{\mathcal {N}}}
\newcommand{\CO}{{\mathcal {O}}}

\newcommand{\CS}{{\mathcal {S}}}
\newcommand{\CT}{{\mathcal {T}}}
\newcommand{\CU}{{\mathcal {U}}}
\newcommand{\CV}{{\mathcal {V}}}
\newcommand{\CW}{{\mathcal {W}}}
\newcommand{\CX}{{\mathcal {X}}}
\newcommand{\CY}{{\mathcal {Y}}}
\newcommand{\CZ}{{\mathcal {Z}}}

\newcommand{\an}{{\mathrm{an}}}
\newcommand{\ran}{{\mathrm{r}{\text{-}\mathrm{an}}}}

\newcommand{\cod}{{\mathrm{cod}}}
\newcommand{\cont}{{\mathrm{cont}}}

\newcommand{\Cl}{{\mathrm{Cl}}}

\newcommand{\Div}{{\mathrm{Div}}}
\newcommand{\Divhat}{{\widehat{\mathrm{Div}}}}
\renewcommand{\div}{{\mathrm{div}}}
\newcommand{\CaCl}{{\mathrm{CaCl}}}
\newcommand{\eqv}{{\mathrm{eqv}}}
\newcommand{\cptf}{{\mathrm{cptf}}}

\newcommand{\id}{{\mathrm{id}}}

\newcommand{\emb}{{\hookrightarrow}}

\newcommand{\Gal}{{\mathrm{Gal}}}

\newcommand{\Hom}{{\mathrm{Hom}}}

\renewcommand{\Im}{{\mathrm{Im}}}

\newcommand{\intb}{{\mathrm{int}}}

\newcommand{\Isom}{{\mathrm{Isom}}}

\newcommand{\Jac}{{\mathrm{Jac}}}

\newcommand{\NS}{{\mathrm{NS}}}

\newcommand{\ord}{{\mathrm{ord}}}

\newcommand{\rank}{{\mathrm{rank}}}

\newcommand{\Pic}{\mathrm{Pic}}
\newcommand{\Picc}{{\mathcal{P}\mathrm{ic}}}
\newcommand{\Pichat}{\widehat{\mathrm{Pic}}}
\newcommand{\Prep}{\mathrm{Prep}}

\renewcommand{\Re}{{\mathrm{Re}}}
\newcommand{\reg}{{\mathrm{reg}}}

\newcommand{\rmod}{{\mathrm{mod}}}

\DeclareMathOperator{\Spec}{Spec}

\newcommand{\tor}{{\mathrm{tor}}}

\newcommand{\tr}{{\mathrm{tr}}}

\newcommand{\vol}{{\mathrm{vol}}}

\newcommand{\sm}{{\mathrm{sm}}}

\newcommand{\red}{{\mathrm{red}}}

\newcommand{\Ly}{{\mathrm{Ly}}}

\newcommand{\wt}{\widetilde}
\newcommand{\wh}{\widehat}
\newcommand{\pp}{\frac{\partial\overline\partial}{\pi i}}
\newcommand{\pair}[1]{\langle {#1} \rangle}

\newcommand{\ds}{\displaystyle}

\newcommand{\ol }{\overline}

\newcommand{\lra}{\longrightarrow}

\newcommand{\iso}{{\overset\sim\lra}}

\newcommand{\fh}{\mathfrak{h}}

\newcommand{\naive}{\mathrm{naive}}
\newcommand{\nef}{\mathrm{nef}}
\newcommand{\snef}{\mathrm{snef}}
\renewcommand{\vert}{\mathrm{vert}}

\newcommand{\bif}{\mathrm{bif}}

\newcommand{\Fal}{\mathrm{Fal}}

\newcommand{\CLL}{\overline{\mathcal L}}
\newcommand{\CMM}{\overline{\mathcal M}}

\newcommand{\CNN}{\overline{\mathcal N}}
\newcommand{\CAA}{\overline{\mathcal A}}
\newcommand{\CBB}{\overline{\mathcal B}}
\newcommand{\CDD}{\overline{\mathcal D}}
\newcommand{\CEE}{\overline{\mathcal E}}
\newcommand{\CFF}{\overline{\mathcal F}}

\newcommand{\OCD}{\overline{\mathcal D}}
\newcommand{\OCE}{\overline{\mathcal E}}

\newcommand{\OL}{{\overline{L}}}
\newcommand{\OM}{{\overline{M}}}
\newcommand{\OD}{{\overline{D}}}
\newcommand{\ON}{{\overline{N}}}
\newcommand{\OH}{{\overline{H}}}
\newcommand{\OK}{{\overline{K}}}

\newcommand{\OA}{{\overline{A}}}
\newcommand{\OB}{{\overline{B}}}
\newcommand{\OQ}{{\overline{Q}}}
\newcommand{\OP}{{\overline{P}}}
\renewcommand{\OE}{{\overline{E}}}

\newcommand{\TL}{{\widetilde{L}}}
\newcommand{\TM}{{\widetilde{M}}}

\newcommand{\CMZ}{\mathcal{M}(\mathbb{Z})}
\newcommand{\CMk}{\mathcal{M}(k)}

\newcommand{\CC}{\mathbb{C}}
\newcommand{\RR}{\mathbb{R}}
\newcommand{\ZZ}{\mathbb{Z}}
\newcommand{\QQ}{\mathbb{Q}}
\newcommand{\FF}{\mathbb{F}}
\newcommand{\PP}{\mathbb{P}}

\newcommand{\kkk}{Let $k$ be either $\ZZ$ or a field. }
\newcommand{\ccc}{Take the uniform terminology in \S\ref{sec uniform}. }

\theoremstyle{theorem}
\newtheorem{thm}{Theorem}[section]
\newtheorem{cor}[thm]{Corollary}
\newtheorem{lem}[thm]{Lemma}
\newtheorem{prop}[thm]{Proposition}
\newtheorem{conj}[thm]{Conjecture}

\theoremstyle{definition}
\newtheorem{defn}[thm]{Definition}

\theoremstyle{remark}
\newtheorem{rmk}[thm]{Remark} 
\newtheorem{problem}[thm]{Problem}

\begin{document}

\title{Adelic Line Bundles on Quasi-projective Varieties}
\author{Xinyi Yuan, Shou-Wu Zhang}

\maketitle

\thispagestyle{empty}
\begin{flushleft}
In memory of Lucien Szpiro (1941--2020).
\end{flushleft}

\tableofcontents



\chapter*{Preface}
\addcontentsline{toc}{chapter}{Preface}

This volume is devoted to a theory of adelic line bundles on quasi-projective varieties over finitely generated fields. 
Our first motivation is to apply Arakelov theory to treat height functions on moduli spaces, such as moduli spaces of smooth projective varieties, vector bundles, dynamical systems, etc, which are usually quasi-projective. 
Our second motivation is to apply Arakelov theory to study some properties of quasi-projective varieties over arbitrary fields since they can be descended to finitely generated fields by the Lefschetz principle.

Our work will be an extension of the theory of adelic line bundles on projective varieties of Zhang \cite{Zha2}, which itself is an extension of the Arakelov theory of hermitian line bundles on integral models of projective varieties initiated by Arakelov \cite{Ara74}, and developed by Faltings \cite{Fal2}, Deligne \cite{Del}, Szpiro \cite{Szp85, Szp90},  and Gillet--Soul\'e \cite{GS1}. These adelic line bundles are defined to be certain limits of hermitian line bundles on integral models. The theory has been applied to N\'eron--Tate heights on abelian varieties or canonical heights on arithmetic dynamical systems. In particular, it has been used as a crucial tool to treat the equidistribution theorem by Szpiro--Ullmo--Zhang \cite{SUZ},
 the Bogomolov conjecture  by 
Ullmo \cite{Ull} and Zhang \cite{Zha4},
and the rigidity theorem of preperiodic points in dynamical systems and the non-archimedean Monge--Amp\'ere equation by Yuan--Zhang \cite{YZ1}.

The adelic line bundles in this volume will be defined as the limits of hermitian line bundles on integral models of the projective compactifications.
Both compactifications and integral models vary during the limit process.   We will define two intersection pairings among these adelic line bundles: Gillet--Soul\'e's intersection number in the absolute setting and the Deligne pairing in the relative setting. We will study relations among intersections, heights, volumes, and positivity. We will prove an equidistribution theorem of small points on quasi-projective varieties over number fields, generalizing the equidistribution theorems of Szpiro--Ullmo--Zhang \cite{SUZ}, Chambert-Loir \cite{CL}, and Yuan \cite{Yua1}.

\medskip

\noindent\textit{Acknowledgments.}
We would like to dedicate this volume to the memory of Professor Lucien Szpiro (1941--2020), who passed away during the preparation of this volume.
His visionary article \cite{Szp90} is the main inspiration for our previous and present works on the Arakelov theory. Shou-Wu Zhang is deeply indebted to him for his mentorship and friendship.

The authors would like to thank 
Jos\'e Ignacio Burgos Gil, 
Alexander Carney, Ziyang Gao, Thomas Gauthier, Niki Myrto Mavraki, Yanshuai Qin,   
and Junyi Xie for their helpful discussions. The authors are also grateful to Yulin Cai, Shiquan Li, Yanshuai Qin, and Junyi Xie for pointing out many mistakes in early versions of the book.

Xinyi Yuan is supported by grant NO. 12250004 from the National Science Foundation
of China and the Xplorer Prize from the New
Cornerstone Science Foundation.
Shou-Wu Zhang is supported by the NSF awards DMS-0970100, DMS-1065839, and DMS-1700883 of the USA.

\vspace{2pt}

\begin{flushright}
Xinyi Yuan, Shou-Wu Zhang\\
October 31, 2024
\end{flushright}

\cleardoublepage



\chapter{Introduction}
The notion of height functions on projective varieties was created by Weil \cite{Wei51} as a measurement for the complexity of solutions to Diophantine equations.  Before its modern name, this notion has its roots in  Euclid's proof of the 
 irrationality of $\sqrt 2$ more than two thousand years ago and Fermat's proof that  $x^4+y^4=z^4$ has no solutions in positive integers
 more than three hundred years ago. It was implicitly used in Weil's proof of Mordell-Weil theorem in \cite{Wei29} in 1928, Siegel's theorem on integral points on curves in 1929, and Northcott's further works on the Northcott property and arithmetic dynamical systems in \cite{Nor49, Nor50} in 1949-1950.
 Since its creation, the theory of heights has been widely used in diophantine geometry, such as the Roth theorem and Schmidt's subspace theorems and the formulation of Birch and Swinnerton-Dyer conjecture.  For more details, see our reviews in \S\ref{app sec height}.

In the 1970s, to translate the proof of Mordell's conjecture from function fields to number fields, Arakelov \cite{Ara74} proposed an intersection theory for arithmetic surfaces. In this theory, the heights are interpreted as degrees of hermitian line bundles on arithmetic curves. Faltings \cite{Fal2} used this interpretation to define a special height function  in his proof of the Mordell conjecture in 1983. After a quick development by Faltings \cite{Fal1}, Deligne \cite{Del}, Szpiro \cite{Szp85}, and Gillet--Soul\'e \cite{GS1}, the Arakelov theory was used again by Vojta \cite{Voj}  for a second proof of the Mordell conjecture and by Faltings \cite{Fal91} for an extension to the Mordell--Lang conjecture for the subvarieties of Abelian varieties. For more details, see our reviews in \S\ref{app sec surface} and \S\ref{app sec intersection}.

In the late 1980s, Szpiro \cite{Szp90} proposed a program to prove the Bogomolov conjecture for small points on curves using arithmetic positivity. The Bogomolov conjecture was eventually reduced to a positivity statement after the development of arithmetic ampleness for adelic line bundles by Zhang \cite{Zha92, Zha93, Zha1, Zha2}, which was an extension of classical Arakelov theory to handle the bad reductions of abelian varieties and algebraic dynamical systems. Shortly after, the Bogomolov conjecture was eventually proved by Ullmo \cite{Ull} and Zhang \cite{Zha4}
using a new ingredient, the equidistribution theorem initiated by Szpiro--Ullmo--Zhang \cite{SUZ}. The arithmetic positivity and equidistribution theorem were further extended by the Yuan \cite{Yua1} to big line bundles and had been widely used to treat problems in arithmetic dynamical systems. For more details, see our reviews in \S\ref{app sec positivity}, \S\ref{app sec adelic},
and \S\ref{app sec equidistribution}.

The goal of this book is to extend the theory of Zhang \cite{Zha2} from projective varieties over number fields to \emph{quasi-projective} varieties over \emph{finitely generated fields}. More precisely, let $F$ be a finitely generated field over $\QQ$ (or a constant field), and let $X$ be a quasi-projective variety over $F$. We introduce a notion of \emph{adelic line bundles} on $X$, consider their intersection theory, study their volumes for effective sections, and introduce heights associated with them. 
The fundamental properties behind these terms are the positivity of adelic line bundles.

An immediate application of our framework is a theory of canonical heights on polarized algebraic dynamical systems over quasi-projective varieties over finitely generated fields. In particular, we introduce N\'eron--Tate heights of abelian varieties over finitely generated fields and extend the arithmetic Hodge index theorem of Faltings \cite{Fal1} and Hriljac \cite{Hri} to this setting. 
Furthermore, we prove an equidistribution theorem of small points on quasi-projective varieties over number fields, generalizing the equidistribution theorems of Szpiro--Ullmo--Zhang \cite{SUZ}, Chambert-Loir \cite{CL}, and Yuan \cite{Yua1}.

The exposition of this book uses a combination of algebraic geometry, complex algebraic geometry, Arakelov theory (cf. \cite{Ara74, GS1}), and Berkovich analytic spaces (cf. \cite{Ber1, Ber2}). 
In the following, we sketch the main constructions and theorems of this book.

\section{Adelic line bundles}

To illustrate the concept quickly, we will take an approach different from the major parts of this book, but it will give equivalent constructions. 

We will use a \emph{uniform terminology}, which will be explained in detail in \S\ref{sec uniform},  to treat both algebraic varieties and arithmetic varieties. 
Namely, we fix a base scheme $\Spec k$, where $k$ is either $\BZ$ or a field. By a {\em variety $X$ over $k$}, we mean an integral scheme $X$ which is finite, flat, and of finite type over $\Spec k$. We say that $X$ is an arithmetic (resp. algebraic ) variety if $k=\BZ$ (resp. $k$ is a field).  We say that $X$ is a  {\em projective variety over $k$} if the structure morphism $X\to \Spec k$ is projective; we say that $X$ is a {\em quasi-projective variety over $k$} if $X$ is an open subscheme of a projective variety over $k$. 
In particular, when $X$ is a quasi-projective variety over $\BZ$, 
we allow finitely many fibers of $X\to \Spec \BZ$ to be empty. 

In the arithmetic situation, we could take a fancier notation $k=\BF_1$, the field of one element,  and thus $\Spec \BZ$ becomes an {\em affine arithmetic curve} over $\BF_1$. 
The curve $\Spec \BZ$ can be further compactified over $\BF_1$ by adding an archimedean place.  So projective varieties over $\ZZ$ can be further {\em compactified} by adding complex varieties over archimedean places. This is the main motivation for Arakelov to introduce hermitian line bundles on projective varieties; see our reviews in \S\ref{app sec surface} and \S\ref{app sec intersection}.  Sometimes, we need to deal with projective varieties over $\BZ[1/N]$ for some positive integer $N$, for example, when we treat abelian varieties and algebraic dynamical systems. In this case, we need to add $p$-adic metrics to the missing primes. This is the main motivation for Zhang to introduce adelic line bundles for projective varieties; see our reviews in \S\ref{app sec adelic}.

The goal of this section is to sketch our theory of adelic line bundles on a quasi-projective variety  $X$ over $k$.
Adelic line bundles are roughly defined as the ``limits" of ``projective models" under our ``boundary topology." 
There are two natural approaches to define these limits. 
The first approach is by an abstract notion of completion by Cauchy sequences (combined with two processes of direct limits), and this is the main approach of this book precisely realized in Chapter \ref{sec adelic}. 
The second approach is to ``put'' all ``projective models'' into the category $\wh\Picc(X^\an)$ of metrized line bundles on a large analytic space $X^\an$. It turns out that 
$\wh\Picc(X^\an)$ is big enough to contain all the limiting line bundles. 
This is essentially treated in Chapter \ref{sec analytic}. 
Each approach has its advantages. 
We will take the second approach in this introduction.

\subsection{Berkovich spaces}
\kkk
Let $X/k$ be a quasi-projective variety. 
There is a natural \emph{Berkovich analytic space} $X^\an$ associated to $X/k$. 
In fact, if $X$ has an open affine cover $\{\Spec A_i\}_i$, then $X^\an=\cup_i \CM(A_i)$, where $\CM(A_i)$ is the set of multiplicative semi-norms $|\cdot|$ on $A_i$; if $k$ is a field, we further require the restriction of $|\cdot|$ to $k$ is trivial.
A \emph{metrized line bundle on $X^\an$} is a pair $\OL=(L,\|\cdot\|)$ consisting of a line bundle $L$ on $X$ and a continuous metric $\|\cdot\|$ of $L$ on $X^\an$. Denote by $\wh \Picc (X^\an)$ the category of metrized line bundles on $X$, in which a morphism between two objects is defined to be an isometry.
There is a forgetful functor
$$\wh \Picc (X^\an)\lra \Picc (X).$$
Here, $\Picc (X)$ denotes the category of line bundles on $X$, in which a morphism between two objects is an isomorphism of line bundles.

\subsection{Model adelic line bundles}
\kkk
Let $X/k$ be a quasi-projective variety. 
Objects of the category $\wh \Picc (X^\an)$ are too general for intersection theory. Instead, we will define a full subcategory 
$\wh\Picc(X/k)$ 
of adelic line bundles in $\wh \Picc (X^\an)$, and a full subcategory 
$\wh\Picc(X/k)_\intb$ 
of integrable adelic line bundles in $\wh \Picc (X^\an)$ for intersection theory.
For this, we will start with model adelic line bundles and take a limit process to extend them to more general notions.

As a convention, we will write tensor products of various line bundles additively, so for example, $mL$ means $L^{\otimes m}$ for $L\in \Picc(X)$ and $m\in \ZZ$.

An object of $\wh \Picc (X^\an)$ with underlying line bundle $L\in \Picc(X)$ is called a {\em model adelic line bundle} if it is induced by a projective model $(\CX, \CLL)$ of $(X, eL)$ over $k$ for some positive integer $e$, where $\CX$ is a projective variety over $k$  with an open immersion $X\emb \CX$, and $\CLL$ is as follows:
\begin{enumerate}
\item if $k$ is a field, then $\CLL=\CL$ is a line bundle on $\CX$ extending $eL$;
\item if $k=\BZ$, then $\CLL=(\CL, \|\cdot\|)$ is a hermitian line bundle on $\CX$, which consists of a line bundle $\CL$ on $\CX$ extending $eL$ and a continuous hermitian metric $\|\cdot\|$ of $\CL(\CC)$ on $\CX(\CC)$. The metric is required to be invariant under the complex conjugate.
\end{enumerate}
Because of the integer $e$ in the definition, $(\CX, e^{-1}\CLL)$ is a projective model of $(X, L)$ in terms of the notion of $\QQ$-line bundles.
Denote by $\wh\Picc (X/k)_\rmod$ the full subcategory of $\wh\Picc (X^\an)$ consisting of {model adelic line bundles} on $X$.

\subsection{Limit process}
\kkk
Let $X/k$ be a quasi-projective variety. 
Choose a projective compactification $X\subset \CX_0$  such that the boundary $\CX_0\setminus X$ is exactly equal to the support of an effective Cartier divisor $\CE_0$ on $\CX_0$. 
If $k$ is a field,  set $\OCE_0=\CE_0$.
If $k=\BZ$ set $\OCE_0=(\CE_0,g_0)$, where $g_0>0$ is a Green function of $\CE_0(\CC)$ on $\CX_0(\CC)$. 
Then $\OCE_0$ induces a Green function $\wt g_0$ of $\CE_0$ on $\CX_0^\an$, which restricts to a continuous function $\wt g_0:X^\an \to \RR_{\geq0}$.

Consider the space $C(X^\an)$ of real-valued continuous functions on $X^\an$. 
It is endowed with a \emph{boundary topology} induced by the extended norm
$$\|f\|_{\wt g_0}:=\sup_{x\in X^\an,\ \wt g_0(x)>0} \frac{|f(x)|}{\wt g_0(x)}.$$
We refer to \cite{Bee} for the basics of extended norms, which are allowed to take values in $[0,\infty]$ but still required to satisfy the triangle inequality.
The boundary topology is independent of the choice of $(\CX_0,\OCE_0)$. 
Moreover, $C(X^\an)$ is complete under the boundary topology.

We say that a sequence 
$\CLL_i=(\CL_i,\|\cdot\|_i)$ in $\wh\Picc (X^\an)$
\emph{converges} to an object $\CLL=(\CL,\|\cdot\|)$ in $\wh\Picc (X^\an)$ if there are isomorphisms $\tau_i:\CL\to \CL_i$ such that the sequence
$-\log (\tau_i^*\|\cdot\|_i/\|\cdot\|)$ converges to 0 in $C(X^\an)$ under the boundary topology.

\subsection{Adelic line bundles} 

There is a notion of nefness of hermitian line bundles on projective arithmetic varieties, and we refer to \S\ref{app subsec ample} for a quick definition.

An object of $\wh \Picc (X^\an)$ is called an \emph{adelic line bundle on $\CU$} if it is isomorphic to the limit of a sequence in $\wh\Picc (X/k)_{\rmod}$.
An adelic line bundle on $X$ is called \emph{strongly nef} if it is isomorphic to the limit of a sequence 
in $\wh\Picc (X/k)_{\rmod}$ induced by projective models $(\CX_i, \CLL_i)$ over $k$ such that $\CLL_i$ is nef on $\CX_i$.
An adelic line bundle $\CLL$ on $X$ is called \emph{nef} if there exists a strongly nef adelic line bundle $\CMM$ on $X$ such that $a\CLL+\CMM$ is strongly nef for all positive integers $a$. An adelic line bundle on $X$ is called \emph{integrable} if it is isometric to 
$\CLL_1- \CLL_2$ for two {strongly nef} adelic line bundle $\CLL_1$ and $\CLL_2$ on $X$. 

Denote by $\wh \Picc (X/k)$ the full subcategory of $\wh \Picc (X^\an)$ consisting of adelic line bundles on $X$.
Denote by $\wh \Picc (X/k)_\nef$ (resp. $\wh \Picc (X/k)_\intb$) the full subcategory of $\wh \Picc (X^\an)$ consisting of nef (resp. integrable) adelic line bundles on $X$. Their objects are called \emph{adelic line bundles} (resp. \emph{nef adelic line bundles}, \emph{integrable adelic line bundles}) on $X/k$.

We can further extend the definition to quasi-projective varieties over finitely generated fields. Namely, let $F$ be a finitely generated field over $k$, i.e., a finitely generated field over the fraction field of $k$. 
Let $X$ be a quasi-projective variety over $F$. Then we define 
$$\wh \Picc (X/k): =\varinjlim_{U\to V} \wh \Picc (U/k), $$
$$\wh \Picc (X/k)_\nef: =\varinjlim_{U\to V} \wh \Picc (U/k)_\nef, $$
$$\wh \Picc (X/k)_\intb :=\varinjlim_{U\to V} \wh \Picc (U/k)_\intb.$$
Here, the limit is over all flat morphisms $U\to V$ of quasi-projective varieties over $k$ whose generic fibers are isomorphic to $X\to \Spec F$. 
Denote by $\wh \Pic(X/k)$ (resp. $\wh \Pic(X/k)_\nef$, $\wh \Pic (X/k)_\intb$, $\wh \Pic (X^\an)$) the \emph{group} of isomorphism classes of objects of
$\wh \Picc (X/k)$ (resp. $\wh \Picc (X/k)_\nef$, $\wh \Picc (X/k)_\intb$, $\wh \Picc (X^\an)$). 
Note that the previous definitions of the analytic terms $X^\an$ and $\wh \Picc (X^\an)$ are actually valid in the current situation. 

As we have seen, our theory of adelic line bundles is valid for both quasi-projective varieties over $k$ and quasi-projective varieties over finitely generated fields over $k$.
We will introduce a natural notion of \emph{essentially quasi-projective varieties over $k$}, which includes both of the above cases.
For simplicity, we will not use this notion in this chapter.

\subsection{Functoriality}

Let $E/F$ be an extension of finitely generated fields over $k$, and  $f: X\to Y$ be an $F$-morphism of quasi-projective varieties $X/E$
and $Y/F$. Then we have a pull-back functor
$$f^*: \wh \Picc (Y/k)\lra \wh\Picc (X/k).$$

When $k=\BZ$, we can also have a base change of a quasi-projective scheme $X/k$ to the generic fiber $X_\BQ/\BQ$. In this case, we denote the functor as 
$$ \wh \Picc (X/k)\lra \wh\Picc (X_\BQ/\BQ),\quad
\OL\longmapsto \wt L.$$
We call $\wt L$ \emph{the geometric part of} $\wh L$. 

Both functors preserve the subcategories of the model (resp. nef, integrable) adelic line bundles.

\section{Intersection theory and heights}

Our intersection theory includes an absolute intersection pairing of Gillet--Soul\'e and a relative intersection pairing that extends the Deligne pairing.

\subsection{Intersection numbers and heights }
\kkk
Let $F$ be a finitely generated field over $k$.
Let $X$ be a quasi-projective variety over $F$.
Then our absolute intersection pairing is a symmetric and multi-linear map
$$
\wh \Pic(X/k)_\intb^{\, d} \longrightarrow \RR,
$$
where $d=\dim X+\dim k+\tr\deg_k F$. Here $\dim k$ denotes the Krull dimension of $k$, and 
$\tr\deg_k F$ denotes the transcendence degree of $F$ over the fraction field of $k$. 
This is the limit version of the intersection theory in algebraic geometry and the arithmetic intersection theory of Gillet--Soul\'e.
See Proposition \ref{intersection1}.

Now let $K$ be a number field if $k=\BZ$; let $K$ be a function field of one variable over $k$ if $k$ is a field.
Let  $X$ be a quasi-projective variety over $K$ of dimension $n$.
Let $\overline L$  be an integrable adelic line bundle on $X$.
Define a \emph{height function}
$$
h_{\overline L}: X(\overline K) \longrightarrow \BR
$$
by
$$
h_{\overline L}(x):
= \frac{ \wh\deg(\overline L|_{x'})}{\deg(x')}.$$
Here $x'$ denotes the closed point of $X$ containing $x$,
$\deg(x')$ denotes the degree of the residue field of $x'$ over $K$,
$\overline L|_{x'}$ denotes the \emph{pull-back} of $\OL$ to 
$\wh\Pic(x'/k)_{\intb}$,
and $\wh\deg:\wh\Pic(x'/k)_{\intb}\to \RR$ is by the intersection theory.

More generally, for any closed $\overline K$-subvariety $Z$ of $X$, define \emph{the  height of $Z$ for  $\overline L$} by
$$
h_{\overline L}(Z):
= \frac{ (\overline L|_{Z' })^{\dim Z+1}}
{(\dim Z+1)(\wt L|_{Z'_K})^{\dim Z} }.$$
Here $Z'$ denotes the image of $Z\to X$ (which is a closed subvariety of $X$ over $K$), and 
$$\OL\longmapsto \overline L|_{Z'}\longmapsto\wt L|_{Z'_K}$$
denotes the image of $\OL$ via the functorial maps
$$\wh\Pic( X/k)_{\intb}\lra \wh\Pic(Z'/k)_{\intb}\lra  \wh\Pic(Z'_K/K)_{\intb},$$
and the self-intersections are as in the above intersection theory.
The height is well-defined only if the denominator is nonzero.

\subsection{Deligne pairing and relative heights}
\kkk
Let $f: X\to Y$ be a projective and flat morphism of relative dimension $n$ between quasi-projective varieties over $k$. 
Assume that $Y$ is normal, which is required in our proof.
\begin{thm} [Theorem \ref{intersection2}]
The Deligne pairing on model adelic line bundles induces a symmetric and multilinear functor
$$\wh \Picc(X)_\intb^{\, n+1} \longrightarrow \wh \Picc (Y)_\intb.
$$
When restricted to nef adelic line bundles, the functor induces a functor
$$\wh \Picc(X)_\nef^{\, n+1} \longrightarrow \wh \Picc (Y)_\nef.$$
Moreover, the functors are compatible with base changes of the form $Y'\to Y$, where $Y'$ is a quasi-projective normal variety over $k$ such that $X'=X\times_YY'$ is integral. 
\end{thm}

In the setting of the theorem, let $F=k(Y)$ be the function field of $Y$, and $X_F\to \Spec F$ the generic fiber of $X\to Y$.
Let $\overline L$  be an object of $\wh \Picc (X)_\intb$.
By this, we can define a \emph{vector-valued} height function
$$
\fh_{\overline L}: X(\overline F) \longrightarrow \wh\Pic (F/k)_{\intb, \QQ}.
$$
Here the group
$$
\wh\Pic (F)_{\intb}:=\varinjlim _U\wh \Pic (U/k)_{\intb},
$$
where $U$ runs through all open subschemes of $Y$.

More generally, for any closed $\overline F$-subvariety $Z$ of $X_F$, define \emph{the vector-valued height of $Z$ for $\overline L$} as
$$
\fh_{\overline L}(Z):
= \frac{ \pair{\overline L|_{Z' }}^{\dim Z+1}}
{(\dim Z+1)(\OL|_{Z'_{F} })^{\dim Z} }\ \in \wh\Pic (F)_{\intb, \QQ}.$$
Here $Z'$ denotes the image of $Z\to X$, $Z'_{F}$ is the generic fiber of $Z'\to Y$, and
$$\OL\longmapsto \overline L|_{Z'}\longmapsto \overline L|_{Z'_{F}}$$
denotes the image of $\OL$ via the functorial maps
$$\wh\Pic( X/k)_{\intb}\lra \wh\Pic(Z'/k)_{\intb}\lra  \Pic(Z'_{F}/F).$$
Note that the first self-intersection is the Deligne pairing, and the second self-intersection is just the degree on the projective variety $Z'_{F}$ in the classical sense.
The height is well-defined only if the denominator is nonzero.

When $F$ is polarized in the sense of Moriwaki \cite{Mor3}, then we can also define the Moriwaki heights. If $K$ is a number field or a finite field, and if $X$ is projective over $F$, we obtain a Northcott property of the Moriwaki heights from that of \cite{Mor3}. 
In general, we obtain the fundamental inequality for the Moriwaki height following the strategy of \cite{Mor3}. 

Hence, part of our height theory extends the previous works of Moriwaki \cite{Mor3, Mor4}. In fact,  \cite{Mor3, Mor4} developed a height theory for projective varieties over finitely generated fields $F$ over $\QQ$, depending on the choice of an arithmetic polarization of $\Spec F$. His motivation was to apply Arakelov geometry to varieties over arbitrary fields (of characteristic 0), and he succeeded in formulating and proving the Bogomolov conjecture in that setting.
His treatment was more on the numerical theory of heights, but ours is more on the geometric theory of adelic line bundles.

\section{Volumes and equidistribution}

As in the projective case, we can define effective sections of adelic line bundles, study their volumes, and prove equidistribution theorems on quasi-projective varieties. 

\subsection{Volumes}
Let $X$ be a quasi-projective variety over $k$.
Let $\OL=(L,\|\cdot\|)$ be an adelic line bundle on $X$.
Define  
$$\wh H^0(X, \overline L):=\{s\in H^0(X, L): \|s(x)\|\leq 1,\ \forall \, x\in X^\an\}.$$
Elements of $\wh H^0(X, \overline L)$ are called \emph{effective sections} of $\overline L$ on $X$.
If $k=\BZ$, denote
$$\wh h^0(X, \OL):=\log\#\wh H^0(X, \overline L);$$
if $k$ is a   field, denote
$$\wh h^0(X, \OL):=\dim_k \wh H^0(X, \overline L).$$
We check that $\wh h^0(X, \overline L)$ is always a finite real number. 
In this setting, we have the following fundamental results. 

\begin{thm} [Theorem \ref{limit}, Theorem \ref{HS}]
Let $X$ be a quasi-projective variety over $k$.
Let $\OL,\OM$ be adelic line bundles on $X$.
Denote $d=\dim X+\dim k$. 
Then the following holds.
\begin{enumerate}
\item The limit 
$$
\wh\vol(X,\OL)
=\lim_{m\to \infty} \frac{d!}{m^{d}}\wh h^0(X, m\OL)
$$
exists. 
\item
If $\OL$ is the limit of a sequence of model adelic line bundles induced by a sequence $\{(\CX_i,\overline \CL_i)\}_{i\geq1}$ of projective models of $(X, L)$ over $k$, then 
$$
\wh\vol(X,\OL)
=\lim_{i\to \infty} \wh\vol(\CX_i,\CLL_i).
$$
\item
If $\OL$ is nef, then
$$
\wh\vol(X,\OL) = \OL^d.
$$ 
\item
If $\OL,\OM$ are nef, then 
$$
\wh\vol(X,\OL-\OM) \geq \OL^d-d\,\OL^{d-1}\OM.
$$
\end{enumerate}
\end{thm}

Part (1) generalizes the classical result of Fujita (cf. \cite[11.4.7]{Laz2}) for line bundles on projective varieties and the result of \cite{Che1, Che2, Yua2} for hermitian line bundles on projective arithmetic varieties.
Part (2) allows us to transfer many previous results in the projective case to the quasi-projective case.
Part (3) generalizes the classical Hilbert--Samuel formula in algebraic geometry and the arithmetic Hilbert--Samuel formula proved by Gillet--Soul\'e \cite{GS2}, Bismut--Vasserot \cite{BV}, and Zhang \cite{Zha1}.
Part (4) generalizes the classical theorem of Siu \cite{Siu} and the arithmetic bigness theorem of Yuan \cite{Yua1}. 

In the setting of the theorem, we say that $\OL$ is \emph{big} 
if $\wh\vol(X,\OL)>0$.
We will see that in this case, we will have nice lower bounds of the height function associated with $\OL$.


\subsection{Height inequality}

Let $K$ be a number field if $k=\BZ$; let $K$ be a function field of one variable over $k$ if $k$ is a field. 
Let $X$ be a quasi-projective variety over $K$.
For an adelic line bundle $\overline L$ on $X/k$, we usually denote by $\wt L$ 
the \emph{geometric part} of $\OL$ on $X/K$, i.e.
 the image of 
$\OL$ under the functorial map $ \wh \Picc (X/k)\to \wh\Picc (X/K)$.

As a quick consequence of the above fundamental results on volumes, we have the following height inequality. 

\begin{thm}[Theorem \ref{height comparison0}] \label{height comparison intro0}
Let $\pi: X\to S$ be a morphism of quasi-projective varieties over $K$. 
Let $\OL\in \wh\Pic(X)$ and $\OM\in \wh\Pic(S)$ be adelic line bundles. 
 If $\OL$ is nef on $X/k$ and $\wt L$ is big on $X/K$, then 
 for any $c>0$, there exist $\epsilon>0$ and a non-empty open subvariety $U$ of $X$ such that 
$$
h_{\OL}(x) \geq \epsilon\, h_{\OM}(\pi(x)) -c, \quad\ \forall\, x\in U(\overline K). 
$$

\end{thm}

We refer to Theorem \ref{height comparison0} for various versions of the height inequality and to Theorem \ref{thm-ht-big} for a partial converse to the height inequality.

In a series of works, Dimitrov--Gao--Habegger \cite{GH, DGH} and K\"uhne \cite{Kuh} proved a uniform Bogomolov conjecture over number fields, and combined the work of Vojta \cite{Voj} on the Mordell conjecture to confirm Mazur's uniform Mordell conjecture.  
A key result of \cite{DGH} is a height inequality in a setting of abelian schemes, which also plays a fundamental role in the further work \cite{Kuh}. Our current height inequality can be viewed as a theoretical version of that of \cite{DGH}. 
We will come back to this connection later, and we refer to the context of Theorem \ref{height comparison} for more details on this connection and for the application of our height inequality to dynamical systems.

In recent work, Yuan \cite{Yua4} has used our theory of adelic line bundles to prove the uniform Bogomolov conjecture over global fields. This gives a new proof of the main results of \cite{DGH, Kuh}, and it works uniformly for both number fields and function fields of any characteristic.
The key ingredient of \cite{Yua4} is to prove the bigness of the admissible canonical bundle of the universal curve using the Deligne pairing and apply the bigness to obtain certain height inequality to control the number of points of small heights. We refer to \S\ref{sec admissible metrics} for a brief introduction to the admissible canonical bundle.

In recent work, Gao--Zhang \cite{GZ24} proved a Northcott property for Gross--Schoen cycles and Ceresa cycles parametrized by 
a non-empty open subset of moduli spaces of curves of genus at least $3$.
One key ingredient in their proof is applying our height inequality to convert the positivity properties of adelic line bundles to a Northcott property.

\subsection{Equidistribution}

One of the most important theorems of this book is an equidistribution theorem for small points of a quasi-projective variety over a number field or a function field of one variable.

\begin{thm}[Theorem \ref{equi3}]
\kkk
Let $K$ be a number field if $k=\ZZ$; let $K$ be the function field of one variable over $k$ if $k$ is a field. Let $X$ be a quasi-projective variety over $K$.
Let $\overline L$ be a nef adelic line bundle on $X$ such that 
$\deg_{\TL}(X)>0$. 
Let $\{x_m\}_m$ be a generic sequence in $X(\overline K)$ such that $\{h_{\OL}(x_m)\}_m$ converges to $h_{\OL}(X)$.
Then the Galois orbit of $\{x_m\}_m$ is equidistributed in $X_{K_v}^\an$ for $d\mu_{\overline L,v}$ for any place $v$ of $K$.
\end{thm}

Here $d\mu_{\overline L,v}$ is a canonical probability measure on $X_{K_v}^\an$, 
defined using the recent theory of Chambert-Loir and Ducros in \cite{CLD} if $v$ is non-archimedean. This generalizes the Monge--Amp\`ere measure and the Chambert-Loir measure from the projective case to the quasi-projective case.

If $k=\ZZ$ and $X$ is projective over $K$, the equidistribution theorem is proved by Szpiro--Ullmo--Zhang \cite{SUZ}, Chambert-Loir \cite{CL}, and Yuan \cite{Yua1}.
Our current theorem still follows the variational principle of the pioneering work \cite{SUZ}, applying our
adelic Hilbert--Samuel formula and adelic bigness theorem.

We can further generalize our equidistribution theorem in two different aspects, which give us an equidistribution theorem (Theorem \ref{equi5}) and an equidistribution conjecture (Conjecture \ref{equi1}).
The equidistribution theorem considers a projective and flat morphism of quasi-projective varieties over a number field or a function field of one variable, and its proof follows a strategy of Moriwaki \cite{Mor3}. 
The equidistribution conjecture considers quasi-projective varieties over finitely generated fields, which is stated as follows.

\begin{conj}[Conjecture \ref{equi1}]
\kkk
Let $F$ be a finitely generated field over $k$. 
Let $v$ be a non-trivial valuation of $F$.
Assume that the restriction of $v$ to $k$ is trivial if $k$ is a field.
Let $X$ be a quasi-projective variety over $F$.
Let $\overline L$ be a nef adelic line bundle on $X$ such that $\deg_{\TL}(X)>0$, where $\wt L$ denotes the image of $\OL$ under the functorial map
$\wh\Pic(X/k) \to \wh\Pic(X/F)$.
Let $\{x_m\}_m$ be a generic sequence of small points in $X(\overline F)$.
Then the Galois orbit of $\{x_m\}_m$ is equidistributed in $X_{F_v}^\an$ for $d\mu_{\overline L,v}$.
\end{conj}

We refer to the context of Conjecture \ref{equi1} for the notion of ``small points'' and the equilibrium measure $d\mu_{\overline L,v}$.

\section{Algebraic dynamics}

Here we apply the theory of adelic line bundles to algebraic dynamics.

\subsection{Algebraic dynamics}
\kkk
Let $S$ be a quasi-projective variety over $k$ with function field $F$. 
Let $(X,f,L)$ be a \emph{polarized algebraic dynamical system} over $S$, i.e. $X$ is a flat and projective integral scheme over $S$, 
$f:X\to X$ is an endomorphism over $S$,
and $L$ is an $f$-ample $\QQ$-line bundle  satisfying $f^*L \simeq qL$ for some rational number $q>1$. 

By Tate's limiting argument, we can construct a canonical adelic $\QQ$-line bundle $\overline L_f\in \wh\Pic(X)_{\QQ,\nef}$ extending $L$ which is $f$-invariant in that $f^*\overline L_f \simeq q\overline L_f$. 
Here $\wh\Pic(X)_{\QQ,\nef}$ denotes the sub-semigroup of 
$\wh\Pic(X)_{ \QQ}$ consisting of positive rational multiples of elements of 
$\wh\Pic(X)_{\nef}$.

For any closed $\overline F$-subvariety $Z$ of $X_F$, we have the canonical height
$$
\fh_{f}(Z)
=\fh_{\overline L_f}(Z):
= \frac{ \left\langle \overline L_f|_{Z'}^{\dim Z+1}\right\rangle }
{(\dim Z+1)\deg_{L}(Z' _F) }\ \in \wh \Pic(F)_{\QQ,\nef}.$$
In particular, we have a height function
$$
\fh_f: X(\overline F) \lra \wh \Pic(F)_{\QQ,\nef}.
$$
Tate's limiting argument also explains these heights. 

The height function $\fh_f$ is $f$-invariant. 
As a consequence, $\fh_f(x)=0$ for a preperiodic point $x\in X(\overline F)$.
In the minimal case that $K$ is a number field or a function field of one variable over a finite field $k$, 
$\fh_f$ satisfies the Northcott property. In this case, 
$\fh_f(x)=0$ for a point $x\in X(\overline F)$ 
implies that $x$ is preperiodic under $f$.

\subsection{Equidistribution of small points}

Our equidistribution conjecture naturally implies an equidistribution conjecture of preperiodic points. 
\begin{conj}[Conjecture \ref{equi2}]
\kkk
Let $F$ be a finitely generated field over $k$. Let $(X, f, L)$ be a polarized dynamical system over $F$.
Let $v$ be a non-trivial valuation of $F$.
Assume that the restriction of $v$ to $K$ is trivial if $k$ is a  field.
Let $\{x_m\}_m$ be a generic sequence of preperiodic points in $X(\overline F)$. 
Then the Galois orbit of $\{x_m\}_m$ is equidistributed in $X_{F_v}^\an$ for the  canonical measure $d\mu_{L, f,v}$.
\end{conj}

As an example of our equidistribution theorem (cf. Theorem \ref{equi3}), we deduce the following equidistribution theorem of small points on non-degenerate subvarieties in 
a family of polarized algebraic dynamical systems.

\begin{thm}[Theorem \ref{equi4}]
Let $S$ be a quasi-projective variety over a number field $K$. 
Let $(X,f, L)$ be a polarized dynamical system over $S$.
Let $Y$ be a non-degenerate closed subvariety of $X$ over $K$. 
Let $\{y_m\}_{m\geq 1}$ be a generic sequence of $Y(\overline K)$ such that 
$h_{\OL_{f}}(y_m) \to 0$.
Then for any place $v$ of $K$,
the Galois orbit of  $\{y_m\}_{m\geq 1}$ is equidistributed on the analytic space $Y_v^\an$
for the  canonical measure $d\mu_{\OL_f|_Y,v}$. 
\end{thm}

In the theorem, a closed subvariety $Y$ of $X$ is called \emph{non-degenerate} 
if $\deg _{\wt L}(Y)>0$. This is equivalent to the property that $\TL|_{Y}$ is \emph{big}.
If $X$ is an abelian scheme and $K$ is a number field, our definition of 
``non-degenerate'' agrees with that of \cite{DGH}, which uses Betti maps in the complex analytic setting.

The theorem generalizes the equidistribution theorem of 
DeMarco--Mavraki \cite{DMM} for families of elliptic curves, and confirms the conjecture (REC) of K\"uhne \cite{Kuh} for abelian schemes.
A weaker version of the theorem for abelian schemes is proved by \cite[Thm. 1]{Kuh} independently and applied to prove a uniform Bogomolov conjecture after the work of Dimitrov--Gao--Habegger \cite{DGH}, as mentioned above. 
The proof of \cite{Kuh} is a limit version of the original proof in \cite{SUZ} and uses a result of Dimitrov--Gao--Habegger \cite{DGH} for uniformity in the limit process.
Inspired by our formulation, Gauthier \cite{Gau} has extended the equidistribution theorem of \cite{Kuh} to more general settings, which has a large overlap with our equidistribution theorem.

\subsection{Heights of points of a non-degenerate subvariety}
\kkk
Let $K$ be a number field if $k=\BZ$ or a function field of one variable if $k$ is a field.
Let $S$ be a quasi-projective variety over $K$. 
Let $(X,f, L)$ be a polarized dynamical system over $S$.
Let $Y$ be a closed subvariety of $X$ over $K$.

Suppose $Y$ is a section of $X\to S$. In that case, our vector-valued height of adelic line bundles generalizes and re-interprets the Tate--Silverman specialization theorem of \cite{Tat, Sil2, Sil3, Sil4}, and the work \cite{DMM} from families of elliptic curves to families of algebraic dynamical systems. See Lemma \ref{specialization} for more details. 

As mentioned above, if $X$ is an abelian scheme and $Y$ is non-degenerate in $X$, there is a height inequality of points of $Y$ by \cite{GH, DGH}, which plays a fundamental role in the treatment of the uniform Mordell conjecture in \cite{DGH, Kuh}. 
In terms of our theory, we have a simple interpretation of the height inequality, which is also valid in families of algebraic dynamical systems. 
As the non-degeneracy is interpreted as the bigness of $\TL|_{Y}$, the height inequality is also interpreted by the bigness of some adelic line bundle.
Applying Theorem \ref{height comparison intro0}(2) to the morphism $Y\to S$ and the adelic line bundle $\OL_f|_Y$ on $Y$, we can have a lower bound of the canonical height of points on $Y$ by Weil heights on $S$. 
See Theorem \ref{height comparison} for more details.

\subsection{Equidistribution of PCF maps}

Let $S$ be a smooth and quasi-projective variety over a number field $K$. Let $X=\BP^1_S$ be the projective line over $S$, and let $f: X\to X$ be a finite morphism over $S$ of degree $d>1$. 
A point $y\in S(\overline K)$ is called \emph{post-critically finite} (PCF) if all the ramification points (i.e. critical points) of $f_y:X_y\to X_y$ are preperiodic under $f_y$.

Denote by $\CM_d$ the moduli space over $K$ of endomorphisms of $\PP^1$ of degree $d$. 
Inside $\CM_d$, there is a closed subvariety corresponding to flexible Latt\'es maps.
By the moduli property, there is a morphism $S\to\CM_d$. 

The main result here is the following equidistribution theorem of Galois orbits of PCF points. 

\begin{thm} [Theorem \ref{equi PCF1}]
Assume that the morphism $S\to \CM_d$ is generically finite and its image is not contained in the flexible Latt\`es locus. 
Let $\{y_m\}_m$ be a generic sequence of PCF points of $S(\overline K)$.
Then the Galois orbit of $\{y_m\}_m$ is equidistributed in $S_{K_v}^\an$ for $d\mu_{\overline M,v}$ for any place $v$ of $K$.
\end{thm}

If $S$ is a family of \emph{polynomial} maps on $\PP^1$, the theorem was previously proved by Favre--Gauthier \cite{FG}.
Their strategy was to reduce the problem to the equidistribution of Yuan \cite{Yua1}.

Now  we explain our proof of the theorem, which will also introduce the key term $\OM$ in the statement.
Denote by $R$ the ramification divisor of the finite morphism $f: X\to X$, viewed as a (possibly non-reduced) closed subscheme in $X$. Then $R$ is finite and flat of degree $2d-2$ over $S$, and the fiber $R_y$ of $R$ above any point $y\in S$ is the ramification divisor of $f_y: X_y\to X_y$.

Let $L$ be a $\QQ$-line bundle on $X$, of degree one on fibers of $X\to S$, such that $f^*L\simeq dL$. Denote by $\OL=\OL_f$ the $f$-invariant extension of $L$ in $\wh\Pic(X)_{\QQ,\nef}$ such that $f^*\OL\simeq d\OL$. 
Denote  
$$
\OM:= N_{R/S}(\OL|_R) \in \wh\Pic(S)_{\QQ,\nef}.
$$
Here the norm map is the Deligne pairing of relative dimension 0. 

Consider the height function
$$
h_\OM:S(\overline K)\lra \RR_{\geq0}.
$$
For any $y\in S(\overline K)$, the height $h_\OM(y)=0$ if and only if $y$ is PCF in $S$. 
Then we are in the situation to apply the previous equidistribution theorem (Theorem \ref{equi3}) to $(S,\OM)$, except that we need to check the condition 
$\deg_{\TM}(S)>0$. 

This requires the bifurcation measure introduced by DeMarco \cite{DeM1, DeM2} and further studied by Bassanelli--Berteloot \cite{BB07}.
The $\deg_{\TM}(S)$ is exactly equal to the total volume of the bifurcation measure on $S_\sigma(\CC)$ for any embedding $\sigma: K\to \CC$. Then $\deg_{\TM}(S)>0$ is eventually equivalent to the condition on $S\to\CM_g$ by the works of \cite{BB07, GOV}.
This proves the theorem and confirms that the equilibrium measure $d\mu_{\overline M,\sigma}$ is a constant multiple of the bifurcation measure
for any embedding $\sigma:K\to \CC$.

Recently, Ji--Xie \cite{JX} proved the dynamical Andre--Oort conjecture for 1-dimensional families, which relies on our theory of adelic line bundles and especially our equidistribution theorem of PCF points.

In the end, we note that the theorem also holds for a family of morphisms on $\PP^n$ with a slightly weaker statement. In particular, the construction of the adelic line bundle $\OM$ works in the same way. 
We refer to Theorem \ref{equi PCF0} for more details.

\subsection{Hodge index theorem on curves}

In the end, we present our generalization of the arithmetic Hodge index theorem of Faltings \cite{Fal1} and Hriljac \cite{Hri} to finitely generated fields. 
We refer to Theorem \ref{thm NT} for a detailed account. 

\kkk
Let $F$ be a finitely generated field over $k$, and let $\pi: X\to \Spec F$ be a smooth,  projective, and geometrically connected curve of genus $g>0$. 
Denote by 
$J=\underline{\Pic}^0_{X/F}$ the Jacobian variety of $X$ and by $\Theta$ 
the symmetric theta divisor on $J$. 
By the dynamical system $(J, [2], \Theta)$, we have a N\'eron--Tate height function 
$$
\hat\fh: \Pic^0(X_{\overline F}) \lra \wh\Pic(F/k)_{\QQ,\nef}.
$$
The height function is quadratic, as in the classical case.

\begin{thm}[Theorem \ref {thm NT}]
\kkk
Let $F$ be a finitely generated field over $k$, and let $\pi: X\to \Spec F$ be a smooth,  projective, and geometrically connected curve.
Let $M$ be a line bundle on $X$ with $\deg M=0$. 
Then there is an adelic line bundle $\overline M_0\in \wh\Pic(X/k)_{\intb,\QQ}$ with underlying line bundle $M$ such that 
$$\pi_*\pair{\overline M_0,\overline V}= 0,\quad\forall\,\overline V\in \wh\Pic(X/k)_{\vert,\QQ}.$$ 
Moreover, for such an adelic line bundle,
$$\pi_*\pair{\overline M_0, \overline M_0}= -2\, \wh\fh(M).$$
\end{thm}

In the theorem, $\wh\Pic(X/k)_{\vert,\QQ}$ is the space of vertical adelic line bundles defined as the kernel of the forgetful map $\wh\Pic(X/k)_{\intb,\QQ}\to \Pic(X)_\QQ$.

\section{Notation and terminology} \label{sec uniform}
We will introduce a uniform system of terminology and notations for both the arithmetic and geometric cases. To achieve this, we need to abuse terminology frequently.

Our base ring $k$ is either $\ZZ$ or an arbitrary field. This is divided into two cases:
\begin{enumerate}
\item (arithmetic case) $k=\ZZ$. In this case, the adelic line bundles will be limits of hermitian line bundles on projective integral schemes over $\ZZ$. This limit process obtains the intersection theory.
\item (geometric case) $k$ is an arbitrary field of arbitrary characteristic. In this case, the adelic line bundles will be the limit of usual line bundles on projective varieties of over $k$. This limit process obtains the intersection theory.
\end{enumerate}

By a \emph{finitely generated field $F$ over $k$}, we mean a field $F$ which is finitely generated over the fraction field of $k$.
For any integral scheme $X$ over $k$, denote by $k(X)$ the \emph{function field} of $X$. 

By a \emph{projective variety} over $k$, we mean  an integral,
 flat, and projective scheme over $k$.
 By a \emph{quasi-projective variety}  over $k$, we mean an open subscheme of projective variety over $k$.
 For a quasi-projective variety $\CU$ over $k$, a \emph{projective model} means a projective  variety $\CX$ over $k$ endowed with an open immersion $\CU\to \CX$ over $k$. 
In the arithmetic case (that $k=\ZZ$), we may also use the terms \emph{quasi-projective arithmetic variety} and \emph{projective arithmetic variety} to emphasize the situation.

In the arithmetic case, for a {projective arithmetic variety} $\CX$ over $\ZZ$, we have the group $\wh\Div(\CX)$ of arithmetic divisors on $\CX$, and the group $\wh\Pic(\CX)$ and the category $\wh\Picc(\CX)$ of hermitian line bundles on $\CX$.

In the geometric case, for a {projective variety} $\CX$ over a field $k$, 
an arithmetic divisor means a Cartier divisor, a hermitian line bundle means a line bundle, 
and we write 
$\wh\Div(\CX)$, $\wh\Pic(\CX)$, $\wh\Picc(\CX)$ 
for 
$\Div(\CX)$, $\Pic(\CX)$, $\Picc(\CX)$.
We take this convention in other similar situations. 

This abuse of notation is only one-way. For example, by $\Div$, $\Pic$, or $\Picc$ in the arithmetic case, we still mean the ones without the Archimedean components. 

Below are a few conventions that are not directly related to the base $k$ but are taken throughout this book. 
\begin{enumerate}
\item Denote $M_\QQ=M\otimes_\ZZ\QQ$ for any abelian group $M$. Take similar conventions for $M_\RR$ and $M_\CC$. 
\item For any field $K$, we fix an algebraic closure $\OK$ of $K$ throughout this book. 
\item Except in \S\ref{sec adelic general} and \S\ref{sec local theory}, all schemes are assumed to be noetherian. 
\item By a \emph{variety} over a field, we mean an integral scheme, separated and of finite type over the field. We do not require it to be geometrically integral. 
\item By a \emph{curve} over a field, we mean a variety over the field of dimension 1. 
\item All \emph{divisors} in this book are \emph{Cartier divisors}, unless otherwise instructed.
\item By a \emph{line bundle} on a scheme, we mean an invertible sheaf on the scheme. 
We often write or mention tensor products of line bundles additively, so $a\CL-b\CM$ means
$\CL^{\otimes a}\otimes \CM^{\otimes (-b)}$.
\item All the categories of (adelic, metrized) line bundles are groupoids, so the morphisms are isomorphisms.
\item A \emph{functor} between two categories may also be called a map or a homomorphism sometimes. 
\end{enumerate}

\chapter{Adelic divisors and adelic line bundles} \label{sec adelic}

In this chapter, we develop a theory of adelic divisors and adelic line bundles on essentially quasi-projective schemes. 
The main ideas of this chapter are explained in the introduction.

\section{Preliminaries on arithmetic varieties} \label{sec hermitian}

In this section, we review some basic notions of arithmetic divisors and hermitian line bundles on projective arithmetic varieties. 
These are standard terminology, and most of them are reviewed in the appendix of this book in slightly different settings.
We also refer to the textbook of Yuan--Guo \cite{YG25} for a first course on Arakelov geometry.

\subsection{Metrics on complex analytic spaces} \label{sec complex}

We refer to  \cite[Chapter II]{Dem2} or \cite{Rem} for detailed introductions to complex analytic varieties. For convenience, all complex analytic varieties in this book are assumed to be reduced and irreducible. 

Let $X$ be a (reduced and irreducible) complex analytic variety. 
The default topology on $X$ is the complex topology unless otherwise instructed. 
The \emph{regular locus} (or equivalently \emph{smooth locus}) $X^\reg$ is a complex manifold, which is open and dense in $X$.
In the following, we introduce metrics of line bundles on $X$ with different types of regularities. 

We take the notion of smooth differential forms following \cite[\S1.1]{Kin}.
For integers $p,q\geq 0$, a \emph{smooth $(p,q)$-form on} $X$ is a smooth $(p,q)$-form $\alpha$ on $X^\reg$ such that for any point $x\in X$, there is an open neighborhood $U$ of $x$ in $X$ and an analytic map $i:U\to M$ to a complex manifold $M$ under which $U$ is a closed analytic subvariety of $M$, such that $\alpha|_{U^\reg}$ can be extended to a smooth $(p,q)$-form on $M$.

The case $p=q=0$ gives the notion of smooth functions. 
Namely, a \emph{smooth function} on $X$ is a continuous function $f: X\to \CC$ such that for any point $x\in X$, there is an open neighborhood $U$ of $x$ in $X$ and an analytic map $i: U\to M$ to a complex manifold $M$ under which $U$ is a closed analytic subvariety of $M$, such that $f|_{U}$ can be extended to an infinitely differentiable function $\tilde f: M\to \CC$.
Note that the smoothness here is stronger than that in \cite{Zha1}. 

With the definition of smooth $(p,q)$-forms, most notions and operations on differential forms and currents on complex manifolds can be extended to complex analytic varieties. 

Let $L$ be a line bundle on $X$. 
By a \emph{continuous metric} (resp. \emph{smooth metric}) of $L$ on $X$, we mean the assignment of a metric $\|\cdot\|$ to the fiber $L(x)$ above every point $x\in X$, which varies continuously (resp. smoothly) in that for any local analytic section $s$ of $L$ defined on an open subset $U$ of $X$, the function $\|s(x)\|^2$ is continuous (resp. smooth) in $x\in U$.

For any continuous metric $\|\cdot\|$ of $L$ on $X$, 
the \emph{Chern current} 
$$
c_1(L,\|\cdot\|):= \frac{1}{\pi i}\partial\overline\partial \log\|s\|+ \delta_{\div(s)}
$$
is a $(1,1)$-current on $X$. Here $s$ is any meromorphic section of $L$ on $X$, and the definition is independent of the choice of $s$.

A continuous metric $\|\cdot\|$ of $L$ on $X$ is called \emph{semipositive} if the \emph{Chern current} is a positive current. Equivalently, for any analytic curve $Y$ of $X$,  and any smooth and compactly supported function $f$, the integration 
$$\int _Y f c_1(L|_Y, \|\cdot\|)\ge 0.$$
As the well-known special case, if $X$ is smooth and the metric $\|\cdot\|$ of $L$ is smooth, the Chern current is represented by the \emph{Chern form} $c_1(L,\|\cdot\|)$ (by abuse of notation).
In this case, $\|\cdot\|$ is semipositive on $X$ if and only if $c_1(L,\|\cdot\|)$ is positive semi-definite as a smooth $(1,1)$-form on $X$.
If $X$ is general and the metric $\|\cdot\|$ is smooth, then $\|\cdot\|$ is semipositive on $X$ if and only if $c_1(L|_{X^\reg},\|\cdot\|)$ is positive semi-definite as a smooth $(1,1)$-form on $X^\reg$.

A continuous metric $\|\cdot\|$ of $L$ on $X$ is called \emph{integrable} if it is the quotient of two semipositive metrics; i.e. there are line bundles $(L_1, \|\cdot\|_1)$ and $(L_2, \|\cdot\|_2)$ endowed with semipositive metrics on $X$ such that $(L, \|\cdot\|)$ is isometric to $(L_1, \|\cdot\|_1)\otimes (L_2, \|\cdot\|_2)^\vee$.

Let $X$ be a complex projective variety. Let $L$ be a line bundle on $X$. 
Then any smooth metric $\|\cdot\|$ of $L$ on $X$ is the quotient of two semipositive metrics of line bundles on $X$.
In fact, if $X$ is smooth, take an ample line bundle $A$ with a positive metric $\|\cdot\|_A$, then $(L,\|\cdot\|)\otimes (A,\|\cdot\|_A)^{\otimes m}$ also have a positive metric for sufficiently large $m$.
If $X$ is singular, take a closed embedding $i:X\to \BP^N$ and set 
$(A,\|\cdot\|_A)=i^*(\CO(1),\|\cdot\|_\mathrm{FS})$, where $\|\cdot\|_\mathrm{FS}$ is the Fubini--Study metric. 
By a local argument, 
$(L,\|\cdot\|)\otimes (A,\|\cdot\|_A)^{\otimes m}$ still has a semipositive metric for sufficiently large $m$.
As a consequence, smooth metrics are integrable.

Note that \cite{Zha2} also has a notion of ``semipositive metrics'' and ``integrable metrics''. We will see that our notion is essentially equivalent to those of the loc. cit. on complex projective varieties.
Let $L$ be a line bundle on a complex projective variety $X$ with a continuous metric $\|\cdot\|$. Then we have the following comparisons:
\begin{enumerate}
\item[(1)] If $L$ is ample, then $\|\cdot\|$ is semipositive (in our sense) if and only if it is semipositive in the sense of \cite{Zha2}.
In fact, by \cite[Cor. C]{CGZ}, any semipositive continuous metric on $L$ is an increasing limit of semipositive smooth metrics on $L$. This limit process is uniform since all the metrics are continuous. 
Then it is semipositive in the sense of \cite{Zha2}. The inverse direction follows from the fact that the decreasing limit of psh functions is again psh. 

\item[(2)] Any semipositive metric $\|\cdot\|$ (without assuming that $L$ is ample) is the quotient of two semipositive continuous metrics on ample line bundles. 
This is trivial by tensoring  $(L,\|\cdot\|)$ by an ample line bundle with a semipositive continuous metric.

\item[(3)] A continuous metric $\|\cdot\|$ is integrable if and only if it is integrable in the sense of \cite{Zha2}. This follows from (1) and (2).
\end{enumerate}
 
As above, let $X$ be a complex projective variety. Denote $n=\dim X$. 
Let  
$(L_1,\|\cdot\|_1),\cdots, (L_n,\|\cdot\|_n)$ be line bundles with integrable metrics on $X$. 
Then there is a Monge--Amp\`ere measure $c_1(L_1,\|\cdot\|_1)\cdots c_1(L_n,\|\cdot\|_n)$ on $X$. 
It is reduced to semipositive metrics by linearity, and then the approximation method of \cite[Thm. 2.1]{BT} (or \cite[Cor. 1.6]{Dem1}).

\subsection{Green functions on complex analytic spaces} \label{sec Green}

Let $X$ be a complex analytic variety. We introduce Green functions of divisors following the above treatment of metrics of line bundles.

Let $D$ be an (analytic) Cartier divisor on $X$ with support $|D|$. 
A \emph{Green function} (resp. \emph{Green function of smooth type}) of $D$ on $X$ is a function 
$g: X\setminus |D| \to \RR$ such that for any meromorphic function $f$ on an open subset $U$ of $X$ satisfying $\div(f)=D|_U$, the function $g+\log |f|$ can be extended to a continuous (resp. smooth) function on $U$.
Sometimes, a Green function is also called \emph{a Green function of continuous type} to emphasize the property. 

Note that the pair $(D,g)$ defines a pair $(\CO(D),\|\cdot\|_g)$ with the metric defined by $\|s_D\|_g=e^{-g}$. Here $s_D$ is the section of $\CO(D)$ corresponding to the meromorphic function $1$ on $X$.

By this correspondence, $g$ is of continuous type (resp. of smooth type) if and only if $\|\cdot\|_g$ is  
continuous (resp. smooth). 
Moreover, we say that $g$ is \emph{semipositive} (resp. \emph{integrable}) if 
$\|\cdot\|_g$ is semipositive (resp. integrable). 

All the definitions and results for metrics and Green functions easily extend to finite disjoint unions of the analytic variety (i.e., the complex projective variety).

\subsection{Hermitian line bundles on arithmetic varieties}

By a \emph{projective arithmetic variety (resp. quasi-projective arithmetic variety)} $\CX$, we mean an integral scheme, projective (resp. quasi-projective) and flat over $\BZ$. 
We usually denote by $\QQ(\CX)$ the function field of $\CX$.

Let $\CX$ be a projective arithmetic variety. 
A  {\em hermitian line bundle on $\CX$} is a pair $\overline \CL=(\CL, \|\cdot\|)$,  
where  $\CL$ is a line bundle (equivalently an invertible sheaf) on $\CX$, and $\|\cdot\|$ is a \emph{continuous metric} of $\CL(\CC)$ on $\CX(\CC)$, invariant under the action of the complex conjugate. 

An \emph{isometry} from a hermitian line bundle  $\overline \CL=(\CL, \|\cdot\|)$ to another hermitian line bundle $\overline \CL'=(\CL', \|\cdot\|')$ is an isomorphism $\CL\to \CL'$ of coherent sheaves compatible with the metrics. 

Denote by $\wh\Pic(\CX)$ the \emph{group of isometry classes} of hermitian line bundles on $\CX$. 
Denote by $\wh\Picc(\CX)$ the \emph{category} of hermitian line bundles on $\CX$, in which the morphisms are isometries of hermitian line bundles. This is a groupoid by definition. 

Note that for a hermitian line bundle $\overline \CL=(\CL, \|\cdot\|)$, we only require the metric to be continuous (instead of smooth). 
This relaxed notion will bring some convenience in approximation later and will also make 
the category of hermitian line bundles on $\CX$ equivalent to the category of adelic line bundles on $\CX/\ZZ$.

For convenience, define $\wh\Picc(\CX)_{\sm}$ (resp. $\wh\Picc(\CX)_{\intb}$)
to be the full subcategory of $\wh\Picc(\CX)$ of hermitian line bundles with smooth metric (resp. integrable metrics). 
Define $\wh\Pic(\CX)_{\sm}$ (resp. $\wh\Pic(\CX)_{\intb}$)
to be the subgroup of $\wh\Pic(\CX)$ similarly.

\subsection{Arithmetic divisors}

Let $\CX$ be a projective arithmetic variety. 
An  {\em arithmetic divisor on $\CX$} is a pair $\OCD=(\CD,g_\CD)$, 
where  $\CD$ is a Cartier divisor on $\CX$, and $g_\CD$ is a \emph{Green function} of $\CD(\CC)$ on $\CX(\CC)$, invariant under the action of the complex conjugate. 
A {\em principal arithmetic divisor on $\CX$} is an arithmetic divisor of the form 
$$\wh\div(f):=(\div(f),-\log|f|)$$ 
for any rational function $f\in \QQ(\CX)^\times$ on $\CX$.

Denote by $\wh\Div(\CX)$ the group of arithmetic divisors on $\CX$, and by 
$\wh\Pr(\CX)$ the group of principal arithmetic divisors on $\CX$. 
Then we have the \emph{arithmetic divisor class group} 
$$
\wh\CaCl(\CX)=\wh\Div(\CX)/\wh\Pr(\CX).
$$

An arithmetic divisor $\OCD=(\CD,g_\CD) \in \Divhat(\CX)$ is \emph{effective}
(resp. \emph{strictly effective}) if
$\CD$ is an effective Cartier divisor on $\CX$ and the Green function $g_\CD\ge 0$
(resp. $g_\CD>0$) on $\CX(\BC)-|\CD(\BC)|$.

There is a canonical map
$$
\wh\Div(\CX)\lra \wh\Pic(\CX), \quad\ \OCD\longmapsto \CO(\OCD),
$$
which induces an isomorphism
$$
\wh\CaCl(\CX)\iso \wh\Pic(\CX).
$$
The inverse image of a hermitian line bundle $\CLL$ is represented by the divisor 
$$
\wh\div(s)=\wh\div_{(\CX, \CLL)}(s):=(\div(s),-\log\|s\|),
$$
where $s$ is any nonzero rational section of $\CL$ on $\CX$.

Similar to hermitian metrics, we only require the Green functions to be of continuous type (instead of smooth type). 
Define $\wh\Div(\CX)_{\sm}$ (resp. $\wh\Div(\CX)_{\intb}$)
to be the subgroup of $\wh\Div(\CX)$ of arithmetic divisors with Green functions of smooth types (resp. integrable Green functions).

Assume now $\CX$ be a projective arithmetic variety.
An arithmetic divisor $\OD$ on $\CX$ is \emph{nef} if the hermitian line bundle $\CO(\OD)$ is nef on $\CX$.
Denote by $\wh\Div(\CX)_\nef$ the sub-semigroup of $\wh\Div(\CX)$ of nef line bundles on $\CX$.

\subsection{Arithmetic intersection numbers}

Let $\CX$ be a projective arithmetic variety as above.  
Let $\CZ$ be a closed integral subscheme of $\CX$ of dimension $d\ge 0$. 
We say that $\CZ$ is horizontal if $\CZ\to \Spec \ZZ$ is surjective; we say that $\CZ$ is vertical if the image of $\CZ\to \Spec \ZZ$ is a closed point.  
Let $\ol \CL_1, \cdots, \ol \CL_d$ be $d$ hermitian line bundles on $\CX$ with integrable metrics. 
Then we define the intersection number $\ol \CL_{1}\cdot \ol \CL_{2}\cdots \ol \CL_{d}\cdot[\CZ]$ by induction on $d$ as follows. 

If $d=0$, then $\CZ$ is a closed point of $\CX$  and thus $\Gamma(\CZ,\CO_\CZ)$ is a finite field. 
Define $[\CZ]=\log \#\Gamma(\CZ,\CO_\CZ)$. 

If $d>0$, we take a nonzero rational section $s_d$ of $\CL_d$ over $\CZ$. 
In terms of Weil divisors, we write $\div (s_d)=\sum _i a_i \CZ_i$ with $a_i\in \BZ$ and $\CZ_i$ integral subschemes of $\CZ$ of dimension $d-1$.
We define 
$$
\ol \CL_{1}\cdots \ol \CL_{d}\cdot[\CZ]=\sum _i a_i\,  \ol \CL_{1}\cdots  \ol \CL_{d-1}\cdot[\CZ_i]
- \int_{\CZ(\BC)}\log \|s_d\| c_1(\ol \CL_1)\cdots c_1(\ol \CL_{d-1}).
$$
Here we take the convention that the integral on the right-hand side is zero if $\CZ(\CC)=\emptyset$, which happens if $\CZ$ is vertical. 
The intersection number is the independence of the choices of $s_d$ and the orderings of $\ol \CL_1,\cdots, \ol \CL_d$.

There is also a similar induction formula for intersection numbers of arithmetic divisors with integrable Green functions. 

If $f:\CX'\to \CX$ is a morphism from another projective arithmetic variety $\CX'$, and $\CZ'$ is a closed integral subscheme of $\CX'$ of dimension $d$ with $f(\CZ')\subset \CZ$, then we have the projection formula 
$$f^*\ol \CL_{1}\cdot f^*\ol \CL_{2}\cdots f^*\ol \CL_{d}\cdot[\CZ']
=\deg(\CZ'/\CZ)\cdot (\ol \CL_{1}\cdot \ol \CL_{2}\cdots \ol \CL_{d}\cdot[\CZ]).$$
Here $\deg(\CZ'/\CZ)$ is the degree between the function fields of these two schemes, which is understood to be 0 if $\CZ'\to \CZ$ is not surjective.

Now we consider a few special cases of intersection numbers. 

If $d\leq 1$ or if $\CZ$ is vertical, the definition is also valid without assuming the hermitian metrics to be integrable. 

If $\CZ$ is vertical, then it is a projective variety over a finite field 
$\FF_p$, and thus we have 
$$\ol \CL_{1}\cdots \ol \CL_{d}\cdot[\CZ]=((\CL_{1}|_{\CZ})\cdot ( \CL_{2}|_{\CZ})\cdots ( \CL_{d}|_{\CZ}))\log p.$$
Here, the intersection number on the right-hand side is the usual intersection number for projective varieties in algebraic geometry, and the hermitian metrics do not play any role here.

If $\CZ=\CX$, we simply write $\ol \CL_{1}\cdots \ol \CL_{d}$ for $\ol \CL_{1}\cdots \ol \CL_{d}\cdot[\CX]$. 
This gives a symmetric and multi-linear map 
$$
\wh\Pic(\CX)^{\dim\CX}\lra \RR. 
$$
By this notation, for general $\CZ$ in $\CX$, we easily have 
$$\ol \CL_{1}\cdot \ol \CL_{2}\cdots \ol \CL_{d}\cdot[\CZ]
=(\ol \CL_{1}|_{\CZ})\cdot (\ol \CL_{2}|_{\CZ})\cdots (\ol \CL_{d}|_{\CZ}).$$

In the case of arithmetic curves that $\CZ=\CX=\Spec O_K$ for a number field $K$,
the intersection number is just the arithmetic degree
$$
\wh\deg(\CLL_1)= \log \#(\CL_1/s_1 O_K)-\sum_{\sigma:K\to \CC}\log \|s_1\|_\sigma,
$$ 
where $s\in \CL_1$ is a nonzero element.

\section{Objects of mixed coefficients} \label{sec convention}

The goal of this section is to introduce notations for divisors and line bundles of mixed coefficients, i.e.
 $\QQ$-line bundles and $\QQ$-divisors 
which are integral to an open subscheme of the ambient scheme. 
These are less standard but will be crucial to define effective sections of adelic line bundles in our theory.

For clarity, in this section, we do not take the uniform terminology in \S\ref{sec uniform} but introduce all the terms case by case.

\subsection{$\QQ$-divisors and $\QQ$-line bundles}

When we say {divisors}, we always mean Cartier divisors unless otherwise specified. 
When we want to distinguish the usual divisors (resp. line bundles) from $\QQ$-divisors (resp. $\QQ$-line bundles), we often say integral divisors (resp. integral line bundles).

Let $\CX$ be a scheme. 
Denote by $\Div(\CX)=H^0(\CX,\CK_X^\times/\CO_X^\times)$
the group of \emph{Cartier divisors} on $\CX$.
Here $\CK_X$ is the sheaf of rational functions on $X$. 
The image of $H^0(\CX,\CK_X^\times)$ in $\Div(\CX)$ is the subgroup of \emph{principle Cartier divisors} on $\CX$, denoted by $\Pr(\CX)$.

The \emph{support} $|\CD|$ of a Cartier divisor $\CD$ on $\CX$ is the complement of the maximal open subscheme of $\CX$ on which $\CD$ is trivial.   
A Cartier divisor $\CD$ on $\CX$ is called \emph{effective} if 
it lies in the image of the semi-group $H^0(\CX,\CO_X/\CO_X^\times)$ in $\Div(\CX)$.

An element of $\Div(\CX)_\QQ=\Div(\CX)\otimes_\ZZ\QQ$ is called a (Cartier)
\emph{$\QQ$-divisor} of $\CX$.
A $\BQ$-divisor $\CD \in \Div(\CX)_\BQ$ is called \emph{effective}  if for some
positive integer $m$, the multiple $m \CD$ is an effective (integral) divisor in $\Div(\CX)$.

Denote by $\Picc(\CX)$ the \emph{category} of line bundles on $\CX$, in which the objects are line bundles (or equivalently invertible sheaves) on $\CX$, and the morphisms are isomorphisms of line bundles.  
Denote by $\Picc(\CX)_\QQ$ the \emph{category} of $\QQ$-line bundles on $\CX$, in which the objects are pairs $(a,\CL)$ (or just written as $a\CL$)
with $a\in \QQ$ and $\CL\in\Picc(\CX)$, and a morphism of two such objects is defined to be
$$\Hom(a\CL,a'\CL'):=\varinjlim_m \Hom(am\CL, a'm\CL'),$$
where $m$ runs through positive integers such that $am$ and $a'm$ are both integers, so that $am\CL$ and $a'm\CL'$ are viewed as integral line bundles, 
and ``$\Hom$'' on the right-hand side represents isomorphisms of integral line bundles. 
For the direct system, for any $m|n$, there is a transition map 
$$\Hom(am\CL, a'm\CL')\lra  \Hom(an\CL, a'n\CL')$$ 
locally given by taking $(n/m)$-th power of an isomorphism. 
The group of isomorphism classes of objects of $\Picc(\CX)_\QQ$
 is isomorphic to $\Pic(\CX)_\QQ=\Pic(\CX)\otimes_\ZZ\QQ$.

Let $a\CL$ be a $\QQ$-line bundle on $\CX$ with $a\in\QQ$ and $\CL\in\Picc(\CX)$. 
A \emph{section} of $a\CL$ on $\CX$ is an element of 
$\ds\Hom(\CO_X,a\CL)=\varinjlim_m \Gamma(\CX,am\CL)$, where 
$m$ runs through positive integers with $am\in \ZZ$.
If $\CX$ is an integral scheme, a \emph{rational section} of $a\CL$ on $\CX$ is an element of $\ds\Hom(\CO_\eta,a\CL_\eta)=\varinjlim_m \Gamma(\eta,am\CL)$, where $\eta$ is the generic point of $\CX$, and
$m$ runs through positive integers with $am\in \ZZ$. 
If $s$ is a section represented by $s_m\in \Gamma(\CX,am\CL)$ or 
a rational section represented by $s_m\in \Gamma(\eta,am\CL)$, then define
$$\div(s):=\frac{1}{m}\div(s_m).$$
This is a $\QQ$-divisor on $\CX$.

If $\CX$ is a projective variety over a field, a $\BQ$-divisor $\CD \in \Div(\CX)_\BQ$ (resp. $\BQ$-line bundle $\CL\in \Picc(\CX)_\QQ$) is called \emph{nef} if for some
positive integer $m$, the multiple $m \CD$ (resp. $m\CL$) is a nef divisor on $\CX$ (resp. nef line bundle on $\CX$) in the usual sense.

\subsection{Arithmetic $\QQ$-divisors and hermitian $\QQ$-line bundles}

The above $\QQ$-notions extend easily to the arithmetic situation.  We sketch them briefly. 

Let $\CX$ be a projective arithmetic variety. 
An element of $\wh\Div(\CX)_\QQ=\wh\Div(\CX)\otimes_\ZZ\QQ$ is called an 
\emph{arithmetic $\QQ$-divisor}.
So an arithmetic $\QQ$-divisor is still represented by a pair $(\CD,g_\CD)$, where $\CD$ is a Cartier $\QQ$-divisor on $\CX$, and $g_\CD$ is a Green function of $\CD$ on $\CX$ defined similarly. 

An arithmetic $\BQ$-divisor $\OCD \in \Divhat(\CX)_\BQ$ is called \emph{effective}
(resp. \emph{strictly effective}) if for some
positive integer $m$, the multiple $m\OCD$ is an effective (resp. strictly effective) (integral) arithmetic divisor in $\Divhat(\CX)$.

Denote by $\wh\Picc(\CX)_\QQ$ the \emph{category} of hermitian $\QQ$-line bundles on $\CX$, in which the objects are pairs $(a,\overline\CL)$ (or just written as $a\overline\CL$)
with $a\in \QQ$ and $\overline\CL\in\wh\Picc(\CX)$, and the morphism of two such objects is defined to be
$$\Hom(a\CLL,a'\CLL'):=\varinjlim_m \Isom(am\CLL, a'm\CLL'),$$
where ``$\Isom$'' represents isometries, and 
$m$ runs through positive integers such that $am$ and $a'm$ are both integers.

If $s$ is a section of $a\CL$ represented by $s_m\in \Gamma(\CX,am\CL)$ or 
a rational section of $a\CL$ represented by $s_m\in \Gamma(\eta,am\CL)$, where $\eta\in\CX$ is the generic point as above, then define
$$
\div(s):=\frac{1}{m}\div(s_m),\quad\
\wh\div(s):=\frac{1}{m}\wh\div(s_m).
$$
These are respectively $\QQ$-divisors and arithmetic $\QQ$-divisors on $\CX$.

An arithmetic $\BQ$-divisor $\OCD \in \wh\Div(\CX)_\BQ$ (resp. hermitian $\BQ$-line bundle $\CLL\in \wh\Picc(\CX)_\QQ$) is called \emph{nef} if for some
positive integer $m$, the multiple $m \OCD$ (resp. $m\CLL$) is a nef arithmetic divisor in $\wh\Div(\CX)$ (resp. nef hermitian line bundle in $\wh\Picc(\CX)$) in the usual sense.

\subsection{$(\QQ,\ZZ)$-divisors}

Let $\CU$ be an open subscheme of an integral scheme $\CX$. 
Define $\Div(\CX,\CU)$ to be the fiber product of the natural map 
$\phi:\Div(\CX)_\QQ\to \Div(\CU)_\QQ$ with the natural map $\psi:\Div(\CU)\to \Div(\CU)_\QQ$; i.e.
$$
\Div(\CX,\CU)=\ker(\phi-\psi:\Div(\CX)_\QQ\oplus \Div(\CU)\to \Div(\CU)_\QQ).
$$
In other words, $\Div(\CX,\CU)$ is the group of pairs $(\CD, \CD')$, where $\CD\in \Div(\CX)_\QQ$ and $\CD'\in \Div(\CU)$ have equal images in 
$\Div(\CU)_\QQ$.

An element $(\CD, \CD')$ of $\Div(\CX,\CU)$ is called a \emph{$(\QQ,\ZZ)$-divisor on $(\CX,\CU)$} or a \emph{$\QQ$-divisor of $\CX$ integral on $\CU$}.
We usually call $\CD$ \emph{the rational part} of $(\CD, \CD')$, and  
call $\CD'$ \emph{the integral part} of $(\CD, \CD')$.

By definition, there are projection maps 
$$
\Div(\CX,\CU)\lra \Div(\CX)_\QQ, \quad\
\Div(\CX,\CU)\lra \Div(\CU).
$$
By abuse of notations, we may abbreviate an element $(\CD, \CD')$ of $\Div(\CX,\CU)$ as $\CD$, and then write $\CD|_\CU$ for $\CD'$, viewed as an integral divisor on $\CU$.

There are canonical maps 
$$
\Div(\CX) \lra \Div(\CX,\CU), \quad \CE\longmapsto (\CE, \CE_\CU) 
$$
and 
$$
\Div(\CU)_\tor \lra \Div(\CX,\CU), \quad \CT\longmapsto (0, \CT).
$$
Here $\Div(\CU)_\tor$ is the subgroup of torsion elements of $\Div(\CU)$.
Then we have a canonical exact sequence 
$$
0\lra \Div(\CU)_\tor\lra \Div(\CX,\CU)\lra \Div(\CX)_\QQ\lra \Div(\CU)_\QQ/\Div(\CU).
$$
Then there is a canonical isomorphism 
$$
\Div(\CX,\CU)_\QQ\iso \Div(\CX)_\QQ.
$$

Take quotient
$$
\CaCl(\CX,\CU):=\Div(\CX,\CU)/\Pr(\CX).
$$
Here $\Pr(\CX)$ is mapped to $\Div(\CX,\CU)$ via $\Div(\CX)\to \Div(\CX,\CU)$. 
Note that $\Pr(\CX)\to \Div(\CX,\CU)$ is not necessarily injective, but the quotient makes sense by group action. 
There are canonical maps
$$
\CaCl(\CX,\CU)\lra \CaCl(\CX)_\QQ, \quad\
\CaCl(\CX,\CU)\lra \CaCl(\CU).
$$

An element of $\Div(\CX,\CU)$ is called \emph{effective} if its images in $\Div(\CX)_\QQ$ and $\Div(\CU)$ are both effective.

If $\CX$ is a projective variety over a field, an element of $\Div(\CX,\CU)$ is called \emph{nef} if its image in $\Div(\CX)_\QQ$ is nef.

\subsection{Arithmetic $(\QQ,\ZZ)$-divisors}

The above mixed notions extend easily to the arithmetic situation.  

Let $\CU$ be an open subscheme of a projective arithmetic variety $\CX$. 
Define $\wh\Div(\CX,\CU)$ to be the fiber product of the natural map 
$\wh\phi:\wh\Div(\CX)_\QQ\to \Div(\CU)_\QQ$ with the natural map $\wh\psi:\Div(\CU)\to \Div(\CU)_\QQ$; i.e.
$$
\wh\Div(\CX,\CU)=\ker(\wh\phi-\wh\psi:\wh\Div(\CX)_\QQ\oplus \Div(\CU)\to \Div(\CU)_\QQ).
$$
In other words, $\wh\Div(\CX,\CU)$ is the group of pairs $(\OCD, \CD')$, where $\OCD\in \wh\Div(\CX)_\QQ$ and $\CD'\in \Div(\CU)$ have equal images in 
$\Div(\CU)_\QQ$.
Note that the second component uses $\Div(\CU)$ (instead of $\wh\Div(\CU)$), so it puts no condition on the Green function. 

An element of $\wh\Div(\CX,\CU)$ is called an \emph{arithmetic $(\QQ,\ZZ)$-divisor on $(\CX,\CU)$} or an \emph{arithmetic $\QQ$-divisor of $\CX$ integral on $\CU$}.
We usually call $\OCD$ \emph{the rational part} of $(\OCD, \CD')$, and  
call $\CD'$ \emph{the integral part} of $(\OCD, \CD')$.

There are projection maps 
$$
\wh\Div(\CX,\CU)\to \wh\Div(\CX)_\QQ,\qquad
\wh\Div(\CX,\CU)\to \Div(\CU)
$$
By abuse of notations, we may abbreviate an element $(\OCD, \CD')$ of $\wh\Div(\CX,\CU)$ as $\OCD$, and then write $\OCD|_\CU$ for $\CD'$, viewed as an integral divisor on $\CU$.

There are canonical maps
$$
\wh\Div(\CX) \lra \wh\Div(\CX,\CU),\quad
\Div(\CU)_\tor \lra \wh\Div(\CX,\CU).
$$
We have a canonical exact sequence 
$$
0\lra \Div(\CU)_\tor\lra \wh\Div(\CX,\CU)\lra \wh\Div(\CX)_\QQ\lra \Div(\CU)_\QQ/\Div(\CU).
$$
Then there is a canonical isomorphism 
$$
\wh\Div(\CX,\CU)_\QQ\lra \wh\Div(\CX)_\QQ.
$$
There is also a canonical injection
$$
\ker(\wh\Div (\CX)_\QQ\to \Div (\CU)_\QQ) \lra \wh\Div (\CX,\CU),\quad\
\OCD\longmapsto (\OCD,0).
$$

Take quotient
$$
\wh\CaCl(\CX,\CU):=\wh\Div(\CX,\CU)/\wh\Pr(\CX).
$$
Here the quotient is via the composition $\wh\Pr(\CX)\to\wh\Div(\CX) \to \wh\Div(\CX,\CU)$.

An element of $\wh\Div(\CX,\CU)$ is called \emph{effective} if its images in $\wh\Div(\CX)_\QQ$ and $\Div(\CU)$ are both effective.

An element of $\wh\Div(\CX,\CU)$ is called \emph{nef} if its image in $\wh\Div(\CX)_\QQ$ is nef.

\section{Essentially quasi-projective schemes} \label{sec models}

This section aims to define some basic terms about arithmetic models and introduce quasi-projective schemes, a special class of schemes on which we can naturally define adelic divisors and adelic line bundles.

\subsection{Pro-open immersions}

A morphism $i:X\to Y$ of schemes is called a  \emph{pro-open immersion}
if it satisfies the following two conditions:
\begin{enumerate}
\item[(i)] $i$ is injective as a map between the underlying topological spaces;
\item[(ii)] $i$ induces isomorphisms between the local rings; i.e. for any point $x\in X$, the induced map $\CO_{Y,i(x)}\to \CO_{X,x}$ is an isomorphism. 
\end{enumerate}

By Raynaud \cite[Prop. 1.1]{Ray1}, pro-open immersions are exactly flat monomorphisms, and they are systematically studied in the loc. cit. We have the following equivalent definitions.

\begin{prop} \label{flat mono}
Let $i: X\to Y$ be a morphism of quasi-compact schemes. 
Then the following are equivalent:
\begin{enumerate}

\item The morphism $i$ is a flat monomorphism; i.e, $i$ is flat, and $\Hom(S, X)\to \Hom(S,Y)$ is injective for any scheme $S$. 

\item The morphism $i$ is a pro-open immersion; i.e. $i$ induces an injection between the underlying spaces and isomorphisms between the local rings. 

\item The map $i:X\to Y$ is a homeomorphism of $X$ to its image $i(X)$ endowed with topology induced from $Y$;  the image $i(X)$ is equal to the intersection of its open neighborhoods in $Y$; the natural morphism $(X,\CO_X)\to (X, i^{-1}\CO_Y)$ is an isomorphism of ringed spaces. Here $i^{-1}\CO_Y$ denotes the pull-back as abelian sheaves.
\end{enumerate}
\end{prop}

\begin{proof}
See \cite[Prop. 1.1, Prop. 1.2]{Ray1}.
\end{proof}

Another property of \cite[Prop. 1.2]{Ray1} is as follows.

\begin{lem} 
Let $i: X\to Y$ be a pro-open morphism of quasi-compact schemes. 
If $Y$ is noetherian, then $X$ is also noetherian. 
\end{lem}

To justify the term ``pro-open,'' note that a pro-open immersion to a scheme $Y$ is given by the projective limit of some system of open subschemes of $Y$;
see \cite[Prop. 2.3]{Ray1}.
We refer to the loc. cit.  for more properties.

\subsection{Essentially quasi-projective schemes}

\kkk
We take the convention in \S\ref{sec uniform} for objects over $k$. 

Recall that, by a \emph{projective variety (resp. quasi-projective variety) over $k$}, we mean an integral scheme, projective (resp. quasi-projective) and flat over $k$. 
We make the following further definitions.

\begin{enumerate}
\item 
A flat integral noetherian scheme $X$ over $k$ is called \emph{essentially quasi-projective over $k$} if there is a pro-open immersion $i: X\to \CX$ over $k$ for some projective variety $\CX$ over $k$. 
\item 
Let $X$ be a flat and essentially quasi-projective integral scheme over $k$.
By a \emph{quasi-projective model} (resp. \emph{projective model}) of  $X$ over $k$, we mean a pro-open immersion $X\to \CU$ (resp. $X\to \CX$) for a quasi-projective variety $\CU$ (resp. projective variety) over $k$.
\end{enumerate}

The following are three important and natural classes of essentially quasi-projective schemes over $k$:
\begin{enumerate}
\item[(a)]  a quasi-projective variety over $k$,
\item[(b)]  a quasi-projective variety $X$ over a finitely generated field $F$ over $k$ (including the case $X=\Spec F$),
\item[(c)] the spectrum of the local ring of a quasi-projective variety over $k$ at a point.
\end{enumerate}
In this book, we are mainly concerned with cases (a) and (b).
If $X$ is in case (a), any pro-open immersion $X\to \CX$ to a projective variety $\CX$ over $k$ is necessarily an open immersion, so the notions of projective models regarding $X$ as a quasi-projective variety and as an essentially quasi-projective variety coincide.
If $X$ is in case (b), its quasi-projective model is not as arbitrary as it seems. Lemma \ref{models} asserts that it essentially comes from the generic fiber of a morphism $\CU\to\CV$ of quasi-projective varieties over $k$.

\subsection{More properties}

The following result describes the pro-open immersion in case (b).

\begin{lem} \label{models}
Let $F$ be a finitely generated field over $k$, $X$ be a quasi-projective variety over $F$, and $i: X\to \CU$ be a quasi-projective model of $X$ over $k$. Then
there is an open subscheme $\CU'$ of $\CU$ containing the image of $X$ together with a flat morphism $\CU'\to \CV$ of quasi-projective varieties over $k$, such that the function field of $\CV$ is isomorphic to $F$, and that 
the generic fiber of $\CU'\to \CV$ is isomorphic to $X\to \Spec F$.

Furthermore, if $X$ is projective over $F$, then we can assume that the flat morphism $\CU'\to \CV$ is projective.
\end{lem}
\begin{proof}
The last statement follows from the quasi-projective case by choosing an open subscheme of $\CV$ since projectivity is an open condition. 

For the quasi-projective case, let $\CV$ be a quasi-projective model of $\Spec F$ over $k$. 
Then the rational map 
$\CU\dashrightarrow \CV$ is defined on an open neighborhood of $X$ in $\CU$. 
Replacing $\CU$ by an open subscheme if necessary, we can assume that the rational map extends to a morphism $\CU\to \CV$. 
Denote by $\eta\in \CV$ the generic point and by $\CU_\eta\to \eta$ the generic fiber of $\CU\to \CV$. 
By the universal property of the fiber product of $\CU\to \CV$ and $\eta\to \CV$, we have a morphism $j:X\to \CU_\eta$ over $F$, whose composition with $\CU_\eta\to \CU$ is exactly $i:X\to \CU$.

The morphism $j: X\to \CU_\eta$ is flat and of finite type, so it is an open map. 
In particular, the image $j(X)$ is open in $\CU_\eta$. 
By Proposition \ref{flat mono}(3), $j:X\to \CU_\eta$ is an open immersion. 
Then the result follows. 
\end{proof}

In the general case, we have the following result about the inverse systems of quasi-projective models. 

\begin{lem} \label{models2} 
Let $X$ be a flat and essentially quasi-projective integral scheme over $k$, and let $\CU$ be a fixed quasi-projective model of $X$ over $k$. 
Then the inverse system of open neighborhoods of $X$ in $\CU$ is cofinal to the inverse system of quasi-projective models of $X$ over $k$.
\end{lem}

\begin{proof}
Let $\CU'$ be a quasi-projective model of $X$ over $k$. 
Then the rational map 
$\CU\dashrightarrow \CU'$ is defined on an open neighborhood $\CU''$ of $X$ in $\CU$. 
Then the system $\{\CU''\}$ is cofinal to the system $\{\CU'\}$.
\end{proof}

\subsection{Effectivity of Cartier divisors}

Over a normal scheme, the effectivity of a Cartier divisor can be checked in terms of the effectivity of the corresponding Weil divisor. Then we have the following result.

\begin{lem}  \label{effectivity normal}
Let $\CX$ be a normal integral scheme, and let $\psi:\CX'\to \CX$ be a birational proper morphism of integral schemes.
Let $\CD$ be a Cartier divisor on $\CX$. 
Then the following are equivalent:
\begin{enumerate}
\item $\CD$ is effective on $\CX$;
\item $\psi^*\CD$ is effective on $\CX'$;  
\item $\CD$ is effective as a $\QQ$-divisor on $\CX$; i.e. $m\CD$ is effective for some positive integer $m$.
\end{enumerate}
\end{lem}
\begin{proof}
The proof is straightforward using Weil divisors, except that to prove that (2) implies (1), we need to replace $\CX'$ by its normalization. 
\end{proof}

Without normality, the situation is more delicate. 
The following result solves the problem for our purpose. 
Recall that for a dominant morphism of integral schemes $V\to W$, we say that  \emph{$W$ is integrally closed in $V$} if the normalization of $W$ in $V$ is isomorphic to $W$.

\begin{lem}  \label{effectivity}
Let $i: X\to \CX$ be a pro-open immersion of integral noetherian schemes. 
Assume that $\CX$ is integrally closed in $X$.
Then a Cartier divisor $\CD$ on $\CX$ is effective if and only if the following two conditions hold simultaneously:
\begin{enumerate}
\item the pull-back $\CD|_X$ is effective on $X$;
\item for any $v\in \CX\setminus X$ of codimension one in $\CX$, the valuation 
$\ord_v(\CD)$ in the discrete valuation ring $\CO_{\CX, v}$ is non-negative.  
\end{enumerate}
 \end{lem}

\begin{proof}
We first claim that the local ring $\CO_{\CX, v}$ in (2) is a discrete valuation ring.
In fact, the base change of $X\to \CX$ by $\Spec \CO_{\CX, v} \to \CX$ is exactly 
$\Spec k(X)\to \Spec \CO_{\CX, v}$. Here $k(X)$ denotes the function field of $X$, which is also the fraction field of $\CO_{\CX, v}$. As a consequence, $\Spec \CO_{\CX, v}$ is integrally closed in $\Spec k(X)$. Then $\CO_{\CX, v}$ is a discrete valuation ring since it has dimension 1.

To prove the lemma, we only need to prove the ``if'' part. 
If $\CX$ is normal, then we can write Cartier divisors in terms of Weil divisors, and the effectivity of them are equivalent. 

In the general case, let $f$ be a local equation of $\CD$ in an affine open subscheme $\CW$ of $\CX$.
Then $f\in k(\CW)^\times$, and we need to prove $f\in \CO(\CW)$. 
By the normal case, $f$ is regular on the normalization $\CW'$ of $\CW$. 
As a consequence, $f$ is integral over $\CO(\CW)$. 
Note that $f\in \CO(\CW\cap X)$ by the assumption that $\CD|_X$ is effective.
By assumption, $\CO(\CW)$ is integrally closed in $\CO(\CW\cap X)$.  
Therefore, $f\in \CO(\CW)$. 
This finishes the proof.
\end{proof}

Now we have the following variant of Lemma \ref{effectivity normal} for non-normal schemes.

\begin{lem}  \label{effectivity birational}
Let $i':X\to \CX'$ and $i: X\to \CX$ be pro-open immersions of integral noetherian schemes, and let $\psi:\CX'\to \CX$ be a birational and proper morphism such that $i\circ \psi=i'$.
Assume that $\CX$ is integrally closed in $X$. 
Let $\CD$ be a Cartier divisor on $\CX$. 
Then the following are equivalent:
\begin{enumerate}
\item $\CD$ is effective on $\CX$;
\item $\psi^*\CD$ is effective on $\CX'$;  
\end{enumerate}
\end{lem}
\begin{proof}
It is trivial that (1) implies (2). 
For the opposite direction, replace $\CX'$ by its normalization in $X$, and apply Lemma \ref{effectivity}. 
\end{proof}

\section{Adelic divisors} \label{sec adelic divisors}

\kkk
Let $X$ be a flat and essentially quasi-projective integral scheme over $k$. The goal of this section is to define adelic divisors on $X$. 
We will start the definition for a quasi-projective variety over $k$, and the general case is obtained as a direct limit over all quasi-projective models.

\subsection{Adelic divisors on a quasi-projective variety}

\kkk
\ccc
Let $\CU$ be a quasi-projective variety over $k$. 

Let $\CX$ be a projective model of $\CU$ over $k$. 
In the spirit of \S\ref{sec uniform}, denote by $\wh\Div(\CX,\CU)$ the group of arithmetic $(\QQ,\ZZ)$-divisors on $(\CX,\CU)$. 
Hence, in the geometric case (that $k$ is a field), 
we take the convention $\wh\Div(\CX,\CU)=\Div(\CX,\CU)$, and an arithmetic 
$(\QQ,\ZZ)$-divisor in this case just means a $(\QQ,\ZZ)$-divisor.
Both cases are introduced in \S\ref{sec convention}.

Projective models $\CX$ of $\CU$ over $k$ form an inverse system. 
Using pull-back morphisms, we can form the direct limits:
$$\wh\Div (\CU/k)_\rmod:=\lim_{\substack{\lra\\ \CX}}\wh\Div(\CX,\CU),$$
$$\wh\Pr (\CU/k)_\rmod:=\lim_{\substack{\lra\\ \CX}}\wh\Pr(\CX).$$
Here the subscript ``$\rmod$'' represents ``model divisors'', as these divisors are defined on single projective models. Now we are going to introduce a topology on $\wh\Div (\CU/k)_\rmod$.

For any $\OCD, \OCE \in \Divhat(\CX)_\BQ$, 
 write $\OCD\geq\OCE$ or $\OCE\leq \OCD$ if $\OCD-\OCE$ is effective. It is a partial order in $\Divhat(\CX)_\BQ$. This induces a partial order in $\Divhat(\CX,\CU)$ by the law that $\OCD\geq \OCE$ or $\OCE\leq \OCD$ in $\Divhat(\CX,\CU)$ if the image of $\OCD- \OCE$ in $\Divhat(\CX)_\BQ$ and the image of $\OCD- \OCE$ in $\Div(\CU)$ are both effective.
By direct limit, we have an induced partial order in $\wh \Div (\CU/k)_{\rmod}$, and we will use the same symbols for it.

In the geometric case (that $k$ is a field), by a \emph{boundary divisor of $\CU/k$}, we mean a pair $(\CX_0,\CE_0)$ consisting of a projective model $\CX_0$ of $\CU$ over $k$ and an effective Cartier divisor $\CE_0$ on $\CX_0$ with support equal to $\CX_0\setminus \CU$.
To see the existence of $(\CX_0,\CE_0)$, take any projective  model 
$\CX'_0$ of $\CU$ over $k$, and blow-up $\CX'_0$ 
along the reduced center $\CX'_0\setminus\CU$. We get a projective  
model $\CX_0$ of $\CU$, and the exceptional divisor of $\CX_0\to \CX_0'$ is a Cartier divisor with support $\CX_0\setminus \CU$. 

In the arithmetic case ($k=\ZZ$), by a \emph{boundary divisor of $\CU/k$}, we mean a pair $(\CX_0,\OCE_0)$ consisting of a projective model $\CX_0$ of $\CU$ over $k$ and a \emph{strictly effective} Cartier divisor $\OCE_0$ on $\CX_0$ such that the support of the finite part $\CE_0$ is equal to $\CX_0\setminus \CU$.
 
To unify the terminology, in the geometric case, write $\OCE_0=\CE_0$ in $\Divhat(\CX_0)=\Div(\CX_0)$. Then in both cases, a boundary divisor is written in the form $(\CX_0,\OCE_0)$, and the following are for both cases. 

For any $r\in \QQ$, view $r\OCE_0$ as an element of $\wh\Div (\CX_0,\CU)$ by setting its integral part in $\Div(\CU)$ to be 0.
Then $r\OCE_0$ is also viewed as an element of
$\wh\Div (\CX)_\rmod$. 

Let $(\CX_0,\OCE_0)$ be a {boundary divisor of $\CU/k$}.
We have a \emph{boundary norm} 
$$\|\cdot\|_{\OCE_0}:\wh\Div (\CU/k)_\rmod
\lra [0,\infty]$$
defined by 
$$
\|\OCD\|_{\OCE_0}:=\inf\{\epsilon\in \BQ_{>0}: \ 
 -\epsilon \OCE_0 \leq
\OCD \leq  \epsilon \OCE_0\}.
$$
Here we take the convention that $\inf(\emptyset)=\infty$.
Note that $\|\cdot\|_{\OCE_0}$ can take value infinity, but it is an \emph{extended norm} in the sense of \cite[Def. 1.1]{Bee}. We refer to \cite{Bee} for more theory on extended norms.
In our situation, we have the following basic properties.

\begin{lem} \label{boundary norm}
The boundary norm $\|\cdot\|_{\OCE_0}$ on $\wh\Div (\CU/k)_\rmod$ satisfies the following properties:
\begin{enumerate}
\item $\|\OCD\|_{\OCE_0}=0$ if and only if $\OCD=0$;
\item $\|\OCD_1+\OCD_2\|_{\OCE_0}\leq \|\OCD_1\|_{\OCE_0}+\|\OCD_2\|_{\OCE_0}$;
\item $\|a\OCD\|_{\OCE_0}\leq |a|\cdot\|\OCD\|_{\OCE_0}$ for any nonzero $a\in \ZZ$. The inequality is strict if and only if both $\CD\neq 0$ and $a\CD=0$ hold in 
$\Div (\CU)$, where $\CD$ denotes the image of $\OCD$ in $\Div (\CU)$.
\end{enumerate}
Moreover, if $(\CX_0',\OCE_0')$ is another boundary divisor, then 
$\|\cdot\|_{\OCE_0'}$ is equivalent to $\|\cdot\|_{\OCE_0}$
in the sense that there is a real number $r>1$ such that 
$$r^{-1}\|\cdot\|_{\OCE_0}\leq \|\cdot\|_{\OCE_0'}\leq r \|\cdot\|_{\OCE_0}.$$
\end{lem}
\begin{proof} 
Note that (2) and (3) are automatic by definition.
For (1), assume that $\|\OCD\|_{\OCE_0}=0$ for some $\OCD$; i.e. $-\epsilon \OCE_0 \leq \OCD \leq  \epsilon \OCE_0$ in $\wh\Div (\CU/k)_\rmod$ for all positive rational numbers $\epsilon$. 
Assume that $\OCD$ comes from $\wh\Div(\CX,\CU)$ for a projective model $\CX$ of $\CU$, and assume that $\CX$ dominates $\CX_0$ and is integrally closed in $\CU$.
By Lemma \ref{effectivity birational}, $-\epsilon \OCE_0 \leq \OCD \leq  \epsilon \OCE_0$ holds in $\wh\Div (\CX,\CU)$ for all positive rational numbers $\epsilon$. 
Then we can conclude $\OCD=0$ by Lemma \ref{effectivity}. 

For the equivalence of the two norms, it suffices to find a rational number $r>1$ such that $r^{-1}\OCE_0\leq \OCE_0' \leq r \OCE_0$ in $\wh\Div (\CU/k)_\rmod$.
In fact, we can find a third projective model $\CY$ of $\CU$ dominating both $\CX_0$ and $\CX_0'$, and
we can further assume that 
$\CY$ is integrally closed in $\CU$. 
Then we only need to treat the inequalities over $\CY$, which is an easy consequence of Lemma \ref{effectivity}. 
\end{proof}

The \emph{boundary topology on $\wh\Div (\CU/k)_\rmod$} is the topology induced by the
boundary norm
$\|\cdot\|_{\OCE_0}$.
Thus, a neighborhood basis at $0$ of the topology is given by
$$
 B(\epsilon, \wh\Div (\CU/k)_\rmod):
 =\{\OCD \in \Divhat(\CU/k)_\rmod: \ 
 -\epsilon \OCE_0 \leq
\OCD \leq  \epsilon \OCE_0\}, \quad \epsilon\in \BQ_{>0}.
$$
By translation, it gives a neighborhood basis at any point.
The topology does not depend on the choice of the boundary divisor 
$(\CX_0,\OCE_0)$ by the lemma. 

Let $\wh \Div  (\CU/k)$ be the \emph{completion} of $\wh \Div  (\CU/k)_\rmod$ for the  boundary topology. 
An element of $\wh \Div  (\CU/k)$ is called an \emph{adelic divisor} (or a \emph{compactified divisor}) of $\CU/k$.
By definition, an adelic divisor is represented by a Cauchy sequence in $\wh \Div  (\CU/k)_\rmod$, i.e. a sequence $\{\OCD_i\}_{i\geq 1}$ in $\wh \Div  (\CU/k)_\rmod$ satisfying the property that there is a sequence $\{\epsilon_i\}_{i\geq 1}$ of positive rational numbers converging to $0$ such that 
$$
 -\epsilon_i \OCE_0 \leq
\OCD_{i'}-\OCD_{i} \leq  \epsilon_i \OCE_0,\quad\ i'\geq i\geq 1.
$$
The sequence $\{\OCD_i\}_{i\geq 1}$ represents $0$ in $\wh \Div  (\CU/k)$ if and only if there is a sequence $\{\delta_i\}_{i\geq 1}$ of positive rational numbers converging to $0$ such that 
$$
 -\delta_i \OCE_0 \leq
 \OCD_i \leq  \delta_i \OCE_0,\quad\ i\geq 1.
$$

Define the \emph{class group of adelic divisors} of $\CU$ to be 
$$\wh\CaCl (\CU/k):=\wh\Div  (\CU/k)/\wh \Pr (\CU/k)_\rmod.$$
The map $\wh\Div(\CX,\CU)\to \Div(\CU)$ induces maps 
$$\wh\Div (\CU/k)_\rmod\lra \Div(\CU), \quad\
\wh\CaCl (\CU/k)_\rmod\lra \CaCl(\CU).$$
We call these maps \emph{restriction maps} or \emph{forgetful maps}.

\begin{rmk} \label{continuous vs smooth}
In the arithmetic case $k=\ZZ$, for the definition 
$$\wh\Div (\CU/k)_\rmod=\lim_{\substack{\lra\\ \CX}}\wh\Div(\CX,\CU),$$
we allow elements of $\wh\Div(\CX,\CU)$ to have Green functions of continuous type instead of Green functions of smooth type. See \S\ref{sec hermitian} for the definitions of these terms. 
However, both choices give the same completion $\wh\Div (\CU/k)$ since continuous functions on $\CX(\CC)$ can be approximated by smooth functions uniformly. 
\end{rmk}

\subsection{Completion of the divisor class group}

Let $\CU$ be a quasi-projective variety over $k$. 
Consider the class group of model divisors:
$$\wh\CaCl (\CU/k)_\rmod= \lim_{\substack{\lra\\ \CX}}
\wh\CaCl(\CX,\CU) \simeq \wh\Div(\CU/k)_\rmod/\wh\Pr (\CU/k)_\rmod.$$
It is endowed with the quotient topology induced by the boundary topology of 
$\wh\Div(\CU/k)_\rmod$.
On the other hand, 
$$
\wh\CaCl (\CU/k)=\wh\Div(\CU/k)/\wh\Pr (\CU/k)_\rmod
$$
is not defined to be the completion of $\wh\CaCl (\CU/k)_\rmod$. 
However, the following result asserts that these two are isomorphic. 

\begin{lem}
The group $\wh\Pr (\CU/k)_\rmod$ is discrete in $\wh\Div(\CU/k)_\rmod$ under the boundary topology. Therefore, $\wh\CaCl (\CU/k)$
is canonically isomorphic to the completion of $\wh\CaCl (\CU/k)_\rmod$.
\end{lem}

\begin{proof}

It suffices to prove the first statement. In the following, we assume the arithmetic case $k=\ZZ$ since the geometric case is similar and easier.

Assume that there is a sequence $\OCD_i$ in $\wh\Pr (\CU/k)_\rmod$ converging to 0. Then there is a sequence $\{\epsilon_i\}_{i\geq 1}$ of positive rational numbers converging to $0$ such that 
$\epsilon_i \OCE_0 \pm \OCD_i \geq 0$ 
in $\wh\Div (\CU/k)_\rmod$ for any $i\geq 1$. 
Assume that $\OCD_i=\wh\div_{\CX_i}(f_i)$ for a projective model $\CX_i$ of $\CU$ and a rational function $f_i\in \QQ(\CX_i)^\times=\QQ(\CU)^\times$.

We first consider the case that $\CU$ and $\CX_i$ are normal for $i\geq0$.
For $i=0$, recall that the projective model $\CX_0$ is the one chosen to define $\OCE_0$. 
By Lemma \ref{effectivity normal}, the relation $\epsilon_i \OCE_0 \pm \OCD_i \geq 0$ 
in $\wh\Div (\CU/k)_\rmod$
is the same as the relation $\epsilon_i \OCE_0 \pm \wh\div_{\CX_0}(f_i) \geq 0$ 
in $\wh\Div (\CX_0)_\QQ$. 
When $\epsilon_i$ is small enough, we must have $\div_{\CX_0}(f_i)=0$.
This implies that $f_i\in \CO(\CX_0)^\times=O_K^\times$.
Here $K$ is the algebraic closure of $\QQ$ in $\QQ(\CX_0)$, and $O_K$ is the ring of algebraic integers in $K$.
In the setting of Dirichlet's unit theorem, the image of $O_K^\times$ in $\RR^{r}$ under the logarithms of archimedean absolute values is discrete.  
Then the relation $\epsilon_i g \pm \log|f_i| \geq  0$ from the Green function implies that $|f_i|_\sigma=1$ for any archimedean place $\sigma$ of $K$ and for sufficiently large $i$. 
Therefore, $f_i$ is a root of unity, and thus $\OCD_i=0$ for such $i$. 
This proves the normal case.

For the general case that $\CX_i$ is not normal, denote by $\CX_i'$ (resp. $\CU'$) the normalization of $\CX_i$ (resp. $\CU$) for all $i\geq 0$. 
Consider the pull-back of the relation $\epsilon_i \OCE_0 \pm \OCD_i \geq 0$ 
to the normalizations. Then the previous case implies that $f_i$ is a root of unity for sufficiently large $i$. 
This implies the image of $\OCD_i=\wh\div_{\CX_i}(f_i)$ in $\wh\Div(\CX_i)_\QQ$ is 0. 
On the other hand, by definition of Cauchy sequences, the integral part $\CD_i|_\CU$ is constant in $\Div(\CU)$. Therefore, the sequence $\OCD_i$ in 
$\wh\Div(\CU/k)_\rmod$, which is a subgroup of $\ds\varinjlim_\CX\wh\Div(\CX)_\QQ \oplus \Div(\CU)$,
is eventually constant.
The proof is complete.
\end{proof}

\subsection{Adelic divisors on essentially quasi-projective schemes}

\kkk
Let $X$ be a flat and essentially quasi-projective integral scheme over $k$ as in \S\ref{sec models}. The set of quasi-projective models $\CU$ of $X$ over $k$ form an inverse system.
Define
\begin{eqnarray*}
\wh\Div(X/k):&=&\lim_{\substack{\lra \\ \CU}}\wh \Div (\CU/k),\\
\wh\CaCl (X/k): &=&\lim_{\substack{\lra \\ \CU}}\wh \CaCl (\CU/k).\end{eqnarray*}
We call elements of $\wh\Div (X/k)$ \emph{adelic divisors} of $X/k$.

If $X$ is a quasi-projective variety over $k$, then $X$ itself is the final object of the inverse system of quasi-projective models of $X$ over $k$. In this case, the definitions in terms of direct limits are compatible with the original ones for quasi-projective varieties.

There are many functorial properties of $\wh\Div(X/k)$ and $\wh\CaCl (X/k)$, which will be introduced together with the theory of adelic line bundles.

\section{Adelic line bundles} \label{sec adelic line bundles}

Now  we define adelic line bundles on essentially quasi-projective schemes to match the above definition of adelic divisors. We will use the notion of hermitian $\QQ$-line bundles in \S\ref{sec convention} and arithmetic models in \S\ref{sec models}. 

Throughout this section, let $k$ be either $\ZZ$ or an arbitrary field.
\ccc

\subsection{Adelic line bundles on a quasi-projective variety}

\kkk
Let $\CU$ be a quasi-projective variety over $k$. 

Let us first introduce a notation for \emph{model adelic divisors of rational maps}.
Let $\CX_1, \CX_2$ be projective models of $\CU$ over $k$. 
Let $\CLL_i$ be a hermitian $\QQ$-line bundle on $\CX_i$ for $i=1,2$.
By a  \emph{rational map} $\ell:\CLL_1\dashrightarrow \CLL_2$ over $\CU$, we mean  an isomorphism $\ell:\CL_1|_\CU\to \CL_2|_\CU$ of $\QQ$-line bundles on $\CU$. 
Let $\CY$ be a projective model of $\CU$ with morphisms
$\tau_i:\CY\to \CX_i$ of projective models of $\CU$.
View $\ell$ as a rational section of $\tau_1^*\CL_1^\vee\otimes \tau_2^*\CL_2$ on $\CY$, so that it defines an arithmetic $\QQ$-divisor 
$\wh\div_{\CY}(\ell)$ on $\CY$ using the metric of $\tau_1^*\CLL_1^\vee\otimes \tau_2^*\CLL_2$. 
Set $\wh\div(\ell)$ to be the image of 
$\wh\div_{\CY}(\ell)$ in $\wh \Div (\CU/k)_{\rmod,\QQ}$.
We also view $\wh\div(\ell)$ as an element of $\wh \Div (\CU/k)_{\rmod}$ by setting the integral part on $\CU$ to be 0. 
The definition of $\wh\div(\ell)$ is independent of the choice of $\CY$.

Let $(\CX_0,\OCE_0)$ be a boundary divisor as in \S\ref{sec adelic divisors}.
Namely, $\CX_0$ is a  projective model of $\CU$ and $\OCE_0=(\CE_0,g_0)$ is a (strictly) effective arithmetic divisor on $\CX_0$ whose finite part $\CE_0$ has support equal to $\CX_0\setminus \CU$.

Define the \emph{category $\wh\Picc (\CU/k)$ of adelic line bundles} on $\CU$ as follows.
An object of $\wh\Picc (\CU/k)$ is a pair
$(\CL, (\CX_i,\overline \CL_i, \ell_{i})_{i\geq 1})$ where:
\begin{enumerate}
\item $\CL$ is an object of $\Picc(\CU)$, i.e. a line bundle on $\CU$;

\item  $\CX_i$ is a projective model of $\CU$ over $k$;

\item  $\overline \CL_i$ is an object of $\wh\Picc(\CX_i)_\QQ$, i.e. a hermitian $\QQ$-line bundle on $\CX_i$;

\item $\ell_i:\CL\to \CL_i|_{\CU}$ is an isomorphism in $\Picc(\CU)_\QQ$.
\end{enumerate}
The sequence is required to satisfy the \emph{Cauchy condition} as follows.
By (4), we obtain an isomorphism $\ell_i\ell_1^{-1}: \CL_1|_{\CU}\to \CL_i|_{\CU}$ of $\QQ$-line bundles, and thus a rational map $\ell_i\ell_1^{-1}: \CLL_1\dashrightarrow \CLL_i$ over $\CU$.
By the above notations, it defines a model divisor $\wh\div(\ell_i \ell_1^{-1})$ 
in $\wh \Div (\CU/k)_\rmod$. 
Then the \emph{Cauchy condition} is that
the sequence $\{\wh \div(\ell_i \ell_1^{-1})\}_{i\geq 1}$ is a 
Cauchy sequence in $\wh\Div(\CU/k)_\rmod$ under the boundary topology.
More precisely, 
there is a sequence $\{\epsilon_i\}_{i\geq 1}$ of positive rational numbers converging to $0$ such that 
$$
-\epsilon_i \OCE_0 \leq \wh \div (\ell_{i'} \ell_1^{-1})-\wh \div (\ell_i \ell_1^{-1}) \leq \epsilon_i \OCE_0
$$
in $\wh\Div(\CU/k)_\rmod$ for any $i'\ge i\geq 1$. 
The relation can also be written as
$$
-\epsilon_i \OCE_0 \leq \wh \div (\ell_{i'} \ell_i^{-1}) \leq \epsilon_i \OCE_0
$$
in $\wh\Div(\CU/k)_\rmod$ for any $i'\ge i\geq 1$. 

For convenience, the object
$(\CL,(\CX_i,\CLL_{i}, \ell_{i})_{i\geq1})$ 
is also called a \emph{Cauchy sequence}
in $\wh\Picc(\CU/k)_\rmod$. 
For simplicity, we may abbreviate $(\CL,(\CX_i,\overline \CL_i, \ell_{i})_{i\geq 1})$ as $(\CL,(\CX_i,\overline \CL_i, \ell_{i}))$ or simply as $(\CL,\CX_i,\overline \CL_i, \ell_{i})$.

A morphism from an object $(\CL, (\CX_i,\overline \CL_i, \ell_{i})_{i\geq 1})$ of $\wh\Picc (\CU/k)$ to another object
$(\CL',(\CX_i',\overline \CL_i', \ell_{i}')_{i\geq 1})$ of $\wh\Picc (\CU/k)$ is an isomorphism  $\iota:\CL\to \CL'$ of the integral line bundles on $\CU$ satisfying the following properties.
As above, the composition $\ell_i'\iota \ell_i^{-1}: \CL_i|_\CU\to \CL_i'|_\CU$ induces a rational map $\ell_i'\iota \ell_i^{-1}:\CLL_i\dashrightarrow \CLL_i'$, and thus defines a model divisor $\wh\div(\ell_i'\iota \ell_i^{-1})$ in $\wh\Div(\CU/k)_\rmod$ whose image in $\Div(\CU)$ is 0. 
Then we require the sequence 
$\{\wh\div(\ell_i'\iota \ell_i^{-1}) \}_{i\geq1}$
of $\wh \Div (\CU/k)_\rmod$ converges to 0 in $\wh \Div (\CU/k)$
under the boundary topology, i.e.
 there is a sequence $\{\delta_i\}_{i\geq 1}$ of positive rational numbers converging to $0$ such that 
$$
 -\delta_i \OCE_0 \leq
\wh\div(\ell_i'\iota \ell_i^{-1})  \leq  \delta_i \OCE_0,\quad\ i\geq 1.
$$
Note that the sequence $\{\wh\div(\ell_i'\iota \ell_i^{-1}) \}_{i\geq1}$ is already a Cauchy sequence by 
$$
\wh\div(\ell_i'\iota \ell_i^{-1})=
\wh\div(\ell_{i}'\ell_1'^{-1})-\wh\div(\ell_{i}\ell_1^{-1})+\wh\div(\iota_1) .
$$

By definition, any morphism in $\wh\Picc (\CU/k)$ is an isomorphism, so $\wh\Picc (\CU/k)$ is a groupoid. This category is equipped with a \emph{tensor product} given by 
$$(\CL,(\CX_i,\overline \CL_i, \ell_{i}))\otimes
(\CL',(\CX_i',\overline \CL_i', \ell_{i}')):=
(\CL\otimes \CL', (\CW_i,\tau_i^*\overline \CL_i\otimes \tau_i'^*\overline \CL_i', \ell_{i}\otimes \ell_{i}')),$$
where $\CW_i$ is the Zariski closure of the image of the diagonal map $\CU\to \CX_i\times_k \CX_i'$, and $\tau_i:\CW_i\to \CX_i$ and $\tau_i':\CW_i\to \CX_i'$ are the two projection maps.
It is also equipped with a \emph{dual} given by 
$$(\CL,(\CX_i,\overline \CL_i, \ell_{i}))^\vee:=
(\CL^\vee,(\CX_i,\overline \CL_i^\vee, \ell_{i}^\vee)).$$
Then the tensor product of an element with its dual is isomorphic to the neutral object
$(\CO_\CU,(\CX_0, \overline\CO_{\CX_0},1))$.

An object of $\wh\Picc (\CU/k)$ is called an \emph{adelic line bundle} (or a \emph{compactified line bundle}) on $\CU$.
Define $\wh\Pic (\CU/k)$ to be the \emph{group} of isomorphism classes of objects of $\wh\Picc (\CU/k)$, where the group operation is the above tensor product. 
We usually write an adelic line bundle in the form $\CLL=(\CL,(\CX_i,\CLL_{i}, \ell_{i})_{i\geq1})$, and call $\CL$ 
\emph{the underlying line bundle} of $\CLL$ on $\CU$.

As in the classical case, we have the following result. 

\begin{prop} \label{isomorphism}
Let $\CU$ be a quasi-projective variety over $k$. Then there is a canonical isomorphism
$$\wh\CaCl(\CU/k)\iso \wh\Pic(\CU/k).$$
\end{prop}

\begin{proof}
The proof is a routine, but we write in detail to familiarize the terminologies here. 
It suffices to define a map 
$$\wh\Div(\CU/k)\lra \wh\Pic(\CU/k)$$
and check that it is surjective with kernel $\wh\Pr(\CU/k)_\rmod$.

To define the map, take an element $\OCD$ of the left-hand side. 
Then $\OCD$ is represented by a Cauchy sequence $\{\OCD_i\}_{i\geq 1}$ in $\wh \Div  (\CU/k)_\rmod$. 
Let $\{\CX_i\}_{i\geq 1}$ be a system of projective models of $\CU$
such that $\OCD_i\in\wh \Div  (\CX_i,\CU)$ for any $i\geq1$.
Note that $\CD_1|_\CU=\CD_i|_\CU$ is an integral divisor on $\CU$.
Set the image of $\OCD$ in $\wh\Pic(\CU/k)$ to be the isomorphism class of the sequence 
$(\CL,(\CX_i,\overline \CL_i, \ell_{i}))$, where
$\CL=\CO(\CD_1|_\CU)$, 
$\overline \CL_i=\CO(\OCD_i)$, and $\ell_{i}:\CL\to \CL_i|_{\CU}$ is the isomorphism induced by the equality $\CD_1|_\CU=\CD_i|_\CU$ in $\Div(\CU)$. 
By definition, $\wh\div(\ell_{i}\ell_1^{-1})=\OCD_i-\OCD_1$. 
Then we see that $(\CL,(\CX_i,\overline \CL_i, \ell_{i}))$ satisfies the Cauchy condition. 
This defines the map. 

Now assume that the above adelic divisor $\OCD$ lies in the kernel of the map. 
It follows that there is an isomorphism from
$(\CO_\CU, (\CX_0,\overline \CO_{\CX_0}, 1))$ to $(\CL,(\CX_i,\overline \CL_i, \ell_{i}))$. 
This includes an isomorphism $\CO_\CU\to \CO(\CD_1|_{\CU})$, which is given by the multiplication by some $f\in \Gamma(\CU,\CO_\CU)^\times$ with $\div(f)=\CD_1|_\CU=0$ on $\CU$. 
The further properties of the isomorphism are equivalent to that $\CD_i$ converges to $-\wh\div(f)$ in $\wh \Div  (\CU/k)_\rmod$. 
This proves that the kernel is $\wh\Pr(\CU/k)_\rmod$.

To see the surjectivity of the map, let $\CLL=(\CL,(\CX_i,\overline \CL_i, \ell_{i}))$ be an adelic line bundle on $\CU$. 
For any rational section $s$ of $\CL$, denote 
$$
\wh\div_{\CLL}(s):=\wh\div_{(\CX_1,\CLL_1)}(s)+\lim_{i\to \infty} 
\wh\div(\ell_i\ell_1^{-1}),
$$
which is an element of $\wh\Div(\CU/k)$. 
This gives a preimage of $\CLL$. Then the map is surjective.
\end{proof}

\subsection{Nef and integrable adelic line bundles}

In \S\ref{sec convention}, we have recalled the notion of nef hermitian line bundles on arithmetic varieties. This notion is generalized to adelic line bundles by the limit process as follows: 

\begin{defn} \label{def nef}
Let $\CU$ be a quasi-projective variety over $k$. 
\begin{enumerate}
\item We say that an adelic line bundle
$\overline \CL\in \wh\Picc(\CU/k)$ is \emph{strongly nef} if it is isomorphic to an object
$(\CL,(\CX_i,\overline \CL_i, \ell_{i}))$ where each $\overline\CL_i$ is nef on $\CX_i$.
\item 
We say that an adelic line bundle
$\overline \CL\in \wh\Picc(\CU/k)$ is \emph{nef} if there exists a strongly nef adelic line bundle 
$\overline \CM\in \wh\Picc(\CU/k)$ such that $a\CLL+\CMM$ is strongly nef for all positive integers $a$. 
\item 
We say that an adelic line bundle
in $\overline \CL\in \wh\Picc(\CU/k)$ is \emph{integrable} if it is 
isomorphic to the difference of two strongly nef ones in  $\wh\Picc(\CU/k)$.
\end{enumerate}
\end{defn}

It is obvious that ``strongly nef'' implies ``nef'', and ``nef'' implies ``integrable''.
Denote by 
$$\wh \Picc (\CU/k)_\snef,\quad
\wh \Picc (\CU/k)_\nef,\quad
\wh \Picc (\CU/k)_\intb$$
respectively the 
\emph{full subcategories} of $\wh \Picc (\CU/k)$ of strongly nef objects, nef objects, and integrable objects.
Denote by 
$$\wh \Pic (\CU/k)_\snef,\quad
\wh \Pic (\CU/k)_\nef,\quad
\wh \Pic (\CU/k)_\intb$$
respectively the 
\emph{subsets} of $\wh \Pic (\CU/k)$ of strongly nef elements, nef elements, and integrable elements.
Then $\wh \Pic(\CU/k)_\snef$ and $\wh \Pic(\CU/k)_\nef$ are semigroups and 
$\wh \Pic (\CU/k)_\intb$ is a group. 

The preimages of 
$$\wh \Pic (\CU/k)_\snef,\quad
\wh \Pic (\CU/k)_\nef,\quad
\wh \Pic (\CU/k)_\intb$$
 in $\wh \Div (\CU/k)$ are denoted by 
$$\wh \Div (\CU/k)_\snef,\quad
\wh \Div (\CU/k)_\nef,\quad
\wh \Div (\CU/k)_\intb$$
respectively. Their elements are respectively called 
\emph{strongly nef} adelic divisors on $\CU/k$,
\emph{nef} adelic divisors on $\CU/k$, and
\emph{integrable} adelic divisors on $\CU/k$.

\subsection{Definition on essentially quasi-projective schemes}

\kkk
Let $X$ be a flat and essentially quasi-projective integral scheme over $k$. Define
\begin{eqnarray*}
\wh\Pic (X/k):&=&\lim_{\substack{\lra \\ \CU}}\wh \Pic (\CU/k),\\
\wh\Picc (X/k): &=&\lim_{\substack{\lra \\ \CU}}\wh \Picc (\CU/k).
\end{eqnarray*}
Here the limits are over all quasi-projective models $\CU$ of $X$ over $k$.
The category $\wh\Picc (X/k)$ defined by the direct limit is understood as follows. 
An object of $\wh\Picc(X/k)$ is a pair $(\CLL, \CU)$, where $\CU$ is a quasi-projective model of $X$ over $k$ and $\CLL$ is an object of $\wh\Picc(\CU/k)$. 
A morphism $(\CLL, \CU)\to (\CLL', \CU')$ between two objects of $\wh\Picc(X/k)$ is an isomorphism $\iota: \CL|_X\to \CL'|_X$ in $\Picc(X)$ satisfying the property that for some quasi-projective model $\CV$ of $X$ over $k$ endowed with open immersions $\psi:\CV\to \CU$ and $\psi':\CV\to \CU'$ extending the identity morphism $X\to X$, the isomorphism $\iota: \CL|_X\to \CL'|_X$ can be extended to an isomorphism $\CL|_\CV\to \CL'|_\CV$ in $\Picc (\CV)$ and induces
an isomorphism $\CLL|_\CV\to \CLL'|_\CV$ in $\wh \Picc (\CV/k)$. 
Here we take the convention $\CL|_X=(\CL|_\CU)|_X$, and
if $\CLL=(\CL, (\CX_i,\overline \CL_i, \ell_{i})_{i\geq 1})$ in $\wh\Picc (\CU/k)$, then $\CLL|_\CV=(\CL|_\CV, (\CX_i,\overline \CL_i, \ell_{i}|_\CV)_{i\geq 1})$ in 
$\wh\Picc (\CV/k)$.

By definition, $\wh\Picc (X/k)$ is a groupoid. 
We call objects of $\wh\Picc (X/k)$ \emph{adelic line bundles} on $X/k$.

As an easy consequence of Lemma \ref{isomorphism}, there is a canonical isomorphism
$$\wh\CaCl(X/k)\iso \wh\Pic(X/k)$$
for any flat and essentially quasi-projective integral scheme $X$ over $k$.

Let $\mathbf P$ represent one of the symbols in $\{\snef, \nef, \intb\}$.
Define
\begin{eqnarray*}
\wh\Div(X/k)_{\mathbf P}:&=&\lim_{\substack{\lra \\ \CU}}\wh \Div (\CU/k)_{\mathbf P},\\
\wh\CaCl(X/k)_{\mathbf P}:&=&\lim_{\substack{\lra \\ \CU}}\wh \CaCl (\CU/k)_{\mathbf P},\\
\wh\Pic (X/k)_{\mathbf P}: &=&\lim_{\substack{\lra \\ \CU}}\wh \Pic (\CU/k)_{\mathbf P},\\
\wh\Picc (X/k)_{\mathbf P}: &=&\lim_{\substack{\lra \\ \CU}}\wh \Picc (\CU/k)_{\mathbf P}.\end{eqnarray*}
Objects of $\wh\Picc (X/k)_{\snef}$ (resp. $\wh\Picc (X/k)_{\nef}$, $\wh\Picc (X/k)_{\intb}$) are called
\emph{strongly nef} (resp. \emph{nef},  \emph{integrable}) adelic line bundles on $X/k$.
Elements of $\wh\Div(X/k)_{\snef}$ (resp. $\wh\Div (X/k)_{\nef}$, $\wh\Div (X/k)_{\intb}$) are called
\emph{strongly nef} (resp. \emph{nef},  \emph{integrable}) adelic divisors on $X/k$.

In special situations, we take the following simplified or alternative notations:
\begin{enumerate}
\item The definitions also work for $X=\Spec F$ for a finitely generated field $F$ over $k$. 
We will write 
$$\wh \Pic (F/k)=\wh \Pic ((\Spec F)/k).$$ 
Apply this similarly to the other groups or categories.
\item If $k$ is minimal, i.e.  $k=\ZZ$ or $k=\FF_p$ for a prime $p$, then we may omit the dependence on $k$ in the groups or categories, as $k$ is determined by $X$ as an abstract scheme. In this case, we will simply write 
$$\wh \Pic (X)=\wh \Pic (X/k),\quad
\wh \Pic (F)=\wh \Pic (F/k).$$
This includes particularly the arithmetic case.
Apply this similarly to the other groups or categories.
\item
If $k$ is a field, we may also write 
$$
\wh \Div (X/k),\quad
\wh \CaCl (X/k),\quad
\wh \Picc (X/k),\quad
\wh \Pic (X/k)
$$
as 
$$
\wt \Div (X/k),\quad
\wt \CaCl (X/k),\quad
\wt \Picc (X/k),\quad
\wt \Pic (X/k).
$$
This is to emphasize that there is no archimedean component involved in the terms.
\end{enumerate}

We can compare our definition with that of Moriwaki \cite{Mor4} in the setting of projective varieties over finitely generated fields. 
Let $F$ be a finitely generated field over $k$, and $X$ be a projective variety over $F$.
Then our adelic line bundle on $X$ comes from an adelic line bundle $(\CL,(\CX_i,\overline \CL_i, \ell_{i})_{i\geq 1})$ on a quasi-projective model $\CU$ of $X$ over $k$. 
We will see that the sequence $\{(\CX_i, \overline \CL_i)\}_{i\geq 1}$ is close to the notion of an adelic sequence in \cite[\S3.1]{Mor4}.
In fact, by Lemma \ref{models}, we can shrink $\CU$ such that there is a projective and flat morphism $\CU\to \CV$ extending $X\to \Spec F$, where $\CV$ is a quasi-projective model of $\Spec F$ over $k$. 
We can further take a boundary divisor $(\CY,\OCE_0)$ of $\CV$ over $k$, and assume that there is a morphism $\CX_i\to \CY$ extending $X\to \Spec F$ for every $i\geq1$.
If $\CLL_i$ is nef, then $\{(\CX_i, \overline \CL_i)\}_{i\geq 1}$ is indeed an adelic sequence in the loc. cit.
Note that the loc. cit. defines a Cauchy condition in terms of both effectivity and intersection numbers, and we define the Cauchy condition purely in terms of effectivity.

\subsection{Forgetful maps} \label{sec forgetful}

Let $\CU$ be a quasi-projective variety over $k$. For any projective model $\CX$ of $\CU$, 
there are forgetful maps
$$
\wh\Div(\CX,\CU) \lra \Div(\CU),\quad\ 
\wh\Picc(\CX,\CU) \lra \Picc(\CU).
$$ 
Taking limits induces forgetful maps
$$
\wh\Div(\CU/k) \lra \Div(\CU),\quad\ 
\wh\Pic(\CU/k) \lra \Pic(\CU), \quad\ 
\wh\Picc(\CU/k) \lra \Picc(\CU).
$$ 
Here the last two maps send an object $\CLL=(\CL,(\CX_i,\overline \CL_i, \ell_{i}))$ to $\CL$. 

Let $X$ be a flat and essentially quasi-projective integral scheme over $k$.
Then the above maps induce forgetful maps 
$$
\wh\Div(X/k) \lra \Div(X),\quad\ 
\wh\Pic(X/k) \lra \Pic(X), \quad\ 
\wh\Picc(X/k) \lra \Picc(X).
$$ 
As a convention, we usually write an object of $\wh\Picc(X/k)$ in the form $\overline L$, where $L$ is understood to be the image of $\overline L$ in $\Picc(X)$.
We often refer $L$ as the \emph{underlying line bundle of} $\overline L$, and refer $\overline L$ as an \emph{adelic extension of} $L$.
We take similar conventions for $\wh\Pic(X/k)$ and $\wh\Div(X/k)$.

\subsection{Functoriality} \label{sec functoriality}

Here we introduce a few functorial maps between the Picard groups and the divisor groups. 
In the following, $\mathbf P$ represents one of the symbols in $\{\text{void}, \intb, \nef,\snef\}$, and take the convention that $\wh\Pic(X/k)_{\mathbf P}$ for ``$\mathbf P=\text{void}$'' means $\wh\Pic(X/k)$.

\medskip
\noindent \emph{Pull-back.} 
\kkk
Let $f: X'\to X$ be a morphism of flat and essentially quasi-projective integral schemes over $k$.
Then there are canonical maps 
$$f^*:\wh\Pic(X/k)_{\mathbf P} \lra \wh\Pic(X'/k)_{\mathbf P},$$
$$f^*:\wh\Picc(X/k)_{\mathbf P} \lra \wh\Picc(X'/k)_{\mathbf P}.$$

In fact, for quasi-projective models $X'\to \CU'$ and $X\to \CU$ over $k$, the rational map 
$\CU'\dashrightarrow \CU$ is defined in an open neighborhood of $X'$ in $\CU'$. Replacing $\CU'$ by that neighborhood if necessary, we obtain a morphism 
$f_\CU:\CU'\to \CU$. 
Then it suffices to define a canonical functor
$\wh\Picc(\CU/k) \to \wh\Picc(\CU'/k).$

Let $(\CL,(\CX_i,\CLL_{i}, \ell_{i}))$ be a Cauchy sequence in $\wh\Picc(\CU/k)_\rmod$.
There is a projective model $\CX_i'$ of $\CU'$ with a morphism $f_i:\CX_i'\to \CX_i$ extending $f_\CU:\CU'\to \CU$.
This can be achieved by taking any projective model $\CX_i'$ of $\CU'$ and blow-up $\CX_i'$ along a suitable center supported on $\CX_i'\setminus\CU'$.
Set the image of $(\CL,(\CX_i,\CLL_{i}, \ell_{i}))$ under $f^*$ to be 
$(f_\CU^*\CL,(\CX_i',f_i^*\CLL_{i}, f_\CU^*\ell_{i}))$.
To prove that the latter is a Cauchy sequence in $\wh\Picc(\CU'/k)_\rmod$, 
we need to compare the boundary topologies.

Let $(\CX_0,\OCE_0)$ (resp. $(\CX_0',\OCE_0')$) be a boundary divisor of $\CU$ (resp. $\CU'$) over $k$. As above, we can further assume that there is a morphism $f_0:\CX_0'\to \CX_0$ extending $\CU'\to \CU$. 
Note that $f_0^*\OCE_0$ is supported in $\CX_0'\setminus \CU'$.
As in our proof that the boundary topology is independent of the choice of the boundary divisor in \S\ref{sec adelic divisors}, there is a rational number $c>0$ such that $f_0^*\OCE_0\leq c\OCE_0'$.
This gives the compatibility of the boundary topologies. 

Hence, we have a functor $\wh\Picc(\CU/k) \to \wh\Picc(\CU'/k)$ and a functor $\wh\Picc(X/k) \to \wh\Picc(X'/k).$ The functor keeps tensor products.

In the above construction, if $f: X'\to X$ is dominant, there is also a canonical map
$$f^*:\wh\Div(X/k)_{\mathbf P} \lra \wh\Div(X'/k)_{\mathbf P}.$$

We will see in Corollary \ref{shrink} that these maps are injective if $X'$ and $X$ are normal and $f$ is birational.

\medskip
\noindent \emph{Varying the base.} 
Let $k'/k$ be a finitely generated extension of fields, and $X$ be an essentially quasi-projective integral scheme over $k'$. 
Then there are canonical maps
$$
\wh\Div(X/k)_{\mathbf P} \lra \wh\Div(X/k')_{\mathbf P},$$
$$\wh\Pic(X/k)_{\mathbf P} \lra \wh\Pic(X/k')_{\mathbf P},$$
$$\wh\Picc(X/k)_{\mathbf P} \lra \wh\Picc(X/k')_{\mathbf P}.$$
To define the maps, note that if $\CU$ (resp. $\CV$) is a quasi-projective model of $X$  (resp. $\Spec k'$) over $k$, then by shrinking $\CU$, we can assume that there is a flat morphism $\CU\to \CV$ extending $X\to \Spec k'$. This is similar to Lemma \ref{models}. 
Then it suffices to define a map $\wh\Div(\CU/k) \to \wh\Div(X/k')$
and it's analog for the line bundles.
By composition, we can further assume that $X$ is isomorphic to the generic fiber of $\CU\to \CV$.

Fix a projective model $\CB$ of $\CV$ over $k$.
For any projective model $\CX$ of $\CU$ over $k$, we can assume that there is a morphism $\CX\to \CB$ extending $\CU\to \CV$ by blowing-up $\CX$. 
Then the generic fiber $\CX_\eta$ of $\CX\to \CB$ is a projective model of $X$ over $k'$.
Finally, the map is induced by the natural map $\Div(\CX) \to \Div(\CX_\eta)$.

If $k'/k$ is a finite extension, the above maps are isomorphisms.

\medskip
\noindent \emph{Base change 1: geometric case.} 
Let $k'/k$ be a finitely generated extension of fields, and $X$ be an essentially quasi-projective integral scheme over $k$. 
Assume that the base change $X_{k'}$ is still integral. 
Then there are canonical maps
$$\wh\Div(X/k)_{\mathbf P} \lra \wh\Div(X_{k'}/k')_{\mathbf P},$$ 
$$\wh\Pic(X/k)_{\mathbf P} \lra \wh\Pic(X_{k'}/k')_{\mathbf P},$$
$$\wh\Picc(X/k)_{\mathbf P} \lra \wh\Picc(X_{k'}/k')_{\mathbf P}.$$
This is induced by the fact that, if $\CU$ is a quasi-projective (resp. projective) model of $X$ over $k$, then 
$\CU_{k'}$ is a quasi-projective (resp. projective) model of $X_{k'}$ over $k'$.
Then the maps are induced by the pull-back maps via the base changes.

\medskip
\noindent \emph{Base change 2: from $\ZZ$ to $\QQ$.} 
Let $X$ be a flat and essentially quasi-projective integral scheme over $\ZZ$.
For any projective  model $\CX$ of $X$ over $\ZZ$, the generic fiber $\CX_\QQ$ is a
projective model of $X_\QQ$ over $\QQ$. 
There are natural maps
$$
\wh\Div(\CX) \lra \Div(\CX_\QQ),\quad\ 
\wh\Pic(\CX) \lra \Pic(\CX_\QQ), \quad\ 
\wh\Picc(\CX) \lra \Picc(\CX_\QQ).
$$ 
These maps induce canonical maps
$$
\wh\Div(X/\ZZ)_{\mathbf P} \lra \wh\Div(X_\QQ/\QQ)_{\mathbf P},
\quad\OD\longmapsto \wt D,$$
$$
\wh\Pic(X/\ZZ)_{\mathbf P} \lra \wh\Pic(X_\QQ/\QQ)_{\mathbf P}, 
\quad\OL\longmapsto \wt L,
$$
$$
\wh\Picc(X/\ZZ)_{\mathbf P} \lra \wh\Picc(X_\QQ/\QQ)_{\mathbf P},
\quad\OL\longmapsto \wt L.
$$ 
We call $\wt D$ (resp. $\wt L$) the \emph{geometric part} of $\OD$ (resp. $\OL$) over $\QQ$.

\medskip
\noindent \emph{Base change 3: from $\ZZ$ to $\FF_p$.} 
Let $X$ be a flat and essentially quasi-projective integral scheme over $\ZZ$.
Let $p$ be a prime number such that the fiber $X_{\FF_p}$ of $X$ over $p$ is integral (and non-empty). 
For any projective model $\CX$ of $X$ over $\ZZ$, the Zariski closure $\CX_{\FF_p}'$ of $X_{\FF_p}$ in $\CX_{\FF_p}$ is a
projective model of $X_{\FF_p}$ over $\FF_p$. 
There are natural maps
$$
\wh\Div(\CX) \lra \Div(\CX_{\FF_p}'),\quad\ 
\wh\Pic(\CX) \lra \Pic(\CX_{\FF_p}'), \quad\ 
\wh\Picc(\CX) \lra \Picc(\CX_{\FF_p}').
$$ 
These maps induce canonical maps
$$
\wh\Div(X/\ZZ)_{\mathbf P} \lra \wh\Div(X_{\FF_p}/{\FF_p})_{\mathbf P},$$
$$
\wh\Pic(X/\ZZ)_{\mathbf P} \lra \wh\Pic(X_{\FF_p}/{\FF_p})_{\mathbf P}, 
$$
$$
\wh\Picc(X/\ZZ)_{\mathbf P} \lra \wh\Picc(X_{\FF_p}/{\FF_p})_{\mathbf P}.
$$

\medskip
\noindent \emph{Mixed situation.} 
Let $F$ be a finitely generated field over $\QQ$, and $X$ be a quasi-projective variety over $F$. Combining the above constructions, we obtain compositions.
$$
\wh\Div(X/\ZZ)_{\mathbf P} \lra \wh\Div(X/\QQ)_{\mathbf P}\lra \wh\Div(X/F)_{\mathbf P},
$$
$$
\wh\Pic(X/\ZZ)_{\mathbf P} \lra \wh\Pic(X/\QQ)_{\mathbf P} \lra 
 \wh\Pic(X/F)_{\mathbf P}.
$$ 
If $X$ is projective over $F$, then the compositions are just the forgetful maps
defined above.
In general, the image of an element of $\wh\Div(X/\ZZ)_{\mathbf P}$
(resp. $\wh\Pic(X/\ZZ)_{\mathbf P}$) in 
$\wh\Div(X/F)_{\mathbf P}$
(resp. $\wh\Pic(X/F)_{\mathbf P}$)
is called the \emph{geometric part} of this element over $F$.

\subsection{Extension to $\QQ$-coefficients} \label{sec rational}
\kkk
Let $X$ be a flat and essentially quasi-projective integral scheme over $k$.
To work with $\QQ$-line bundles on $X$, we write 
$$\wh\Pic(X/k)_{ \QQ}=\wh\Pic(X/k)\otimes_\ZZ\QQ,\qquad
\wh\Pic(X/k)_{\intb, \QQ}=\wh\Pic(X/k)_{\intb}\otimes_\ZZ\QQ.$$
We further write $\wh\Pic(X/k)_{\QQ,\snef}$ (resp. $\wh\Pic(X/k)_{\QQ,\nef}$) as the sub-semigroup of 
$\wh\Pic(X/k)_{ \QQ}$ consisting of \emph{positive} rational multiples of elements of 
$\wh\Pic(X/k)_{\snef}$ (resp. $\wh\Pic(X/k)_{\nef}$).
Extend the notations to $\wh\Div$ and $\wh\CaCl$ similarly. 

Let $\CU$ be a quasi-projective variety over $k$. Then we can interpret elements of the above groups directly in terms of Cauchy sequences.
In fact, by the isomorphism $\wh\Div(\CX,\CU)_\QQ\to \wh\Div(\CX)_\QQ$,
the group $\wh\Div(\CU/k)_{ \QQ}$ is simply isomorphic to the completion of 
$$
\wh\Div(\CU/k)_{\rmod, \QQ}=\varinjlim_\CX \wh\Div(\CX)_{\QQ}
$$
for the boundary topology defined similarly and 
$$
\wh\CaCl(\CU/k)_{ \QQ}\simeq 
\wh\Div(\CU/k)_{\rmod, \QQ}/ \wh\Pr(\CU/k)_{\QQ}.
$$

On the other hand, $\wh\Pic(\CU/k)_{ \QQ}$ is the group of isomorphism classes of
objects of a \emph{category} $\wh\Picc(\CU/k)_{ \QQ}$ defined as follows.
An object of $\wh\Picc(\CU/k)_{ \QQ}$ is a
sequence $(\CL,(\CX_i,\overline \CL_i, \ell_{i}))$, whose definition is similar to the integral case, except that  
$\CL$ is an object of $\Picc(\CU)_\QQ$ (instead of $\Picc(\CU)$). 
A morphism from an object $(\CL,(\CX_i,\overline \CL_i, \ell_{i}))$ to another object 
$(\CL', (\CX_i',\overline \CL_i', \ell_{i}'))$ is also similar to the integral case, except that it is given by an isomorphism $\CL\to \CL'$ of $\QQ$-line bundles (instead of integral line bundles) on $\CU$ which induces an isometry between the two objects in a similar sense.

Note that we do not derive $\wh\Picc(\CU/k)_{ \QQ}$ from the category $\wh\Picc(\CU/k)$ by multiplying a rational number to the objects, but they give equivalent categories. We choose the current definition for its simplicity. 

These descriptions can be used to define $\wh\Div(\CU/k)_{ \QQ}$, $\wh\Picc(\CU/k)_{ \QQ}$ and $\wh\Pic(\CU/k)_{ \QQ}$ without 
introducing their integral versions
$\wh\Div(\CU/k)$, $\wh\Picc(\CU/k)$ and $\wh\Pic(\CU/k)$ first. 
Moreover, we can even define the integral versions from the $\QQ$-versions, which can serve as a slightly different approach to the theory. 
For example, $\wh\Div(\CU/k)$ can be defined by the canonical exact sequence
$$
0\lra \wh\Div(\CU/k)
\lra \wh\Div(\CU/k)_{ \QQ} \oplus \Div(\CU)
\lra \Div(\CU)_{\QQ}.
$$
Here the last arrow sends a pair $(\OCD,\CD')$ to $\OCD|_\CU-\CD'$. 
Here $\OCD\mapsto\OCD|_\CU$ is the forgetful map as in \S\ref{sec forgetful}.

\section{Examples and constructions} \label{sec example}

Here we present some natural examples or constructions of adelic line bundles. 
Many constructions will be introduced in detail in the further sections, but we sketch them here to give readers more motivation for our theory of adelic line bundles.

\subsection{Arithmetic curves}

Let $K$ be a number field. 
Here we compute the group $\wh\Pic(K)$ of adelic line bundles.
Note that our definition in this case agrees with that in Zhang \cite{Zha2}.

\begin{lem}\label{arithmetic curve}
Let $\CU$ be an open subscheme of $\CX=\Spec O_K$. 
Denote 
$$
\wh\Pic^0(\CU) =\ker(\deg: \wh\Pic(\CU) \to \RR),
$$
$$
\wh\Pic^0(K) =\ker(\deg: \wh\Pic(K) \to \RR).
$$
There is a canonical exact sequence
$$
0\lra (\CO(\CU)^\times/\mu_K)\otimes_{\ZZ}(\RR/\ZZ) \lra \wh\Pic^0(\CU) \lra \Pic(\CU) \lra 0
$$ 
and a canonical isomorphism
$$
(K^\times/\mu_K)\otimes_{\ZZ}(\RR/\ZZ) \lra \wh\Pic^0(K).
$$ 
Here $\mu_K$ is the group of roots of unity in $K$.
\end{lem}

\begin{proof}

It suffices to prove the results for $\CU$. 
Denote $\CE=\CX\setminus \CU$, 
endowed with the reduced scheme structure.
Denote by $|\CX|, |\CU|, |\CE|$ the set of closed points of the corresponding schemes. 
Recall that
$$
\wh\Pic(\CU)
=\wh\Div(\CU) /\wh\Pr(\CU)_\rmod.$$ 
Note that $\CX$ is the only normal projective model of $\CU$.
We simply have 
$$
\wh\Div(\CU)_\rmod=\wh\Div(\CX,\CU),\quad\
\wh\Pr(\CU)_\rmod=\wh\Pr(\CX).
$$
Explicitly, 
$$
\wh\Div(\CX,\CU)=(\oplus_{v\in |\CU|} \ZZ )\oplus \QQ^{|\CE|}\oplus \RR^{ M_\infty},\quad\
\wh\Pr(\CX)= \wh\div_\CX(K^\times).
$$
Here 
$$
\wh\div_\CX: K^\times\lra \wh\Div(\CX)
$$ 
is the map of taking principal divisors. 

Take the arithmetic divisor 
$$
\OCE_0 = (\CE,1)=\sum_{v\in |\CE|\cup M_{\infty}} [v]. 
$$ 
It defines the boundary topology on $\wh\Div(\CX,\CU)$. 
The completion gives
$$
\wh\Div(\CU)=(\oplus_{v\in |\CU|} \ZZ)\oplus \RR^{|\CE|\cup M_\infty}.$$

The restriction map $\wh\Div(\CU)\to \Div(\CU)$ induces a canonical exact sequence
$$
0\lra  \wh\Pr(\CX)\cap \RR^{|\CE|\cup M_\infty}
\lra \RR^{|\CE|\cup M_\infty}  \lra \wh\Pic(\CU) \lra \Pic(\CU) \lra 0.
$$
It is easy to have 
$$
\wh\Pr(\CX)\cap \RR^{|\CE|\cup M_\infty}
\simeq \wh\div_\CX(\CO(\CU)^\times)
\simeq \CO(\CU)^\times/\mu_K.
$$
By Dirichlet's unit theorem, 
$\wh\div_\CX(\CO(\CU)^\times)$
is a full lattice of the hyperplane
$$
(\RR^{|\CE|\cup M_\infty})_0
:=\ker(\wh\deg:\RR^{|\CE|\cup M_\infty}\to \RR).
$$
This gives an isomorphism
$$
\wh\div_\CX(\CO(\CU)^\times)\otimes_\ZZ\RR = (\RR^{|\CE|\cup M_\infty})_0.
$$
The result follows.
\end{proof}

\subsection{Families of algebraic dynamical systems}

This example is our major motivation for introducing the theory of adelic line bundles. 
\kkk
Let the base $S$ be either one of the following:
\begin{enumerate}
\item[(a)] a quasi-projective variety over $k$,
\item[(b)] a quasi-projective variety over a field $F$ which is finitely generated field over $k$. 
\end{enumerate}
Let $(X, f, L)$ be a \emph{polarized dynamical system} over $S$; i.e.
\begin{enumerate}
\item $X$ is a flat projective scheme over $S$;
\item $f:X\to X$ is a morphism over $S$;
\item $L\in \Pic(X)_\QQ$ is a $\QQ$-line bundle, relatively ample over $S$, such that
$f^*L=qL$ from some rational number $q>1$.
\end{enumerate}
By Tate's limiting argument, there is a nef adelic line bundle $\overline L_f\in \wh\Pic(X/k)_{ \QQ}$ extending $L$ such that $f^*\overline L_f=q\overline L_f$.  

The construction is similar to the case that $S$ is the spectrum of a number field in \cite{Zha2}. We refer to Theorem \ref{invariant metric} for the current case. 

There are many naturally polarized dynamical systems over bases $S$ of positive dimensions. For example, this happens if $S$ is a moduli space, which includes the moduli space of endomorphisms of $\BP^N$ of fixed degree $d>1$ and a fine moduli space of polarized abelian varieties over $k$.

\subsection{Hodge bundles for Faltings heights} \label{sec hodge bundle}

In the proof of the Mordell conjecture by Faltings \cite{Fal2},  Faltings heights of abelian varieties over number fields are interpreted in terms of the height function associated with the Hodge bundle on certain moduli space of abelian varieties to deduce the Northcott property of Faltings heights. It is known that the Hodge bundle and the Faltings metric do not form a hermitian line bundle (on compactifications of the moduli space) due to the singularity of the metric at the boundary. Still, they do form an adelic line bundle in the framework of this book.
Here we describe the situation briefly and state the result for general families of abelian varieties instead of just moduli spaces.

The base ring here is $k=\ZZ$. 
Let $S$ be either one of the following:
\begin{enumerate}
\item a flat and quasi-projective integral scheme over $\ZZ$,
\item a quasi-projective variety over $\QQ$. 
\end{enumerate}
Our results hold for essentially quasi-projective $S$ by pull-back, but we restrict to the current cases for simplicity. 
Let $\pi: X\to S$ be a principally polarized abelian scheme of relative dimension $g$. 
Denote by 
$e:S\to X$ its identity section.

The \emph{Hodge bundle} on $S$ is the line bundle  
$$
\omega(S) = e^* \Omega_{X/S}^g\simeq \pi_*\Omega_{X/S}^g.$$
The \emph{Faltings metric} of $\omega(S)$ on $S(\CC)$ is defined by 
$$
\|\alpha\|_\Fal^2= \frac{1}{2^g} \left|\int_{X_s(\CC)} \alpha\wedge\bar \alpha\right|
= \frac{i^{g^2}}{2^g} \int_{X_s(\CC)} \alpha\wedge\bar \alpha
$$
for any point $s\in S(\CC)$ and any element $\alpha$ of the fiber
$$
\omega(S)(s)
=e_s^*\Omega_{X_s/ \CC}^g
\simeq
\Gamma(X_s, \Omega_{X_s/ \CC}^g).
$$

Then we have a metrized line bundle 
$(\omega(S), \|\cdot\|_\Fal)$
on $S$.
If $S$ is not projective over $\ZZ$, then it is not a hermitian line bundle in our strict sense. In general, it does not extend to a hermitian line bundle on a projective model of $S$ over $\ZZ$ due to the logarithmic singularity of the metric at the boundary. 
However, we will see that $(\omega(S), \|\cdot\|_\Fal)$ extends canonically to an adelic line bundle $\overline{\omega(S)}$ on $S/\ZZ$, and the height function associated to $\overline{\omega(S)}$ exactly computes the stable Faltings heights of the abelian varieties on fibers of $X\to S$. 

The precise statement and proof require the analytification and the height function of adelic line bundles, which will be introduced in the future sections, so we postpone the treatment to Theorem \ref{hodge bundle}.

\subsection{Hyperbolic metrics on families of curves} \label{sec hyperbolic metric}

The construction here is similar to the above setting of Hodge bundles, except that we work on families of curves instead of families of abelian varieties. The construction here is previously studied by Wolpert \cite{Wol} and Zhang \cite{Zha5}.

The base ring here is still $k=\ZZ$. 
Let $S$ be a flat and quasi-projective normal integral scheme over $\ZZ$ or $\QQ$. 
Let $\pi: X\to S$ be a smooth projective morphism whose fibers are geometrically integral curves of genus $g>1$. 

Consider the relative dualizing sheaf $\omega_{X/S}$ on $X$. 
The \emph{hyperbolic metric} $\|\cdot\|_\mathrm{hyp}$ of $\omega_{X/S}$ on $X(\CC)$ is defined as follows.
Endow the differential sheaf $\omega_{\CH}$ of the upper half plane $\CH$ the hyperbolic metric (also called the Petersson metric or the Poincare metric) normalized by
$$\|d z\|_\mathrm{hyp}=2\Im(z).$$ 
For any point $s\in S(\CC)$, via the universal covering $\CH\to X_s$, the hyperbolic metric of $\omega_{\CH}$ on $\CH$ descends to a hyperbolic metric of the line bundle $\omega_{X_s}$ on the Riemann surface $X_s.$
By the canonical isomorphism $\omega_{X/S}|_{X_s}\simeq \omega_{X_s}$ for varying $s\in S(\CC)$, 
we obtain the hyperbolic metric $\|\cdot\|_\mathrm{hyp}$ of $\omega_{X/S}$ on $X(\CC)$.  

Hence, we have a metrized line bundle 
$(\omega_{X/S}, \|\cdot\|_\mathrm{hyp})$
on $X$.
Our conclusion is that $(\omega_{X/S}, \|\cdot\|_\mathrm{hyp})$ extends canonically to an adelic line bundle 
$\overline\omega_{X/S,\mathrm{hyp}}$ on $X/\ZZ$. 
This result will be stated and briefly proved in Theorem \ref{hyperbolic metric}, as it is very similar to Theorem \ref{hodge bundle} for the Hodge bundle.

\subsection{Admissible metrics on families of curves} \label{sec admissible metrics}

In the above, have just extended $(\omega_{X/S}, \|\cdot\|_\mathrm{hyp})$ to an adelic line bundle. However, in Arakelov geometry, the hyperbolic metrics of curves are less used than the Arakelov metrics of curves. The goal here is to sketch a similar result for Arakelov metrics by Yuan \cite{Yua4}.

Recall that in Arakelov's original work \cite{Ara74}, on a compact Riemann surface $C$ of genus $g\geq 1$, there is an Arakelov metric $\|\cdot\|_\mathrm{Ar}$ on the canonical bundle $\omega_C$ defined to have a simple arithmetic adjunction formula. We refer to \cite[\S A.1]{Yua4} for a quick construction of the metric. 

Now  let us return to the global setting. 
The base ring here is still $k=\ZZ$. 
Let $S$ be a flat and quasi-projective integral scheme over $\ZZ$ or $\QQ$. 
Let $\pi: X\to S$ be a smooth projective morphism whose fibers are geometrically integral curves of genus $g\geq 1$.
Note that $g=1$ is allowed here.  

Consider the relative dualizing sheaf $\omega_{X/S}$ on $X$. 
By the canonical isomorphism $\omega_{X/S}|_{X_s}\simeq \omega_{X_s}$ for varying $s\in S(\CC)$, the Arakelov metric on $\omega_{X_s}$ patches together to form the  \emph{Arakelov metric} $\|\cdot\|_\mathrm{Ar}$ of $\omega_{X/S}$ on $X(\CC)$.  

Hence, we have a metrized line bundle 
$(\omega_{X/S}, \|\cdot\|_\mathrm{Ar})$
on $X$.
By \cite[Thm. 2.3]{Yua4},  $(\omega_{X/S}, \|\cdot\|_\mathrm{Ar})$ extends  to a canonically defined adelic line bundle 
$\overline\omega_{X/S,a}$ on $X/\ZZ$. 
The adelic line bundle 
$\overline\omega_{X/S,a}$
 has natural adjunction properties, and is called the \emph{admissible canonical bundle} of $X/S/\ZZ$.
Yuan \cite{Yua4} proves that $\overline\omega_{X/S, a}$ is nef and big on $X/\ZZ$ if the moduli map $S\to M_g$ is generically finite and applies this bigness to give a new proof of the uniform Bogomolov conjecture and the uniform Mordell conjecture previously proved by Vojta \cite{Voj}, Dimitrov--Gao--Habegger \cite{DGH}, and K\"uhne \cite{Kuh}.

In the case $S=\Spec K$ for a number field $K$, the admissible canonical bundle $\overline\omega_{X/S,a}$ is previously constructed by Zhang \cite{Zha3}.
We also refer to \S\ref{app subsec admissible} for a sketch of the construction of \cite{Zha3}.
The construction of \cite{Yua4} is a family version of \cite{Zha3}.

\subsection{Line bundles on Zariski--Riemann spaces}   \label{sec ZR}

Our model adelic line bundles can be realized on some generalized Zariski--Riemann space as introduced by Temkin \cite{Tem}. Here we recall the definitions and connections briefly. The treatment here will not be used in this book elsewhere.  

To illustrate the idea, we restrict to the geometric case that $k$ is a field. 
Let $X$ be an essentially quasi-projective integral scheme over $k$, as defined in \S\ref{sec models}. 
Define \emph{the Zariski--Riemann space associated to $X$} to be the ringed space
$$
\wt X=\varprojlim_\CX \CX,
$$
where the limit is over all projective models $\CX$ of $X$ over $k$. 
In the limit process, the underlying space $\wt X$ is endowed with the limit topology, i.e. the coarsest topology, so that all the projections $p_{\CX}:\wt X\to \CX$ to projective models $\CX$ of $X$ are continuous. 
The structure sheaf $\CO_{\wt X}$ is defined to be the direct limit of 
$p_{\CX}^{-1}\CO_\CX$ over all projective models $\CX$ of $X$. 

If $X=\Spec F$ for a finitely generated field $F$ over $k$, the space $\wt X$ is exactly the classical Zariski--Riemann space introduced by Zariski. 
In the general situation, the space $\wt X$ is the relative Zariski--Riemann space $\mathrm{RZ}_{X}(\CX_0)$ introduced by \cite[\S B2]{Tem}, here $\CX_0$ is a fixed projective model of $X$ over $k$. 

By definition, there is a canonical morphism $X\to\wt X$ induced by the morphism $X\to \CX$. Another key property is that $\wt X$ is quasi-compact; 
see \cite[Prop. B.2.3]{Tem}.

Since $\wt X$ is a ringed space, coherent sheaves, invertible sheaves, and Cartier divisors are defined on $\wt X$.
Then we can still define line bundles on $\wt X$ to be invertible sheaves. 
By the quasi-compactness of $\wt X$, we can prove that the natural maps
$$
\varinjlim_\CX \Div(\CX) \lra \Div(\wt X),
$$
$$
\varinjlim_\CX \Pic(\CX) \lra \Pic(\wt X)
$$ 
are isomorphisms.
If $X=\Spec F$, the first limit is the group of \emph{Carter b-divisors} introduced by Shokurov \cite{Sho} in the minimal model program. 

To connect to our adelic divisors, we see that 
$$
\Div(X)_{\rmod,\QQ}
=\varinjlim_\CX \Div(\CX)_\QQ
$$
is canonically isomorphic to $\Div(\wt X)_\QQ$.
Then the group $\wh\Div(X/k)_\QQ$ is a suitable completion of $\Div(\wt X)_\QQ$. 
Similarly, the group $\wh\Div(X/k)$ is a suitable completion of the mixed divisor group 
$$
\Div(\wt X,X):=\ker(\Div(\wt X)_\QQ\oplus \Div(X)\to \Div(X)_\QQ). 
$$
Here the arrow sends $(\wt D, D)$ to $\wt D|_X-D$ as before.

\section{Definitions over more general bases}  \label{sec adelic general}

The above theory of $\wh\Div(X/k)$ and $\wh\Picc(X/k)$, when $k$ is either $\ZZ$ or a field, covers all the global situations we are interested in, but it does not include the local situation that $k$ comes from a local field. 
Moreover, the definitions can also be generalized by replacing $k$ with a general Dedekind scheme. 
The goal of this section is to sketch the treatment in these situations, which also includes the function field case.  

The exposition here is very similar to the previous case. Still, we use it sparingly throughout this book to avoid the extra burden of terminology and potential confusion of cases. 
The setup here is only restricted to this section and partly to \S\ref{sec local theory} and \S\ref{sec local pairing}.

\subsection{Valuations}

By a \emph{valuation} over a field $K$, we mean a map $|\cdot |:K\to \RR$ satisfying the following properties:
\begin{enumerate}
\item (positivity) $|0|=0$, and $|a|>0$ for any $a\in K^\times$. 
\item (triangle inequality) $|a+b|\leq |a|+|b|$ for any $a,b\in K$. 
\item (multiplicativity) $|ab|= |a|\cdot |b|$ for any $a,b\in K$. 
\end{enumerate}
The valuation is \emph{trivial} if $|a|=1$ for any $a\in K^\times$. 
The valuation is \emph{archimedean} (resp. \emph{non-archimedean}) if $|n|$ is unbounded (resp. bounded) for all $n\in \ZZ$ viewed as elements of $K$ under the natural map $\ZZ\to K$.
If $|\cdot|$ is non-archimedean, the \emph{valuation ring} of $K$ is $O_K:=\{x\in K:|x|\leq 1\}$. 

By a  \emph{non-archimedean field $(K,|\cdot|)$}, we mean a field $K$ endowed with a complete non-archimedean non-trivial valuation $|\cdot|$.

\subsection{Base valued schemes} 

Recall that an integral domain is called a \emph{Pr\"ufer domain} if all of its local rings are valuation rings. This is a classical term widely studied in commutative algebra. We refer to \cite[p. 558, Ch. VII, \S2, Ex. 12]{Bou} for 14 equivalent definitions of Pr\"ufer domains, and refer to \cite{BG} for more properties and history of the concept. 

It is easy to see that a Pr\"ufer domain is a Dedekind domain if and only if it is noetherian. 
Thus, Pr\"ufer domains can be viewed as a non-noetherian generalization of Dedekind domains, and thus, many nice properties of Dedekind domains also hold for Pr\"ufer domains.
For example, a module over a Pr\"ufer domain is flat if and only if it is torsion-free, which can be checked by taking localizations.

A quasi-compact integral scheme is called a \emph{Pr\"ufer scheme} if all of its local rings are valuation rings.
We introduce this concept to include the following three important classes:
\begin{enumerate}
\item $\Spec k$ for a field $k$;
\item a Dedekind scheme, i.e. a regular and integral noetherian scheme of dimension 1;
\item $\Spec O_K$, where $K$ is a non-archimedean field and $O_K$ is the valuation ring. 
\end{enumerate}
As a consequence of the flatness property, a reduced scheme of $X$ over a 
Pr\"ufer scheme $B$ is flat over $B$ if and only if every irreducible component of $X$ has a Zariski dense image in $B$.

By a \emph{base valued scheme}, we mean a pair $\OB=(B,\Sigma)$ consisting of a Pr\"ufer scheme $B$ and a subset $\Sigma$ of $\Hom(K,\CC)$, where $K$ denotes the function field of $B$. The set $\Sigma$ is allowed to be empty, in which case we get a scheme $\OB=B$. 

Note that every $\sigma\in \Sigma$ induces an archimedean valuation $|\cdot|_\sigma$ over $K$. We may also think $\OB$ as $(B,\{|\cdot|_\sigma\}_{\sigma\in \Sigma})$, but note that $|\cdot|_\sigma=|\cdot|_{\sigma'}$ if and only if $\sigma'=c\circ\sigma$ for the complex conjugate $c:\CC\to \CC$.

We introduce this definition to include the following important and natural types of base schemes:
\begin{enumerate}
\item (geometric case) $\Spec k$, where $k$ is a field;
\item (number field case) $(\Spec O_K, \Hom(K,\CC))$, where $K$ is a number field;
\item (function field case) a projective and geometrically integral regular curve $B$ over a field $k$;
\item (archimedean case) $(\Spec \RR, i_\mathrm{st})$ or $(\Spec \CC, \mathrm{id})$, where 
$i_\mathrm{st}:\RR\to \CC$ is the standard injection and $\mathrm{id}:\CC\to \CC$ is the identity map;
\item (non-archimedean case) $\Spec O_K$, where $K$ is a non-archimedean field and $O_K$ is the valuation ring;
\item (Dedekind case) a Dedekind scheme $B$.
\end{enumerate}
Note that $\Sigma=\emptyset$ in cases (1), (3), (6), and $\Sigma$ is finite in all the cases.
We usually write $K$ for the function field of $B$. 

Case (1) with any $k$ and case (2) with $K=\QQ$ are exactly our original case $k=\ZZ$ or $k$ is a field in \S\ref{sec uniform}. 
The Dedekind case includes the function field case, but we list the function case separately for its importance. 

Let $B$ be a Pr\"ufer scheme. 
By an \emph{arithmetic variety} over $B$, we mean an integral scheme that is flat, separated, and of finite type over $B$.
By a \emph{quasi-projective arithmetic variety} (resp. \emph{projective arithmetic variety}) over $B$, we mean an arithmetic variety over $B$ which is quasi-projective (resp. projective) over $B$.
For a quasi-projective arithmetic variety $\CU$ over $B$, a \emph{projective model} means a projective arithmetic variety $\CX$ over $B$ endowed with an open immersion $\CU\to \CX$ over $B$. 

As in \S\ref{sec models}, a flat integral scheme $X$ over $B$ is \emph{essentially quasi-projective over $B$} if there is a \emph{pro-open} immersion $X\to \CX$ to a projective arithmetic variety $\CX$ over $B$. 
A \emph{quasi-projective model} (resp. \emph{projective model}) of $X$ means a quasi-projective (resp. projective) arithmetic variety $\CU$ over $B$ endowed with a pro-open immersion $X\to \CU$ over $B$.

\subsection{Model adelic divisors and adelic line bundles}

Let $\OB=(B,\Sigma)$ be a base valued scheme. 
Let $\CX$ be a projective arithmetic variety over $B$.
We define arithmetic divisors and hermitian line bundles on $\CX$ as follows.

An \emph{arithmetic divisor} on $\CX$ is a pair $(\CD,g_\CD)$, where $\CD$ is a Cartier divisor on $\CX$, and $g_\CD:\CX_{\Sigma}(\CC)\setminus |\CD_{\Sigma}(\CC)|\to \RR$ is a Green function of continuous type of $\CD_\Sigma(\CC)$ on $\CX_\Sigma(\CC)$ as in \S\ref{sec Green}. 
Here
$\CX_{\Sigma}(\CC):=\coprod_{\sigma\in\Sigma}\CX_{\sigma}(\CC)$
is a projective analytic variety, and the 
Cartier divisor $\CD_\Sigma(\CC)$ on $\CX_\Sigma(\CC)$ is defined similarly.
By restriction, we have a Green function $g_{\CD,\sigma}:\CX_{\sigma}(\CC)\setminus |\CD_{\sigma}(\CC)|\to \RR$ of $\CD_{\sigma}(\CC)$ on $\CX_{\sigma}(\CC)$. Then we can also think of $g_\CD$ as a collection of $g_{\CD,\sigma}$ over $\sigma\in\Sigma$.

The Green function $g_\CD$ is further required to be \emph{invariant under the complex conjugate $c:\CC\to \CC$} in the sense that for any $\sigma\in \Sigma$ such that $\bar \sigma=c\circ\sigma\in \Sigma$, we require $g_{\CD,\bar\sigma}=g_{\CD,\sigma}\circ c$.

The divisor $(\CD,g_\CD)$ is \emph{effective} (resp. \emph{strictly effective}) if $\CD$ is an effective Cartier divisor on $\CX$, and $g_\CD\geq 0$ (resp. $g_\CD>0$). 

A \emph{principal arithmetic divisors} on $\CX$ is an arithmetic divisor of the form 
$$\wh\div(f):=(\div(f), -\log|f|)$$ 
for some $f\in K(\CX)^\times$. 

A \emph{hermitian line bundle} on $\CX$ is a pair $(\CL,\|\cdot\|)$, where $\CL$ is a line bundle on $\CX$, and $\|\cdot\|$ is a continuous hermitian metric of 
$\CL_\Sigma(\CC)$ on $\CX_\Sigma(\CC)$ as in \S\ref{sec complex}. 
As above, the metric $\|\cdot\|$ is equivalent to a collection of continuous metrics $\|\cdot\|_\sigma$ of the line bundle $\CL_\sigma(\CC)$ on $\CX_\sigma(\CC)$ over $\sigma\in \Sigma$.
The metric is also required to be \emph{invariant under the complex conjugate} in the above sense. 

Now we have the following groups (or category)
$$
\wh\Div(\CX), \quad
\wh\Pr(\CX), \quad
\wh\Pic(\CX), \quad
\wh \Picc(\CX), \quad
\wh\Picc(\CX)_\QQ.
$$
Here $\wh\Div(\CX)$ (resp. $\wh\Pr(\CX)$) is the group of
arithmetic divisors (resp. principal arithmetic divisors) on $\CX$.
And $\wh\Picc(\CX)$ (resp. $\wh\Pic(\CX)$) is the category (resp. group) of
{hermitian line bundles} on $\CX$ under isometry.
The category $\wh\Picc(\CX)_\QQ$ is defined from $\wh\Picc(\CX)$ similar to that in \S\ref{sec convention}. 

Define
$$
\wh\CaCl(\CX):=\wh\Div(\CX)/\wh\Pr(\CX).
$$
Then there is a canonical isomorphism
$$
\wh\CaCl(\CX)
\lra
\wh \Pic(\CX).
$$

If $\Sigma=\emptyset$, and thus $\OB=B$ is a scheme, then the above groups are just the usual ones
$$
\Div(\CX), \quad
\Pr(\CX), \quad
\CaCl(\CX),\quad
\Picc(\CX), \quad
\Pic(\CX).
$$

Let $\CU$ be an open subscheme of $\CX$.
As in \S\ref{sec convention},  
we also have the groups of objects of $(\QQ,\ZZ)$-coefficients:
$$
\wh\Div(\CX,\CU), \quad
\wh\CaCl(\CX,\CU).
$$ 
For example, $\wh\Div(\CX,\CU)$ is the fiber product of the natural map 
$\wh\Div(\CX)_\QQ\to \Div(\CU)_\QQ$ with the natural map $\Div(\CU)\to \Div(\CU)_\QQ$, whose elements are pairs $(\OCD, \CD')$, called \emph{arithmetic $(\QQ,\ZZ)$-divisors on $(\CX,\CU)$},
where $\OCD\in \wh\Div(\CX)_\QQ$ and $\CD'\in \Div(\CU)$ have equal images in 
$\Div(\CU)_\QQ$.

An element of $\wh\Div(\CX,\CU)$ is called \emph{effective} if its images in $\wh\Div(\CX)_\QQ$ and $\Div(\CU)$ are both effective.
The effectivity induces a \emph{partial order} on $\wh\Div(\CX,\CU)$ as before.

\subsection{Adelic divisors on a quasi-projective variety}

Let $\CU$ be a quasi-projective arithmetic variety over $B$.
Using pull-back morphisms, define 
$$\wh\Div (\CU/\OB)_\rmod:=\lim_{\substack{\lra\\ \CX}}\wh\Div(\CX,\CU),
\qquad \wh\Pr (\CU/\OB)_\rmod:=\lim_{\substack{\lra\\ \CX}}\wh\Pr(\CX).$$
Here the limits are over projective models $\CX$ of $\CU$ over $B$.

The direct limit is filtered, i.e. for any two projective models $\CX_1, \CX_2$ of $\CU$ over $B$, there is a third projective model $\CY$ of $\CU$ over $B$ dominating $\CX_1, \CX_2$ in the sense that there are morphisms $\CY\to \CX_1$ and $\CY\to \CX_2$ of projective models of $\CU$ over $B$.
It suffices to take $\CY$ to be the Zariski closure of the image of the composition
$\CU \to \CU\times_{B}\CU \to \CX_1\times_{B}\CX_2$,
where the first map is the diagonal map. 
Here $\CY$ is flat over $B$ by the Pr\"ufer property.

The partial order in $\wh \Div (\CX,\CU)$ induces a partial order in $\wh\Div (\CU/\OB)_\rmod$ by the limit process.

Let $(\CX_0,\OCE_0)$ be \emph{a boundary divisor of $\CU$ over $B$}; that is, $\CX_0$ is a projective model of $\CU$ over $B$, $\OCE_0$ is a strictly effective arithmetic divisor on $\CX_0$ with support $|\CE_0|=\CX_0\setminus \CU$. 
This gives an \emph{extended norm} 
$$\|\cdot\|_{\OCE_0}:\wh\Div (\CU/\OB)_\rmod
\lra [0,\infty]$$
defined by 
$$
\|\OCD\|_{\OCE_0}:=\inf\{\epsilon\in \BQ_{>0}: \ 
 -\epsilon \OCE_0 \leq
\OCD \leq  \epsilon \OCE_0\}.
$$
Here the inequalities are again defined by the effectivity of divisors,
and we take the convention that $\inf(\emptyset)=\infty$.
The \emph{boundary topology on $\wh\Div (\CU/\OB)_\rmod$} is the topology over
$\wh\Div (\CU/\OB)_\rmod$ induced by the
extended norm
$\|\cdot\|_{\OCE_0}$.
Thus, a neighborhood basis at $0$ of the topology is given by
$$
 B(\epsilon, \wh\Div (\CU/\OB)_\rmod):
 =\{\OCD \in \Divhat(\CU/\OB)_\rmod: \ 
 -\epsilon \OCE_0 \leq
\OCD \leq  \epsilon \OCE_0\}, \quad \epsilon\in \BQ_{>0}.
$$
By translation, it gives a neighborhood basis at any point.

Let $\wh \Div  (\CU/\OB)$ be the \emph{completion} of $\wh \Div  (\CU/\OB)_\rmod$ for the  boundary topology. 
An element of $\wh \Div  (\CU/\OB)$ is called an \emph{adelic divisor} (or a \emph{compactified divisor}) on $\CU$.
Define the \emph{class group of adelic divisors} of $\CU$ to be 
$$\wh\CaCl (\CU/\OB):=\wh\Div  (\CU/\OB)/\wh \Pr (\CU/\OB)_\rmod.$$

\subsection{Adelic line bundles on a quasi-projective variety}

Let $\CU$ be a quasi-projective variety over $B$. Let $(\CX_0,\OCE_0)$ be as above.
Define the \emph{category $\wh\Picc (\CU/\OB)$ of adelic line bundles} on $\CU$ as follows.
An object of $\wh\Picc (\CU/\OB)$ is a pair
$(\CL, (\CX_i,\overline \CL_i, \ell_{i})_{i\geq 1})$ where:
\begin{enumerate}
\item $\CL$ is an object of $\Picc(\CU)$, i.e. a line bundle on $\CU$;

\item  $\CX_i$ is a projective model of $\CU$ over $B$;

\item  $\overline \CL_i$ is an object of $\wh\Picc(\CX_i)_\QQ$, i.e. a hermitian $\QQ$-line bundle on $\CX_i$;

\item $\ell_i:\CL\to \CL_i|_{\CU}$ is an isomorphism in $\Picc(\CU)_\QQ$.
\end{enumerate}
Similar to \S\ref{sec adelic line bundles},  the sequence is required to satisfy the \emph{Cauchy condition} that
the sequence $\{\wh \div(\ell_i \ell_1^{-1})\}_{i\geq 1}$ is a 
Cauchy sequence in $\wh\Div(\CU/\OB)_\rmod$ under the boundary topology.

A morphism from an object $(\CL, (\CX_i,\overline \CL_i, \ell_{i})_{i\geq 1})$ of $\wh\Picc (\CU/\OB)$ to another object
$(\CL',(\CX_i',\overline \CL_i', \ell_{i}')_{i\geq 1})$ of $\wh\Picc (\CU/\OB)$ is an isomorphism  $\iota:\CL\to \CL'$ of the integral line bundles on $\CU$ satisfying the following properties.
Denote by $\iota_1:\overline \CL_1\dashrightarrow \overline \CL_1'$ the rational map on $\CU$ induced by $\iota$, which induces an element $\wh\div(\iota_1)$
of $\wh \Div (\CU/\OB)_\rmod$.
Then we require that the sequence 
$\{\wh\div(\ell_{i}'\ell_1'^{-1})-\wh\div(\ell_{i}\ell_1^{-1})+\wh\div(\iota_1) \}_{i\geq1}$
of $\wh \Div (\CU/\OB)_\rmod$ converges to 0 in $\wh \Div (\CU/\OB)$
under the boundary topology.

An object of $\wh\Picc (\CU/\OB)$ is called an \emph{adelic line bundle} (or a \emph{compactified line bundle}) on $\CU$.
Define $\wh\Pic (\CU/\OB)$ to be the \emph{group} of isomorphism classes of objects of $\wh\Picc (\CU/\OB)$. 
As before,  there is a canonical isomorphism
$$\wh\CaCl(\CU/\OB)\lra \wh\Pic(\CU/\OB).$$

\subsection{Definitions on essentially quasi-projective schemes}

Let $\OB=(B,\Sigma)$ be a base valued scheme. 
Let $X$ be a flat and essentially quasi-projective integral scheme over $B$.
Define
$$\wh \Div (X/\OB):=\lim _{\substack{\lra \\ \CU}}\wh \Div (\CU/\OB),$$
$$\wh \CaCl (X/\OB):=\lim _{\substack{\lra \\ \CU}}\wh \CaCl (\CU/\OB),$$
$$\wh \Picc (X/\OB):=\lim _{\substack{\lra \\ \CU}}\wh \Picc (\CU/\OB),$$
$$\wh \Pic (X/\OB):=\lim _{\substack{\lra \\ \CU}}\wh \Pic (\CU/\OB).$$
An element of $\wh \Div (X/\OB)$ is called an \emph{adelic divisor on $X/\OB$}.
An object of $\wh \Picc (X/\OB)$ is called an \emph{adelic line bundle on $X/\OB$}.

We take the following alternative notations:
\begin{enumerate}
\item
If $\Sigma=\emptyset$ and thus $\OB=B$ is a Pr\"ufer scheme, we may also write 
$$
\wh \Div (X/\OB),\quad
\wh \CaCl (X/\OB),\quad
\wh \Picc (X/\OB),\quad
\wh \Pic (X/\OB)
$$
as 
$$
\wt \Div (X/B),\quad
\wt \CaCl (X/B),\quad
\wt \Picc (X/B),\quad
\wt \Pic (X/B).
$$
This is to emphasize that there is no archimedean component involved in the terms.
If $B=\Spec R$ is affine, they are further written as
$$
\wt \Div (X/R),\quad
\wt \CaCl (X/R),\quad
\wt \Picc (X/R),\quad
\wt \Pic (X/R).
$$

\item
If $\OB=(\Spec O_K, \Hom(K,\CC))$ is in the arithmetic case, we may also write  
$$
\wh \Div (X/\OB),\quad
\wh \CaCl (X/\OB),\quad
\wh \Picc (X/\OB),\quad
\wh \Pic (X/\OB)
$$
as 
$$
\wh \Div (X/O_K),\quad
\wh \CaCl (X/O_K),\quad
\wh \Picc (X/O_K),\quad
\wh \Pic (X/O_K).
$$
We take a similar notation for the archimedean case $\OB=(\Spec \RR, i_\mathrm{st})$ or $\OB=(\Spec \CC, \mathrm{id})$.
\end{enumerate}

To compare the notations with the original setting ($k=\ZZ$ or a field), we have the following. 
\begin{enumerate}
\item[(a)]
If $k$ is a field, the current term $\wh\Div(X/B)$ with $B=\Spec k$ is the same as the original term $\wh\Div(X/k)$. 
\item[(b)]
If $k=\ZZ$, the current term $\wh\Div(X/\OB)$ with $\OB=(\Spec \ZZ,\Hom(\QQ,\CC))$ is the same as the original term $\wh\Div(X/\ZZ)$. 
On the other hand, the term $\wt\Div(X/\ZZ)$ removes the Green functions from the arithmetic case. Then we have naturally forgetful maps
$$\wh\Div(X/\ZZ)\lra \wt\Div(X/\ZZ),\qquad
\wh\Picc(X/\ZZ)\lra \wt\Picc(X/\ZZ).$$
They are surjective (or essentially surjective).
\item[(c)]
If $K$ is a number field, for any flat and essentially quasi-projective integral scheme $X$ over $O_K$, we have canonical isomorphisms
$$\wh\Div(X/O_K)\lra \wh\Div(X/\ZZ),\qquad
\wh\Picc(X/O_K)\lra \wh\Picc(X/\ZZ).$$
This follows from the fact that a scheme over $O_K$ is projective (resp. flat) over $O_K$ if and only if it is projective (resp. flat) over $\ZZ$.
Therefore, our original approach essentially includes this case.
\end{enumerate}

For any base valued scheme $\OB=(B,\Sigma)$, there are canonical forgetful maps 
$$
\wh\Div(X/\OB)\lra \wt\Div(X/B) \lra \Div(X), 
$$
and 
$$
\wh\Picc(X/\OB) \lra \wt\Picc(X/B) \lra \Picc(X).
$$ 
These are induced by the forgetful functor
$$
\wh\Picc(\CU/\OB) \lra \wt\Picc(\CU/B) \lra \Picc(\CU),$$
given by
$$ 
(\CL,(\CX_i, \CLL_i, \ell_{i}))\longmapsto (\CL,(\CX_i, \CL_i, \ell_{i}))\longmapsto \CL.
$$ 
As a convention, our notation for the three objects is usually denoted by
$$ 
\OL\longmapsto \TL \longmapsto L.
$$ 
We often refer $L$ as the \emph{underlying line bundle of} $\OL$ and $\TL$, and refer $\OL$ as an \emph{adelic extension of} $L$.

\subsection{The theory over function fields}

Arakelov geometry is analogous to algebraic geometry over fields, which is why this book uses uniform terminology. 
However, Arakelov geometry is more analogous to algebraic geometry over a projective curve. 
Here we explore this analog briefly. 

Let $k$ be a field and $B$ be a projective and regular curve over $k$. 
Denote by $K=k(B)$ the function field of $B$. 
In the above perspective, the counterpart of the arithmetic object $\wh\Div(\cdot /\ZZ)$ should be the geometric object $\wt\Div(\cdot /B)$.

Let $X$ be flat and essentially a quasi-projective integral scheme over $B$.
Then $X$ is also essentially quasi-projective over $k$.
We claim that there are canonical isomorphisms
$$
\wt\Div(X/B) \lra \wt\Div(X/k),
$$
$$
\wt\CaCl(X/B) \lra \wt\CaCl(X/k),
$$
$$
\wt\Pic(X/B) \lra \wt\Pic(X/k),
$$
$$
\wt\Picc(X/B) \lra \wt\Picc(X/k).
$$
In this sense, we do not lose much in our original setup by considering the objects over the absolute base field $k$.

To see the isomorphisms, note that any quasi-projective (resp. projective) model of $X$ over $B$ is a quasi-projective (resp. projective) model of $X$ over $k$. 
Moreover, for any quasi-projective (resp. projective) model $\CU$ of $X$ over $k$, 
the rational map $\CU\dasharrow B$ is defined along $X$ and can be turned to a morphism by shrinking $\CU$ (resp. blowing-up $\CU$ along a center disjoint from $X$). 
Therefore, the inverse systems of quasi-projective (resp. projective) models of $X$ over $B$ are cofinal to the inverse system of quasi-projective (resp. projective) models of $X$ over $k$.

\chapter{Interpretation by Berkovich spaces} \label{sec analytic}

\kkk
Let $X$ be a flat and essentially quasi-projective integral scheme over $k$, as defined in \S\ref{sec models}. 
In \S \ref{sec adelic line bundles}, we have introduced the category $\wh\Picc(X/k)$ of adelic line bundles on $X$. The goal of this chapter is to introduce a category 
$\wh\Picc(X^\an)$ of metrized line bundles on the Berkovich analytic space $X^\an$ associated to $X$, and study the analytification functor from $\wh\Picc(X/k)$ to $\wh\Picc(X^\an)$. 
The analytification functor is fully faithful and thus provides a convenient interpretation of adelic line bundles. 
This generalizes the work of \cite{Zha2} for projective varieties over number fields.
We refer to \S\ref{app sec adelic} and \S\ref{app sec equidistribution} for the projective case.

\section{Berkovich spaces}   \label{sec Berkovich}

In this section, we review definitions and some basic properties of Berkovich spaces.
In the end, we introduce a density result that will be useful in analytification of adelic divisors and adelic line bundles.

\subsection{Generality on Berkovich spaces}

Berkovich spaces are best known as analytic spaces associated with varieties over non-archimedean fields, whose foundation was introduced by Berkovich \cite{Ber1}. 
By Berkovich \cite[\S1]{Ber2}, the base fields are relaxed to be Banach rings, and the old construction works similarly. 
In the following, we recall the construction of \cite[\S1]{Ber2} to adapt our setting so that the schemes are not required to be of finite type. 

Let $k$ be a commutative Banach ring with unity 1. 
Let $X$ be a scheme over $k$. 
In the following, we recall the definition and basic properties of the \emph{Berkovich space $X^\an$ associated to} $X$, which is more rigorously written as $(X/k)^\an$ to emphasize the dependence on $k$.
\begin{enumerate}
\item \textbf{Affine case.}
If $X=\Spec A$, then $X^\an$ is defined to be the space $\CM(A)=\CM(A/k)$ of multiplicative
semi-norms on $A$ whose restriction to $k$ is bounded by $|\cdot|_\mathrm{Ban}$. For each $x\in \CM(A)$, 
denote its corresponding semi-norm on $A$ by $|\cdot|_x:A\to \RR$.
For any $f\in A$, write $|f|_x$ as $|f(x)|$, which gives a real-valued function $|f|$ on $\CM(A)$.
The topology on $\CM(A)$ is the weakest one such that the function $|f|:\CM(A)\to\RR$ is continuous for all $f\in A$. 

\item \textbf{General case.}
If $X$ is covered by an affine open cover $\{\Spec A_i\}_i$, then $X^\an$ is defined to be the union of $\CM(A_i)$, glued canonically. 
The topology of $X^\an$ is the weakest one such that each $\CM(A_i)$ is an open subspace of $X^\an$.

\item \textbf{Residue field.}
For each $x\in \CM(A)$, the corresponding semi-norm $|\cdot|_x$ induces a norm on the integral domain $A/\ker(|\cdot|_x)$. The completion of the fraction field of  $A/\ker(|\cdot|_x)$ is called the \emph{residue field} of $x$ and denoted by $H_x$. Denote by $|\cdot|$ the valuation (multiplicative norm) on $H_x$ induced by $|\cdot|_x$. 
Then $|\cdot|_x:A\to \RR$ is equal to the composition 
$$A\lra H_x\overset{|\cdot|}{\lra} \BR.$$
We write the first map as $f\mapsto f(x)$, which is compatible with the convention $|f|_x=|f(x)|$.
The notation $H_x$ generalizes to any scheme $X$ over $k$.

\item \textbf{Contraction.}
There is a canonical contraction map $\kappa: X^\an\to X$. It suffices to describe it in the case $X=\Spec A$. 
For each $x\in \CM(A)$, 
the kernel of the map $|\cdot|_x:A\to \RR$ is a prime ideal of $A$, and thus defines an element $\kappa(x)\in \Spec A$.

\item \textbf{Injection.}
Assume that for any $x\in \Spec k$, the semi-norm $|\cdot|_{x,0}$ on $k$, induced by the trivial norm on the residue field $k/x$, is bounded by $|\cdot|_\mathrm{Ban}$.
This gives a natural injection $\iota:\Spec k\to \CM(k)$ by sending $x$ to $|\cdot|_{x,0}$.

Under this assumption, there is a natural injection $\iota: X\to X^\an$
defined similarly. It suffices to describe it in the case $X=\Spec A$. 
For any $x\in \Spec A$, still denote by $|\cdot|_{x,0}$ the semi-norm 
on $A$ induced by the trivial norm on the residue field $A/x$. 
Then $\iota:X\to X^\an$ sends $x$ to $|\cdot|_{x,0}$.
The $\kappa\circ\iota$ is the identity map on $X$.

\item \textbf{Functoriality.}
Any morphism $f:X\to Y$ over $k$ induces a continuous map $f^\an:X^\an\to Y^\an$. 
For any point $v\in Y^\an$, the \emph{fiber} $X_v^\an=(f^\an)^{-1}(v)$, defined as a subspace of $X^\an$, is canonically homeomorphic to the Berkovich space $(X_{H_v}/H_v)^\an$.
More generally, for any subset $T\subset Y^\an$, denote by $X_T^\an$ the preimage of $T$, viewed as a subspace of $X^\an$. 
This notation automatically applies to the case $Y=\Spec k$ and $Y^\an=\CM(k)$.
\end{enumerate}
By \cite[Lem. 1.1, Lem. 1.2]{Ber2}, we have the following basic topological properties: 
\begin{enumerate}
\item If $X$ is separated and of finite type over $k$, then $X^\an$ is Hausdorff.
\item If $X$ is of finite type over $k$, then $X^\an$ is locally compact.
\item If $X$ is projective over $k$, then $X^\an$ is compact.
\end{enumerate}

It is well-known that Berkovich spaces also include complex analytic spaces coming from algebraic varieties. 
\begin{enumerate}
\item If $k=\CC$ with the standard absolute value and 
$X$ is of finite type over $\CC$, then $X^\an$ is homeomorphic to the analytic space $X(\CC)$.
\item If $k=\RR$ with the standard absolute value and 
$X$ is of finite type over $\RR$, then $X^\an$ is homeomorphic to the quotient of the analytic space $X(\CC)$ by the action of the complex conjugate.
\end{enumerate}

In general, we have a decomposition
$$
X^\an=X^\an[\infty]\cup X^\an[\mathrm{f}],  
$$
where $X^\an[\infty]$ is the subset of all archimedean semi-norms in $X^\an$, and $X^\an[\mathrm{f}]$ is the subset of all non-archimedean semi-norms  in $X^\an$.

\subsection{Our choice of base ring}

\kkk
Similar to \S\ref{sec uniform}, we introduce a uniform terminology for these two cases. 
Endow $k$ with a norm $|\cdot|_\mathrm{Ban}$ as follows.
If $k=\ZZ$, $|\cdot|_\mathrm{Ban}$ is the usual archimedean absolute value $|\cdot|_\infty$;
if $k$ is a field, $|\cdot|_\mathrm{Ban}$ is the trivial valuation $|\cdot|_0$. 
This makes $k$ into a Banach ring.

Concerning our special situation, we have the following  results and notations:
\begin{enumerate}
\item
If $k$ is a field, then $\CM(k)$ has only one element $v_0=|\cdot|_0$ by definition.
In this case, if $X$ is a finite type over $k$, then $X\mapsto X^\an$ is just the analytification functor constructed in \cite[\S3.5]{Ber1}. 

\item
If $k$ is a field, and $X$ is a projective regular curve over $k$, then $X^\an$ is the union of the closed line segments $\{|\cdot|_{v}^t:0\leq t\leq \infty\}$ for all closed points $v\in X$, by identifying $|\cdot|_{v}^0$ with the trivial norm $|\cdot |_0$ for all $v\in X$ as one point. 
Here $|\cdot|_{v}$ denotes the normalized valuation $\exp(-\ord_v)$.
The space $\CM(k(X))$ for the function field $k(X)$ is exactly the subspace of $X^\an$ obtained by removing the subset $\{|\cdot|_{v}^\infty:v\in X \text{ closed}\}$.
\item
In the arithmetic case ($k=\ZZ$), the space $\CM(\ZZ)$ is compact and path-connected. As described in \cite[1.4.3]{Ber1}, it is the union of the closed line segment 
$$[0,1]_\infty:=\{|\cdot|_{\infty}^t: 0\leq t\leq 1\},$$ 
and the closed-line segments 
$$[0,\infty]_p:=\{|\cdot|_{p}^t:0\leq t\leq \infty\}$$ 
for all finite primes $p$, by identifying the endpoints $|\cdot|_{\infty}^0$ and $|\cdot|_{p}^0$ for all finite primes $p$ with the trivial norm $|\cdot |_0$ of $\ZZ$. Here $|\cdot|_{\infty}$ and $|\cdot|_{p}$ denote the usual normalized valuations. 
The canonical injection $\iota:\Spec \ZZ\to \CM(\ZZ)$ sends the generic point to the trivial norm $|\cdot |_0$, and sends a prime $p$ to $|\cdot|_{p}^\infty$, the semi-norm of $\ZZ$ induced by the trivial norm of $\FF_p$.

\qquad The space $\CM(\QQ)$, defined by viewing $\QQ$ as a ring over the Banach ring $\ZZ$,  is exactly the subspace of $\CM(\ZZ)$ obtained by removing the subset $\{|\cdot|_{p}^\infty:p<\infty\}$. There is a very similar description for number fields.

\qquad For convenience, denote 
$$v_0=|\cdot|_0,\quad
v_\infty=|\cdot|_{\infty},\quad
v_\infty^t=|\cdot|_{\infty}^t,\quad
v_p=|\cdot|_{p},\quad
v_p^t=|\cdot|_{p}^t.$$
We may also write $\infty$ and $p$ for $v_\infty$ and $v_p$, viewed as points of $\CMZ$.
For convenience, denote by 
$$(0,1]_\infty,\quad (0,1)_\infty,\quad (0,\infty]_p,\quad [0,\infty)_p,\quad (0,\infty)_p$$ 
the sub-intervals of the line segments obtained by removing one or two endpoints;
for example,
$$(0,\infty)_p:=\{|\cdot|_{p}^t:0< t< \infty\}.$$ 
\item
In the arithmetic case ($k=\ZZ$), there is a structure map $X^\an \to \CM(\ZZ)$. 
This gives disjoint unions.
$$
X^\an=\bigcup_{v\in \CM(\ZZ)} X_v^\an,
$$
where $X_v^\an$ is the fiber of $X^\an$ above $v$. 
The most distinguished fibers are
$$
X_{\infty}^\an=X_{v_\infty}^\an=X_{\RR}^\an,\qquad
X_{p}^\an=X_{v_p}^\an=X_{\QQ_p}^\an.
$$
According to the structure of $\CM(\ZZ)$, 
we can further write $X^\an$ as a disjoint union of the following subspaces:
\begin{enumerate}
\item[(i)] $X_{v_0}^\an=(X_{\QQ}/\QQ)^\an$ under the trivial norm of $\QQ$; 
\item[(ii)] $X_{v_p^\infty}^\an=(X_{\FF_p}/\FF_p)^\an$ under the trivial norm of $\FF_p$ for finite primes $p$;
\item[(iii)] $X_{(0,\infty)_p}^\an$, homeomorphic to $X_{\QQ_p}^\an\times (0,\infty)$ for finite primes $p$; 
\item[(iv)] $X_{(0,1]_\infty}^\an$, homeomorphic to $X_{\RR}^\an\times (0,1]$. 
\end{enumerate}

\item 
In both cases (that $k$ is $\ZZ$ or a field), if $X$ is connected and of finite type over $k$, then $X^\an$ is path-connected. 
We can assume that $X$ is normal by passing to its normalization.

\qquad We first treat the geometric case that $k$ is a field. By blowing up $X$, there is a flat morphism $X\to C$ to a connected regular curve $C$ over $k$.
We can further assume that the fibers of $X\to C$ are connected by taking the integral closure of $C$ in $X$. 
The fibers of $X^\an\to C^\an$ are path-connected by induction and by the well-known case of non-trivial valuation fields.
There are finitely many (connected) closed curves $C_1,\cdots, C_n$ in $X$ 
such that $\Im(C_1\to C),\cdots, \Im(C_n\to C)$  is a Zariski open cover of $C$. 
Note that $C_1^\an,\cdots, C_n^\an$ are path-connected by example (2) above, so $X^\an$ is path-connected. 

\qquad In the arithmetic case, denote by $O_K$ the integral closure of $\ZZ$ in $\CO_X$. 
The fibers of $X^\an\to \CM(O_K)$ are path-connected.
For any finite extension $K'$ of $K$ and any open subscheme $C'$ of $\Spec O_{K'}$, the space $C'^\an$ is connected by an explicit description similar to (3).
So $X^\an$ is path-connected as in the geometric case. 
\end{enumerate}

Note that any multiplicative semi-norm on $\ZZ$ is bounded by the standard Archimedean absolute value and thus belongs to $\CMZ$. 
The space $X^\an$ in the arithmetic space is the ``largest'' Berkovich space associated with $X$ defined in terms of multiplicative semi-norms.

\kkk
Let $X$ be a proper scheme over $k$.
There is a \emph{specialization map} (or \emph{reduction map})
$$
r: X^\an \lra X
$$
defined as follows. 
For any point $x\in X^\an[\mathrm{f}]$, recall that $H_x$ is the (complete) residue field of $x$ in $X^\an$, denote by $R_x$ the valuation ring of $H_x$, and denote by $m_x$ the maximal ideal of $R_x$.
As $X$ is proper over $k$, the valuative criterion gives a unique $k$-morphism
$\Spec R_x\to X$ extending the $k$-morphism $\Spec H_x\to X$ associated to $x$.
Define $r(x)$ to be the image of the unique closed point of $\Spec R_x$ in $X$. 

For any point $x\in X^\an[\infty]$, we still have a morphism $\Spec H_x\to X$. Here $H_x$ is isomorphic to either $\RR$ or $\CC$. Define $r(x)$ to be the image of $\Spec H_x$ in $X$.

\subsection{Density result}

We are interested in $X^\an$ for an essentially quasi-projective scheme $X$ over $k$.
The following result asserts that the Berkovich spaces of essentially quasi-projective schemes do not lose "many points" from those of their quasi-projective models. 
In the arithmetic case, the space $X^\an$ is somehow determined by its fibers above the non-trivial absolute values of $\QQ$.

\begin{lem}\label{density}
\kkk
Let $X$ be a flat and essentially quasi-projective integral scheme over $k$.
\begin{enumerate}
\item
Let $X\to \CU$ be a quasi-projective model of $X$ over $k$. 
Then the induced map $X^\an\to \CU^\an$ is continuous, injective, and with a dense image.
Moreover, the set of $v\in X^\an$ corresponding to discrete or archimedean valuations of $H_v$ is dense in $\CU^\an$. 
\item 
If $k=\ZZ$, then $X^\an\setminus X^\an_{\iota(\Spec \ZZ)}$ is dense in $X^\an$. Here 
$\iota:\Spec\ZZ\to \CMZ$ is the canonical injection whose image consists of $v_0$ and $v_p^\infty$ for finite primes $p$. 
\end{enumerate}

\end{lem}
\begin{proof}
We first prove (1). Only the density is not automatic from the definitions. 
Denote by $F$ the function field of $X$, which is also the function field of $\CU$. There is a composition of injections $(\Spec F)^\an\to X^\an \to \CU^\an$. 
It suffices to prove that the set of discrete or archimedean $v\in (\Spec F)^\an$ is dense in $\CU^\an$.

We first prove that $(\Spec F)^\an[\infty]$ is dense in $\CU^\an[\infty]$. 
Assume  $k=\ZZ$ to have a non-trivial statement. 
Recall that we have  
$$(\Spec \ZZ)^\an[\infty]=\{v_\infty^t=|\cdot|_{\infty}^t: 0< t\leq 1\}\simeq (0,1].$$ 
As before, denote by $(\Spec F)^\an_{v_\infty^t}$ and $\CU^\an_{v_\infty^t}$ the fibers of 
$(\Spec F)^\an$ and $\CU^\an$ above $v_\infty^t\in \CM(\ZZ)$. 
It suffices to prove that $(\Spec F)^\an_{v_\infty^t}$ is dense in $\CU^\an_{v_\infty^t}$ for every $t\in (0,1]$. 
By Ostrowski's theorem, any  $v\in \CU^\an_{v_\infty^t}$ is induced by the valuation $v_\infty^t$ of $\CC$ via
a morphism $\Spec \CC\to \CU$. It follows that we have a natural surjection 
$\CU(\CC)\to \CU^\an_{v_\infty^t}$, and it induces a homeomorphism 
$\CU(\CC)/\Gal(\CC/\RR)\to \CU^\an_{v_\infty^t}$. 
An element of $\CU(\CC)$ gives an element of $(\Spec F)^\an_{v_\infty^t}$ in this process if and only if the image of the corresponding morphism $\Spec \CC\to \CU$ is the generic point $\Spec F$ of $\CU$. 
Denote by $\CU(\CC)_\mathrm{gen}$ the subspace of such elements of $\CU(\CC)$. 
This induces a homeomorphism $\CU(\CC)_\mathrm{gen}/\Gal(\CC/\RR)\to (\Spec F)^\an_{v_\infty^t}$.
It is reduced to prove that $\CU(\CC)_\mathrm{gen}$ is dense in $\CU(\CC)$. 
Note that $\CU(\CC)\setminus \CU(\CC)_\mathrm{gen}$ is the union of $V(\CC)$ over all Zariski closed $V\subsetneq \CU_\QQ$. Then $\CU(\CC)\setminus \CU(\CC)_\mathrm{gen}$ is a countable union of proper Zariski closed subsets of $\CU(\CC)$. 
This implies that $\CU(\CC)_\mathrm{gen}$ is dense in $\CU(\CC)$.
In fact, as the smooth locus of $\CU(\CC)$ is dense in  $\CU(\CC)$, 
we can assume that $\CU(\CC)$ is smooth, and by cover $\CU(\CC)$ by open balls, we see that $\CU(\CC)\setminus \CU(\CC)_\mathrm{gen}$ has Lebesgue measure 0 over each ball. 

Now  we treat the non-archimedean points (where $k=\ZZ$ or $k$ is a field). 
We can assume that $\CU$ is projective over $k$ by passing to a projective model.
A point $\xi\in \CU^\an[\mathrm{f}]$ is called \emph{divisorial} if there is a birational morphism $\CU'\to \CU$ from a normal integral scheme of finite type over $k$, together with a prime Weil divisor $D\subset \CU'$, such that $|\cdot|_\xi=\exp(-t\,\ord_D)$ for some constant $t>0$. 
Note that $\xi$ is discrete and actually lies in $(\Spec F)^\an[\mathrm{f}]$. 
It suffices to prove that the set of divisorial points is dense in $\CU^\an[\mathrm{f}]$ for all projective variety $\CU$ over $k$. 

Note that the analogous statement for a non-archimedean field $k$ (with a non-trivial valuation) is a well-known result. For example, in this case, Berkovich's theory implies that the analytic space $\CU^\an$ has a topological basis consisting of strictly $k$-affinoid domains, and any $k$-affinoid domain has a (non-empty) Shilov boundary. By \cite[A.3, A.6, A.9]{GM}, any point in the Shilov boundary of a strictly $k$-affinoid domain is actually a divisorial point. 

Return to the original $(\CU,k)$. We will prove the density by induction on the dimension of $\CU$. Assume that $\CU$ is normal by passing to its normalization.
The case of dimension one essentially follows from the explicit descriptions, noting that the space $\CM(O_F)$ for a number field $F$ has a description analogous to that of $\CM(\ZZ)$.
Assume that $\dim \CU>1$ and that the density statement holds for all lower dimensions. 

We first prove the case $k=\ZZ$ (for the density of divisorial points). Let $F$ be the algebraic closure of $\QQ$ in $\QQ(\CU)$, so that $\CU$ has geometrically connected fibers over $\Spec O_F$. 
Write $\CU^\an$ as unions of fibers $\CU^\an_v$ above $v\in \CM(O_F)$.
By induction and by the case of non-archimedean fields, it suffices to prove that, for any $v\in \CM(O_F)\setminus \CM(F)$ corresponding to the trivial norm of the residue field $\FF_\wp$ for some prime ideal $\wp$ of $O_F$, and any irreducible component $U$ (endowed with the reduced scheme structure) of 
the special fiber $\CU_{\FF_\wp}$, any divisorial point $\xi\in (U/\FF_\wp)^\an$ lies in the closure of the divisorial points of $\CU^\an$.

To illustrate the key idea, we will first treat the nice case that $U=\CU_{\FF_\wp}$ is smooth over $\FF_p$, $\xi\in (U/\FF_\wp)^\an$ is a divisorial point corresponding to a prime divisor $D\subset U$, and that there is a section $S$ of $\CU_{O_{F_v}}$ over $O_{F_v}$ such that the point $u=S_{\FF_\wp}$ is regular in all $D, U, \CU$. 
In that case, we can find a local system of parameters $\varpi,x_1,\cdots, x_d\in \CO_{\CU,u}$, such that locally at $u\in \CU$, $U$ is defined by the ideal $(\varpi)$, and $D$ is defined by the ideal $(\varpi,x_1)$. Here we require $\varpi\in O_F$ to be a generator of $\wp$.
Any $f\in \CO_{\CU,u}$ can be uniquely written as a power series
$$
f= \sum_{n\geq0} a_n x_1^n,\quad\ a_n\in O_{F_\wp}[[x_2,\cdots, x_d]].
$$
Let $0<r<1$ and $t>0$ be real numbers.
Set 
$$
|f|_{y_t}= \max_{n \geq 0} \big(|a_n|_\wp^{t}r^n\big).
$$
Here $|\cdot|_\wp=\exp(-\ord_\wp)$. 
We can check that $y_t\in \CU^\an$ with 
$$
\lim_{t\to \infty}|f|_{y_t}=r^{\min\{n:|a_n|_\wp=1\}}=r^{\ord_{x_1}(f\, \mathrm{mod}\, \wp)}.
$$
As $r$ varies, the right-hand side gives exactly all the divisorial points 
of $(U/\FF_\wp)^\an$ corresponding to the divisor $D\subset U$.
This finishes the nice case.

Now  we need to make a few operations and replacements to convert the general case $(\CU_{O_{F_\wp}}, \wp,U,\xi,D)$ to the nice case. Note that the situation is local at $\wp$ over $O_F$, so we are only concerned with $\CU_{O_{F_\wp}}$ instead of $\CU$.  
The process is mostly geometric but a little tedious.

First, we convert to the case that $U$ is normal and $\xi\in (U/\FF_\wp)^\an$ is a divisorial point corresponding to a prime divisor $D\subset U$. 
Assume that $\xi\in (U/\FF_\wp)^\an$ is a divisorial point given by a prime divisor $D'$ on a normal projective scheme $U'$ over $\FF_\wp$ with a birational morphism $U'\to U$. 
Then $ U'to U$ is obtained by blowing up $U$ along a closed subscheme $Z$ of $U$. 
Let $\CU'\to \CU$ be the blowing-up of $\CU$ along $Z$, the strict transform of $U$ is exactly $U'$. 
Replace $(\CU,U)$ by $(\CU',U')$.

Second, there is a finite and flat morphism 
$\CU_{O_{F_\wp}}\to \BP^d_{O_{F_\wp}}$ with $d=\dim\CU$.
Take any ample line bundle $\CL$ on $\CU_{O_{F_\wp}}$. 
Denote by $\CL_0$ its pull-back to $\CU_{\FF_\wp}$. 
There is a positive integer $m$, such that $\Gamma(\CU_{\FF_\wp},m\CL_0)$ contains a base-point-free subspace $V$ of dimension $d+1$ such that the corresponding morphism $\CU_{O_{F_\wp}}\to \BP(V)$ is finite.
This is done by the classical argument of embedding to a projective space of a high dimension and projecting to hyperplanes successively.
For $m$ large enough, the map 
$\Gamma(\CU_{O_{F_\wp}},m\CL)\to \Gamma(\CU_{\FF_\wp},m\CL_0)$ is surjective, and thus we can lift $V$ to an $O_{F_\wp}$-submodule 
$\wt V$ of $\Gamma(\CU_{O_{F_\wp}},m\CL)$ of the same rank. 
The morphism $\CU_{O_{F_\wp}}\to \BP(\wt V)$
satisfies the requirement. 

Third, in the morphism $\CU_{O_{F_\wp}}\to \BP^d_{O_{F_\wp}}$, denote by $U_0$ and $D_0$ the images of $U$ and $D$ respectively. Note that $U_0=\BP^d_{\FF_\wp}$. 
We claim that the result for $(\BP^d_{O_{F_\wp}}, U_0, D_0, \xi_0)$ implies that for 
$(\CU_{O_{F_\wp}}, U, D, \xi)$. Here $\xi_0$ is the image of $\xi$ in $U_0^\an$. 
In fact, denote by $\CU'\to \BP^d_{O_{F_\wp}}$ the Galois closure of 
$\CU_{O_{F_\wp}}\to \BP^d_{O_{F_\wp}}$, defined to be the normalization of the Galois closure of the function fields.
Denote by $\xi_1',\cdots, \xi_m'$ the preimage of $\xi$ in $(\CU'/O_{F_\wp})^\an$. 
If the result for $(\BP^d_{O_{F_\wp}}, U_0, D_0, \xi_0)$ holds, then by taking preimages of divisorial point, one of $\xi_i'$ lies in the closure of divisorial points of 
$(\CU'/O_{F_\wp})^\an$. This also uses compactness of $(\CU'/O_{F_\wp})^\an$.
As the Galois group acts transitively on 
$\xi_1',\cdots, \xi_m'$, any $\xi_i'$ lies in the closure of divisorial points of $(\CU'/O_{F_\wp})^\an$.
Taking images in $\CU^\an$, we see that $\xi$ lies in the closure of divisorial points $\CU^\an$.

Replace $(\CU_{O_{F_\wp}}, U, D, \xi)$ by $(\BP^d_{O_{F_\wp}}, U_0, D_0, \xi_0)$. We have converted the problem to the case that $U=\CU_{\FF_\wp}$ is smooth over $\FF_\wp$ and $\xi\in (U/\FF_\wp)^\an$ is a divisorial point corresponding to a prime divisor $D\subset U$. 

Fourth, it remains to construct a section $S$ of $\CU_{O_{F_\wp}}$ over 
$O_{F_\wp}$ passing through a regular point of $D$.
This is easily done by the base change by a finite unramified extension of $F_\wp$.

Therefore, our proof of (1) for the arithmetic case $k=\ZZ$ is complete. 
If $k$ is a field, by blowing-up $\CU$, there is a fibration $\CU\to \BP_k^1$.
Then the above induction argument still works.

Now we prove (2). 
Let $\CU$ be a normal projective model of $X$ as in (1).
Still consider the composition $(\Spec F)^\an\to X^\an \to \CU^\an$. 
The topologies of $(\Spec F)^\an$ and $X^\an$ are the same as the subspace topologies induced from $\CU^\an$. 
It suffices to prove that $(\Spec F)^\an\setminus (\Spec F)^\an_{\iota(\Spec \ZZ)}$ is dense in $\CU^\an$.
In (1), we have proved that $(\Spec F)^\an$ is dense in $\CU^\an$, but the proof actually gives the statement that $(\Spec F)^\an\setminus (\Spec F)^\an_{\iota(\Spec \ZZ)}$ is dense in 
$\CU^\an\setminus \CU^\an_{v_0}$.
Here $v_0$ denotes the trivial norm of $\ZZ$.
Therefore, it suffices to prove that any divisorial point $\xi$ of $\CU^\an_{v_0}$ lies in the closure of $\cup_{p<\infty} \CU^\an_{\FF_p}$
in $\CU^\an$.

As above, we have a geometrically connected morphism $\CU\to \Spec O_F$. 
Then the divisorial point $\xi$ lies in $\CU^\an_{v_0}=(\CU_F/F)^\an$ under the trivial norm of $F$.
By replacing $\CU$ with a blowing-up if necessary, we can assume that $\xi$ corresponds to a prime divisor $D$ of $\CU_F$.
Denote by $\CD$ the Zariski closure of $D$ in $\CU$.
For all but finitely many prime ideal $\wp$ of $O_F$, the reduction $\CU_{\FF_\wp}$ is normal and $\CD_{\FF_\wp}$ is a prime divisor of $\CU_{\FF_\wp}$.
For any rational function $f$ of $\CU$, for all but finitely many prime ideal $\wp$ of $O_F$, the specializations above $\wp$ of the irreducible components of $\div(f|_{\CU_F})$ are irreducible and distinct, which give $\ord_D(f|_{\CU_F})=\ord_{\CD_{\FF_\wp}}(f|_{\CU_{\FF_\wp}})$. 
This proves that $\exp(-\ord_D)$ is the limit of $\exp(-\ord_{\CD_{\FF_\wp}})$ as $\wp$ varies. 
Therefore, $x$ lies in the closure of divisorial points of $(\Spec F)^\an\setminus (\Spec F)^\an_{v_0}$ in $\CU^\an$.
This proves (2).
\end{proof}

\section{Arithmetic divisors and metrized line bundles}
\label{sec analytic objects}

In this section, we introduce arithmetic divisors and metrized line bundles on Berkovich spaces, which are analytic counterparts of the adelic divisors in  \S\ref{sec adelic divisors} and the adelic line bundles in \S\ref{sec adelic line bundles}.

\subsection{Arithmetic divisors}

Let $k$ be a commutative Banach ring, which is also an integral domain. 
Let $X$ be an integral scheme over $k$.
Let $X^\an=(X/k)^\an$ be the Berkovich space defined above.
 
Let $D$ be a Cartier divisor on $X$. By a {\em Green function} of the divisor $D$ on $X^\an$, we mean a continuous function $g: X^\an\setminus |D|^\an \to \RR$ with logarithmic singularity along $D$ in the sense that, for any rational function $f$ on a Zariski open subset $U$ of $X$ satisfying $\div(f)=D|_U$, the function $g+\log |f|$ can be extended to a continuous function on $U^\an$.

The pair $\overline D=(D, g)$ is called an {\em arithmetic divisor} on $X^\an$.  
An {arithmetic divisor} is called {\em effective} if $D$ is an effective Cartier divisor on $X$ and $g\geq 0$ on $X^\an\setminus |D|^\an$.
An {arithmetic divisor} is called {\em principal} if it is of the form 
$$\wh\div_{X^\an}(f):=(\div(f), -\log |f|)$$ 
for some nonzero rational function $f$ on $X$.

An arithmetic divisor $\OD$ or its Green function $g$ is called 
\emph{norm-equivariant} if for any points $x,x_1\in X^\an\setminus |D|^\an$ satisfying $|\cdot |_x=|\cdot |_{x_1}^t$ for some $0\leq t< \infty$ locally on $\CO_X$, we have 
$g(x)=t\, g(x_1)$.
By definition, principal arithmetic divisors are norm-equivariant.

Denote by $\wh \Div (X^\an)$ the group of arithmetic divisors on $X^\an$, by $\wh \Pr (X^\an)$ the group of {principal arithmetic divisors} on $X^\an$, and by $\wh \Div (X^\an)_\eqv$ the group of norm-equivariant arithmetic divisors on $X^\an$.
Denote the class group of arithmetic divisors as
$$\wh \CaCl (X^\an):=\wh \Div (X^\an)/\wh \Pr (X^\an),$$
$$\wh \CaCl (X^\an)_\eqv:=\wh \Div (X^\an)_\eqv/\wh \Pr (X^\an).$$

Notice that for any arithmetic divisor $\overline D=(D, g)$ on $X^\an$, 
the algebraic part $D$ is a Cartier divisor on $X$ (instead of $X^\an$), and $g$ is a function on $X^\an$. 
We take this ad hoc definition to avoid defining general Cartier divisors on $X^\an$ by the lack of a good theory of analytic functions on $X^\an$.

\subsection{Metrized line bundles}
 
Let $k$ be a commutative Banach ring, which is also an integral domain. 
Let $X$ be an integral scheme over $k$.
Let $X^\an=(X/k)^\an$ be the Berkovich space defined above.

Let $L$ be a line bundle on $X$. 
At each point $x\in X^\an$, denote by $\bar x$ the image of $x$ in $X$. 
The \emph{fiber} $L^\an(x)$ of $L$ at $x$ is defined to be the $H_x$-line $L(\bar x)\otimes_{k(\bar x)} H_x$, or equivalently the completion of the fiber $L(\bar x)$ of $L$ on $\bar x$ for the  semi-norm $|\cdot|_x$.
By a \emph{metric} $\|\cdot \|$ of $L$ on $X^\an$ we mean a continuous metric on $\coprod_{x\in X^\an} L^\an(x)$ compatible with the semi-norms
 on $\CO_X$. More precisely, to each point $x\in X^\an$, we assign a norm $\|\cdot \|_x$ on the $H_x$-line
 $L^\an(x)$ which is compatible with the norm $|\cdot |_x$ of $H_x$ in the sense that
 $$\|f\ell\|_x=|f|_x\cdot \|\ell\|_x, \qquad f\in H_x, \quad \ell \in L^\an(x).$$
 We always assume that the metric $\|\cdot \|$ on $L$ is {\em continuous} in the sense that, for any section $\ell$ of $L$ on a Zariski open subset $U$ of $X$, the function $\|\ell (x)\|=\|\ell (x)\|_x$ is continuous in $x\in U^\an$.
 
The pair $(L,\|\cdot\|)$ above is called a \emph{metrized line bundle} on $X^\an$.
An \emph{isometry} from a metrized line bundles $(L,\|\cdot\|)$ to another one $(L',\|\cdot\|')$ is an isomorphism $i:L\to L'$ of line bundles on $X$ such that $\|\cdot\|=i^*\|\cdot\|'$.

A metrized line bundle $\OL=(L,\|\cdot\|)$ or its metric $\|\cdot\|$ is called 
\emph{norm-equivariant} if for any rational section $s$ of $L$ on $X$, and any points $x,x_1\in X^\an\setminus |\div(s)|^\an$ satisfying $|\cdot |_x=|\cdot |_{x_1}^t$ for some $0\leq t< \infty$ locally on $\CO_X$, we have 
$\|s \|_x=\|s \|_{x_1}^t$.

Denote by $\wh \Picc (X^\an)$ the \emph{category} of metrized line bundles on $X^\an$, where the morphisms are isometries. 
Denote by $\wh \Pic (X^\an)$ the \emph{group} of isometry classes of metrized line bundles on $X^\an$. 
Denote by $\wh \Picc (X^\an)_\eqv$ (resp. $\wh \Pic (X^\an)_\eqv$) the full subcategory (resp. the subgroup) of norm-equivariant line bundles in $\wh \Picc (X^\an)$ (resp. $\wh \Pic (X^\an)$).

Similar to $\wh\Div(X^\an)$, elements of $\wh\Pic (X^\an)$ are of the form $(L,\|\cdot\|)$, where $L$ is a line bundle on $X$ (instead of $X^\an$) and $\|\cdot\|$ is a metric on $X^\an$. 
We have forgetful maps
 $$\wh\Picc (X^\an )\lra \Picc (X),\quad\
 \wh\Pic (X^\an )\lra \Pic (X).$$
The fibers of the second map are homogeneous spaces of the group of metrics on $\CO_X$.

There are canonical isomorphisms 
$$\wh \CaCl (X^\an) \lra \wh \Pic (X^\an),$$
$$\wh \CaCl (X^\an)_\eqv \lra \wh \Pic (X^\an)_\eqv.$$
In fact, given any arithmetic divisor, $(D, g)$ on $X^\an$, the term $e^{-g/2}$ defines a metric on $\CO (D)$, and thus we obtain a metrized line bundle on $X^\an$.
Conversely, for any metrized line bundle $(L,\|\cdot \|)$ on $X^\an$, if $s$ is a rational section of $L$, then 
$$\wh\div_{X^\an}(s):=(\div(s),-\log \|s\|)$$ 
defines an arithmetic divisor on $X^\an$. 
Both processes keep the properties of being norm-equivariant.

In the case $k=\ZZ$, a norm-equivariant Green function or a norm-equivariant metric on a line bundle on $X^\an$ is uniquely determined by its restriction to the disjoint union of the distinguished fibers $X_{v}^\an=X_{\QQ_v}^\an$ over all places $v\leq \infty$. This follows from Lemma \ref{density}(2). 
Later on, all of Green functions and metrics in our consideration will be 
norm-equivariant.

\section{Analytification of adelic divisors} \label{sec analytification divisors}

The adelic divisors in  \S\ref{sec adelic divisors} induce norm-equivariant arithmetic divisors on Berkovich spaces. The goal of this section is to study this analytification process. The main result is as follows:

\begin{prop}  \label{injection1}
\kkk
Let $X$ be a flat and essentially quasi-projective integral scheme over $k$.
There are canonical injective maps 
$$\wh \Div (X/k)\lra \wh \Div (X^\an)_\eqv,$$
$$\wh \CaCl (X/k)\lra \wh \CaCl (X^\an)_\eqv.$$ 
\end{prop}

\begin{rmk}
In a recent work, Song \cite{Son24} proves the the maps are also surjective if $X$ is quasi-projective over $k$. We refer to the loc. cit. for more details, while we will only treat the injectivity in this book.
\end{rmk}

For an adelic divisor $\OD$ on $X/k$ with underlying divisor $D$, its image under the first map takes the form $(D, \wt g_\OD)$, where $\wt g_\OD$ is a Green function of $D$ on $X^\an$.
We call $\wt g_\OD$ the \emph{total Green function} of $\OD$ on $X^\an$.  
In the following, we will construct the maps and prove the injectivity in the order of projective case, quasi-projective case, and essentially quasi-projective case.

\subsection{Projective case} 

\kkk
\ccc

Let $\CX$ be a projective variety over $k$. 
Then there is a canonical map 
$$\wh\Div (\CX)\lra \wh \Div (\CX^\an)_\eqv.$$
In the following, for any $\OCD=(\CD,g)\in \wh\Div (\CX)$, we will introduce a Green function $\wt g$ of $\CD$ on $\CX^\an$, and define the map by $(\CD,g)\mapsto (\CD,\wt g)$.

We will define $\wt g$ according to the decomposition $\CX^\an=\CX^\an[\mathrm{f}]\cup \CX^\an[\infty]$, and then check the continuity. 

For any point $x\in \CX^\an[\mathrm{f}]$, recall that there is a specialization map $r:\CX^\an[\mathrm{f}]\to \CX$ by the properness of $\CX$ over $k$.
Let $\CU$ be a Zariski open subscheme of $\CX$ containing $r(x)$ such that $\CD|_\CU$ is defined by a single equation $f\in k(\CU)^\times$ on $\CU$.
By $r(x)\in \CU$, the image of $\Spec R_x\to X$ lies in $\CU$ and thus $x\in \CU^\an$.
Define $\wt g(x)=-\log|f(x)|$. 
This definition is independent of the choice of $(\CU,f)$.

It is easy to define $\wt g$ on $\CX^\an[\infty]$ (in the arithmetic case). 
In fact, the Green function $g$ on $\CX(\CC)$ descends to the fiber $\CX^\an_{\infty}=\CX_\RR^\an$. This gives the definition of $\wt g$ on $\CX^\an_{\infty}$. 
It extends to $\CX^\an[\infty]$ by requiring $\wt g$ to be norm-equivariant. 
In fact, for any point $x\in \CX^\an[\infty]$, there is a unique point $x_1 \in \CX^\an_{\infty}$ such that $|\cdot |_x=|\cdot |_{x_1}^t$ for some $0<t<1$, and then we set $\wt g(x)=t\,\wt g(x_1)$. 

Now we prove that $\wt g$ is indeed a Green function; i.e. $\wt g$ is continuous on $\CX^\an\setminus |\CD|^\an$, and has logarithmic singularity along $\CD$.

\begin{lem} \label{green's functions}
The function $\wt g$ is a Green function of $\CD$ on $\CX^\an$.
\end{lem}

\begin{proof}
We first note that the continuity of $\wt g$ on $\CX^\an\setminus |\CD|^\an$ (for all such $\OCD$) implies that $\wt g$ has logarithmic singularity along $\CD$.
Let $f$ be a local equation of $\CD$ on an open subscheme $\CU$ of $\CX$.
Assume that the continuity holds for the arithmetic divisor $(\CD-\div_\CX(f),g+\log|f|_\infty)$, i.e. 
$\wt g+\log|f|$ is continuous on $\CX^\an\setminus |\CD-\div(f)|^\an$, 
Then $\wt g$ has the correct logarithmic singularity on $\CU^\an$. 
Vary $(\CU,f)$ to cover $\CX$. 

Now we prove the continuity of $\wt g$ on $\CX^\an[\mathrm{f}]\setminus |\CD|^\an$.
Let $r:\CX^\an\to \CX$ be the specialization map.
Let $\{x_m\}_{m\geq1}$ be a sequence in $\CX^\an[\mathrm{f}]\setminus |\CD|^\an$ converging to 
$x\in \CX^\an[\mathrm{f}]\setminus |\CD|^\an$. 
We need to prove that $\wt g(x_m)$ converges to $\wt g(x)$.
Let $\CU_1,\cdots, \CU_n$ be an open cover of $\CX$ such that, for any $i=1,\cdots, n$, 
$\CU_i$ contains $r(x)$ and $\CD$ is defined by a single equation $f_i$ on $\CU_i$. 

To see the existence of the open cover, by quasi-compactness, it suffices to prove that for any point $y\in \CX$, there is an open neighborhood of $\{r(x),y\}$ in $\CX$ such that $\CD$ is principal on $\CU$.
Note that $\CD$ is principal on $\CU$ if and only if the line bundle $\CO(\CD)$ is trivial in $\Pic(\CU)$.
We can further assume that $\CD$ is very ample by writing $\CD$ as the difference of two very ample Cartier divisors on $\CX$. 
Then there is an embedding $\CX\to \BP^N_k$ using global sections of $\CO(\CD)$. 
For any hyperplane $\CH$ of $\BP^N_k$, the line bundle $\CO(\CD)$ is trivial in $\Pic(\CX\setminus \CH)$, since $\CO_{\BP^N_k}(1)$ is trivial in $\Pic(\BP^N_k\setminus \CH)$. Now it suffices to choose a hyperplane $\CH$ disjoint with $\{r(x),y\}$. 
This is easy if $k$ is infinite. 
If $k$ is finite, it is also easy if $N$ is large. 

Now we have the open cover $\CU_1,\cdots, \CU_n$.
Denote by $I_i$ the set of $m\geq 1$ such that $r(x_m)\in \CU_i$. 
Then we have $I_1\cup\cdots \cup I_n=\{1,2,\cdots\}$. 
It suffices to prove $\lim_{m\in I_i}\wt g(x_m)=\wt g(x)$ for each $i=1,\cdots,n$.
By definition, $r(x)\in \CU$ implies that the image of the closed point of $\Spec R_x\to \CX$ lies in $\CU$, where $R_x\subset H_x$ is the valuation ring. 
This implies that the image of $\Spec R_x\to \CX$ lies in $\CU$, and thus $x\in \CU^\an$. 
Thus $\wt g(x)=-\log|f(x)|$ by definition. 
Similarly, $m\in I_i$ implies $x_m\in \CU^\an$ and $\wt g(x_m)=-\log|f(x_m)|$.
It follows that $\lim_{m\in I_i}\wt g(x_m)=\wt g(x)$. 
This proves that $\wt g$ is continuous on $\CX^\an[\mathrm{f}]\setminus |\CD|^\an$.

If $k=\ZZ$, we need to make extra arguments to extend the continuity of $\wt g$ from 
$\CX^\an[\mathrm{f}]\setminus |\CD|^\an$ to $\CX^\an\setminus |\CD|^\an$.
By definition, $\wt g$ is continuous on $\CX^\an[\infty]\setminus |\CD|^\an$.
It remains to prove that $\wt g$ is continuous when $\CX^\an[\infty]\setminus |\CD|^\an$ approaches $\CX^\an_{v_0}\setminus |\CD|^\an$, where $v_0\in \CMZ$ is the trivial norm of $\ZZ$. Namely, let $\{x_m\}_{m\geq1}$ be a sequence in $\CX^\an[\infty]\setminus |\CD|^\an$ converging to 
a point $x$ in $\CX^\an_{v_0}\setminus |\CD|^\an$. We need to prove that $\wt g(x_m)$ converges to $\wt g(x)$.

The canonical homeomorphism $\CX^\an[\infty]\to \CX^\an_{\infty}\times (0,1]$
induces a projection $\pi:\CX^\an[\infty]\to \CX^\an_{\infty}$.
Here $\CX^\an_{\infty}= \CX_{\RR}^\an=\CX(\CC)/\Gal(\CC/\RR)$ is compact.
To prove $\lim_{m\geq1}\wt g(x_m)=\wt g(x)$, by proof by contradiction, it suffices to prove 
$\lim_{m\in I}\wt g(x_m)=\wt g(x)$ for all subsequences $I$ of $\{1,2,\cdots\}$ such that  
$\{\pi(x_m)\}_{m\in I}$ converges in $\CX^\an_{\infty}$.

For such a subsequence $I$, denote by $z=\lim_{m\in I}\pi(x_m)$ in 
$\CX^\an_{\infty}$.
There is an open neighborhood $\CU$ of $\{r(x),r(z)\}$ in $\CX$ such that $\CD$ is defined by a single equation $f$ on $\CU$. 
The existence of $\CU$ has already been proved in the above non-archimedean case.  

Similarly, the condition $r(x),r(z)\in\CU$ implies $x,z\in \CU^\an$. 
This holds for $x$ as in the above non-archimedean case and holds for $z$ since $r(z)$ is the image of $\Spec H_z\to \CX$. 

By removing finitely many elements of $I$, we can assume that $x_m$ lies in $\CU^\an$ for every $m\in I$. 
Then $-\log |f|(x_m)$ for $m\in I$ converges to $-\log |f|(x)=\wt g(x)$.
Denote $h(y)=\wt g(y)+\log |f|(y)$, as a function on $\CU^\an$. 
It suffices to prove $\lim_{m\in I}h(x_m)=0$.

Note that $h$ is norm-equivariant on $\CU^\an$.
For $m\in I$, denote by $t_m$ the image of $x_m$ under the canonical projection 
$\CX^\an[\infty]\to (0,1]_\infty=(0,1]$. We have $\lim_{m\in I} t_m\to 0$.
It follows that 
$\lim_{m\in I}h(x_m)=\lim_{m\in I}t_m \, h(\pi(x_m))=0$
as $\lim_{m\in I}\pi(x_m)=z$ in $\CU^\an_{\infty}$. 
This finishes the proof.
\end{proof}

The following effectivity result will be very useful in proving the injectivity of the analytification map. 

\begin{lem}  \label{effectivity1}
\kkk
Let $\CX$ be a projective variety over $k$ and let $i: X\to \CX$ be a pro-open immersion.  
Let $\OCD=(\CD,g)$ be an arithmetic divisor on $\CX$, and denote by $\tilde g$ the
Green function of $\CD$ on $\CX^\an$ induced by $\OCD$. 
Assume one of the following two conditions:
\begin{enumerate}
\item $\CX$ is normal;
\item the scheme $\CX$ is integrally closed in $X$, and the Cartier divisor $\CD|_X$ is effective on $X$. 
\end{enumerate}
Then $\OCD$ is effective if and only if $\tilde g\geq 0$ on $\CX^\an\setminus |\CD|^\an$.
\end{lem}

\begin{proof}
Note that (1) is a special case of (2) with $X$ equal to the generic point of $\CX$, but we list it separately for its independent importance.
It suffices to prove the ``if'' part.  
Assuming $\tilde g\geq 0$,   we need to prove that $\CD$ is effective. 
This is an analytic version of Lemma \ref{effectivity}.
It suffices to prove that for any $v\in \CX\setminus X$ of codimension one in $\CX$, the valuation $\ord_v(\CD)$ in the local ring $\CO_{\CX, v}$ is non-negative.  
Consider the divisorial point $\xi=\exp(-\ord_v)$ of $\CX^\an$. 
Let $f$ be a local equation of $\CD$ in an open neighborhood of $v$ in $\CX$.
By definition, 
$$\wt g(\xi)=-\log |f(\xi)|=-\log (\exp(-\ord_vf))=\ord_vf=\ord_v(\CD).$$ 
It follows that $\ord_v(\CD)\geq 0$.
\end{proof}

\subsection{Quasi-projective case} 

\kkk
Let $\CU$ be a quasi-projective variety over $k$ and $\CX$ be a projective model of $\CU$.
The analytification map 
$$\wh\Div (\CX)\lra \wh \Div (\CX^\an)_\eqv$$
defined above induces a map 
$$\wh\Div (\CX,\CU)\lra \wh \Div (\CU^\an)_\eqv,$$ 
which sends $\OCD=(\CD,g_\CD)$ to 
$\OCD^\an:=(\CD|_{\CU}, \tilde g_{\CD})$. Here $\CD|_{\CU}$ is an integral part of $\OCD$, which is an integral Cartier divisor on $\CU$. 
By direct limit, the map gives a map
$$\wh\Div (\CU/k)_\rmod\lra \wh \Div (\CU^\an)_\eqv.$$
In the following, we prove that the map can be extended to adelic divisors of quasi-projective varieties by taking limits.

\begin{proof}[Proof of Proposition \ref{injection1}: quasi-projective case]

We need to define and prove the injectivity of
$$\wh\Div (\CU/k) \lra \wh \Div (\CU^\an)_\eqv,$$
$$\wh\CaCl (\CU/k) \lra \wh \CaCl (\CU^\an)_\eqv.$$
The injectivity of the first map implies that of the second map.
In fact, if 
$\wh\Div (\CU/k) \to \wh \Div (\CU^\an)_\eqv$
is defined and injective, 
then the map $\wh\Pr (\CU/k)_\rmod \to \wh \Pr (\CU^\an)$
is also injective.
Thus $\wh\Pr (\CU/k)_\rmod \to \wh \Pr (\CU^\an)$ is bijective as both groups are quotients of $k(\CU)^\times$.
The quotients give a well-defined and injective map $\wh\CaCl (\CU/k) \to \wh \CaCl (\CU^\an)_\eqv$.

To treat the first map, we will extend the map
$$\wh\Div (\CU/k)_\rmod \lra \wh \Div (\CU^\an)_\eqv$$
to a map 
$$\wh\Div (\CU/k) \lra \wh \Div (\CU^\an)_\eqv$$
by continuity.
Recall that the left-hand side is endowed with the boundary topology using $\OCE_0$; similarly, we endow the right-hand side with the boundary topology using the divisor $\OCE_0^\an$. 
Here 
$$\OCE_0^\an=(\CE_0|_{U}, \tilde g_0)=(0,\tilde g_0)$$ 
is the image of $\OCE_0=(\CE_0, g_0)$ in 
$\wh \Div (\CU^\an)_\eqv$. 

Note that that the map $\wh\Div (\CU/k)_\rmod \to \wh \Div (\CU^\an)_\eqv$
keeps the partial order of effectivity, so it sends Cauchy sequences to Cauchy sequences. 
To prove that the map is well-defined, it suffices to prove that $\wh \Div (\CU^\an)_\eqv$ is complete under the boundary topology. 
Let $\{(\CD_i,\tilde g_i)\}_{i\geq1}$ be a Cauchy sequence in $\wh \Div (\CU^\an)_\eqv$. 
Then we have $\CD_i=\CD_1$ for all $i$, and there is a sequence $\{\epsilon_i\}_i$ of positive rational numbers converging to 0 such that
$$
-\epsilon_i \tilde g_0 \leq \tilde g_i-\tilde g_j \leq \epsilon_i \tilde g_0,\quad \forall\, j\geq i\geq1.
$$
Note that $\tilde g_0$ is continuous on $\CU^\an$ and thus bounded on any compact subset of $\CU^\an$.
Then $\{\tilde g_i-\tilde g_1\}_i$ is uniformly convergent (to a continuous function) on any compact subset of $\CU^\an$.
As $\CU^\an$ is locally compact, the sequence $\{\tilde g_i-\tilde g_1\}_i$ is pointwise convergent to a continuous function on $\CU^\an$.
Then $\tilde g_i=\tilde g_1+(\tilde g_i-\tilde g_1)$ converges to a Green function of $\CD_1$ on $\CU^\an$.
This gives the limit of $\{(\CD_i,\tilde g_i)\}_{i\geq1}$ in $\wh\Div (\CU^\an)_\eqv$.
Therefore, $\wh\Div (\CU^\an)_\eqv$ is complete, and the first map of the proposition is well-defined.

In the definition
$$\wh\Div (\CU/k)_\rmod=\lim_{\substack{\lra\\ \CX}}\wh\Div(\CX,\CU),$$
we can replace each $\CX$ by its normalization in $\CU$, so that $\CX$ is integrally closed in $\CU$.
By Lemma \ref{effectivity1}, an element of $\wh\Div (\CU/k)_\rmod$ is effective if and only if its image in $\wh\Div (\CU^\an)_\eqv$ is effective.
As a consequence, for any sequence $\{\OCD_i\}_i$ of $\wh\Div (\CU/k)_\rmod$, if the image of $\{\OCD_i\}_i$ in $\wh\Div (\CU^\an)_\eqv$ is a Cauchy sequence equivalent to 0, then $\{\OCD_i\}_i$ is a Cauchy sequence in $\wh\Div (\CU/k)_\rmod$ equivalent to 0.
This proves the injectivity. 
\end{proof}

\subsection{Essentially quasi-projective case}

\kkk
Let $X$ be a flat and essentially quasi-projective integral scheme over $k$.
Recall 
$$
\wh\Div(X/k)=\lim_{\substack{\lra \\ \CU}}\wh \Div (\CU/k),$$
$$
\wh\CaCl(X/k)=\lim_{\substack{\lra \\ \CU}}\wh \CaCl (\CU/k).
$$
Here the limits are over quasi-projective models $\CU$ of $X$ over $k$. 

Note that we have already had an injection
$$\wh \Div (\CU/k)\lra \wh \Div (\CU^\an)_\eqv$$
for quasi-projective models $\CU$ of $X$.
Its direct limit gives an injection
$$
\wh \Div (X/k)\lra \varinjlim_\CU \wh \Div (\CU^\an)_\eqv.
$$
Composing with the map
$$
\varinjlim_\CU \wh \Div (\CU^\an)_\eqv\lra \wh \Div (X^\an)_\eqv,
$$
we get a map 
$$
\wh \Div (X/k)\lra \wh \Div (X^\an)_\eqv. 
$$
This is the map in Proposition \ref{injection1}. 
Now  we are ready to finish the proof of the proposition. 

\begin{proof}[Proof of Proposition \ref{injection1}: essentially quasi-projective case]
Similar to the quasi-projective case, it suffices to prove the injectivity of
$$
\wh \Div (X/k)\lra \wh \Div (X^\an)_\eqv. 
$$
By the above composition, it suffices to prove that the map
$$
\varinjlim_\CU \wh \Div (\CU^\an)_\eqv\lra \wh \Div (X^\an)_\eqv
$$
is injective.
By Lemma \ref{density}, $X^\an\to \CU^\an$ is injective with a dense image. 
Then it suffices to prove that the map
$$
\Phi:\varinjlim_\CU \Div (\CU)\lra \Div (X)
$$
is injective.

Fix a quasi-projective model $\CU_0$ of $X$.
By Lemma \ref{models2}, in the above limits, we can take 
$\{\CU\}$ to be the inverse system of open subschemes of $\CU_0$ containing $X$.
If $\CD$ is an element in the kernel of $\Phi$, then we can assume that $\CD$ lies in 
$\Div (\CU)$ for some $\CU$. At any point $x\in X$, $\CD$ is defined by a single equation $f$ in a neighborhood of $x$ in $\CU$. By assumption, $f$ is invertible in $\CO_{\CU,x}$, so $f$ is invertible on a neighborhood $\CV_x$ of $x$ in $\CU$, or equivalently $\CD$ is 0 on $\CV_x$. Taking unions of $\CV_x$ for all $x\in X$, we see that $\CD$ is 0 on an open neighborhood of $X$ in $\CU_0$.
Thus $\CD=0$. This proves the injectivity of $\Phi$.
The proof is complete.
\end{proof}

\section{Analytification of adelic line bundles} \label{sec analytification line bundles}

The adelic line bundles in \S\ref{sec adelic line bundles} induce norm-equivariant metrized line bundles on Berkovich spaces. The goal of this section is to study this analytification process. 
The main result is as follows:

\begin{prop}  \label{injection2}
\kkk
Let $X$ be a flat and essentially quasi-projective integral scheme over $k$.
There is a canonical, fully faithful functor 
$$\wh \Picc (X/k)\lra \wh\Picc (X^\an)_\eqv,$$
which induces an injective map
$$\wh\Pic (X/k) \lra \wh \Pic (X^\an)_\eqv.$$
\end{prop}

Most of the process is parallel and implied by the analytification of adelic divisors in \S\ref{sec analytification divisors}. 
We will include it for the sake of readers.
As in the case of adelic divisors, we will construct the maps and prove the injectivity in the order of projective case, quasi-projective case, and essentially quasi-projective case.
In the end, we will consider the canonical measures on Berkovich spaces induced by this process.

\subsection{Projective case} 

\kkk
\ccc
Let $\CX$ be a projective variety over $k$. 
Then there is a canonical functor 
$$\wh\Picc (\CX)\lra \wh \Picc (\CX^\an)_\eqv$$
and a canonical map 
$$\wh\Pic (\CX)\lra \wh \Pic (\CX^\an)_\eqv.$$
This is very similar to that construction in \S\ref{sec analytification divisors}. It is a consequence of the latter for choosing a rational section $s$ of a line bundle on $\CX$ and converting metrics $\|s\|$ to Green functions $-\log\|s\|$.

For importance, we sketch the definition here. 
Let $\CLL$ be a hermitian line bundle on $\CX$.
We need to define a metric of $\CL$ on $\CX^\an$. 
The metric of the fibers of $\CL$ on $\CX^\an_\infty=\CX^\an_\RR$ are given by the original hermitian metric, and it extends to $\CX^\an[\infty]$ by norm-equivariance.
For the metric of $\CL$ at a point $x\in \CX^\an[\mathrm{f}]$,
let $\phi_x^\circ:\Spec R_x\to \CX$ be the $k$-morphism extending the $k$-morphism 
$\phi_x:\Spec H_x\to \CX$ under the valuative criterion. 
Then $(\phi_x^\circ)^*\CL$ is a free module over $R_x$ of rank 1.
Let $s_x$ be the basis of this free module. 
Define the metric of $\CL(x)=\phi_x^*\CL$ by setting $\|s_x\|=1$.
The continuity of the metric is a consequence of Lemma \ref{green's functions}.

\subsection{Quasi-projective case} 

\kkk
Let $\CU$ be a quasi-projective variety over $k$.
We are going to have a canonical functor
$$\wh \Picc (\CU/k)\lra \wh\Picc (\CU^\an)_\eqv$$
and a canonical map 
$$\wh \Pic (\CU/k)\lra \wh\Pic (\CU^\an)_\eqv.$$

The functor is described as follows. 
Recall from \S\ref{sec adelic line bundles} that an object of $\wh\Picc (\CU/k)$ is a sequence $\CLL=(\CL,(\CX_i,\overline \CL_i, \ell_{i})_{i\geq1})$. Resume the other notations for this sequence in \S\ref{sec adelic line bundles}. 
Note that each $\overline\CL_i$ induces a metric $\|\cdot\|_i^*$ of $\CL_{i}$ on $\CX_i^\an$.
By the isomorphism $\ell_{i}:\CL\to \CL_i|_\CU$, and by restriction, we get a metric $\|\cdot\|_i$ of $\CL$ on $\CU^\an$.  
We will see that the Cauchy condition implies that these metrics converge pointwise to a continuous metric $\|\cdot\|$ of $\CL$ on $\CU^\an$. 
Then $\CLL^\an:=(\CL,\|\cdot\|)$ defines an element of 
$\wh \Picc(\CU^\an)_\eqv$, which is the desired image of the functor. 

By the above idea, we prove Proposition \ref{injection2} for quasi-projective varieties.

\begin{proof}[Proof of Proposition \ref{injection2}: quasi-projective case]

We need to prove that the above construction gives a functor 
$$\wh \Picc (\CU/k)\lra \wh\Picc (\CU^\an)_\eqv,$$
and prove that the functor is fully faithful. 
This is more or less a consequence of the quasi-projective case of Proposition \ref{injection1}. We will write some parts of the proof and convert some other parts to the proposition. 

Resume the above notations in the construction of the functor.
Let $\CLL=(\CL,(\CX_i,\overline \CL_i, \ell_{i})_{i\geq1})$ be an object of $\wh \Picc (\CU/k)$.
Denote by $\|\cdot\|_i$ the metric of $\CL$ on $\CU^\an$ induced by $(\CX_i,\overline \CL_i)$.

We first check that the metric $\|\cdot\|_i$ converges pointwise to a metric of $\|\cdot\|_i$ of $\CL$ on $\CU^\an$. 
This is very similar to the quasi-projective case of Proposition \ref{injection1}. 
In fact, from \S\ref{sec adelic line bundles}, the Cauchy condition means that 
there is a sequence $\{\epsilon_j\}_{j\geq 1}$ of positive rational numbers converging to $0$ such that in $\wh\Div (\CU/k)_\rmod$,
$$
-\epsilon_j \OCE_0 \leq \wh \div (\ell_{i}\ell_1^{-1})-\wh \div (\ell_{j}\ell_1^{-1}) \leq \epsilon_j \OCE_0, \quad\ i\ge j\geq 1.
$$
This implies 
$$
-\epsilon_j \tilde g_0 \leq  \log(\|\cdot\|_i/\|\cdot\|_j)
\leq  \epsilon_j \tilde g_0, \quad\ i\ge j\geq 1.
$$
Here $\tilde g_0$ is the Green function of $\CE_0$ on $\CX_0^\an$ induced by $\OCE_0$, which is continuous on $\CU^\an$.
Write $f_i=\log(\|\cdot\|_i/\|\cdot\|_1)$ as a continuous function on $\CU^\an$. 
Then the above condition gives 
$$
-\epsilon_j \tilde g_0 \leq f_i-f_j\leq \epsilon_j \tilde g_0, \quad\ i\ge j\geq 1.
$$
As in the proof of Proposition \ref{injection1}, since $\CU^\an$ is locally compact,
 $f_i$ converges pointwise to a continuous function $f$ on $\CU^\an$. 
As a consequence, $\|\cdot\|_i$ converges pointwise to a continuous metric $\|\cdot\|$, and 
$$
-\epsilon_j \tilde g_0 \leq  \log(\|\cdot\|/\|\cdot\|_j)
\leq  \epsilon_j \tilde g_0, \quad\ j\geq 1.
$$
This gives the functor image $\overline \CL^\an=(\CL,\|\cdot\|)$.

To check that it is indeed a fully faithful functor, let 
$\CLL'$ be another 
object of $\wh \Picc (\CU/k)$ with image $\CLL'^\an$ in $\wh\Picc (\CU^\an)_\eqv$. 
We need to prove that there is a canonical isomorphism 
$$
\Hom(\CLL',\CLL)
\lra \Hom(\CLL'^\an,\CLL^\an).
$$
This is equivalent to a canonical isomorphism
$$
\Hom(\ol\CO_{\CX_0}, \CLL'^\vee\otimes\CLL)
\lra \Hom(\ol\CO_{\CU},(\CLL'^\an)^\vee\otimes\CLL^\an).
$$
Here $\ol\CO_{\CX_0}=(\CO_{\CU},(\CX_0, \ol\CO_{\CX_0}, 1))$
and $\ol\CO_{\CU}=(\CO_{\CU},\|\cdot\|_0)$ are the neural elements, where $\|\cdot\|_0$ is defined by $\|1\|_0=1$.

Replacing $\CLL'^\vee\otimes\CLL$ by $\CLL$, it suffices to prove that there is a canonical isomorphism
$$
\Phi:\Hom(\ol\CO_{\CX_0}, \CLL)
\lra \Hom(\ol\CO_{\CU}, \CLL^\an).
$$
Write $\CLL=(\CL,(\CX_i,\overline \CL_i, \ell_{i})_{i\geq1})$ as above.

Elements of both sides of $\Phi$ are represented by regular sections $s$ of $\CL$ everywhere non-vanishing on ${\CU}$.
Such a section $s$ gives an element of the right-hand side of $\|s\|=1$ on $\CU^\an$.
The section $s$ gives an element of the left-hand side of the Cauchy sequence 
$\{\wh\div(\ell_i\ell_1^{-1})+\wh\div_{(\CX_1,\CLL_1)}(s)  \}_{i\geq1}$
of $\wh \Div (\CU/k)_\rmod$ converges to 0 in $\wh \Div (\CU/k)$
under the boundary topology. These two conditions on $s$ are equivalent since
$$\wh\Div (\CU/k) \lra \wh \Div (\CU^\an)_\eqv$$
is injective by Proposition \ref{injection1}.
The proof is complete.
\end{proof}

\subsection{Essentially quasi-projective case}

\kkk
Let $X$ be a flat and essentially quasi-projective integral scheme over $k$.
Recall that 
$$ \wh\Pic (X/k) =\lim_{\substack{\lra \\ \CU}}\wh \Pic (\CU/k),$$
$$ \wh\Picc (X/k) =\lim_{\substack{\lra \\ \CU}}\wh \Picc (\CU/k).$$
Here the limits are over quasi-projective models $\CU$ of $X$ over $k$. 
Note that we have a fully faithful functor
$$\wh \Picc (\CU/k)\lra \wh \Picc (\CU^\an)_\eqv$$
for quasi-projective models $\CU$ of $X$.
Its direct limit gives a  fully faithful functor
$$
\wh \Picc (X/k)\lra \varinjlim_\CU \wh \Picc (\CU^\an)_\eqv.
$$
Composing with the functor
$$
\varinjlim_\CU \wh \Picc (\CU^\an)_\eqv\lra \wh \Picc (X^\an)_\eqv,
$$
we get a functor 
$$
\wh \Picc (X/k)\lra \wh \Picc (X^\an)_\eqv. 
$$
This is the functor in Proposition \ref{injection2}. 
Now  we are ready to finish the proof of the proposition.

\begin{proof}[Proof of Proposition \ref{injection2}: essentially quasi-projective case]
It suffices to prove that the functor
$$
\varinjlim_\CU \wh \Picc(\CU^\an)_\eqv\lra \wh \Picc (X^\an)_\eqv
$$
is fully faithful. 
Similarly, by Lemma \ref{density}, $X^\an\to \CU^\an$ is injective with a dense image, so it suffices to prove that the functor
$$
\Psi:\varinjlim_\CU \Picc(\CU)\lra \Picc (X)
$$
is fully faithful. 

Fix a quasi-projective model $\CU_0$ of $X$.
By Lemma \ref{models2}, in the above limits, we can take 
$\{\CU\}$ to be the inverse system of open subschemes of $\CU_0$ containing $X$.

To prove that the functor $\Psi$ is fully faithful, it suffices to prove that for any line bundles $\CL,\CL'$ on some open neighborhood of $X$ in $\CU_0$, the canonical map
$$
\varinjlim_{\CU}\Hom(\CL|_{\CU}, \CL'|_{\CU}) \lra \Hom(\CL|_{X}, \CL'|_{X})
$$
is an isomorphism. 
The map is isomorphic to
$$
\varinjlim_{\CU}\Gamma(\CU, \CL^\vee\otimes \CL') \lra \Gamma(X,\CL^\vee\otimes \CL').
$$
The injectivity is clear as both sides are subgroups of rational sections of $\CL^\vee\otimes \CL'$ on $X$. 
For the surjectivity, it suffices to prove that if a rational section $s$ of $\CL^\vee\otimes \CL'$ is regular and nowhere vanishing on $X$, then it is regular and nowhere vanishing on a neighborhood of $X$ in $\CU_0$.
In fact, for any $x\in X$, as $s$ is regular and non-vanishing at $x$, it is so at an open neighborhood $\CV_x$ of $x$ in $\CU_0$. Take unions of $\CV_x$ for all $x\in X$. It gives an open neighborhood of $X$ in $\CU_0$ satisfying the requirement. 
This finishes the proof.
\end{proof}

\subsection{Consequence on shrinking the underlying scheme} 

A quick consequence of Proposition \ref{injection1} and Proposition \ref{injection2}
is the following injectivity result. 

\begin{cor}  \label{shrink}
\kkk
Let $f: X\to Y$ be a morphism of flat and essentially quasi-projective integral schemes over $k$. 
Assume that $X$ and $Y$ are normal, and $f$ induces an isomorphism $k(Y)\to k(X)$ of the function fields. 
Then the canonical maps 
$$\wh\Div (Y/k) \lra \wh \Div (X/k),$$
$$\wh\Pic (Y/k) \lra \wh \Pic (X/k),$$
$$\wh\Div (Y^\an) \lra \wh \Div (X^\an),$$
$$\wh\Pic (Y^\an) \lra \wh \Pic (X^\an)$$
are injective.
\end{cor}
\begin{proof}
By Proposition \ref{injection1} and Proposition \ref{injection2}, it suffices to prove the last two maps are injective.

We claim that the third map implies the injectivity of the fourth map.
In fact, by the isomorphism between $\wh\Pic$ and $\wh\CaCl$, the fourth map is isomorphic to the canonical map
$$\wh\CaCl (Y^\an) \lra \wh \CaCl (X^\an).$$
As $k(Y)\to k(X)$ is an isomorphism, the canonical map 
$$\wh\Pr(Y^\an) \lra \wh \Pr(X^\an)$$
is surjective.
Then the injectivity of the third map implies that of the fourth. 

Now  we prove the injectivity of the third map
$$\wh\Div (Y^\an) \lra \wh \Div (X^\an).$$
Assume that an arithmetic divisor $(D,g_D)$ on $Y^\an$ lies in the kernel of this map.
Then $g_D$ is zero on $(\Spec k(Y))^\an$. 
By Lemma \ref{density}, $g_D$ is zero on $Y^\an$. 
Note that $g_D$ has logarithmic singularity along $|D|^\an$ in $Y^\an$. 
This implies that $|D|^\an$ is empty, and thus $D=0$. 
It finishes the proof. 
\end{proof}

\section{Restricted analytic spaces} \label{sec restricted}

\kkk
Let $X$ be a flat and essentially quasi-projective integral scheme over $k$.
The Berkovich space $X^\an=(X/k)^\an$ is intrinsic and functorial.
Moreover, the analytification map 
$$\wh\Pic (X/k) \lra \wh \Pic ((X/k)^\an)_\eqv$$
is functorial in the sense that it is compatible with the 
functoriality maps listed in \S \ref{sec functoriality}. 

However, a disadvantage is that the space is too large and too abstract to work on, mainly because it contains ``too many" redundant points. 
The goal here is to consider a smaller subspace $(X/k)^\ran$ of $(X/k)^\an$, as the union of some distinguished fibers, which is sufficient for many applications.
It turns out that metrized line bundles on $(X/k)^\ran$ are very close to the adelic line bundles of Zhang \cite{Zha2} for projective varieties over number fields reviewed in \S\ref{app sec adelic}.
In this case, Theorem \ref{image} asserts that the adelic line bundles of \cite{Zha2} are equivalent to our current notions.
We will first introduce the arithmetic case and give a sketch of the geometric case.

\subsection{Arithmetic case}

Recall that $(\Spec \ZZ)^\an=\CM(\ZZ)$ is the set of all multiplicative semi-norms of $\ZZ$. 
Define $(\Spec \ZZ)^{\ran}$ to be the subspace of $(\Spec \ZZ)^{\an}$ of non-trivial standard absolute values $|\cdot |_v$ of $\ZZ$. 
Hence $(\Spec \ZZ)^{\ran}$ is bijective to $M_\QQ=\{\infty, 2, 3,5,7,\cdots\}$ and endowed with the discrete topology. 

Let $X$ be a scheme over $\ZZ$. 
There is a structure map $X^\an \to (\Spec \ZZ)^{\an}$. 
Define the \emph{restricted analytic space $X^\ran =(X/\ZZ)^\ran$ associated to} $X/\ZZ$ to be the preimage of $(\Spec \ZZ)^\ran$ under the map $X^\an \to (\Spec \ZZ)^{\an}$. It follows that
$$
X^\ran =\coprod_{v\in (\Spec \ZZ)^\ran} X_v^\an,
$$
where $X_v^\an$ is the fiber of $X^\an$ above $v$.
Then $X_v^\an$ is canonically homeomorphic to $X_{\QQ_v}^\an=(X_{\QQ_v}/\QQ_v)^\an$,
the Berkovich space associated to $X_{\QQ_v}$ over the complete field $\QQ_v$.
The topology on $X^\ran$ is induced by the disjoint union so that 
each $X_v^\an$ is both open and closed in $X^\ran$.

Define an \emph{arithmetic divisor} on $X^\ran$ to be a pair $(D,g_D)$, where $D$ is a Cartier divisor on $X$, and $g_D$ is a \emph{Green function} of $D$ on $X^\ran$, i.e. a continuous function $g:X^\ran\setminus |D|^\ran \to \RR$ with logarithmic singularity along $D$ in the sense that, for any rational function $f$ on a Zariski open subset $U$ of $X$ satisfying $\div(f)=D|_U$, the function $g+\log |f|$ can be extended to a continuous function on $U^\ran$.

An {arithmetic divisor} on $X^\ran$ is called {\em principal} if it is of the form 
$$\wh\div_{X^\ran}(f):=(\div(f), -\log |f|)$$ 
for some nonzero rational function $f$ on $X$.

Denote by $\wh\Div(X^\ran)$ (resp. $\wh\Pr(X^\ran)$) the \emph{group} of arithmetic divisors (resp. principal arithmetic divisor) on $X$. 
Define 
$$
\wh\CaCl(X^\ran):=\wh\Div(X^\ran)/\wh\Pr(X^\ran). 
$$

Define a \emph{metrized line bundle} on $X^\ran$ to be a pair $(L,\|\cdot\|)$, where $L$ is a line bundle on $X$, and $\|\cdot\|$ is a continuous metric of fibers of $L$ on $X^\ran$. 
This is similar to the original case; we omit the details.

Denote by $\wh\Pic(X^\ran)$ the \emph{group} of isometry classes of metrized line bundle on $X$, and $\wh\Picc(X^\ran)$ the \emph{category} of metrized line bundle on $X$ whose morphisms are isometries.
There is a canonical isomorphism
$$
\wh\CaCl(X^\ran)\lra \wh\Pic(X^\ran). 
$$

Finally, the motivation for restricted analytic spaces is as follows. 

\begin{prop}  \label{injection3}
Let $X$ be a flat and essentially quasi-projective integral scheme over $\ZZ$.
There are canonical injective maps
$$\wh\Div (X) \lra \wh \Div (X^\ran),$$
$$\wh\Pic (X) \lra \wh \Pic (X^\ran),$$
and a canonical fully faithful functor 
$$\wh \Picc (X)\lra \wh\Picc (X^\ran).$$
\end{prop}
\begin{proof}
This is a consequence of Proposition \ref{injection1} and Proposition \ref{injection2}.
For example, the first map is obtained as the composition 
$$\wh\Div (X) \lra \wh \Div (X^\an)_\eqv \lra \wh \Div (X^\ran).$$
Here the first arrow is injective by Proposition \ref{injection1}. 
The second arrow is injective since a norm-equivalent Green function on $X^\an$ is determined by its restriction to $X^\an\setminus X^\an_{\iota(\Spec \ZZ)}$ by Lemma \ref{density}, and the latter is determined by its restriction to $X^\ran$ by norm-equivariance. 
\end{proof}

\subsection{Comparison in the projective case}

Let us compare the current theory with the original theory of adelic line bundles on projective varieties over a number of fields of Zhang \cite{Zha2}. 
We refer to \S\ref{app sec adelic} for an overview of old theory, and we will follow the terminology of that section.

Let $X$ be a projective variety over a number field $K$. Denote by $M_K$ the set of places of $K$. 
Recall that in \S\ref{app sec adelic}, an $M_K$-metrized line bundle $\OL=(L, (\|\cdot\|_v)_v)$ on $X$ consists of a line bundle $L$ on $X$, and an $M_K$-metric $(\|\cdot\|_v)_v$ of $L$ on $X_{M_K}=\coprod_{v\in M_K} X_{K_v}$, i.e. a collection 
 of {m-continuous} $K_v$-metrics $\|\cdot\|_v$  of $L_{K_v}$ on $X_{K_v}$ over all places $v$ of $K$. 
Here ``m-continuous metrics'' are introduced in \S \ref{subsec-metrics} as uniform limits of  metrics induced by projective models in the non-archimedean case. 
An adelic line bundle on $X$ in the sense of \cite{Zha2} is
a metrized line bundle on $X_{M_K}$ satisfying the coherence condition, i.e. over all but finitely many places, the metric is induced by a single projective integral model over an open subscheme of $\Spec O_K$. 

To avoid confusion, denote by $\wh\Picc(X_{M_K})^*$ the category of $M_K$-metrized line bundles on $X$, 
and denote by $\wh\Picc(X)^*$ the category of adelic line bundles on $X$ introduced in \S\ref{app subsec adelic}. Note that we use the notations $\wh\Picc(X_{M_K})^*$ and $\wh\Picc(X)^*$ instead of $\wh\Picc(X_{M_K})$ and $\wh\Picc(X)$ to distinguish them from our current set of notations. 

To compare it with the current theory, view $X$ as a projective variety over
$\QQ$, which is generally not geometrically integral.

Denote by $\wh\Picc(X^\ran)_\mathrm{coh}$ the full subcategory of $\wh\Picc(X^\ran)$ whose objects are adelic line bundles on $X$ satisfying the coherence condition. 
Here the \emph{coherence condition} of a metrized line bundle $\OL$ on $X^\ran$ is the existence of an open subscheme $\CV$ of $\Spec\ZZ$, a projective and flat morphism $\CU\to \CV$ whose generic fiber is isomorphic to $X\to \Spec \QQ$, and a line bundle $\CL$ on $\CU$ endowed with an isomorphism $\CL_\QQ\to L$ over $X$, such that the metric of $\OL$ on $X^\an_p\subset X^\ran$ is equal to the metric of $L$ on $X^\an_p$ induced by $\CL$ for all primes $p\in \CV$.

The following result asserts that we have equivalences  
$$\wh\Picc(X/\ZZ)\lra \wh\Picc(X^\ran)_\mathrm{coh}\lra \wh\Picc(X)^*$$ 
of categories. 
Therefore, we have three equivalent definitions of adelic line bundles on $X$ in the projective case.

\begin{thm}  \label{image}
Let $X$ be a projective variety over a number field $K$.
Then there is a canonical equivalence 
$$ \wh\Picc(X^\ran)\lra \wh\Picc (X_{M_K})^*,$$ 
which induces equivalences 
$$\wh\Picc(X/\ZZ)\lra \wh\Picc(X^\ran)_\mathrm{coh}\lra \wh\Picc(X)^*.$$ 
\end{thm}
\begin{proof}
We first reduce the problem to the case $K=\QQ$. 
View $X$ as a projective variety over $\QQ$, and denote this variety by $X'$ to avoid confusion. 
For any place $v$ of $\QQ$, there is a canonical isomorphism  
$$
K\otimes_\QQ \QQ_v \simeq \prod_{w|v} K_w,
$$
where the product is over all places $w$ of $K$ above $v$.
This induces canonical isomorphisms
$$
X' \otimes_\QQ \QQ_v \simeq \coprod_{w|v} X\otimes_K K_w, \qquad
X'(\ol \QQ_v)=\prod_{w|v} X(\OK_w).
$$
As a consequence, we have a canonical equivalence
$$
\wh\Picc (X_{M_K})^*\lra \wh\Picc (X'_{M_\QQ})^*, \qquad
\wh\Picc (X)^* \lra \wh\Picc (X')^*.$$
In other words, the categories are essentially independent of the number field $K$. 
Then we can and will assume that $K=\QQ$ in the following, which will simplify notations significantly. 

Next, we treat the functor $\wh\Picc(X^\ran)\to \wh\Picc (X_{M_\QQ})^*$. 
It suffices to make a canonical bijection, for a fixed line bundle $L_{\QQ_v}$ on $X_{\QQ_v}$ at a place $v$ of $\QQ$,  from the set of continuous metrics of $L_{\QQ_v}$ on $X_v^\an$ to the set of m-continuous $\QQ_v$-metrics of $L_{\QQ_v}$ on $X_{\QQ_v}$. Recall that ``m-continuous metrics'' are introduced in \S \ref{subsec-metrics} as uniform limits of model metrics. 
The bijection is automatic for $v=\infty$, so we assume $v\neq\infty$. 
Denote by $|X_{\QQ_v}|_0$ the set of closed points of $X_{\QQ_v}$. 
We have an injection $|X_{\QQ_v}|_0\to X_v^\an$ with a dense image. 
The canonical maps $X(\ol\QQ_v) \to |X_{\QQ_v}|_0\to X_v^\an$ induce maps between the fibers of $L_{\QQ_v}$ in the opposite directions, and thus transfer 
 metrics of $L_{\QQ_v}$ on $X_v^\an$ to $\QQ_v$-metrics of $L_{\QQ_v}$ on $X_{\QQ_v}$. 
It remains to check that continuous metrics of $L_{\QQ_v}$ on $X_v^\an$ 
 correspond exactly to m-continuous $\QQ_v$-metrics of $L_{\QQ_v}$ on $X_{\QQ_v}$. 

The situation is very similar to that in \S\ref{app subsec measure}.
In fact, a projective model of $(X_{\QQ_v}, L_{\QQ_v})$ induces a model metric on each side which corresponds to each other.  
Taking quotients of the metrics by a fixed model metric of $L_{\QQ_v}$, it is reduced to the case that $L_{\QQ_v}$ is the trivial line bundle, and then the process
$\|\cdot\|\mapsto -\log \|1\|$ transfers metrics to functions. 
Thus it is reduced to make a canonical bijection from the set of continuous functions on $X_v^\an$ to the set of m-continuous functions on $|X_{\QQ_v}|_0$. 
Note that m-continuous functions are defined as uniform limits of model functions. 
This bijection is a consequence of Gubler \cite[Theorem 7.12]{Gub3} (cf. \cite[Lem. 3.5]{Yua1}), which asserts that the space of model functions on $X^\an$ is dense in the space of continuous functions on $X^\an$ under the uniform topology.

Hence, we have defined $\wh\Picc(X^\ran)\to \wh\Picc (X_{M_\QQ})^*$ and proved that it is an equivalence. 
Moreover, it induces an equivalence 
$\wh\Picc(X^\ran)_\mathrm{coh}\to \wh\Picc (X)^*$. 
 
By Proposition \ref{injection3}, we have a fully faithful  functor
 $\wh\Picc(X)\to \wh\Picc(X^\ran)$. 
It remains to prove that its essential image is 
$\wh\Picc (X^\ran)_\mathrm{coh}$.

By definition, 
$$
\wh \Picc (X)= \varinjlim_\CU \wh \Picc (\CU),
$$
where the limit is over quasi-projective models $\CU$ of $X$ over $\ZZ$.
Replacing $\CU$ by an open subscheme if necessary, we can assume that there is a projective and flat morphism $\CU\to \CV$ for some open subscheme 
$\CV$ of $\Spec \ZZ$. 
It is reduced to characterize the essential image of the functor
$\wh \Picc (\CU)\to \wh\Picc (X^\ran).$

Let $\CLL=(\CL,(\CX_i,\CLL_i,\ell_i)_{i\geq1})$ be an object of $\wh \Picc (\CU)$, which has underlying line bundle $\CL$ on $\CU$. Denote by $\OL=(L, (\|\cdot\|_v)_{v})$ the image of 
$\CLL$ in $\wh\Picc (X^\ran)$, with underlying line bundle $L=\CL|_X$ on $X$.  
Denote by $\|\cdot\|_{i,v}$ (resp. $\|\cdot\|_{v}^\circ$) the metric of $L$ on $X_v^\an$ induced by $(\CX_i,\CLL_i)$ (resp. $(\CU,\CL)$) for $v\leq \infty$ (resp. $v\in\CV$).
By definition, $\|\cdot\|_v$ is the pointwise limit of $\|\cdot\|_{i,v}$ on $X_v^\an$.

As the base change of $(\CX_i,\CLL_i)$ to $\CV$ is isomorphic to $(\CU,\CL)$, they 
induce the same metrics of $L$ at any closed point $v\in \CV$. 
It follows that $\|\cdot\|_{i,v}=\|\cdot\|_{v}^\circ=\|\cdot\|_{v}$ for $v\in \CV$. 
This proves that the metric of $\OL$ satisfies the coherence condition over $\CV$. 
Therefore, the essential image of $\wh\Picc (X)\to \wh\Picc (X^\ran)$ lies in 
$\wh\Picc (X^\ran)_\mathrm{coh}$.

We claim that the convergence of $\|\cdot\|_{i,v}$ to $\|\cdot\|_v$ is the uniform convergence on $X_v^\an$ for any place $v\leq \infty$.
In fact, for the sake of the boundary topology, take any projective model $\CX_0$ of $\CU$, and take the arithmetic divisor 
$\ol\CE_0=(\CE_0, 1)$ 
over $\Spec\ZZ$, where $\CE_0=(\Spec\ZZ)\setminus\CV$ is endowed with the reduced structure. 
Take $\ol\CF_0$ to be the pull-back of $\ol\CE_0$ to $\CX_0$.
Use $(\CX_0,\ol\CF_0)$ to define the boundary topology of $\wh \Picc (\CU)_\rmod$. 
There is also a boundary topology of $\wh \Picc (X^\ran)$ defined by $\ol\CF_0^\an$.
Note that the Green function $\wt g$ of $\ol\CF_0^\an$ is 0 on $X_v^\an$ for any $v\in\CV$ and a positive constant on $X_v^\an$ for $v\notin \CV$ (including $v=\infty$). 
As a consequence, the convergence of $\|\cdot\|_{i,v}$ to $\|\cdot\|_v$ is the uniform convergence.

It remains to prove that $\wh\Picc (X)\to \wh\Picc (X^\ran)_\mathrm{coh}$ is essentially surjective. 
Given any metrized line bundle $\OL=(L, (\|\cdot\|_v)_{v})$ on $X^\ran$ satisfying the coherence condition over an open subscheme $\CV$ of $\Spec \ZZ$, we will prove that $\OL$ is isomorphic to the image of some adelic line bundle $\CLL=\{(\CX_i,\CLL_i,\ell_i)\}_{i\geq1}$ on some quasi-projective model $\CU$ of $X$.
In fact, by the coherence condition, there is a model $(\CU\to \CV,\, \CL)$ of 
$(X\to\Spec\QQ,\, L)$ inducing the metric of $\OL$ above $\CV$.
For any $v\notin \CV$, the metric $\|\cdot\|_v$ is continuous. 
By \cite[Thm. 7.12]{Gub3} again, 
$\|\cdot\|_v$ is a uniform limit of model metrics 
$(\|\cdot\|_{i,v})_{i\geq1}$ if $v$ is finite.
For any $i\geq1$, we can find an integral model $(\CX_i,\CLL_i)$ of $(X,L)$ which 
extends $(\CU,\CL)$ and induces the metric $(\|\cdot\|_{i,v})_v$ of $L$. 
This gives the adelic line bundle $\CLL$.
It finishes the proof.
\end{proof}

Let $X$ be a projective variety over a number field $K$.
One can also check that nefness in $\wh\Picc(X)$ (defined in \S\ref{app sec adelic}), nefness in $\wh\Picc(X/\ZZ)$,
and strong nefness in $\wh\Picc(X/\ZZ)$ are equivalent. 
Then the notions of integrable adelic line bundles in these two categories are also equivalent.

\begin{rmk}
The restricted analytic space $X^\ran$ seems very artificial, but it has a functorial interpretation in terms of the Gelfand spectrum. 
Let $X$ be a projective variety over a number field $K$.
We call an adelic line bundle $\OL$ on $X$ \emph{vertical} if the underlying line bundle $L$ is isomorphic to the trivial line bundle $\CO_X$. 
Then the space $V(X)=\bigoplus_v C(X_v^\an)$ exactly consists of the functions $(-\log\|1\|_v)_v$ for all vertical adelic line bundles $\OL$ on $X$. 
This can be viewed as a space of  ``vertical adelic divisors" on $X$. 
Moreover, $V(X)$ has a natural ring structure. 
It is complete under the natural topology, where a sequence $f_i=(f_{i,v})_v$  converges to $f=(f_v)_v$
if there is a finite set  $S$ of places of $\QQ$ such that $f_{i,v}=f_v$=0 for all $v\notin S$, and that $f_{i,v}$ converges to $f_v$ for places $v$. 
Then we have
  $$X^\ran=\Hom _{\cont}(V(X), \BR), \qquad V(X)=C_c (X^\ran).$$
Here ``$\Hom$" denotes the set of continuous ring homomorphisms, and the first isomorphism is by Gelfand spectrum. 
\end{rmk}

\subsection{Function field case}  \label{sec function field analytic}

Let $k$ be any field. Let $X$ be a scheme over $k$.
Recall that the analytic space $X^\an=(X/k)^\an$ is (Zariski locally) given by multiplicative semi-norms trivial over $k$.
To define a restricted subspace of $X^\an$ as in the arithmetic case, we need extra data to get a global field. 
This fits the setting at the end of \S\ref{sec adelic general}.

Let $B$ be a projective regular curve over $k$. 
Denote the function field by $K=k(B)$. 
Any closed point $v\in B$ gives a normalized absolute value $|\cdot |_v=\exp(-\ord_v)$ of $K$. 
Define $B^{\ran}=(B/k)^{\ran}$ to be the subspace of $B^{\an}=(B/k)^{\an}$ of non-trivial normalized absolute values of $K$. 
Therefore, $B^\ran$ is bijective to the set of closed points of $B$ and endowed with the discrete topology.

Let $X$ be a scheme over $B$ (instead of just over $k$). 
There is a natural map $X^\an \to B^{\an}$. 
Define the \emph{restricted analytic space 
$$X^\ran=(X/B)^\ran=(X/B/k)^\ran$$ 
associated to} $X/B/k$ to be the preimage of $B^\ran$ under the map $X^\an \to B^{\an}$. It follows that
$$
X^\ran =\coprod_{v\in B^\ran} X_v^\an,
$$
where $X_v^\an$ is the fiber of $X^\an$ above $v$.
Then $X_v^\an$ is canonically homeomorphic to $X_{K_v}^\an=(X_{K_v}/K_v)^\an$,
the Berkovich space associated to $X_{K_v}$ over the complete field $K_v$.
The topology on $X^\ran$ is induced by the disjoint union so that 
each $X_v^\an$ is both open and closed in $X^\ran$.

Similar to the arithmetic case, we can define arithmetic divisors and metrized line bundles
over $X^\ran$. Then Proposition \ref{injection3} also holds for any flat and essentially quasi-projective integral scheme $X$ over $B$. 
More precisely, there are canonical injective homomorphisms
$$ \wh\Div (X/k) \lra \wh \Div (X^\ran),$$
$$\wh\Pic (X/k) \lra \wh \Pic (X^\ran),$$
and a canonical, fully faithful functor 
$$\wh \Picc (X/k)\lra \wh\Picc (X^\ran).$$
Recall that from \S\ref{sec adelic general}, we also have canonical isomorphisms
$$
\wh\Div(X/B) \lra \wh\Div(X/k),
$$
$$
\wh\Pic(X/B) \lra \wh\Pic(X/k),
$$
$$
\wh\Picc(X/B) \lra \wh\Picc(X/k).
$$

In the end, we remark that the theory depends on the structure $X/B/k$.
In general, if we are only given $X/k$, then we may use some geometric operations to construct the curve $B$ in the middle. For example, if $X$ is quasi-projective over $k$ of dimension at least 2, take $K=k(t)$ for some transcendental element $t\in k(X)$, which gives a rational map $X\dashrightarrow B$ with $B=\BP^1_k$, and then blow-up $X$ to get a morphism $X\to B$.

\section{Local theory} \label{sec local theory}

Let $X$ be a quasi-projective variety over $\QQ$.
We want to know more about the essential image of the functor 
$$\wh \Picc (X)\lra \wh\Picc (X^\ran).$$
If $X$ is projective,  Theorem \ref{image} gives a satisfactory answer. 
If $X$ is not projective, such a task might be impossible, but the situation simplifies when restricted to fibers of $X^\ran$, i.e. considering the essential image of  
$$\wh \Picc (X)\lra \wh\Picc (X_v^\an)$$
for any place $v$ of $\QQ$. 
This is the motivation of the theory in this section. 

In this section, we will first study the completion process of adelic divisors on Berkovich spaces over complete fields and introduce the Chambert-Loir measure in this situation. 
Most of this section is over complete fields, except that Theorem \ref{singularity global} gives a sufficient condition for a metrized line bundle over $\ZZ$ to be an adelic line bundle over $\ZZ$. 

\subsection{The analytification functor}

Let $K$ be a field complete for a non-trivial valuation $|\cdot|$. 
If $K$ is non-archimedean, denote by $O_K$ the valuation ring of $K$.
If $K$ is archimedean ($K=\CC$ or $\RR$), write $O_K=K$ for convenience.

Let $U$ be a quasi-projective variety over $K$.
By \S\ref{sec adelic general}, we have introduced the groups
$$
\wh\Div(U/\OB),\quad
\wh\CaCl(U/\OB),\quad
\wh\Picc(U/\OB),\quad
\wh\Pic(U/\OB).
$$
Here we understand that the base valued scheme $\OB$ to be $\Spec O_K$ in the non-archimedean case, and to be $(\Spec \RR, i_\mathrm{st})$ or 
$(\Spec \CC, \mathrm{id})$ in the archimedean case.
By abuse of notations, we will write the groups uniformly by
$$
\wh\Div(U/O_K),\quad
\wh\CaCl(U/O_K),\quad
\wh\Picc(U/O_K),\quad
\wh\Pic(U/O_K).
$$

Let $U^\an$ be the (usual) Berkovich analytic space associated with $U$ over $K$.
As in \S\ref{sec analytic objects}, an \emph{arithmetic divisor} on $U^\an$ is pair 
$\OD=(D,g_D)$, where $D$ is a Cartier divisor on $U$, and 
 $g_D:U^\an\setminus |D|^\an\to \RR$ is a Green function of continuous type of $D$ on $U^\an$. 
Similarly, a \emph{metrized line bundle} on $U^\an$ is pair 
$\OL=(L,\|\cdot\|)$, where $L$ is a line bundle on $U$, and $\|\cdot\|$ is a continuous metric of $L$ on $U^\an$. 

Therefore, we have the following groups
$$
\wh\Div(U^\an),\quad
\wh\Pr(U^\an),\quad
\wh\CaCl(U^\an),\quad
\wh\Picc(U^\an),\quad
\wh\Pic(U^\an).
$$
Here $\wh\Div(U^\an)$ (resp. $\wh\Pr(U^\an)$) is the group of
arithmetic divisors (resp. principal arithmetic divisors) on $U^\an$.
And $\wh\Picc(U^\an)$ (resp. $\wh\Pic(U^\an)$) is the category (resp. group) of
{metrized line bundles} on $U^\an$ under isometry.
The group
$$\wh\CaCl(U^\an):=\wh\Div(U^\an)/\wh\Pr(U^\an)$$
is canonically isomorphic to
$\wh \Pic(U^\an).$

The local counterparts of Proposition \ref{injection1} and Proposition \ref{injection2} are as follows. The proof is similar to and easier than the global case, so we omit them. 

\begin{prop}  \label{injection local}
Let $K$ be a field complete for a non-trivial valuation. 
Let $U$ be a quasi-projective variety over $K$.
There are canonical injective maps 
$$\wh \Div (U/O_K)\lra \wh \Div (U^\an),$$
$$\wh \CaCl (U/O_K)\lra \wh \CaCl (U^\an),$$ 
$$\wh \Picc (U/O_K)\lra \wh\Picc (U^\an),$$
$$\wh\Pic (U/O_K) \lra \wh \Pic (U^\an).$$
\end{prop}

In the following, we are going to study the images of the analytification functors.
Denote 
$$\wh \Div (U^\an)_\cptf=\Im(\wh \Div (U/O_K)\to \wh \Div (U^\an)),$$
$$\wh \Picc (U^\an)_\cptf=\Im(\wh \Picc (U/O_K)\to \wh \Picc (U^\an)),$$
$$\wh \Pic (U^\an)_\cptf=\Im(\wh \Picc (U/O_K)\to \wh \Pic (U^\an)).$$
They are compactifications of 
$$\wh \Div (U^\an)_\rmod=\Im(\wh \Div (U/O_K)_\rmod\to \wh \Div (U^\an)),$$
$$\wh \Picc (U^\an)_\rmod=\Im(\wh \Picc (U/O_K)_\rmod\to \wh \Picc (U^\an)),$$
$$\wh \Pic (U^\an)_\rmod=\Im(\wh \Picc (U/O_K)_\rmod\to \wh \Pic (U^\an)).$$
We will first describe the compactification process directly on $U^\an$.

\subsection{Compactified arithmetic divisors}

Let $K$ be a field complete for a non-trivial valuation $|\cdot|$. 
Let $U$ be a quasi-projective variety over $K$.
Recall that \emph{projective model} of $U$ over $O_K$ is a flat and projective integral scheme $\CX$ over $O_K$ together with an open immersion $U\to \CX_K$. 

If $K$ is non-archimedean, an \emph{arithmetic model} (or \emph{integral model}) of a Cartier divisor $D$ of $U$ is a pair $(\CX,\CD)$, where $\CX$ is a projective model of $U$ over $O_K$, and $\CD$ is a Cartier $\QQ$-divisor on $\CX$ extending $D$ in that $\CD$ and $D$ have the same image in $\Div(U)_\QQ$.

If $K$ is archimedean, an \emph{arithmetic model} of a Cartier divisor $D$ of $U$ is a pair $(\CX,\CD)$, where $\CX$ is a projective model of $U$ over $O_K=K$, and $\CD=(\wt D,g_{\wt D})$ consisting of a $\QQ$-divisor $\wt D$ on $\CX$ extending $D$ and a  Green function $g_{\wt D}:\CX^\an \setminus |\wt D|^\an\to \RR$ of continuous type of $\wt D$ on $\CX^\an$.

In both cases, the arithmetic model $(\CX,\CD)$ induces a Green function 
$g_\CD$ of $\CD_K$ on $\CX_K^\an$, and thus a Green function $g_{\CD}|_{U^\an}$ of $D$ on $U^\an$ by restriction.
The process is essentially the same as the global case described in \S\ref{sec analytification divisors}. 
The Green function $g_{\CD}|_{U^\an}$ is called a \emph{model Green function}, and the arithmetic divisor 
$(D,g_{\CD}|_{U^\an})$ is called a \emph{model arithmetic divisor} on $U^\an$.
The model Green function or the model arithmetic divisor is called \emph{nef} (or \emph{semipositive}) if either $\CD$ is nef on $\CX$ in the non-archimedean case, or $\CD$ has a semipositive Chern current on $\CX(\CC)$ in the archimedean case (cf. \S\ref{sec hermitian}).

By definition, the image
$$\wh \Div (U^\an)_\rmod=\Im(\wh \Div (U/O_K)_\rmod\to \wh \Div (U^\an))$$
is the group of all model arithmetic divisors on $U^\an$.
It is a natural subgroup of $\wh\Div(U^\an)$. 
Now  we endow it with a boundary topology as in \S\ref{sec adelic divisors}.

By a \emph{boundary divisor} of $U$ over $O_K$, we mean an arithmetic model 
$(\CX_0,\CE_0)$ over $O_K$ of the divisor 0 on $U$  such that 
the support of $\CE_{0,K}$ on $\CX_{0,K}$ is exactly $\CX_{0,K}\setminus U$, and such that the induced Green function $g_{\CE_0}>0$ on $\CX_{0,K}^\an$. 
Then $(\CX_0,\CE_0)$ induces an arithmetic divisor $\OE_0=(0,g_0)$ on $U^\an$. Here $g_0=g_{\CE_0}|_{U^\an}$ is a continuous function on $U^\an$.
Moreover, $g_{\CE_0}$ has a strictly positive lower bound on $X_0^\an$.
We call $\OE_0=(0,g_0)$ a \emph{boundary divisor} of $U^\an$.

Now have a \emph{boundary norm} 
$$\|\cdot\|_{\OE_0}:\wh\Div (U^\an) \lra [0,\infty]$$
defined by 
$$
\|\OD\|_{\OE_0}:=\inf\{\epsilon\in \BQ_{>0}: \ 
 -\epsilon \OE_0 \leq
\OD \leq  \epsilon \OE_0\}.
$$
Here we take the convention that $\inf(\emptyset)=\infty$.
Then $\|\cdot\|_{\OE_0}$ is an extended norm. 
Now we have a \emph{boundary topology on $\wh\Div (U^\an)$} induced by the boundary norm, for which a neighborhood basis at $0$ is formed by
$$
 B(\epsilon, \wh\Div (U^\an)):=\{\OD \in \Divhat(U^\an): \ 
 -\epsilon \OE_0 \leq
\OD \leq  \epsilon \OE_0\}, \quad \epsilon\in \BQ_{>0}.
$$
Here ``$\leq$'' is still given by effectivity. 
By translation, it gives a neighborhood basis at any point.
The topology does not depend on the choice of $\OE_0$.

By a similar method, we have a boundary topology over
$\wh\Div (U^\an)_\rmod$, which is the same as the subspace topology induced from $\wh\Div (U^\an)$. 

By construction, 
$$\wh \Div (U^\an)_\cptf=\Im(\wh \Div (U/O_K)\to \wh \Div (U^\an))$$
is equal to the {completion} of $\wh \Div  (U^\an)_\rmod$ for the  boundary topology. 
An element of $\wh \Div  (U^\an)_\cptf$ is called a \emph{compactified divisor} on $U^\an$.
A compactified divisor or its Green function is called \emph{strongly nef} (or \emph{strongly semipositive}) if it is a limit of nef model arithmetic divisors under the boundary topology. 
A compactified divisor $\OD$ of $\wh \Div  (U^\an)_\cptf$ or its Green function is called \emph{nef} (or \emph{semipositive}) if there exists a strongly nef element $\OD_0$ of $\wh \Div  (U^\an)_\cptf$ such that $a\OD+\OD_0$ is strongly nef for all positive integers $a$. 

\begin{lem} \label{local subspace}
The space $\wh \Div  (U^\an)$ is complete for the boundary topology and contains $\wh \Div  (U^\an)_\cptf$ as a subspace. 
\end{lem}
\begin{proof}
The second statement is by definition, while the first statement is similar to and easier than Lemma \ref{local subspace2} below. We omit the proof.
\end{proof}

\subsection{Singularity of Green functions}

It turns out that there is a surprisingly explicit description of $\wh \Div  (U^\an)_\cptf$, which is determined by the space of Green functions in it. 
For that purpose, we start with the following spaces of real-valued functions on $U^\an$.
\begin{enumerate}
\item $C(U^\an)$ denotes the space of real-valued continuous functions on $U^\an$;
\item $G(U^\an)$ denotes the space of Green functions on $U^\an$ associated to Cartier divisors of $U$.
\end{enumerate}

By definition, there is a natural injection.
$$
\wh\Div(U^\an) \lra \Div(U)\oplus G(U^\an)
$$ 
and a canonical exact sequence
$$
0\lra C(U^\an) \lra \wh\Div(U^\an) \lra \Div(U) \lra 0.
$$

In terms of the boundary divisor $(0,g_0)$ with $g_0=g_{\CE_0}|_{U^\an}$, we have a boundary topology on $C(U^\an)$ and $G(U^\an)$. The topologies are compatible under inclusion. 
For example, the \emph{boundary topology on $G(U^\an)$} is 
induced by the \emph{boundary norm on $G(U^\an)$} given by 
$$
\|g\|_{g_0}:=\|g/g_0\|_{\sup}=\sup\{|g(x)/g_0(x)|:x\in U^\an\}.
$$
Then a neighborhood basis at $0$ to be formed by
$$
 B(\epsilon, G(U^\an)):=\{g \in G(U^\an): \ 
 -\epsilon g_0 < g <  \epsilon g_0\}, \quad \epsilon\in \BQ_{>0}.
$$
Here the inequalities are understood to hold pointwise away from the loci of the logarithmic singularities.
By translation, it gives a neighborhood basis at any point.

We have the following basic result, which is essentially contained in our previous treatments. 

\begin{lem} \label{local subspace2}
The space $C(U^\an)$ is complete for the boundary topology. 
If $U$ is normal, the space $G(U^\an)$ is complete for the boundary topology. 
\end{lem}
\begin{proof}
This is similar to the quasi-projective model case of the proof of Proposition \ref{injection1}. 
We only treat $G(U^\an)$, as $C(U^\an)$ is similar. 
In fact, let $\{f_i\}_{i\geq 1}$ be a Cauchy sequence in 
$G(U^\an)$. 
Then there is a sequence $\{\epsilon_j\}_{j\geq 1}$ of positive rational numbers converging to $0$ such that 
$$
 -\epsilon_j g_0 <
f_i-f_j<  \epsilon_j g_0,\quad\ i\geq j\geq 1.
$$ 
Note that $g_0$ is continuous on $U^\an$, so $h_i=f_i-f_1$ is bounded on any compact subset of $U^\an$. 
Note that $h_i=f_i-f_1$ has logarithmic singularity along a Cartier divisor 
$D_i$ on $U$. 
By assumption, $U$ is normal, and we can view $D_i$ as a Weil divisor of $U$.
The boundedness of $h_i$ implies $D_i=0$ and thus implies that $h_i$ is continuous on $U^\an$. 

Hence, $\{h_i\}_{i\geq1}$ is a sequence of continuous functions on $U^\an$.
By the boundary norm, $\{h_i/g_0\}_{i\geq1}$ is uniformly convergent, and thus the limit is a continuous function. 
Then $\{h_i\}_{i\geq1}$ pointwise converges to a continuous function $h$ on $U^\an$. 
Then $f_1+h$ is the limit of $\{f_i\}_{i\geq 1}$ in $G(U^\an)$. 
\end{proof}

In order to study $\wh\Div(U^\an)_\cptf$, we introduce the following spaces.
\begin{enumerate}
\item[(3)] $C(U^\an)_\rmod$ denotes the space of \emph{model functions} on $U^\an$, i.e. model Green functions induced by a pair $(\CX,\CD)$, where $\CX$ is a projective model of $X$ over $O_K$, and $\CD$ is an arithmetic $\QQ$-divisor on $\CX$ such that the generic fiber $\CD_K=0$ on $\CX_K$ (instead of just on $U$);
\item[(4)] $G(U^\an)_\rmod$ denotes the space of model Green functions on $U^\an$ associated to Cartier divisors of $U$.
\item[(5)] $C(U^\an)_\cptf$ denotes the completion of $C(U^\an)_\rmod$ for the boundary topology;
\item[(6)] $G(U^\an)_\cptf$ denotes the completion of $G(U^\an)_\rmod$ for the boundary topology.
\end{enumerate}
As $C(U^\an)$ and $G(U^\an)$ are complete, we have inclusions 
$$
C(U^\an)_\rmod \lra C(U^\an)_\cptf \lra C(U^\an),
$$
$$
G(U^\an)_\rmod \lra G(U^\an)_\cptf \lra G(U^\an).
$$

By the direct limit defining $\wh\Div(U^\an)_\rmod$ (commuting with exact sequences), we have a canonical exact sequence
$$
0\lra C(U^\an)_\rmod \lra \wh\Div(U^\an)_\rmod \lra \wt\Div(U/K)_\rmod \lra 0.
$$
Here 
$$
\wt\Div(U/K)_\rmod=\wh\Div(U/K)_\rmod=\varinjlim_{\CX} \Div(X,U),
$$
where the limit is over all projective models $X$ of $U$ over $K$.
We use $\wt\Div$ instead of $\wh\Div$ to avoid confusion. 
We further have a canonical injection
$$
\wh\Div(U^\an)_\rmod \lra \wt\Div(U/K)_\rmod \oplus G(U^\an)_\rmod.
$$

Taking completions, we have a canonical injection 
$$
\wh\Div(U^\an)_\cptf \lra \wt\Div(U/K) \oplus G(U^\an)_\cptf,
$$
and a sequence
$$
0\lra C(U^\an)_\cptf \lra \wh\Div(U^\an)_\cptf \lra \wt\Div(U/K) \lra 0.
$$
Our main result below claims that the sequence is exact and gives an explicit description of $C(U^\an)_\cptf$. 

Recall that the boundary divisor $(0,g_0)$ on $U^\an$ induced by the boundary divisor $(\CX_0,\CE_0)$ on $U$. 
We further denote $X_0=\CX_{0,K}$ and $E_0=\CE_{0,K}$. 

\begin{thm} \label{singularity}
Let $K$ be a field complete for a non-trivial valuation. 
Let $U$ be a quasi-projective variety over $K$. The following is true:
\begin{enumerate}
\item
The canonical sequence 
$$
0\lra C(U^\an)_\cptf \lra \wh\Div(U^\an)_\cptf \lra \wt\Div(U/K) \lra 0
$$
is exact.
\item
For any projective model $X$ of $U$ over $K$, 
denote by $C(U^\an,X^\an)$ the space of continuous functions $h:X^\an\to\RR$
supported on $U^\an$; i.e. $h(x)=0$ for any $x\in X^\an\setminus U^\an$. 
Denote by $C(U^\an)_0$ the image of the injection
$$C(U^\an,X^\an)\lra C(U^\an), \quad h\longmapsto h|_{U^\an}.$$
Then $C(U^\an)_0$ is independent of the choice of $X$ as a projective model of $U$ over $K$. Moreover,
$$
C(U^\an)_\cptf=g_0\cdot C(U^\an)_0
=\{g_0 h: h\in C(U^\an)_0 \}.
$$
\end{enumerate}
\end{thm}

\subsection{Proof of Theorem \ref{singularity}}

The proof of Theorem \ref{singularity} is long as it contains many different parts. We include a detailed proof in the following. 

We first prove that the space $C(U^\an)_0$ in Theorem \ref{singularity}(2) is independent of the choice of 
$X$. 
Namely, if $X'$ is another projective model of $U$ over $K$, then 
$$
\Im(C(U^\an,X^\an)\to C(U^\an))
=\Im(C(U^\an,X'^\an)\to C(U^\an)).$$
We can assume that there is a birational morphism $\pi: X'\to X$ extending the identity map of $U$.
By pull-back via $\pi^\an:X'^\an\to X^\an$, we have an inclusion 
$$
\Im(C(U^\an,X^\an)\to C(U^\an))
\subset \Im(C(U^\an,X'^\an)\to C(U^\an)).$$
It suffices to prove the inverse direction. 
For any $h'\in C(U^\an, X'^\an)$, we want to descend it to the left-hand side. 
Define $h:X^\an\to\RR$ by setting $h|_{U^\an}=h'|_{U^\an}$ and 
$h|_{X^\an \setminus U^\an}=0$. 
Then the pull-back of $h$ via $\pi^\an:X'^\an\to X^\an$ is exactly $h'$. 
It suffices to prove that $h$ is continuous on $X^\an$.
Since $X'^\an$ and $X^\an$ are both Hausdorff and compact, their closed subsets are the same as compact subsets, so $\pi^\an$ is a closed map. Then the continuity of $h'$ implies that of $h$ by the basic result listed in Lemma \ref{closed map}.
This proves the independence on $X$.

Now we prove the second statement of Theorem \ref{singularity}(2), i.e. 
$C(U^\an)_\cptf=g_0\cdot C(U^\an)_0$. 
For any projective model $X$ of $U$ over $K$ dominating $X_0$, the image of 
the natural map $C(X^\an) \to C(U^\an)$ is contained in $g_0\cdot C(U^\an)_0$.
As a consequence, we have a composition of injections 
$$C(U^\an)_\rmod \lra \varinjlim_{X}C(X^\an) \lra g_0\cdot C(U^\an)_0.$$
Here the limit is over all projective models $X$ of $U$ over $K$.
To prove the result, it suffices to prove that $g_0\cdot C(U^\an)_0$ is complete and that $C(U^\an)_\rmod$ is dense in $g_0\cdot C(U^\an)_0$ under the boundary topology. 

It is easy to prove that $g_0\cdot C(U^\an)_0$ is complete. 
In fact, under the bijection $C(U^\an)_0\to g_0\cdot C(U^\an)_0$,
the boundary topology on $g_0\cdot C(U^\an)_0$ corresponds to the uniform topology on $C(U^\an)_0$, which also corresponds to the uniform topology on $C(U^\an,X_0^\an)\subset C(X_0^\an)$. 
The last space is complete. 

Now we prove that $C(U^\an)_\rmod$ is dense in $g_0\cdot C(U^\an)_0$. 
Note that $C(X_0^\an)_\rmod$ is dense in $C(X_0^\an)$ under the uniform topology.
This is already used in the proof of Theorem \ref{image}, as a theorem of Gubler (cf. \cite[Thm. 7.12]{Gub3} and \cite[Lem. 3.5]{Yua1}).
As the uniform topology is stronger than the boundary topology, 
we see that $C(X_0^\an)$ lies in the closure of $C(U^\an)_\rmod$ in 
$g_0\cdot C(U^\an)_0$.
Thus, it is reduced to prove that 
$C(X_0^\an)$ is dense in $g_0\cdot C(U^\an,X_0^\an)$ under the boundary topology.

Let $f=g_0 h$ be an element of $g_0\cdot C(U^\an,X_0^\an)$.
Define 
$$f_n:=\min\{g_0,n\}\cdot h.$$
One checks that 
$\min\{g_0,n\} \in C(X_0^\an)$ and thus $f_n\in C(X_0^\an)$. 
Denote 
$$
Z_n=\{x\in X_0^\an: g_0 (x)\geq n\}, \quad
\epsilon_n=\max\{|h(x)|:x\in Z_n\}.
$$
Note that $\{Z_n\}_n$ decreases to $|E_0|^\an$, so $\epsilon_n$ decreases to 0.
Then we have 
$$|f-f_n|=\max\{0,g_0 -n\}\cdot |h| \leq \epsilon_n g_0.$$
Thus $\{f_n\}_n$ converges to $f$. 
This proves Theorem \ref{singularity}(2).

Now we prove Theorem \ref{singularity}(1), i.e. the exactness of 
$$
0\lra C(U^\an)_\cptf \lra \wh\Div(U^\an)_\cptf \lra \wt\Div(U/K) \lra 0.
$$
We first prove the exactness in the middle. 
Let $\OD$ be an element in the kernel of $\wh\Div(U^\an)_\cptf \to \wt\Div(U/K)$. 
So $\OD$ is the limit of a sequence $\OD_i=(D_i,g_i)$ (with $i\geq1$) in $\wh\Div(U^\an)_\rmod$ with $\lim_i D_i=0$ in $\wt\Div(U/K)$.
We need to prove that $g_\infty:=\lim_i g_i$ lies in $C(U^\an)_\cptf$.
Denote $h_i=g_i/g_0$, viewed as a continuous function on $U^\an$.
It suffices to prove that $h_\infty:=\lim_i h_i$ defines an element of $C(U^\an)_0$ naturally.
We will use the following properties.
\begin{enumerate}
\item[(a)]
The sequence $h_i$ converges uniformly to $h_\infty$ on $U^\an$.
\item[(b)]
There is a compact subset $W_i$ of $U^\an$ for each $i\geq1$, such that
$$
\|h_i\|_{U^\an\setminus W_i,\sup}:=\sup\{|h_i(x)|: x\in U^\an\setminus W_i \}
$$
converges to $0$ as $i\to\infty$. 
\end{enumerate}
Property (b) comes from the condition $\lim_i D_i=0$ in $\wt\Div(U/K)$.
In fact, the condition gives $-\epsilon_i E_0\leq D_i \leq \epsilon_i E_0$ with $\epsilon_i\to 0$. 
In terms of Green functions, this implies that
$\epsilon_i g_0\pm g_i$ is bounded below on $U^\an$. 
As $\epsilon_i g_0$ goes to infinity along the boundary of $U^\an$, we see that 
$2\epsilon_i g_0\pm g_i=\epsilon_i g_0+(\epsilon_i g_0\pm g_i) \geq 0$
in a neighborhood of $\CX_{i,K}^\an\setminus U^\an$ in $\CX_{i,K}^\an$, where $(\CX_{i},\CD_i)$ is a projective model of $(X,0)$ over $O_K$ inducing $\OD_i$. 

We claim that (a) and (b) imply that $h_\infty=\lim_i h_i$ lies in $C(U^\an)_0$.
This is basic in topology. 
In fact, the function $h_\infty$ lies in $C(U^\an)$ by (a). 
It suffices to prove that $h_\infty$ converges to $0$ along the boundary $X_0^\an\setminus U^\an$. 
Assume the contrary. Then there is a sequence $\{x_j\}_{j\geq1}$ in $U^\an$ converging to a point $x_\infty\in X_0^\an\setminus U^\an$ such that $|h_\infty(x_j)|>c$ for a constant $c>0$. 
By (a), we can assume that there is $i_0$ such that $|h_i(x_j)|>c/2$ for all $i\geq i_0$ and $j\geq 1$. 
This implies that 
$\|h_i\|_{U^\an\setminus W_i,\sup}\geq c/2$
for all $i\geq i_0$, which contradicts to (b). 
This proves the exactness in the middle.

It remains to prove the right exactness in Theorem \ref{singularity}(1), i.e. the surjectivity of 
$\wh\Div(U^\an)_\cptf \to \wt\Div(U/K)$.
Let $\wt D$ be an element of $\wt\Div(U/K)$, represented by a Cauchy sequence 
$\{D_i\}_i$ in $\wt\Div(U/K)_\rmod$. 
We need to find a preimage of $\wt D$ in $\wh\Div(U^\an)_\cptf$. 
There is a sequence $\epsilon_i$ of positive rational numbers such that. 
$$
-\epsilon_i E_0\leq  D_i-D_{i+1}\leq \epsilon_i E_0.
$$
By the Cauchy property, replacing $\{D_i\}_i$ by a subsequence if necessary, we can assume that $\sum_{i\geq1} \epsilon_i$ converges.
We claim that for any $i\geq 1$, there is a model Green function $g_i$ of $D_i$ on $U^\an$ such that the sequence $\{g_i\}_i$ satisfies
$$
-\epsilon_i g_0\leq  g_i-g_{i+1}\leq \epsilon_i g_0,\quad i\geq1.
$$
If the claim holds, the sequence $\{(D_i, g_i)\}_i$ is a Cauchy sequence, and represents a preimage of $\wt D$ in $\wh\Div(U^\an)_\cptf$. 

It remains to prove the claim. 
We will construct $g_i$ inductively. 
Assume that $g_1,\cdots, g_i$ is constructed, and we need to construct $g_{i+1}$ satisfying the requirement. 
Let $g_{i+1}'$ be a model Green function of $D_{i+1}$ on $U^\an$. 
Assume that $(D_i,g_i)$ and $(D_{i+1},g_{i+1})$ can be realized as a model arithmetic divisors of mixed coefficients on $(X_{i+1},U)$ for some projective model $X_{i+1}$ of $U$ over $K$. 
Set $g_{i+1}=g_{i+1}'-f$ for $f\in C(X_{i+1}^\an)_\rmod$.
It suffices to find $f\in C(X_{i+1}^\an)_\rmod$ satisfying 
$$
-\epsilon_i g_0\leq  g_i-g_{i+1}'+f\leq \epsilon_i g_0.
$$
As before, $C(X_{i+1}^\an)_\rmod$ is dense in $C(X_{i+1}^\an)$ under uniform convergence. 
So it suffices to find $f\in C(X_{i+1}^\an)$
satisfying
$$
-\epsilon_i (g_0-c_0)\leq  g_i-g_{i+1}'+f\leq \epsilon_i (g_0-c_0),
$$
where $c_0>0$ is a constant with $g_0>c_0$ on $U^\an$.
The condition is equivalent to
$$
g_{i+1}'-g_i-\epsilon_i (g_0-c_0)\leq  f\leq g_{i+1}'-g_i+\epsilon_i (g_0-c_0).
$$
This is an inequality of Green functions on $X_{i+1}^\an$ corresponding to the divisor relation
$$
D_{i+1}-D_i-\epsilon_i E_0\leq  0\leq D_{i+1}-D_i+\epsilon_i E_0.
$$
By the first inequality of divisors, $g_{i+1}'-g_i-\epsilon_i (g_0-c_0)$ has a finite upper bound $c$ on $X_{i+1}^\an$. 
Then we can take 
$$f=\min\{c,\ g_{i+1}'-g_i+\epsilon_i (g_0-c_0)\},$$ 
which is continuous by the second inequality of divisors. 
This finishes the proof of Theorem \ref{singularity}.

In the above proof, the following basic result was used. 
We list it separately since it will be used again later. 
\begin{lem} \label{closed map}
Let $\pi:M\to N$ be a surjective, closed, and continuous map of topological spaces.
Let $f:N\to \RR$ be a map, and $\pi^*f=f\circ \pi:M\to \RR$ the pull-back. 
Then $f$ is continuous if and only if $\pi^*f$ is continuous. 
\end{lem}
\begin{proof}
For the ``if'' part, prove that inverse images of closed sets under $f$ are closed.
\end{proof}

\subsection{Global version of Theorem \ref{singularity}}

As a dilation, we introduce a global version of Theorem \ref{singularity}, which gives a quick sufficient condition for a metrized line bundle to arise from an adelic line bundle. 
For simplicity, we state it in terms of adelic divisors, and we only put singularities at archimedean places. Still, we notice that the result holds if we allow similar singularities at infinitely many places of $\QQ$. 

\begin{thm} \label{singularity global}
Let $\CU$ be a quasi-projective arithmetic variety. 
Let $(\CX_0,\CEE_0)$  be a boundary divisor of $\CU/\ZZ$ with $\CEE_0=(\CE_0, g_0)$. 
Let $(\CD, g_\CD)$ be a pair consisting of a divisor $\CD$ on $\CX$ and a continuous function $g_\CD:\CX(\CC)\setminus (|\CD(\CC)|\cup |\CE_0(\CC)|) \to \RR$.
Assume that $g_\CD$ is a Green function of $\CD$ on $\CX(\CC)$ with $o(g_0)$-singularity; i.e. for any Green function $g_\CD':\CX(\CC)\setminus |\CD(\CC)| \to \RR$ of continuous type of $\CD$ on $\CX(\CC)$, the difference $g_\CD-g_\CD'$  extends to a continuous function on $\CX(\CC)\setminus  |\CE_0(\CC)|$ and grows as $o(g_0)$ along $|\CE_0(\CC)|$. 
Then $(\CD, g_\CD)$ extends to a unique adelic divisor $\CDD'\in \wh\Div(\CU/\ZZ)$ in the sense that 
 the underlying divisor of $\CDD'$ is $\CD|_\CU$, and that the Green function of $\CD|_\CU$ on $\CU^\ran$ induced by 
$\CDD'$ (via Proposition \ref{injection3}) 
is the same as 
the Green function of $\CD|_\CU$ on $\CU^\ran$ induced by 
$(\CD, g_\CD)$. 
\end{thm}

\begin{proof}
The uniqueness follows from the injectivity of the analytification map in Proposition \ref{injection3}.
For the existence, consider the pair $(0, f)$ on $\CX_0$, where $0$ is the zero element of $\Div(\CX)$, and $f=g_\CD-g_\CD'$ is continuous on $\CX(\CC)\setminus  |\CE_0(\CC)|$ and grows as $o(g_0)$ along $|\CE_0(\CC)|$. 
It suffices to prove that $(0, f)$ extends to a unique adelic divisor on $\CU/\ZZ$. 
This is a global version of Theorem \ref{singularity}. 
The proof of the theorem holds for the global $(0, f)$ (over $\ZZ$) since the only problem appears at the archimedean places. We omit the details here.
\end{proof}

A major feature of the theorem is that the extra singularity of $g_\CD$ grows as $o(g_0)$ along the boundary, which is a very broad type of singularity. 
For a natural example of a compactified divisor on $\CC$ with a Green function of singularity $O(\log g_0)$, see \S\ref{sec hodge} for the Hodge bundle of the moduli space of principally polarized abelian varieties. 

In the literature, there are many theories and results on Green functions with singularities. 
We first note that Burgos--Kramer--K\"uhn \cite{BKK} has introduced a general arithmetic intersection theory of arithmetic Chow cycles with pre-log-log currents. This theory treats particular singularities of type $O(\log g_0)$ and thus includes the above example of Hodge bundles. 
We also refer to Bost \cite{Bos} and Moriwaki \cite{Mor1} for an arithmetic intersection of arithmetic Chow cycles with $L_1^2$-currents. 
To compare with our theory, all these references treat intersection theory on a fixed projective arithmetic variety and focus on singularities of Green currents. In contrast, we treat the intersection theory of suitable limits of hermitian line bundles (which corresponds to arithmetic Chow cycles of co-dimension one). In the limit process, our underlying line bundles also vary.

\subsection{Compactified metrics}

Here we briefly introduce the corresponding notions of compactified line bundles. 
Resume the above local notations.
Namely, let $K$ be a field complete for non-trivial absolute value 
$|\cdot|$.
Set $O_K$ to be $K$ in the Archimedean case and to be the valuation ring in the non-archimedean case.
Let $U$ be a quasi-projective variety over $K$.

Recall from Proposition \ref{injection local} that there are canonical injective maps 
$$\wh \Picc (U/O_K)\lra \wh\Picc (U^\an),$$
$$\wh\Pic (U/O_K) \lra \wh \Pic (U^\an).$$
Denote 
$$\wh \Picc (U^\an)_\cptf=\Im(\wh \Picc (U/O_K)\to \wh \Picc (U^\an)),$$
$$\wh \Pic (U^\an)_\cptf=\Im(\wh \Picc (U/O_K)\to \wh \Pic (U^\an)).$$
They are compactifications of 
$$\wh \Picc (U^\an)_\rmod=\Im(\wh \Picc (U/O_K)_\rmod\to \wh \Picc (U^\an)),$$
$$\wh \Pic (U^\an)_\rmod=\Im(\wh \Picc (U/O_K)_\rmod\to \wh \Pic (U^\an)).$$
As in the case of arithmetic divisors, we are going to describe these groups or categories directly on $U^\an$.

If $K$ is non-archimedean, an \emph{arithmetic model} (or \emph{integral model}) of a line bundle $L$ on $U$ is a pair $(\CX,\CL)$, where $\CX$ is a projective model of $U$ over $O_K$, and $\CL$ is a $\QQ$-line bundle on $\CX$ extending $L$.

If $K$ is archimedean, an \emph{arithmetic model} of a line bundle $L$ on $U$ is a pair $(\CX,\CL)$, where $\CX$ is a projective model of $U$ over $O_K=K$, and $\CL=(\wt L, \|\cdot\|_{\wt L})$ consisting of a $\QQ$-line bundle $\wt L$ on $\CX$ extending $L$ and a continuous metric $\|\cdot\|_{\wt L}$ of $\wt L$ on $\CX^\an$.

In both cases, the arithmetic model $(\CX,\CL)$ induces a metric 
$\|\cdot\|_\CL$ of $\CL_K$ on $\CX_K^\an$, and thus a metric of $L$ on $U^\an$ by restriction.
The process is essentially the same as the global case described in \S\ref{sec analytification line bundles}. 
The metric $\|\cdot\|_\CL$ of $L$ on $U^\an$ is called a \emph{model metric}, and the metrized line bundle
$(L, \|\cdot\|_\CL)$ on $U^\an$ is called a \emph{model metrized line bundle}.

The model metric or the model metrized line bundle is called \emph{nef} (or \emph{semipositive}) if either $\CL$ is nef on $\CX$ in the non-archimedean case or the metric of $\CL$ is semipositive on $\CX(\CC)$ (cf. \S\ref{sec hermitian}).

Recall that we have boundary topologies on $\wh\Div(X^\an)$ and $C(X^\an)$
in terms of $g_0=g_{\CE_0}|_{U^\an}$ obtained by the choice of a pair $(\CX_0,\CE_0)$.

A metrized line bundle $\OL=(L, \|\cdot\|)$ on $U^\an$ or its metric is called \emph{compactified} if there is a sequence of model metrics 
$\{\|\cdot\|_i\}_{i\geq1}$ of $L$ on $U^\an$, such that the continuous function $\log(\|\cdot\|_i/\|\cdot\|)$ on $U^\an$ converges to $0$ under the boundary topology on $C(X^\an)$. 

The metrized line bundle $\OL$ or its metric $\|\cdot\|$
is said to be \emph{strongly nef} (or \emph{strongly semipositive}) if there exists such a sequence such that every model metric $\|\cdot\|_i$ is nef. 
The metrized line bundle $\OL$ or its metric $\|\cdot\|$
is said to be \emph{nef} (or \emph{semipositive}) if there exists a strongly nef metrized line bundle $\OM$  such that $a\OL+\OM$ is strongly nef for all positive integers $a$.
The metrized line bundle $\OL$ or its metric $\|\cdot\|$
is said to be \emph{integrable} if $\OL$ is isometric to the difference of two strongly nef metrized line bundles. 

Finally, our result is as follows:
\begin{enumerate}
\item $\wh\Pic(U^\an)_\rmod$ (resp. $\wh\Pic(U^\an)_\cptf$) is the subgroup of 
$\wh\Pic(U^\an)$ consisting of model metrized (resp. compactified)  line bundles on $U^\an$.
\item
 $\wh\Picc(U^\an)_\rmod$ (resp. $\wh\Picc(U^\an)_\cptf$) is equivalent to the full subcategory of $\wh\Picc(U^\an)$ consisting of model (resp. compactified) metrized line bundles on $U^\an$.
\end{enumerate}

\subsection{Chambert-Loir measures} \label{sec measure}

Let $U$ be a quasi-projective variety over a complete field $K$ with a non-trivial valuation as above. Denote $n=\dim U$.
Let $\OL_1,\OL_2,\cdots, \OL_n$ be strongly nef (compactified) metrized line bundles on $U^\an$.
We will see that there is a canonical Radon measure  
$c_1(\OL_1)c_1(\OL_2)\cdots c_1(\OL_n)$ on the Berkovich space $U^\an$, which generalizes the Monge--Amp\`ere measure in the complex case. 
We will call this measure \emph{the Chambert-Loir measure}.

If $K$ is archimedean, this is treated by classical analysis in \cite[Thm. 2.1]{BT} (or \cite[Cor. 1.6]{Dem1}).  
If $K$ is non-archimedean and $U$ is projective, 
this is constructed by Chambert-Loir \cite{CL} when $K$ has a dense and countable subfield and extended to general $K$ by Gubler \cite{Gub1}. 
If $K$ is non-archimedean and $U$ is quasi-projective, we will follow the theory of 
Chambert-Loir and Ducros in \cite{CLD}, a vast generalization of the construction of \cite{CL} via a local analytic approach. 

In the following, assume that $K$ is non-archimedean and that $U$ is quasi-projective over $K$.
We are going to apply \cite[Cor. 5.6.5]{CLD} to $\OL_1,\OL_2,\cdots, \OL_n$. 
Before that, we claim that the metric of a strongly nef metrized line bundle 
 satisfies the condition of
\cite[Cor. 5.6.5]{CLD}; i.e. it is locally psh-approachable on $U^\an$
in the sense of \cite[6.3.1, Def. 5.6.3, Def. 5.5.1]{CLD}.

In fact, let $\OL=(L,\|\cdot\|)$ be a strongly nef metrized line bundle on $U^\an$.
By definition, the metric $\|\cdot\|$ is the limit of a sequence of nef model metrics 
$\|\cdot\|_i$ under the boundary topology of $C(U^\an)$.
Since $U^\an$ is locally compact, the convergence is locally uniform as in Lemma \ref{local subspace2}. 
Therefore, it suffices to prove the nef model case.
So, we assume that the metric $\|\cdot\|$ is a nef model metric, and we need to prove that it is locally psh-approachable. 
This is a consequence of \cite[Cor. 6.3.4]{CLD}, since the metric $\|\cdot\|$ 
is induced by a nef line bundle $\CL$ on a projective model $\CX$ of $U$ over $O_K$.
Note that the loc. cit. is only stated for the ample case but can be extended to the nef case. In fact, take any ample line bundle $\CM$ on $\CX$, which induces a metric $\|\cdot\|_\CM$ of $M=\CM|_U$ on $U^\an$. 
For any local sections $s$ and $t$ of $L$ and $M$ regular and everywhere 
non-vanishing on a Zariski open set $W$ of $U$, the function 
$-\log \|s\|-\epsilon\log \|t\|_\CM$ is globally psh-approachable on $W^\an$ for any positive rational numbers $\epsilon>0$. 
As $\epsilon\to 0$, the function converges to $-\log \|s\|$, which is uniform on any compact subset of $W^\an$.
This proves that $\|\cdot\|$ is locally psh-approachable and finishes the quasi-projective case.

Finally, by \cite[Cor. 5.6.5]{CLD}, there is a canonical  measure 
$$
c_1(\OL_1)\cdots c_1(\OL_n)
=d'd''(-\log \|\cdot\|_{1}) \wedge \cdots \wedge d'd''(-\log \|\cdot\|_{n})
$$
over $U^\an$.
Here the right-hand side is understood as follows. 
If $W$ is a Zariski open subset of $X$ and 
$s$ is a regular and everywhere non-vanishing section of $L$ on $W$, 
then we set $d'd'' (-\log \|\cdot\|_{i})=d'd''(- \log \|s\|_{i})$ on $W^\an$. 
This is independent of the choice of $s$ by the Poincar\'e-Lelong formula in
\cite[Thm. 4.6.5]{CLD}.

The measure is defined by a weak convergence process. 
We describe it as follows. 
For any $i=1,\cdots, n$, the metric of $\OL_i$ is the limit of model metrics induced by (projective) arithmetic models $(\CX_{i,j}, \CL_{i,j})$ of $(U,L)$ over $O_K$.
We can assume that $\CX_{i,j}$ is independent of $i$ and write it as $\CX_j$.
Denote $X_j=\CX_{j,K}$, which is a projective model of $U$ over $K$. 
Denote by $\OL_{i,j}=(L_{i,j}, \|\cdot\|_{i,j})$ the metrized line bundle on the compact space $X_j^\an$, induced by the model $(\CX_{i,j}, \CL_{i,j})$. 
Denote by $C_c(U^\an)$ the space of real-valued, continuous, and compactly supported function on $U^\an$.
Then the construction gives, for any $f\in C_c(U^\an)$, 
$$
\int_{U^\an}f c_1(\OL_1)\cdots c_1(\OL_n)
= \lim_{j\to \infty} \int_{X_j^\an} f c_1(\OL_{1,j})\cdots c_1(\OL_{n,j}).
$$
As $X_j$ is projective over $K$, the right-hand side is equal to the integration defined by global intersection numbers by \cite{CL, Gub1}. 

It is worth noting that by \cite[Cor. 4.2]{CT}, the integral of $c_1(\OL_{1,j})\cdots c_1(\OL_{n,j})$ on any Zariski closed subset of $X_j^\an$ of positive codimension is 0.

\subsection{Application to finitely generated fields} \label{sec app global}

\kkk
\ccc
Let $F$ be a finitely generated field over $k$.
Let $v$ be a point of $(\Spec F)^\an=\CM(F/k)$ that is \emph{not} the trivial valuation over $F$. 
Denote by $F_v$ the completion of $F$ for $v$. 
It can be either archimedean or non-archimedean.

Let $X$ be a \emph{quasi-projective} variety of dimension $n$ over $F$.
Let $\overline L$ be a \emph{strongly nef} adelic line bundle on $X$ with an underlying line bundle $L$ on $X$. 
Let $X_v^\an$ be the Berkovich space associated to the variety $X_{F_v}$ over the complete field $F_v$, which is the fiber of $X^\an\to (\Spec F)^\an$ above $v$. 
By Proposition \ref{injection2}, $\overline L$ induces a metric $\|\cdot\|$ of $L$ on $X^\an$, which restricts to a $F_v$-metric $\|\cdot\|_v$ of $L$ on $X_v^\an$.  

\begin{lem} \label{psh}
Assume that $X$ is quasi-projective over $F$ and that $\OL$ is strongly nef on $X$. 
Then the metric $\|\cdot\|_v$ of $L$ on $X_v^\an$ induced by $\OL$ is strongly nef.
\end{lem}

With the lemma, there is a Chambert-Loir measure 
$$
c_1(\OL)_v^n:=c_1(L_{F_v}, \|\cdot\|_v)^n
$$
over the Berkovich space $X_v^\an$ for any point $v\in \CM(F/k)$ which is non-trivial over $F$.
This measure will be used in our equidistribution conjectures and theorems. 

By multi-linearity, for any integrable line bundles $\overline L_1,\cdots, \OL_n$ on $X$, there is a (signed)
Chambert-Loir measure 
$$
c_1(\OL_1)_v\cdots c_1(\OL_n)_v:=c_1(L_{1, F_v}, \|\cdot\|_v)\cdots c_1(L_{n, F_v}, \|\cdot\|_v)
$$
over the Berkovich space $X_v^\an$ for any point $v\in \CM(F/k)$ which is non-trivial over $F$.
Now we prove the lemma. 

\begin{proof}[Proof of Lemma \ref{psh}]
We only treat the case that $v$ is non-archimedean since the archimedean case is 
easier. 

Assume that $\OL$ is represented by a Cauchy sequence 
$\CLL=(\CL,(\CX_i,\CLL_{i}, \ell_{i})_{i\geq1})$ in $\wh\Picc(\CU/k)_\rmod$. Here $\CU$ is a quasi-projective model of $X$, and each $\CLL_{i}$ is nef on $\CX_i$.
By Lemma \ref{models}, we can assume that $\CU$ is equipped with a flat morphism $\CU\to \CV$ to a quasi-projective variety $\CV$, whose generic fiber is isomorphic to $X\to \Spec F$.
Let $\CS$ be a fixed projective model of $\CV$.
By blowing-up $\CX_i$ if necessary, we can assume that $\CU\to \CV$ extends to a morphism $\CX_i\to \CS$.

The point $v\in (\Spec F)^\an \subset \CS^\an$ has a residue field $F_v$ and a valuation ring $R_v\subset F_v$. 
By the valuative criterion, the morphism $\Spec F_v \to \CS$ extends to a morphism
$\Spec R_v \to \CS$. 
The base change of $\CX_{i}\to\CS$ gives a morphism $\CX_{i,R_v}\to \Spec R_v$ whose generic fiber contains $X_{F_v}$ as an open subvariety.
Denote by $\CX_{i,R_v}'$ the Zariski closure of $X_{F_v}$ in $\CX_{i,R_v}$, 
so that $\CX_{i,R_v}'$ is the unique irreducible component of $\CX_{i,R_v}$ flat over $R_v$. 
By pull-back, we get a sequence of $\QQ$-line bundles 
$\CL_i|_{\CX_{i,R_v}'}$ on $\CX_{i,R_v}'$, which induces a sequence of model metrics 
$\|\cdot\|_i$ of $L$ on $X_v^\an$.
The limit of these metrics is exactly the desired metric $\|\cdot\|_v$. 
Moreover, the convergence of $\{\|\cdot\|_i\}_i$ to $\|\cdot\|_v$ is for the boundary topology, as we can see in the proof of Proposition \ref{injection2} for quasi-projective varieties.
\end{proof}

\begin{rmk}
By the lemma, if $X$ is projective over $F$, then $\|\cdot\|_v$ is semipositive in the sense that it is a uniform limit of metrics induced by nef models. In this case, we can also use the construction of \cite{CL, Gub1} to define the measure. 
\end{rmk}

\chapter{Intersection theory}

In this chapter, we develop an intersection theory of integrable adelic line bundles.
There are two types of intersection pairings. The first type gives an absolute intersection number; the second type is an intersection pairing in a relative setting in terms of the Deligne pairing. 
While the absolute intersection number is easy to obtain, the construction of the relative intersection pairing takes most of this chapter.

\section{Intersection theory}

In this section, we state both intersection pairings, prove the existence of the absolute version, and leave the proof of the relative version to the rest of this section.

\subsection{Absolute intersection numbers}

In algebraic geometry, for a projective variety $\CX$ of dimension $d$ over a base field, 
there is an intersection pairing $\Pic (\CX)^{d}\to \ZZ$.

In Arakelov geometry, there is an intersection of hermitian line bundles by Deligne \cite{Del} and Gillet--Soul\'e \cite{GS1}.  
Namely, for each projective variety $\CX$ of absolute dimension $d$ over $\ZZ$, 
there is an intersection pairing $\wh\Pic (\CX)_\sm^{d}\to \BR$, which was extended to a pairing $\wh\Pic (\CX)_\intb^{d}\to \BR$ as recalled in \S\ref{sec hermitian}.
See also \S\ref{app subsec intersection} for the smooth case.

We are going to extend these pairings to adelic line bundles.
As in the case of \cite{Zha2,Mor4}, we cannot expect the intersection to be defined for all adelic line bundles, but we require all but one adelic line bundles to be integrable.

\begin{prop} \label{intersection1}
\kkk
For any flat and essentially quasi-projective integral scheme $X$ over $k$, the intersection pairing above extends to a canonical multi-linear homomorphism
$$
\wh \Pic (X/k)\times \wh \Pic (X/k)_\intb^{d-1} \longrightarrow \BR.
$$
Here $d$ is the absolute dimension of a quasi-projective model of $X$ over $k$.
The homomorphism is symmetric in the last $d-1$ variable;
if the first variable is also in $\wh \Pic (X/k)_\intb$, then the homomorphism is symmetric in all $d$ variables. 

Moreover, if $\OL_1,\cdots, \OL_d$ are nef adelic line bundles on $X$, then their intersection number $\OL_1\cdot \OL_2 \cdots \OL_d\geq 0$.
\end{prop}

\begin{proof}
We only define the intersection number, and omit proofs of the other properties. 
It suffices to treat the case that $X=\CU$ is a quasi-projective variety over $k$. 
We first define the intersection pairing 
$\wh \Pic (\CU/k)_\intb^{\, d} \to \BR.$
By linearity, it suffices to define
$\overline\CD_1\cdot \overline\CD_2 \cdots \overline\CD_{d}$ for any
$\overline\CD_1,\cdots, \overline\CD_{d}\in\wh \Div (\CU/k)_{\QQ,\snef}$.
Here we have switched to adelic $\QQ$-divisors for simplicity of notations.

Let $(\CX_0,\OCE_0)$ be a boundary divisor of $\CU$ over $k$, and we will use it to define the boundary topology of $\wh \Div (\CU/k)_{\rmod,\QQ}$.
We can further assume that $\OCE_0$ is nef, which is possible by simply replacing 
$\OCE_0$ by a nef arithmetic divisor $\OCE_0'$ on $\CX_0$ satisfying $\CE_0'\geq \CE_0$ and replacing $\CU$ by $\CX_0\setminus |\CE_0'|$.

For $j=1,\cdots,d$, assume that $\overline\CD_j$ is represented by a Cauchy sequence
$\{(\CX_i,\overline\CD_{j,i})\}_{i\geq1}$, where each $\overline\CD_{j,i}$ is a nef arithmetic $\QQ$-divisor on the projective model $\CX_i$ dominating $\CX_0$.
Here we assume that the model $\CX_i$ is independent of $j$, which is always possible.
There is a sequence $\{\epsilon_i\}_{i\geq1}$ of positive rational numbers converging to 0 such that 
$$-\epsilon_i\OCE_0\leq 
\CDD_{j,i'}-\CDD_{j,i}
\leq \epsilon_i\OCE_0, \qquad i'>i$$
for any $j=1,\cdots, d$.

For any subset $J\subset \{1, \cdots, d\}$, consider the intersection number
$$\alpha_{J, i}:=\OCE_0 ^{d-|J|}\prod_{j\in J}\overline \CD_{j, i}.$$
We will prove by induction that $\{\alpha_{J, i}\}_{i\geq1}$ is a Cauchy sequence and thus convergent in $\RR$.
When $J$ is the full set, the limit of the Cauchy sequence gives our definition of 
$\overline\CD_1\cdot \overline\CD_2 \cdots \overline\CD_{d}$.

There is nothing to prove if  $J$ is the empty set.
Assume the claim is true for any $|J|<r$ for some $r>0$.
We need to prove the result for any $J$ with $|J|=r$. 
Without loss of generality, assume $J=\{1,2,\cdots, r\}$. 
Then
\begin{align*}
\alpha_{J, i'}-\alpha_{J, i}
=&\ \OCE_0^{d-r}\overline\CD_{1,i'}\cdots \overline\CD_{r,i'}-\OCE_0^{d-r} \overline\CD_{1,i}\cdots \overline\CD_{r,i}\\
\leq &\ \OCE_0^{d-r}(\overline\CD_{1,i}+\epsilon_i \OCE_0)\cdots (\overline\CD_{r,i}+\epsilon_i \OCE_0)-\OCE_0^{d-r} \overline\CD_{1,i}\cdots \overline\CD_{r,i}\\
= &\ \sum_{J'\subsetneq J} \epsilon_i^{r-|J'|}\alpha_{J', i}.
\end{align*}
Similarly,  
\begin{align*}
\alpha_{J, i}-\alpha_{J, i'}
\leq &\ \OCE_0^{d-r}(\overline\CD_{1,i'}+\epsilon_i \OCE_0)\cdots (\overline\CD_{r,i'}+\epsilon_i \OCE_0)-\OCE_0^{d-r} \overline\CD_{1,i'}\cdots \overline\CD_{r,i'}\\
= &\ \sum_{J'\subsetneq J} \epsilon_i^{r-|J'|}\alpha_{J', i'}.
\end{align*}
It follows that $\{\alpha_{J, i}\}_i$ is a Cauchy sequence.

Therefore, we have defined the intersection number
$\overline\CD_1\cdot \overline\CD_2 \cdots \overline\CD_{d}$ for any
$\overline\CD_1,\cdots, \overline\CD_{d}\in\wh \Div (\CU/k)_{\QQ,\snef}$.
It is independent of the choice of the Cauchy sequence $\{(\CX_i,\overline\CD_{j,i})\}_{i\geq1}$ for each $j$, since we can merge any two different Cauchy sequences into a single one. 

Now we extend the definition of the intersection number to $\CDD_1 \in\wh \Div (\CU/k)_\QQ$ and $\overline\CD_2,\cdots, \overline\CD_{d}\in\wh \Div (\CU/k)_{\QQ,\snef}$.
We need to further approximate $\CDD_1$. 
Take the above notation for the Cauchy sequence
$\{(\CX_i,\overline\CD_{1,i})\}_{i\geq1}$ and the relation 
$$-\epsilon_i\OCE_0\leq 
\CDD_{1,i'}-\CDD_{1,i}
\leq \epsilon_i\OCE_0, \qquad i'>i.$$
Note that $\CDD_{1,i}$ is an arithmetic $\QQ$-divisor on $\CX_i$, which is not assumed to be nef any more, but we can assume that $\CDD_{1,i}$ has a Green function of smooth type. It follows that
the intersection number $\beta_i=\overline\CD_{1,i}\cdot \overline\CD_2 \cdots \overline\CD_{d}$ is already defined. 
It remains to prove that $\{\beta_i\}_{i\geq 1}$ is a Cauchy sequence. 
In fact, we simply have
$$
\beta_{i}-\beta_{i'}=(\overline\CD_{1,i}-\overline\CD_{1,i'})\cdot \overline\CD_2 \cdots \overline\CD_{d}
\leq \epsilon_i\,\OCE_0\cdot \overline\CD_2 \cdots \overline\CD_{d}
$$
and 
$$
\beta_{i}-\beta_{i'}=(\overline\CD_{1,i}-\overline\CD_{1,i'})\cdot \overline\CD_2 \cdots \overline\CD_{d}
\geq -\epsilon_i\,\OCE_0\cdot \overline\CD_2 \cdots \overline\CD_{d}.
$$
Here we have used the fact that the intersection number of an effective arithmetic $\QQ$-divisor with 
$\overline\CD_2 \cdots \overline\CD_{d}$ is non-negative.
This finishes the proof.
\end{proof}

A basic property of the intersection number is the following projection formula. 

\begin{prop} [projection formula]
\kkk
Let $f:X'\to X$ be a morphism of flat and essentially quasi-projective integral schemes over $k$.
Assume that the absolute dimensions of quasi-projective models of $X'$ and $X$ over $k$ are all equal to $d$.
Let $\OL_1,\cdots, \OL_d$ be integrable adelic line bundles on $X$. 
Then  
$$
f^*\OL_1\cdot f^*\OL_2 \cdots f^*\OL_d=\deg(f) \, (\OL_1\cdot \OL_2 \cdots \OL_d).
$$
Here if $f$ is dominant in that it maps the generic point of $X'$ to the generic point of $X$, then $\deg(f)$ is the degree of the extension between the function fields; otherwise, we take the convention $\deg(f)=0$.
\end{prop}

\begin{proof}
By the limit process, it is reduced to the well-known formula in the projective case.
\end{proof}

\subsection{Deligne Pairing: main theorem}

Let $f: X\to Y$ be a projective and flat morphism of noetherian schemes of pure relative dimension $n$. The Deligne pairing is a multi-linear functor
$$\Picc (X)^{n+1}\lra \Picc (Y), \quad
(L_1,\cdots,L_{n+1})
\longmapsto \pair{L_1, \cdots, L_{n+1}}.
$$
The functor refines the intersection of the Chern classes of the line bundles. It satisfies many natural functorial properties, including the base change property, the multi-linearity, the symmetry, and the induction formula.

For a brief history of pairing, the case $n=0$ is just the norm functor $N_{X/Y}$. 
Deligne \cite{Del} constructed the functor for $n=1$ and speculated a similar pairing for general $n$. 
Deligne's major motivation is to formulate an arithmetic Riemann--Roch theorem for families of curves. We will not need this formulation here, but we refer interested readers to \S\ref{app sec surface} for a sketch on it.  
For general $n$, the pairing was constructed by Elkik \cite{Elk1} for any $f$ which is projective, flat, and further Cohen--Macaulay, and by Munoz Garcia \cite{MG} for any $f$ which is projective, equi-dimensional and of finite Tor-dimension (which implies the projective and flat case). Moreover, Ducrot \cite{Duc} had a different treatment of the projective and flat case.

If $X$ and $Y$ are smooth varieties over $\CC$ and $f$ is smooth, and if $L_1,\cdots,L_{n+1}$ are endowed with smooth hermitian metrics, then the metrics transfer to a canonical smooth hermitian metric on 
$\pair{L_1, \cdots, L_{n+1}}$, as constructed by Deligne \cite{Del} and Elkik \cite{Elk2}. 
As we will prove later, the metric construction can be generalized to the projective and flat case, and in this case, the Deligne pairing transfers continuous metrics to continuous metrics. 

Our goal is to extend the Deligne pairing to adelic line bundles. This section's main result is as follows.

\begin{thm} \label{intersection2}
\kkk
Let $Y$ be a flat and essentially quasi-projective integral scheme over $k$.
Let $f: X\to Y$ be a projective and flat morphism of relative dimension $n$.
Assume that $X$ is integral and $Y$ is normal.
Then the Deligne pairing induces a symmetric and multilinear functor
$$\wh \Picc(X/k)_\intb^{\, n+1} \longrightarrow \wh \Picc (Y/k)_\intb.
$$
When restricted to strongly nef or nef adelic line bundles, the functor induces functors
$$\wh \Picc(X/k)_\snef^{\, n+1} \longrightarrow \wh \Picc (Y/k)_\snef,$$
$$\wh \Picc(X/k)_\nef^{\, n+1} \longrightarrow \wh \Picc (Y/k)_\nef.$$
Moreover, the maps are compatible with base changes of the form $Y'\to Y$, where $Y'$ is any normal integral scheme, flat and essentially quasi-projective over $k$, such that $X'=X\times_YY'$ is integral. 
\end{thm}

The proof of this theorem will take up the rest of this chapter. 
After some preparations about metrics of the Deligne pairings and basic properties in the model case, the proof of the theorem will be given in \S \ref{sec pairing proof}.

\section{Metrics of the Deligne pairing: statements}

The goal of this section is two-fold. 
First, we review the treatment of the Deligne pairing of \cite{MG} to set up a framework for our treatment. 
Second, we state some results on natural metrics of the Deligne pairing from metrics of the original line bundles, which generalizes the result of \cite{Del, Elk2} from the smooth case to the general case.  
Note that the treatments of \cite{Zha3,Mor2} on the metrics have gaps due to misinterpretations of the definition of the canonical section $\pair{s_1, \cdots, s_{n+1}}$
of $\pair{L_1, \cdots, L_{n+1}}$ on $Y$.

\subsection{Deligne Pairing: review}

Here we recall some results of the Deligne pairing in \cite{MG}. Our main interest is the Deligne pairing for projective and flat morphisms. Still, it seems inevitable to treat non-flat morphisms of finite Tor-dimension if we want to pass to generic hyperplane sections by an induction formula. 
Therefore, we will follow the generality of \cite{MG} to treat morphisms of finite Tor-dimension.

Recall that a morphism $f:X\to Y$ of noetherian schemes is \emph{of pure relative dimension $n$} if for every $y\in Y$, every irreducible component of $X_y$ (if non-empty) has dimension $n$.

Recall that a morphism $f:X\to Y$ of noetherian schemes is \emph{of finite Tor-dimension} if one of the following two equivalent conditions holds:
\begin{itemize}
\item[(a)] there is an integer $d_0$ such that $\mathrm{Tor}_d^{B}(A,M)=0$ for any $d>d_0$, for any affine open subscheme $\Spec A$ of $X$ whose image under $f$ lies in an affine open subscheme $\Spec B$ of $Y$, and for any $B$-module $M$.
\item[(b)] there is an integer $d_0$ such that $\mathrm{Tor}_d^{\CO_{Y,y}}(\CO_{X,x},M)=0$ for any $d>d_0$,
 for any point $x\in X$ with $y=f(x)\in Y$, and for any $\CO_{Y,y}$-module $M$.
\end{itemize}
See \cite[III, \S3, Def. 3.2, Prop. 3.3]{SGA6} for more information. 
Note that this holds automatically if $f$ is flat or $Y$ is regular. 
Moreover, if $f:X\to Y$ is {of finite Tor-dimension}, and $Z$ is an effective Cartier divisor of $X$, then $Z\to Y$ is also of finite Tor-dimension. 

Let $f:X\to Y$ be a projective morphism of noetherian schemes of finite Tor-dimension and pure relative dimension $n$. 
Let $s_1,\cdots, s_{n+1}$ be global sections of $L_1,\cdots, L_{n+1}$ on $X$ respectively.
For any $i=1,\cdots, n+1$, denote by $Z_i=\div(s_1)\cap \cdots \cap \div(s_i)$ the schematic intersection in $X$. Set $Z_0=X$ for convenience.
Following \cite[Def. 4.3.2]{MG}, we say that the sequence $(s_1,\cdots, s_{n+1})$ is \emph{strongly regular} if the following conditions hold:
\begin{itemize}
\item[(1)] for any $i=1,\cdots, n+1$, the section $s_i$ is not a zero-divisor on $Z_{i-1}$ in the sense that the morphism $\CO_{Z_{i-1}}\to L_i|_{Z_{i-1}}$ induced by $s_i$ is injective;
\item[(2)] for any $i=1,\cdots, n$, the scheme $Z_i$ is purely of relative dimension $n-i$ over $Y$.
\end{itemize}
If (1) and (2) hold for $i=1,\cdots, n$, then we say that the sequence $(s_1,\cdots, s_{n})$ is \emph{strongly regular}.
The condition (1) is symmetric in $s_1,\cdots, s_{n+1}$ by a basic property of regular sequences in local rings. 
Note that the notion of \emph{strongly regular} is stronger than the notion of \emph{very regular} in \cite[Def. 3.2.1, Def. 3.2.4]{MG}, and is more convenient in applications.

The following existence of strongly regular sequence will be frequently used in our treatment. 
If $L_1,\cdots, L_{n+1}$ are $f$-ample on $X$, then there is a finite Zariski open cover $V$ of $Y$ and a positive integer $m$, such that the base change 
$(L_1^{\otimes m})_V,\cdots, (L_{n+1}^{\otimes m})_V$ has a strongly regular sequence of sections for the morphism $f_V:X_V\to V$. 
In fact, for any closed point $y\in Y$, let $V_y$ be an affine open neighborhood of $y$ in $Y$.
Then we can find a global section $s_1$ of
$(L_1^{\otimes m})_{V_y}$ on $X_{V_y}$ for some positive integer $m$ such that
$s_1$ is non-vanishing at any associated point of $X_{V_y}$ or $X_y$. 
This can be guaranteed by requiring $s_1$ to be non-vanishing at a prescribed closed point of every irreducible component and every embedded component of $X_{V_y}$ and $X_y$. 
Then $s_1$ is not a zero-divisor on $X_{V_y}$ or $X_y$, and thus $\div(s_1)\cap X_y$ is pure of dimension $n-1$. 
By semi-continuity of dimensions of fibers (cf. \cite[IV-3, Cor, 13.1.5]{EGA}),  $\div(s_1)$ is of pure relative dimension $n-1$ over a neighborhood of $y$ in $V_y$. 
Replace $V_y$ with this open neighborhood. 
By induction, this gives a strongly regular sequence.

Let $(s_1,\cdots, s_{n+1})$ be a strongly regular sequence of sections of $(L_1,\cdots, L_{n+1})$ on $X$. 
By \cite[Prop. 3.2.6]{MG},  
there is a canonical global section $\pair{s_1, \cdots, s_{n+1}}$
of $\pair{L_1, \cdots, L_{n+1}}$ on $Y$. 
There is a canonical isomorphism 
$$
r:\pair{L_1, \cdots, L_{n+1}} \lra N_{Z_n/Y}(L_{n+1}).
$$
The global section $s_{n+1}$ gives a global section $N_{Z_n/Y}(s_{n+1})$ of $N_{Z_n/Y}(L_{n+1})$. 
Set 
$$\pair{s_1, \cdots, s_{n+1}}=r^{-1}(N_{Z_n/Y}(s_{n+1})).$$
Note that $Z_n\to Y$ is finite but not necessarily flat over $Y$, the norm functor 
$$N_{Z_n/Y}:\Picc(Z_n)\lra \Picc(Y)$$ 
is defined in \cite[\S1.2]{MG} as a natural generalization of the finite and flat case.
The section $\pair{s_1, \cdots, s_{n+1}}$
of $\pair{L_1, \cdots, L_{n+1}}$ behaves well if switching the orders of $L_1, \cdots, L_{n+1}$; see \cite[Thm. 3.4.2]{MG}. 

This essentially gives construction of $\pair{L_1, \cdots, L_{n+1}}$ for relatively ample line bundles $L_1, \cdots, L_{n+1}$ on $X$. By linearity, it generalizes to arbitrary line bundles $L_1, \cdots, L_{n+1}$ on $X$. 

As a convention, the Deligne pairing
$\pair{L_1, \cdots, L_{n+1}}$
for the morphism $f:X\to Y$ will also be written as 
$$
f_*\pair{L_1, \cdots, L_{n+1}},\quad
\pair{L_1, \cdots, L_{n+1}}_{X/Y},\quad 
\pair{L_1, \cdots, L_{n+1}}_{X}.
$$
This may be used when we vary $f: X\to Y$ to avoid confusion.
We take this convention for all similar pairings introduced later.

\subsection{Deligne Pairing: metric at a point} \label{sec metric at point}

Let $Y=\Spec \CC$ and $f:X\to Y$ be a projective morphism of pure relative dimension $n$.
Let $\OL_1,\cdots, \OL_{n+1}$ be line bundles on $X$, endowed with integrable metrics. 
The goal is to endow a metric of $\pair{L_1, \cdots, L_{n+1}}$ on $Y$ in this general setting. 
Note that $\pair{L_1, \cdots, L_{n+1}}$ is just a 1-dimensional complex vector space. 

As we do not assume that $X$ is integral,  we need to extend the definition of metrics and integrations to this setting. 
Denote by $X_{\red}$ the reduced structure of $X$. 
Denote by $X_1,\cdots, X_r$ the irreducible components of $X_\red$, endowed with the reduced structures. For $i=1,\cdots, r$, denote by $\eta_i$ the generic point of $X_i$. 

Define \emph{the multiplicity of $X_i$ in $X$} to be 
$$
\delta(X_i)=\delta(X_i,X)= \mathrm{length}_{\CO_{X,\eta_i}}(\CO_{X,\eta_i}). 
$$
See \cite[\S9.1, Def. 3]{BLR} for example.
For integrations, we define
$$
\int_X \alpha:=\sum_{i=1}^r \delta(X_i)\int_{X_i} \alpha|_{X_i}
$$
in reasonable settings to be used later.
For example, if $X$ is a finite scheme over $\CC$ (so $n=0$), then for any function $\alpha:X_\red\to \RR$, we take the convention
$$
\int_X \alpha=\sum_{i=1}^r \delta(X_i) \alpha(X_i).
$$

Most notions in \S\ref{sec complex}-\ref{sec Green} can be generalized to the current setting. 
By a \emph{continuous function} on $X$, we mean a continuous function on $X_\red$. 
By a \emph{smooth function} on $X$, mean a continuous function $g: X_\red\to \RR$, such that for any closed point $x\in X$, there is an open subscheme $U$ of $X$ containing $x$ together with a closed immersion $U\to M$ to a complex manifold $M$ such that $g|_{U}$ can be extended to a smooth function on $M$.
Let $L$ be a line bundle on $X$. 
By a \emph{continuous metric} of $L$ on $X$, we mean a continuous metric of $L|_{X_\red}$ on $X_\red$. 
By a \emph{smooth metric} $L$ on $X$, we mean a continuous metric $\|\cdot\|$ of 
$L|_{X_\red}$ on $X_\red$ such that $\|s\|^2$ is a smooth function for any local section $s$ of $X$, which is not a zero-divisor Zariski locally. 
Define \emph{Chern currents}, \emph{semipositive metrics}, 
\emph{integrable metrics} similarly. 
In terms of integration, we essentially only care about the pull-back of these terms to $X_1,\cdots, X_n$.

Let $\OL_1,\cdots, \OL_{n+1}$ be line bundles on $X$, endowed with integrable metrics. Then  $\pair{L_1, \cdots, L_{n+1}}$ is a 1-dimensional complex vector space. 
We endow a metric of $\pair{L_1, \cdots, L_{n+1}}$ as follows. 

We assume that all $L_i$ are very ample by linearity. 
For any nonzero section $s_1$ of $L_1$ on $X$, which is a regular sequence in that $s_1$ is not a zero-divisor Zariski locally on $X$, 
we have a natural isomorphism 
$$
[s_1]:\pair{L_1, \cdots, L_{n+1}} \lra \pair{L_2, \cdots, L_{n+1}}_{Z_1}.
$$
Here $Z_1=\div(s_1)$, and the right-hand side is the Deligne pairing for the morphism $Z_1\to Y$.
Define the norm of the map $[s_1]$ by 
$$
\log\|[s_1]\|=-\int_{X} \log \|s_1\| c_1(\OL_{2})\cdots c_1(\OL_{n+1}).
$$
This defines the metric of $\pair{L_1, \cdots, L_{n+1}}$ by induction on $\dim X$.

We have a few remarks to justify this definition. 
First, the integral is a sum of integrals on $X_i$ with weight $\delta(X_i)$, so we only need to consider the pull-back of the measure $c_1(\OL_{2})\cdots c_1(\OL_{n+1})$ to $X_i$.
Second, the measure $c_1(\OL_{2})\cdots c_1(\OL_{n+1})$ (over $X_i$) is defined by \cite[Thm. 2.1]{BT} (or \cite[Cor. 1.6]{Dem1}), and the integral on the right-hand side is convergent by \cite[Thm. 4.1]{CT}. 
Third, there is a Stokes formula as follows. 

\begin{lem}[Stokes formula] \label{Stokes}
 If $(s_1,s_2)$ is a strongly regular sequence of sections of $(L_1,L_2)$ on $X$. Then  
\begin{multline*}
\int_{X} \log \|s_1\| c_1(\OL_{2})c_1(\OL_{3})\cdots c_1(\OL_{n+1})
-\int_{X} \log \|s_2\| c_1(\OL_{1})c_1(\OL_{3})\cdots c_1(\OL_{n+1})\\
=\int_{\div(s_2)} \log \|s_1\| c_1(\OL_{3})\cdots c_1(\OL_{n+1})
-\int_{\div(s_1)} \log \|s_2\| c_1(\OL_{3})\cdots c_1(\OL_{n+1}).
\end{multline*}
\end{lem}
\begin{proof}
This is a generalization of \cite[I.1.3]{Elk2}, and we only sketch a proof. If $X$ is integral, the formula holds for integrable metrics by a regularization process or as an easy consequence of \cite[Thm. 4.1]{CT}.
If $X$ is not integral, by the case of integral schemes, it suffices to check that for any irreducible component $V$ of $\div(s_1)$, endowed with the reduced structure, 
$$
\delta(V,\div(s_1))=\sum_{i=1}^r\delta(V,\div(s_1|_{X_i}))\delta(X_i,X).
$$  
This is a consequence of \cite[\S 9.1, Lem. 6]{BLR} by setting $A=\CO_{X,\eta_V}$, $M=A$ and $a$ to be a defining equation of $V$ in $A$. Here $\eta_V$ denotes the generic point of $V$. 
\end{proof}

With the Stokes formula, as in \cite[Thm. I.1.1(c)]{Elk2}, we can prove that the definition of the metric is independent of the choices of the induction process, and the Deligne pairing with the metric is symmetric and multi-linear. 

In a single formula, if $(s_1,\cdots, s_{n+1})$ is a strongly regular sequence of sections of $(L_1,\cdots, L_{n+1})$ on $X$, then the metric of $\pair{L_1, \cdots, L_{n+1}}$ is given by
$$
-\log \|\pair{s_1, \cdots, s_{n+1}}\|
=-\sum_{i=1}^{n+1} \int_{Z_{i-1}} \log \|s_i\| c_1(\OL_{i+1})\cdots c_1(\OL_{n+1}).
$$
If $X$ is integral, this is exactly the local intersection number 
$$
\wh\div(s_1)\cdot \wh\div(s_2)\cdots \wh\div(s_{n+1}).
$$
See \cite[\S2]{CT} or \cite[Appendix 1]{YZ1} for basic properties of the local intersection number.

\subsection{Relation to integral schemes}

The following result converts the Deligne pairing of non-integral schemes to those of their irreducible components. It can substitute for the above treatment of non-integral schemes and will also be used later.

\begin{lem} \label{non-integral}
Let $Y$ be either the spectrum of a field or an integral Dedekind scheme. 
Let $f:X\to Y$ be a projective and flat morphism of pure relative dimension $n$.
Denote by $X_1,\cdots X_r$ the irreducible components of $X$, endowed with the reduced structures. 
Assume that for each $i=1,\cdots, r$, the morphism $X_i\to Y$ is smooth at the generic point of $X_i$. 
Let $L_1, \cdots, L_{n+1}$ be line bundles on $X$. 
Then there is a canonical isomorphism 
$$
\pair{L_1, \cdots, L_{n+1}} \lra \otimes_{i=1}^r\pair{L_1|_{X_i}, \cdots, L_{n+1}|_{X_i}}^{\otimes \delta(X_i)}.
$$

Moreover, if $Y=\Spec\CC$, $L_1, \cdots, L_{n+1}$ are endowed with integrable metrics on $X$, and both sides are endowed with the induced metrics, then the isomorphism is an isometry.
\end{lem}

\begin{proof}

We will only prove the first statement, as the second statement can be checked through the same process. 

Denote by $\eta_i$ the generic point of $X_i$, and denote by $\wt X_i$ the schematic closure of $\eta_i$ in $X$. 
Then we have a birational morphism 
$$
\coprod_{i=1}^r \wt X_i \lra X.
$$
Apply \cite[Thm. 5.3.1]{MG} to this morphism, we have a canonical isomorphism
$$
\pair{L_1, \cdots, L_{n+1}} \lra \otimes_{i=1}^r\pair{L_1|_{\wt X_i}, \cdots, L_{n+1}|_{\wt X_i}}.
$$
Therefore, it suffices to establish for each $i$ a canonical isomorphism 
$$
\pair{L_1|_{\wt X_i}, \cdots, L_{n+1}|_{\wt X_i}} \lra \pair{L_1|_{X_i}, \cdots, L_{n+1}|_{X_i}}^{\otimes \delta(X_i)}.
$$

Let $\widetilde U$ be an affine open subscheme of $\widetilde X_i$ such that
the reduced structure $U=(\widetilde U)_\red$ is smooth over $Y$.
By the infinitesimal lifting theorem (cf. \cite[\S2.2, Prop. 6]{BLR}), the identity morphism $U\to U$ can be lifted to a morphism $\phi:\widetilde U\to U$ over $Y$. 
Replacing $\wt U$ by an open subscheme if necessary, we can further assume that $\phi:\widetilde U\to U$ is flat. 
Note that $\phi:\widetilde U\to U$ is finite automatically. 
With the morphism $\phi$, computing multiplicity in terms of depths gives
$$
\delta(X_i)
=\delta(U, \widetilde U)
=\deg(\phi).$$

The morphism $\wt U\to U$ gives a rational map $\wt X_i\dasharrow X_i$. 
By blowing up $\wt X_i$, the rational map becomes a morphism $\wt X_i'\dasharrow X_i$. 
Apply the Raynaud--Gruson flattening theorem in \cite[Thm. 5.2.2]{RG}. 
We can further blow up $\wt X_i'$ and $X_i$ to change the rational map into a flat morphism $\psi:\wt Z\to Z$.
Here $\wt U$ and $U$ are respectively open subschemes of $\wt Z$ and $Z$, and $\wt Z\to Z$ extends the morphism $\wt U\to U$.
We can assume that $Z$ is normal by taking the base change of $\wt Z\to Z$ by the normalization of $Z$. The morphism $(\wt Z)_\red\to Z$ is finite, birational, and equi-dimensional, so it must be an isomorphism. 
Note that the blowing-up does not affect the Deligne pairings by \cite[Thm. 5.3.1]{MG}.

Now  it suffices to establish a canonical isomorphism 
$$
\tau:\pair{\psi^*M_1,\cdots, \psi^*M_{n+1}}
\lra \pair{M_1,\cdots, M_{n+1}}^{\otimes \deg(\psi)}
$$
for line bundles $M_1, \cdots, M_{n+1}$ on $Z$. 
Here $\psi:\wt Z\to Z$ is the finite and flat morphism. 
Write the isomorphism in the form 
$$
\pair{\psi^*M_1,\cdots, \psi^*M_{n+1}}
\lra \pair{M_1,\cdots, M_n, N_{\wt Z/Z}(\psi^*M_{n+1})}
$$
This isomorphism follows from the projection formula of \cite[Prop. 5.2.3.b]{MG}.

Now  we have established the desired isomorphism. We can also check that the isomorphism is independent of the choice of $\wt U\to U$. In fact, the morphism $\tau$ sends the section $\pair{s_1,\cdots, s_{n+1}}$ to the section $\pair{s_1|_Z,\cdots, s_{n+1}|_Z}^{\otimes \delta(Z)}$, for any strongly regular sequence 
$(s_1,\cdots, s_{n+1})$ of sections of $(\psi^*M_1,\cdots, \psi^*M_{n+1})$ on $\wt Z$ which can be descended to a strongly regular sequence of sections of $(M_1,\cdots, M_{n+1})$ on $Z$. 
Then we can check the independence by comparing different strongly regular sequences. 
\end{proof}

\subsection{Deligne Pairing: metrics in a family}

Let $f:X\to Y$ be a projective morphism of quasi-projective varieties over $\CC$ of finite Tor-dimension and pure relative dimension $n$. Let $\OL_1,\cdots, \OL_{n+1}$ be line bundles on $X$, endowed with integrable metrics. The goal is to endow a natural metric of $\pair{L_1, \cdots, L_{n+1}}$ on $Y$ in this general setting.

For any closed point $y\in Y$, we have a canonical metric $\|\cdot\|_{X_y}$ of $\pair{L_{1,y}, \cdots, L_{n+1,y}}$ at $y$. This is just the above construction applied to $f_y:X_y\to y$.
By the canonical isomorphism 
$$\pair{L_{1}, \cdots, L_{n+1}}_y\lra \pair{L_{1,y}, \cdots, L_{n+1,y}},$$
we get a natural metric of the left-hand side. 
Varying $y$, this gives a ``metric'' of $\pair{L_1, \cdots, L_{n+1}}$ on $Y$. 
Denote this metric by 
$\|\cdot\|_{X/Y,\mathrm{fibral}}$, to indicate that it is fiberwise defined. 
The metric is not a priori continuous. 
The main result of this subsection asserts that it is indeed continuous if $f$ is flat and can be ``modified'' to a continuous one if $Y$ is normal.

\begin{thm} \label{metric}
Let $f:X\to Y$ be a projective morphism of quasi-projective varieties over $\CC$ of finite Tor-dimension and pure relative dimension $n$. 
Assume that either $f$ is flat or $Y$ is normal. 
Let $\OL_1,\cdots, \OL_{n+1}$ be line bundles on $X$, endowed with integrable metrics. Then there is a continuous integrable metric $\|\cdot\|_{X/Y}$ of $\pair{L_1, \cdots, L_{n+1}}$ on $Y$ satisfying the following properties.
\begin{enumerate}
\item Let $V$ be the maximal open subscheme of $Y$ such that $X_V$ is flat over $V$. Then the metric $\|\cdot\|_{X/Y}$ is equal to the metric $\|\cdot\|_{X/Y,\mathrm{fibral}}$ at all fibers of $\pair{L_1, \cdots, L_{n+1}}$ above $V$. 

\item The metric $\|\cdot\|_{X/Y}$ is compatible with base changes by morphisms $Y'\to Y$ of quasi-projective varieties such that the image of $Y'$ intersects $V$ and that $X\times Y'\to Y'$ has finite Tor-dimension.

\item The metric $\|\cdot\|_{X/Y}$ is symmetric and multi-linear in the components $L_1, \cdots, L_{n+1}$. 

\item The Chern current 
$$c_1(\pair{L_1, \cdots, L_{n+1}}, \|\cdot\|_{X/Y})=f_*(c_1(\OL_1)\cdots c_1(\OL_{n+1}))$$ 
as $(1,1)$-currents on $Y$.

\item If the metrics of $\OL_1,\cdots, \OL_{n+1}$ are semipositive, then the metric 
$\|\cdot\|_{X/Y}$ of $\pair{L_1, \cdots, L_{n+1}}$ is also semipositive.
\end{enumerate}

\end{thm}

If $f$ is flat, then we have $\|\cdot\|_{X/Y}=\|\cdot\|_{X/Y,\mathrm{fibral}}$ everywhere, so $\|\cdot\|_{X/Y}$ is continuous. 
In this case, part (2) holds for any base change $Y'\to Y$. 

In general, $\|\cdot\|_{X/Y}$ is determined by $\|\cdot\|_{X/Y,\mathrm{fibral}}$ by continuity, but it may happen that they are not equal. 

In the theorem, by continuity, (1) determines the metric uniquely and implies (2) and (3). It is also easy to see that (4) implies (5). 
Thus, the task is to prove that the metric $\|\cdot\|_{X/Y}$ determined by (1) exists and also satisfies (4).
The proof of these two parts will be given in the next section.

\section{Metrics of the Deligne pairing: proofs}

The goal of this section is to prove Theorem \ref{metric}.
The idea is to apply Stoll and King's classical analytic results to treat the continuity of relative integrals.

\subsection{Continuity of relative integral}

As a preparation to prove Theorem \ref{metric}, 
we first convert classical results of Stoll \cite{Sto1,Sto2} and King \cite{Kin} into the following statement.

\begin{thm} \label{integral}
Let $f:X\to Y$ be a projective morphism of quasi-projective varieties over $\CC$ of pure relative dimension $n$. 
Let $\alpha$ be a continuous differential $(n,n)$-form on $X$. 
Denote by $I_{X/Y}:Y(\CC)\to \RR$ the function defined by
$$I_{X/Y}(y)=\int_{X_y} \alpha,\qquad y\in Y(\CC).$$
The following is true:
\begin{enumerate}
\item If $f$ is flat, then $I_{X/Y}$ is continuous on $Y(\CC)$.
\item If $Y$ is normal, there is a unique continuous function $\tilde I_{X/Y}:Y(\CC)\to \RR$ such that 
$\tilde I_{X/Y}(y)=I_{X/Y}(y)$ for all $y\in Y(\CC)$ over which $X$ is flat.
\end{enumerate}
\end{thm}

For a singular complex variety, there are notions of continuous differential forms and smooth differential forms in \cite[\S1.1]{Kin}. Some of these are recalled in \S\ref{sec hermitian}. 

Recall that by definition, the integration 
$$I_{X/Y}(y)=\int_{X_y} \alpha=\sum_{i=1}^r \delta(W_i) \int_{W_i} \alpha|_{W_i},$$
where $W_1,\cdots, W_r$ are irreducible components of $X_y$ endowed with reduced structures, and $\delta(W_i)$ is the multiplicity of $W_i$ in $X_y$ introduced in last subsection. 

We first recall some results of Stoll \cite{Sto1, Sto2}.
Let $f:X\to Y$ and $\alpha$ be as in the above proposition. 
Assume, furthermore, that $X$ is regular and $Y$ is normal. 
Then \cite[Thm. 3.9]{Sto2} asserts that the integral 
$$
\ds I_{X/Y}^*(y)=\int_{(X_y)_\red} \nu_f \alpha
$$
defines a continuous function of $y\in Y(\CC)$. 
Here the multiplicity function $\nu_f:(X_y)_\red(\CC)\to \ZZ$ is defined in 
\cite[p. 17, p. 48]{Sto1}. 
Instead of reviewing the definitions of the multiplicity function, 
we first state the following result, which is sufficient for our application, and then we review some details on the multiplicity function in the proof. 
 
\begin{lem}  \label{multiplicity function}
Let $f:X\to Y$ be a projective morphism of smooth varieties over $\CC$ of pure relative dimension $n$. Let $y\in Y(\CC)$ be a closed point. 
Then the following holds:
\begin{enumerate}
\item  For any smooth (closed) point $x$ of $(X_{y})_\red$, we have 
$\nu_f(x)=\delta(W(x),X_y)$. 
Here $W(x)$ is the irreducible component of $(X_{y})_\red$ containing $x$. 
\item For any continuous differential $(n,n)$-form on $X$, we have 
$$
\int_{(X_y)_\red} \nu_f \alpha
=\int_{X_y} \alpha.
$$
\end{enumerate}
\end{lem}
\begin{proof}
Note that (1) implies (2), since it implies 
$\nu_f(x)=\delta(W(x),X_y)$ for $x$ outside a subset of $(X_{y})_\red$ of measure 0. 

Now we prove the case $n=0$ of (1). 
For the purpose later, we will prove the following slightly more general statement:

\emph{Let $f:X\to Y$ be a morphism of varieties over $\CC$ with $\dim X=\dim Y$. Let $x\in X$ and $y\in Y$ be closed points with $f(x)=y$.
Assume that $X$ is smooth at $x$ and that $Y$ is smooth at $y$. 
Assume that $x$ is an isolated point of $X_y$, i.e. $\{x\}$ is a connected component of $X_y$. Then 
$$\nu_f(x)=\dim_\CC (\CO_{X,x}/m_y\CO_{X,x}).$$ 
Here $m_x$ (resp. $m_y$) denotes the maximal ideal of $\CO_{X,x}$ (resp. $\CO_{Y,y}$). }

For a brief definition of $\nu_f(x)$, recall that there is an open neighborhood $U$ of $x$ under the analytic topology such that $U\to f(U)$ is proper and $f^{-1}(y)\cap U=\{x\}$. Then $\nu_f(x)$ is the degree of $U\to f(U)$, i.e. the common order of $f^{-1}(f(z))$ for any $z\in U\setminus R$, where $R$ is an analytic subset of $U$ of positive codimension.
See also \cite[Chap. 3, Def. 3.12]{Mum2}. 

Denote by $\CO_{X,x}^\an$ (resp. $\CO_{Y,y}^\an$) the local ring of germs of analytic functions at a point $x\in X(\CC)$  (resp. $y\in Y(\CC)$). 
By \cite[Appendix to Chap. 6, Thm. A.8]{Mum2}, the formula of Weil gives
$$
\nu_f(x)=\rank_{\CO_{Y,y}^\an} \CO_{X,x}^\an.
$$
By convention, the rank of an $R$-module $M$ for an integral domain $R$ means the dimension of the base change of $M$ to the fraction field of $R$.

Note that $\CO_{X,x}^\an$ is a finite module over $\CO_{Y,y}^\an$ by \cite[Appendix to Chap. 6, Prop. A.7]{Mum2}. 
As a consequence, 
$$
\nu_f(x)=\rank_{\wh\CO_{Y,y}} \wh\CO_{X,x}.
$$
Here $\wh \CO_{X,x}$ (resp. $\wh \CO_{Y,y}$) is the completion of $\CO_{X,x}^\an$ (resp. $\CO_{Y,y}^\an$), which is canonically isomorphic to the completion of 
$\CO_{X,x}$ (resp. $\CO_{Y,y}$).

Note that $\CO_{X,x}$ is flat over $\CO_{Y,y}$ by the miracle flatness 
(cf. \cite[Thm. 23.1]{Mat}). 
It follows that $\wh\CO_{X,x}$ is flat (and finite) over $\wh\CO_{Y,y}$. 
It follows that 
$$\nu_f(x)=\dim_\CC (\wh\CO_{X,x}/m_y \wh\CO_{X,x})
=\dim_\CC (\CO_{X,x}/m_y \CO_{X,x}).$$ 
This proves the case $n=0$.

Now we prove (1) for $n>0$. The idea is to reduce it to the case $n=0$.
Assume $n>0$. 
Fix an irreducible component $W$ of $(X_y)_\red$.
Denote by $U$ an affine open subscheme of $W$, which is smooth over $\Spec(\CC)$. 
Denote by $\widetilde U$ the unique open subscheme of $X_y$ supported on $U$.
Then we have the reduced structure $U=(\widetilde U)_\red$.
By the infinitesimal lifting theorem (cf. \cite[\S2.2, Prop. 6]{BLR}), the identity morphism $U\to U$ can be lifted to a morphism $\phi:\widetilde U\to U$. 
Replacing $\wt U$ by an open subscheme if necessary, we can further assume that $\phi:\widetilde U\to U$ is flat. 
Note that $\phi:\widetilde U\to U$ is finite automatically. 
By the morphism $\phi$, a little argument gives
$$
\delta(W, X_y)
=\delta(U, \widetilde U)
=\deg(\phi).$$
This technique to treat the multiplicity is also used in the proof of Lemma \ref{non-integral}. 

We are going to prove $\nu_f(x)=\delta(W, X_y)$
for any closed point $x\in U$. 
This extends to all closed points of $W$ that are smooth in $(X_y)_\red$ by 
\cite[Thm. 5.6]{Sto1} about the global multiplicity function. 

Let $x\in U$ be any closed point. 
Let $t_1,\cdots, t_n\in \CO_{U,x}$ be a coordinate system, i.e. a minimal set of generators of the maximal ideal of the regular local ring $\CO_{U,x}$.
For $i=1,\cdots, n$, denote by $\wt t_i=\phi^*t_i\in \CO_{\wt U,x}$  the pull-back via the morphism $\phi:\widetilde U\to U$.
Denote by $t_i^*$ a lifting of $\wt t_i$ in $\CO_{X,x}$.
Then $t_1^*,\cdots, t_n^*$ are defined on an open neighborhood $W$ of $x$ in $X$.
Finally, denote by $Z$ the closed subscheme of $W$ defined by the equations $t_1^*,\cdots, t_n^*$. 
The base charge of $\phi:\widetilde U\to U$ gives a finite and flat morphism 
$\Spec (\CO_{Z\cap X_y,x})\to x$ of the same degree. 
It follows that
$$
\delta(x, Z\cap X_y)
=\deg(\phi)
=\delta(W, X_y).$$

On the other hand, by \cite[Thm. 5.5]{Sto1}, $\nu_f(x)=\nu_{f|_Z}(x)$.
By the case $n=0$ we have just proved, we further have
$\nu_{f|_Z}(x)=\delta(x, Z\cap X_y)$. 
Thus $\nu_f(x)=\delta(W, X_y).$
This finishes the proof.
\end{proof}

Now we can prove Theorem \ref{integral}.

\begin{proof}[Proof of Theorem \ref{integral}]

In (1), by Lemma \ref{closed map}, we can take a normalization and take the base change, so we will assume that $Y$ is also normal in (1).

We will start the proof with (2) and then move to (1). 
Let $f:X\to Y$ be as in (2) so that $Y$ is normal. 
By \cite[Thm. 3.3.2]{Kin}, there is a continuous function $\tilde I_{X/Y}:Y(\CC)\to \RR$ 
representing the current $f_*\alpha$.  
Recall that we also have functions $I_{X/Y}:Y(\CC)\to \RR$ and $I_{X/Y}^*:Y(\CC)\to \RR$
defined by 
$$
I_{X/Y}(y)=\int_{X_y} \alpha,\qquad
I_{X/Y}^*(y)=\int_{(X_y)_\red} \nu_f \alpha.
$$
Here we will only need  $\ds I_{X/Y}^*(y)$ for the case that $X$ and $Y$ are smooth. 
We are going to compare $\tilde I_{X/Y}$,  $I_{X/Y}$ and $I_{X/Y}^*$.

Denote by $\psi:X'\to X$ a generic desingularization of $X$. 
Then there is a Zariski open and dense subset $V_0$ of $Y$ such that $V_0$ is regular, $X$ is flat over $V_0$, and $X'$ is smooth over $V_0$. 
Shrinking $V_0$ if necessary, we can assume that for any point $y\in V_0$, the morphism $X'_y\to X_y$ is a birational morphism of reduced schemes. 
Then for any $y\in V_0(\CC)$,
$$
I_{X/Y}(y)=\int_{X_y} \alpha=\int_{X'_y} \psi^*\alpha=\int_{X'_y} \nu_{f'}\, \psi^*\alpha
=I_{X'/Y}^*(y).
$$
Here the third equality inequality follows easily from Lemma \ref{multiplicity function}. 
By \cite[Thm. 3.9]{Sto2}, $I_{X/Y}^*(y)=I_{X'/Y}^*(y)$ is continuous in $y\in V_0(\CC)$. 
This is also easy to prove directly since $X'_{V_0}$ is diffeomorphic to a constant family over $V_0$ by Ehresmann's fibration theorem.
By continuity, we have  
$I_{X/Y}(y)=\wt I_{X/Y}(y)$ for any $y\in V_0(\CC)$.

Now let $V$ be the maximal open subscheme of $Y$ such that $X_V$ is flat over $V$.
We need to prove $I_{X/Y}(y)=\wt I_{X/Y}(y)$ for any $y\in V$.
Note that $V$ contains $V_0$.
It suffices to treat the case $Y=V$; i.e. $f:X\to Y$ is flat.
This is case (1).

By taking a desingularization of $Y$ and taking the base change of $f$ accordingly, 
we can assume that $Y$ is smooth over $\CC$. This uses Lemma \ref{closed map} again.
Fix a point $y\in Y(\CC)$. 
Take a smooth curve $C\subset Y$ passing through $y$ and intersecting $V_0$.
This can be done by successively applying Bertini's theorem. 
Consider the base change $g: Z\to C$ of $f:X\to Y$ by $C\to Y$. 
Then $g$ is projective and flat. 
Since $C$ intersects $V_0$, the generic fiber of $g$ is integral. 
Then the flatness of $g$ implies that $Z$ is integral. 

Note that we need to prove $I_{X/Y}(y)=\wt I_{X/Y}(y)$ for all $y\in C(\CC)$. 
As they are equal for $y\in C(\CC)\cap V_0(\CC)$, it suffices to prove that 
$I_{X/Y}(y)$ is continuous in $y\in C(\CC)$. 
Since $I_{X/Y}(y)=I_{Z/C}(y)$ for all $y\in C(\CC)$, we only need to consider everything for the fibration $g:Z\to C$.

If $Z$ is smooth over $\CC$, this is a consequence of 
\cite[Thm. 3.9]{Sto2} and Lemma \ref{multiplicity function}. 
Otherwise, we need to take a resolution of singularity and check that $I_{Z/C}(y)$ does not change in this process.
The advantage of $\dim C=1$ is that the resolution of singularity does not violate the flatness of $Z$ over $C$. 

By Hironaka's theorem, there is a birational and projective morphism $Z'\to Z$ from a projective and smooth variety $Z'$ over $\CC$. 
We need to check that $I_{Z/C}(y)=I_{Z'/C}(y)$ for any $y\in C(\CC)$. 
Let $W$ be an irreducible component of $(Z_y)_\red$. 
Denote by $W_1',\cdots W_a'$ the irreducible components of $(Z'_y)_\red$ mapping surjectively to $W$. 
To prove $I_{Z/C}(y)=I_{Z'/C}(y)$, by pull-back of integrals, it suffices to prove 
$$
\delta(W,Z_y)=\sum_{i=1}^a \delta(W_i',Z'_y) \deg(W_i'/W). 
$$

Take a finite morphism $Z\to \PP^n_C$ over $C$, which exists by replacing $C$
by a Zariski open cover. 
The construction is similar to the construction of the morphism $\CU_{O_{F_\wp}}\to \BP^d_{O_{F_\wp}}$ in the proof of Lemma \ref{density}, so we will not repeat it here.
Denote $Z_0=\PP^n_C$ in the following.

Denote by $\eta_0$ (resp. $\eta, \eta_i'$) the generic point of $Z_{0,y}$ (resp. $W$ and $W_i'$). 
Denote by $\CO_{C,y}$, $\CO_{Z_0,\eta_0}$, $\CO_{Z,\eta}$ the local rings. 
Denote by $\CO_{Z',\eta}$ the base change $\CO_{Z'} \otimes_{\CO_{Z}} \CO_{Z,\eta}$, which is the semi-local ring of $Z'$ at the points $\eta_1',\cdots, \eta_a'$. 
All these rings are integral domains of dimension 1. 
Moreover, $\CO_{C,y}$ and $\CO_{Z_0,\eta_0}$ are discrete valuation rings. 
Then $\CO_{Z',\eta}$ and $\CO_{Z,\eta}$ are finite and flat over $\CO_{Z_0,\eta_0}$.

The inclusion $\CO_{Z,\eta} \to \CO_{Z',\eta}$ gives the same fraction fields since it comes from the birational morphism $Z'\to Z$. 
As a consequence, $\CO_{Z,\eta}$ and $\CO_{Z',\eta}$ have the same rank over $\CO_{Z_0,\eta_0}$.
Computing the degrees between the fibers above $y$, we have  
$$
\deg(\Spec \CO_{Z,\eta}/\Spec \CO_{Z_0,\eta_0})=\delta(W) \deg(\eta/\eta_0)
$$
and 
$$
\deg(\Spec \CO_{Z',\eta}/\Spec \CO_{Z_0,\eta_0})=\sum_{i=1}^a \delta(W_i') \deg(\eta_i'/\eta_0).
$$
The equality of these two degrees gives the desired result. 
The proof of Theorem \ref{integral} is complete.
\end{proof}

\subsection{Deligne Pairing: patching metrics}

Now we prove Theorem \ref{metric}. The major task is to prove part (1) of the theorem. 
Note that we have two cases: $f$ is flat, or $Y$ is normal. 
These correspond to the two cases of Theorem \ref{integral}. 

For convenience, denote 
$$
\pair{\OL_1,\cdots, \OL_{n+1}}_{\mathrm{fibral}}
=\pair{\OL_1,\cdots, \OL_{n+1}}_{X/Y,\mathrm{fibral}}
=(\pair{L_1,\cdots, L_{n+1}}, \|\cdot\|_{X/Y,\mathrm{fibral}})
$$ 
and 
$$
\pair{\OL_1,\cdots, \OL_{n+1}}
=\pair{\OL_1,\cdots, \OL_{n+1}}_{X/Y}
=(\pair{L_1,\cdots, L_{n+1}}, \|\cdot\|_{X/Y})
$$ 
in the following, the second metric is the continuous one to be constructed.

\medskip\noindent \textbf{Find a smooth metric.}
By multi-linearity, we can assume that $L_1, \cdots, L_{n+1}$ are all isomorphic to the same $f$-ample line bundle $L$ on $X$; see \cite[\S1, Step 2]{Mor2} for the argument for this reduction process. Of course, the metrics of $\OL_i$ are allowed to be very different. 

We first claim that, up to passing to a Zariski open cover of $Y$, there exists 
a smooth metric $\|\cdot\|$ of $L$, such that the induced metric 
$\|\cdot\|_{X/Y, \mathrm{fibral}}$ of $\pair{L, \cdots, L}$
is also smooth. 

Replacing $Y$ by a Zariski open cover and replacing $L$ by a tensor power if necessary, we can assume that there is a finite morphism $\psi:X\to \PP^n_Y$ over $Y$ such that $\psi^*\CO_{\PP^n_Y}(1)\simeq L$.
The construction is similar to the construction of the morphism $\CU_{O_{F_\wp}}\to \BP^d_{O_{F_\wp}}$ in the proof of Lemma \ref{density}, so we will not repeat here. 

Denote $\OM_0=(\CO_{\PP^n_\CC}(1), \|\cdot\|_\mathrm{FS})$ with the Fubini-Study metric $\|\cdot\|_\mathrm{FS}$ on $\PP^n_\CC$. 
Denote $\OM=p^* \OM_0$, where $p:\PP^n_Y\to \PP^n_\CC$ is the projection. 
Denote $\OL=\psi^*\OM$, or equivalently $\OL=(L, \|\cdot\|)$ with $\|\cdot\|=(p\circ \psi)^*\|\cdot\|_\mathrm{FS}$. 

By the base change $q:Y\to \Spec\CC$, we have a canonical isometry 
$$
q^*\pair{\OM_0,\cdots, \OM_0}_{\PP^n_\CC/\CC,\mathrm{fibral}}
\lra 
\pair{\OM,\cdots, \OM}_{\PP^n_Y/Y,\mathrm{fibral}}.
$$
As a consequence, the right-hand side is isomorphic to the trivial bundle $\CO_Y$ with a constant metric.

There is also a natural isometry 
$$
\psi^*\pair{\OM,\cdots, \OM}_{\PP^n_Y/Y,\mathrm{fibral}}
\lra
\pair{\OL,\cdots, \OL}_{X/Y,\mathrm{fibral}}.
$$
The functionality gives a natural isomorphism of the underlying line bundles and a natural isometry of the fibers, which are compatible.

As a consequence, the metric $\|\cdot\|_{X/Y, \mathrm{fibral}}$ of
$\pair{L,\cdots, L}_{X/Y}$ is smooth. 
This gives the requirement.

\medskip\noindent \textbf{Compare the metrics.}
Consider the identity map 
$$ \gamma: \pair{L,\cdots, L}\lra\pair{L_1,\cdots, L_{n+1}}.$$
We first prove Theorem \ref{metric}(1) in the case that $f$ is flat. 
Then it suffices to prove that the norm $\|\gamma\|$ of $\gamma$ under the fibral metrics is continuous on $Y$ in this case. 

For $i>1$, denote $f_i=-\log (\|\cdot\|_i/\|\cdot\|)$, which is a continuous function on $X$. Write $\gamma$ as the composition of 
$$ \gamma_i: \pair{L_1,\cdots L_{i-1}, L,\cdots, L} \lra\pair{L_1,\cdots L_{i}, L,\cdots, L}.$$
for $i=1,\cdots, n+1$. 
The norm of $\gamma$ at any $y\in Y(\CC)$ is given by 
$$
-\log \|\gamma\|(y)=
\sum_{i=1}^{n+1}
\int_{X_y} f_i c_1(\OL_1) c_1(\OL_2) \cdots c_1(\OL_{i-1}) c_1(\OL)^{n+1-i}.
$$
Denote $d=\partial+\bar\partial$ and $d^c=(\partial-\bar\partial)/(2\pi i)$. Note that some literature normalizes $d^c$ by a denominator $4\pi i$ instead of $2\pi i$. 
By $c_1(\OL_j)=c_1(\OL)+ dd^cf_j$, we see that 
$-\log \|\gamma\|(y)$
is a linear combination of 
$$
\int_{X_y} f_i (\wedge^{j\in J}dd^c f_j) \wedge c_1(\OL)^{n-|J|}.
$$
Here $i\in \{1,\cdots,n+1\}$ and $J\subset \{1,\cdots,i-1,i+1,\cdots, n+1\}$.

We are going to prove that for any $i=1,\cdots, n+1$, and for
any integrable functions $f_1,\cdots, f_i$ on $X$, the function 
$$
y\longmapsto \int_{X_y} f_1 (dd^c f_2)\wedge \cdots \wedge (dd^c f_i) \wedge c_1(\OL)^{n+1-i}
$$
is continuous in $y\in Y(\CC)$. 
Here an integrable function $f$ on $X$ is a continuous function such that the trivial bundle $\CO_X$ with the metric defined by $\|1\|=e^{-f}$ is integrable. 

If $f_2,\cdots, f_i$ are all smooth, the continuity is given by Theorem \ref{integral}. 
In general, the strategy is to approximate them by smooth functions. 
For any $j=2,\cdots, i$, by the Stokes formula, 
\begin{eqnarray*}
&&\int_{X_y} f_1 (dd^c f_2)\wedge \cdots \wedge (dd^c f_i) \wedge c_1(\OL)^{n+1-i}\\
&=&\int_{X_y} f_j (dd^c f_1) \wedge\cdots \wedge (dd^c f_{j-1})\wedge (dd^c f_{j+1})\wedge\cdots \wedge (dd^c f_i) \wedge c_1(\OL)^{n+1-i}.
\end{eqnarray*}
This is an easier version of Lemma \ref{Stokes}. 
Over any compact subset of $Y$, $f_j$ is a uniform limit of smooth functions on $X$.
Looking at the second integral, it suffices to prove the same statement, assuming that $f_j$ is smooth. 
By this method, we can assume that all $f_2,\cdots, f_i$ are all smooth. 
This proves the continuity for flat $f$. 

In the case that $Y$ is normal (but $f$ is not necessarily smooth), let $V$ be the maximal open subscheme of $Y$ 
over which $X$ is flat. 
Then we have already proved that all the relative integrals above are continuous on $V$, and it suffices to prove that they can extended to continuous functions on $Y$. 
This is proved in the same way by Theorem \ref{integral}(2).

\medskip\noindent \textbf{The Chern current.}
Once we have part (1) of Theorem \ref{metric}, it is easy to obtain part (4) of the theorem. The goal is to prove 
$$c_1(\pair{L_1, \cdots, L_{n+1}}, \|\cdot\|_{X/Y})=f_*(c_1(\OL_1)\cdots c_1(\OL_{n+1}))$$ 
as $(1,1)$-currents on $Y$.
Recall that for any metrized line bundle $(M,\|\cdot\|)$ on $Y$, the Chern current
$$
c_1(M,\|\cdot\|)= dd^c(-\log\|s\|)+ \delta_{\div(s)}
$$
for any rational section $s$ of $M$. 

Similar to the above, it suffices to prove the formula when all $L_i$ are isomorphic to a single $L$.
In the above, we have the identity map 
$$ \gamma: \pair{L,\cdots, L}\lra\pair{L_1,\cdots, L_{n+1}}.$$
Then 
$$
c_1(\pair{\OL_1,\cdots, \OL_{n+1}})
=c_1(\pair{\OL,\cdots, \OL})+ dd^c(-\log\|\gamma\|).
$$
Here if $f$ is not flat, then $\|\cdot\|_{X/Y}$ is not necessarily equal to $\|\cdot\|_{X/Y,\mathrm{fibral}}$ at a subvariety of $Y$ of positive codimension. Still, the ambiguity can be ignored in the sense of currents. 

Note that the identity 
$$c_1(\pair{\OL, \cdots, \OL})=f_*(c_1(\OL)^n)$$ 
holds as both sides are 0 since $\OL$ is constructed from a constant family. 
Considering the expression of $\log\|\gamma\|$ above in terms of the function
$$
F(y)=\int_{X_y} f_1 (dd^c f_2)\wedge\cdots\wedge (dd^c f_i) \wedge c_1(\OL)^{n+1-i}, \qquad y\in Y(\CC). 
$$
It suffices to prove that 
$$
dd^c F
=f_*( (dd^c f_1)\wedge\cdots\wedge (dd^c f_i) \wedge c_1(\OL)^{n+1-i}).
$$

Denote $d=\dim Y$. 
For any smooth and compactly supported $(d-1,d-1)$-form $\alpha$ on $Y$, we have by definition
$$\pair{dd^cF,\ \alpha}
=\int_Y F dd^c\alpha.
$$
By the expression of $F$, the right-hand side is equal to 
$$\int_X f_1 (dd^c f_2)\wedge\cdots\wedge (dd^c f_i) \wedge c_1(\OL)^{n+1-i}\wedge f^*dd^c\alpha.$$
By the Stokes formula, this becomes
$$\int_X (dd^c f_1)\wedge (dd^c f_2)\wedge\cdots\wedge (dd^c f_i) \wedge c_1(\OL)^{n+1-i}\wedge f^*\alpha,$$
which is exactly
$$
\pair{f_*( (dd^c f_1)\wedge\cdots\wedge (dd^c f_i) \wedge c_1(\OL)^{n+1-i}),\
\alpha}.
$$
As $\alpha$ is an arbitrary test form, this finishes the proof.

\section{Positivity of the Deligne pairing}

In this section, we consider the Deligne pairing for projective varieties in both the geometric case and the arithmetic case.
We will focus on some positive results later. 
For simplicity, we will only focus on the flat case. 
For clarity, we do not use uniform terminology here.

\subsection{Geometric case}

The following easy result asserts that Deligne pairing sends nef (resp. ample) line bundles to nef (resp. ample) line bundles. 
Nakayama \ proved it cite[Cor. 4.6]{Nak}, but we provide a more direct proof here. 

\begin{lem} \label{pairing nef1}
Let $f:X\to Y$ be a flat morphism of relative dimension $n$ of projective varieties of over a field $k$. Let $L_1,\cdots, L_{n+1}$ be line bundles on $X$. Then the following are true:
\begin{enumerate}
\item If $\dim Y=1$, then 
$$
\deg(\pair{L_1,\cdots, L_{n+1}})= L_1\cdot L_2\cdots L_{n+1}.
$$
\item If $L_1,\cdots, L_{n+1}$ are nef, then $\pair{L_1,\cdots, L_{n+1}}$ is nef. 
\item If $L_1,\cdots, L_{n+1}$ are ample, then $\pair{L_1,\cdots, L_{n+1}}$ is ample. 
\end{enumerate}
\end{lem}

\begin{proof}

For (1), we can assume that $Y$ is regular by taking its normalization (and taking the corresponding base change of $X\to Y$). 
We can assume that $X$ is normal by taking its normalization and applying \cite[Thm. 5.3.1]{MG}. 
By linearity, we can assume that $L_1,\cdots, L_{n+1}$ are very ample on $X$.
The intersection number of $(L_1,\cdots, L_{n+1})$ on $X$ is equal to 
$\deg(\div(s_1)\cap \cdots \cap\div(s_{n+1}))$ for a strongly regular sequence $(s_1,\cdots, s_{n+1})$ of sections of $(L_1,\cdots, L_{n+1})$ on $X$.
Many Bertini-type results guarantee the existence of strongly regular sequences in the current situation. The quickest one is the Bertini-type theorem of Seidenberg \cite{Sei} for normal varieties.
Then (1) holds essentially by definition.

To prove (2), it suffices to prove that 
$\pair{L_1,\cdots, L_{n+1}}$ has a non-negative degree on any closed integral curves $C$ in $Y$.
Take the base change of $X\to Y$ by $C\to Y$.
It suffices to compute the degree of the Deligne pairing for 
$X_C\to C$. 
If $X_C$ is integral, this is just (1).
If $X_C$ is not integral, take a finite and flat base change $C'\to C$ for some regular projective curve $C'$ so that the reduced structure of $X_{C'}$ is smooth over $C'$ at all the generic points of $X_{C'}$. 
Then we can apply Lemma \ref{non-integral} to convert to the integral case in (1).

To prove (3), assume that $L_1,\cdots, L_{n+1}$ are ample on $X$. 
Let $L$ be an ample line bundle on $Y$. 
By \cite[Prop. 5.2.1]{MG}, $\pair{f^*L, L_2,\cdots, L_{n+1}}$ is a positive multiple of $L$ and thus is ample.
If necessary, we can replace $L_1$ by a multiple, and we can assume that $L_1- f*L$ is ample on $X$. 
Now we have 
$$
\pair{L_1, L_2,\cdots, L_{n+1}}
\simeq 
\pair{f^*L, L_2,\cdots, L_{n+1}}
+ 
\pair{L_1-f^*L, L_2,\cdots, L_{n+1}}.
$$
The two terms on the right-hand sides are respectively ample and nef, so the left-hand side is ample. This finishes the proof. 
\end{proof}

Next, we introduce a mixed pairing between Cartier divisors and line bundles in a suitable situation and consider its effectiveness. 

Let $f:X\to Y$ be a projective and flat morphism of integral noetherian schemes of pure relative dimension $n$.
Let $L_1, \cdots, L_{n}$ be line bundles on $X$. 
Let $D$ be a Cartier divisor on $X$, and $\CO(D)$ be the line bundle associated with $D$. 
Let $V$ be a dense and open subvariety of $Y$, and denote by $U\to V$ the base change of $X\to Y$ by $V\to Y$. 
Assume that $D|_U$ is the trivial divisor on $U$, which gives a canonical isomorphism $\CO_U\to \CO(D)|_U$. 
There is a canonical isomorphism
$$
\pair{\CO(D),L_1,\cdots, L_n}|_V
\lra \pair{\CO(D)|_U,L_1|_U,\cdots, L_n|_U},
$$
and canonical isomorphisms 
$$
\pair{\CO(D)|_U,L_1|_U,\cdots, L_n|_U}
\lra \pair{\CO_U,L_1|_U,\cdots, L_n|_U}
\lra \CO_V.$$
Here the last map is a special case of \cite[Prop. 5.2.1.a]{MG}.
Thus, we have a canonical isomorphism 
$$
\CO_V\lra \pair{\CO(D),L_1,\cdots, L_n}|_V.
$$
This defines a rational map. 
$$
\CO_Y\dasharrow \pair{\CO(D),L_1,\cdots, L_n}
$$
and thus a rational section $s$ of $\pair{\CO(D),L_1,\cdots, L_n}$. 
Define our \emph{mixed Deligne pairing} by
$$
\pair{D,L_1,\cdots, L_n}:=\div(s),
$$
which is a Cartier divisor on $Y$, supported on $Y\setminus V$. 
Note that $\pair{D,L_1,\cdots, L_n}$ is multi-linear in $L_1,\cdots, L_n$. 

The following result concerns the pairing's effectivity, which is compatible with the general fact that the intersection number of an effective divisor with nef divisors is non-negative.

\begin{lem} \label{pairing effective1}
Let $f:X\to Y$ be a flat morphism of relative dimension $n$ of projective varieties over a field $k$. 
Let $L_1, \cdots, L_{n}$ be line bundles on $X$. 
Let $D$ be a Cartier divisor on $X$.
Let $V$ be a dense and open subvariety of $Y$, and denote by $U\to V$ the base change of $X\to Y$ by $V\to Y$. 
Assume that $D|_U$ is the trivial divisor on $U$.
Then the following are true:

\begin{enumerate}
\item If $D=f^*D_0$ for a Cartier divisor $D_0$ on $Y$, then 
$$\pair{D,L_1,\cdots, L_n}=(L_{1,\eta}\cdot L_{2,\eta} \cdots L_{n,\eta})\, D_0,$$
where $(L_{1,\eta}\cdot L_{2,\eta} \cdots L_{n,\eta})$ is the intersection numbers of $L_1, \cdots, L_{n}$ on the generic fiber of $f:X\to Y$.
\item If $Y$ is normal, $D$ is effective, and $L_1, \cdots, L_{n}$ are nef,  
then $\pair{D,L_1,\cdots, L_n}$ is effective on $Y$. 
\end{enumerate}

\end{lem}

\begin{proof}

Note that (1) is a consequence of \cite[Prop. 5.2.1.a]{MG}. 
For (2), we first assume that $L_1,\cdots, L_n$ are ample on $X$. 
As $Y$ is normal, by passing to Weil divisors, 
$\pair{D,L_1,\cdots, L_n}$ is effective on $Y$ if and only if some positive multiple of it is effective. 
Thus, we can replace $L_1,\cdots, L_n$ by positive multiples if necessary. 
Therefore, passing to a Zariski open cover of $Y$, we can find a strongly regular sequence $(s_1,\cdots, s_{n})$ of sections of $(L_1,\cdots, L_{n})$ on $X$. 
By the induction formula, this reduces the problem to $Z_n=\div(s_1)\cap \cdots \cap \div(s_n)$. 
Then the effectivity follows since the norm map from $Z_n$ to $Y$ sends global sections to global sections by \cite[Prop. 1.2.4(4)]{MG}. 
This proves the ample case.

Now we consider the case that $L_1,\cdots, L_n$ are nef on $X$. 
Let $A$ be an ample line bundle on $X$. Then we have proved that 
$$
D_m=\pair{D, mL_1+A,\cdots, mL_n+A}
$$
is effective for all positive integers $m$. 
Note that $D_m$ is a linear combination of 
the finitely many prime divisors of $Y$ supported on $Y\setminus V$. 
Then 
$$D=\lim_{m\to\infty}m^{-n}D_m$$
is effective.
\end{proof}

\subsection{Arithmetic case}

Now  we consider the arithmetic analogs of the above results. 
Let $f:\CX\to \CY$ be a flat morphism of relative dimension $n$ of projective arithmetic varieties (over $\ZZ$). 
Let $\CLL_1,\cdots, \CLL_{n+1}$ be hermitian line bundles with integrable metrics on $\CX$. 
Define their Deligne pairing
$$
\pair{\CLL_1,\cdots, \CLL_{n+1}}:= (\pair{\CL_1,\cdots, \CL_{n+1}}, \|\cdot\|_{X/Y})
$$
Here the metric on the right-hand side is given by Theorem \ref{metric}. 
This defines a functor
$$
\wh\Picc(\CX)_\intb^{n+1} \lra \wh\Picc(\CY)_\intb.
$$

We say that a hermitian line bundle $\CLL$ on a projective variety $\pi:\CX\to \Spec \ZZ$ is \emph{arithmetically positive} if the following holds:
\begin{enumerate}
\item the generic fiber $\CL_\QQ$ is ample on $\CX_\QQ$; 
\item there exist a hermitian line bundle $\CNN$ on $\Spec \ZZ$ with $\wh\deg(\CNN)>0$ such that $\CLL-\pi^*\CNN$ is nef on $\CX$.
\end{enumerate}
This is equivalent to the definition of the same notion in \S\ref{app subsec ample}. 
Recall that Theorem \ref{app NM} states
Zhang's arithmetic Nakai--Moishezon theorem holds for arithmetically positive hermitian line bundles; see also  
\cite[Cor. 4.8]{Zha1} and \cite[Cor. 5.1]{Mor7}.

\begin{lem} \label{pairing nef2}
Let $f:\CX\to \CY$ be a flat morphism of relative dimension $n$ of projective arithmetic varieties (over $\ZZ$). 
Let $\CLL_1,\cdots, \CLL_{n+1}$ be hermitian line bundles with integrable metrics on $\CX$. 
Then the following are true:
\begin{enumerate}
\item If $\dim \CY=1$, then 
$$
\wh\deg(\pair{\CLL_1,\cdots, \CLL_{n+1}})= \CLL_1\cdot \CLL_2\cdots \CLL_{n+1}.
$$
\item If $\CLL_1,\cdots, \CLL_{n+1}$ are nef, then $\pair{\CLL_1,\cdots, \CLL_{n+1}}$ is nef. 
\item If $\CLL_1,\cdots, \CLL_{n+1}$ are arithmetically positive, then $\pair{\CLL_1,\cdots, \CLL_{n+1}}$ is arithmetically positive. 
\end{enumerate}
\end{lem}

\begin{proof}
By Theorem \ref{metric}, the Deligne pairing of semipositive metrics is semipositive.
The other parts of the proof are similar to that of Lemma \ref{pairing nef1}. We omit it here. 
\end{proof}

Now  we introduce the arithmetic counterpart of Lemma \ref{pairing effective1}.
The situation is more or less included in the geometric case, except that there is an extra metric involved.

Let $f:\CX\to \CY$ be a projective and flat morphism of projective arithmetic varieties of pure relative dimension $n$.
Let $\CLL_1, \cdots, \CLL_{n}$ be hermitian line bundles on $\CX$ with integrable metrics. 
Let $\OCD$ be an arithmetic divisor on $\CX$, with an integrable Green function, and $\CO(\OCD)$ be the hermitian line bundle associated with $\OCD$. 
Let $\CV$ be a dense and open subvariety of $\CY$, and denote by $\CU\to \CV$ the base change of $\CX\to \CY$ by $\CV\to \CY$. 
Assume that $\CD|_\CU$ is the trivial divisor on $\CU$. 
As in the geometric case, we have a rational map 
$$
\CO_\CY\dasharrow \pair{\CO(\CD),\CL_1,\cdots, \CL_n}
$$
and thus a rational section $s$ of $\pair{\CO(\CD),\CL_1,\cdots, \CL_n}$. 
Define our \emph{mixed Deligne pairing} by
$$
\pair{\OCD,\CLL_1,\cdots, \CLL_n}:=\wh\div(s)
=(\div(s),-\log\|s\|),
$$
which is an arithmetic divisor on $\CY$.
The Green function uses the canonical metric of the Deligne pairing, which is simply given by 
$$
-\log \|s\|
=\int_{\CX(\CC)} g_{\OCD} \, c_1(\CLL_1)\cdots c_1(\CLL_n).
$$
As in the geometric case, $\pair{\OCD,\CLL_1,\cdots, \CLL_n}$ is multi-linear in $\CLL_1,\cdots, \CLL_n$.

The following effectivity result is the arithmetic version of Lemma \ref{pairing effective1}. 

\begin{lem} \label{pairing effective2}
Let $f:\CX\to \CY$ be a flat morphism of relative dimension $n$ of projective arithmetic varieties (over $\ZZ$). 
Let $\CLL_1, \cdots, \CLL_{n}$ be hermitian line bundles with integrable metrics on $\CX$. 
Let $\OCD$ be an arithmetic divisor on $\CX$ with an integrable Green function.
Let $\CV$ be a dense and open subvariety of $\CY$, and denote by $\CU\to \CV$ the base change of $\CX\to \CY$ by $\CV\to \CY$. 
Assume that $\CD|_\CU$ is the trivial divisor on $\CU$.
Then the following are true:

\begin{enumerate}
\item If $\OCD=f^*\OCD_0$ for an arithmetic divisor $\OCD_0$ on $\CY$, then 
$$\pair{\OCD,\CLL_1,\cdots, \CLL_n}=(\CL_{1,\eta}\cdot \CL_{2,\eta} \cdots \CL_{n,\eta})\, \OCD_0,$$
where $(\CL_{1,\eta}\cdot \CL_{2,\eta} \cdots \CL_{n,\eta})$ is the intersection numbers of $\CL_1, \cdots, \CL_{n}$ on the generic fiber of $f:\CX\to \CY$.
\item If $\CY$ is normal, $\OCD$ is effective, and $\CLL_1, \cdots, \CLL_{n}$ are nef,  
then $\pair{\OCD,\CLL_1,\cdots, \CLL_n}$ is effective on $\CY$. 
\end{enumerate}

\end{lem}

\begin{proof}
This is similar to Lemma \ref{pairing effective1}. 
In (2), the Green function  
$$
-\log \|s\|
=\int_{\CX(\CC)} g_{\OCD} \, c_1(\CLL_1)\cdots c_1(\CLL_n)
$$
is positive, since the current $c_1(\CLL_1)\cdots c_1(\CLL_n)$ is positive by the nefness of $\CLL_1, \cdots, \CLL_n$. 
\end{proof}

\section{Deligne pairing of adelic line bundles}
\label{sec pairing proof}

Now  we are ready to prove Theorem \ref{intersection2}. 
With the preparation in the previous sections, the proof here is similar to that of Proposition \ref{intersection1}. 

\begin{proof}[Proof of Theorem \ref{intersection2}]

Note that 
$$
\wh \Picc(X/k)_\intb =\varinjlim_{\CU\to\CV}  \wh \Picc(\CU/k)_\intb,
$$
where the direct limit is over all quasi-projective models $\CU\to \CV$ of $X\to Y$, i.e. projective and flat morphisms $\CU\to \CV$ extending $X\to Y$, where $\CU$ and $\CV$ are quasi-projective models of $X$ and $Y$ over $k$. 
Similar to Lemma \ref{models}, for any quasi-projective models $\CU$ and $\CV$ of $X$ and $Y$, the rational map $\CU\dashrightarrow\CV$ can be turned into a projective and flat morphism by shrinking $\CU$ and $\CV$ suitably. 
We can further assume that $\CV$ is normal.

Therefore, it suffices to prove the results for projective and flat morphisms $f:\CU \to \CV$ of quasi-projective varieties $\CU,\CV$ over $k$, where $\CV$ is assumed to be normal. 
We only need to define the functor 
$$\wh \Picc(\CU/k)_\snef^{n+1} \longrightarrow \wh \Picc (\CV/k)_\snef.$$

The functor is extended to integrable adelic line bundles by linearity. 
To extend it to nef adelic line bundles, it suffices to check that 
if $\overline\CL_1,\cdots,\overline\CL_{n+1}$ are nef on $\CU/k$, then 
$\CMM=\pair{\overline\CL_1,\cdots,\overline\CL_{n+1}}\in \wh \Picc (\CV/k)_\intb$ is nef 
on $\CV/k$.
There is a strongly nef adelic line bundle $\CNN$ on $\CU/k$ such that 
$\overline\CL_1+a\CNN,\cdots, \overline\CL_{n+1}+a\CNN$ are strongly nef for all positive rational numbers $a$.
It follows that 
$$\CMM_a=\pair{\overline\CL_1+a\CNN,\cdots,\overline\CL_{n+1}+a\CNN}$$
is strongly nef all positive rational numbers $a$.
Expanding it in terms of powers of $a$, we see that 
$$\CMM_a=\CMM+a \CNN_1+\cdots+a^{n+1} \CNN_{n+1}$$
is strongly nef
for integrable adelic line bundles $\CNN_1,\cdots, \CNN_{n+1}$ on $\CV/k$.
By integrability, there is a strongly nef adelic line bundle $\overline\CK$ such that 
$\overline\CK-\CNN_1,\cdots, \overline\CK-\CNN_{n+1}$ are strongly nef. 
As a consequence, $\CMM+(a+\cdots+a^{n+1})\overline\CK$ is strongly nef.
This implies that $\CMM$ is nef.

Now  we construct the functor for strongly nef adelic line bundles.
For the sake of the boundary topology, let $(\CY_0,\OCE_0)$ be a boundary divisor of $\CV$ over $k$.
Assume that there is a projective model $\CX_0$ of $\CU$ with a morphism 
$f_0:\CX_0\to \CY_0$ extending $f:\CU\to \CV$. 
Then $(\CX_0,f^*\OCE_0)$ is a boundary divisor of $\CU$ over $k$.

Let $\overline\CL_1,\cdots,\overline\CL_{n+1}$ be objects of $\wh\Picc (\CU/k)_\snef$.
For each $j=1,\cdots,n+1$, suppose that $\overline\CL_j$ is represented by a Cauchy sequence
$(\CL_j,(\CX_i,\overline\CL_{j,i},\ell_{j,i})_{i\geq1})$ with each $\overline\CL_{j,i}$ nef on a projective model $\CX_i$ of $\CU$ over $k$.
Here we assume that the integral model $\CX_i$ is independent of $j$, which is always possible.
For any $i\geq1$, assume that there is a projective model $\CY_i$ of $\CV$ with a morphism $f_i:\CX_i\to \CY_i$ extending $f:\CU\to \CV$.  
We assume that for each $i'>i\geq0$, we have morphisms $\CX_{i'}\to \CX_i$ and $\CY_{i'}\to \CY_i$ extending the identity maps of $\CU$ and $\CV$. 

Apply the Raynaud--Gruson flattening theorem in \cite[Thm. 5.2.2]{RG}. After blowing up $\CY_i$ and replacing $\CX_i$ by its pure transform, we can assume that $f_i: \CX_i \to \CY_i$ is flat for any $i\geq0$.
By the Deligne pairing, we have a line bundle 
$$
\CM=\pair{\CL_{1}, \CL_{2},  \cdots,  \CL_{n+1}}
$$
on $\CV$, and a
hermitian $\QQ$-line bundle 
$$
\CMM_i=\pair{\overline\CL_{1,i}, \overline\CL_{2,i},  \cdots,  \overline\CL_{n+1,i}}
$$
on $\CY_i$ for any $i\geq1$. 
The isomorphism $\ell_{j,i}:\CL_j\to \CL_{j,i}|_\CU$ induces an isomorphism 
$m_i:\CM\to \CM_i|_\CV$ of $\QQ$-line bundles on $\CV$.
By Lemma \ref{pairing nef1} and Lemma \ref{pairing nef2},  
each $\CMM_i$ is nef on $\CY_i$. 

To prove the theorem, we will define the Deligne pairing $\pair{\CLL_{1}, \CLL_{2},  \cdots,  \CLL_{n+1}}$
to be 
$$\CMM=(\CM, (\CY_i, \CMM_i,m_{i})_{i\geq1}).$$ 
For that, we need to check that $(\CM, (\CY_i, \CMM_i,m_{i})_{i\geq1})$ is indeed a Cauchy sequence in $\wh\Picc(\CV)_\rmod$.
Then it suffices to prove that $\{\wh\div(m_{i}m_1^{-1})\}_i$ is a Cauchy sequence in $\wh\Div(\CV)_\rmod$.

For any $j=1,\cdots, n+1$, 
by the Cauchy condition, 
$$
-\epsilon_i f_0^*\OCE_0\leq \wh\div(\ell_{j,i'}\ell_{j,i}^{-1}) 
\leq \epsilon_i f_0^*\OCE_0, \qquad 1\leq i\leq i'.
$$
Here $\{\epsilon_i\}_{i\geq 1}$ is a sequence of rational numbers converging to zero.

We claim that for any $i<i'$,
$$
-\epsilon_i \deg(\CU_\eta) \OCE_0\leq \wh\div(m_{i'}m_{i}^{-1}) 
\leq \epsilon_i \deg(\CU_\eta) \OCE_0
$$
in $\wh\Div(\CV)_\rmod$.
Here
$$
\deg(\CU_\eta)=\sum_{j=1}^{n+1} \deg(\CU_\eta)_j,
$$
with
$$\deg(\CU_\eta)_j= \deg (\CL_{1,\eta}\cdot\CL_{2,\eta}\cdots \CL_{j-1,\eta}\cdot \CL_{j+1,\eta} \cdots \CL_{n+1,\eta}),
$$
where $\CU_\eta\to \eta$ is the generic fiber of $f:\CU\to \CV$, and $\CL_{j,\eta}$ is the restriction of $\CL_{j}$ to $\CU_\eta$. 

The situation is similar to the proof of Proposition \ref{intersection1}. 
Note that the isomorphism $m_{i'}\circ m_{i}^{-1}:\CM_{i}|_\CV\to \CM_{i'}|_\CV$ is induced by the isomorphism $\ell_{j,i'}\circ\ell_{j,i}^{-1}:\CL_{j,i}|_\CU\to \CL_{j,i'}|_\CU$ for $j=1,\cdots, n+1$ via the construction of the Deligne pairing. 

In the following, for simplicity of notations, view line bundles on $\CX_i$ as line bundles on $\CX_{i'}$ via pull-back by abuse of notations. Apply similar conventions to $\CY_{i}$ and $\CY_{i'}$. 

Write the rational map $m_{i'}\circ m_{i}^{-1}:\CM_i\dashrightarrow \CM_{i'}$ as a composition of the rational maps 
$$
t_j:\pair{\overline\CL_{1,i'}, \cdots, \overline\CL_{j-1,i'}, \overline\CL_{j,i},\cdots,  \overline\CL_{n+1,i}}
\dashrightarrow 
\pair{\overline\CL_{1,i'}, \cdots, \overline\CL_{j,i'}, \overline\CL_{j+1,i},\cdots,  \overline\CL_{n+1,i}}
$$
for $j=1,\cdots, n+1$,
which are induced by the natural isomorphisms on $\CU$.
View $t_j$ as a rational section of 
$$
\pair{\overline\CL_{1,i'}, \cdots, \overline\CL_{j,i'}, \overline\CL_{j+1,i},\cdots,  \overline\CL_{n+1,i}}
-\pair{\overline\CL_{1,i'}, \cdots, \overline\CL_{j-1,i'}, \overline\CL_{j,i},\cdots,  \overline\CL_{n+1,i}},$$
which is canonically isomorphic to 
$$
\CNN_j=\pair{\overline\CL_{1,i'}, \cdots, \overline\CL_{j-1,i'}, \overline\CL_{j,i'}-\overline\CL_{j,i}, \overline\CL_{j+1,i},\cdots,  \overline\CL_{n+1,i}}$$
over $\CY_{i'}$.
It suffices to prove 
$$
-\epsilon_i \deg(\CU_\eta)_j \OCE_0\leq \wh\div(t_j) 
\leq \epsilon_i \deg(\CU_\eta)_j \OCE_0
$$
in $\wh\Div(\CV)_\rmod$.

The line bundle, $\CNN_j$, fits the framework of Lemma \ref{pairing effective1} and  Lemma \ref{pairing effective2}.
In terms of the mixed Deligne pairing, we exactly have 
$$
\wh\div(t_j) 
=\pair{ \wh\div(\ell_{j,i'}\ell_{j,i}^{-1}),\,
\overline\CL_{1,i'}, \cdots, \overline\CL_{j-1,i'}, \overline\CL_{j+1,i},\cdots,  \overline\CL_{n+1,i}}.
$$
Apply Lemma \ref{pairing effective1} and  Lemma \ref{pairing effective2}. We get 
$$
\wh\div(t_j) 
\leq \pair{ \epsilon_i f_0^*\OCE_0, \,
\overline\CL_{1,i'}, \cdots, \overline\CL_{j-1,i'}, \overline\CL_{j+1,i},\cdots,  \overline\CL_{n+1,i}}
=\epsilon_i \deg(\CU_\eta)_j \OCE_0
$$
by the Cauchy condition
$$
-\epsilon_i f_0^*\OCE_0\leq \wh\div(\ell_{j,i'}\ell_{j,i}^{-1}) 
\leq \epsilon_i f_0^*\OCE_0.
$$
Similarly, we have 
$$
\wh\div(t_j) 
\geq -\epsilon_i \deg(\CU_\eta)_j \OCE_0.
$$
It finishes the proof.
\end{proof}

\section{More functorialities of the pairing} \label{sec local pairing}

In Theorem \ref{intersection2}, we have listed that the Deligne pairing is compatible with base change. In this section, we list two more natural properties.
The first one is the behavior of the pairing in some situations under compositions, and the second one is a non-archimedean local version of the pairing.

\subsection{Functoriality properties} \label{sec local pairing}

We first present various projection formulas on Deligne pairings. 
To avoid confusion, for a morphism $\psi:X\to Y$, we use write $\psi_*\pair{\cdots}$
for the Deligne pairing for this morphism. 

\begin{lem} \label{basic2}
\kkk
Let $\psi:X\to Y$ be a projective morphism of relative dimension $r$ over $k$, and 
$\pi:Y\to S$ be a projective morphism of relative dimension $m$ over $k$.
Here $X, Y, S$ are quasi-projective and flat integral schemes over $k$. 
Assume that $\pi:Y\to S$ and $\pi\circ\psi:X\to S$ are flat, and assume that $S$ is normal. 
Let $\OL_1,\cdots, \OL_{r+1}$ be integrable adelic line bundles on $X$,  $\OM_1,\cdots, \OM_{m+1}$ be integrable adelic line bundles on $Y$, and $\ON_1$ be an integrable adelic line bundles on $S$.
\begin{enumerate}
\item
Assume that $\psi:X\to Y$ is flat, and assume that $Y$ is normal. 
Then there is a canonical isomorphism  
\begin{align*}
& (\pi\circ\psi)_*\pair{\OL_1,\cdots, \OL_{r+1}, \psi^*\OM_1,\cdots, \psi^*\OM_m}\\
\lra&\pi_*\pair{ \psi_*\pair{\OL_1,\cdots, \OL_{r+1}}, \OM_1,\cdots, \OM_m}.
\end{align*}

\item
There is a canonical isomorphism  
$$
\pi_*\pair{\OM_1,\cdots, \OM_{m},\pi^*\ON_1}
\lra e\, \ON_1.
$$
Here $e$ is the intersection number of the underlying line bundles of $\OM_1,\cdots, \OM_m$ on the generic fiber of $\pi:Y\to S$. 
\item
There is a canonical isomorphism  
$$
(\pi\circ\psi)_*\pair{\OL_1,\cdots, \OL_{r}, \psi^*\OM_1,\cdots, \psi^*\OM_{m+1}} \lra 
d\,\pi_*\pair{ \OM_1,\cdots, \OM_{m+1}}.
$$
Here $d$ is the intersection number of the underlying line bundles of $\OL_1,\cdots, \OL_r$ on the generic fiber of $\psi:X\to Y$ if $\psi$ is surjective;
set $d=0$ if $\psi$ is not surjective. 

\item 
If $r=0$, then there is a canonical isomorphism  
$$
(\pi\circ\psi)_*\pair{\psi^*\OM_1,\cdots, \psi^*\OM_{m+1}} \lra 
\deg(\psi)\,\pi_*\pair{ \OM_1,\cdots, \OM_{m+1}}.
$$
Here $\deg(\psi)$ is degree of the extension between the function fields of $X$ and $Y$ induced by $\psi:X\to Y$ if $\psi$ is surjective; set $\deg(\psi)=0$ if $\psi$ is not surjective. 
\end{enumerate}

\end{lem}

\begin{proof}
We first prove (1). 
See \cite[Prop. 5.2.3.b]{MG} for the isomorphism of the underlying line bundles. 
By taking the limit, this already implies the result for the geometric case that $k$ is a field. 

If $k=\ZZ$, we need an extra argument to check the compatibility of the hermitian metrics. 
By Theorem \ref{metric}, metrics of Deligne pairings are fiberwise defined, so it suffices to check the equality of the metrics assuming that $S=\Spec \CC$, and $\OL_i,\OM_j$ are metrized line bundles on the complex varieties $X,Y$.
Induct on $m=\dim Y$. 
By linearity, assume that $M_m$ is very ample on $Y$, and take a section $s\in \Gamma(Y, M_m)$ such that $Y'=\div(s)$ is integral, and $X'=X\times_YY'$ is also integral. 
Denote by $\psi':X'\to Y'$ and $\pi':Y'\to S$ the morphisms.
Denote $\OL_i=\OL_i|_{X'}$ and $\OM_j'=\OM_j|_{Y'}$. 
By induction, we have an isometry  
\begin{align*}
& (\pi'\circ\psi')_*\pair{\OL_1',\cdots, \OL_{r+1}', \psi^*\OM_1',\cdots, \psi^*\OM_{m-1}'} \\ \lra &\,
\pi'_*\pair{ \psi'_*\pair{\OL_1',\cdots, \OL_{r+1}'}, \OM_1',\cdots, \OM_{m-1}'}.
\end{align*}
It suffices to check that the changes in the metrics of both sides are equal. 
By \S\ref{sec metric at point}, this amounts to check
\begin{align*}
& \int_X \log\|\psi^*s\| c_1(\OL_1)\cdots c_1(\OL_{r+1})c_1(\psi^*\OM_1)\cdots c_1(\psi^*\OM_{m-1})\\
=&\int_Y \log\|s\| c_1(\psi_*\pair{\OL_1,\cdots, \OL_{r+1}}) c_1(\OM_1)\cdots c_1(\OM_{m-1}).
\end{align*}
This follows from Theorem \ref{metric}(4).
It proves (1).

Note that (2) is the special case of (3) when $\pi$ is the identity map on $Y$, and (4) is also a special case of (3). We also have two quick proofs of (2). 
First, \cite[Prop. 5.2.1.a]{MG} gives an isomorphism of the underlying line bundles in (2), and then we can extend it to the adelic case as in the proof of (1).
Alternatively, write $\ON_1=\CO(\ol D)$ for some adelic divisor $\ol D$ on $S$. 
Then the result follows from limit versions of Lemma \ref{pairing effective1}(1) and Lemma \ref{pairing effective2}(1).

To prove (3), the key is to establish a canonical isomorphism 
$$
(\pi\circ\psi)_*\pair{L_1,\cdots, L_{r}, \psi^*M_1,\cdots, \psi^*M_{m+1}} \lra 
d\,\pi_*\pair{ M_1,\cdots, M_{m+1}}
$$
of the underlying line bundles since the extension of this to the adelic case is similar to that in (1). 
For the isomorphism of the underlying line bundles, we can assume that all $M_i, L_j$ are very ample by linearity. By passing to a finite Zariski open cover of $S$, we can find a section $t_{m+1}$ of $M_{m+1}$ such that $\div(t_{m+1})$ is integral and flat over $S$. This reduces the problem from $(X,Y, S)$ to $(X_{\div(t)}, \div(t), S)$, and thus eventually we can assume that $m=0$. 
Then $\pi:Y\to S$ is finite and flat, and by passing to a Zariski open cover of $S$,
we can take a section of $L_r$ to reduce $\dim X$, and eventually we can also assume that $r=0$. 
The case $m=r=0$ can be checked by an easy relation of the norm maps. 
This finishes the proof.
\end{proof}

\subsection{Local theory} \label{sec local pairing}

In this subsection, we are going to consider the Deligne paring over a non-archimedean field and treat the metrics in this setting.

Let $K$ be a non-archimedean field with a discrete valuation, and let $O_K$ be the valuation ring. 
Let $X$ be a \emph{projective} variety of dimension $n$ over $K$. 
Recall that in \S\ref{sec adelic general}, we have defined $\wh\Picc(X/O_K)$ as the completion of 
$$
\wh\Picc(X/O_K)_\rmod=\varinjlim_{\CX} \Picc(\CX,X)
$$
along the boundary topology. Here the limit is over projective models $\CX$ of $X$ over $O_K$.

Similar to the global case, we can introduce the category $\wh\Picc(X/O_K)_\snef$ (resp. $\wh\Picc(X/O_K)_\nef$, $\wh\Picc(X/O_K)_\intb$)
of \emph{strongly nef} (resp. \emph{nef}, \emph{integrable}) objects of $\wh\Picc(X/O_K)$.
A line bundle on a projective model, $\CX$ of $X$ over $O_K$, is \emph{nef} if it has a non-negative degree on every projective and integral curve in the special fiber of $\CX\to \Spec O_K$. 
An adelic line bundle on $X$ is \emph{strongly nef} if it is the limit under the boundary topology of model adelic line bundles induced by nef line bundles on projective models of $X$. 
An adelic line bundle, $\OL$ on $X$, is \emph{nef} if there exists a strongly nef adelic line bundle $\OM$ on $X$ such that $a\OL+\OM$ is strongly nef for all positive integers $a$.
An adelic line bundle on $X$ is \emph{integrable} if it is isomorphic to the difference of two strongly nef adelic line bundles on $X$. 

By continuity, the Deligne pairing
$$
\Picc(\CX)^{n+1}\lra \Picc(O_K)
$$
extends to a canonical pairing
$$
\wh\Picc(X/O_K)_\intb^{n+1}\lra \wh\Picc(K/O_K).
$$
Note that $\wh\Picc(K/O_K)=\wh\Picc(\Spec K/O_K)$ is equivalent to the category of triples $(L, \CL, \ell)$ with 
$L\in \Picc(K)$, $\CL\in \Picc(O_K)_\QQ$, and $\ell: L\to \CL|_{\Spec K}$ an isomorphism in $\Picc(K)_\QQ$. 
If $O_K$ is a discrete valuation ring, the proof of this extension is similar to Theorem \ref{intersection2}, and it is easier without the Archimedean metrics. 
If $O_K$ is not a discrete valuation ring, an extra ingredient of the extension is from Xia \cite{Xia}, who extends the Deligne pairing of \cite{Del, Elk1, MG, Duc} to non-noetherian schemes. 
The idea of \cite[Prop. 3.7]{Xia} is that a projective and flat morphism of (possibly non-noetherian) schemes can be Zariski locally descended to a projective and flat morphism of noetherian schemes.

On the other hand, in Proposition \ref{injection local}, we have a canonical fully faithful functor
$$\wh \Picc (X/O_K)\lra \wh\Picc (X^\an).$$
We have essential images
$$\wh \Picc (X^\an)_\cptf=\Im(\wh \Picc (X/O_K)\to \wh \Picc (X^\an)),$$
$$\wh \Picc (X^\an)_\snef=\Im(\wh \Picc (X/O_K)_\snef\to \wh \Picc (X^\an)),$$
$$\wh \Picc (X^\an)_\nef=\Im(\wh \Picc (X/O_K)_\nef\to \wh \Picc (X^\an)),$$
$$\wh \Picc (X^\an)_\intb=\Im(\wh \Picc (X/O_K)_\intb\to \wh \Picc (X^\an)).$$
Note that in \S\ref{sec local theory}, we have also given direct descriptions of compactified (resp. strongly nef, nef, integrable) metrized line bundles on $X^\an$. 

Parallel to the archimedean setting in \S\ref{sec metric at point}, we use integration to define a Deligne pairing 
$$
\wh\Picc(X^\an)_\intb^{n+1}\lra \wh\Picc(K^\an).
$$
In fact, let $\OL_1,\cdots, \OL_{n+1}$ be integrable metrized line bundles on $X^\an$, we endow a metric of the 1-dimension $K$-space $\pair{L_1, \cdots, L_{n+1}}$ as follows. 

We assume that all $L_i$ are very ample by linearity. 
Let $s_1$ be a nonzero section of $L_1$ on $X$.
We have a natural isomorphism 
$$
[s_1]:\pair{L_1, \cdots, L_{n+1}} \lra \pair{L_2, \cdots, L_{n+1}}_{Z_1/K}.
$$
Here $Z_1=\div(s_1)$ and the right-hand side is the Deligne pairing for the morphism $Z_1\to \Spec K$.
Define the norm of the map $[s_1]$ by 
$$
\log\|[s_1]\|=-\int_{X^\an} \log \|s_1\| c_1(\OL_{2})\cdots c_1(\OL_{n+1}).
$$
Here the right-hand side uses the Chambert-Loir measure.
This defines the metric of $\pair{L_1, \cdots, L_{n+1}}$ by induction on $\dim X$.
Note that if $Z_1$ is not integral, we can use Lemma \ref{non-integral} to convert it to the pairings from its irreducible components.
As in the Archimedean case, the definition is independent of the choice of $s_1$ by \cite[Thm. 4.1]{CT}, and the pairing is symmetric and multi-linear. 

Similar to the Archimedean case, in a single formula, if $(s_1,\cdots, s_{n+1})$ is a strongly regular sequence of sections of $(L_1,\cdots, L_{n+1})$ on $X$, then the metric of $\pair{L_1, \cdots, L_{n+1}}$ is given by
$$
-\log \|\pair{s_1, \cdots, s_{n+1}}\|
=-\sum_{i=1}^{n+1} \int_{Z_{i-1}^\an} \log \|s_i\| c_1(\OL_{i+1})\cdots c_1(\OL_{n+1}).
$$
This is exactly the local intersection number 
$$
\wh\div(s_1)\cdot \wh\div(s_2)\cdots \wh\div(s_{n+1}).
$$
See \cite[\S2]{CT} or \cite[Appendix 1]{YZ1} for basic properties of the local intersection number.

Finally, we have the following result, which asserts that the two pairings are compatible.
\begin{thm} \label{local compatibility}
Let $K$ be a non-archimedean field, and let $O_K$ be its valuation ring. 
Let $X$ be a projective variety of dimension $n$ over $K$. 
The Deligne pairings 
$$
\wh\Picc(X/O_K)_\intb^{n+1}\lra \wh\Picc(K/O_K)
$$
and 
$$
\wh\Picc(X^\an)_\intb^{n+1}\lra \wh\Picc(K^\an)
$$
are compatible with the analytification functors
$$
\wh\Picc(X/O_K)_\intb\lra \wh\Picc(X^\an)_\intb
$$
and
$$
\wh\Picc(K/O_K)\lra \wh\Picc(K^\an).
$$
\end{thm}

\begin{proof}
It suffices to prove the model case. Namely, let $\CL_1,\cdots, \CL_{n+1}$ be line bundles on a projective model $\CX$ of $X$ over $O_K$, with generic fibers $L_1,\cdots, L_{n+1}$ on $X$. Then the metric of $\pair{L_1,\cdots, L_{n+1}}$ on $X^\an$ induced by $\pair{\CL_1,\cdots, \CL_{n+1}}$ is equal to the one defined by the integrals. 

To prove the model case, assume that $\CX$ is normal by taking a normalization. 
Assume that all $\CL_i$ are very ample by linearity. 
Let $s_1$ be a nonzero section of $\CL_1$ on $\CX$ such that 
$\CZ_1=\div(s_1)$ is flat over $O_K$. 
Then we have a natural isomorphism 
$$
[s_1]:\pair{\CL_1, \cdots, \CL_{n+1}} \lra \pair{\CL_2, \cdots, \CL_{n+1}}_{\CZ_1/O_K}.
$$
Thus $\pair{\CL_1, \cdots, \CL_{n+1}}$ and 
$\pair{\CL_2, \cdots, \CL_{n+1}}_{\CZ_1/O_K}$
induce compatible metrics on the line bundles $\pair{L_1, \cdots, L_{n+1}}$ and $\pair{L_2, \cdots, L_{n+1}}_{\CZ_{1,K}/K}.$
By induction, it suffices to prove that the analytic term
$$
\log\|[s_1]\|=-\int_{X^\an} \log \|s_1\| c_1(\OL_{2})\cdots c_1(\OL_{n+1})
$$
vanishes.
By definition, the Chambert-Loir measure on the right-hand side is supported on the divisorial points of $X^\an$ corresponding to the irreducible components of the special fiber of $\CX$. 
On the other hand, $\|s_1\|=1$ at these divisorial points by the assumption that 
$\div(s_1)$ is flat over $O_K$. 
Then the integral vanishes.
\end{proof}

\chapter{Volumes and heights}

In this chapter, we are going to study effective sections of adelic line bundles, volumes of adelic line bundles, heights of algebraic points and subvarieties, and equidistribution of small points. 
It turns out that many definitions and results for hermitian line bundles can be extended 
to the current situation.

As before, we will treat the geometric case and the arithmetic case uniformly, taking the uniform terminology in \S\ref{sec uniform}.

\section{Effective sections of adelic line bundles} \label{sec HS}

The goal of this section is to introduce effective sections of adelic line bundles and derive some basic finiteness properties.

\subsection{Effective adelic divisors}

Effective sections of adelic line bundles are defined in terms of effective adelic divisors, so we start with the following definition.

\begin{defn} \label{def effectivity}
\kkk
\begin{enumerate}
\item 
Let $\CU$ be a quasi-projective variety over $k$. 
An adelic divisor $\OCD$ in $\wh\Div(\CU/k)$ is called \emph{effective} if it can be represented by a Cauchy sequence of effective divisors in $\wh\Div(\CU/k)_\rmod$. 
\item
Let $X$ be a flat and essentially quasi-projective integral scheme over $k$. 
An adelic divisor $\OD$ in $\wh\Div(X/k)$ is called \emph{effective} if it is the image of an effective adelic divisor of $\wh\Div(\CU/k)$ for some quasi-projective model $\CU$ of $X$. 
\end{enumerate}
\end{defn}

As before, we will use $\geq$ and $\leq$ to denote the partial order on $\wh\Div(X/k)$ induced by effectivity.

Note that the above definition is very similar to the definition of strong nefness.
The following is the justification for this definition.

\begin{lem} \label{effectivity2}
\kkk
Let $X$ be a flat and essentially quasi-projective integral scheme over $k$.
Then an adelic divisor $\OD\in \wh\Div(X/k)$ is effective if and only if its image 
$\OD^\an$ in $\wh\Div(X^\an)$ is effective.

If furthermore $X$ is normal, then $\OD$ is effective if and only the total Green function $\wt g_{\OD}$ induced by $\OD$ is non-negative on $X^\an\setminus |D|^\an$.
\end{lem}

\begin{proof}
It suffices to prove the case that $X=\CU$ is a quasi-projective variety over $k$, and we write  $\OCD$ for $\OD$ as a convention.
For the first statement, assume that $\OCD^\an$ is effective, and we need to prove that $\OCD$ is effective. 
Assume that $\OCD$ is represented by a sequence $\{\OCD_i\}_{i\geq1}$ in 
$\wh\Div(\CU/k)_\rmod$. 
By definition, there is a sequence $\{\epsilon_j\}_{j\geq 1}$ of positive rational numbers converging to $0$ such that 
$$
 -\epsilon_j \OCE_0 \leq
\OCD_i-\OCD_j \leq  \epsilon_j \OCE_0,\quad\ i\geq j\geq 1.
$$
This implies 
$$
 -\epsilon_j \tilde g_{\CE_0} \leq
\tilde g_{\CD_i}-\tilde g_{\CD_j}\leq  \epsilon_j \wt g_{\CE_0},\quad\ i\geq j\geq 1.
$$
Here $\wt g_\bullet$ denotes the corresponding Green function on $\CU^\an$.
Set $i\to \infty$. It gives 
$$
\tilde g_{\CD_j}+ \epsilon_j \tilde g_{\CE_0} \geq 
\tilde g_{\CD} \geq 0,\quad\  j\geq 1.
$$
By Lemma \ref{effectivity1}, $\OCD_j+ \epsilon_j \OCE_0$ is effective in 
$\wh\Div(\CU/k)_\rmod$.
Note that $\OCD$ is also represented by the Cauchy sequence $\{\OCD_i+\epsilon_i \OCE_0\}_{i\geq1}$ in 
$\wh\Div(\CU/k)_\rmod$. 
Then it is effective.

For the second statement, it suffices to prove that $\wt g_{\OCD}\geq 0$ implies $\CD\geq0$. This can be proved as in the proof of Corollary \ref{shrink}.  
\end{proof}

\subsection{Effective sections of adelic line bundles}

\kkk
Let $X$ be a flat and essentially quasi-projective integral scheme over $k$, and let $\OL$ be an adelic line bundle on $X$. For any nonzero rational section $s$ of $L$ on $X$, there is an arithmetic divisor $\wh\div(s)$, defined as an element of $\wh\Div(X/k)$. 
It suffices to define this for any quasi-projective model $\CU$ of $X$. This is in Lemma \ref{isomorphism}. 
Namely, if $\CLL=(\CL,(\CX_i,\overline \CL_i, \ell_{i})_{i\geq1})$ is an adelic line bundle on $\CU$, and $s$ is a nonzero rational section of $\CL$ on $\CU$, then
$$
\wh\div(s)=\wh\div_{(\CX_1,\CLL_1)}(s)+\lim_{i\to \infty} \wh\div(\ell_i\ell_1^{-1})
$$
in $\wh\Div(\CU/k)$. 

Now  we are ready to introduce the key definitions. Note that the term $\wh h^0$ below is defined to be either finite real numbers or infinity. 

\begin{defn}\label{def effective sections}
\kkk 
Let $X$ be a flat and essentially quasi-projective integral scheme over $k$.
\begin{enumerate}
\item
Let $\OL$ be an adelic line bundle on $X$ with underlying line bundle $L$ on $X$.  
Define 
$$\wh H^0(X, \OL):=\{s\in H^0(X, L): \wh\div(s)\geq 0\}.$$
Here the partial order is in $\wh\Div(X/k)$. 
Elements of $\wh H^0(X, \overline L)$ are called \emph{effective sections} of $\overline L$ on $X$. 
If $k=\ZZ$, denote 
$$\wh h^0(X, \OL):=\log\#\wh H^0(X, \OL);$$ 
if $k$ is a field, denote 
$$\wh h^0(X, \OL):=\dim_k 
\wh H^0(X, \OL).$$
We say that $\OL$ if \emph{effective} if $\wh h^0(X, \OL)>0$. 

\item
Let $\overline L$ be a metrized line bundle on $X^\an$ with underlying line bundle $L$.
For any $s\in H^0(X, L)$ and any $v\in \CMk$, define the supremum norms
$$
\|s\|_{\sup}:= \sup_{x\in X^\an} \|s(x)\|, $$
$$
\|s\|_{v,\sup}:= \sup_{x\in X_v^\an} \|s(x)\|.
$$
Both ``$\sup$'' are allowed to be infinity.
Define 
$$\wh H^0(X, \overline L):=\{s\in H^0(X, L): \|s\|_{\sup}\leq 1\}.$$
Elements of $\wh H^0(X, \overline L)$ are called \emph{effective sections} of $\overline L$ on $X$. 
If $k=\ZZ$, denote 
$$\wh h^0(X, \OL):=\log\#\wh H^0(X, \OL);$$ 
if $k$ is a field, denote 
$$\wh h^0(X, \OL):=\dim_k \wh H^0(X, \OL).$$
\end{enumerate}
\end{defn}

The definitions in (1) and (2) are compatible. Namely, for any adelic line bundle $\OL$ on $X$, which induces a metrized line bundle $\OL^\an$ on $X^\an$,
 the canonical map 
$$
\wh H^0(X, \OL) \lra \wh H^0(X, \OL^\an)
$$
is bijective.
This follows from Lemma \ref{effectivity2}.

If $k$ is a field, it is easy to see that $\wh H^0(X, \OL)$ is a vector space over $k$ in the setting of (2), and thus the same holds in the setting of (1). So the dimension $\wh h^0(X, \OL)$ is well-defined (as a finite number or infinity).

In both the arithmetic case and the geometric case,
we will prove that the number $\wh h^0(X, \OL)$ in (1) is finite.
The proof is easy, but we will postpone it till Lemma of {finiteness3} to set up a framework to bound sections of adelic line bundles.

In terms of arithmetic divisors, the definitions are written more easily. 
For example, if $\OD$ is an adelic divisor on $X$, 
then 
$$\wh H^0(X, \CO(\OD))=\{f\in k(X)^\times: \wh\div(f)+ \OD\geq 0\}\cup\{0\}.$$
Because of this, we may work on adelic divisors instead of adelic line bundles. 
For simplicity, we will denote 
$$\wh H^0(X, \OD)=\wh H^0(X, \CO(\OD)),\qquad
\wh h^0(X, \OD)=\wh h^0(X, \CO(\OD)).$$

\subsection{Effective sections of arithmetic $\QQ$-divisors} \label{sec alternative sections}

For the purpose later, we generalize the definition of effective sections to arithmetic $\QQ$-divisors. 

\kkk
Let $\CX$ be a projective variety over $k$.
Let $\OCD$ be an arithmetic $\QQ$-divisor on $\CX$.
Denote 
$$\wh H^0(\CX, \OCD)'=\{f\in k(\CX)^\times: \wh\div(f)+ \OCD\geq 0
\text{ in } \wh\Div(\CX)_\QQ
\}\cup \{0\}.$$
If $k=\ZZ$, denote 
$$\wh h^0(\CX, \OCD)':=\log\#\wh H^0(\CX, \CO(\OCD))';$$ 
if $k$ is a field, and if $\CX$ is normal, denote 
$$\wh h^0(\CX, \OCD)':=\dim_k H^0(\CX, \CO(\OCD))'.$$

Note that if $k$ is a field and if $\CX$ is normal, $\wh H^0(\CX, \OCD)'$ is a $k$-vector space, since $\wh\div(f)+ \OCD\geq 0$ is equivalent to
$\ord_v(f)\geq -\ord_v (\OCD)$ for all codimension one points $v\in\CX$. 
However, if $\CX$ is not normal, $\wh H^0(\CX, \OCD)'$ might fail to be a group.

If $\OCD$ is integral, recall the usual set of effective sections defined by
$$\wh H^0(\CX, \OCD):=\{f\in k(\CX)^\times: \wh\div(f)+ \OCD\geq 0
\text{ in } \wh\Div(\CX)
\}\cup \{0\}.$$
There is a canonical injection 
$$\wh H^0(\CX, \OCD)\to\wh H^0(\CX, \OCD)',$$
which might fail to be bijective due to the difference of the effectivity relations in $\wh\Div(\CX)$ and $\wh\Div(\CX)_\QQ$. 
However, if $\CX$ is normal, then it is bijective by Lemma \ref{effectivity normal}.

In the definitions above, $\wh h^0(\CX, \OCD)'$ is always finite.
In fact, the finiteness holds if $\CX$ is normal and $\OCD$ is integral by $\wh h^0(\CX, \OCD)'=\wh h^0(\CX, \OCD)$. 
The normality condition can be removed by passing to the normalization, and the integrality condition can be removed by bounding $\OCD$ by an integral arithmetic divisor on $\CX$ under the relation ``$\leq$''. 
In fact, the following technical lemma guarantees the existence of such integral arithmetic divisors. 

\begin{lem}  \label{divisor relations}
Let $\CX$ be a projective variety over $k$, and let $\CU$ be an open subscheme of $\CX$.
Then the following are true. 
\begin{enumerate}
\item Let $\CDD\in \wh\Div(\CX)_\QQ$ be an arithmetic $\QQ$-divisor on $\CX$. 
Then there are $\OCD_1, \OCD_2\in \wh\Div(\CX)$ satisfying 
 $\OCD_1\leq 0\leq \OCD_2$ in $\wh\Div(\CX)$ and $\OCD_1\leq\OCD\leq\OCD_2$ in $\wh\Div(\CX)_\QQ$.

\item
Let $\OCD\in  \wh\Div(\CX,\CU)$ be an arithmetic $(\QQ,\ZZ)$-divisor on $(\CX,\CU)$. 
Then there are $\OCD_1, \OCD_2\in \wh\Div(\CX)$
satisfying 
 $\OCD_1\leq 0\leq \OCD_2$ in $\wh\Div(\CX)$
and $\OCD_1\leq\OCD\leq\OCD_2$ in $\wh\Div(\CX,\CU)$. 
\end{enumerate}
\end{lem} 

\begin{proof}
We claim that \emph{for any divisor $\CE$ on $\CU$, there is an effective divisor $\CE'$ on $\CX$ such that $\CE'|_\CU\geq \CE$ in $\Div(\CU)$.}

We first prove the claim assuming that $\CE$ is effective on $\CU$. 
In this case, the zero locus of $\CE$ defines a closed subscheme $\CZ$ of $\CU$. 
Extend $\CZ$ to a closed subscheme $\wt\CZ$ of $\CX$ by taking the schematic closure. 
Take a very ample line bundle $\CA$ on $\CX$ with a nonzero global section $s$ vanishing along $\wt\CZ$; i.e. $s$ lies in the kernel of $H^0(\CX,\CA)\to H^0(\wt\CZ,\CA|_{\wt\CZ})$. 
Then we can set $\CE'=\div(s)$ on $\CX$.

Next, we prove the claim for general $\CE$. 
In fact, take a finite affine open cover $\{\CU_i\}_i$ of $\CU$ such that 
$\CE|_{\CU_i}=\div(f_i)$ for some $f_i\in k(\CU_i)^\times$ for every $i$. 
Write $f_i=f_i'/f_i''$ with $f_i', f_i''\in \Gamma(\CU_i, \CO_{\CU_i})$. 
Set $\CE_i=\div(f_i')\in \Div(\CU_i)$. 
Then $\CE_i\geq 0$ and $\CE_i\geq \CE|_{\CU_i}$ in $\Div(\CU_i)$. 
By step (1), there is an effective divisor $\CE_i'$ on $\CX$ such that 
$\CE_i'|_{\CU_i}\geq \CE_i$ in $\Div(\CU_i)$.
Then $\CE'=\sum_i \CE_i'$ is an effective divisor on $\CX$ such that 
$\CE'|_{\CU_i}\geq \CE_i'|_{\CU_i}\geq \CE_i\geq \CE|_{\CU_i}$ in $\Div(\CU_i)$, and thus 
$\CE'|_\CU\geq \CE$ in $\Div(\CU)$.
This proves the claim. 

Now we can prove part (1). 
It suffices to construct $\OCD_2$, since $\OCD_1$ can be obtained by considering $-\OCD$.
Write $\OCD=a \OCE$ for $a\in\QQ_{>0}$ and $\OCE\in \wh\Div(\CX)$. 
By the case $\CU=\CX$ in the claim, there is an effective divisor $\CE'\in \Div(\CX)$ with $\CE'\geq \CE$ in  $\Div(\CX)$. 
Taking a suitable Green function for $\CE'$ (in the arithmetic case), 
we get an effective arithmetic divisor $\OCE'\in \wh\Div(\CX)$ with $\OCE'\geq \OCE$ in 
$\wh\Div(\CX)$. 
Then we take $\OCD_2=a' \OCE'$ for an integer $a'\geq a$. 

It remains to prove part (2). 
Again, it suffices to construct $\OCD_2$. 
We will find effective divisors $\OCD_2', \OCD_2''\in \wh\Div(\CX)$ with $\OCD\leq\OCD_2'$ in $\wh\Div(\CX)_\QQ$ and 
$\CD|_\CU\leq\CD_2''|_\CU$ in $\wh\Div(\CU)$.  
Then we can just set $\OCD_2=\OCD_2'+\OCD_2''$.
Note that $\OCD_2'$ is already constructed in (1). 
For $\OCD_2''$, the claim gives its underlying divisor, and it suffices to choose a suitable Green function (in the arithmetic case) to make it effective.
This finishes the proof. 
\end{proof}

\subsection{Model case}

\kkk
Let $\CX$ be a projective variety over $k$.
Let $\CU$ be an open subscheme of $\CX$.
Let $\OCD$ be an arithmetic $(\QQ,\ZZ)$-divisor on $(\CX,\CU)$. 
Definition \ref{def effective sections} gives
$$\wh H^0(\CU, \OCD)=\{f\in k(\CU)^\times: \wh\div(f)+ \OCD\geq 0 
\text{ in } \wh\Div(\CU/k)
\}\cup \{0\}.$$
Note that the partial order is taken in $\wh\Div(\CU/k)_\rmod$.
Here $\wh\div(f)\in \wh\Div(\CX)$ is viewed as an element of $\wh\Div(\CX,\CU)$
via the canonical map $\wh\Div(\CX)\to\wh\Div(\CX,\CU)$.

If $\CX$ is normal, using the rational part of $\OCD$, we have a well-defined $\wh H^0(\CX, \OCD)'$ in the above, which might be different from $\wh H^0(\CU, \OCD)$, as they use different effectivity relations. The following result gives some inequalities between these different notions. 

\begin{lem} \label{finiteness2}
Let $X$ be a flat and essentially quasi-projective integral scheme over $k$. Let $\CU$ be a quasi-projective model of $X$, and let $\CX$ be a projective model of $\CU$.
Let $\OCD$ be an arithmetic $(\QQ,\ZZ)$-divisor on $(\CX,\CU)$. 
\begin{enumerate}
\item
There is a canonical injection
$$
\wh H^0(X, \OCD) \lra \wh H^0(\CX', \pi^*\OCD)'.
$$
Here $\pi:\CX'\to \CX$ is the normalization of $\CX$.
\item
If $\OCD$ is the image of an integral arithmetic divisor $\OCD^*\in \wh\Div(\CX)$ in $\wh\Div(\CX,\CU)$, then there is a canonical injection
$$
\wh H^0(\CX, \OCD^*) \lra \wh H^0(X, \OCD).
$$
\item $\wh h^0(X, \OCD)$ is always finite.
\end{enumerate}
\end{lem} 
\begin{proof}
Part (2) is trivial. 
Part (3) is a direct consequence of (1).
For (1), denoted by $X'$, the generic point of $\CX'$.
The canonical map 
$$\wh H^0(X, \OCD) \lra \wh H^0(X', \pi^*\OCD)$$
is injective, 
so it suffices to prove that the canonical injection 
$$
\wh H^0(\CX', \pi^*\OCD)'\lra 
\wh H^0(X', \pi^*\OCD)
$$
is bijective.
It suffices to note that for any non-empty open subscheme $\CU'$ of $\CX'$, and for any $\OCE\in \wh\Div(\CX',\CU')$, the relations $\OCE\geq 0$, viewed in 
$$\wh\Div(\CX',\CU'),\quad \wh\Div(\CX')_\QQ, \quad \wh\Div(X'/k),$$ 
are all equivalent. This is a consequence of Lemma \ref{effectivity normal}, Lemma \ref{effectivity1}, and Lemma \ref{effectivity2} by converting effectivity to positivity of Green functions. 
\end{proof}

\begin{rmk}
If $X$ is normal, then the injection in (1) is an isomorphism.  
\end{rmk}

\subsection{Adelic case}

Now  we can easily obtain the finiteness of $\wh h^0$ in Definition \ref{def effective sections}(1).

\begin{lem} \label{finiteness3}
\kkk
Let $\OD$ be an adelic divisor on a flat and essentially quasi-projective integral scheme $X$ over $k$. Then the following are true.
\begin{enumerate}
\item
There is a model adelic divisor $\OD'$ on $\CU$, induced by an effective and nef arithmetic divisor on a projective model of $\CU/k$, such that 
$$
-\OD'\leq \OD\leq \OD'
$$
in $\wh\Div(\CU/k)$.
\item $\wh h^0(X, \OCD)$ is always finite.
\end{enumerate}
\end{lem} 
\begin{proof}
Part (1) implies part (2)  by Lemma \ref{finiteness2}(3).
For part (1), assume that $\OD$ is represented by a Cauchy sequence
$\{\OCD_i\}_{i\geq1}$ in $\wh\Div(\CU/k)_\rmod$ for a quasi-projective model $\CU$ of $X$.
The Cauchy condition implies that for some rational number $\epsilon_1>0$,
$$-\epsilon_1 \OCE_0 \leq \OCD_i-\OCD_1 \leq \epsilon_1 \OCE_0, \quad \forall i>1.$$ 
The limit gives
$$-\epsilon_1 \OCE_0 \leq \OD-\OCD_1 \leq \epsilon_1 \OCE_0.$$ 
This gives a model adelic divisors $\OD'$ on $\CU$ such that 
$-\OD'\leq \OD\leq \OD'$. 
Assume that $\OD'$ is defined on a projective model $\CX$ of $X$. 
We can find a nef and effective arithmetic divisor $\OD'$ on $\CX$ such that 
$-\OD''\leq \OD'\leq \OD''$. 
This finishes the proof.
\end{proof}

\section{Volumes of adelic line bundles} \label{sec HS}

The goal of this section is to extend many fundamental properties on volumes of hermitian line bundles to adelic line bundles, including the arithmetic Hilbert--Samuel formula, the arithmetic bigness theorems, the Fujita approximation theorem, the log-concavity theorem, and continuity of the volume function.
The key to these extensions is that volumes of hermitian line bundles naturally approximate volumes of adelic line bundles.

\subsection{Volumes on arithmetic varieties}

We refer to \S\ref{app sec positivity} for an overview of many major theorems on positivity and volumes of hermitian line bundles. In the following, we will still repeat some crucial ones in both the geometric case and the arithmetic case for the purpose here.

\kkk
Let $\CX$ be a projective variety over $k$ of absolute dimension $d$.
For any hermitian line bundles $\CLL$ on $\CX$ (with continuous metrics), denote the volume
$$
\wh\vol(\CX,\CLL)
:=\lim_{m\to \infty} \frac{d!}{m^d}\wh h^0(\CX, m\CLL).
$$
The limit defining the volume always exists.
In the geometric case, this is a result of Fujita (cf. \cite[11.4.7]{Laz2}).
In the arithmetic case, this is originally proved by Chen \cite{Che1}, and Yuan \cite{Yua2} gives a different proof. 
We need the following basic properties of the volume function.

\begin{enumerate} 
\item In the arithmetic case, if there is a sequence $\{\CLL_i\}_{i\geq1}$ of hermitian line bundles on $\CX$ with underlying line bundles $\CL_i= \CL$ such that the metrics of $\CLL_i$ converges to the metric of $\CLL$ uniformly, then $\wh\vol(\CX,\CLL_i)$ converges to $\wh\vol(\CX,\CLL)$. This is a direct consequence of \cite[Prop. 2.1]{YZt}. 
Hence, in many situations, we can easily extend the results from smooth metrics to continuous metrics. 
\item In both cases, the volume function is homogeneous in that 
$\wh\vol(\CX,m\CLL)=m^d\wh\vol(\CX,\CLL)$ for any positive integer $m$.
Therefore, the definition of $\wh\vol$ extends to all hermitian $\QQ$-line bundles by homogeneity. 
\item In both cases, the volume function is a birational invariant. 
Namely, if $\pi:\CX'\to \CX$ is a birational morphism of projective varieties over $k$, then $\wh\vol(\CX',\pi^*\CLL)=\wh\vol(\CX,\CLL)$.
The geometric case is proved in \cite[Prop. 2.2.43]{Laz1}, and the arithmetic case is proved by Moriwaki \cite[Thm. 4.3]{Mor5}.
\end{enumerate}

The arithmetic Hilbert--Samuel formula asserts that, for any nef hermitian line bundles $\CLL$ on $\CX$, 
$$
\wh\vol(\CX,\CLL) = \CLL^d.
$$
Here the right-hand side denotes the arithmetic self-intersection number.
In the geometric case, this is the classical Hilbert--Samuel formula in algebraic geometry (cf. \cite[Cor. 1.4.41]{Laz1}).
Now  we briefly describe the history of the formula for the arithmetic case.
If $\CLL$ is ample in the sense of Zhang \cite{Zha1}. The formula is a consequence of the arithmetic Riemann--Roch theorem of Gillet--Soul\'e \cite{GS2}, an estimate of analytic torsions of Bismut--Vasserot \cite{BV}, and the arithmetic Nakai--Moishezon theorem of  Zhang \cite{Zha1}.
See \cite[Corollary 2.7]{Yua1} for more details of the implications.
The formula was further extended to the nef case with continuous metrics by Moriwaki \cite{Mor5, Mor6}.

We will also need a bigness theorem,  which asserts that for hermitian line bundles $\CLL, \CMM$ on $\CX$ such that $\CLL$ and $\CMM$ are ample, 
$$
\wh \vol(\CX, \CLL-\CMM) \geq  \CLL^d-d\,\CLL^{d-1}\CMM. 
$$
In the geometric case, this is a theorem of Siu \cite{Siu}.
In the arithmetic case, this is the main theorem of Yuan \cite{Yua1}. 
This extends to the case that $\CLL$ and $\CMM$ are nef. Fix an ample hermitian line bundle $\CAA$ on $\CX$. For any positive rational number $\epsilon>0$, apply the result to 
$(\CLL+ \epsilon\CAA,\ \CMM+ \epsilon\CAA)$. Then set $\epsilon\to 0$.

\subsection{Main theorems on volumes}

Now  we are ready to state our generalization of the theorems to adelic line bundles. 

\kkk
Let $X$ be a flat and essentially quasi-projective integral scheme over $k$.
Let $\OL$ be an adelic line bundle on $X$.
Define 
$$
\wh\vol(X,\OL)
:=\lim_{m\to \infty} \frac{d!}{m^{d}}\wh h^0(X, m\OL).
$$
Here $d$ is the absolute dimension of a quasi-projective model of $X$ over $k$.
We will prove that the limit exists.

An adelic line bundle $\OL$ on $X$ is said to be \emph{big} if $\wh\vol(X,\OL)>0$. 
Many results on big hermitian line bundles can be generalized to the current setting.

We will freely use the fact that the volume function is increasing under effectivity. More precisely, 
If $\OL'$ is another adelic line bundles on $X$ with $\wh h^0(X,\OL'-\OL)>0$, then 
$\wh\vol(X,\OL)\leq \wh\vol(X,\OL')$. Taking a nonzero section 
$s\in \wh H^0(X,\OL'-\OL)$, multiplication by $s$ induces an injection $\wh H^0(X,\OL)
\to \wh H^0(X,\OL')$ and thus an inequality $\wh h^0(X,\OL)
\leq \wh h^0(X,\OL')$. 

Our first main result asserts that the limit defining $\wh\vol(X,\OL)$ exists. 

\begin{thm}  \label{limit}
\kkk
Let $X$ be a flat and essentially quasi-projective integral scheme over $k$.
Let $\OL$ be an adelic line bundle on $X$.
\begin{enumerate}
\item The limit 
$$
\wh\vol(X,\OL)
=\lim_{m\to \infty} \frac{d!}{m^{d}}\wh h^0(X, m\OL)
$$
exists. Here $d$ is the absolute dimension of a quasi-projective model of $X$ over $k$. 
\item
If $\OL$ is represented by an adelic line bundle 
$(\CL,(\CX_i,\overline \CL_i, \ell_{i})_{i\geq1})$ on $\CU$ for a quasi-projective model $\CU$ of $X$ over $k$, 
then 
$$
\wh\vol(X,\OL)
=\lim_{i\to \infty} \wh\vol(\CX_i,\CLL_i).
$$
\end{enumerate}
\end{thm}

On the right-hand side, $\wh\vol(\CX_i,\CLL_i)$ is the volume of $\CLL_i$ as a hermitian $\QQ$-line bundle on $\CX_i$, defined by homogeneity.
By the theorem, the definition of $\wh\vol(X,\OL)$ extends to adelic $\QQ$-line bundles on $X$ by homogeneity.

The proof of Theorem \ref{limit} will take up most of the rest of this section. 
Let us first note that by the theorem, the arithmetic Hilbert--Samuel formula and the arithmetic bigness theorem can be generalized to the following theorem. 

\begin{thm}  \label{HS}
\kkk
Let $X$ be a flat and essentially quasi-projective integral scheme over $k$.
Denote by $d$ the absolute dimension of a quasi-projective model of $X$ over $k$. 
\begin{enumerate}
\item
Let $\OL$ be a nef adelic line bundles on $X$. Then
$$
\wh\vol(X,\OL) = \OL^d.
$$ 
\item
Let $\OL,\OM$ be nef adelic line bundles on $X$. Then 
$$
\wh\vol(X,\OL-\OM) \geq \OL^d-d\,\OL^{d-1}\OM.
$$
\end{enumerate}
\end{thm}

It is immediate that Theorem \ref{HS} holds for strongly nef adelic line bundles as a limit version of its counterpart on projective (arithmetic) varieties by Theorem \ref{limit}. 
The theorem will be further extended to nef adelic line bundles by the continuity of the volume function in Theorem \ref{volume continuity} below.

In application, the above theorem is usually combined with the following basic result.  

\begin{prop}  \label{change norm}
\kkk
Let $X$ be a flat and essentially quasi-projective integral scheme over $k$.
Denote by $d$ the absolute dimension of a quasi-projective model of $X$ over $k$. 
Let $\OL$ be an adelic line bundle on $X$. 
\begin{enumerate}
\item
If $k=\ZZ$, let $\ON\in \wh\Pic(\ZZ)$ be a hermitian line bundle with 
$\wh\deg(\ON)>0$. 
Then for any positive rational number $c$, 
$$
\wh\vol(X,\OL-c\ON) \geq \wh\vol(X,\OL)-d\,c\,\wh\deg(\ON)\, \wh\vol(X_\QQ,\TL).
$$
Here $\ON$ is viewed as an adelic line bundle on $X$ via pull-back, and 
$\TL$ is the image of $\OL$ under the canonical map 
$\wh\Pic(X/\ZZ)\to \wh\Pic(X_\QQ/\QQ)$.
\item
If $k$ is a field, assume that there is a projective and regular curve $B$ over $k$ together with a flat $k$-morphism $X\to B$. 
Let $N\in \Pic(B)$ be a line bundle with $\deg(N)>0$. 
Then for any positive rational number $c$, 
$$
\wh\vol(X,\OL-cN) \geq \wh\vol(X,\OL)-d\,c\, \deg(N)\, \wh\vol(X_K,\TL).
$$
Here $N$ is viewed as an adelic line bundle on $X$ via pull-back, $K$ is the function field of $B$, and $\TL$ is the image of $\OL$ under the canonical composition
$$\wh\Pic(X/k)\lra \wh\Pic(X_K/k)\lra \wh\Pic(X_K/K).$$
\end{enumerate}
\end{prop}

\begin{proof}
By Theorem \ref{limit}, the problem is reduced to the case that $X$ is projective over $k$, and $\OL$ is a hermitian line bundle on $X$.
Then the result is more or less well-known, and one easily checks that the result depends only on 
$c\,\wh\deg(\ON)$ (or $c\deg(N)$). 
The arithmetic case is implied by \cite[Prop. 4.2(2)]{Mor5}.
The geometric case can be proved by assuming that $N$ is linearly equivalent to a closed point $P\in B$ and applying the exact sequence
$$
0\lra H^0(X, aL-bN)\lra H^0(X, aL-(b-1)N)\lra H^0(X_P, (aL-(b-1)N)|_{X_P})
$$
to count the dimensions for $a\geq b\geq 1$. 
\end{proof}

\subsection{Volumes of model adelic divisors} 

For the proof of Theorem \ref{limit}, we need a slightly generalized limiting expression about volumes of arithmetic $\QQ$-divisors. 

\kkk
Let $\CX$ be a projective variety over $k$ of absolute dimension $d$.
Let $\OCD$ be an arithmetic $\QQ$-divisor on $\CX$.
Recall that in \S\ref{sec alternative sections} we have introduced  
$$\wh H^0(\CX, \OCD)'=\{f\in k(\CX)^\times: \wh\div(f)+ \OCD\geq 0
\text{ in } \wh\Div(\CX)_\QQ
\}\cup \{0\}.$$
If $k=\ZZ$, we have 
$$\wh h^0(\CX, \OCD)'=\log\#\wh H^0(\CX, \CO(\OCD))';$$ 
if $k$ is a field and $\CX$ is normal, we have 
$$\wh h^0(\CX, \OCD)'=\dim_k H^0(\CX, \CO(\OCD))'.$$

On the other hand, we have an extended definition of  
$$\wh\vol(\CX,\OCD)=\wh\vol(\CX,\CO(\OCD))$$ 
from integral divisors to $\QQ$-divisors by homogeneity. 
Namely, let $a$ be a positive integer such that $a\OCD$ can be realized as an integral arithmetic divisor $\OCD^*$ on $\CX$.
Then
$$\wh\vol(\CX,\OCD):=\frac{1}{a^d}\wh\vol(\CX,\OCD^*)
=\frac{1}{a^d}\lim_{m\to \infty} \frac{d!}{m^d}\wh h^0(\CX, m\OCD^*).
$$ 
It turns out that we have the following compatibility.

\begin{lem} \label{volume rational}
Let $\pi:\CX'\to \CX$ be the normalization morphism. 
Let $\OCD\in \wh\Div(\CX)_\QQ$ and $\OCE\in \wh\Div(\CX')_\QQ$  
be arithmetic $\QQ$-divisors. Then
$$
\wh\vol(\CX,\OCD)
=\wh\vol(\CX',\pi^*\OCD)
=\lim_{m\to \infty} \frac{d!}{m^d}\wh h^0(\CX', m\, \pi^*\OCD+\OCE)'.
$$
\end{lem}

\begin{proof}
The first equality follows from the birational invariance of the geometric case in \cite[Prop. 2.2.43]{Laz1} and the arithmetic case in \cite[Thm. 4.3]{Mor5}. 
The second equality holds for arithmetic $\RR$-divisors by Moriwaki \cite[Thm. 5.2.2(1)]{Mor6} (in the arithmetic case).
For our purpose of arithmetic $\QQ$-divisors, the situation is much easier. We sketch a proof, which will be used later. 

For the second equality, it suffices to assume that $\CX=\CX'$ is normal. 
By Lemma \ref{divisor relations}, there are arithmetic divisors $\OCD_1$ and $\OCD_2$ on $\CX$ with $\OCD_1\leq\OCD\leq\OCD_2$ and $\OCD_1\leq\OCE\leq\OCD_2$
in $\wh\Div(\CX)_\QQ$. 
Let $a$ be a positive integer such that $a\OCD$ can be realized as an integral arithmetic divisor on $\CX$.
For any $r=0,\cdots, a-1$, we have
$$
\wh h^0(\CX, ma\OCD+(r+1)\OCD_1)
\leq  \wh h^0(\CX, (am+r)\OCD+\OCE)'
\leq \wh h^0(\CX, ma\OCD+(r+1)\OCD_2).
$$
On the other hand, we have
$$
\lim_{m\to \infty} \frac{d!}{m^d}\wh h^0(\CX, ma\OCD+(r+1)\OCD_j)
=\wh\vol(\CX,a\OCD), \quad j=1,2.
$$
In the arithmetic case, the extra term $(r+1)\OCD_j$ does not change the limit by \cite[Thm. 4.4]{Mor5}. The corresponding result also holds in the geometric case, and we omit them. 
\end{proof}

Now  we introduce a compatibility result on volumes of model adelic divisors, which will be used in the proof of Theorem \ref{limit}.
The setting is similar to that of Lemma \ref{finiteness2}. 

\begin{lem}  \label{finiteness8}
Let $X$ be a flat and essentially quasi-projective integral scheme over $k$. Let $\CU$ be a quasi-projective model of $X$, and let $\CX$ be a projective model of $\CU$.
Let $\OCD$ be an arithmetic $(\QQ,\ZZ)$-divisor on $(\CX,\CU)$. 
Denote by $d$  the absolute dimension of $\CX$.
Then the limit
$$
\wh\vol(X,\OCD)=\lim_{m\to \infty} \frac{d!}{m^d}\wh h^0(X, m\OCD)
$$
exists and equals $\wh\vol(\CX,\OCD)$.
\end{lem} 
\begin{proof}
By Lemma \ref{divisor relations}, 
there are $\OCD_1, \OCD_2\in \wh\Div(\CX)$
with $\OCD_1\leq\OCD\leq\OCD_2$ in $\wh\Div(\CX,\CU)$. 
By the method of the proof of Lemma \ref{volume rational}, it suffices to prove that for any $\OCD,\OCE\in \wh\Div(\CX)$,
$$
\lim_{m\to \infty} \frac{d!}{m^d}\wh h^0(X, m\OCD+\OCE)
=\wh\vol(\CX,\OCD).
$$
Assume $\OCE=0$ for simplicity of notations since the general case is similar.

Denote by $\pi:\CX'\to \CX$ the normalization. 
By Lemma \ref{finiteness2}(1) and Lemma \ref{volume rational}, 
$$
\limsup_{m\to \infty} \frac{d!}{m^d}\wh h^0(X, m\OCD)
\leq \limsup_{m\to \infty} \frac{d!}{m^d}\wh h^0(\CX',m\,\pi^*\OCD)'
=\wh\vol(\CX,\OCD).
$$
On the other hand, Lemma \ref{finiteness2}(2) implies 
$$
\liminf_{m\to \infty} \frac{d!}{m^d}\wh h^0(X, m\OCD)
\geq \liminf_{m\to \infty} \frac{d!}{m^d}\wh h^0(\CX, m\OCD)
=\wh\vol(\CX,\OCD).
$$
This finishes the proof. 
\end{proof}

\subsection{Proof of Theorem \ref{limit}}

Now we are ready to prove Theorem \ref{limit}.
It is easier to write the proof in terms of divisors, so we reformulate the problem as follows.

Let $\OD$ be an adelic divisor on $X$.
Assume that $\OD$ is represented by $\OCD\in \wh\Div(\CU/k)$ for a quasi-projective model $\CU$ of $X$, and that $\OCD$ is a Cauchy sequence
$\{\OCD_i\}_{i\geq1}$ in $\wh\Div(\CU/k)_\rmod$.
The goal is to prove that the limit 
$$
\lim_{m\to \infty} \frac{d!}{m^{d}}\wh h^0(X, m \OD)
$$
and the limit
$$
\lim_{i\to \infty} \wh\vol(\CX_i, \OCD_i)
$$
exist and are equal. Here we write $\OD$ instead of $\CO(\OD)$ in the notations for $\wh h^0$ and $\wh\vol$ and take similar conventions in the following.

For convenience, denote 
$$
\wh\vol(X, \OD)_-=\liminf_{m\to \infty} \frac{d!}{m^{d}}\wh h^0(X, m \OD),
$$
$$
\wh\vol(X, \OD)_+=\limsup_{m\to \infty} \frac{d!}{m^{d}}\wh h^0(X, m \OD).
$$

The Cauchy condition implies that there is a sequence $\{\epsilon_j\}_{j\geq 1}$ of positive rational numbers converging to $0$ such that 
$$
 \OCD_j -\epsilon_j \OCE_0 \leq
\OCD_i\leq  \OCD_j +\epsilon_j \OCE_0,\quad\ i\geq j\geq 1.
$$
Here $(\CX_0,\OCE_0)$ is a boundary divisor, and we have assumed that there is a morphism $\CX_i\to \CX_0$ extending the identity map on $\CU$ for any $i\geq1$.
The effectivity relations hold in $\wh\Div(\CU/k)_\rmod$, but we can actually assume that it holds in $\wh\Div(\CX_j)_\QQ$ by replacing each $\CX_i$ by a projective model of $\CU$ dominating $\CX_i$.

Set $i\to \infty$. This gives
$$
\OCD_j  -\epsilon_j \OCE_0 \leq
\OCD \leq \OCD_j + \epsilon_j \OCE_0,\quad\  j\geq 1.
$$
The effectivity relations hold in $\wh\Div(\CU/k)$.

Take $\wh h^0$ in the above inequality. 
We have
$$
\wh h^0(X,\OCD_j  -\epsilon_j \OCE_0) \leq
\wh h^0(X,\OCD) \leq \wh h^0(X,\OCD_j + \epsilon_j \OCE_0)
$$
It follows that 
$$
\wh\vol(X, \OD)_- 
\geq 
\liminf_{m\to \infty} \frac{d!}{m^{d}}
 \wh h^0(X,m(\OCD_j  -\epsilon_j \OCE_0))
=
\wh\vol(\CX_i, \OCD_i -  \epsilon_i \OCE_0),
$$
Here equality follows from Lemma \ref{finiteness8}. 
Here by abuse of notations, $\OCE_0$ denotes its pull-back to $\CX_i$. 
Similarly, we have
$$
\wh\vol(X, \OD)_+ \leq 
\limsup_{m\to \infty} \frac{d!}{m^{d}}
 \wh h^0(X,m(\OCD_j  +\epsilon_j \OCE_0))
 =\wh\vol(\CX_i, \OCD_i +  \epsilon_i \OCE_0).
$$
Combining them together, we have
$$
 \wh\vol(\OCD_i -  \epsilon_i \OCE_0)
\leq \wh\vol(X, \OD)_- 
\leq \wh\vol(X, \OD)_+ \leq 
 \wh\vol(\OCD_i +  \epsilon_i \OCE_0).
$$
Here we also omit the dependence on $\CX_i$ of the volume function, noting that the volume function on projective varieties is a birational invariant.

We first consider the case
$$\liminf_{i\to \infty} \wh\vol(\OCD_i + \epsilon_i \OCE_0)=0.$$ 
Then we have $\wh\vol(X, \OD)_+ \leq 0$, and thus $\wh\vol(X,\OCD)=0$. 
Since $\OCD_i\leq \OCD_j + \epsilon_j \OCE_0$ for any $i>j$, 
we have $\ds\lim_{i\to \infty} \wh\vol(\OCD_i)=0$.
It follows that the theorem holds in this case.

Now we consider the major case  
$$
\liminf_{i\to \infty} \wh\vol(\OCD_i + \epsilon_i \OCE_0)>0.
$$
By removing finitely many terms, we assume that 
$\wh\vol(\OCD_i + \epsilon_i \OCE_0)>0$
for every $i$.
It suffices to prove 
$$
\lim_{i\to \infty} \big(\wh\vol(\OCD_i + \epsilon_i \OCE_0)-\wh\vol(\OCD_i -  \epsilon_i \OCE_0)\big)=0.
$$

We are going to apply Fujita's approximation theorem proved in \cite{Fuj}, and its arithmetic counterpart proved independently by Yuan \cite{Yua2} and Chen \cite{Che2}. 
With the current notations, the theorems assert that, for any $\delta>0$, there is a birational morphism $\psi:\CX_i'\to \CX_i$ of projective varieties over $k$, together with an ample arithmetic $\QQ$-divisor $\CFF_i$ on $\CX_i'$ such that 
$$\wh\vol(\CFF_i)>\wh\vol(\OCD_i+ \epsilon_i \OCE_0)-\delta$$ 
and
such that $\psi^*(\OCD_i+ \epsilon_i \OCE_0)-\CFF_i$ is an effective arithmetic $\QQ$-divisor on $\CX_i$. 
Then we have 
$$
\wh\vol(\OCD_i - \epsilon_i \OCE_0)
\geq \wh\vol(\CFF_i - 2\epsilon_i \OCE_0).
$$

Now  we are going to apply Siu's theorem and Yuan's arithmetic version, which we recalled before. 
Write $\OCE_0=\CAA-\CBB$ for nef divisors $\CAA$ and $\CBB$ on $\CX_0$. 
Then
$$
\wh\vol(\CFF_i - 2\epsilon_i \OCE_0)
=\wh\vol(\CFF_i +2\epsilon_i\CBB- 2\epsilon_i \CAA)
\geq (\CFF_i +2\epsilon_i\CBB)^d-2d \epsilon_i(\CFF_i +2\epsilon_i\CBB)^{d-1}  \CAA.
$$

We need a uniform upper bound on $\CFF_i$. 
We claim that there is a nef arithmetic $\QQ$-divisor $\CNN$ on $\CX_1$ 
such that $\CNN\geq \CFF_i$ in $\wh\Div(\CX_i')_\QQ$ for any $i$.
In fact, by the Cauchy condition, 
$$\CFF_i\leq \OCD_i+ \epsilon_i \OCE_0\leq 
\OCD_1 +  \epsilon_i \OCE_0+\epsilon_1 \OCE_0$$
 for any $i$. 
Then it is easy to find $\CNN$ to bound $\OCD_1 +  \epsilon_i \OCE_0+\epsilon_1 \OCE_0$. 
See also Lemma \ref{finiteness3}(1). 

With the uniform bound $\CNN$, we have 
$$
(\CFF_i +2\epsilon_i\CBB)^{d-1}  \CAA
\leq (\CFF_i +2\epsilon_i\CBB)^{d-2} (\CNN +2\epsilon_i\CBB)^{d-1}  \CAA
\leq \cdots
\leq (\CNN +2\epsilon_i\CBB)^{d-1}  \CAA.
$$
It follows that 
$$
\wh\vol(\OCD_i - \epsilon_i \OCE_0)
\geq \CFF_i^d- 2\epsilon_id(\CNN +2\epsilon_i\CBB)^{d-1}  \CAA.
$$
Set $\delta\to 0$, so that $\CFF_i^d\to \wh\vol(\OCD_i+\epsilon_i \OCE_0)$.
The bound becomes 
$$
\wh\vol(\OCD_i - \epsilon_i \OCE_0)
\geq \wh\vol(\OCD_i+\epsilon_i \OCE_0)
- 2\epsilon_id(\CNN +2\epsilon_i\CBB)^{d-1}  \CAA.
$$
As a consequence, we have 
$$
\lim_{i\to \infty} \big(\wh\vol(\OCD_i+\epsilon_i \OCE_0)-\wh\vol(\OCD_i -  \epsilon_i \OCE_0)\big)=0.
$$
This proves the theorem.

\subsection{More properties of the volume function}

Here we generalize some other fundamental properties of volumes of (hermitian) line bundles to adelic line bundles. 
These results for hermitian line bundles and adelic line bundles in the sense of \cite{Zha2} are listed in Theorem \ref{app volume3} and Theorem \ref{app volume4}.
The first result is the log-concavity property.

\begin{thm}[log-convavity] \label{log-concavity}
\kkk
Let $X$ be a flat and essentially quasi-projective integral scheme over $k$.
Let $\OL_1,\OL_2$ be two effective adelic line bundles on $X$. Then 
$$
\wh\vol(\OL_1+\OL_2)^{1/d}
\geq \wh\vol(\OL_1)^{1/d}
+\wh\vol(\OL_2)^{1/d}.
$$
Here $d$ is the absolute dimension of a quasi-projective model of $X$ over $k$.
\end{thm}
\begin{proof}
The result is easy if $\OL_1$ or $\OL_2$ is not big. 
Assume that both $\OL_1$ and $\OL_2$ are big. 
Apply Theorem \ref{limit}. The problem is converted to the model case. 
Then the geometric case is the classical result in \cite[Thm. 11.4.9]{Laz2}, and the arithmetic case is \cite[Thm. B]{Yua2}.
\end{proof}

A morphism between two flat and essentially quasi-projective integral schemes over $k$ is called \emph{birational} if it induces an isomorphism between the function fields.
The next result says that the volume function is also a birational invariant.

\begin{thm}[birational invariance] 
\kkk
Let $\pi: X'\to X$ be a birational morphism of essentially quasi-projective integral schemes over $k$.
Let $\OL$ be an adelic line bundle on $X$. Then 
$$
\wh\vol(X',\pi^*\OL)
= \wh\vol(X,\OL).
$$
\end{thm}
\begin{proof}
This is the adelic version of the geometric case in \cite[Prop. 2.2.43]{Laz1} and the arithmetic case in \cite[Thm. 4.3]{Mor5}.
For the current case, it suffices to check the case that $X'$ is the generic point of $X$. Then the result follows from Theorem \ref{limit}(2). 
\end{proof}

We also have the Fujita approximation theorem for adelic line bundles. 
There are many slightly different notions of ampleness of hermitian line bundles; we take the notion of ``arithmetically positive'' introduced right before Lemma \ref{pairing nef2}. By abuse of terminology, ``arithmetically positive'' on a projective variety over a field means ``ample''.

\begin{thm}[Fujita approximation] \label{fujita}
\kkk
Let $X$ be a flat and essentially quasi-projective integral scheme over $k$.
Let $\OL$ be a big adelic $\QQ$-line bundle on $X$. Then for any $\epsilon>0$,
there exist a flat and essentially quasi-projective integral scheme $X'$ over $k$, a birational projective morphism $\pi:X'\to X$ over $k$, a projective model $\CX'$ of $X'$ over $k$, 
and an arithmetically positive hermitian $\QQ$-line bundle $\CAA$ on $\CX'$ such that $ \wh\vol(\CX',\CAA)\geq\wh\vol(X,\OL)-\epsilon$ and
$\pi^*\OL-\CAA$ is effective in $\wh\Pic(X'/k)_\QQ$. 
\end{thm}
\begin{proof}
Apply Theorem \ref{limit}. The problem is converted to the original Fujita approximation theorem proved in \cite{Fuj}, and its arithmetic counterpart proved independently by Yuan \cite{Yua2} and Chen \cite{Che2}. 
\end{proof}

Now  we consider a continuity property of the volume function 
$\wh\vol:\wh\Pic(X/k)_\QQ\to \RR$. 
Recall that in the projective case, the volume function has very nice continuity properties by Lazarsfeld \cite[Thm. 2.2.44]{Laz1} for the geometric case and by Moriwaki \cite{Mor5} for the arithmetic case. 
The following result generalizes these two, but our proof is different from theirs.

\begin{thm}[continuity] \label{volume continuity}
\kkk
Let $X$ be a flat and essentially quasi-projective integral scheme over $k$.
Let $\OL,\OM_1,\cdots, \OM_r$ be adelic $\QQ$-line bundles on $X$. Then 
$$
\lim_{t_1,\cdots, t_r\to 0}\wh\vol(\OL+t_1\OM_1+\cdots+t_r \OM_r)
=\wh\vol(\OL).
$$
Here $t_1,\cdots, t_r$ are rational numbers converging to $0$.
\end{thm}
\begin{proof}

For convenience, a model adelic line bundle is called \emph{nef} if it is induced by a nef hermitian line bundle on a projective model. 
We will apply Theorem \ref{HS} for nef model adelic line bundles. 
We already know that Theorem \ref{HS} holds for strongly nef adelic line bundles by Theorem \ref{limit}. 

By Lemma \ref{finiteness3}(1), there is a nef and effective model adelic line bundle $\OM_i'$ on $X$ such that $\OM_i'\pm \OM_i$ are effective for any $i$.
Denote $\OM=\OM_1'+\cdots+\OM_r'$, which is still a nef and effective model adelic line bundle.
It suffices to prove 
$$
\lim_{t\to 0}\wh\vol(\OL+t\OM)
=\wh\vol(\OL).
$$

We first treat the case that $\OL$ is big. 
If $t<0$, denote $t'=-t$. 
Apply the Fujita approximation theorem in Theorem \ref{fujita}.
By replacing $X$ by some $X'$ with a birational morphism $X'\to X$ if necessary, we can assume that $\OL\geq \OA$ for some nef model adelic line bundle 
$\OA$ on $X$ with $\wh\vol(\OA)\to \wh\vol(\OL)$. 
Then the bigness result in Theorem \ref{HS}(2) gives
$$
\wh\vol(\OL-t'\OM)
\geq \wh\vol(\OA-t'\OM)
\geq \OA^d-dt'\OA^{d-1}\OM.
$$
Here $d$ is the absolute dimension of a quasi-projective model of $X$ over $k$.

We can bound $\OA^{d-1}\OM$ from above as in the proof of Theorem \ref{limit}. 
In fact, by Lemma \ref{finiteness3}(1) again, we can find a nef model adelic line bundle $\ON$ on $X$ such that $\OL\leq \ON$. 
This implies $\OA\leq \ON$. It follows that 
$$
\OA^{d-1}\OM\leq \ON^{d-1}\OM
$$
is bounded as $\OA$ varies.
Then the above lower bound of $\wh\vol(\OL-t'\OM)$ gives
$$
\liminf_{t'\to 0+}\wh\vol(\OL-t'\OM)
\geq \wh\vol(\OL).
$$
This proves the case $t<0$.

If $t>0$, by the log-concavity theorem in Theorem \ref{log-concavity}, 
$$
\wh\vol(\OL+t\OM)^{1/d}
\leq 2\,\wh\vol(\OL)^{1/d}
-\wh\vol(\OL-t\OM)^{1/d}.
$$
Then the result follows from the case $t<0$.
This idea of applying the log-concavity theorem is inspired by 
 Chen \cite{Che3}. 

Now  we treat the case that $\OL$ is not big. 
Then $\wh\vol(\OL)=0$, and we need to prove 
$\wh\vol(\OL+t\OM)$ converges to $0$. 
Assume that it is not true, and thus there is a constant $c>0$ such that 
$\wh\vol(\OL+t\OM)>c$ for all rational number $t>0$. 
Apply the Fujita approximation theorem in Theorem \ref{fujita} again.
There is a birational morphism $X_t\to X$ such that $\OL+t\OM\geq \OA_t$ on $X_t$ for some nef model adelic line bundle 
$\OA_t$ on $X_t$ with $\wh\vol(\OA_t)>c/2$. 
Then the bigness result in Theorem \ref{HS}(2) gives
$$
\wh\vol(\OL)
\geq \wh\vol(\OA_t-t\OM)
\geq \OA_t^d-dt\OA_t^{d-1}\OM.
$$

We can bound $\OA_t^{d-1}\OM$ by the above method. 
In fact, as above, we have a nef model adelic line bundle $\ON$ on $X$ such that $\OL\leq \ON$. 
This implies $\OA_t\leq \ON+t\OM$. It follows that 
$$
\OA_t^{d-1}\OM\leq (\ON+t\OM)^{d-1}\OM
$$
is bounded as $t\to 0$.
As a consequence, 
$$
\wh\vol(\OL)
\geq c/2-O(t), \quad t\to 0.
$$
This contradicts to $\wh\vol(\OL)=0$.
The proof is complete.
\end{proof}

In the end, we present a basic result that asserts that the bigness of the geometric part is close to the bigness of the whole adelic line bundle. 

\begin{lem}\label{bigness lemma}
\kkk
If $k=\ZZ$, let $K$ be a number field;
if $k$ is a field, let $K$ be a function field of one variable over $k$.
Let $X$ be a quasi-projective variety over $K$. 
Let $\ON\in \wh\Pic(K/k)$ be an adelic line bundle with $\wh\deg(\ON)>0$, viewed as an adelic line bundle on $X/k$ via pull-back. 
Let $\OL\in \wh\Pic(X/k)$ be an adelic line bundle on $X/k$. 
Assume that the image $\wt L$ of $\OL$ under the canonical map $\wh\Pic(X/k)\to \wh\Pic(X/K)$ is big on $X/K$. 
Then for sufficiently large rational number $c$, the adelic line bundle $\OL+c\ON$ is big on $X/k$.  
\end{lem}

\begin{proof}
We only consider the arithmetic case $k=\ZZ$, since the geometric case is similar. 
Let $\CU$ be a quasi-projective model of $X$ over $\ZZ$ such that $\OL$ actually lies in 
$\wh \Pic(\CU/\ZZ)$.
It is more convenient to work on adelic divisors, so we take an adelic divisor
$\OCD\in\wh \Div(\CU/\ZZ)$ linearly equivalent to $\OL$.

We need the general fact that any adelic divisor is represented by an increasing Cauchy sequence, which is proved seriously in  
Lemma \ref{increasing sequence} below. 
Then $\OCD$ is the limit of the increasing sequence
$\{\OCD_i\}_{i\geq 1}$ in $\wh \Div  (\CU/\ZZ)_\rmod$.
Then the image $\wt \CD$ of $\OCD$ in $\wh\Div(X/K)$ 
is the limit of the increasing sequence
$\{\CD_{i,\QQ}\}_{i\geq 1}$ in $\wh \Div  (\CU_\QQ/\QQ)_\rmod$. 
Note that $\wh\vol(\CD_{i,\QQ})$ converges to $\wh\vol(\wt \CD)>0$. 
Then there is an $i$ such that $\wh\vol(\CD_{i,\QQ})>0$. 
Note that 
$$\wh\vol(\OCD+c\ON)=\lim_{j\to\infty}\wh\vol(\OCD_{j}+c\ON) \geq \wh\vol(\OCD_{i}+c\ON)$$ 
by the increasing property of the sequence. 
It suffices to prove that $\OCD_{i}+c\ON$ is big
for sufficiently large rational number $c$, under the condition that $\CD_{i,\QQ}$ is big. 
This reduces the problem from the adelic divisor $\OCD$ to the arithmetic divisor $\OCD_i$.

The arithmetic case is well-known to experts. In fact, by linear equivalence, we can reduce the problem to the case that $\ON$ is represented by (the pull-back of) the arithmetic divisor $(0,1)$ 
on $\Spec \ZZ$ with underlying divisor $0\in \Div(\ZZ)$ and with Green function $1$ at the archimedean place. 
Then the arithmetic case of the result we need is an easy consequence of \cite[Cor. 2.4(1)(4)]{Yua1}.
\end{proof}

In the above proof, we have used the following basic but useful lemma, which will be used again later. 

\begin{lem} \label{increasing sequence}
\kkk
Let $\CDD$ be an adelic divisor on a quasi-projective variety $\CU$ over $k$. 
Then $\OCD$ is represented by a increasing Cauchy sequence, i.e. a Cauchy sequence 
$\{\OCD_i\}_{i\geq 1}$ in $\wh \Div  (\CU/k)_\rmod$
such that $\OCD_{i'}\geq \OCD_i$ for any $i'\geq i$.
Similarly, $\OCD$ is represented by a decreasing Cauchy sequence. 
\end{lem}

\begin{proof}
It suffices to prove the existence of the increasing sequence, since the decreasing sequence can be obtained by taking the increasing sequence for $-\CDD$. 
By definition, $\OCD$ is represented by a Cauchy sequence in $\wh \Div  (\CU/k)_\rmod$, i.e. a sequence $\{\OCD_i\}_{i\geq 1}$ in $\wh \Div  (\CU/k)_\rmod$ satisfying the property that there is a sequence $\{\epsilon_i\}_{i\geq 1}$ of positive rational numbers converging to $0$ such that 
$$
 -\epsilon_i \OCE_0 \leq
\OCD_{i'}-\OCD_{i} \leq  \epsilon_i \OCE_0,\quad\ i'\geq i\geq 1.
$$
Replacing $\{\OCD_i\}_{i\geq 1}$ by a subsequence if necessary, we can assume  
$\epsilon_{i+1}\leq \epsilon_i/2$ for every $i\geq 1$.  
Now $\OCD$ is represented by the Cauchy sequence
 $\{\OCD_i-2\epsilon_i \OCE_0\}_{i\geq 1}$, which is increasing by 
$$
(\OCD_{i+1}-2\epsilon_{i+1} \OCE_0)-(\OCD_i-2\epsilon_i \OCE_0)
\geq  -\epsilon_i \OCE_0-2\epsilon_{i+1} \OCE_0+2\epsilon_i \OCE_0
\geq 0.
$$
\end{proof}

\subsection{Pseudo-effective adelic line bundles}
\label{subsec pseudo}

We have already encountered positivity notions including ampleness, nefness, bigness, and effectivity. In this subsection, we introduce one more positivity notion called pseudo-effectivity, which is the weakest one among all these natural positivity notions. 

Sometimes, it is more convenient to work with adelic $\QQ$-line bundles or adelic $\QQ$-divisors.
Recall that an adelic $\QQ$-line bundles is nef (resp. big, effective) if some positive integer multiple of it is a nef (resp. big, effective) adelic line bundle. A similar statement holds for adelic $\QQ$-divisors.

\kkk
Let $X$ be a flat and essentially quasi-projective integral scheme over $k$. 

\begin{defn} \label{def pseudo}
Let $\OL$ be an adelic line bundle or an adelic $\QQ$-line bundle on $X/k$. 
We say that $\OL$ is \emph{pseudo-effective} on $X/k$ 
if for any big adelic $\QQ$-line bundle $\OM$ on $X/k$, the sum $\OL+\OM$ is a
big  adelic $\QQ$-line bundle on $X/k$.
\end{defn}

By linear equivalence, we have notions of pseudo-effective adelic divisors
and pseudo-effective adelic $\QQ$-divisors.

\begin{lem} \label{pseudo lemma}
Let $\OL$ be an adelic $\QQ$-line bundle  on $X/k$.
Then the following three statements are equivalent:
\begin{enumerate}
\item $\OL$ is {pseudo-effective};
\item there exists an adelic $\QQ$-line bundle $\OM_0$ on $X/k$ such that for every positive rational number $\epsilon$, the sum $\OL+\epsilon\OM_0$ is effective;
\item for any adelic $\QQ$-line bundle $\OM$ on $X/k$, we have
 $\wh\vol(\OL+\OM)\geq \wh\vol(\OM).$ 
\end{enumerate}
Moreover, the following statements are true:
\begin{enumerate}
\item[(a)] if $\OL$ is nef, then $\OL$ is pseudo-effective;
\item[(b)] if $\OL$ is pseudo-effective, then $\OL\cdot \ON_1\cdots \ON_{d-1}\geq 0$ for any nef adelic $\QQ$-line bundles $\ON_1,\cdots, \ON_{d-1}$ on $X/k$.
Here $d$ is the dimension of a projective model of $X$ over $k$.
\end{enumerate}
\end{lem}

\begin{proof}
Note that (1)$\Rightarrow$(2) and (3)$\Rightarrow$(1) are trivial. 
Now we prove (2)$\Rightarrow$(3). 
Let $\OM_0$ and $\OM$ be as in (2) and (3). 
If (2) holds, then for $\epsilon>0$, 
$$
\wh\vol(\OL+\OM) \geq \wh\vol(\OM-\epsilon \OM_0)  \to \wh\vol(\OM),\qquad \epsilon\to 0.
$$
Here the convergence follows from the continuity of the volume function in Theorem \ref{volume continuity}.

For (a), it suffices to choose $\OM_0$ to be nef and big in (2).
For (b), by (2), it is reduced to the case that $\OM$ is effective, which follows from the corresponding statement on projective varieties over $k$. 
\end{proof}

In $\wh\Div(X)_\QQ$, effective adelic $\QQ$-divisors are exactly limits of effective model adelic $\QQ$-divisors under the boundary topology. 
In $\wh\Pic(X)_\QQ$, pseudo-effective adelic $\QQ$-line bundles are exactly limits of effective adelic $\QQ$-line bundles under the ``linear topology''. 
Here on a (possibly infinite-dimensional) vector space $V$ over $\QQ$, a closed subset of $V$ under the linear topology is of the form $S\cap W$, where $W$ is a finite-dimensional $\QQ$-subspace of $V$, and $S$ is a closed subset of $W\otimes_{\QQ}\RR\simeq \RR^{\dim W}$ under the Euclidean topology. 
We can also compare the situation with the definition of nef adelic line bundles in Definition \ref{def nef}.

In algebraic geometry, the BDPP criterion proved by Boucksom--Demailly--Paun--Peternell
\cite{BDPP} (cf. \cite[Thm. 11.4.19]{Laz2})
asserts that the cone of pseudo-effective divisors is dual to the cone of mobile curves, where mobile curves can be realized as push-forward's of intersections of nef divisors under birational morphisms.
The arithmetic analogue of this criterion is proved by Ikoma \cite[Thm. 6.4]{Iko}. 
Now we generalize the criterion to the adelic setting as follows.

\begin{thm}[BDPP criterion] \label{BDPP criterion}
 \kkk
Let $X$ be a flat and essentially quasi-projective integral scheme over $k$. 
Let $d$ be the dimension of a projective model of $X$ over $k$.
Let $\OL$ be an adelic $\QQ$-line bundle on $X/k$.
Then the following statements are equivalent:
\begin{enumerate}
\item $\OL$ is {pseudo-effective};
\item for any proper birational morphism $f:X'\to X$, where $X'$ is a flat and essentially quasi-projective integral scheme over $k$, and for any 
nef adelic $\QQ$-line bundle $\ON$ on $X'/k$, the intersection number 
 $f^*\OL\cdot \ON^{d-1}\geq 0$.  
 \end{enumerate}
\end{thm}

\begin{proof}
Note that (1)$\Rightarrow$(2) follows from Lemma \ref{pseudo lemma}(b). 
It suffices to prove (2)$\Rightarrow$(1). We can assume that $X=\CU$ is quasi-projective over $k$, and that $\OL$ is linearly equivalent to an adelic divisor $\OCD$ on $\CU/k$. 
By Lemma \ref{increasing sequence},  $\OCD$ 
 is the limit of a decreasing sequence
$\{\OCD_i\}_{i\geq 1}$ in $\wh \Div  (\CU/k)_\rmod$. 
We have $\OCD_i\geq \OCD$, as the limit of $\OCD_i\geq \OCD_{i'}$ by $i'\to \infty$. 
Therefore, (2) implies 
 $$f^*\CDD_i\cdot \ON^{d-1}\geq f^*\CDD\cdot \ON^{d-1}\geq 0.$$  
Note that $\CDD_i$ is a model divisor. 
By the theorems of \cite{BDPP} and \cite{Iko}, $\CDD_i$ is pseudo-effective. 
This implies that $\CDD$ is pseudo-effective by Lemma \ref{pseudo lemma}(3).
\end{proof}

\section{Heights on quasi-projective varieties}
In this section, we introduce heights on quasi-projective varieties over finitely generated fields.
If the varieties are projective, we can define vector-valued heights, which refines the Moriwaki height in \cite{Mor3, Mor4} in the arithmetic case. 
Note that Moret-Bailly introduced vector-valued heights \cite{MB} for different purposes.
The Northcott property in the projective case is deduced from that of Moriwaki, and the fundamental inequality is extended to the current case following an idea of Moriwaki.

\subsection{Vector-valued heights}

\kkk
\ccc
Let $F$ be a finitely generated field over $k$ and $X$ be a quasi-projective variety over $F$.
Let $\overline L$ be an element of $\wh \Pic (X/k)_{\intb,\QQ}$, i.e. an integrable adelic $\QQ$-line bundle.

For any point $x\in X(\overline F)$, define the \emph{vector-valued height}
$$
\fh_{\overline L}(x)
:=\fh_{\overline L}(x')
:= \frac{ \left\langle \overline L|_{x'}\right\rangle }
{\deg(x') }\ \in\, \wh \Pic(F/k)_{\intb,\QQ}.$$
Here $x'$ denotes the closed point of $X$ containing $x$,
$\deg(x')$ is the degree of the residue field of $x'$ over $F$,
$\overline L|_{x'}$ denotes the \emph{pull-back} of $\OL$ in 
$\wh\Pic(x'/k)_{\intb,\QQ}$,
and $\langle\overline L|_{x'}\rangle$ is the image of $\overline L|_{x'}$
under the norm map $\wh\Pic(x'/k)_{\intb,\QQ} \to \wh\Pic(F/k)_{\intb,\QQ}$, which is the  Deligne pairing of relative dimension $0$ in Theorem \ref{intersection2}.

Therefore, we have a height function 
$$
\fh_{\overline L}: X(\overline F) \longrightarrow \wh \Pic(F/k)_{\intb,\QQ}.
$$
Note that we do not require $X$ to be projective here.

We can generalize the definition to high-dimensional projective subvarieties.
Let $Z$ be a closed projective $\overline F$-subvariety of $X$, i.e. 
a closed subvariety of $X_{\overline F}$ which is \emph{projective} over $\overline F$. 
Define \emph{the vector-valued height of $Z$ for $\overline L$} as
$$
\fh_{\overline L}(Z)
:=\fh_{\overline L}(Z')
:= \frac{ \left\langle (\overline L|_{Z'})^{\dim Z+1}\right\rangle }
{(\dim Z+1)\deg_{L}(Z') }\ \in \wh \Pic(F/k)_{\intb,\QQ}.$$
Here $Z'$ denotes the image of $Z\to X$,
$\overline L|_{Z'}$ denotes the \emph{pull-back} in $\wh\Pic(Z'/k)_\intb$,
and the self-intersection is as in Theorem \ref{intersection2}.
As we do not require $X$ to be projective or $L$ to be ample, the height $\fh_{\overline L}(Z)$ is only well-defined if $Z$ is projective and $\deg_{L}(Z')\neq 0$. 

The following are some special situations:
\begin{enumerate}
\item If $\OL$ is nef on $X$, then the height $\fh_{\overline L}(Z)$ is also nef if it is defined.
\item If $X$ is projective over $F$ and $L$ is ample on $X$, then the degree $\deg_{L}(Z')$ is well-defined and positive for all closed subvarieties of $X_{\overline F}$. This gives a function 
$$
\fh_{\overline L}:  |X_{\overline F}| \longrightarrow \wh \Pic(F/k)_{\intb,\QQ}.
$$
Here $|X_{\overline F}| $ denotes the set of closed subvarieties of $X_{\overline F}$.
\item
If $k=\ZZ$, let $F$ be a number field;
if $k$ is a field, let $F$ be a function field of one variable over $k$.
There is a degree map 
$$\wh\deg: \wh \Pic(F/k)_{\QQ}\longrightarrow  \BR.$$ 
This follows from Lemma \ref{arithmetic curve} in the number field case and limits of degrees of divisors on curves in the function field case. In both cases,
$$h_\OL(Z):=\wh\deg\,\overline\fh_{\overline L}(Z)$$ 
generalizes the height function of Zhang \cite{Zha2} (from the projective case to the quasi-projective case).

\qquad If furthermore $X$ is projective over $F$, by the case $\dim Z=0$, we have a height function $h_\OL:X(\overline F)\to \RR$, which is a Weil height associated to $L\in \Pic(X)_\QQ$.  
\end{enumerate}

\subsection{High-dimensional base}

The above definition works well because any scheme over $\Spec F$ is automatically flat over $\Spec F$, which is required in the Deligne pairing. This still works well if we change $\Spec F$ to a Dedekind scheme, but if we change $\Spec F$ to a high-dimensional base, we easily lose this convenience, and thus, the definition only works in some special cases.

\kkk
Let $B$ be a flat and essentially quasi-projective integral scheme over $k$, and $X$ be a flat and quasi-projective integral scheme over $B$.
Let $\overline L$ be an element of $\wh \Pic (X/k)_{\intb,\QQ}$.
Then we have a vector-valued height function 
$$
\fh_{\overline L}: X(B) \longrightarrow \wh \Pic(B/k)_{\intb,\QQ}, 
\quad x\longmapsto x^*\OL.
$$

In general, for any integral subscheme $Y$ of $X$ which is projective and flat over $B$,  we can still define a vector-valued height of $Y$ for $\overline L$ 
in terms of the Deligne pairing. We omit it here since we will not use it in the current book.

\subsection{Moriwaki heights}

Let $(k, F,X,\overline L)$ be as above. Namely, 
$k$ is either $\ZZ$ or a field,
$F$ is a field finitely generated over $k$, $X$ is a quasi-projective variety over $F$, and $\overline L$ is an element of $\wh \Pic (X/k)_{\intb,\QQ}$.
Denote by $d$ the absolute dimension of a quasi-projective model of $F$ over $k$,
and denote by $n$ the dimension of $X$.

Let $\overline H_1,\cdots, \overline H_{d-1}$ be any $d-1$ elements in $\wh \Pic(F/k)_{\intb,\QQ}$. For any point $x\in X(\overline F)$, define the \emph{Moriwaki height of $x$ for $\overline L$ and 
$(\overline H_1,\cdots, \overline H_{d-1})$} by
$$
h_{\overline L}^{\overline H_1,\cdots, \overline H_{d-1}}(x):
=h_{\overline L}^{\overline H_1,\cdots, \overline H_{d-1}}(x'):
=\frac{ \overline L|_{x'} \cdot \overline H_1\cdots \overline H_{d-1}}
{\deg(x') }
\, \in\, \RR.$$
Here $x'$ and $\overline L|_{x'}$ are as above, and the intersection  number is taken in $\wh\Pic(x'/k)_{\intb,\QQ}$,  
as defined by Proposition \ref{intersection1}, where $\overline H_1,\cdots, \overline H_{d-1}$ are viewed as elements of $\wh\Pic(x'/k)_{\intb,\QQ}$ via pull-back.
This gives a height function 
$$
h_{\overline L}^{\overline H_1,\cdots, \overline H_{d-1}}: X(\overline F) \longrightarrow \RR.
$$

We can also generalize the definition to high dimensions. 
For any closed $\overline F$-subvariety $Z$ of $X$,  \emph{the Moriwaki height of $Z$ for $\overline L$ and 
$(\overline H_1,\cdots, \overline H_{d-1})$} is
$$h_{\overline L}^{\overline H_1,\cdots, \overline H_{d-1}}(Z)
:=h_{\overline L}^{\overline H_1,\cdots, \overline H_{d-1}}(Z')
:=\frac{  (\overline L|_{Z'})^{\dim Z+1} \cdot \overline H_1\cdots \overline H_{d-1}}
{(\dim Z+1)\deg_{\wt L}(Z'/F) }.$$
Here 
$Z'$ and $\overline L|_{Z'}$ are as above, and the intersection  number is taken in $\wh\Pic(Z'/k)_{\intb,\QQ}$,  
as defined by Proposition \ref{intersection1}, where $\overline H_1,\cdots, \overline H_{d-1}$ are viewed as elements of 
$\wh\Pic(Z'/k)_{\intb,\QQ}$ via pull-back.
The term $\wt L$ denotes the image of $\OL$ under the canonical map
$$
\wh\Pic(X/k)_{\intb,\QQ} \lra \wh\Pic(X/F)_{\intb,\QQ}
$$
introduced at the end of \S\ref{sec functoriality}, and 
$$\deg_{\wt L}(Z'/F):=(\wt L|_{Z'})^{\dim Z}\,\in\,\RR$$ 
is the self-intersection number of $\wt L|_{Z'}$ in $\wh\Pic(Z'/F)_{\intb,\QQ}$ defined in Proposition \ref{intersection1}.

The height $h_{\overline L}^{\overline H_1,\cdots, \overline H_{d-1}}(Z)$ is only well-defined if $\deg_{\wt L}(Z'/F)\neq 0$.
If $Z$ is projective, then we have $\wt L|_{Z'}=L|_{Z'}$, and thus $\deg_{\wt L}(Z'/F)=\deg_{L}(Z')$.

The vector-valued height refines the Moriwaki height by the simple formula
$$
h_{\overline L}^{\overline H_1,\cdots, \overline H_{d-1}}(Z)
=\fh_{\overline L}(Z)\cdot \overline H_1\cdots \overline H_{d-1}
$$
as long as the right-hand side is well-defined. 
We introduce the following simplified notations. 
\begin{enumerate}
\item
For any $\overline H\in \wh \Pic(F/k)_{\intb,\QQ}$, denote  
$h_{\overline L}^{\overline H}=h_{\overline L}^{\overline H,\cdots, \overline H}$, where the right-hand has $d-1$ copies of $\OH$. 

\item
If $F$ is a number field, then $d=1$ and thus $h_{\overline L}^{\overline H}$ is independent of $\OH$, so we just write $h_{\overline L}=h_{\overline L}^{\overline H}$. 
In this case, we simply have 
$$h_{\overline L}(Z)
=h_{\overline L}(Z')
=\frac{  (\overline L|_{Z'})^{\dim Z+1} }{(\dim Z+1)\deg_{\wt L}(Z'/F)}
\,\in\,\RR.$$
A similar convention holds for function fields of one variable. 
\end{enumerate}

In the arithmetic case ($k=\ZZ$), if $X$ is projective over $F$, and both $\overline L$ and $(\overline H_1,\cdots, \overline H_{d-1})$ are realized on some projective model
$\CX\to \CS$ of $X\to \Spec F$, then $h_{\overline L}^{\overline H_1,\cdots, \overline H_{d-1}}$ is exactly the height function introduced in \cite{Mor3}. In \cite{Mor4}, Moriwaki generalizes the definition to the case that an adelic sequence gives $\overline L$. 

Let us briefly compare our adelic line bundles with the adelic sequence in \cite{Mor4}. 
Roughly speaking, the adelic sequence in the loc. cit. are more numerical since they use intersection numbers to define their topology, while our adelic line bundles use effectivity to define their topology. 
Then our notion includes more restrictive objects and allows coarser equivalence relations. 
These two notions are similar in the definition of absolute intersection numbers, but our notion has the advantage of having Deligne pairings, effective sections, and volumes.

\subsection{Northcott property in the projective case}

In the projective case, we have the following Northcott property of Moriwaki Heights, which generalizes \cite[Prop. 3.3.7(4)]{Mor3}. 

\begin{thm}[Northcott property]
Let $k$ be either $\ZZ$ or a finite field. 
Let $F$ be a finitely generated field over $k$, and let $d$ be the absolute dimension of a quasi-projective model of $F$ over $k$.
Let $X$ be a projective variety over $F$.
Let $\overline L$  be an element in $\wh \Pic (X/k)_{\intb,\QQ}$  with an ample generic fiber $L$. 
Let $\overline H_1,\cdots, \overline H_{d-1}$ be nef and big elements in 
$\wh \Pic (F/k)_{\intb,\QQ}$.
Then for any
$D\in \BR$ and $A\in \BR$,
the set
$$
\{ x\in X(\overline F):\ \deg(x)<D,\ h_{\overline L}^{\overline H_1,\cdots, \overline H_{d-1}}(x)<A \}
$$
is finite. 
\end{thm}

\begin{proof}
We only treat the arithmetic case $k=\ZZ$, since the geometric case is similar. 
If $\overline L, \overline H_1,\cdots, \overline H_{d-1}$ are model adelic divisors, this follows from \cite[Prop. 3.3.7(4)]{Mor3}. 
We will extend it to the current generality by replacing $\overline L, \overline H_1,\cdots, \overline H_{d-1}$ successively by more general adelic line bundles. 

Let $\OL'$ be any integrable adelic line bundle on $X$ with underlying line bundle $L'= L$. 
To replace $\OL$ by $\OL'$, it suffices to check that 
$h_{\overline L}^{\overline H_1,\cdots, \overline H_{d-1}}-h_{\overline L'}^{\overline H_1,\cdots, \overline H_{d-1}}$ is a bounded function on $X(\overline F)$.
Assume that $\OL$ and $\OL'$ lie in $\wh\Picc(\CU/\ZZ)_{\QQ}$ for a quasi-projective model $\CU$ of $X$ over $\ZZ$ with a projective and flat morphism $f:\CU\to \CV$ to a quasi-projective model $\CV$ of $F$ over $\ZZ$.
Let $(\CY_0,\OCE_0)$ be a boundary divisor for $\CV$.
By $L'= L$,  the difference $\OL-\OL'$ is represented by a Cauchy sequence $\OCD=\{\OCD_i\}_{i\geq1}$ in $\wh\Div(\CU)_{\rmod, \QQ}$ with $\CD_i|_X=0$. 
As in the proof of Lemma \ref{finiteness3}, the Cauchy condition implies  
$$\OCD_1-\epsilon f^*\OCE_0\leq \OCD \leq \OCD_1+\epsilon f^*\OCE_0$$
for some positive rational number $\epsilon$. 
Thus $-f^*\OCD_1\leq \OCD\leq f^*\OCD_1$ 
for some $\OCD_1\in \wh\Div(\CV)_{\rmod, \QQ}$.
It follows that 
$$
|h_{\overline L}^{\overline H_1,\cdots, \overline H_{d-1}}-h_{\overline L'}^{\overline H_1,\cdots, \overline H_{d-1}}|
\leq \OCD_1\cdot \overline H_1 \cdots \overline H_{d-1}.
$$

Now we replace $\overline H_1,\cdots, \overline H_{d-1}$ successively by more general line bundles. 
By symmetry, it suffices to do that for $\OH_1$. 
Let $\OH_1'$ be a nef and big element in $\wh \Pic (F/\ZZ)_{\intb,\QQ}$.
To replace $\OH_1$ by $\OH_1'$, it suffices to have 
$h_{\overline L}^{\overline H_1',\cdots, \overline H_{d-1}}\geq c\, h_{\overline L}^{\overline H_1,\cdots, \overline H_{d-1}}$ 
for some positive rational number $c$.
We can further assume that $\OL$ is nef. 
Then it suffices to prove that $\OH_1'-c\,\OH_1$ is effective for some positive rational number $c$.
By Theorem \ref{HS}, 
$$
\wh\vol(\OH_1'-c\,\OH_1)\geq \OH_1^d-d\,c\, \OH_1^{d-1}\OH_1'
$$
is positive if $c$ is sufficiently small. 
This finishes the proof.
\end{proof}

\subsection{Fundamental inequality for Moriwaki heights}

The fundamental inequality, a part of the theorem of successive minima of Zhang \cite{Zha1,Zha2} as reviewed in Theorem \ref{app sm2}, is generalized to projective varieties over finitely generated fields by Moriwaki \cite{Mor3}. Now  we further generalize the result to quasi-projective varieties over finitely generated fields.

We first introduce the Moriwaki condition on polarizations of finitely generated fields. 
\kkk
\ccc
Let $F$ be a finitely generated field over $k$.
Denote by $d$ the absolute dimension of a quasi-projective model of $F$ over $k$.
Let $\overline H\in \wh \Pic(F/k)_{\QQ}$ be an adelic $\QQ$-line bundle. 

If $k=\ZZ$ and $d>1$, we say that $\overline H$ satisfies \emph{the Moriwaki condition} if $\OH$ is nef on $F/\ZZ$, the arithmetic top self-intersection number 
$\wh\deg_{\OH}(F/\ZZ)=\overline H^{d}=0$, and the geometric top self-intersection number 
$\deg_{\wt H}(F/\QQ)=\wt H^{d-1}>0$. 
Here $\wt H$ is the image of $\OH$ under the canonical map 
$$\wh\Pic(F/\ZZ)_{\QQ} \lra \wh\Pic(F/\QQ)_{\QQ}.$$

If $k$ is a field  and $d>1$, we say that $\overline H$ satisfies \emph{the Moriwaki condition} if $\OH$ is nef on $F/k$, the geometric top self-intersection number 
$\deg_{\OH}(F/k)=\overline H^{d}=0$, and the geometric top self-intersection number 
$\deg_{\wt H}(F/K)=\wt H^{d-1}>0$ for some extension $K$ of $k$ in $F$ of transcendental degree 1, where $\wt H$ is the image of $\OH$ under the canonical map 
$$\wh\Pic(F/k)_{\QQ} \lra \wh\Pic(F/K)_{\QQ}.$$
Note that the definition depends on the choice of $K$.

We will assume that $d\geq 1$ and take the convention that the Moriwaki condition is automatically satisfied if $d=1$. 
Now  we are ready to state the theorem. 

\begin{thm}[fundamental inequality] \label{minima}
\kkk
Let $F$ be a finitely generated field over $k$.
Assume that $F$ is an infinite extension of $k$ if $k$ is a field. 
Let $\overline H$ be an element of $\wh\Pic (F/k)_{\QQ}$ satisfying the Moriwaki condition. 
Let $X$ be a quasi-projective variety over $F$.
Let $\overline L$  be a nef element in $\wh \Pic (X/k)_{\QQ}$ such that its image $\wt L$
in $\wh \Pic (X/F)_{\QQ}$ is big. 
Then
$$e_1^{\overline H}(X,\overline L) \geq h_{\overline L}^{\overline H}(X).$$
\end{thm}


Recall that the essential minimum
$$
e_1^{\overline H}(X,\overline L) = \sup_{U\subset X} \inf_{x\in U(\overline F)}  h_{\overline L}^{\overline H}(x),
$$
where the supremum is taken over all Zariski open subschemes $U$ of $X$. 
Recall that $\wt L$ denotes the image of $\OL$ under the map
$$\wh\Pic(X/k)_{\QQ} \lra \wh\Pic(X/F)_{\QQ}.$$

The theorem is a part of \cite[Cor. 5.2]{Mor3} if $k=\ZZ$, $X$ is projective over $F$, and $\OL, \OH$ are model adelic line bundles. 
The general case here is proved similarly, while the new ingredient is our results on volumes of adelic line bundles in Theorem \ref{HS}.

\begin{proof}[Proof of Theorem \ref{minima}]

Let $\ON\in \wh\Pic(B/k)_\QQ$ be an element of degree 1. 
Here is $k=\ZZ$, then $K=\QQ$ and $B=\Spec\ZZ$; 
if $k$ is a field, then $K$ is the function field of one variable over $k$ in $F$ defining the Moriwaki condition, and $B$ is the unique projective regular curve over $k$ with function field $K$. 

View $\ON$ as elements of $\wh\Pic(F/k)_\QQ$ and $\wh\Pic(X/k)_\QQ$ by pull-back.
Denote $\OL'=\OL-c\ON$ with $c\in\QQ$. 
Note that 
$$
e_1^{\overline H}(X,\overline L) 
-e_1^{\overline H}(X,\overline L') 
= h_{\overline L}^{\overline H}(X)
- h_{\overline L'}^{\overline H}(X)
=c\, \wt H^{d-1}.
$$
Thus it suffices to prove that for any $c\in\QQ$ such that $h_{\overline L'}^{\overline H}(X)> 0$, we also have $e_1^{\overline H}(X,\overline L') \geq 0$. 

By the assumption $\overline H^{d}=0$, we see that both $e_1^{\overline H}(X,\overline L')$ and $h_{\overline L'}^{\overline H}(X)$ remain the same if we replace 
$\OL'$ by $\OL'+m\OH$ for some positive rational number $m$. 
Note that we always have 
$$
\wh\vol(\OL'+m\OH) 
\geq (\OL'+m\OH)^{d+n}.
$$
This follows from Theorem \ref{HS}(1) if $c\leq 0$, and follows 
from Theorem \ref{HS}(1) and Proposition \ref{change norm} if $c\geq 0$.

By the assumption $\OH^d=0$, 
$$
(\OL'+m\OH)^{d+n}={d+n \choose d-1} \OL'^{n+1}\OH^{d-1} m^{d-1}+O(m^{d-2}),\quad m\to \infty.
$$
Therefore, there is a positive integer $m$ such that $\wh\vol(\OL'+m\OH)>0$. 
By definition, there is a positive integer $N>0$ such that $N(\OL'+m\OH)$ is an integral adelic line bundle with an effective section $s$ on $X$.
For any point $x\in X(\overline F)$ outside the support $|\div(s)|$, we have 
$$
h_{\OL'}^\OH(x)
=h_{\OL'+m\OH}^\OH(x)
=\frac{1}{N\deg(x)}\wh\div(s)|_{x'} \cdot \pi^*\OH^{d-1}
\geq 0.
$$
Here the intersection is on the closed point $x'\in X$ corresponding to $x$.
This finishes the proof. 
\end{proof}

\subsection{Fundamental inequality over global fields}

Here we provide a different proof for Theorem \ref{minima} for the special but rather important case that $F$ is a global field. We can see more clearly the role of the small sections from this proof. 
Moreover, we also include an inequality in the opposite direction, which is a weak version of Zhang's theorem of successive minima (cf. Theorem \ref{app sm2}) in the current setting.

\begin{thm}[fundamental inequality: global field case] \label{minima2}
\kkk
If $k=\ZZ$, let $K$ be a number field;
if $k$ is a field, let $K$ be a function field of one variable over $k$.
Let $X$ be a quasi-projective variety of dimension $n$ over $K$.
Let $\overline L$  be a nef element in $\wh \Pic (X/k)_\QQ$ such that its image $\wt L$ in 
$\wh \Pic (X/K)_\QQ$ is big.  Then
$$e_1(X,\overline L) \geq h_{\overline L}(X) 
\geq \frac{1}{n+1} e_1(X,\overline L).$$
\end{thm}

The first inequality is a consequence of the arithmetic Hilbert--Samuel formula in Theorem \ref{HS} by the following result.

\begin{lem}  \label{minima3}
Let $X/K/k$ be as in Theorem \ref{minima2}.
Let $\overline L$  be an element of $\wh \Pic (X/k)$, and let
$\wt L$  be its image in $\wh \Pic (X/K)$. 
\begin{enumerate}
\item
If $k=\ZZ$,  for any positive integer $m$ such that $\wh h^0(X,m\OL)>0$,
$$
e_1(X,\overline L) \geq 
\frac{\wh h^0(X,m\OL)}{m\, \wh h^0(X,m\TL)}-\frac2m[K:\QQ]
$$
if the right-hand side is strictly positive.

\item
If $k$ is a field,  for any positive integer $m$ such that $\wh h^0(X,m\OL)>0$,
$$
e_1(X,\overline L) \geq 
\frac{\wh h^0(X,m\OL)}{m\, \wh h^0(X,m\TL)}-\frac1m [k':k]
$$
if the right-hand side is strictly positive.
Here $k'$ denotes the algebraic closure of $k$ in $K$. 

\item
In both cases,
$$
e_1(X,\overline L) \geq 
\frac{\wh\vol(\OL)}{(n+1) \wh\vol(\TL)}
$$
if both $\wh\vol(\OL)$ and $\wh\vol(\TL)$ are strictly positive.
\end{enumerate}
\end{lem}

\begin{proof}
Note that (3) is the limit of (1) and (2) as $m\to\infty$.

Let us first prove (1).  
Recall that $\wh H^0(X,m\TL)$ is a vector space of dimension $\wh h^0(X,m\TL)$ over $K$, which contains the finite set $\wh H^0(X,m\OL)$.
By Definition \ref{def effective sections}, we have a $v$-adic norm $\|\cdot\|_{v,\sup}$ on $\wh H^0(X,m\TL)$ for any place $v$ of $\QQ$. 
We claim that there is a nonzero element $s\in \wh H^0(X,m\OL)$ such that 
$$
-\log \|s\|_\infty \geq \frac{\wh h^0(X,m\OL)}{[K:\QQ]\wh h^0(X,m\TL)}- 2.
$$

The claim follows from a basic result in the classical geometry of lattices.
For example, we can apply \cite[Prop. 2.1(1)]{YZt}.
To match the notations, denote by $M$ the $\ZZ$-submodule of 
$\wh H^0(X,m\TL)$ generated by $\wh H^0(X,m\OL)$.
Then $\OM=(M, \|\cdot\|_{\infty,\sup})$ is a normed $\ZZ$-module in the sense of the loc. cit..
Denote $r=\rank\, M$, which is at most $[K:\QQ]\wh h^0(X,m\TL)$.
Note that 
$$\alpha
:= \frac1r \wh h^0(\OM)- 2
\geq \frac{\wh h^0(X,m\OL)}{[K:\QQ]\wh h^0(X,m\TL)}- 2
>0.
$$
By the first inequality of \cite[Prop. 2.1]{YZt}, 
we have 
$$\wh h^0(\OM(-\alpha))>
\wh h^0(\OM)-r\alpha-r\log 3
>0.$$ 
Then there is a nonzero element $s\in \wh H^0(\OM(-\alpha))$ satisfying 
$-\log \|s\|_{\infty,\sup}\geq \alpha$.  
This proves the claim. 

With the section $s\in \wh H^0(X,m\OL)$, for any $x\in X(\OK)$ not contained in $\div_X(s)$, we have
$$
m\,h_\OL(x)=\frac{1}{\deg_K(x)} \wh\deg(m\OL|_{x'})
=\frac{1}{\deg_K(x)} \sum_v \sum_{y\in x'_v} (-\log \|s(y)\|_v^{\deg_{\QQ_v}(y)}).
$$
The first summation is over all places $v$ of $\QQ$, but the point $x\in X(\OK)$, so there are a lot of Galois orbits in the above.
Namely, $x'$ is the closed point of $X$ corresponding to $x$, 
$x'_v$ is the image of $x'\times_\QQ \QQ_v$ in $X_{\QQ_v}$, 
which is a finite set of closed points of $X_{\QQ_v}$. 
Any $y\in x'_v$ is also viewed as a classical point of $X_v^\an$. 
Then we have 
$$m\,h_\OL(x)\geq [K:\QQ]\alpha,$$ 
since $\|s\|_v\leq 1$ for any finite $v$ and $-\log \|s\|_{\infty,\sup}\geq \alpha$.
It follows that 
$$
e_1(X,\overline L) \geq \frac1m[K:\QQ]\alpha.
$$
This proves (1). 

The proof of (2) is similar to that of (1). 
By replacing $k$ by $k'$, we can assume that $k$ is algebraically closed in $K$.
Denote by $B$ the unique projective regular curve over $k$ with function field $K$.
Take a rational point $v_0\in B(k)$, which exists by passing to a finite extension of $k$. 
Similar to (1), it suffices to prove
that there is a nonzero element $s\in \wh H^0(X,m\OL)$ such that 
$$
-\log \|s\|_{v_0} \geq \frac{\wh h^0(X,m\OL)}{\wh h^0(X,m\TL)}- 1.
$$

In fact, the finite-dimensional $K$-space $\wh H^0(X,m\TL)$ can be viewed as a vector bundle on the generic point $\eta\in B$, and thus a quasi-coherent sheaf on $B$ via push-forward. 
Denote by $\CM$ the $\CO_B$-submodule of 
$\wh H^0(X,m\TL)$ generated by $\wh H^0(X,m\OL)$.
Then $\CM$ is a vector bundle on $B$. Denote $r=\rank\, \CM$.
As the geometric version of \cite[Prop. 2.1]{YZt}, for any positive integer $\alpha$, we have 
$$
h^0(B,\CM(-\alpha v_0)) \geq h^0(B,\CM)-r \alpha.
$$
To see the truth of this inequality, it suffices to consider the case $\alpha=1$, 
but then it follows from the exact sequence 
$$
0\lra H^0(B,\CM(- v_0)) \lra H^0(B,\CM) \lra H^0(v_0,\CM|_{v_0}). 
$$
Once we have the geometric inequality, by taking $\alpha$ to be the integral part of 
$\wh H^0(X,m\OL)/r$, 
we conclude that there is a nonzero element $s\in \wh H^0(X,m\OL)$ such that 
$$
-\log \|s\|_{v_0} \geq \alpha \geq \frac{\wh h^0(X,m\OL)}{\wh h^0(X,m\TL)}- 1.
$$
The remaining part is similar to (1).
\end{proof}

Now we prove the second inequality of Theorem \ref{minima2}. 

\begin{proof}[Proof of Theorem \ref{minima2}]
We have already proved the first inequality of the theorem. 
Now we prove the second inequality of the theorem, i.e. 
$$ h_{\overline L}(X) 
\geq \frac{1}{n+1} e_1(X,\overline L).$$
Assume that $\OL$ is linearly equivalent to an adelic divisor 
$\CDD$ on a quasi-projective model $\CU$ of $X$ over $O_K$.
The proof is done by a series of reductions. For convenience, we write the proof in terms of adelic $\QQ$-divisors. 

First, we reduce to the case that $\CDD$ is strongly nef on $\CU$. 
In fact, there is a strongly nef adelic divisor $\CAA$ on $\CU$ such that $\CDD+\epsilon \CAA$ is strongly nef for any rational number $\epsilon>0$. 
If the inequality holds for all  $\CDD+\epsilon \CAA$, then it holds for $\CDD$ by taking $\epsilon\to 0$, since 
$e_1(X,\CDD+\epsilon \CAA) \geq e_1(X,\CDD)$. 
We assume that $\CDD$ is strongly nef in the following.

Second, we reduce to the case that $\CDD^{n+1}>0$.
This follows a similar process as in the first step by taking $\CAA$ to the be pull-back of some $\CNN\in \wh\Div(B/k)$ with $\wh\deg(\CNN)>0$. 
Here $B=\Spec \ZZ$ if $k=\ZZ$; $B$ is the unique projective regular curve over $k$ with function field $K$ if $k$ is a field.
We assume that $\CDD^{n+1}>0$ in the following.

Third, we reduce to the case that $\CU$ is projective over $k$.
In fact, since $\CDD$ is strongly nef, it is a limit of nef model divisors $\CDD_i$ on $\CU$.
Then there is a sequence $\epsilon_i$ of positive rational numbers converging to $0$ such that 
$$
\OCD\leq \OCD_{i} +  \epsilon_i \OCE_0 \leq \OCD_{i} +  \epsilon_i \OCE.
$$
Here $(\CX_0, \CEE_0)$ is a boundary divisor of $\CU$, 
and $\CEE$ is a fixed nef arithmetic divisor on $\CX_0$ with $\OCE\geq \OCE_0$. 
Let $t>0$ be a positive rational number to be determined later. 
By Theorem \ref{HS}, 
$$
\wh\vol(t\CDD-\CEE) \geq t^{n+1} \CDD^{n+1}-(n+1)t^n\,\CDD^n\CEE.
$$
As $\CDD^{n+1}>0$, there exists $t>0$ such that $\wh\vol(t\CDD-\CEE)>0$.
It follows that 
$$
\CEE\leq_\mathrm{lin} t\CDD, \qquad
\OCD\leq_\mathrm{lin} \OCD_{i} +  \epsilon_i t \OCD, \qquad
(1-\epsilon_i t)\OCD\leq_\mathrm{lin} \OCD_{i}.
$$
Here $\leq_\mathrm{lin}$ denotes the relation in $\wh\CaCl(\CU/k)_\QQ$ (instead of in 
$\wh\Div(\CU/k)_\QQ$) induced by effectivity. 
Now if the inequality holds for $\OCD_{i}$, then 
$$
 h_{\CDD_i}(X) 
 \geq \frac{1}{n+1} e_1(X,  \OCD_{i}) \geq 
 \frac{1 -  \epsilon_i t}{n+1}e_1(X, \OCD).
$$
This implies the inequality for $\CDD$ by taking $i\to \infty$. 

Finally, we prove the case that $\CU=\CX$ is projective over $k$, that $\CLL$ is nef on 
$\CX$, and that $\CL_K$ is big on $\CX_K$.  
By the method in the first step again, we can assume that $\CL_K$ is ample on $\CX_K$.  
In the arithmetic case, 
 the inequality follows from Zhang's original theorem of successive minima in Theorem \ref{app sm1}, by using the fact $e_1(X, \CLL)=e_1(\CX_K, \CLL)$ and the fact 
 $e_i(\CX_K, \CLL)\geq 0$ for nef $\CLL$.
In the geometric case, Zhang's theorem was proved by 
Gubler \cite[Lem. 4.1, Prop. 4.3]{Gub07b}.
This finishes the proof.
\end{proof}

\subsection{Essential minimum and pseudo-effectivity}

We still consider the situation of global fields. 
 \kkk
If $k=\ZZ$, let $K$ be a number field;
if $k$ is a field, let $K$ be a function field of one variable over $k$.

Recall that we have introduced the notion of pseudo-effective line bundles in 
\S\ref{subsec pseudo}. 
It turns out that this notion is closely related to positivity of the essential minimum. 
In fact, we have the following quasi-projective version of Yuan's conjecture stated in \cite[Conj. 1.3]{QY}. 

\begin{conj}[pseudo-effectivity] \label{pseudo conj}
Let $X$ be a quasi-projective variety over $K$.
Let $\OL$ be an adelic line bundle on $X$. 
Then $\OL$ is pseudo-effective if and only if $e_1(X,\OL)\geq 0$. 
\end{conj}

This conjecture looks surprising at the beginning, but ``most of it'' has been proved. In fact, if $X$ is projective over $K$, 
Balla\"y \cite[Thm. 1.1]{Bal} proved the ``if'' part assuming that the metrics of $\OL$ at archimedean places are semipositive, and proved the ``only if'' part assuming that $L$ is big.
Yuan's conjecture is inspired by Balla\"y's work and some ideas from Diophantine geometry. 
More recently, Qu--Yin applied their arithmetic Demailly approximation to remove Balla\"y's assumption on the semipositivity of metrics, and deduced the ``if'' part of the full conjecture; see \cite[Thm. 1.4, Thm. 1.7, Thm. 1.8]{QY}. 
Both works on the ``if'' part are based on the BDPP criterion and the arithmetic version of Ikoma \cite[Thm. 6.4]{Iko} (cf. Theorem \ref{BDPP criterion}).

The following theorem re-organizes the results of Balla\"y \cite{Bal} and Qu--Yin \cite{QY} in our current setting. 

\begin{thm} \label{pseudo thm}
Let $X$ be a quasi-projective variety of over $K$.
Let $\OL$ be an adelic line bundle on $X$. 
Then the following are true. 
\begin{enumerate}
\item If $e_1(X,\OL)\geq 0$, then $\OL$ is pseudo-effective. 
\item If $\OL$ is pseudo-effective, and if $\wt L$ is big, then $e_1(X,\OL)\geq 0$. 
Here $\wt L$ denotes the image of $\OL$ in $\wh \Pic (X/K)$. 
\end{enumerate}
\end{thm}

\begin{proof}
The arithmetic case of (1) is just \cite[Thm. 1.8]{QY}, and the geometric case is similar. Alternatively, it is easy to reduce it to the projective case in \cite[Thm. 1.4, Thm. 1.7]{QY}
by applying Lemma \ref{increasing sequence} to get a decreasing sequence of model 
adelic line bundles converging to $\OL$. 

For (2), we apply the method of \cite{Bal} directly to the quasi-projective case.
In fact, since $\wt L$ is big, by Lemma \ref{bigness lemma}, 
there is an adelic line bundle $\ON\in \wh\Pic(K/k)$ such that
$\OL+\ON$ is big on $X/k$.  
Since $\OL$ is pseudo-effective, for any rational number $\epsilon>0$, 
$\OL+\epsilon (\OL+\ON)$ is big. 
It follows that $e_1(\OL+\epsilon (\OL+\ON), X)> 0$ by Lemma \ref{minima3}.  
Then we have 
$$
0\leq e_1(\OL+\epsilon (\OL+\ON), X)=e_1((1+\epsilon)\OL, X)+\epsilon\, \wh\deg(\ON)
=(1+\epsilon) e_1(\OL, X)+\epsilon\, \wh\deg(\ON).
$$
Setting $\epsilon\to 0$, we have $e_1(\OL, X)\geq0$.
\end{proof}

\subsection{The height inequality}

In the end, we present the following height inequality, which is a general form of the height inequality in Theorem \ref{height comparison}. It holds over finitely generated fields in a suitable sense, but we restrict it to global fields for simplicity. 

\begin{thm}[height inequality]\label{height comparison0}
\kkk
If $k=\ZZ$, let $K$ be a number field;
if $k$ is a field, let $K$ be a function field of one variable over $k$.
Let $\pi:X\to S$ be a morphism of quasi-projective varieties over $K$. 
Let $\OL\in \wh\Pic(X/k)$ and $\OM\in \wh\Pic(S/k)$ be adelic line bundles. 
Denote by $\TL$ the image of $\OL$ in $\wh\Pic(X/K)$.
\begin{enumerate}
\item
If $\OL$ is big on $X$, then there exist $\epsilon>0$ and a non-empty open subvariety $U$ of $X$ such that 
$$
h_{\OL}(x) \geq \epsilon\, h_{\OM}(\pi(x)), \quad\ \forall\, x\in U(\overline K). 
$$
\item If $\OL$ is nef on $X/k$, and $\TL$ is big on $X/K$, then 
 for any $c>0$, there exist $\epsilon>0$ and a non-empty open subvariety $U$ of $X$ such that 
$$
h_{\OL}(x) \geq \epsilon\, h_{\OM}(\pi(x)) -c, \quad\ \forall\, x\in U(\overline K). 
$$
\item If $\TL$ is big on $X/K$, then 
there exist $c>0$,  $\epsilon>0$, and a non-empty open subvariety $U$ of $X$ such that 
$$
h_{\OL}(x) \geq \epsilon\, h_{\OM}(\pi(x)) -c, \quad\ \forall\, x\in U(\overline K). 
$$
\end{enumerate}
\end{thm}

\begin{proof}
We only write the proofs in the arithmetic case $k=\ZZ$, since the geometric case is similar.

We first see that (1) implies (2).
In fact, if $(\OL,\TL)$ is as in (2), 
let $\ON\in \wh\Pic(O_K)$ be a hermitian line bundle with $\wh\deg(\ON)=1$, and view $\ON$ as an adelic line bundle on $X$ by pull-back.  
Denote $\OL'=\OL+c\ON$ for a rational number $c>0$. 
It follows that 
$$
\OL'^{d}=(\OL+c\ON)^d=\OL^d+dc \TL^{d-1}>0.
$$
Here $d=\dim X+1$. 
Then $\OL'$ is big, and we can apply (1) to $(\OL',\OM)$.
This gives (2) by the 
simple relation 
$$
h_{\OL'}(x)
=h_{\OL}(x)+c.
$$

Now we see that (1) implies (3).
In fact, if $(\OL,\TL)$ is as in (3), we still denote
$\OL'=\OL+c\ON$ for a rational number $c>0$.  
By Lemma \ref{bigness lemma}, 
 $\OL'$ is big on $X/\ZZ$ for sufficiently large rational number $c>0$. 
We still apply (1) to $(\OL',\OM)$.

Now  we prove (1). The key is that there exists a rational number $\epsilon>0$ such that $\wh\vol(\OL-\epsilon \pi^*\OM)>0$. 
If $\OL$ and $\OM$ are nef, this is a consequence of Theorem \ref{HS}, which asserts
$$
\wh\vol(\OL-\epsilon \pi^*\OM) \geq 
\OL^d-d\epsilon\, \OL^{d-1}\cdot \pi^*\OM.
$$
In general, by the continuity of the volume function in Theorem \ref{volume continuity}, 
$$
\lim_{\epsilon\to 0} \wh\vol(\OL-\epsilon \pi^*\OM) = \wh\vol(\OL).
$$
Then there still exists such an $\epsilon$.

Consequently, there is a positive integer $m$ and a nonzero effective section $s$ of 
$m(\OL-\epsilon \pi^*\OM)$ on $X$. 
We claim that this implies 
$$
h_{\OL-\epsilon \pi^*\OM}(x)\geq 0,\quad
\forall x\in X(\overline K),\ s(x)\neq 0.
$$
Then the result is followed by the simple relation.
$$
h_{\OL-\epsilon \pi^*\OM}(x)
=h_{\OM}(x)
-\epsilon h_{\OM}(\pi(x)).
$$

For the claim, the reason is already in the proof of Lemma \ref{minima3}.
Alternatively, denote by $x'$ the closed point of $X$ corresponding to $x$, 
Then the pull-back of $\OL-\epsilon \pi^*\OM$ to $x'$ gives an adelic line bundle on $x'$ with a nonzero effective section and thus is linearly equivalent to an effective adelic divisor on $x'$. 
This effective adelic divisor can be written as a limit of effective model divisors on a quasi-projective model of $x'$, and thus the degree is non-negative.
\end{proof}

In the following, we want to prove a partial converse to 
Theorem \ref{height comparison0}(3).
\begin{thm}\label{thm-ht-big} 
\kkk
If $k=\ZZ$, let $K$ be a number field;
if $k$ is a field, let $K$ be a function field of one variable over $k$.
Let $X$ be a quasi-projective variety over $K$. 
Let $\OL, \OM\in \wh\Pic(X/k)$ be adelic line bundles, and let $\wt L, \wt M\in \wh\Pic (X/K)$ be their images in $\wh\Pic (X/K)$.
Assume that $\wt M$ is big. Then $\wt L$ is big if and only if  
there exist $c>0$,  $\epsilon>0$, and a non-empty open subvariety $U$ of $X$ such that 
$$
h_{\OL}(x) \geq \epsilon\, h_{\OM}(x) -c, \quad\ \forall\, x\in U(\overline K). 
$$
\end{thm}

\begin{proof} 
The ``only if" part follows from Theorem \ref{height comparison0}(3). 
We will prove the ``if" part. 
Let $\ON\in \wh\Pic(K/k)_\QQ$ be an adelic $\QQ$-line bundle with $\wh\deg(\ON)=1$, and view $\ON$ as an adelic $\QQ$-line bundle on $X/k$ via pull-back.  
We have 
$$
h_{\OL}(x) \geq \epsilon\, h_{\OM}(x) -c
\quad \Longleftrightarrow \quad
h_{\OL-\epsilon \OM+c\ON}(x) \geq 0.
$$
Thus the ``if'' condition implies 
the essential minimum $e_1(\OL-\epsilon \OM+c\ON,X)\geq 0$. 
By Theorem \ref{pseudo thm}(1), which is due to Qu--Yin \cite{QY},  
$\OL-\epsilon \OM+c\ON$ is pseudo-effective on $X/k$.
By Lemma \ref{pseudo lemma}(2),  pseudo-effectivity is preserved by the map 
$\wh\Pic(X/k) \to  \wh\Pic(X/K)$. 
Apply this to $\OL-\epsilon \OM+c\ON$. 
We see that 
$\wt L-\epsilon \wt M$ is pseudo-effective on $X/K$.
Then $\wt L=(\wt L-\epsilon \wt M)+\epsilon \wt M$ is big on $X/K$.
\end{proof}

\section{Equidistribution: conjectures and theorems}\label{sec equi}

In this section, we formulate an equidistribution conjecture and prove two equidistribution theorems for small points. More precisely, we have the following:
\begin{enumerate}
\item Theorem \ref{equi3}, an equidistribution theorem for quasi-projective varieties over number fields or function fields of one variable;
\item Conjecture \ref{equi1}, an equidistribution conjecture for quasi-projective varieties over finitely generated fields;
\item Theorem \ref{equi5}, an equidistribution theorem for morphisms between quasi-projective varieties over number fields or function fields of one variable. 
\end{enumerate}
All these statements generalize the equidistribution theorems of Szpiro--Ullmo--Zhang \cite{SUZ}, Chambert-Loir \cite{CL}, and Yuan \cite{Yua1} for projective varieties over number fields. We refer to Theorem \ref{app adelic equidistribution} for the equidistribution theorem of Yuan \cite{Yua1}. 
Theorem \ref{equi5} also generalizes an equidistribution theorem of Moriwaki \cite{Mor3}.
Conjecture \ref{equi1} generalizes Theorem \ref{equi3} by changing the base fields; 
Theorem \ref{equi5} generalizes Theorem \ref{equi3} by changing it to a relative version. 

Our main ingredient is the extension of the arithmetic Hilbert--Samuel formula and Yuan's bigness theorem to quasi-projective varieties in Theorem \ref{HS}, with which we can apply the variational principle of Szpiro--Ullmo--Zhang to the current quasi-projective situation.

\subsection{Small points}

We will first state the equidistribution conjecture (Conjecture \ref{equi1}) and then prove the two equidistribution theorems. 
We start with some definitions.

\kkk
Let $X$ be a quasi-projective variety of dimension $n$ over a finitely generated field $F$ over $k$. Let $d$ be the dimension of any quasi-projective model of $F$ over $k$.
Let $\overline L$ be a nef adelic line bundle on $X$.   
Recall that we have a Moriwaki height function
$$
h_{\overline L}^{\overline H}:  X(\overline F)\lra \RR_{\geq0}
$$
for any polarization $\OH\in \wh\Pic(F/k)_{\nef}$.

Denote by $\wt L$ the image of $\OL$ under the canonical map
$\wh\Pic(X/k)_{\nef} \to \wh\Pic(X/F)_{\nef}$
introduced in \S\ref{sec functoriality}. 
Assume that the self-intersection number (defined in Proposition \ref{intersection1})
$$\deg_{\wt L}(X/F)=\wt L^{\dim X}>0.$$  
Then we have a well-defined Moriwaki height  
$$h_{\overline L}^{\overline H}(X)
=\frac{ \overline L^{n+1} \cdot \overline H^{d-1}}
{(n+1)\deg_{\wt L}(X/F) }.$$

A sequence $\{x_m\}_{m\geq 1}$ in $X(\overline F)$
is said to be \emph{generic} if any closed subvariety $Y\subsetneqq X$ contains only finitely many terms of the sequence.

Let $\overline H\in \Pic (F/k)_{\QQ,\nef}$ be a polarization.
A sequence $\{x_m\}$ in $X(\overline F)$ is said to be {\em directionally small for $\OH$}, or just \emph{$h_{\overline L}^{\OH}$-small}, if $h_{\overline L}^{\OH} (x_m)$ converges to $h_{\overline L}^{\OH} (X)$.

A sequence $\{x_m\}$ in $X(\overline F)$ is said to be {\em small}, or just 
\emph{$\fh_{\overline L}$-small}, if it is {$h_{\overline L}^{\OH}$-small} for any polarization $\overline H\in \Pic (F/k)_{\QQ,\nef}.$

If $X$ is projective, so that $\fh_{\overline L}(X)$ is well-defined. 
Then the sequence is {$\fh_{\overline L}$-small} if and only if  
$\fh_{\overline L}(x_m)$ converges to $\fh_{\overline L}(X)$ \emph{numerically} in $\Pic (F/k)_\mathrm{int,\QQ}$ in the sense that
$$\lim_{m\to\infty}\fh_{\overline L}(x_m)\cdot \overline H_1\cdot\overline H_2\cdots \overline H_{d-1}
= \fh_{\overline L}(X)\cdot \overline H_1\cdot\overline H_2\cdots \overline H_{d-1}$$
for any $\overline H_1, \cdots, \overline H_{d-1}\in \Pic (F/k)_{\intb,\QQ}.$

If $k=\ZZ$ and $F$ is a number field, both $\wh\deg\,\fh_{\overline  L}$ and $h_{\overline  L}^{\overline  H}$ are equal to the usual height function $h_{\overline  L}$. Then both smallness is equivalent to the usual one given by $h_{\overline  L}(x_m)\to h_{\overline  L}(X)$.

\subsection{Equilibrium measure}

Resume the above notations for $(k,F,X,\OL)$.

Let $v\in \CM(F/k)$ be a point corresponding to a non-trivial valuation of $F$; 
i.e. a non-trivial multiplicative norm $|\cdot|_v:F\to \BR$. It could be either Archimedean or non-archimedean.
Denote by $F_v$ the completion of $F$ for $v$. 
Then we have the Berkovich space $X_v^\an$ 
associated to the variety $X_{F_v}$ over the complete field $F_v$.

There is an \emph{equilibrium measure} 
$$d\mu_{\overline L,v}:=\frac{1}{\deg_{\wt L}(X/F)}c_1(\OL)_v^n$$
over the analytic space $X_v^\an$.
If $v$ is Archimedean, this is classical. 
If $v$ is non-archimedean, this is defined in terms of the Chambert-Loir measure developed by Chambert-Loir and Ducros in \cite{CLD}. 
We refer to \S\ref{sec app global} for the precise definition.

\subsection{Equidistribution conjecture over finitely generated fields}

For each point $x\in X(\overline F)$, we have the measure
$$\mu_{x,v}:=\frac{1}{\deg(x)}\delta_{x'_v}$$ on $X_v^\an$.
Here $\delta_{x'_v}$ is the Dirac measure for the Galois orbit $x'_v$ of $x$ in $X_v^\an$. More precisely, $x'$ is the closed point of $X$ corresponding to $x$, 
and $x'_v$ is the image of $x'\times_F F_v$ in $X_{F_v}$, 
viewed as a finite set of classical points of $X_v^\an$.

We say the Galois orbit of a sequence $\{x_m\}_{m\geq 1}$ of points of $X(\overline F)$ is \emph{equidistributed in $X_v^\an$ for} $d\mu_{\overline L,v}$ if the weak convergence
$$\mu_{x_m,v}\lra d\mu_{\overline L,v}$$
holds on $X_v^\an$.
Namely, 
$$\int_{X_v^\an} f\, \mu_{x_m,v}\lra \int_{X_v^\an} f\, d\mu_{\overline L,v}$$
for any $f\in C_c(X_v^\an)$. Here $C_c(X_v^\an)$ is the space of real-valued continuous and compactly supported functions on $X_v^\an$.

Finally, we are ready to state our equidistribution conjecture. Recall that $\wt L$ denotes the image of $\OL$ under the map
$\wh\Pic(X/k)_{\nef} \to \wh\Pic(X/F)_{\nef}$.

\begin{conj}[equidistribution over finitely generated fields] \label{equi1}
\kkk
Let $F$ be a finitely generated field over $k$. 
Let $v$ be a non-trivial valuation of $F$.
Assume that the restriction of $v$ to $k$ is trivial if $k$ is a field.
Let $X$ be a quasi-projective variety over $F$.
Let $\overline L$ be a nef adelic line bundle on $X/k$ such that $\deg_{\wt L}(X/F)>0$. 
Let $\{x_m\}_m$ be a generic sequence of small points in $X(\overline F)$.
Then the Galois orbit of $\{x_m\}_m$ is equidistributed in $X_v^\an$ for $d\mu_{\overline L,v}$.
\end{conj}

In the arithmetic case ($k=\ZZ$), if $F$ is a number field and $X$ is projective, the conjecture is fully known previously. 
The pioneering work of Szpiro--Ullmo--Zhang \cite{SUZ} proved the equidistribution for number fields $F$ and archimedean places $v$ assuming pointwise positivity of the Chern form $c_1( L, \|\cdot\|_v)$.
Their work was extended to non-archimedean places $v$ by Chambert-Loir \cite{CL}. 
Yuan proved the full case of number fields with $L$ ample \cite{Yua1} by developing a bigness theorem for the difference of ample hermitian line bundles.
The proof of \cite{Yua1} works by replacing the ampleness of $L$ by the positivity $\deg_{L}(X)>0$. For more history of this subject, we refer to \cite[\S6.3]{Yua3}. 

The above arguments were also generalized to the geometric case. 
In that case, if $X$ is projective over $F$, and the valuation $v$ of $F$ comes from a prime divisor of a projective model of $F$ over $k$, the conjecture was proved independently by Faber \cite{Fab} when the transcendental degree of $F/k$ is one and by Gubler \cite{Gub2} for general transcendental degrees.

In Theorem \ref{equi3} below, we will prove the conjecture for any quasi-projective $X$ and for any number field $F$ or function field $F$ of one variable. 
However, the conjecture seems widely open if $F$ has a positive transcendental degree over $\ZZ$.

\begin{rmk}
In the conjecture, we have assumed that $v$ is a non-trivial valuation of $F$. 
Nonetheless, if $v$ is the trivial valuation of $F$, a similar equidistribution theorem on $X_v^\an$ was proved by \cite[Cor. 5.6]{Xie}. Here the equilibrium measure on $X_v^\an$ is the Dirac measure supported at the point corresponding to the trivial valuation of the function field of $X_{F_v}$. 
\end{rmk}

\subsection{Equidistribution theorem over number fields}

The goal here is to prove the following theorem, which asserts that Conjecture \ref{equi1} holds if $F$ is a number field or a function field of one variable.
It is a consequence of the variational principle of \cite{SUZ, Yua1} and the bigness result in Theorem \ref{HS}. 

\begin{thm}[equidistribution over number fields] \label{equi3}
\kkk
Let $K$ be a number field if $k=\ZZ$; let $K$ be the function field of one variable over $k$ if $k$ is a field. 
Let $X$ be a quasi-projective variety over $K$.
Let $\overline L$ be a nef adelic line bundle on $X/k$ such that $\deg_{\wt L}(X/K)>0$. 
Let $\{x_m\}_m$ be a generic sequence in $X(\overline K)$ such that $\{h_{\OL}(x_m)\}_m$ converges to $h_{\OL}(X)$.
Then the Galois orbit of $\{x_m\}_m$ is equidistributed in $X_v^\an$ for $d\mu_{\overline L,v}$ for any place $v$ of $K$.
\end{thm}

\begin{proof}
We only write the proof for the arithmetic case $k=\ZZ$, since the geometric case is similar. 
Apply the variational principle of \cite{SUZ, Yua1} to Theorem \ref{HS}.
The process is standard at the beginning, and then there will be a new situation due to quasi-projectivity. 

The conditions and the result do not change if replacing $\OL$ by $\OL+\pi^*\ON$ for an element $\ON\in \Pichat(K)_\intb$ with $\wh\deg(\ON)>0$. 
Here $\pi^*: \Pichat(K)_\intb\to \Pichat(X)_\intb$ is the pull-back map. 
As a consequence, we can assume $\OL^{n+1}>0$.
Here we denote $n=\dim X$. 

Let $\OM$ be an element in the kernel of the map
$\wh\Pic(X)_{\intb}\to \wh\Pic(X/K)_{\intb}$. 
Let $\epsilon$ be a nonzero rational number.
By Lemma \ref{minima3}, 
$$
e_1(X,\overline L+\epsilon\OM) \geq 
\frac{\wh\vol(\OL+\epsilon \OM)}{(n+1) \wh\vol(\TL)},
$$
if both $\wh\vol(\OL+\epsilon \OM)$ and $\wh\vol(\TL)$ are strictly positive.

Now  it is straightforward to apply Theorem \ref{HS}. 
In fact, by writing $\OM$ as the difference of two nef adelic line bundles, 
Theorem \ref{HS} implies
$$
\wh\vol(\OL+\epsilon \OM) \geq 
\OL^{n+1}+\epsilon (n+1)\OL^n\OM+O(\epsilon^2).
$$
By the assumption $\OL^{n+1}>0$, the right-hand side is strictly positive if $|\epsilon|$ is sufficiently small.
Theorem \ref{HS} also implies the geometric volume
$$
\wh\vol(\TL) = \TL^n=\deg_\TL(X)
$$
which is assumed to be strictly positive.
It follows that 
$$
e_1(X,\overline L+\epsilon\OM)
 \geq 
\frac{\OL^{n+1}+\epsilon (n+1)\OL^n\OM}{(n+1) \deg_\TL(X)}+O(\epsilon^2).
$$

Apply the inequality to the generic sequence $\{x_m\}_m$.
We have
$$
\liminf_{m\to \infty} h_{\OL+\epsilon \OM}(x_m) \geq 
\frac{\OL^{n+1}+\epsilon (n+1)\OL^n\OM}{(n+1) \deg_\TL(X)}+O(\epsilon^2).
$$
By assumption, 
$$
\lim_{m\to \infty} h_{\OL}(x_m) = h_\OL(X)=
\frac{\OL^{n+1}}{(n+1) \deg_\TL(X)}.
$$
Then the inequality implies 
$$
\liminf_{m\to \infty} \epsilon h_{\OM}(x_m) \geq 
\epsilon \frac{\OL^n\OM}{\deg_\TL(X)}+O(\epsilon^2).
$$

If $\epsilon>0$, the above implies 
$$
\liminf_{m\to \infty} h_{\OM}(x_m) \geq 
 \frac{\OL^n\OM}{\deg_\TL(X)}+O(\epsilon).
$$
If $\epsilon<0$, the above implies 
$$
\limsup_{m\to \infty} h_{\OM}(x_m) \leq 
 \frac{\OL^n\OM}{\deg_\TL(X)}+O(|\epsilon|).
$$
Set $\epsilon\to 0$ in each case. We obtain
$$
\lim_{m\to \infty} h_{\OM}(x_m) =  \frac{\OL^n\OM}{\deg_\TL(X)}.
$$

We are going to deduce the equidistribution theorem on $X_v^\an$ from the above limit identity.
Assume that $\OL\in \wh\Pic(\CU)_{\QQ,\nef}$ for a quasi-projective model $\CU$ of $X$ over $\ZZ$, and assume that $\OL$ is represented by a Cauchy sequence 
$(\CL,(\CX_i,\overline \CL_i, \ell_{i})_{i\geq1})$ in $\wh\Picc(\CU)_{\rmod,\QQ}$.
Here $\CX_i$ is a projective model of $\CU$, and $\CLL_i$ is a hermitian $\QQ$-line bundle on $\CX_i$. 
Assume that there is a morphism $\psi_i:\CX_i\to \CX_1$ extending the identity morphism of $\CU$.
Denote $X_i=\CX_{i,\QQ}$, which contains $X$ as an open subvariety.
 
Let $\CX_1'$ be another projective model of $X_1$ over $\ZZ$.
Let $\CMM$ be a hermitian $\QQ$-line bundle on $\CX_1'$, with a fixed isomorphism $\CM_\QQ\to \CO_{X_1}$. Then it induces a metric $\|\cdot \|_w$ of $\CO_{X_1}$ on 
$X_{1,w}^\an$ for any place $w$ of $K$. 
Assume that the metric $\| 1\|_w=1$ for any place $w\neq v$ of $K$.
Denote $f=-\log\|1\|_v$, which is continuous on $X_{1,v}^\an$. 
By definition, 
$$
h_{\CMM}(x_m) = \int_{X_v^\an} f \mu_{x_m,v},
$$
and
$$
\OL^n\CMM 
= \lim_{i\to\infty}  \CLL_i^n\CMM
= \lim_{i\to\infty} 
\int_{X_{i,v}^\an} f c_1(\CLL_i)_v^n.
$$
Then the above result gives a limit identity
$$
\lim_{m\to\infty}
\int_{X_v^\an} f \mu_{x_m,v}
= \frac{1}{\deg_\TL(X)} \lim_{i\to\infty} 
\int_{X_{i,v}^\an} f\, c_1(\CLL_i)_v^n.
$$
Here $f$ is viewed as a function on $X_{i,v}^\an$
by the pull-back induced by $\psi_{i,\QQ}:X_i\to X_1$.

Now we going to vary $f=-\log\|1\|_v$, which is a model function on $X_{1,v}^\an$ associated to $(\CX_1',\CMM)$.
By Gubler's density theorem (cf. \cite[Thm. 7.12]{Gub3} and \cite[Lem. 3.5]{Yua1}), 
the space of all such model functions $f$ is dense in 
$C(X_{1,v}^\an)$ under the topology of uniform convergence. 
Note that  
$$
\lim_{i\to\infty} \int_{X_{i,v}^\an}  c_1(\CLL_i)_v^n
= \lim_{i\to\infty} (\CL_{i,\QQ})^n
= \TL^n= \deg_\TL(X).
$$
Therefore, the limit identity also holds for any $f\in C(X_{1,v}^\an)$.

Finally, assume $f\in C_c(X_{v}^\an)$, viewed as an element of $C(X_{i,v}^\an)$ by the open immersion $X\to X_i$. Then
$$
\lim_{i\to\infty} 
\int_{X_{i,v}^\an} f c_1(\CLL_i)_v^n
=\lim_{i\to\infty} 
\int_{X_{v}^\an} f\, c_1(\CLL_i)_v^n|_{X_v^\an}
=\int_{X_{v}^\an} f c_1(\OL)_v^n.
$$
Here the last equality follows from the definition of $c_1(\OL)_v^n$ in \S\ref{sec local theory}, based on the theory of \cite{CLD} in the non-archimedean case. Therefore, the limit identity becomes
$$
\lim_{m\to\infty}
\int_{X_v^\an} f \mu_{x_m,v}
= \frac{1}{\deg_\TL(X)} \int_{X_{v}^\an} f c_1(\OL)_v^n.
$$
This proves the equidistribution theorem. 
\end{proof}

\subsection{Total volume}

In the following, we prove that the equilibrium measure $d\mu_{\overline L,v}$ in Theorem \ref{equi3} is indeed a probability measure. Our proof uses the global intersection theory to bounded local integrals, so it only works over number fields 
and function fields of one variable. 
We refer to Gauthier--Vigny \cite[Thm. B]{GV} for a complex approach of such a result in the dynamical setting.

\begin{lem} \label{total volume}
\kkk
Let $K$ be a number field if $k=\ZZ$; let $K$ be the function field of one variable over $k$ if $k$ is a field. 
Let $X$ be a quasi-projective variety of dimension $n$ over $K$.
Let $\overline L_1,\cdots, \overline L_n$ be integrable adelic line bundles on $X/k$, and let $\wt L_1,\cdots, \wt L_n$ be their images under the map $\wh\Pic(X/k)\to \wh\Pic(X/K)$.
Then for any place $v$ of $K$,
$$\int_{X_v^\an} c_1(\OL_1)_v\cdots c_1(\OL_n)_v=\wt L_1\cdot \wt L_2\cdots \wt L_n.$$ 
\end{lem}
\begin{proof}
By multi-linearity, it suffices to assume that all $\overline L_1,\cdots, \overline L_n$ are strongly nef.
By multi-linearity again, it suffices to assume that all $\overline L_1,\cdots, \overline L_n$ are isomorphic to the same adelic line bundle $\OL$ on $X$.

Assume that $\OL$ is represented by a Cauchy sequence 
$\CLL=(\CL,(\CX_i,\overline \CL_i, \ell_{i})_{i\geq1})$ in $\wh\Picc(\CU)_\rmod$. Here $\CU$ is a quasi-projective model of $X$ over $k$, and each $\CLL_{i}$ is nef on $\CX_i$.
We further assume that for each $i\geq1$, there is a morphism
$\phi_i:\CX_i\to \CX_0$ extending the identity morphism of $\CU$.
Here $(\CX_0,\OCE_0)$ is a boundary divisor. 
Denote $X_i=\CX_{i,\QQ}$, which is a projective model of $X$ over $K$. 

The weak convergence formula gives, for any $f\in C_c(X_v^\an)$, 
$$
\int_{X_v^\an}f c_1(\OL)_v^n
= \lim_{i\to \infty} \int_{X_{i,v}^\an} f c_1(\CLL_{i})_v^n.
$$
See \S\ref{sec measure} for more details.
As $X_i$ is projective over $K$, the right-hand side is equal to the integration defined by global intersection numbers by \cite{CL, Gub1}. 
It suffices to extend this formula to the case $f=1$.

Denote by $\wt g_v\geq 0$ the Green function of $\OCE_0$ on $X_{0,v}^\an$. 
For any $m\geq1$, define a continuous  and compactly supported function $f_m:X_{0,v}^\an\to \RR$ by 
\begin{enumerate}
\item $f_m(x)=1$ if $\wt g_v (x) \leq m$;
\item $f_m(x)=m+1-\wt g_v (x)$ if $m\leq \wt g_v (x)\leq m+1$;
\item $f_m(x)=0$ if $\wt g_v (x) \geq m+1$.
\end{enumerate}
Then $f_m$ increases to the constant function one on $X_v^\an$. We have
$$
\int_{X_v^\an} c_1(\OL)_v^n
= \lim_{m\to \infty} \int_{X_v^\an} f_m c_1(\OL)_v^n
= \lim_{m\to \infty} \lim_{i\to \infty} \int_{X_{i,v}^\an} f_m c_1(\CLL_{i})_v^n.
$$
The first equality follows from Lebesgue's monotone convergence theorem, and the second equality holds by viewing $f_m$ as an element of $C_c(X_v^\an)$.
Then it suffices to prove 
$$
\lim_{m\to \infty} \lim_{i\to \infty} \int_{X_{i,v}^\an} (1-f_m) c_1(\CLL_{i})_v^n=0.
$$
By definition, $0\leq 1-f_m\leq \wt g_v/m$ everywhere on $X_v^\an$.
Therefore, it suffices to prove $\ds\int_{X_{i,v}^\an} \phi_i^*\wt g_v\, c_1(\CLL_{i})_v^n$
is bounded above as $i$ varies.

By the global intersection formula of Chambert-Loir and Thuillier in \cite[Thm. 1.4]{CT}, 
$$
 \CLL_i^n\cdot \phi_i^*\OCD
= (\CLL_i|_{\CH_i})^n
+\sum_{v} \int_{X_{i,v}^\an} \phi_i^*\wt g_v\, c_1(\CLL_{i})_v^n.
$$
Here $\CH_i$ is the horizontal part of $\phi_i^*\CD$, as an effective divisor on $\CX_i$.
As $\CLL_i$ is nef, $\OCD$ is effective, and $\wt g_0\geq0$, every term on the right-hand side is non-negative.
This gives 
$$
\int_{X_{i,v}^\an} \phi_i^*\wt g_v\, c_1(\CLL_{i})_v^n\leq  \CLL_i^n\cdot \phi_i^*\OCD.
$$
The right-hand sides converges to $\OL^n\cdot \OCD$ by Proposition \ref{intersection1}. This finishes the proof.
\end{proof}

\begin{rmk}
In Lemma \ref{total volume}, the result is local at the place $v$, but the condition assumes that the adelic line bundles come from a global field. This global assumption seems strange but gives us the convenience of bounding local integrals by global intersection numbers. In recent work, Guo \cite{Guo}  proves the volume formula by local methods and thus removes our global assumption. 
\end{rmk}

\subsection{Equidistribution theorem in the relative situation}

Inspired by an original idea of Moriwaki \cite{Mor3}, we generalize Theorem \ref{equi3} to equidistribution of directionally small points in the relative situation. 
The statement is closely related to the fundamental inequality in Theorem \ref{minima}. 
The key is still the variational principle of \cite{SUZ, Yua1} and the bigness result in Theorem \ref{HS}. 
The theorem is as follows. 

\begin{thm}[equidistribution in the relative case] \label{equi5}
\kkk
Let $K$ be a number field if $k=\ZZ$; let $K$ be the function field of one variable over $k$ if $k$ is a field. 
Let $\pi:U\to V$ be a flat morphism of relative dimension $n$ of quasi-projective varieties over $K$.
Denote $d=\dim V+1$. 
Let $X\to \Spec F$ be the generic fiber of $U\to V$.

Let $\overline H$ be an element of $\wh\Pic (V/k)_{\nef}$ satisfying the Moriwaki condition that $\overline H$ is nef, $\OH^{d}=0$ and $\wt H^{d-1}>0$.
Here $\wt H$ is the image of $\OH$ in $\wh\Pic (V/K)_{\nef}$. 

Let $\overline L$  be an element of $\wh\Pic (U/k)_{\nef}$ such that 
$\deg_{\wt L}(X/F)>0$. 
Here $\wt L$ is the image of $\OL$ under the canonical composition 
$$\wh\Pic (U/k)_{\nef}\lra \wh\Pic (X/k)_{\nef}
\lra \wh\Pic (X/F)_{\nef}.$$ 

Let $\{x_m\}_m$ be a generic sequence in $X(\overline F)$ such that 
$h_{\overline L}^{\overline H}(x_m)$ converges to $h_{\overline L}^{\overline H}(X)$.
Then for any place $v$ of $K$, there is a weak convergence
$$
\frac{1}{\deg(x_m)} \delta_{\Delta(x_m),v}\, c_1(\pi^*\OH)_v^{d-1}
\lra \frac{1}{\deg_{\wt L}(X/F)} c_1(\OL)_v^{n}c_1(\pi^*\OH)_v^{d-1}
$$
of measures on $U_v^\an$. 
Here $\Delta(x_m)\subset U$ denotes the Zariski closure of the image of $x_m$ in $U$, and 
$\delta_{\Delta(x_m),v}$ denotes the Dirac current of $\Delta(x_m)_{K_v}^\an$ in $U_v^\an$. 
\end{thm}

In the theorem, the weak convergence means that 
$$
\frac{1}{\deg(x_m)} \int_{\Delta(x_m)_{K_v}^\an}  f\, c_1(\pi^*\OH)_v^{d-1}
\lra \frac{1}{\deg_{\wt L}(X/F)} \int_{U_v^\an}  f\, c_1(\OL)_v^{n}c_1(\pi^*\OH)_v^{d-1}
$$
for any continuous and compactly supported function $f:U_v^\an\to \RR$. 
Here the measures are defined in \S\ref{sec app global}.

The prototype of the theorem is \cite[Thm. 6.1]{Mor3}, which proves  
the equidistribution at archimedean places with the additional assumption that $U\to V$ is projective and the metric of $\OL$ is smooth and strictly positive (at archimedean places).

\begin{proof}
The proof is a hybrid of the proofs of Theorem \ref{equi3} and Theorem \ref{minima}. 
As in the proof of Theorem \ref{equi3}, let $\OM$ be an element in the kernel of the map $\wh\Pic(X/k)_{\intb}\to \wh\Pic(X/K)_{\intb}$. 
Let $\epsilon$ be a nonzero rational number.
The key is the claim that 
$$
e_1^{\overline H}(X,\overline L+\epsilon\OM) \geq 
h_{\overline L+\epsilon\OM}^{\overline H}(X)+O(\epsilon^2),\ \quad \epsilon \to 0.
$$

Let us first see how the claim implies the equidistribution theorem, following the proof of Theorem \ref{equi3}. 
The claim gives
$$
 \liminf_{m\to \infty} \epsilon\, h_{ \OM}(x_m) \geq 
\epsilon\frac{  \OM\cdot \OL^n\cdot \OH^{d-1}}{ \deg_\TL(X)}+O(\epsilon^2).
$$
Then this implies
$$
 \lim_{m\to \infty}  h_{ \OM}(x_m) = 
\frac{  \OM\cdot \OL^n\cdot \OH^{d-1}}{ \deg_\TL(X)}.
$$
As in the proof of Theorem \ref{equi3}, it further implies the equidistribution theorem by taking $\OM$ to be the trivial bundle on $U$ with metrics given by model functions. 

Now  we prove the claim. 
The proof is similar to that of Theorem \ref{minima}, but more delicate due to the extra term $\epsilon \OM$. 
As in that proof, let $\ON\in \wh\Pic(B/k)_\QQ$ be an element of degree 1. 
Here is $k=\ZZ$, then $B=\Spec O_K$ and further assume that $\ON$ comes from the pull-back of $\wh\Pic(\ZZ)_\QQ$; 
if $k$ is a field, then $B$ is the unique projective and regular curve over $k$ with function field $K$. 

We make two convenient assumptions.
First, assume that $\epsilon>0$, which can be achieved by replacing $\OM$ by $-\OM$ if necessary. 
Second, assume that $\OL^{n+1}\OH^{d-1}>0$. 
This can be achieved by replacing $\OL$ with $\OL+\ON$, which does not affect the inequality we want to prove.

Denote $\OL'=\OL-c\ON$ with $c\in\QQ$. 
We still have 
$$
e_1^{\overline H}(X,\overline L+\epsilon\OM) 
-e_1^{\overline H}(X,\overline L'+\epsilon\OM) 
= h_{\overline L+\epsilon\OM}^{\overline H}(X)
- h_{\overline L'+\epsilon\OM}^{\overline H}(X)
=c\, \wt H^{d-1}.
$$
It suffices to prove that 
$$
e_1^{\overline H}(X,\overline L'+\epsilon\OM) 
- h_{\overline L'+\epsilon\OM}^{\overline H}(X)
\geq O(\epsilon^2)
$$
for some rational number $c$.

Write $\OM=\OA-\OB$ for nef adelic line bundles $\OA,\OB$ on $X$. 
Denote by $\OL_K, \OM_K, \OA_K, \OB_K$ the images of $\OL, \OM, \OA, \OB$ in $\wh\Pic(X/K)$. 
Note that $\OM_K=0$ by assumption. 
In the following, take  
$$
c=c(\epsilon)=\frac{(\OL+\epsilon\OA)^{n+1}\cdot \OH^{d-1}
-(n+1)(\OL+\epsilon\OA)^{n}\cdot\epsilon \OB \cdot \OH^{d-1}}{
(n+1)(\TL^{n})  (\wt H^{d-1})}
-\delta(\epsilon),
$$
where $\delta:\QQ_{>0}\to \RR$ is a fixed function such that $0< \delta(\epsilon) <\epsilon^2$ and such that $c(\epsilon)$ is always a rational number.
For this choice of $c$, we will check that 
$$
 h_{\overline L'+\epsilon\OM}^{\overline H}(X)
= O(\epsilon^2)
$$
and 
$$
e_1^{\overline H}(X,\overline L'+\epsilon\OM) >0.
$$
This implies the claim.

By definition, it is easy to have
$$
c(\epsilon)=\frac{(\OL+\epsilon\OM)^{n+1}\cdot \OH^{d-1}}{
(n+1)(\TL^{n})  (\wt H^{d-1})}
+O(\epsilon^2).
$$
This implies  
$$
 h_{\overline L'+\epsilon\OM}^{\overline H}(X)
= h_{\overline L+\epsilon\OM}^{\overline H}(X)-c\, \wt H^{d-1}
= O(\epsilon^2).
$$

It remains to prove
$$
e_1^{\overline H}(X,\overline L'+\epsilon\OM) >0.
$$
The assumption $\overline H^{d}=0$ still implies that
$$
e_1^{\overline H}(X,\overline L'+\epsilon\OM)
=e_1^{\overline H}(X,\overline L'+\epsilon\OM+m\OH)
$$ 
for all rational numbers $m$. 
Then it suffices to prove 
$$
e_1^{\overline H}(X,\overline L'+\epsilon\OM+m\OH)\geq 0
$$
for some $m$. 
As in the proof of Theorem \ref{minima}, it suffices to prove 
$$\wh\vol(\OL'+\epsilon\OM+m\OH)>0$$ 
for sufficiently large $m$.

Now we estimate $\wh\vol(\OL'+\epsilon\OM+m\OH)$ for $m>0$.
Note that $c(\epsilon)\geq 0$ when $\epsilon$ is sufficiently small, due to the assumption 
$\OL^{n+1}\OH^{d-1}>0$. 
By Proposition \ref{change norm},
$$
 \wh\vol(\OL'+\epsilon\OM+m\OH) 
\geq  \wh\vol(\OL+\epsilon\OM+m\OH)-(d+n)\, c(\epsilon)\, \wh\vol(\OL_K+m\wt H).
$$
Write 
$$\OL+\epsilon\OM+m\OH=(\OL+\epsilon\OA+m\OH)-(\epsilon \OB),
$$
and apply Theorem \ref{HS} to the above terms. 
We have 
\begin{eqnarray*}
&& \wh\vol(\OL'+\epsilon\OM+m\OH) \\
&\geq & (\OL+\epsilon\OA+m\OH)^{d+n}-(d+n) (\OL+\epsilon\OA+m\OH)^{d+n-1}\cdot \epsilon \OB \\
&&-(d+n)\, c(\epsilon)(\OL_K+m\wt H)^{d+n-1}.
\end{eqnarray*}

By the assumption $\OH^d=0$, the right-hand side is a polynomial in $m$ of degree at most $d-1$, and the coefficient of $m^{d-1}$ is equal to 
\begin{eqnarray*}
&& {d+n \choose d-1} (\OL+\epsilon\OA)^{n+1}\OH^{d-1} 
-(d+n){d+n-1 \choose d-1} (\OL+\epsilon\OA)^{n}\cdot \OH^{d-1} \cdot \epsilon \OB\\
&-& (d+n)\, c(\epsilon)\, {d+n-1 \choose d-1} (\wt L^{n})(\wt H^{d-1}).
\end{eqnarray*}
This is exactly the product of $\ds {d+n \choose d-1}$ with 
$$
 (\OL+\epsilon\OA)^{n+1}\OH^{d-1} 
-(n+1)(\OL+\epsilon\OA)^{n}\cdot \OH^{d-1} \cdot \epsilon \OB
- (n+1)\, c(\epsilon)(\wt L^{n})(\wt H^{d-1}).$$
It is strictly positive by the definition of $c(\epsilon)$. 
This finishes the proof. 
\end{proof}

\begin{rmk}
Note that
Conjecture \ref{equi1} can be viewed as a fiberwise version of Theorem \ref{equi5}. 
We expect that Theorem \ref{equi5} implies Conjecture \ref{equi1}, while the obstruction is some complicated regularization processes. 
\end{rmk}

\section{The Hodge bundle} \label{sec hodge}

In \S\ref{sec hodge bundle}, we have introduced the example of Hodge bundles and mentioned that it naturally defines an adelic line bundle. 
The goal of this section is to state the result precisely and give proof. 
We will also sketch a proof of a similar result for canonical bundles of families of curves endowed with the hyperbolic metrics mentioned in \S\ref{sec hyperbolic metric}.

\subsection{Hodge bundle for a general family} 

Recall from \S\ref{sec hodge bundle} that $S$ is a flat and quasi-projective integral scheme over $\ZZ$ or $\QQ$, and $\pi:X\to S$ is a principally polarized abelian scheme of relative dimension $g$.
Recall that
$\omega(S) = e^* \Omega_{X/S}^g$
is the Hodge bundle on $S$, and the Faltings metric $\|\cdot\|_\Fal$ of $\omega(S)$ on $S(\CC)$ is defined by integration.
Our precise theorem is as follows.

\begin{thm}\label{hodge bundle}
There is a canonically defined adelic line bundle $\overline{\omega(S)}$ on $S/\ZZ$ which extends the pair $(\omega(S), \|\cdot\|_\Fal)$. Moreover,
$$
h_{\overline{\omega(S)}}(s)=h_\Fal(X_s),\quad \forall\, s\in S(\ol\QQ).
$$
\end{thm}

Here we explain some of the terms of the theorem. 
First, that $\overline{\omega(S)}$ extends $(\omega(S), \|\cdot\|_\Fal)$
means that the underlying line bundle of $\overline{\omega(S)}$ is $\omega(S)$, and that the metric of $\omega(S)$ on $S(\CC)$ induced by 
$\overline{\omega(S)}$ (via Proposition \ref{injection3}) 
is equal to $\|\cdot\|_\Fal$.

Second, by restriction, $\overline{\omega(S)}$ induces an adelic line bundle on 
$S_{\QQ}$, and thus defines a height function 
$h_{\overline{\omega(S)}}:S_{\QQ}(\ol\QQ)\to \RR$. 

Third, the stable Faltings height $h_\Fal(X_s)$ of the abelian variety 
$X_s$ over $\ol\QQ$ associated to $y$ is defined as follows.
Note that $X_s$ descends to an abelian variety $G$ over a number field $K$ with semi-abelian reduction. 
Then we define the stable Faltings height by
$$
h_\Fal(X_s)=\frac{1}{[K:\QQ]}\wh\deg(\omega_\CG, \|\cdot\|_{\Fal}).
$$
Here $\omega_\CG=e_\CG^*\Omega_{\CG/O_K}^g$ is the Hodge bundle of the N\'eron model $\CG$ of $G$ over $O_K$, where $e_\CG:\Spec\, O_K\to \CG$ is the identity section, 
and $\|\cdot\|_{\Fal}$ is the Faltings metric of $\omega_\CG$ defined by  
$$
\|\alpha\|_{ \Fal}^2= \frac{i^{g^2}}{2^g} \int_{G_\sigma(\CC)} \alpha\wedge\bar \alpha
$$
for any embedding $\sigma:K\to\CC$ and any element $\alpha$ of 
$$
\omega_\CG\otimes_\sigma \CC
\simeq
\Gamma(G_\sigma(\CC), \Omega_{G_\sigma(\CC)/ \CC}^g).
$$
The definition is independent of the choice of $(G,K)$.

\subsection{Hodge bundles for moduli spaces} 

Theorem \ref{hodge bundle} is implied by a similar result for the minimal compactification of the coarse moduli scheme of abelian varieties. To introduce it, we will start with many constructions by Faltings--Chai \cite{FC}.
We will eventually only work on schemes, but the construction is easier to describe in terms of stacks.

Denote by $\CA_g$ the moduli stack of principally polarized abelian varieties over $\ZZ$. It is a smooth Deligne--Mumford stack over $\ZZ$, endowed with a universal abelian scheme $\CX_g\to \CA_g$.
Denote by $\CA_g'$ the coarse moduli scheme of $\CA_g$, which is a flat and quasi-projective integral scheme over $\ZZ$. 

By \cite[IV, Thm. 5.7]{FC}, there is a toroidal compactification 
$\CA_{g}^\mathrm{tor}$ of $\CA_{g}$ (by choosing a suitable combinatorial datum), which is a proper Deligne--Mumford stack over $\ZZ$ containing $\CA_{g}$ as an open and dense substack. Moreover, the universal abelian scheme $\CX_g\to \CA_{g}$ extends to a semi-abelian scheme
$\CX_g^\mathrm{tor}\to \CA_{g}^\mathrm{tor}$.

In terms of the universal abelian scheme (resp. semi-abelian scheme), we have a Hodge bundle $\omega(\CA_{g})$ on $\CA_{g}$ (resp. $\omega(\CA_{g}^\mathrm{tor})$ on $\CA_{g}^\mathrm{tor}$) defined similarly to the Hodge bundle $\omega(S)$ on $S$.

By \cite[V, Thm. 2.3]{FC}, there is a minimal compactification 
$\CA_{g}^*$ of the coarse moduli scheme $\CA_g'$. It is a normal projective scheme over $\ZZ$ defined by contracting $\CA_{g}^\tor$ via linear systems associated to
$\omega(\CA_{g}^\tor)$.
As a consequence, the Hodge bundle $\omega(\CA_{g}^\tor)$  
descends to a $\QQ$-line bundle $\omega(\CA_{g}^*)$ on $\CA_{g}^*$. 
In fact,  $\omega(\CA_{g}^*)$ is just $m^{-1} \CL$ in 
the notation of \cite[V, Thm. 2.3]{FC}, so it is indeed a $\QQ$-line bundle. 
Denote by $\omega(\CA_{g}')$ the restriction of
$\omega(\CA_{g}^*)$ to $\CA_{g}'$.

Note that $(\CA_{g}^*, \omega(\CA_{g}^*))$ is constructed by choosing a toroidal compactification, but the final result does not depend on the choices.

Since $\CA_g'$ is the coarse moduli scheme, any point $y\in \CA_{g}'(\CC)$ corresponds to a complex abelian variety $G$. Then the fiber 
$\omega(\CA_{g}^*)(y)^{\otimes m}$ is canonically isomorphic to the $m$-th tensor power of the Hodge bundle of $G/\CC$. 
Then the integration on $G(\CC)$ as before induces a Faltings metric $\|\cdot\|_\Fal$ of $\omega(\CA_{g}^*)(y)$.
Varying $y$, we obtain a Faltings metric $\|\cdot\|_\Fal$ of $\omega(\CA_{g}^*)$ on $\CA_{g}'(\CC)$.

Consider pair $(\omega(\CA_{g}^*), \|\cdot\|_\Fal)$.
It is similar to the original pair $(\omega(S), \|\cdot\|_\Fal)$, but it has the huge advantage that $\CA_{g}^*$ is projective over $\ZZ$. 
In particular, $(\omega(\CA_{g}^*), \|\cdot\|_\Fal)$
induces a metrized line bundle 
$\ol{\omega(\CA_{g}')}^\ran$ in $\Picc(\CA_{g}'^\ran)_\QQ$ with underlying $\QQ$-line bundle $\omega(\CA_{g}')$. 

The following is an analog of Theorem \ref{hodge bundle}, which is still based on the analytification functor in Proposition \ref{injection3}.

\begin{thm}\label{hodge bundle2}
The metrized $\QQ$-line bundle 
$\ol{\omega(\CA_{g}')}^\ran$ in $\wh\Picc(\CA_{g}'^\ran)_\QQ$
is the image of a unique adelic $\QQ$-line bundle 
$\ol{\omega(\CA_{g}')}$ in $\wh\Picc(\CA_{g}'/\ZZ)_\QQ$ via the analytification functor. 
Moreover, for any $y\in \CA_{g}'(\ol\QQ)$ corresponding to an abelian variety $G$ over $\ol\QQ$, we have
$h_{\overline{\omega(\CA_{g}')}}(y)=h_\Fal(G).$
\end{thm}

\begin{proof}

Let $\CA_{g}^{**}\to \CA_{g}^*$ be the blowing-up of $\CA_{g}^*$ along the boundary $\CA_{g}^*\setminus \CA_{g}'$. 
The exceptional divisor $\CE$ can be extended to a boundary divisor $(\CE,g_\CE)$ on $\CA_{g}^{**}$. 
Let $\omega(\CA_{g}^{**})$ be the pull-back of $\omega(\CA_{g}^*)$ to $\CA_{g}^{**}$. 
It suffices to consider the pair $(\omega(\CA_{g}^{**}), \|\cdot\|_\Fal)$.

By \cite[V, Def. 4.2, Rem. 4.3, Prop. 4.5]{FC}, the metric $\|\cdot\|_\Fal$ has logarithmic singularities along the boundary 
$\CE(\CC)$. 
Namely, take any hermitian metric $\|\cdot\|'$ of $\omega(\CA_{g}^{**})$ on $\CA_{g}^{**}(\CC)$. 
Denote 
$$
f=\log(\|\cdot\|_\Fal/\|\cdot\|'), 
$$
which is a continuous function on $\CA_{g}'(\CC)$. 
Then the logarithmic singularity means that 
$$
|f|
< c\log g_{\CE}
$$
over $\CA_{g}^{**}(\CC)$
for some constant $c>0$.

By Theorem \ref{singularity global}, the pair 
$(0, f)$ lies in $\wh\Div(\CA_{g}')$.
Therefore, we have proved that the metrized $\QQ$-line bundle 
$\ol{\omega(\CA_{g}')}^\ran$ comes from an adelic $\QQ$-line bundle 
$\ol{\omega(\CA_{g}')}$ on $\CA_{g}'$. 

It remains to prove the identity
$h_{\overline{\omega(\CA_{g}')}}(y)=h_\Fal(G)$
for $y\in \CA_{g}'(\ol\QQ)$.
We can assume that $y$ corresponds to a point $y:\Spec K\to \CA_{g}$ for a number field $K$. 
This induces a point $y:\Spec K\to \CA_g^\mathrm{tor}$ on the proper stack $\CA_g^\mathrm{tor}$ over $\ZZ$.
By the properness and the valuative criterion, by enlarging $K$ if necessary, we can assume that $y:\Spec K\to \CA_g^\mathrm{tor}$ extends to a morphism $\tilde y:\Spec O_K\to \CA_g^\mathrm{tor}$. 
Via the universal semi-abelian scheme $\CX_g^\tor\to \CA_g^\mathrm{tor}$, we obtain a semi-abelian scheme $\CG= \tilde y^*\CX_g^\tor$ over $O_K$. 
The generic fiber $\CG_K$ is a descent of the abelian variety $G$ to $K$.
By this, we see that $h_\Fal(G)$ is equal to $\wh\deg(\tilde  y^*( \omega(\CA_g^\mathrm{tor}), \|\cdot\|_\Fal))/[K:\QQ]$.
Here the Faltings metric $\|\cdot\|_\Fal$ of $\omega(\CA_g^\mathrm{tor})$ on 
$\CA_g(\CC)$ is defined by integration as before.
By the compatibility of the Hodge bundles, 
$( \omega(\CA_g^\mathrm{tor}), \|\cdot\|_\Fal)$ can be changed to 
$( \omega(\CA_g^{*}), \|\cdot\|_\Fal)$.
This finishes the proof.
\end{proof}

Once we have Theorem \ref{hodge bundle2}, the proof of Theorem \ref{hodge bundle} is immediate. The family $X\to S$ induces a moduli morphism $S\to \CA_g$. Composing with the canonical morphism $\CA_g\to \CA_g'$, we obtain a morphism $S\to \CA_g'$.
Then $\ol{\omega(S)}$ is just the pull-back of $\ol{\omega(\CA_{g}')}$ to $S$. 
This pull-back is a priori only an adelic $\QQ$-line bundle. Still, it is uniquely realized as an adelic line bundle since the underlying line bundle $\omega(S)$ is an integral line bundle on $S$.

\subsection{Hyperbolic metrics on families of curves} 

Recall from \S\ref{sec hyperbolic metric} that $S$ is a flat and quasi-projective normal integral scheme over $\ZZ$ or $\QQ$, and $\pi:X\to S$ is a smooth projective morphism whose fibers are geometrically integral curves of genus $g>1$. 
Recall that the relative dualizing sheaf 
$\omega_{X/S}$ has a hyperbolic metric $\|\cdot\|_\mathrm{hyp}$ on $X(\CC)$ defined by fiberwise universal covering.
Our precise theorem is as follows.

\begin{thm}\label{hyperbolic metric}
There is a canonically defined adelic line bundle $\overline\omega_{X/S,\mathrm{hyp}}$ on $X/\ZZ$ which extends the pair $(\omega_{X/S}, \|\cdot\|_\mathrm{hyp})$.
\end{thm}

The meaning of the term ``extend'' is similar to that in Theorem \ref{hodge bundle}. 
The proof is also similar, so we only sketch it in the following. 

First, we can assume that $S$ is quasi-projective over $\ZZ$. In fact, 
if $S$ is quasi-projective over $\QQ$, by taking suitable integral models of $X\to S$ over $\ZZ$, we can convert it to the case that $S$ is quasi-projective over $\ZZ$. 

Second, if $S'\to S$ is a finite and flat morphism for an integral scheme $S'$ over $\ZZ$,  
then the result holds for $X\to S$ if and only if it holds for the base change $X'\to S'$ of $X\to S$ by $S'\to S$. 
We can recover $\overline\omega_{X/S,\mathrm{hyp}}$ (up to a multiple) from $\overline\omega_{X'/S',\mathrm{hyp}}$ via the Deligne pairing (or equivalently the norm map) $\Picc(X'/\ZZ)\to \Picc(X/\ZZ)$. 

Third, by replacing $S$ by a suitable $S'$ as above, we can assume that $\pi:X\to S$ has a
\emph{stable compactification} $\pi^c:X^c \to S^c$, i.e. $S^c$ (resp. $X^c$) is a projective model of $S$ (resp. $X$) over $\ZZ$, $\pi^c:X^c \to S^c$ is a projective and flat morphism extending $\pi:X\to S$, and every  fiber of $\pi^c:X^c \to S^c$ is a stable curve in the sense of Deligne--Mumford \cite{DM}. 
The existence of $S'$ is explained in \cite[\S3.1.4]{Yua4}. 

Fourth, let $\pi^c:X^c \to S^c$ be a stable compactification as above. Then the relative dualizing sheaf $\omega_{X^c/S^c}$ is a line bundle on $X^c$ extending $\omega_{X/S}$. 
By Wolpert \cite{Wol}, the metric $\|\cdot\|_\mathrm{hyp}$ of $\omega_{X/S}$ on $X(\CC)$ extends to a continuous metric $\|\cdot\|_\mathrm{hyp}^c$ of $\omega_{X^c/S^c}$ on $(X^c\setminus N)(\CC)$, where $N\subset X^c$ is the Zariski closed subset of nodes of fibers of $X^c\to S^c$.  
The restriction of $\|\cdot\|_\mathrm{hyp}^c$ to each non-compact fiber of $(X^c\setminus N)(\CC) \to S^c(\CC)$ is still a hyperbolic metric (normalized suitably). 
Moreover, the metric $\|\cdot\|_\mathrm{hyp}^c$ has logarithmic singularity along $N(\CC)$; more rigorously, this logarithmic singularity should be understood after pull-back via the blowing-up of $X^c$ along $N$. 

Finally, as in the proof of Theorem \ref{hodge bundle}, by  Theorem \ref{singularity global}, the logarithmic singularity of the metric implies that $(\omega_{X^c/S^c},\|\cdot\|_\mathrm{hyp}^c)$ extends to an adelic line bundle on $(X^c\setminus N)/\ZZ$.

\chapter{Algebraic dynamics} 

In this chapter, we first develop a theory of admissible adelic line bundles for polarized algebraic dynamical systems over finitely generated fields, following the idea of \cite{Zha2, YZ1}. Then we generalize the arithmetic Hodge index theorem of Faltings \cite{Fal1} and Hriljac \cite{Hri} to projective curves over finitely generated fields. 

To work with adelic $\QQ$-line bundles on flat and essentially quasi-projective integral schemes $X$ over $k$, we recall the definitions of  
$\wh\Pic(X/k)_{\QQ}$,
 $\wh\Pic(X/k)_{\intb, \QQ}$ and
$\wh\Pic(X/k)_{\QQ,\nef}$ 
in \S\ref{sec rational}. 
Recall the categories $\wh\Picc(X/k)_{\QQ}$,
 $\wh\Picc(X/k)_{\intb, \QQ}$ and
$\wh\Picc(X/k)_{\QQ,\nef}$ defined similarly.

\section{Invariant adelic line bundles} \label{sec dynamics}

Let $(X, f, L)$ be a \emph{polarized dynamical system} over an integral scheme $S$, i.e.
\begin{enumerate}
\item $X$ is an integral scheme projective and flat over $S$;
\item $f:X\to X$ is a morphism over $S$;
\item $L\in \Picc(X)_\QQ$ is a $\QQ$-line bundle on $X$, relatively ample over $S$, such that
$f^*L\simeq qL$ for some rational number $q>1$.
\end{enumerate}
We refer to \cite[\S1.7]{Laz1} for relative ampleness. In particular, \cite[Thm. 1.7.8]{Laz1} asserts that a line bundle on $X$ is relatively ample over $S$ if and only if it is ample on every fiber of $X$ over $S$. 

If $S$ is the spectrum of a number field, Zhang \cite{Zha2} applied Tate's limiting argument to construct a nef adelic $\QQ$-line bundle $\overline L_f$ extending $L$ and with $f^*\overline L_f\simeq q\overline L_f$.  
The goal here is to generalize the result to finitely generated fields or even 
essentially quasi-projective schemes $S$.

\subsection{Invariant adelic line bundle}

\kkk
\ccc
Let $S$ be a flat and essentially quasi-projective integral scheme over $k$.
Let $(X, f, L)$ be a {polarized dynamical system} over $S$.
Fix an isomorphism $f^*L\to qL$ with $q>1$ by assumption.

Choose a projective model $\pi:\CX\rightarrow\CS$ of $X\to S$, i.e. a projective model $\CS$ of $S$ over $k$ and a flat morphism 
$\pi:\CX\rightarrow\CS$ of projective varieties over $k$ whose base change by $S\to \CS$ is isomorphic to $X\to S$. 
Choose a hermitian $\QQ$-line bundle
$\overline\CL=(\CL, \|\cdot\| )$ on $\CX$ such that $(\CX_S, \CL_S)\simeq(X,L)$.

For each positive integer $i$, consider the composition $X \stackrel{f^i}{\rightarrow} X \rightarrow \CX$. Denote the normalization of the composition by $f_i: \CX_i\to \CX$, and denote the induced map to $\CS$ by $\pi_i: \CX_i\to \CS$.
Denote $\overline \CL_i =q^{-i} f_i^* \overline\CL$, which lies in $\wh\Picc(\CX_i)_{\BQ}$.

The sequence $\{(\CX_i, \overline \CL_i)\}_{i\geq 1}$ is an adelic sequence in the sense of Moriwaki \cite[\S3.1]{Mor4}.
In our setting, we will complete the datum to an adelic line bundle 
$\OL_f=(\CL_\CV,(\CX_i, \overline \CL_i,\ell_i)_{i\geq 1})$ for a quasi-projective model $\CU$ of $X$ over $k$.

In fact, there is an open subscheme $\CV$ of $\CS$ containing $S$, such that $\CU=\CX_\CV$ is projective and flat over $\CV$, and that
$f:X\to X$ extends to a morphism $f_\CV:\CU\to \CU$ and such that the isomorphism $f^*L\to qL$ extends to an isomorphism 
 $f_\CV^*\CL_\CV\to q\CL_\CV$ in $\Picc(\CU)_\QQ$.
By the construction, we make identifications $\CX_{i,\CV}=\CX_{\CV}=\CU$ and $\CL_i|_\CU=\CL_{i,\CV}$. 

Start with the isomorphism 
$$\ell:\CL_\CV\lra q^{-1}f_\CV^*\CL_\CV$$ 
in $\Picc(\CU)_\QQ$.
By applying $q^{-1}f_\CV^*$ to $\ell$ successively, we obtain canonical isomorphisms
$$
\CL_\CV\lra q^{-1}f_\CV^*\CL_\CV \lra q^{-2}(f_\CV^*)^2\CL_\CV
\lra \cdots \lra q^{-i}(f_\CV^*)^i\CL_\CV
$$
in $\Picc(\CU)_\QQ$.
This induces an isomorphism 
$$
\ell_i: \CL_\CV\lra \CL_{i,\CV}
$$
in $\Picc(\CU)_\QQ$ by the identification 
$\CL_{i,\CV}= q^{-i}(f_\CV^*)^i\CL_\CV$. 
Then we have introduced every term in $(\CL_\CV,(\CX_i, \overline \CL_i,\ell_i)_{i\geq 1})$. 

Note that if $S$ is already a quasi-projective variety over $k$, then we can simply take $(\CU,\CV)=(X,S)$. This is the essential case of the result. 

\begin{thm} \label{invariant metric}
\kkk
Let $S$ be a flat and essentially quasi-projective integral scheme over $k$.
Let $(X, f, L)$ be a {polarized dynamical system} over $S$.
Fix an isomorphism $f^*L\to qL$ in $\Picc(X)_\QQ$ with $q>1$.

The above sequence $(\CL_\CV,(\CX_i, \overline \CL_i,\ell_i)_{i\geq 1})$ converges in $\wh \Picc(\CU/k)_{\QQ}$, and thus defines an object $\overline L_f$ of $\wh \Picc(X/k)_{\QQ}$. 
The adelic line bundle $\OL_f$ is uniquely determined by $(S,X,f,L)/k$ and $f^*L\to qL$ up to isomorphism, and satisfies the following properties.
\begin{enumerate}
\item
$\OL_f$ is $f$-invariant in the sense that $f^* \overline L_f\simeq q\overline L_f$ in $\wh \Picc(X/k)_{\QQ}$.

\item
$\OL_f$ is nef in $\wh \Picc(X/k)_{\QQ}$. If $S$  has an affine quasi-projective model over $k$, then $\OL_f$ is strongly nef in
$\wh \Picc(X/k)_{\QQ}$.

\item
If furthermore $L\in \Picc(X)$ (instead of $\Picc(X)_\QQ$) and $q\in\ZZ_{>1}$ with $f^*L\simeq qL$ in $\Picc(X)$, then all the results hold in $\wh \Picc(X/k)$ (instead of $\wh \Picc(X/k)_{\QQ}$).
\end{enumerate}

\end{thm}

\begin{proof}
We first prove the existence of the limit.
By blowing up $\CS$ along $\CS\setminus \CV$ if necessary, we can assume that there is a boundary divisor $(\CS,\OCE_0)$ of $\CV$. 
Then we get a boundary divisor $(\CX,\pi^*\OCE_0)$ of $\CU$.

View the isomorphism
$\ell:\CL_\CV\to q^{-1}f_\CV^*\CL_\CV$
as a rational map 
$\CLL\dashrightarrow \CLL_1$.
This defines a model adelic divisor $\wh\div(\ell)$ in $\wh\Div(\CU/k)_{\rmod,\QQ}$
whose image in $\Div(\CU)$ is 0. 
Then there exists $r>0$ such that 
$$
-r \pi^* \OCE_0\leq \wh\div(\ell) \leq r \pi^* \OCE_0
$$
holds in $\wh\Div(\CU)_{\rmod,\QQ}$. 
The existence of $r$ can be seen in the comparison of the boundary norms in the proof of Lemma \ref{boundary norm}.

By construction, the isomorphism 
$\ell_{i+1}\ell_i^{-1}:\CL_{i,\CV}\to \CL_{i+1,\CV}$
is obtained from $\ell:\CL_\CV\to q^{-1}f_\CV^*\CL_\CV$
by applying $(q^{-1}f_\CV^*)^i$. 
Accordingly, the rational map  
$\ell_{i+1}\ell_i^{-1}:\CLL_{i}\dashrightarrow \CLL_{i+1}$
is obtained from the rational map $\ell:\CLL \dashrightarrow \CLL_1$
by ``applying'' $(q^{-1}f^*)^i$. 
The situation can be conveniently described by the analytification functor in Proposition \ref{injection2} or the Zariski--Riemann space in \S\ref{sec ZR}. Still, we give a precise description in terms of projective models of $\CU$ as follows.

Write $\CX_0=\CX$ and $\CLL_0=\CLL$ for convenience.
There are projective models $\CY_1$ and $\CY_{i+1}$ of $\CU$ over $k$, together with morphisms
$$
\tau_1:\CY_1\to \CX_1,\quad 
\tau_1':\CY_1\to \CX_0,\quad
\tau_{i+1}:\CY_{i+1}\to \CX_i,\quad 
\tau_{i+1}':\CY_{i+1}\to \CX_{i+1}$$
extending the identity morphism $\CU\to\CU$,
and a morphism
$$
g_i:\CY_{i+1}\to \CY_1$$
extending the morphism $f_\CV^{i}:\CU\to\CU$.
Then the rational map $\ell:\CLL_0 \dashrightarrow \CLL_1$ is realized as a rational map $\ell':\tau_1'^*\CLL_0 \to \tau_1^*\CLL_1$ over $\CY_1$;
the rational map $\ell_{i+1}\ell_i^{-1}:\CLL_{i}\dashrightarrow \CLL_{i+1}$ is realized as a rational map
$(\ell_{i+1}\ell_i^{-1})': \tau_{i+1}'^*\CLL_i \to \tau_{i+1}^*\CLL_{i+1}$ over $\CY_{i+1}$.
The second rational map, including its source and its target, is obtained by applying $q^{-i}g_i^*$ to the first rational map via $g_i:\CY_{i+1}\to \CY_1$. 
As a consequence, we have
$$
\wh\div((\ell_{i+1}\ell_i^{-1})')= q^{-i} g_i^*\wh\div(\ell')
$$
in $\wh\Div(\CY_{i+1})_{\QQ}$.

Denote by $\pi_{1}':\CY_{1}\to \CS$ and $\pi_{i+1}':\CY_{i+1}\to \CS$ the structure morphisms. 
Note that $g_{i}^*\pi_{1}'^*\OCE_0=\pi_{i+1}'^*\OCE_0$ is equal to $\pi^*\OCE_0$ in 
$\wh\Div(\CU)_{\rmod,\QQ}$. We obtain
$$
-\frac{r}{q^{i}} \pi^* \OCE_0\leq \wh\div(\ell_{i+1}\ell_i^{-1}) \leq \frac{r}{q^{i}} \pi^* \OCE_0
$$
holds in $\wh\Div(\CU)_{\rmod,\QQ}$. 
As a consequence, $\{\wh\div(\ell_{i}\ell_1^{-1})\}_{i\geq1}$ is a Cauchy sequence in $\wh\Div(\CU)_{\rmod,\QQ}$.

This finishes the existence of the limit. 
The independence of the limit on the auxiliary data can be proved similarly, so we omit it.
It remains to treat the nefness of $\OL_f$ on $X$. 

At the beginning of the construction, if we can choose $(\CX,\overline\CL)$ such that $\CLL$ is nef on $\CX$, then every $\CLL_i$ is nef on $\CX_i$ by pull-back, and thus $\OL_f$ is strongly nef by definition. 
This happens if $S$ has an affine quasi-projective model $\CV$ over $k$.
In fact, in this case, we can assume that $S$ is an open subscheme of $\CV$, and then the relative ampleness of $L$ on $S$ implies the ampleness of $L$ on $X$, so we can choose $(\CX,\overline\CL)$ such that $\CLL$ is nef.

However, such $(\CX,\overline\CL)$ might not exist in general, and we will have to make a slightly weaker choice. Namely, we claim that there is a projective model $\pi:\CX\rightarrow\CS$ of $X\to S$ over $k$, together with a hermitian $\QQ$-line bundle
$\overline\CL$ on $\CX$ extending $L$ and a nef hermitian $\QQ$-line bundle $\CMM$ over $\CS$, such that $\CLL'=\CLL+\pi^*\CMM$ is nef on $\CX$.

To prove the claim, by taking a sufficiently small quasi-projective model of $X\to S$ over $k$, we can assume that $S$ is quasi-projective over $k$. 
Since $L$ is relatively ample, there is an ample line bundle $M$ on $S$ such that 
$L+\pi^*M$ is ample on $X$.
Take a tensor power of $L+\pi^*M$, use it to embed $X$ into $\PP^N_k$, and take the Zariski closure of $X$. 
Then $L+\pi^*M$ extends to an ample $\QQ$-line bundle $\CL'$ on a projective model $\CX$ of $X$ over $k$. Extend $\CL$ to a nef hermitian line bundle $\CLL'$ on $\CX$. 
Similarly, using a tensor power of $M$ to embed $S$ into $\PP^{N'}_k$ and taking the Zariski closure, we have a projective model $\CS$ of $S$ such that $M$ extends to a nef hermitian line bundle $\CMM$ on $\CS$.
The rational map $\CX\dashrightarrow \CS$ extends to a morphism $\pi:\CX\to \CS$ by blowing-up $\CX$, and we can further assume that $\CX\to\CS$ is flat by the Raynaud--Gruson flattening theorem in \cite[Thm. 5.2.2]{RG}. 
Finally, we set $\CLL=\CLL'-\pi^*\CMM$. This proves the claim.

Now  we prove that $\OL_f$ is nef.
Let $\displaystyle\overline L_f= (\CL_\CV,(\CX_i, \overline \CL_i,\ell_i)_{i\geq 1})$
be constructed using the new pair $(\CX,\CLL)$ as in the claim. 
Note that $\CLL'=\CLL+\pi^*\CMM$ is nef on $\CX$. 
Then $\CLL_i+q^{-i}\pi_i^*\CMM=q^{-i} f_i^* \overline\CL'$ is nef on $\CX_i$ for any $i\geq1$.
It follows that for any positive integer $a$, the line bundle $\CLL_{i+a}+q^{-a}\pi^*\CMM$ is nef for any $i\geq 1$. 
View $\OL_f$ as the limit of 
$(\CL_\CV,(\CX_{i+a}, \overline \CL_{i+a},\ell_{i+a})_{i\geq 1})$. 
We see that $\OL_f+q^{-a}\pi^*\CMM$ is strongly nef. 
This proves that $\OL_f$ is nef.  
\end{proof}

For the uniqueness of $\OL_f$, we have the following result. For convenience of applications, we do not require $L$ to be ample.

\begin{thm} \label{invariant metric2}
\kkk
Let $X$ and $S$ be flat and essentially quasi-projective integral schemes over $k$. 
Let $\pi:X\to S$ be a projective and flat morphism with geometrically connected fibers.
Let $f:X\to X$ be a morphism over $S$. 
Let $L\in \Pic(X)_\QQ$ be an element such that
$f^*L= qL$ in $\Pic(X)_\QQ$ for some rational number $q>1$.
The following is true:
\begin{enumerate}
\item
There exists a unique preimage $\OL$ of $L$ under the map $\wh\Pic(X)_\QQ\to \Pic(X)_\QQ$ such that 
$f^*\OL= q\OL$ in $\wh\Pic(X)_\QQ$.
\item
If $f':X\to X$ is a morphism over $k$ such that $f'f=ff'$ and that $f'^*L= q'L$ in $\Pic(X)_\QQ$ for some rational number $q'\neq0$,
then the adelic line bundle $\OL$ defined in (1) satisfies $f'^*\OL= q'\OL$ in $\wh\Pic(X)_\QQ$.
\end{enumerate}

\end{thm}
\begin{proof}
Note that (1) implies (2).
In fact, $f'f=ff'$ implies $f^*(f'^*\OL)= qf'^*\OL$. Then $\OL'=q'^{-1}f'^*\OL$ is an extension of $L$ with $f^*\OL'= q\OL'$. By the uniqueness in (1), we have $\OL'=\OL$.
This proves (2).

For (1), the existence of $\OL$ is similar to Theorem \ref{invariant metric}.
For the uniqueness, we can assume that $L=\CO_X$ is the trivial line bundle. 

By Proposition \ref{injection2}, there is a canonical injection
$$\wh \Pic (X/k) _{\QQ} \lra \wh\Pic (X^\an)_\QQ.$$
As $L=\CO_X$, the image of $\overline L$ in $\wh\Pic (X^\an)_\QQ$ is represented by an element $(0, g)$ of $\wh \Div (X^\an) _{\QQ}$, where the underlying divisor is $0$ on $X$, and the Green function $g$ is actually a continuous function on $X^\an$. 
The condition $f^*\OL= q\OL$ implies in $\wh \Div (X^\an)$ 
$$m(0,f^*g-qg) =(\div(\alpha),-\log |\alpha|), \quad \alpha\in k(X)^\times,\ m\in \ZZ, \ m\neq 0.$$ 

This implies $\div(\alpha)=0$ on $X$, and thus 
$\alpha$ lies in $\Gamma(X,\CO_X^\times)=\Gamma(S,\CO_S^\times)$. 
As a result, the difference
$$
f^*g-qg = -\frac{1}{m}\log |\alpha|
$$
is constant on every fiber of $X^\an \to S^\an$.

Let $v\in S^\an$ be a point with residue field $H_v$.
The fiber $X^\an_v$  of $X^\an$ above $v$ is exactly the Berkovich space associated to $X_{H_v}$ over $H_v$.
We have that $f^*g-qg=c_{v}$ is constant on $X^\an_v$.
Denote by $g_{\max}$ and $g_{\min}$ the global maximal value and the global minimum value of the continuous function $g$ on the compact space $X^\an_v$.  
Note that $f:X^\an_v \to X^\an_v$ is surjective. 
The relation $f^*g=qg+c_v$ gives 
$g_{\max}=qg_{\max}+c_v$ and thus $g_{\max}=-c_v$. 
Similarly, $g_{\min}=-c_v$. 
This forces $g_{\max}=g_{\min}$ and thus $g$ is constant on $X^\an_v$.

As a consequence, $f^*g=g$ on $X^\an$. The original equation gives
$$m(1-q)(0,g) =(\div(\alpha),-\log |\alpha|).$$ 
Then $\overline L$ is 0  in $\wh \Pic (X/k) _{\QQ}$. 
This finishes the proof. 
\end{proof}

\subsection{Abelian schemes}

The most important example of the above construction is for abelian schemes. 
In this case, we can prove that the adelic line bundles $\OL_f$ in Theorem \ref{invariant metric2} is integrable (without assuming that $L$ is relatively ample.)

\begin{thm} \label{invariant metric3}
\kkk
Let $S$ be a flat and essentially quasi-projective integral scheme over $k$. 
Let $\pi:X\to S$ be an abelian scheme with the identity section $e:S\to X$. 
Let $L$ be a line bundle on $X$ with a rigidification, i.e. an isomorphism $e^*L\to \CO_S$. 
Assume that $[-1]^*L\simeq \epsilon L$ for some $\epsilon\in \{\pm 1\}$. 

Then there is an adelic line bundle $\OL$ on $X$ extending $L$ satisfying
$[2]^*\OL\simeq 4\OL$ for $\epsilon=1$
and $[2]^*L\simeq 2L$ for $\epsilon=-1$. The adelic line bundle $\OL$ is uniquely determined by the rigidification. 

Moreover, $\OL$ is always integrable. For any integer $m$,   
$[m]^*L\simeq m^2L$ if $\epsilon=1$;
and $[m]^*L\simeq mL$ if $\epsilon=-1$. 
\end{thm}

\begin{proof}
Set $i=2$ for the symmetric case $\epsilon=1$, and $i=1$ for the anti-symmetric case $\epsilon=-1$. 
Note that $[-1]^*L\simeq \epsilon L$ implies that $[m]^*\OL\simeq m^i\OL$.
We first see that $[m]^*\OL- m^i\OL$ is trivial on fibers of $\pi:X\to S$, and thus is isomorphic to $\pi^*M$ for some $M\in \Pic(S)$. 
But $M$ is trivial by the rigidification. 

The rigidification determines a unique choice of an isomorphism $[2]^*\OL\to 2^i\OL$.
Apply Theorem \ref{invariant metric2}(1) to the dynamical system $(X, [2], L)$ over $S$.
We obtain a unique adelic line bundle $\OL$ on $X$ extending $L$ such that 
$[2]^*\OL\simeq 2^i\OL$. 
Moreover, Theorem \ref{invariant metric2}(2) implies $[m]^*\OL\simeq m^i\OL$. 

It remains to prove that $\OL$ is integrable.  
In the case $\epsilon=1$, if $L$ is relatively ample, then $\OL$ is nef by Theorem \ref{invariant metric}.
In the case $\epsilon=1$ for general $L$, we can write it as the difference of two relatively ample line bundles with rigidification, and then the integrability still follows. 

Assume $\epsilon=-1$ in the following. Let $X^\vee\to S$ be the dual abelian scheme of $X\to S$. 
Let $P$ be the Poincare line bundle on $X\times_S X^\vee$, with a rigidification along the identity section of $X\times_S X^\vee\to S$. 
Then $L$ corresponds to a section $\sigma:S\to X^\vee$ in the sense that 
$L\simeq (\id, \sigma\circ\pi)^*P$. Here $(\id,\sigma\circ \pi)$ is the composition 
$X\to X\times_S S\stackrel{(\id,\sigma)}{\to} X\times_S X^\vee$.

For any $m\in \ZZ$, denote by 
$$[m]:X\times_S X^\vee\lra X\times_S X^\vee$$ the (total) multiplication of the abelian scheme $X\times_S X^\vee$ over $S$, and denote by 
$$[m]':X\times_S X^\vee\lra X\times_S X^\vee$$ 
the (partial) multiplication of the abelian scheme $X\times_S X^\vee$ on $X^\vee$. 
By the universal property, $[-1]^*P\simeq P$ and $[-1]'^*P\simeq -P$.
Then there is a unique adelic line bundle $\OP$ in $\wh\Picc(X\times_S X^\vee/k)$ extending $P$ with $[2]^*\OP\simeq 4\OP$. 
It further gives $[2]'^*\OP\simeq 2\OP$ by Theorem \ref{invariant metric2}(2).
Moreover, $\OP$ is integrable by the case $\epsilon=1$. 

Finally, under $L\simeq (\id, \sigma\circ\pi)^*P$, we have 
$\OL':=(\id, \sigma\circ\pi)^*\OP$ extends $L$ and satisfies $[2]^*\OL'\simeq 2\OL'$. 
It follows that $\OL\simeq \OL'$ by the uniqueness. Finally, $\OL$ is integrable since so is $\OP$.
This finishes the proof. 
\end{proof}

\subsection{Canonical height}

\kkk
Let $F$ be a finitely generated field over $k$.  
Let  $(X,f,L)$ be a polarized dynamical system over $F$.

By Theorem \ref{invariant metric}, there is an $f$-invariant line bundle $\overline L_f$ in $\wh\Pic(X/k)_{\QQ,\nef}$.
For any closed $\overline F$-subvariety $Z$ of $X$, define \emph{the vector-valued canonical height of $Z$} as
$$
\fh_{f}(Z)=\fh_{L,f}(Z):= \fh_{\overline L_f}(Z)\in \wh \Pic(F/k)_{\intb,\QQ}.$$
It gives a map $\fh_{f}:|X_{\overline F}|\to \wh \Pic(F/k)_{\intb,\QQ}$. 

We can also define the canonical height by Tate's limiting argument:
$$ \fh_{f}(Z)=
\lim_{m\rightarrow\infty} \frac{1}{q^m} \fh_{(\CX,\overline\CL)}(f^m(Z)).
$$
Here $(\CX,\overline\CL)$ is any initial model of $(X,L)$ as in the construction of $\overline L_f$ above.
Then one can check that it is convergent in $\wh \Pic(F/k)$ and compatible with the previous definition.

\begin{prop} \label{basicproperties}
Let $Z$ be a closed subvariety of $X$. Then the following are true:
\begin{enumerate}
\item The height $\fh_{f}(Z)$ lies in $\wh \Pic(F/k)_{\QQ,\nef}$.
\item The height is $f$-invariant in the sense that $\fh_{f}(f(Z))=q\ \fh_{f}(Z)$.
\item If $Z$ is preperiodic under $f$, the height $\fh_{f}(Z)= 0$ in $\wh \Pic(F/k)_\intb$. 
\end{enumerate}
\end{prop}

\begin{proof}
Since $\overline L_f$ is nef, the height $\fh_{f}(Z)$ is nef.
The formula $\fh_{f}(f(Z))=q  \fh_{f}(Z)$ follows from the projection formula in Lemma \ref{basic2}(4) and the invariance of $\overline L_f$. Thus $\fh_{f}(Z)= 0$ if $Z$ is preperiodic under $f$. 
\end{proof}

By choosing adelic line bundles $\overline H_1,\cdots, \overline H_{d-1}\in \wh \Picc(F/k)_{\QQ,\nef}$, we can form \emph{the canonical Moriwaki height}
$$
h_{f}^{\overline H_1,\cdots, \overline H_{d-1}}(Z):=\fh_{f}(Z)\cdot \overline H_1 \cdots \overline H_{d-1}.
$$
It is a non-negative real number.

\subsection{N\'eron--Tate height}
\kkk
Let $F$ be a finitely generated field over $k$.
Let $X$ be an abelian variety over $F$, $f=[2]$ be the multiplication by 2, and $L$ be any symmetric and ample line bundle on $X$. Then the canonical height 
$$\wh \fh_L= \fh_{L,[2]}: X(\overline F)\lra \wh \Pic(F/k)_{\QQ,\nef},$$
as a generalization of the N\'eron--Tate height, 
is quadratic in the sense that 
$$
\pair{x,y}_L:= \wh \fh_L(x+y)-\wh \fh_L(x)-\wh \fh_L(y)
$$
gives a bilinear map  
$$X(\overline F)\times X(\overline F)\lra \wh \Pic(F/k)_{\intb,\QQ}.$$
It can be proved by the theorem of the cube as in the classical case over number fields. 
We refer to \cite[\S3.3]{Ser} for the classical case and omit the proof in the current case.

\subsection{Equidistribution conjecture of preperiodic points}

Since all preperiodic points of a polarized dynamical system have height $0$, Conjecture \ref{equi1} implies the following conjecture. 

\begin{conj}[equidistribution of preperiodic points]  \label{equi2} 
\kkk
Let $F$ be a finitely generated field over $k$.
Let $v$ be a non-trivial valuation of $F$.
Assume that the restriction of $v$ to $k$ is trivial if $k$ is a field.
Let $(X,f,L)$ be a polarized dynamical system over $F$. Let $\{x_m\}_m$ be a generic sequence of preperiodic points in $X(\overline F)$. 
Then the Galois orbit of $\{x_m\}_m$ is equidistributed in $X_v^\an$ for the   measure $d\mu_{L, f,v}$.
\end{conj}

Here $X_v^\an$ is the Berkovich space associated with $X_{F_v}$, where $F_v$ is the completion of $F$ for $v$. The equilibrium measure is the Chambert-Loir measure
$$d\mu_{L,f,v}=\frac{1}{\deg_{L}(X)}c_1(L, \|\cdot\|_{f,v})^{\dim X}$$
over the analytic space $X_v^\an$, 
where $\|\cdot\|_{f,v}$ is an $f$-invariant metric of $L$ on $X_v^\an$ obtained by Tate's limiting argument. 

One can also formulate the consequence of Theorem \ref{equi5} for preperiodic points. We omit it here.

\section{Heights of points on a subvariety}
\label{sec subvariety}

Let $S$ be a quasi-projective variety over a number field $K$. 
Let $(X,f,L)$ be a polarized dynamical system over $S$.
Let $\OL_f \in \wh\Pic(X)_{\QQ,\nef}$ be the $f$-invariant extension of $L$.
Let $Y$ be a closed subvariety of $X$. 
The goal of this section is to explore the properties of the height function 
$$
h_{\OL_f}: Y(\overline K) \lra \RR.
$$
We consider two special cases.
If $Y$ is a section, then we have a specialization theorem. 
If $Y$ is non-degenerate, then we have an equidistribution theorem. 

The following exposition also works over function fields of one variable, but we restrict to number fields for simplicity.

\subsection{Height of specialization}
\label{sec specialization}

Now  we consider the variation of the height of a section specializing in an algebraic family of algebraic dynamical systems.

Let $S$ be a quasi-projective variety over a number field $K$. 
Let $(X,f,L)$ be a polarized dynamical system over $S$.
Let $\OL_f \in \wh\Pic(X)_{\QQ,\nef}$ be the $f$-invariant extension of $L$.
Let $i:S\to X$ be a section of $\pi:X\to S$. 
Denote the vector-valued height
$$\OM:=\fh_{\OL_f}(i)=i^*\OL_f \in \wh\Pic(S)_{\QQ,\nef}.$$
This gives a height function
$$
h_\OM: S(\overline K) \lra \RR.
$$

For any point $s\in S(\overline K)$, denote by $s'$ the closed point of $S$ corresponding to $s$. 
Then $i(s)\in X(\overline K)$ is actually a point on the polarized dynamical system $(X_{s'}, f_{s'}, L_{s'})$ over $s'$. 
Denote by $h_{\OL_{s', f_{s'}}}(i(s))$ the canonical height of $i(s)$ for the  polarized dynamical system $(X_{s'}, f_{s'}, L_{s'})$ over $s'$.
Now  we have the following identity. 

\begin{lem} [specialization] \label{specialization} 
For any point $s\in S(\overline K)$,
$$
h_{\OL_{s', f_{s'}}}(i(s))=h_{\OM}(s).
$$
Therefore, $i(s)$ is preperiodic under $f$ if and only 
$h_{\OM}(s)=0$.
\end{lem}

\begin{proof}
By definition, the $f_{s'}$-invariant extension of $L_{s'}$ on $X_{s'}$ is exactly $\OL_{f}|_{X_{s'}}$.
Then $h_{\OM}(s)$ is the normalized degree of the pull-back of $\OL_f$ via the composition $s' \to S \to X$, and $h_{\OL_{s', f_{s'}}}(i(s))$ is the normalized degree of the pull-back of $\OL_f$ via the composition $s' \to X_{s'} \to X$. 
Then both terms are equal to $h_{\OL_{f}}(i(s))$. 
\end{proof}

If $X$ is a family of elliptic curves over a smooth curve $S$ over $K$, a similar height identity was obtained by DeMarco--Mavraki \cite[Thm. 1.1]{DMM}. Their approach was very different, and their result was stronger in this case. They proved that there is an adelic line bundle $\OM'$ on the unique smooth projective model $S'$ of $S$ over $K$ such that $h_{\OM'}(s)= h_{\OL_{f}}(i(s))$ for any $s\in S(\overline K)$. In other words, their result implies that our $\OM$
lies in the image of $\wh\Pic(S')_{\QQ,\nef}\to \wh\Pic(S)_{\QQ,\nef}$.

As the work of \cite{DMM} is a refinement of the specialization theorem of Tate \cite{Tat} and Silverman \cite{Sil2, Sil3, Sil4} for elliptic surfaces, our height identity can be viewed as a generalization and new interpretation of the specialization theorem for families of algebraic dynamic systems.

In the setting of Lemma \ref{specialization}, if we know certain bigness property of $\OM$, then we may apply some height inequality in
Theorem \ref{height comparison0} to conclude that the height function $h_{\OM}$ is ``big''.

\subsection{Non-degenerate subvarieties}
\label{sec non-degenerate}

Let $S$ be a quasi-projective variety over a number field $K$. 
Let $(X,f,L)$ be a polarized dynamical system over $S$.
Let $Y$ be a closed subvariety of $X$. 

Let $\OL_f \in \wh\Pic(X)_{\QQ,\nef}$ be the $f$-invariant extension of $L$.
Denote by 
$$\OM:=\OL_f|_Y$$ 
the image of $\OL_f$ under 
the pull-back map 
$$\wh\Pic(X/\ZZ)_{\QQ,\nef}\lra \wh\Pic(Y/\ZZ)_{\QQ,\nef},$$
and denote by $\wt M$ the image of $\OM$ under 
the canonical composition 
$$\wh\Pic(Y/\ZZ)_{\QQ,\nef}\lra \wh\Pic(Y/\QQ)_{\QQ,\nef}\lra \wh\Pic(Y/K)_{\QQ,\nef}.$$ 
Note that the last arrow is an isomorphism. 
We refer to \S\ref{sec functoriality} for the definitions of these maps. 
By nefness, both self-intersection numbers 
$$\wh\deg_{\ol M}(Y)=\ol M^{\dim Y+1},\qquad
\deg_{\wt M}(Y)=\wt M^{\dim Y}$$ 
are non-negative.

We say that $Y$ is \emph{non-degenerate} in $X$ if $\deg_{\wt M}(Y)>0$.
As $\wt M$ is nef on $Y$,  the condition is equivalent to that $\wt M$ is big on $Y$.
Another related result is Lemma \ref{total volume}, which asserts that, for any embedding $\sigma:K\to \CC$,
$$\deg_{\wt M}(Y)=\int_{Y_\sigma(\CC)} c_1(\OL_f)_\sigma^{\dim Y}.$$ 
So $Y$ is {non-degenerate} if and only if 
the measure $c_1(\OL_f)_\sigma^{\dim Y}|_{Y_\sigma(\CC)}$ is nonzero on $Y_\sigma(\CC)$. 
The same result holds over non-archimedean places.

If $X\to S$ is an abelian scheme over a smooth variety $S$ over $K$,  in terms of Tate's limiting argument, $c_1(\OL_f)_\sigma$ defines a semipositive smooth $(1,1)$-form on $X_\sigma(\CC)$. 
In particular, it is the Betti form as defined in \cite[\S2]{CGHX}. 
By \cite[Prop. 2.2]{DGH}, $c_1(\OL_f)_\sigma^{\dim Y}$ is non-zero on $Y_\sigma(\CC)$ if and only if the Betti map $Y_\sigma(\CC)_V\to (\RR/\ZZ)^{2g}$ has a full rank at some point of $Y_\sigma(\CC)_V$ for some simply connected open subset of $S_\sigma(\CC)$. Strictly speaking, the Betti form in \cite[Prop. 2.2]{DGH} is the one that comes from a principal polarization (instead of a general $L$). Still, Betti forms of any two relatively ample line bundles can bound each other by positive constant multiples. 
Therefore, our definition of ``non-degenerate'' agrees with that of the loc. cit., and generalizes to families of algebraic dynamical systems.

Now  we have the following theorem, which generalizes \cite[Thm. 1.4]{GH} and \cite[Thm. 1.6]{DGH} from abelian schemes to dynamical systems. 
Our proof follows the idea of \cite{DGH} but is simplified significantly by our new notion of adelic line bundles.

\begin{thm}[height inequality]\label{height comparison}
Let $S$ be a quasi-projective variety over a number field $K$. 
Let $(X,f,L)$ be a polarized dynamical system over $S$.
Let $Y$ be a non-degenerate closed subvariety of $X$ over $K$.
Let $\OB\in \wh\Pic(S)_\QQ$ be an adelic $\QQ$-line bundle on $S$. 
Then for any $c>0$, there exist $\epsilon>0$ and a non-empty open subvariety $U$ of $Y$ such that 
$$
h_{\OL_f}(y) \geq \epsilon\, h_{\OB}(\pi(y)) -c, \quad\ \forall\, y\in U(\overline K). 
$$
Here $\pi:X\to S$ denotes the structure morphism. 
\end{thm}

\begin{proof}
Apply Theorem \ref{height comparison0}(2) to the morphism $Y\to S$ and the adelic line bundles $\OL_f|_Y$ and $\OM$. 
\end{proof}

\subsection{Equidistribution theorem over non-degenerate subvarieties}

Restricted to the setting of non-degenerate subvarieties, 
we get a special example of Theorem \ref{equi3}.

\begin{thm}[equidistribution over non-degenerate subvarieties] \label{equi4}
Let $S$ be a quasi-projective variety over a number field $K$. 
Let $(X,f,L)$ be a polarized dynamical system over $S$.
Let $Y$ be a non-degenerate closed subvariety of $X$ over $K$. 
Let $\{y_m\}_{m\geq 1}$ be a generic sequence of $Y(\overline K)$ such that 
$h_{\OL_{f}}(y_m) \to 0$.
Then for any place $v$ of $K$,
the Galois orbit of  $\{y_m\}_{m\geq 1}$ is equidistributed over the analytic space $Y_v^\an$
for the  canonical measure $d\mu_{\OL_f|_Y,v}$. 
\end{thm}

The theorem generalizes \cite[Cor. 1.2]{DMM}, which treats the family of elliptic curves described above. 
If $X\to S$ is an abelian scheme, the theorem confirms the conjecture (REC) of K\"uhne \cite{Kuh}, and our proof is independent of the slightly weaker version in \cite[Thm. 1]{Kuh}. The proof of \cite{Kuh} is a limit version of the original proof in \cite{SUZ} and uses a result of Dimitrov--Gao--Habegger \cite{DGH} for uniformity in the limit process.

Note that the existence of the sequence $\{y_m\}_{m\geq 1}$ implies 
$h_{\OM}(Y)=0$ and thus $\wh\deg_{\ol M}(Y)=0$, as a consequence of Theorem \ref{minima2}.
In the following, we make some remarks on the existence of $Y$
satisfying the condition of the theorem.

First, the non-degeneracy of $Y$ is easy to check if $\dim Y=\dim S=1$.
In fact, in this case, it becomes $\deg(\wt M)>0$, and 
$\deg(\wt M)$ is exactly the canonical height $\wh h(Y_\eta)$ of the closed point $Y_\eta$ for the  polarized dynamical system $(X_\eta,f_\eta,L_\eta)$ over the generic point $\eta=\Spec K(S)$ of $S$. 
For example, if $X$ is a family of abelian varieties over $S$ with trivial $K(S)/K$-trace, then $\wh h(Y_\eta)=0$ if and only if $Y_\eta$ is torsion in $X_\eta(\eta)$. See \cite[Thm. 9.15]{Con} for example.

For an abelian scheme $X\to S$ of relative dimension $g$ with a high-dimensional base $S$, there are natural generalizations of the above situation by Andr\'e--Corvaja--Zannier \cite{ACZ20} and Gao \cite{Gao1}.  
Namely, by \cite[Thm. 2.3.1, Prop. 2.1.1]{ACZ20} and \cite[Thm. 9.1]{Gao1}, a closed subvariety $Y$ of $X$ is non-degenerate and contains a Zariski dense set of torsion points if the following conditions hold:
\begin{enumerate}
\item $\dim S=g$;
\item  the morphism from $S$ to the moduli space of abelian varieties of dimension $g$ (with a polarization of degree equal to $\deg(L_\eta)/g!$) is generically finite;
\item $X$ is simple over the algebraic closure of the function field of $S$;
\item $Y$ is a non-torsion section of $X\to S$.
\end{enumerate}
Contrary to the case $g=1$, the result does not hold for general $g$ if we change condition (3) to the statement that $X$ has a trivial $K(S)/K$-trace. We refer to Gao \cite[Thm. 1.4(ii)]{Gao2} for a counter-example for $g=4$.

\section{Equidistribution of PCF maps} 
\label{sec PCF}

In this section, we consider the equidistribution of post-critically finite endomorphisms on $\PP^n$ as another application of the equidistribution theorem (Theorem \ref{equi3}). The equidistribution fits perfectly to the setting of the dynamical Andre--Oort conjecture of Baker--DeMarco \cite[Conj. 1.10]{BD}. Our treatment plays a crucial role in the recent solution of the dynamical Andre--Oort conjecture for 1-dimensional families by Ji--Xie \cite{JX}.  

We will only write the case of number fields, though some of the results also hold over function fields of one variable. 

\subsection{Post-critically finite maps}

Let $f:\PP^n\to \PP^n$ be a finite separable morphism over a field. 
Assume that its algebraic degree $d$ (defined by $f^*\CO(1)\simeq \CO(d)$) is strictly larger than 1. Denote by $R(f)$ the ramification divisor (or critical locus) of $f$ in 
$\PP^n_k$, whose definition will be recalled below in the family version. 
The morphism $f$ is said to be \emph{post-critically finite} (PCF) if every irreducible component of $R(f)$ (with reduced structure) under $f$ is preperiodic. 

Let $S$ be a smooth and quasi-projective variety over a number field $K$. Let $X=\BP^n_S$ be the projective space over $S$, and let $f:X\to X$ be a finite morphism over $S$ of algebraic degree $d>1$ (over the fibers above $S$). 
A point $y\in S(\overline K)$ is called \emph{post-critically finite} (PCF) if the morphism $f_y:X_y\to X_y$ is post-critically finite.

The main result here is the construction of a natural adelic line bundle $\OM$ over $S$ and equidistribution theorems of Galois orbits of PCF points. 

\subsection{The adelic line bundle $\OM$}

Let $S$ and $f:X\to X$ be as above. 
Namely, $S$ is a smooth and quasi-projective variety over a number field $K$, $X=\BP^n_S$, and $f:X\to X$ is a finite morphism over $S$ of algebraic degree $d>1$. 

Denote by $\pi:X\to S$ the structure morphism.
The canonical morphism $f^*\omega_{X/S}\to \omega_{X/S}$ induces
a global section $\delta f$ of $\omega_f=\omega_{X/S}\otimes f^*\omega_{X/S}^\vee$ on $X$.
The \emph{ramification divisor}  $R=R(f)$ of the finite morphism $f:X\to X$
is defined to be the divisor of the section $\delta f$. It is also viewed as a (possibly non-reduced) closed subscheme in $X$.
By definition, we have a canonical isomorphism 
$\omega_f\simeq \CO(R)$. 
We have the following basic result. 

\begin{lem}
The scheme $R(f)$ and every irreducible component of it (with the reduced structure) are projective and flat of relative dimension $n-1$ over $S$. 
The fiber $R(f)_y$ of $R(f)$ above any point $y\in S$ is equal to the ramification divisor $R(f_y)$ of $f_y:X_y\to X_y$.
\end{lem}
\begin{proof}
Since the canonical map $f_y^*\omega_{X_y/y}\to\omega_{X_y/y}$ is the base change of 
$f^*\omega_{X/S}\to \omega_{X/S}$ via $y\to S$, we see that $\delta(f_y)$ is the base change of $\delta f$, and $R(f_y)$ is the base change of $R(f)$ via $y\to S$. 
Then $R(f)$ is of pure relative dimension $n-1$ over $S$.
Since $R(f)$ is a Cartier divisor on $X$, it is Cohen--Macaulay over $S$.
By the miracle flatness (cf. \cite[Thm. 23.1]{Mat}), the morphism $R(f)\to S$ is flat. 
Similarly, any irreducible component of $R$ is flat over $S$.
\end{proof}

Let $L$ be a $\QQ$-line bundle on $X$, isomorphic to $\CO(1)$ on fibers of $S$, such that $f^*L\simeq dL$. 
There is a unique class in $\Pic(X)_\QQ$ satisfying these requirements.
In fact, we can set $L=\CO_{\BP^1_S}(1)\otimes \pi^*N$ for a suitable $\QQ$-line bundle $N$ on $S$. 
Then $f^*L\simeq dL$ becomes 
$f^*\CO_{\BP^1_S}(1)-\CO_{\BP^1_S}(d)=(d-1)\pi^*N$. 
Note that $f^*\CO_{\BP^1_S}(1)-\CO_{\BP^1_S}(d)$ is trivial on fibers of $X\to S$, and thus lies in $\pi^*\Pic(S)$. 
The equality determines the class $N\in \Pic(S)_\QQ$ uniquely.

Denote by $\OL=\OL_f$ the nef $f$-invariant extension of $L$ in $\wh\Pic(X)_\QQ$ such that $f^*\OL\simeq d\OL$, as constructed in Theorem \ref{invariant metric}.
Recall that the ramification divisor $R$ is projective and flat of pure relative dimension $n-1$ over $S$.
Define  
$$
\OM:=\pair{\OL|_R}_{R/S}^n= \pair{\OL|_R,\cdots, \OL|_R}_{R/S} \in \wh\Pic(S)_\QQ.
$$
Here the Deligne pairing is as in Theorem \ref{intersection2}.
Since the theorem requires $R$ to be integral,  we need to extend the definition if $R$ is not integral. 
In fact, write $R=\sum_{i=1}^r m_iR_i$ in terms of distinct prime divisors $R_1,\cdots, R_r$ of $X$, and interpret the definition by 
$$
\OM= \sum_{i=1}^r m_i \pair{\OL|_{R_i}}^n_{R_i/S} \in \wh\Pic(S)_\QQ.
$$ 
In all cases, $\OM$ is a nef adelic $\QQ$-line bundle on $S$. 

If $n=1$, then $R$ is finite and flat over $S$, so
$$
\OM= N_{R/S}(\OL|_R) \in \wh\Pic(S)_\QQ
$$
is given by the norm map.

As before, denote by 
$$\OL\longmapsto \wt L\longmapsto L,\quad
\OM\longmapsto \wt M\longmapsto M$$
the images of $\OL$ and $\OM$ under the maps 
$$\wh\Pic(X)_\QQ \lra \wh\Pic(X/K)_\QQ \lra \Pic(X)_\QQ,\quad
\wh\Pic(S)_\QQ \lra \wh\Pic(S/K)_\QQ \lra \Pic(S)_\QQ.$$

\subsection{The height function}

Consider the height function 
$$
h_\OM:S(\overline K)\lra \RR.
$$
It detects PCF points using the following result: 

\begin{lem} \label{height 0}
Let $y\in S(\overline K)$ be a point. The following is true:
\begin{enumerate}
\item $h_\OM(y)\geq 0$. 
\item If $y$ is PCF in $S$, then $h_\OM(y)=0$. 
\item If $n=1$, then $y$ is PCF in $S$ if and only if $h_\OM(y)=0$. 
\end{enumerate}
\end{lem}
\begin{proof}
Part (1) holds since $\OM$ is nef. 
For (2) and (3), for convenience, assume that $y$ is a closed point of $S$ instead of an algebraic point. 
By Theorem \ref{intersection2}, the Deligne pairing is compatible with base change $y\to S$. 
It follows that 
$$\OM|_{y}=\pair{\OL|_R}^n|_{y}=\pair{\OL|_{R_y}}^n
=\sum_{i} m_{R_{y,i}} \pair{\OL|_{R_{y,i}}}^n.$$
Here $R_y=\sum_i m_{R_{y,i}} R_{y,i}$ is the decomposition into prime divisors in $X_y$. 
Then we have 
$$
\wh\deg(\OM|_{y})=\sum_{i} m_{R_{y,i}} \OL|_{R_{y,i}}^n.
$$
In terms of heights, we have  
$$
h_\OM(y)= \sum_{i} m_{y,i}' h_{\OL}(R_{y,i}),
$$
Here $m_{y,i}'=m_{R_{y,i}} n \deg_{L_y}(R_{y,i})/\deg(y)$ is strictly positive.

Then $h_\OM(y)=0$ if and only if $h_{\OL}(R_{y,i})=0$ for every irreducible component $R_{y,i}$ of $R_y$. 
This gives (2) immediately. 
For (3), $R_{y,i}$ is a closed point, and thus $h_{\OL}(R_{y,i})=0$ further implies $R_{y,i}$ is preperiodic. 
\end{proof}

\begin{problem}
We raise the question of whether Lemma \ref{height 0}(3) holds for $n\geq2$.
This amounts to ask: for a finite morphism $f:\PP^n_K\to \PP^n_K$ of algebraic degree $d>1$ over a number field $K$, if every irreducible component of the ramification divisor $R(f)$ has canonical height 0, does it follow that every irreducible component of $R(f)$ is preperiodic?
This is the dynamical Manin--Mumford conjecture for $R(f)$ under the dynamical system $f:\PP^n_K\to \PP^n_K$. 
We refer to Ghioca--Tucker--Zhang \cite{GTZ} for various versions and examples of 
the dynamical Manin--Mumford conjecture.
\end{problem}

\subsection{The equidistribution theorem}

With the nef adelic line bundle $\OM$ over $S$, we have the following equidistribution theorem, which is a direct consequence of Theorem \ref{equi3}. 

\begin{thm} [equidistribution: PCF maps on projective space] \label{equi PCF0}
Let $S$ be a smooth and quasi-projective variety over a number field $K$. Let $X=\BP^n_S$ be the projective space over $S$, and let $f:X\to X$ be a finite morphism over $S$ of algebraic degree $d>1$. 
Assume that $\deg_{\TM}(S)>0$. 
Let $\{y_m\}_m$ be a generic sequence of PCF points of $S(\overline K)$.
Then the Galois orbit of $\{y_m\}_m$ is equidistributed in $S_v^\an$ for $d\mu_{\overline M,v}$ for any place $v$ of $K$.
\end{thm}

Note that the existence of a generic sequence of PCF points implies $h_{\OM}(S)=0$.
This follows from the fundamental inequality. 
$$
\liminf_{y\in S(\overline K)} h_{\OM}(y) \geq h_{\OM}(S)
$$
proved in Theorem \ref{minima2}.

The condition $\deg_{\TM}(S)>0$ (equivalent to the bigness of $\TM$) seems very hard to check in general. 
However, in the case $n=1$, it is equivalent to a very clean condition in terms of the moduli space of endomorphisms. As we will see later, this equivalence is obtained as a combination of some geometric constructions and technical, analytic arguments building on a history of stability analysis. 

To describe the condition, denote by $\CM_d^n$ the moduli space over $K$ of endomorphisms $\PP^n$ of algebraic degree $d$. 
The moduli space was constructed using Mumford's geometric invariant theory by the works of Silverman \cite{Sil5}, Levy \cite{Lev}, and Petsche--Szpiro--Tepper \cite{PST}.

If $n=1$, there is a special type of PCF morphisms $\PP^1\to \PP^1$, called the flexible Latt\`es maps, which are descended from multiplication morphisms of elliptic curves. 
We refer to Silverman \cite[\S6.5]{Sil6} for the basics of the flexible Latt\`es maps. 
In $\CM_d^1$, there is a distinguished closed subvariety, called the flexible Latt\`es locus, parametrizing the flexible Latt\`es maps on $\PP^1$. 
The flexible Latt\`es locus is empty if $d$ is not a perfect square and has dimension one if $d$ is a perfect square. 

Return to the dynamical system $f:X\to X$ for $X=\PP^1_S$.
Recall  
$$\OM= N_{R/S}(\OL_f|_R) \in \wh\Pic(S)_\QQ.$$
By the moduli property, there is a morphism $S\to \CM_d^1$.
Finally, the main result here is the following variant of Theorem \ref{equi PCF0}. 

\begin{thm} [equidistribution: PCF maps on projective line] \label{equi PCF1}
Let $S$ be a smooth and quasi-projective variety over a number field $K$. Let $X=\BP^1_S$ be the projective line over $S$, and let $f:X\to X$ be a finite morphism over $S$ of algebraic degree $d>1$.  
Assume that the morphism $S\to \CM_d^1$ is generically finite and its image is not contained in the flexible Latt\`es locus. 
Let $\{y_m\}_m$ be a generic sequence of PCF points of $S(\overline K)$.
Then the Galois orbit of $\{y_m\}_m$ is equidistributed in $S_v^\an$ for $d\mu_{\overline M,v}$ for any place $v$ of $K$.
\end{thm}

If $S$ is a family of polynomial maps on $\PP^1$, the theorem was previously proved by Favre--Gauthier \cite{FG}.
Their strategy is to reduce the problem to the equidistribution of Yuan \cite{Yua1}, which works for polynomial maps but not for rational maps. 

As a dilation, we remark that the nef adelic line bundle $\OL_f$ for general $n$ is strongly nef. For this, it suffices to treat the case that $S$ is the moduli space $\CM_d^n$. 
Note that the corresponding moduli space over $\ZZ$ is affine (cf. \cite[Thm. 1.1]{Lev}).
Then $\OL_f$ is strongly nef by Theorem \ref{invariant metric}.

\subsection{The Lyapunov exponent}

Let us recall the classical Lyapunov exponent in the current setting. For completeness, we will include both the Archimedean case and the non-archimedean case. 

Let $K$ be a complete field with a non-trivial absolute value $|\cdot|$. 
Let $f: \PP^n_K\to \PP^n_K$ be a finite and separable morphism of algebraic degree $d>1$.
Recall that the ramification divisor $R(f)=\div(\delta f)$, where 
$\delta f$ is the section of $\omega_f=\omega_{\PP^n_K}\otimes f^*\omega_{\PP^n_K}^\vee$ induced by the canonical morphism $f^*\omega_{\PP^n_K} \to \omega_{\PP^n_K}$.
This definition gives a canonical isomorphism $\omega_f\simeq \CO(R(f))$.

Fix a continuous metric $\|\cdot\|_0$ of $\omega_{\PP^n_K}$ on the analytic space 
$\PP_K^{n,\an}$ over the valued field $(K,|\cdot|)$, and take the metric 
 $f^*\|\cdot\|_0$ of $f^*\omega_{\PP^n_K}$ on $\PP_K^{n,\an}$. 
Then we have the quotient metric $\|\cdot\|_1$ of
$\omega_f=\omega_{\PP^n_K}\otimes f^*\omega_{\PP^n_K}^\vee$ on $\PP_K^{n,\an}$. 
The function 
$-\log\|\delta f\|_1$ is a Green function of $R(f)$ on $\PP_K^{n,\an}$. 
The \emph{Lyapunov exponent} of $f$ is defined by 
$$
\mathrm{Ly}(f)=\int_{\PP_K^{n,\an}} \log \|\delta f\|_1\, d\mu_{f}.
$$
Here $d\mu_{f}=c_1(\ol{\CO(1)}_{f})^n$ is the $f$-invariant probability measure on 
$\PP_K^{n,\an}$. 
The definition is independent of the choice of the metric $\|\cdot\|_0$ of $\omega_{\PP^n_K}$. 
In fact, if $\|\cdot\|_0'=\|\cdot\|_0e^{h}$ is a different choice for a continuous function $h$ on 
$\PP_K^{n,\an}$, then the integral of 
$$
 \log \|\delta f\|_1'- \log \|\delta f\|_1=h-f^*h
$$
for $\mu_{f}$ is 0 by $f_*d\mu_{f}=d\mu_{f}$. 

It is convenient to choose $\|\cdot\|_0$ to be an $f$-invariant metric of 
$\omega_{\PP^n_K}\simeq \CO(-n-1)$ on $\PP_K^{n,\an}$. This metric is unique up to constant multiples, but then the induced metric on $\omega_f$ does not depend on the constant multiple, and we will denote this metric by $\|\cdot\|_{f}$. 
Then the {Lyapunov exponent} is just
$$
\mathrm{Ly}(f)=\int_{\PP_K^{n,\an}} \log \|\delta f\|_{f}\,d \mu_{f}.
$$

Now we assume that $K$ is a number field, and that $f: \PP^n_K\to \PP^n_K$ is still a finite morphism of algebraic degree $d>1$.
Choose $\|\cdot\|_0=(\|\cdot\|_{0,v})_v$ to be an $f$-invariant adelic metric of  
$\omega_{\PP^n_K}\simeq \CO(-n-1)$ on $\PP_K^{n,\an}$ in the classical sense of \cite{Zha2}.
This induces an $f$-invariant adelic metric $\|\cdot\|_f=(\|\cdot\|_{f,v})_v$ of 
$\omega_f=\omega_{\PP^n_K}\otimes f^*\omega_{\PP^n_K}^\vee\simeq \CO((n+1)(d-1))$ on $\PP_K^{n,\an}$. 
Write $\ol\omega_f=(\omega_f, \|\cdot\|_f)$ for the adelic line bundle. 
Via these metrics, the arithmetic intersection number
$$
\ol\omega_f\cdot \ol{\CO(1)}_f ^n
=(n+1)(d-1) \ol{\CO(1)}_f ^{n+1}=0.
$$
On the other hand, we can apply \cite[Thm. 1.4]{CT} to the section $\delta f$ of $\omega_f$ to compute the arithmetic intersection number. 
It gives
$$
\ol\omega_f\cdot \ol{\CO(1)}_f ^n
=\big( \ol{\CO(1)}_f|_{R(f)}\big) ^n-\sum_v \int_{\PP_{K_v}^{n,\an}} \log \|\delta f\|_{f}\, 
c_1(\ol{\CO(1)}_{f})_v ^n.
$$
This gives the height formula
$$
\big( \ol{\CO(1)}_f|_{R(f)}\big) ^n
=\sum_v \mathrm{Ly}(f_{K_v}).
$$

\subsection{The bifurcation measure}

Return to the situation of Theorem \ref{equi PCF0}, where $f:X\to X$ and
$X=\PP_S^n$ are over a number field $K$. 

It turns out that the equilibrium measure at a complex place $v$ in the theorem is exactly the probability measure associated with the bifurcation measure. The bifurcation $(1,1)$-current was first introduced by DeMarco \cite{DeM1, DeM2} for $n=1$, and the higher forms (including the bifurcation measure) of the bifurcation current for general $n$ were introduced by Bassanelli--Berteloot \cite[\S5]{BB07}. The goal here is to explore this relation in our adelic setting, which implies
the identity of the measures in both the Archimedean setting and the non-archimedean setting. The exposition here is a family version of the above height formula in terms of the Lyapunov exponents.

Let $v$ be a place of $K$.  
The Lyapunov exponent defines a function 
$$\Ly_v:S_v^\an\lra \RR, \quad
y\longmapsto \Ly(f_y).$$ 
If $v$ is archimedean, the pull-back of $\Ly_v$ to $S_v(\CC)$ is continuous and psh;
if $v$ is non-archimedean, then $\Ly_v$ is locally psh-approachable on $S_v^\an$
in the sense of \cite[6.3.1, Def. 5.6.3, Def. 5.5.1]{CLD}.
The archimedean case of this statement can be derived from \cite{DeM1, DeM2, BB07}, and we will present an approach including both cases.
The \emph{bifurcation measure} of $(X,f)$ over $S_v^\an$ is defined to be the Monge--Amp\`ere measure
$$d\mu_{\bif,v}=(dd^c \Ly_v)^{\dim S}.$$ 

Now we have the following description of the equilibrium measure in the setting of 
Theorem \ref{equi PCF0}. In the non-archimedean case, currents are understood in the sense of Chambert-Loir--Ducros \cite{CLD}. 

\begin{thm}[bifurcation measure] \label{bif}
Let $v$ be a place of $K$. 
As $(1,1)$-currents on $S_v^\an$,
$$
c_1(\OM)_v=dd^c \Ly_v.
$$ 
As measures on $S_v^\an$,
$$
c_1(\OM)_v^{\dim S}=(dd^c \Ly_v)^{\dim S},
$$
and 
$$d\mu_{\overline M,v}=\frac{1}{\deg_{\wt M}(S)}(dd^c \Ly_v)^{\dim S}.$$
\end{thm}

Note that the third equality follows from the second one. In fact, by  Lemma \ref{total volume}, we have 
$$
\deg_{\wt M}(S)=\int_{S_v^\an} c_1(\OM)_\sigma^{\dim S}
=\int_{S_v^\an} (dd^c \Ly_v)^{\dim S}.
$$
This also implies that the integral on the right-hand side is independent of $v$.

The complex version of Theorem \ref{bif} is essentially  
Bassanelli--Berteloot \cite[Cor. 4.6]{BB07}. 
The following theorem is an adelic treatment of the situation, which asserts that the Lyapunov exponents for all places $v$ can be glued together to form an adelic divisor.

\begin{thm} \label{adelic bif}
The following is true.
\begin{enumerate}
\item
There is a unique adelic divisor $\OD_\bif$ on $S$ with underlying divisor $0$ whose total   Green function $\wt g_{\OD_\bif}:S^\an\to \RR$ satisfies $\wt g_{\OD_\bif}|_{S_v^\an}=\Ly_v$ on $S_v^\an$ for every place $v$ of $K$.  
\item
For the above adelic divisor $\OD_\bif$, we have $\OM=\CO(\OD_\bif)$ in $\wh\Pic(S)_\QQ$.
\end{enumerate}
\end{thm}

It is easy to see that Theorem \ref{adelic bif} implies Theorem \ref{bif}. It also implies the continuity and reasonable psh properties of $\Ly_v:S_v^\an\to \RR$, since $\OM$ is nef by construction.

\begin{proof}[Proof of Theorem \ref{adelic bif}]

The uniqueness in (1) follows from Proposition \ref{injection3}.
The major part of the proof follows from an adelic version of the above construction to derive the height formula in terms of $\sum_v \Ly(f_{K_v})$ on single dynamical systems. 

Recall that $\OL=\OL_f\in \wh\Picc(X)_\QQ$ is a nef $f$-invariant extension of $L$ with $f^*\OL\simeq d\OL$, and $\OM \in \wh\Picc(S)_\QQ$ is the Deligne pairing  
$\pair{\OL|_R}^n$. 
Here we assume that $\OL$ and $\OM$ are adelic $\QQ$-line bundles instead of just isomorphism classes.

Recall that $R=\div(\delta f)$ is the divisor of the canonical global section $\delta f$ of 
$\omega_f=\omega_{X/S}- f^*\omega_{X/S}$ on $X$.
Here we write the operation of line bundles additively again.  
By comparing the fibers above $S$, there is an isomorphism 
$\tau_0:\omega_{X/S}\to -(n+1)L+\pi^*N$ for some $\QQ$-line bundle $N\in \Picc(S)_\QQ$.
This induces an isomorphism $\tau_1: \omega_f \to (n+1)(d-1)L$, which does not depend on the choice of $\tau_0$. Different choices of $\tau_0$ (for fixed $L, N$) are up to multiples by elements of $\Gamma(S,\CO_S^\times)$, and these elements are killed in the definition of $\tau_1$. 

Denote by $\ol\omega_f$ the adelic line bundle on $X$ with underlying line bundle 
$\omega_f$, such that $\tau_1: \omega_f \to (n+1)(d-1)L$ induces an isomorphism 
$\ol \omega_f \to (n+1)(d-1)\ol L$.  
Note that the extension $\ol\omega_f$ of $\omega_f$ is unique up to unique isomorphism. 

With these extensions, the Deligne pairing
$$
\pair{\ol\omega_f, \OL ^n}_{X/S}
=(n+1)(d-1) \pair{\OL}_{X/S}^{n+1}=0
$$
in $\Pic(S)_\QQ$. 
Here the last equality is similar to Proposition \ref{basicproperties}(3), as a consequence of the projection formula in Lemma \ref{basic2}(4) and the invariant property $f^*\OL\simeq d\OL$. 

On the other hand, using the section $\delta f$ of $\omega_f$ to compute the Deligne pairing, we have a canonical isomorphism 
$$
\pair{\omega_f, L ^n}_{X/S} \lra \pair{L|_R}_{R/S} ^n=M.
$$
This gives a canonical section $t_\bif$ of $M-\pair{\omega_f, L ^n}_{X/S}$, and thus an adelic $\QQ$-divisor $\OD_\bif=\wh\div(t_\bif)$ on $S$ for $\OM-\pair{\ol\omega_f, \OL ^n}_{X/S}$. 
The underlying divisor $D_\bif=0$ on $S$ by definition. 
By construction, we have 
$$
\OM\simeq \pair{\ol\omega_f, \OL ^n}_{X/S}+ \CO(\OD_\bif) \simeq  \CO(\OD_\bif)
$$
in $\Picc(S)_\QQ$. 

Now  we can compute the total Green function $\wt g_{\OD_\bif}$ on $S_v^\an$. 
If $v$ is archimedean, for any $y\in S_v(\CC)$,  the fiberwise formula in \S\ref{sec metric at point} gives
$$
\wt g_{\OD_\bif}(y)=-\log\|t_\bif\|(y)= \int_{X_y^\an} \log\| \delta f\|_y \, c_1(\OL)_y^n
=\Ly_v(y).
$$
The formula also holds for non-archimedean $v$ by Theorem \ref{local compatibility}.
This proves Theorem \ref{adelic bif}.
\end{proof}

\subsection{Bigness problem}

In the case $n=1$, to deduce Theorem \ref{equi PCF1} from Theorem \ref{equi PCF0}, 
it suffices to prove that if the morphism $S\to \CM_d^1$ is generically finite and its image is not contained in the flexible Latt\`es locus, then $\wt M$ is big on $S/K$. By Theorem \ref{bif},  
it suffices to check that the total volume of $\mu_{\bif,\sigma}$ is strictly positive in this case.  
The positivity is proved by \cite[Prop. 6.3]{BB07} and \cite[Lem. 6.8]{GOV}.
This finishes our proof of Theorem \ref{equi PCF1}.

\begin{rmk}
In the case $n=1$, the height function
$h_\OM:S(\overline K)\to \RR$
is equal to the critical height considered by Ingram \cite{Ing} and
Gauthier--Okuyama--Vigny \cite{GOV}.
This leads to a new proof of \cite[Thm. 1]{Ing}.  
If the morphism $S\to \CM_d^1$ is generically finite. Its image does not intersect the flexible Latt\`es locus; we know from the above argument that $\wt M$ is big on $S/K$. Then we can apply Theorem \ref{height comparison0}(2) to bound $h_\OM$ by a usual Weil height (in both directions) outside a Zariski closed subset $S_1$ of $S$. 
Apply the argument to irreducible components of $S_1$ repeatedly. We eventually cover every point of $S$. 
\end{rmk}

In the case $n=1$, it is well-known that the set of PCF points is Zariski dense in 
$\CM_d^1$. See \cite[Thm. A]{DeM3}, for example. 

In the case $n>1$, the situation is very different.
In fact, by the work of Ingram--Ramadas--Silverman \cite{IRS}, PCF points 
in $\CM_d^n$ are expected to be very sparse in some sense. 
As in \cite[Question 5]{IRS}, we do not know if the set of PCF points in $\CM_d^n$ is Zariski dense.
Then we raise the following question.

\begin{problem}
Assume $n\geq2$, $S=\CM_d^n$, and $f:\PP^n_S\to \PP^n_S$ is the universal family. 
Is $\TM$ big on $S$? 
Is $\OM$ big on $S$? 
\end{problem}

The bigness of $\TM$ is equivalent to $\deg_{\TM}(S)>0$, which is a condition of the equidistribution theorem. 
The bigness of $\OM$ is equivalent to $\wh\deg_{\OM}(S)>0$, which becomes 
$h_{\OM}(S)>0$ assuming $\deg_{\TM}(S)>0$.
It is further related to the existence of a generic and small sequence for $h_\OM$ considering Theorem \ref{minima2}. 
In particular, if $\OM$ is big, then the set of PCF points in $\CM_d^n$ is not Zariski dense.

\section{Admissible extensions of line bundles}

Let $(X, f, L)$ be a {polarized dynamical system} over a finitely generated field $F$ over $\QQ$. 
Assume that $X$ is normal.  
We have already constructed an adelic line bundle 
$\overline L_f\in \wh\Pic(X/k)_{\QQ,\nef}$ extending $L$ and with $f^*\overline L_f=q\overline L_f$.  
Following the idea of \cite{YZ1}, we can construct an admissible extension 
in $\wh\Pic(X/k)_{\QQ,\intb}$
for any line bundle $M\in \Pic(X)_\QQ$. 

Our exposition could be clearer, and we refer to \cite[\S4.3]{YZ1} for the common arguments, but we will explain the difference in the current case. 
Moreover, we will only restrict to the arithmetic case ($k=\ZZ$) and refer to \cite{Car1, Car2} for the counterparts in the geometric case, where the extra argument is to treat the contribution of the $F/k$-image of $\underline{\Pic}^0_{X/F}$.

\subsection{Semisimplicity}

The pull-back map $f^*$ preserves the  exact sequence
$$0\lra \Pic ^0(X)\lra \Pic (X)\lra \NS (X)\lra 0.$$
We refer to \cite[Appendix 1]{YZ1} for a list of properties of this sequence. 
In particular, $\NS(X)$ is a finitely generated $\ZZ$-module. 
By the Lang--N\'eron theorem (cf. \cite[Thm. 2.1]{Con}), $\Pic^0(X)$ is also a finitely generated $\ZZ$-module, since it is the Mordell--Weil group of the Picard variety representing the functor $\underline{\Pic}^0_{X/F}$ over the finitely generated field $F$. 
The counterpart of \cite[Theorem 4.7]{YZ1} is as follows. 

\begin{thm} \label{semisimplicity}
Let $(X, f, L)$ be a {polarized dynamical system} over a finitely generated field $F$ over $\QQ$. Assume that $X$ is normal. 
\begin{enumerate}
\item The operator $f^*$ is semisimple on 
$\Pic^0(X)_\BC$ (resp. $\NS (X)_\BC$) with eigenvalues of absolute values $q^{1/2}$ (resp. $q$). 
\item The operator $f^*$ is semisimple on 
$\Pic(X)_\BC$ with eigenvalues of absolute values $q^{1/2}$ or $q$. 
\end{enumerate} 
\end{thm}

\begin{proof}
The proof is similar to its counterpart. The only difference is some extra work to prove that $f^*$ is semisimple on $\Pic^0(X)_\BC$ with eigenvalues of absolute values $q^{1/2}$. We describe it briefly here. 

As before, $(X,f, L)$ extends to a dynamical system $(U,f, L_V)$ over a smooth quasi-projective variety $V$ over $\QQ$ with function field $F$. Here $U\to V$ is a projective and flat morphism with generic fiber $X\to \Spec F$, $f:U\to U$ is a $V$-morphism extending $f:X\to X$, and $L_V$ is a $\QQ$-line bundle on $U$, relatively ample over $V$, and with $f^*L_V=qL_V$. We can further assume that all the fibers of $U\to V$ are normal. 
We claim that there is a closed point $v\in V$ such that the reduction map 
$\Pic^0(X)_\CC\to \Pic^0(U_v)_\CC$ is injective. 
If this holds, then the result follows from its counterpart over number fields.

Note that the Picard functor $\underline\Pic_{U/V}$ is representable by a group scheme by \cite[\S8.2, Thm. 1]{BLR}. Its relative identity component $\underline\Pic_{U/V}^0$ is an abelian scheme over $V$ by \cite[Thm. 9.5.4]{Kle}. 
Then the injectivity is a consequence of the specialization theorem of Wazir \cite{Waz}, which is a generalization of the specialization theorem of Silverman \cite{Sil1} using the Moriwaki height. 
\end{proof}

By the theorem above, the exact sequence
$$0\lra \Pic ^0(X)_\BQ\lra \Pic(X)_\BQ\lra \NS (X)_\BQ\lra 0.$$
has an $f^*$-equivariant splitting
$$\ell_f:  \NS (X)_\BQ\lra \Pic(X)_\BQ.$$
Denote by $\Pic_f(X)_\QQ$ the image of $\ell_f$.

We say an element of $\Pic(X)_\BQ$ is \emph{$f$-pure of weight $1$ (resp. $f$-pure of weight $2$)} if it lies in $\Pic^0(X)_\BQ$ (resp. $\Pic_f(X)_\QQ$).

\subsection{Admissible extensions}

The action $f^*: \wh\Pic(X)_{\QQ}\to \wh\Pic(X)_{\QQ}$ is compatible with the action $f^*: \Pic(X)_\BQ\to \Pic(X)_\BQ$.
The goal is to study the spectral theory of this action.
The following result is the generalization of \cite[Thm. 4.9]{YZ1}. 

\begin{thm} \label{admissible metric}
Let $(X, f, L)$ be a {polarized dynamical system} over a finitely generated field $F$ over $\QQ$. Assume that $X$ is normal. The projection 
$$
\wh\Pic(X)_{\QQ}\lra \Pic(X)_\QQ
$$
has a unique section 
$$
M\longmapsto \overline M_f
$$
as $f^*$-modules. 
The image $\overline M_f$ is always integrable. 
If $M\in \Pic_f(X)_\BQ$ is ample, then $\overline M_f$ is nef. 
\end{thm}

We call $\OM_f$ \emph{the $f$-admissible extension of $M$} in 
$\wh\Pic(X)_{\QQ}$. An adelic line bundle in $\wh\Pic(X)_{\QQ}$ which is isomorphic to some $\OM_f$ is called \emph{$f$-admissible}. 

Note that the theorem for abelian schemes is actually
Theorem \ref{invariant metric3}. 
In fact, any $\QQ$-line bundle $L$ on an abelian scheme $X$ can be written as the sum of the symmetric $\QQ$-line bundle $(L+[-1]^*L)/2$ with the anti-symmetric $\QQ$-line bundle $(L-[-1]^*L)/2$.

As in the case of number fields, we also have the following result as the counterpart of \cite[Cor. 4.11]{YZ1}.

\begin{cor} 
For $M\in \Pic(X)_\BQ$, the following are true:
\begin{enumerate}
\item If $f^*M=\lambda M$ for some $\lambda\in \QQ$, then 
$f^*\overline M_f=\lambda \overline M_f$ in $\wh\Pic(X)_{\QQ}$.
\item For any $x\in \Prep(f)$, one has $\overline M_f|_{x'}=0$ in $\wh\Pic(x')_{\QQ}$.
Here $x'$ is the closed point of $X$ corresponding to $x$. Hence, the height function $\fh_{\overline M_f}$ is zero on $\Prep(f)$. 
\end{enumerate}
\end{cor}

Now we sketch a proof of Theorem \ref{admissible metric}, following the line of that of \cite[Thm. 4.9]{YZ1}.

\begin{proof}[Proof of Theorem \ref{admissible metric}]
Assume that $X$ is geometrically connected over $F$, which can be achieved by replacing $F$ by its algebraic closure in $F(X)$. 
Let $\CV$ be a quasi-projective model of $\Spec F$ over $\ZZ$, and $(\CU,f,\CL)$ be a polarized dynamical system over $\CV$ whose generic fiber is the polarized dynamical system $(X,f, L)$ over $\Spec F$.

\medskip\noindent\textit{Step 1.} 
We claim that there is an affine open subscheme $\CV'$ of $\CV$ such that the canonical map $\Pic(\CU_{\CV'})\to \Pic(X)$ is an isomorphism. This is a well-known fact, but we provide proof due to a lack of precise reference. 
\begin{enumerate}
\item There is an open subscheme $\CV'$ of $\CV$ such that $\CV'$ is regular and $\CU_{\CV'}\to \CV'$ has geometrically connected fibers. 
\item We can assume that $\Pic(\CV')$ is trivial by \cite[Chap. 2, Cor. 7.7]{Lan83}. 
Then $\Pic(\CV'')$ is trivial for any open subscheme $\CV''$ of $\CV'$ since $\Pic(\CV')\to \Pic(\CV'')$ is surjective by passing to Weil divisors, where the key is that $\CV'$ is regular. 
\item The canonical map $\Pic(\CU_{\CV''})\to \Pic(X)$ is injective for any open subscheme $\CV''$ of 
$\CV'$.
It suffices to prove that $\CaCl(\CU_{\CV''})\to \CaCl(X)$ is injective for the class groups of Cartier divisors. 
Then it suffices to prove that $\Cl(\CU_{\CV''})\to \Cl(X)$ is injective for the class groups of Weil divisors. 
If a Weil divisor of $\CU_{\CV''}$ is trivial on $X$, then it is vertical in the sense that it is the pull-back of a Weil divisor from $\CV''$, which is linearly equivalent to 0 by $\Pic(\CV'')=0$.
\item The canonical map $\ds\varinjlim_{\CV''}\Pic(\CU_{\CV''})\to \Pic(X)$ is an isomorphism by \cite[IV-3, Thm. 8.5.2]{EGA}. 
\item By (3) and (4), $\Pic(\CU_{\CV''})\to \Pic(X)$ is an isomorphism for sufficiently small open subscheme $\CV''$ of $\CV'$, since $\Pic(X)$ is finitely generated.
\end{enumerate}

Therefore, we can assume that the canonical map $\Pic(\CU)\to \Pic(X)$ is an isomorphism by replacing $\CV$ with a sufficiently small affine open subscheme.

Let $\pi:\CX\to \CS$ be a projective model of $\CU\to \CV$; i.e. $\CX$ and $\CS$ are projective models of $\CU$ and $\CV$ respectively, and $\CX\to \CS$ is a morphism extending $\CU\to \CV$. 
We can further assume that there is a strictly effective arithmetic divisor $\OCE_0$ on $\CS$, whose finite part has support equal to $\CS\setminus\CV$. 
Use the boundary divisor $(\CX,\pi^*\OCE_0)$ to define the boundary topology of $\wh\Div(\CU)_{ \QQ}$.

\medskip\noindent\textit{Step 2.} 
Consider the exact sequence
$$
0\lra \wh\Pic(\CU)_{\vert,\QQ} \lra \wh\Pic(\CU)_{\QQ} \lra \Pic(\CU)_{\QQ}\lra 0.
$$
Here $\wh\Pic(\CU)_{\vert,\QQ}$ is defined by the left exactness. 
For the right exactness, it suffices to prove that any effective Cartier divisor $\CD$ on $\CU$ can be extended to a projective model $\CX'$ of $\CU$.
This is easy by setting $\CX'$ to be the blowing-up of $\CX$ along the Zariski closure of $\CD$ in $\CX$.

Denote by $R(t)$ the characteristic polynomial of $f^*$ on the finite-dimensional vector space $\Pic(\CU)_{\QQ}=\Pic(X)_{\QQ}$.
We claim that 
$$R(f^*):\wh\Pic(\CU)_{\vert,\QQ}\lra \wh\Pic(\CU)_{\vert,\QQ}$$ 
is surjective.

Contrary to the proof of \cite[Thm. 4.9]{YZ1}, we do not use the interpretation of the metrics on Berkovich analytic spaces as in Proposition \ref{injection3} since it would be hard to control the convergence in terms of the boundary topology. 

Define $\wh\Div(\CU)_{\vert,\QQ}$ by the left exactness of
$$
0\lra \wh\Div(\CU)_{\vert,\QQ} \lra \wh\Div(\CU)_{\QQ} \lra \Div(\CU)_{\QQ}\lra 0.
$$
In terms of Proposition \ref{isomorphism}, there is a canonical surjection
$$
\wh\Div(\CU)_{\vert,\QQ}\lra \wh\Pic(\CU)_{\vert,\QQ}. 
$$
It suffices to prove that 
$$R(f^*):\wh\Div(\CU)_{\vert,\QQ}\lra \wh\Div(\CU)_{\vert,\QQ}$$ 
is surjective.

Take the Taylor expansion at $t=0$ by
$$
\frac{1}{R(t)}=\sum_{m=0}^\infty a_m t^m, \quad a_m\in \QQ.
$$
By Theorem \ref{semisimplicity}, the roots of the polynomial $R(t)$ have absolute values equal to $q$ or $q^{1/2}$.
Using partial fractions to expand $1/R(t)$,
there is a polynomial $Q(t)$ of rational coefficients such that 
$$
|a_m|\leq Q(m)q^{-m/2}, \ \forall m. 
$$
Denote 
$$
S_i(t)=\sum_{m=0}^i a_m t^m, \quad i\geq1.
$$

To prove the surjectivity, take any $\OCD\in \wh\Div(\CU)_{\vert,\QQ}$. We claim that 
the sequence $\{S_i(f^*)\OCD\}_i$ converges in $\wh\Div(\CU)_{\vert,\QQ}$. 
If so, then the limit gives an inverse image of $\OCD$ under $R(f^*)$. 

For the convergence, note that there is a positive rational constant $c$
such that 
$$-c\, \pi^*\OCE_0\leq \OCD\leq c\, \pi^*\OCE_0.$$
This holds automatically if $\OCD$ lies in the kernel of 
$\wh\Div(\CU)_{\rmod,\QQ} \to \Div(\CU)_{\QQ}$. 
In general, $\OCD$ is a limit of such elements, but then the Cauchy condition of $\OCD$ gives the constant $c$.

For any $i>j\geq 1$, we have 
$$
S_i(f^*)\OCD-S_j(f^*)\OCD
=\sum_{m=j+1}^i a_m\, (f^*)^m \OCD
\leq 
c\sum_{m=j+1}^i |a_m|\, \pi^*\OCE_0.
$$
We similarly have 
$$
S_i(f^*)\OCD-S_j(f^*)\OCD
\geq
-c\sum_{m=j+1}^i |a_m|\, \pi^*\OCE_0.
$$
By the bound of $a_m$, we see that $\{S_i(f^*)\OCD\}_i$ converges in $\wh\Div(\CU)_{\vert,\QQ}$.

\medskip\noindent\textit{Step 3.} 
The remaining part of the proof is almost identical to that of \cite[Thm. 4.9]{YZ1}. 
In fact, for any $M\in \Pic(\CU)_{\QQ}$, take any extension $\OM^0$ of $M$ in $\wh\Pic(\CU)_{\intb,\QQ}$, and set
$$
\OM_f= \OM^0-R(f^*)|_{\wh\Pic(\CU)_{\vert,\QQ}}^{-1} (R(f^*)\OM^0).
$$

The proof of the nefness of $\OM_f$ under the ampleness of $M$ on $\CU$, though lengthy, is similar to that in \cite[Thm. 4.9]{YZ1}, so we omit it. 
\end{proof}

\section{N\'eron-Tate height on a curve}  \label{section NT}

When $X$ is a projective curve over a finitely generated field, we present a theorem (Theorem \ref{thm NT}) which interprets the intersection numbers in terms of the N\'eron--Tate height. It generalizes the Hodge index theorem of Faltings \cite{Fal1} and Hriljac \cite{Hri} (cf. Theorem \ref{app hodge index1}) to finitely generated fields.

\subsection{The arithmetic Hodge index theorem}

\kkk
\ccc
Let $F$ be a finitely generated field over $k$, and let $\pi: X\to \Spec F$ be a smooth,  projective, and geometrically connected curve of genus $g>0$. 
We first introduce the canonical height function 
$$
\hat\fh: \Pic^0(X_{\overline F}) \lra \wh\Pic(F/k)_{\QQ,\nef}.
$$

Denote by $J=\underline{\Pic}^0_{X/F}$ the Jacobian variety of $X$. 
Denote by $\Theta$ the symmetric line bundle on $J$ associated with the theta divisor. Namely, choose a point $x_0\in X(\overline F)$ and denote by $j: X_{\overline F}\hookrightarrow J_{\overline F}$ the embedding $x\mapsto [x-x_0]$. 
Denote by $\theta$ the image of the composition $X_{\overline F}^{g-1}\hookrightarrow J_{\overline F}^{g-1}\to J_{\overline F}$. The second map is the sum under the group law. Then $\theta$ is a divisor of $J_{\overline F}$. Denote by $\Theta$ the line bundle on $J_{\overline F}$ associated to $\theta+[-1]^*\theta$. 
The isomorphism class of $\Theta$ does not depend on the choice of $x_0$, so it is Galois invariant and descends to a line bundle on $J$. 
See \cite[\S5.6]{Ser} for more details about the construction.

By the symmetric and ample line bundle $\Theta$ on $J$, we have the canonical height
$$
\hat\fh_{\Theta}: J(\overline F) \lra \wh\Pic(F/k)_{\QQ,\nef}.
$$
By convention, we set 
$$\ds\hat\fh=\frac12\hat\fh_{\Theta}.$$

The goal of this section is to prove the following extension of the arithmetic Hodge index theorem of Faltings \cite{Fal1} and Hriljac \cite{Hri} to finitely generated fields.

\begin{thm}[arithmetic Hodge index theorem] \label{thm NT}
\kkk
Let $F$ be a finitely generated field over $k$, and let $\pi: X\to \Spec F$ be a smooth,  projective, and geometrically connected curve.
Let $M$ be a line bundle on $X$ with $\deg M=0$. 
Then there is an adelic line bundle $\overline M_0\in \wh\Pic(X/k)_{\intb,\QQ}$ with underlying line bundle $M$ such that 
$$\pi_*\pair{\overline M_0,\overline V}= 0,\quad\forall\,\overline V\in \wh\Pic(X/k)_{\vert,\QQ}.$$ 
Moreover, for such an adelic line bundle,
$$\pi_*\pair{\overline M_0, \overline M_0}= -2\, \wh\fh(M)$$
in $\wh\Pic(F/k)_{\QQ}$.
\end{thm}

In the theorem, $\pi_*\pair{\cdot,\cdot}$ denotes the Deligne pairing 
$$\wh\Pic(X/k)_{\intb,\QQ} \times \wh\Pic(X/k)_{\intb,\QQ} \lra \wh\Pic(F/k)_{\intb,\QQ}$$
introduced in Theorem \ref{intersection2}. 
And
$\wh\Pic(X/k)_{\vert,\QQ}$ denotes the kernel of the forgetful map $\wh\Pic(X/k)_{\intb,\QQ}\to \Pic(X)_\QQ$.

If we fix a polarization of $F/k$ and intersect both sides of the equality with the polarization, then we obtain equality about the Moriwaki heights. Moriwaki \ proved this cite[Thm. B]{Mor2}. 

\begin{rmk}
We will see in \cite{YZ2} that the extension $\overline M_0$ is unique up to translation by $\pi^*\wh \Pic(F/k)_\intb$.
\end{rmk}

\subsection{The universal adelic line bundle}

Now  we construct the extension $\overline M_0$ in Theorem \ref{thm NT}. It is written almost the same as the number field case. We include it here briefly. For basic geometric results on abelian varieties and Jacobian varieties, we refer to Mumford \cite{Mum1} and Serre \cite{Ser}. 

Denote by $p_1: X\times J\to X$ and $p_2: X\times J\to J$ the projections. 
Via $p_1$, we view $X\times J$ as an abelian scheme on $X$. 
Denote by $[m]_X:X\times J\to X\times J$ the multiplication by an integer $m$ 
as abelian schemes on $X$, i.e. the map sending $(x,y)$ to $(x,my)$.

We claim that there is a universal line bundle $Q\in \Pic(X\times J)_\QQ$ satisfying the following properties:
\begin{enumerate}
\item For any $\alpha\in J(\overline F)$, the $\QQ$-line bundle $Q|_{X \times \alpha}$ on $X \times \alpha= X_{\overline F}$ is equal to $\alpha$ in $\Pic^0(X_{\overline F})_\QQ$.
\item For any integer $m$, $[m]_X^*Q=mQ$ in $\Pic(X\times J)_\QQ$.  
\end{enumerate}
The line bundle $Q$ is unique up to translation by $p_2^*\Pic^0(J)_\QQ$. 

Let $\alpha_0$ be a line bundle on $X$ of degree $d>0$.
Denote the canonical morphism 
$$i_0:X\lra J,\quad x\longmapsto dx-\alpha_0.$$ 
Denote by 
$$(i_0,\id):X\times J\lra J\times J$$ 
the natural morphism. 
Set 
$$Q=\frac{1}{d}(i_0,\id)^*P,$$ 
where $P$ is the Poincar\'e line bundle on $J\times J$. 

If there is a line bundle on $X$ of degree $1$, we can choose $Q$ to be an integral line bundle on $X\times J$. 
If $X(F)$ is non-empty, take $x_0\in X(F)$ and use it to define $i_0:X\to J$.
Then $Q$ is an integral line bundle on $X\times J$ such that $Q_{x_0\times J}=0$ and that for any $\alpha\in J(\overline F)$, the line bundle $Q|_{X \times \alpha}$ on $X \times \alpha= X_{\overline F}$ is equal to $\alpha$ in $\Pic^0(X_{\overline F})$.
These properties determine $Q$ uniquely. 

With the universal line bundle $Q\in \Pic(X\times J)_\QQ$, by Theorem \ref{invariant metric3}, there is a unique extension $\overline Q\in \wh\Pic(X\times_F J/k)_{\intb,\QQ}$ of $Q$ such that $[2]_X^*\overline Q=2\overline Q$. 

Let $\alpha$ be the point of $J(F)$ represented by the line bundle $M\in \Pic^0(X)$. 
Set 
$$\overline M_0:= \overline Q|_{X\times \alpha}\in \wh\Pic(X/k)_{\intb,\QQ}.$$ 
We need to prove that $\overline M_0$ satisfies the requirement of Theorem \ref{thm NT}; i.e. 
$$\pi_*\pair{\overline M_0, \overline V}= 0,\quad\forall\,\overline V\in \wh\Pic(X/k)_{\vert,\QQ}.$$ 
Consider the adelic line bundle
$$
\overline R:=p_{2,*}\pair{\overline Q, p_1^* \overline V}
$$
in $\wh\Pic(J/k)_{\intb,\QQ}$. Note that $\overline R$ is universal in the sense that the pull-back of $\overline R$
via $\alpha:\Spec(F) \to J$ is exactly $\pi_*\pair{\overline M_0, \overline V}$. 
Thus it suffices to prove that the adelic line bundle $\overline R=0$ in $\wh\Pic(J/k)_{\intb,\QQ}$.

This is a consequence of Theorem \ref{invariant metric2} by noting the following two properties:
\begin{enumerate}
\item  the underlying line bundle $R=0$ in $\Pic(X\times J)_\QQ$, as a consequence of the underlying line bundle $V=0$;
\item  $[2]_X^*\overline R=2\overline R$ in $\wh\Pic(J/k)_{\QQ}$ by the dynamical property of $\overline Q$;
\end{enumerate}

\subsection{The height equality}
It remains to prove 
$$\pi_*\pair{\overline M_0, \overline M_0}= -2\, \wh\fh(M).$$
Replacing the field $F$ by a finite extension if necessary, we can assume that $X(F)$ is non-empty. 
We first express the left-hand side as a height function. 

Take $x_0\in X(F)$. Use $x_0$ to define the canonical embedding $i_0:X\to J$, and identity $X$ as a subvariety of $J$.
As before, let $Q$ be the restriction of the Poincar\'e line bundle $P$ from $J\times J$ to $X\times J$. 

Note that $P$ is symmetric on $J\times J$. 
Thus $[2]^*P=4P$ and we can extend it to $[2]^*\OP=4\OP$ for some $\overline P\in\wh\Pic(J\times_F J/k)_{\intb}$ by Theorem \ref{invariant metric3}.
We claim that $\overline Q=\overline P|_{X\times J}$ in $\wh\Pic(X\times_F J/k)_{\QQ}$.

Note that $[2]: J\times J\to J\times J$ is multiplication by two on both components, while $[2]_X: X\times J\to X\times J$ is only the multiplication by two on the second component. 
Denote by $[2]_2:J\times J\to J\times J$ the multiplication by two on the second component.
By Theorem \ref{invariant metric2}(2), $[2]^*\OP=4\OP$ implies $[2]_2^*\OP=2\OP$.
This argument was used in the proof of Theorem \ref{invariant metric3}.
This implies $\overline Q=\overline P|_{X\times J}$ by the uniqueness of $\OQ$ in Theorem \ref{invariant metric2}(1). 
All these qualities are viewed as isomorphism classes of adelic $\QQ$-line bundles.

\begin{lem}
For any $\alpha,\beta\in J(F)$, we have 
$$
\pi_*\pair{\overline P_\alpha, \overline P_\beta} = \fh_{\overline P}(\alpha,\beta).
$$ 
Here $\overline P_\alpha=\overline P|_{X\times \alpha}$ and $\overline P_\beta=\overline P|_{X\times \beta}$ are viewed as adelic line bundles on $X$. 
\end{lem}
\begin{proof}
Note both sides are bilinear in $(\alpha,\beta)$. 
We can assume that $\alpha$ represents the divisor $x-x_0$ on $X$. 
Then $\alpha=j(x)$. 
Here we assume $x\in X(F)$ by replacing $F$ with a finite extension if necessary. 
Then we have 
$$
\pi_*\pair{\overline P_\alpha, \overline P_\beta}
=\pi_*\pair{\hat x-\hat x_0, \overline P_{\beta}}
=\pi_*\pair{\hat x, \overline P_{\beta}}-\pi_*\pair{\hat x_0, \overline P_{\beta}}.
$$
Here $\hat x$ and $\hat x_0$ are any extensions of $x$ and $x_0$ in $\wh\Pic(X/k)_{\intb,\QQ}$. 
Note that $\overline P_{\beta}$ has zero intersection with any vertical classes. The above becomes 
$$
\pi_*(\overline P|_{x\times \beta}) 
-\pi_*(\overline P|_{x_0\times \beta}) 
=\pi_*(\overline P|_{x\times \beta})
=\fh_{\overline P}(\alpha,\beta).
$$
\end{proof}

Now  we are ready to prove
$$\pi_*\pair{\overline M_0,\overline M_0}= -2\, \wh\fh(M).$$
By the lemma, it suffices to prove 
$$\fh_{\overline P}(\alpha,\alpha)
=-\fh_{\overline \Theta}(\alpha),\quad \forall\ \alpha\in J(F).$$ 
It is well known that the Poincar\'e bundle on $J\times J$ has the expression 
$$
P=   p_1^* \theta + p_2^*\theta-m^* \theta. 
$$
Here $ m,  p_1,  p_2: J\times J\to J$
denotes the addition law and the projections. 
It induces
$$
2P=   p_1^* \Theta + p_2^*\Theta-m^* \Theta. 
$$
We use $\Theta$ because it is also symmetric. It follows that
$$
2\overline P=   p_1^* \overline\Theta + p_2^*\overline\Theta-m^*\overline \Theta. 
$$
Computing heights using the identity we have 
$$
2\fh_{\overline P}(\alpha,\alpha)
=\fh_{\overline \Theta}(\alpha)
+\fh_{\overline \Theta}(\alpha)
-\fh_{\overline \Theta}(2\alpha)
=-2\fh_{\overline \Theta}(\alpha).
$$
This finishes the proof of Theorem \ref{thm NT}.

\subsection{High-dimensional bases}

The above setting treats $X\to \Spec F$ for a finitely generated field $F$ over $k$. We can replace $\Spec F$ by an essentially quasi-projective scheme $S$ over $k$. Still, due to the flatness problem, we have to restrict the vector-valued height of sections of the relative Jacobian scheme. 

\kkk
Let $S$ be a normal integral scheme, flat and essentially quasi-projective over $k$.
Let $\pi: X\to S$ be a projective and smooth morphism whose fibers are smooth and geometrically connected curves.
Denote by $\Pic^0(X/S)$ the group of line bundles on $X$ with degree 0 on the fibers of $X\to S$. 
We first introduce a canonical height function 
$$
\hat\fh: \Pic^0(X/S) \lra \wh\Pic(S/k)_{\QQ,\nef}.
$$
This is obtained as a slight generalization of the case $S=\Spec F$ and is thus compatible with the latter.

Denote by $J=\Pic^0_{X/S}$ the Jacobian scheme of $X$ over $S$.
For the basics of Jacobian schemes, we refer to \cite[Chap. 6]{MFK}.
By \cite[\S6.1, Prop. 6.9]{MFK}, there is a canonical principal polarization 
$\lambda_1:J\to J^\vee$ over $S$. 
By the construction of \cite[\S6.2, Prop. 6.10]{MFK}, there is a symmetric line bundle 
$\Theta$ on $J$ such that the polarization $\lambda_\Theta:J\to J^\vee$ corresponding to $\Theta$ is exactly twice of $\lambda_1:J\to J^\vee$. 
To relate it to our previous case of fields, $\Theta$ recovers that on each fiber of $J\to S$. 
We can uniquely determine $\Theta$ by the rigidification $e^*\Theta\simeq \CO_S$ for the identity section $e: S\to J$.

Finally, by the symmetric and relatively ample line bundle $\Theta$ on $J$, we have 
a unique extension $\overline \Theta$ of $\Theta$ in $\wh\Pic(J/k)_{\QQ,\nef}$ such that $[2]^*\overline \Theta=4\overline \Theta$. 
Then we have the vector-valued height function
$$
\hat\fh_{\overline \Theta}: J(S) \lra \wh\Pic(S/k)_{\QQ,\nef}.
$$
By convention, we set 
$$\ds\hat\fh=\frac12\hat\fh_{\overline \Theta}.$$
By the canonical map $\Pic^0(X/S)\to J(S)$, we obtain
$$
\hat\fh: \Pic^0(X/S) \lra \wh\Pic(S/k)_{\QQ,\nef}.
$$
As in the classical case, this height function is also quadratic.

As before, we have a universal $\QQ$-line bundle $Q\in \Pic(X\times_S J)_\QQ$ satisfying the following properties:
\begin{enumerate}
\item For any base change $X'\to S'$ of $X\to S$, and for any $\alpha\in \Pic(X')$ of degree 0 on fibers of $X'\to S'$, 
 the pull-back of the $\QQ$-line bundle $Q$ via
 $(\id,\alpha):X \times_S S'\to X\times_SJ$ is equal to $\alpha$ in $\Pic(X')_\QQ$.
\item For any integer $m$, $[m]_X^*Q=mQ$ in $\Pic(X\times_S J)_\QQ$.  
\end{enumerate}
The line bundle $Q$ is unique up to translation by $p_2^*\Pic(J)_\QQ$. 
There is a unique extension $\overline Q\in \wh\Pic(X\times_S J/k)_{\intb,\QQ}$ of $Q$ such that $[2]_X^*\overline Q=2\overline Q$. 

The following is a variant of Theorem \ref{thm NT} over high-dimensional bases.

\begin{thm} \label{thm NT2}
\kkk
Let $S$ be a normal integral scheme, flat and essentially quasi-projective over $k$.
Let $\pi: X\to S$ be a projective and smooth morphism of relative dimension 1 with geometrically connected fibers.
Let $M$ be a line bundle on $X$ with degree 0 on fibers of $X\to S$.
Then there is an adelic line bundle $\overline M_0\in \wh\Pic(X/k)_{\intb,\QQ}$ with underlying line bundle $M$ such that 
$$\pi_*\pair{\overline M_0,\overline V}= 0,\quad\forall\,\overline V\in \wh\Pic(X/k)_{\vert,\QQ}.$$ 
Moreover, for such an adelic line bundle,
$$\pi_*\pair{\overline M_0, \overline M_0}= -2\, \wh\fh(M)$$
in $\wh\Pic(S/k)_{\QQ}$.
\end{thm}

\begin{proof}
The existence of $\OM_0$ can be obtained by generalizing the construction by the universal line bundle $\overline Q$ in Theorem \ref{thm NT} to the general base $S$. 
Namely, let $i:S\to J$ be the morphism corresponding to $[M]\in J(S)$, i.e.
by the morphism $(\id,i\circ\pi):X\to  X\times_S J$, the pull-back
$(\id,i\circ\pi)^*Q= M+ \pi^*N$ for some $N\in\Pic(S)_\QQ$. 
Then we set 
$$\OM_0=(\id,i\circ\pi)^*\overline Q- \pi^*\ON$$
for any $\ON\in\Pic(S/k)_{\intb,\QQ}$ extending $N$.
For the height identities are consequences of Theorem \ref{thm NT}, as $\wh\Pic(S/k)_{\QQ}\to \wh\Pic(F/k)_{\QQ}$ is injective by Corollary \ref{shrink}. 
Here $F=k(S)$ is the function field. 
\end{proof}

The following universal Hodge index theorem is essentially equivalent to  
Theorem \ref{thm NT2}. 

\begin{cor} \label{thm NT3}
By the second projection $p_2:X\times_SJ\to J$,
$$p_{2*}\pair{\overline Q, \overline Q}= -\overline \Theta$$
in $\wh\Pic(J/k)_{\QQ}$.
\end{cor}

\begin{proof}
Let $S'$ be a normal integral scheme, flat and essentially quasi-projective over $k$, endowed with a $k$-morphism $S'\to S$.
Denote by $\pi':X'\to S'$ the base change of $\pi:X\to S$ by $S'\to S$. 
Let $i: S'\to J$ be a morphism over $S$. 
Consider the morphism $(\id,i):X\times_S S'\to  X\times_S J$.
Apply Theorem \ref{thm NT2} to $\pi':X'\to S'$ and the line bundle $M=(\id,i)^*Q$ in $\Pic(X')_\QQ$.
We obtain
$$\pi_{*}'\pair{(\id,i)^*\overline Q, (\id,i)^*\overline Q}= -i^*\overline \Theta.$$
Set $S'=J$ and set $i: S'\to J$ to be the identity morphism. 
\end{proof}

\subsection{High-dimensional fibers}

For completeness, we state an arithmetic Hodge index theorem for adelic line bundles on projective varieties over finitely generated fields proved by Yuan--Zhang \cite{YZ1, YZ2} and Carney \cite{Car1, Car2}. 
It can be viewed as a generalization of Theorem \ref{thm NT} from projective curves to projective varieties.
We refer to Theorem \ref{app hodge index3} for a statement of the main theorem of  Yuan--Zhang \cite{YZ1}.

To state the theorem, we start with the following general positivity notions. 

\begin{defn} \label{def-positivity0} 
\kkk
Let $X$ be a flat and essentially quasi-projective integral scheme over $k$. 
Denote by $d$ the dimension of projective models of $X$ over $k$, and assume $d\geq 1$. 
Let $\overline L, \OM$ be adelic line bundles or adelic $\QQ$-line bundles on $X/k$. 
\begin{enumerate}
\item We say that $\overline L$ is \emph{numerically trivial} if 
 $\overline L\cdot \overline N_1\cdots \overline N_{d-1}= 0$ for any $\overline N_1, \cdots, \overline N_{d-1}$ in $\wh \Pic (X/k)_{\intb}$. 

\item We say that $\OM$ is \emph{$\overline L$-bounded} if  there is a rational number $\epsilon >0$ such that both 
$\overline L+ \epsilon \OM$ and $\overline L- \epsilon \OM$ are nef. 
\end{enumerate}
\end{defn}

To compare with the property ``numerically trivial'', we recall the property ``pseudo-effective'' in \S\ref{subsec pseudo}. 
In particular, if $\overline L$ is pseudo-effective, 
the top intersection number $\overline L\cdot \overline N_1\cdots \overline N_{d-1}\geq 0$ for any $\overline N_1, \cdots, \overline N_{d-1}$ in $\wh \Pic (X/k)_{\nef}$. 
Moreover, the BDPP criterion in Theorem \ref{BDPP criterion} gives an inverse of this statement by passing to modifications of $X$. 

The following notion is specific to our setting of the arithmetic Hodge index theorem. 

\begin{defn} \label{def-positivity1} 
\kkk
Let $F$ be a finitely generated field over $k$. 
Let $X$ be a geometrically integral projective variety over $F$.
Let $\overline L$ be an adelic line bundle or an adelic $\QQ$-line bundle on $X/k$.
Define the notion $\overline L\gg 0$ in the following two cases.
\begin{enumerate} 
\item 
If $k=\ZZ$, we write $\overline L\gg 0$ if $L$ is ample, and $\overline L- \overline N$ is nef for some $\overline N\in \wh\Pic(\ZZ)$ with $\wh\deg (\overline N)>0$. 
Here the adelic line bundle $\overline N$ is viewed as an element of $\wh \Pic (X/k)_{\intb}$ by the natural pull-back map. 
\item 
If $k$ is a field and $F/k$ has transcendence degree $d\geq 1$, we write $\overline L\gg 0$ if $L$ is ample, and $\overline L- \overline N$ is nef for some $\overline N\in \wh\Pic(k_1/k)$ with $\wh\deg (\overline N)>0$. 
Here $k_1/k$ is an intermediate extension of $F/k$ of transcendence degree $1$, and 
the adelic line bundle $\overline N$ is viewed as an element of $\wh \Pic (X/k)_{\intb}$ by the natural pull-back map. 
\end{enumerate}
\end{defn}

Finally, the arithmetic Hodge index theorem is as follows.

\begin{thm}[\cite{YZ1,YZ2, Car1,Car2}] \label{hodge high dim} 
\kkk
Let $F$ be a finitely generated field over $k$. 
If $k$ is a field, assume that 
the transcendence degree of $F$ over $k$ is at least one and that $k$ is algebraically closed in $F$. 
Let $\pi: X\to \Spec F$ be a geometrically integral and geometrically normal projective variety of dimension $n\geq 1$. 
Let $\overline M$ be an integrable adelic $\QQ$-line bundle on $X$, and $\overline L_1, \cdots, \overline L_{n-1}$ be $n-1$ nef adelic $\QQ$-line bundles on $X$ where each $L_i$ is big on $X$. 
Assume $M\cdot L_1\cdots L_{n-1}=0$ on $X$. 
Then 
$$-\pi_*\langle\overline M, \OM,  \overline L_1, \cdots,  \overline L_{n-1}\rangle
$$
is pseudo-effective
in $\wh \Pic(F/k)_{\intb,\QQ}$.

Moreover, if $\overline L_i\gg 0$, and $\overline M$ is $\overline L_i$-bounded for each $i$, then
$$\pi_*\langle\overline M, \OM,  \overline L_1, \cdots,  \overline L_{n-1}\rangle$$
is numerically trivial in $\wh \Pic(F/k)_{\intb,\QQ}$
if and only if one of the following two cases holds: 
\begin{enumerate}
\item 
 $k=\ZZ$ and $\overline M\in \pi^* \wh \Pic(F/k)_{\intb,\QQ}$;
\item 
 $k$ is a field and $\overline M\in \pi^* \wh \Pic(F/k)_{\intb,\QQ}+(\tr_{K/k} \Pic^0(X))_\QQ$, 
 where $\tr_{K/k} \Pic^0(X)$ denotes the image of the natural composition
 $$
 \Pic^0(A_0)\lra \Pic(A_0)\stackrel{\sim}{\lra} \wh\Pic(A_0/k) \lra \wh\Pic(A/k) \lra \wh\Pic(X/k).
 $$
Here 
$A=(\underline{\Pic}^0_{X/F,\red})^\vee$ denotes the Albanese variety of $X$ over $F$, 
and $A_0$ denotes Chow's $(F/k)$-image of $A$.  
The last two arrows are the pull-back maps via the canonical map $A\to A_0$ and an Albanese map $X\to A$. 
\end{enumerate}
\end{thm} 

If $k=\ZZ$ and $F$ is a number field, the theorem was proved in \cite{YZ1} (cf. Theorem \ref{app hodge index3});
if $k=\ZZ$ and $F$ is general, the theorem was proved in \cite{YZ2};
if $k$ is a field and $F$ has transcendence degree one over $k$, the theorem was proved in \cite{Car1}; If $k$ is a field and $F$ is general, the theorem was proved in \cite{Car2}. 
Note that the cases of \cite{YZ1, Car1} require only the original theory of adelic line bundles of Zhang \cite{Zha2}. 
We refer to the original references for proofs of the theorem and the applications to algebraic dynamical systems.

\cleardoublepage
\appendix

\chapter[Heights and intersection theory]{Review on heights and arithmetic intersection theory}\label{ch-htint}



One major motivation of Arakelov geometry is to apply geometric ideas, usually inspired by algebraic geometry and complex geometry, to study Diophantine problems.
In this consideration, heights of algebraic points or subvarieties on projective varieties are interpreted as arithmetic intersection numbers. 

In this appendix, we review some major theories and important theorems of Arakelov geometry with an emphasis on applications to heights.
The first intention here is to provide background materials for our new theory of adelic line bundles on quasi-projective varieties. 
The second intention here is to provide a historical account of Arakelov geometry for readers who are not familiar with the theory. 
Due to limited space here, our exposition is very sketchy;
due to the limited ability of the authors, our exposition might miss many important parts of the theory;
due to various different mathematics settings, 
the terminology of this appendix might be slightly different from the main body of this book.
Our exposition will roughly follow the chronological order of materials appearing in history, with minor adjustments according to mathematics subjects. 
Roughly speaking,  \S\ref{app sec height}, \S\ref{app sec surface}, and \S\ref{app sec intersection} give historical accounts of definitions of heights and arithmetic intersection numbers;
\S\ref{app sec adelic}, \S\ref{app sec positivity}, and \S\ref{app sec equidistribution}
give an overview of the more recent theory of arithmetic positivity and its applications to heights of algebraic points.

In particular, \S\ref{app sec adelic} reviews Zhang's old theory of adelic line bundles on projective varieties. 
This is used in \S\ref{sec restricted} for a comparison of the old theory with our new theory, which will be helpful to understand the new theory. 
Because of this, \S\ref{app sec adelic} is written with more rigourous definitions. 
 
Most materials of this appendix, if used in the main body of this book, will be reviewed or quoted in the relevant part of the main body.
For readers unfamiliar with Arakelov geometry, it might be helpful to read this appendix before reading the main body of this book; 
for readers familiar with Arakelov geometry, it is natural to read the main body without reading through this appendix.

For a first course in Arakelov geometry, we refer to the textbook of Yuan--Guo \cite{YG25}.

\section{Heights and degrees}\label{app sec height}

In this section, we review Weil heights for projective varieties,  
N\'eron--Tate heights for abelian varieties,  Call--Silverman heights for polarized dynamical systems, and Faltings heights of abelian varieties. We also 
introduce their relations to  arithmetic degrees of hermitian line bundles on arithmetic curves. Our basic references are  Lang \cite{Lan83}, Serre \cite{Ser},  Bombieri--Gubler \cite{BG06}, and Szpiro \cite[\S1, \S3]{Szp85}. For updated details, see  Yuan--Guo \cite[\S 2, \S3, \S4]{YG25}.

\subsection{Naive heights}

Let us start with heights of algebraic points of projective varieties originally introduced by Weil \cite{Wei51}.

For each place $p$ of $\QQ$ in $M_\QQ= \{\infty, 2, 3, 5, 7,\cdots\}$, let $\ol\BQ_p$ denote an algebraic closure of the local field $\BQ_p$.  
Endow $\ol\BQ_p$ with the $p$-adic norm $|\cdot|_p$ so that $|p|_p=1/p$ for $p\ne \infty$ and that $|\cdot|_\infty$ is the usual absolute induced from $\BC=\ol\BQ_\infty$. 

For a point $x=[x_0, \cdots, x_n]\in \BP^n (\ol\BQ)$, the naive height of $x$ is defined by 
$$h_\naive(x):=\frac 1{[K:\BQ]}\sum _{p\in M_\QQ} \sum _{\sigma: K\to \ol \BQ_p}
\log \max\{ |\sigma (x_0)|_p,\cdots, |\sigma (x_n)|_p\},$$
where $K$ is a finite extension of $\BQ$ containing all $x_i$, $p$ runs through the set of places of $\BQ$, and $\sigma$ runs through the set of embeddings from $K$ to $\ol \BQ_p$. This definition does not depend on the choice of the homogenous coordinate $[x_0, \cdots, x_n]$ and the field $K$.
We can also express this height in terms of places $v$ of $K$ by
$$h_\naive(x)=\frac 1{[K:\BQ]} \sum _{v\in M_K}
\log \max\{ |x_0|_v,\cdots, |x_n|_v\},$$
where $M_K$ denotes the set of places of $K$, and  
the $v$-adic norm is normalized by $|x|_v=|\sigma (x)|_p^{[K_v: \BQ_p]}$ for a place $v$ of $K$ above a place $p$ of $\QQ$ and for any continuous homomorphism $\sigma: K_v\to \ol \BQ_p$.

In the simplest case $x\in \BP^n (\BQ)$, taking $x_0,\cdots, x_n$ to be coprime integers,  we have 
$$h_\naive(x)=\log\max\{ |x_0|_\infty,\cdots, |x_n|_\infty\}.$$
Then $h_\naive (x)$ is essentially $1/(\log 10)$ times the maximum of the numbers of digits of $|x_i|_\infty$. Thus the naive height function is indeed a function to measure the arithmetic complexity of rational points in $\BP^n$.

\subsection{Height machine}
As a convention, we write tensor products of line bundles additively throughout the book. 

The above naive height function induces height functions on every closed subvariety $X$ of $\BP^n$ defined over a number field $K$. 
Up to a bounded function on $X(\ol K)$, this height function depends only on the line bundle $L:=\CO(1)|_X$, not on the choice of sections of $L$ defining the embedding.
 Thus, we have a well-defined map $L\mapsto h_L$ from the set of very ample line bundles to the space of functions on $X(\ol K)$ 
 modulo the subspace of bounded functions. There are two important properties: this map is additive and injective.
Hence, we can extend this map to the vector space $\Pic (X)_\RR=\Pic (X)\otimes_\ZZ\RR$ to obtain a homomorphism
  $$\Pic (X)_\RR\lra \frac{\{\text{$\BR$-valued functions on $X(\ol K)$}\}}{\{\text{Bounded $\BR$-valued functions on $X(\ol K)$}\}},\qquad L\longmapsto h_L.$$
This homomorphism is usually called the height machine, and every function $h_L:X(\OK)\to \RR$ representing the image of $L\in \Pic (X)_\RR$ under the homomorphism is called a Weil height function associated to $L$. 
 
In 1949, Northcott \cite{Nor49} proved a beautiful finiteness property for any Weil height function $h_L$ associated to an ample line bundle $L$ on $X$. 
Namely,  for any positive integers $D$ and $H$, the set of points $x\in X(\ol K)$ with 
$\deg x\le D$ and $h_L(x)\le H$ is finite. As a consequence, he
proved  a weak Mordell--Weil property for an endomorphism $f: X\to X$ polarized by an ample line bundle $L$  in the sense that $f^*L\simeq qL$ for some integer $q>1$.
We call a point $x\in X(\ol K)$ preperiodic if there are integers $m>n\ge 0$  such that $f^m(x)=f^n(x)$.
Then the weak Mordell--Weil property of Northcott \cite{Nor50} asserts that the set of preperiodic points in $X(K)$ is finite.

\subsection{N\'eron--Tate heights and Call--Silverman heights}
In the height machine, for a fixed line bundle $L$, the Weil height $h_L:X(\OK)\to \RR$ is only unique up to bounded functions, and there is usually no natural choice of $h_L$. 
However, in the 1960s, when $X=A$ was an abelian variety defined over a number field $K$, a normalization of the height machine is found by N\'eron in \cite{Ner65}  and by Tate in an unpublished note.
 N\'eron's definition is by a decomposition into a sum of normalized local heights,
  and Tate's definition is a simple limiting argument as follows. 
  
First of all, we have a decomposition 
  $$\Pic (A)_\BQ=\Pic (A)_\BQ^+\oplus \Pic(A) _\BQ ^-$$
 into $(\pm 1)$-eigenspaces according to the action by the pull-back via $[-1]:A\to A$.
Then for any element $L^+ \in \Pic (A)_\BQ^+$ (resp. $L^- \in \Pic (A)_\BQ^-$) with a Weil height function $h_{L^+}$ (resp. $h_{L^-}$), we define a normalized height in the class of $h_{L^+}$ (resp. $h_{L^-}$) 
  by
  $$\wh h_{L^+}(x):=\lim _{n\to \infty} \frac {h(2^n x)}{4^n }, \qquad \wh h_{L^-}(x):=\lim _{n\to \infty} \frac {h(2^n x)}{2^n }.$$
  This is Tate's limiting argument.
  It can be shown that $\wh h_{L^+}$ is quadratic on $A(\ol K)$ and $\wh h_{L^-}$ is linear on $A(\ol K)$. In particular, they both vanish on the subgroup $A(\ol K)_\tor$ of torsion points.
  Moreover, if $L^+$ is ample, then $\wh h_{L^+}$ is positive definite on $A(\ol K)/A(\ol K)_\tor$. The normalized height is now called the canonical height or the N\'eron--Tate height. It is an important ingredient for the formulation of the refined Birch and Swinnerton-Dyer conjecture by Tate \cite{Tat66}.

  In 1993, Call--Silverman \cite{CS93} generalized  Tate's limiting argument to any endomorphism $f: X\to X$ on a projective variety $X$ over a number field $K$ and any line bundle $L$ on $X$ so that $f^*L\simeq qL $ for some integer $q>1$. In this case, they defined the canonical height by
  $$\wh h_{L,f}(x):=\lim _{n\to \infty} \frac {h(f^n (x))}{q^n }, \qquad x\in X(\ol K).$$
 
 Moreover, if $L$ is ample (and thus gives a polarization of $f$), then we can do a similar construction to
  give a canonical height function $\wh h_M:X(\ol K)\to \RR$ for each $M\in \Pic (X)$ so that 
$$  \wh h_M (f (x))=\wh h_{f^*M}(x), \qquad \forall\, M \in \Pic (X), \quad \forall\, x\in X(\ol K).
$$
In fact, the polarization $L$  implies that the action $f^*:\Pic(X)_\BC\to \Pic(X)_\BC$ is semisimple and all of its eigenvalues have absolute values $q^{1/2}$ or $q$.
  Thus, we may apply Tate's limiting argument to eigenvectors in $\Pic (X)_\BC$.
  For more details, see Yuan--Zhang \cite{YZ1}.
  
  \subsection{Hermitian line bundles on arithmetic curves}
  
  Let $K$ be a number field with the ring $O_K$ of integers.
  Then we have an arithmetic curve $\Spec O_K$. 
  This curve can be compactified by adding the set of archimedean places $v\mid \infty$. Let $\CL$ be a line bundle  on $\Spec O_K$, which can be identified as a locally free $O_K$-module of rank 1. The line bundle $\CL$  can be compactified by adding a $v$-adic norm on $\CL_v=\CL\otimes _{O_K} K_v$ for each archimedean place $v$, i.e. a nonzero map
$$\|\cdot\|_v: \CL_v\lra \BR_{\geq 0}$$
satisfying
$$\|a\ell\|_v=|a|_v \|\ell\|_v, \qquad \forall\, a\in K_v, \, \ell \in \CL_v.$$
We call the pair $\ol \CL=(\CL, \|\cdot\|)$ a hermitian line bundle on $\Spec O_K$.
The arithmetic degree of $\CLL$ is defined by
$$\wh\deg(\ol\CL):=\log  \# (\CL/\ell O_K)-\sum_{v\mid \infty }\log\|\ell\|_v,$$
where $\ell\in \CL$ is a non-zero section. 

Denote by $\wh\Pic (O_K)$  the isometry group of hermitian line bundles on $\Spec O_K$.
Taking arithmetic degrees gives a homomorphism $\wh\deg: \wh \Pic (O_K)\to \BR$.
Then we have an exact sequence
$$0\lra (\BR ^r)^0/\log O_K^\times\lra \wh\Pic(O_K)\lra \Pic (O_K)\times \BR\lra 0,$$
where $r$ is the number of archimedean places of $K$,  
$(\BR ^r)^0$ is the space of $a=(a_v)\in \BR ^r$ with $\sum_v a_v=0$,
the first nonzero map
sends $a=(a_v)\in \BR ^r$ to  the hermitian line bundle $\CO(a)=(O_K, \|\cdot\| )$ 
with metric given by $\|1\|_v=e^{-a_v}$, the second nonzero map is the combination of the forgetful map and the  degree map.

\subsection{Faltings heights}
As one application, we introduce the modular heights for abelian varieties over 
number  fields defined by Faltings \cite{Fal2} 
 in his proof of the Mordell conjecture in 1983. 
 Let $A$ be an abelian variety of dimension $g$ over a number field $K$.
By Grothendieck's semistable reduction theorem, there is 
 a finite extension $K'$ over $K$ so that the base change $A_{K'}=A\otimes _KK'$ has a semi-abelian model $\CA'$ over $\CO_{K'}$.
 Then we have the Hodge bundle of invariant $g$-forms by
 $$\underline\omega _{\CA'/O_{K'}}:=\wedge ^g e^*\Omega _{\CA'/O_{K'}},$$
 where $e: \Spec O_K'\to \CA'$ is the zero section. This bundle has a Faltings metric, at each archimedean place $v$ of $K'$, defined by
 $$\|\alpha \|_v^2=\frac 1{2^g}\int _{A_v(\CC)}|\alpha \wedge \bar \alpha|.$$
 Thus, we obtain a metrized line bundle $\widehat{\underline \omega}_{\CA'/O_{K'}}$. 
 Then we define the stable Faltings height by
 $$h_{\Fal}(A):=\frac 1{[K':K]}\deg (\widehat{\underline \omega}_{\CA'/O_{K'}}).$$
This height does not depend on the choice of $K'$.

\subsection{Riemann--Roch and positivity}
As in classical algebraic geometry, one expects   a Riemann--Roch theorem for the hermitian line bundle $\ol \CL$ on $\Spec O_K$ to count the number of effective sections, or more precisely, to estimate
$$\wh h^0(\ol \CL):=\log \#\{\ell \in \CL: \|\ell\|_v\le 1,  \forall\, v\mid \infty\}.$$
By Minkowski's  theorems, it is reduced to estimate the covolume of the lattice $\CL$ in the
Euclidean space $\CL_\RR=\CL\otimes_{\ZZ} \BR$ with a norm defined by $\|\cdot\|_v$.
Thus, we define the Euler characteristic of $\ol \CL$ by
$$\chi (\ol \CL):=-\log \vol (\CL_\RR/\CL),$$
where we normalize the Haar measure of $\CL_\RR=\Pi_{v|\infty}\CL_v$ so that the product of unit balls of $\CL_v$ has measure 1. 
Then, the Riemann--Roch theorem has the classical form:
$$\chi (\ol \CL)=\wh\deg (\ol \CL)+\chi (O_K), \qquad \chi(O_K)=\log (2^{r_1}\pi ^{r_2}|D_K|^{-1/2}),$$
where $r_1, r_2$ are respectively the numbers of real and complex places of $K$, and $D_K$ is the discriminant of $K$.

By Minkowski's first theorem, 
$$\wh h^0(\ol \CL)\ge \chi (\ol \CL)-[K:\BQ]\log 2.$$
 Thus $\wh h^0(\ol \CL)\ne 0$ if 
 $\wh\deg(\ol \CL)\ge \log ((4/\pi)^{r_2}|D_K|^{1/2})$. Moreover, if $\wh\deg(\ol \CL)>0$, then by Minkowski's first and second theorems, for any hermitian line bundle $\ol \CM$, we have the asymptotic formula 
$$\wh h^0(n\ol \CL+\ol \CM)=n\,\wh\deg (\ol \CL)+o(n), \qquad n\to \infty.$$

In \cite{Szp85},  Szpiro applied this terminology to give a new treatment of classical results in algebraic number theory, such as the finiteness of the class group, the Dirichlet unit theorem,
the simple connectedness of $\Spec \BZ$, and Hermite's finiteness theorem on number fields with bounded discriminants.  

\subsection{Adelic line bundles on arithmetic curves}
As an obvious generalization, we can replace $\Spec O_K$ by its generic point $\Spec K$, and consider a line bundle $L$ on $\Spec K$ with an adelic metric 
$\|\cdot\|_\BA$. Equivalently, for each place $v$ of $K$, we have a $K_v$-metric $\|\cdot\|_v$ on $L_v=L\otimes _K K_v$ such that 
for any nonzero $\ell\in L$, the metric $\|\ell\|_v=1$ for all but finitely many $v$. 
The pair $\ol L=(L, \|\cdot\|)$ form an adelic line bundle on $\Spec K$. 
In this way, we have a simple formula for the arithmetic degree by 
$$
\wh\deg (\ol L):=-\sum_{v\in M_K} \log \|\ell\|_v.
$$
There is a natural injection $\wh \Pic (O_K)\to \wh \Pic (K)$ by sending a hermitian line bundle $\ol \CL=(\CL, \|\cdot\|)$ to the adelic line bundle 
$(\CL_K, \|\cdot\|_\BA)$ whose $K_v$-metric $\|\cdot\|_v$ on $\CL_{K}\otimes _K {K_v}$ at any non-archimedean place $v$ is given by the unit ball $\CL\otimes _{O_K} O_{K_v}$.
The injection is compatible with the arithmetic degrees.

\subsection{Hermitian line bundles on arithmetic varieties}

  Let $K$ be a number field.
 By an 
 \emph{arithmetic variety} $\CX$ over $O_K$, we mean an integral scheme $\CX$ endowed with a projective and flat morphism $\pi: \CX\to \Spec O_K$. 
Note that we require arithmetic varieties to be projective in this appendix, which is different from the convention of the main body of this book.  
 
 In this way, we call $\CX$ an integral model of the projective variety $X=\CX_K$. As $\Spec O_K$ can be compactified by adding archimedean places, $\CX$ 
 can compactified by adding the complex space $\CX_\infty:=\coprod _{v\mid\infty}\CX(\ol K_v)$. For the sake of notations, it is more convenient to use the complex space 
$\CX(\CC)=\coprod _{\sigma:K\to \CC}\CX_\sigma(\CC)$, where $\CX(\CC)=\Hom_\ZZ (\Spec \CC, \CX)$ and  $\CX_\sigma(\CC)=\Hom_{(O_K,\sigma)} (\Spec \CC, \CX)$.
 
By a hermitian line bundle $\CX$, we mean a pair $\ol \CL=(\CL, \|\cdot\|)$ consisting of a line bundle $\CL$ on $\CX$ and a continuous metric $\|\cdot\|$ of the corresponding line bundle $\CL(\CC)$ on $\CX(\CC)$ which is invariant under the complex conjugate. 
Compatible with the main body of this book, the metrics of hermitian line bundles are only required to be continuous;
if we require the metric to be smooth, we will mention it specifically.
 
Now $\CLL$ defines a height function 
 $h_{\CLL}:\CX(\ol K)\to\RR$ on $X=\CX_K$. 
 In fact,
 for any algebraic point $x\in X(\ol K)$, we have a morphism $\bar x: \Spec O_{K(P)}\to \CX$ so that $\bar x_K$ is the closed point of $X$ defined by $x$.
 We then define the height of $x$ for $\ol \CL$ by 
 $$h_{\ol \CL}(x):=\frac {\wh\deg (\bar x^* \ol \CL)}{[K(P): \BQ]}.$$
 It is clear that this formalism defines a $\BR$-linear homomorphism
$$
  \wh\Pic (\CX)_\BR\lra \text{\{$\BR$-valued functions on $X(\ol K)$\}},\qquad \ol \CL\longmapsto h_{\ol \CL}.
$$
Here $\wh \Pic (\CX)$ denotes the group of isometry classes of hermitian line bundles on $\CX$.
We claim that the homomorphism is compatible with the  height machine via the natural map 
$\wh\Pic (\CX)_\RR\to \Pic (X)_\RR$, or equivalently that $h_{\CLL}$ is a Weil height on $X$ for $L=\CL_K$. 

To prove the claim, we first note that it depends only on $(X,L)$; i.e. for two different such pairs $(\CX,\CLL)$ and $(\CX',\CLL')$ over $O_K$ with generic fibers $(\CX_K,\CL_K)$ and $(\CX'_K,\CL'_K)$ isomorphic to $(X,L)$, the difference $h_{\ol \CL}-h_{\ol \CL'}$ is bounded on $X(\ol K)$. 
Second, we 
take $\CX=\BP^n_\ZZ$ and $\CL=\CO(1)$ with a hermitian metric defined by 
  $$\|\ell\|_\infty (P)=\frac {|a_0t_0+\cdots+ a_n t_n|}{\max\{ |t_0|, \cdots, |t_n|\}}, \qquad  \quad P=[t_0, \cdots, t_n],$$
where the global section $\ell$ of $\CL(\CC)$ on $\CX(\CC)$ is represented by the homogeneous polynomial $a_0x_0+\cdots+ a_n x_n$ with $a_i\in \CC$. 
  The resulting height function $h_{\CLL}$ in this case is exactly the naive height function 
  $h_{\naive}:\PP^n(\ol\QQ)\to \RR$. 
 This proves the claim for $X=\PP^n_\QQ$, and thus for general $(X,L)$ over $K$ by embedding into projective spaces over $K$.

\subsection{Faltings heights via hermitian line bundles}

Now we deal with the Faltings heights for abelian varieties. Fix a positive integer $g$, let $S_g$ denote the moduli space of principally polarized abelian varieties with some fine level structure, and let $\CS$ be a toroidal compactification of $S_g$ over $\ZZ$. Then $\CS$ is an arithmetic variety over $\BZ$ with a universal semi-abelian scheme $\CA\to \CS$. 

We have the Hodge bundle $\omega =\wedge ^g e^*\Omega _{\CA/\CS}$
 where $e: \CS\to \CA$ is the zero section. Over the open part $S_g(\BC)$ of $\CS(\CC)$, we have a Faltings metric $\|\cdot\|_\Fal$ on $\omega$ defined by glueing the Faltings metrics on fibers of $\omega$ on $S_g(\BC)$. 
By the above procedure, the resulting metrized line bundle 
 $\bar \omega$ defines a height function $h_{\bar\omega}:S_g(\ol \BQ)\to \RR$, which glues  Faltings heights in that $h_{\bar\omega}(s)=h_\Fal(\CA_s)$ for any $s\in S_g(\ol \BQ)$. 

The Faltings metric on $\omega$ is smooth on $S_g(\BC)$, but not continuous along the boundary $\CS(\CC)\setminus S_g(\BC)$.  
In fact, Faltings proved that 
 the metric on $\omega $ has a log-log singularity in the sense that for every point $z_0\in \CS(\CC)\setminus S_g(\BC)$, there is an open neighborhood $U\subset \CS(\BC)$ of $z_0$ such that the boundary $Z=U\setminus S_g(\BC)$ is defined by equations $f_1, \cdots, f_k$ on $U$, that the restriction $\omega|_U$ is generated by  a section $\alpha$, and that there is a constant $C>0$ satisfying
 $$\Big|\log \|\alpha(z)\|_\Fal\Big|\le C\log \Big(-\log \big(\sum_{i=1}^k  |f_i(z)|\big)\Big), \qquad \forall\, z\in U.$$
As a consequence, Faltings proves that the height function $h_{\bar\omega}:S_g(\ol \BQ)\to \RR$ satisfies the Northcott property. 
We will see that the metrized line bundle $\bar\omega$ on $S_g$ is a canonical example of an adelic line bundle on a quasi-projective varieties in our theory. 
 
In the end, let us try to treat the N\'eron--Tate height function, or more generally the Call--Silverman height function
 in terms of hermitian line bundles. Then, we immediately realize that we could not do it directly as an algebraic dynamical system over a number field does not necessarily extend to a (projective) arithmetic variety.
 To cover Arakelov theory for these heights, we need to introduce adelic metrics.
 The key point is to define continuous metrics on algebraic points over non-archimedean fields.
 This will be the main topic in \S\ref{app sec adelic}.

\section{Arithmetic surfaces}\label{app sec surface}

In this section, we describe  intersection theories on arithmetic surfaces based on the works of Arakelov \cite{Ara74}, Faltings \cite{Fal2}, Deligne \cite{Del}, Szpiro \cite{Szp90},
and Zhang \cite{Zha92, Zha93}. See also Szpiro \cite[\S2]{Szp85}  and Yuan--Guo \cite[\S5]{YG25}.
We will only deal with arithmetic varieties over $\BZ$. 
Of course, most results hold for a function field of one variable over a field.

\subsection{Deligne pairings}
Let $S$ be a noetherian scheme, and let $\pi: X\to S$ be a flat and proper family of curves which is a local complete intersection. 
Let $\Picc (X)$ be the groupoid of line bundles on $X$; i.e. $\Picc (X)$ is the category whose objects are line bundles on $X$ and whose morphisms are isomorphisms of line bundles on $X$. 
The Deligne pairing, first introduced by Deligne \cite{Del}, is a canonical bilinear functor 
$$\Picc(X)\times \Picc(X)\lra \Picc(S), \qquad (L, M)\longmapsto \pair{L, M}.$$
At a point $s\in S$, the fiber $\pair{L, M}(s)$ is generated by symbols $\pair{\ell, m}$ indexed by nonzero rational sections $\ell$ and $m$ of $L|_{X_s}$ and $M|_{X_s}$ with 
$|\div(\ell)|\cap |\div(m)|=\emptyset$,  such that for any rational function $f$ on $X$ with $|\div(f)|\cap |\div(m)|=\emptyset$,
$$\pair{f\ell, m}=f(\div (m))\,\pair{\ell, m}.$$
Here we take the multiplicative convention $f(\sum_i a_i x_i)=\prod_i f(x_i)^{a_i}$ for values of $f$ on divisors. 

Moreover, let $\omega=\omega_{X/S}$ denote the relative dualizing sheaf of $\pi$. Then we have   a functorial Riemann--Roch isomorphism
\begin{equation}\label{eq-RR}
2\det R\pi_* L\iso \pair{L, L-\omega}+2\det R\pi_* \CO_X.
\end{equation}
This isomorphism is unique up to multiplication by $\pm 1$.

If $\pi$ is smooth, then there is a functional Noether isomorphism 
$$\alpha: \pair{\omega, \omega}\iso 12\det R\pi_*\omega.
$$
If $\pi$ is only smooth at every generic point of $S$, then the above isomorphism defines a Cartier divisor $\Delta$ on $S$, called the discriminant of $\pi$, so that $\alpha$ induces an isomorphism 
\begin{equation}\label{eq-noet}
\alpha: \pair{\omega, \omega}\otimes \CO_S(\Delta) \iso 12\det R\pi_*\omega.
\end{equation}

If $\pi$ is a fibration from a smooth projective surface $X$ to a smooth projective curve $S$ over a field, then taking degrees gives us a Riemann--Roch formula
$$\deg(\det R\pi_*L)=\frac 12 L^2-\frac 12 L\cdot \omega+\frac 1{12}(\omega ^2+\deg \Delta).$$
This induces the classical Riemann--Roch theorem for $X$ using the Riemann--Roch theorems on $S$ and on the generic fiber of $X\to S$.

If $X$ and $S$ are smooth complex varieties, then the Deligne pairing can be extended to hermitian line bundles with smooth metrics to get a canonical bilinear pairing
$$\wh\Picc(X)_\sm\times \wh \Picc(X)_\sm\lra \wh \Picc(S), \qquad (\ol L, \ol M)\longmapsto \pair{\ol L, \ol M}.$$
Here $\wh\Picc(X)_\sm$ denotes the category of hermitian line bundles on 
$X$ with smooth metrics. 
The metric of $\pair{\ol L, \ol M}$ at a point $s\in S$ is given by 
$$
\log \|\pair{\ell, m}\|= \int _{X_s} \log \|\ell\| c_1(\ol M)+\log \|m\| (\div (\ell)).
$$
Here we take the additive convention $\log \|m\|(\sum_i a_i x_i)=\sum_i {a_i} \log \|m(x_i)\|$ for values of $\log \|m\|$ on divisors. 
The metric on $\pair{\ol L, \ol M}$ is continuous with the curvature form given by
$$c_1(\pair{\ol L, \ol M})=\int _{X/S} c_1(\ol L) c_1(\ol M).$$

As a dilation, let us note that Deligne \cite{Del} speculated a similar pairing for suitable projective morphisms $\pi:X\to S$ of general relative dimension $n\geq 0$. Then we need construct a line bundle on $S$ from $n+1$ line bundles on $X$.  
This problem was settled by by Elkik \cite{Elk1} assuming that $\pi$ is projective, flat, and Cohen--Macaulay, and by Munoz Garcia \cite{MG} assuming that $\pi$ is projective, equi-dimensional and of finite Tor-dimension (which implies the projective and flat case). 
Moreover, if $\pi$ is a projective and smooth morphism of complex varieties, Elkik \cite{Elk2} also treated Deligne pairings of hermitian line bundles with smooth metrics.

\subsection{Quillen metric}
Let $X$ be a compact Riemann surface whose canonical sheaf  $\omega$ is endowed with a smooth hermitian metric. We write the data as $\ol X=(X, \bar\omega)$.
Let $\ol L=(L, \|\cdot\|)$ be a hermitian line bundle with a smooth metric on $X$.  We want to put a norm on $\det H^*(X, L)$. 
Notice that  the cohomology $H^*(X, L)$ can be computed by the Dolbeault complex
$$\Omega ^{0, 0} (L)\overset {\bar \partial} \lra \Omega^{0, 1} (L).$$
The metric on $\omega$  defines a volume form $\ds d\mu=\frac i2 \alpha \wedge \bar \alpha$ on $X$, where $\alpha$ is a $(1,0)$-form with norm 1 under the metric of $\omega$.
Thus, we have $L^2$-norms on the space of smooth forms. 
The Quillen metric on  
$$\det H^*(X, L):=(\det H^0(X, L))\otimes (\det H^1(X, L))^\vee$$
is formally defined 
by the formula
$$\det H^* (\ol X, \ol L)_Q=(\det \Omega^{0,0} (\ol L)_{L^2})\otimes (\det \Omega^{0,1} (\ol L)_{L^2})^\vee,$$
provided a reasonable regularization of the right-hand side. 

Notice that the $L^2$-norms on smooth forms induce $L^2$-norms on the cohomology groups by the orthogonal decompositions
$$\Omega ^{0,0}(\ol L)_{L^2}=H^0(X, \ol L)_{L^2} \oplus (\ker \bar\partial )^\perp,$$
$$\Omega^{0,1}(\ol L)_{L^2}=\Im\, \bar\partial \oplus H^1(X,\ol L)_{L^2}.$$
Thus, for an element $\ell$ of $\det H^* (X, \ol L)$, the Quillen metric is defined by 
$$\|\ell\|_Q=\|\ell \|_{L^2} \cdot \|\det \bar\partial'\|^{-1},$$
where $\bar\partial': (\ker\bar\partial)^\perp\to \Im\, \bar\partial$ is the restriction of $\bar\partial$.
As a  convention of metrized line bundles twisted by a number, we may write 
$$
\det H^* (\ol X, \ol L)_Q=\det H^* (X, \ol L)_{L^2} (\tau (\ol L))
$$
with $\tau (\ol L)$ Quillen's {\em analytic torsion} formally given by
$$\tau(\ol L)=\log  \|\det \bar\partial'\|.$$

To define Quillen's analytic torsion,  we consider the adjoint operator $\bar\partial ^*$ of $\bar\partial $ and form the Laplacian operator 
$\Delta:=\bar\partial ^*\bar\partial$ on $\Omega^{0,0} (\ol L)$. Then we have an orthonormal basis  $\{\varphi_i\}_i$
formed by eigenvectors $\varphi_i$ of $\Delta$ with
$$\Delta \varphi_i=\lambda _i  \varphi_i.$$
The space $(\ker\bar\partial)^\perp$ is spanned by those $\varphi_i$ with $\lambda _i>0$.
By definition, it is easy to see that $\bar\partial \varphi_i$ are orthogonal to each other. 
Thus we have the formal expression
$$\|\det \bar\partial'\|^2=\prod_{\lambda _i>0} \|\bar\partial \varphi_i\|^2
=\prod_{\lambda _i>0} \lambda _i.$$
To regularize the last infinite product, Quillen uses the zeta function
$$\zeta _{\ol L} (t)=\sum _{\lambda _i>0} \lambda _i^{-t}.$$
This function is convergent if $\Re(t)$ is sufficiently large and has a meromorphic continuation to the complex plane.
Taking derivative at $t=0$, we get formal equalities
$$\zeta'_{\ol L}(0)=-\sum _{\lambda_i>0}\log \lambda _i=-2\log \|\det \bar\partial'\|.$$
Therefore, Quillen's analytic torsion is rigorously defined by 
$$\tau (\ol L):=-\frac 12 \zeta _{\ol L}'(0).$$

\subsection{Deligne's Riemann--Roch Theorem}
Let $X$ be a compact Riemann surface, and let  
$\ol X=(X, \bar\omega)$ be a pair as above.
Deligne proves that the Riemann--Roch isomorphism  and the Noether isomorphism are  isometries for the Quillen metrics:
$$2\det H^*(\ol X, \ol L)_Q \iso \pair{\ol L, \ol L-\bar \omega}+2\det H^* (\ol X, \CO_X)_Q,$$
$$
\pair{\bar \omega, \bar\omega}\iso 12\det H^*(X, \CO_X)_Q.
$$

Write $\tau (X)=\tau (\CO_{X})$. Then $ \tau (X)$ measures how far the K\"ahler manifold $\ol X$ becomes singular. In fact, $-12\tau (X)$ is the right archimedean analogue of the discriminant for $X$ over a non-archimedean field.
 Thus, we modify Deligne's determinant by 
\begin{equation}\label{eq-normQ}
\det H^*(\ol X, \ol L):=\det H^*(\ol X, \ol L)_Q(-\tau (X)).
\end{equation}
Then we obtain 
\begin{equation}\label{eq-RRQ}
2\det H^*(\ol X, \ol L) \iso \pair{\ol L, \ol L-\bar \omega}+2\det H^*(\ol X, \CO_X)_{L^2},\end{equation}
\begin{equation}\label{eq-noetQ}
\pair{\bar \omega, \bar\omega}(-12 \tau(X))\iso 12\det H^*(\ol X, \CO_X)_{L^2}
\end{equation}

Finally, we notice that the metric on the bundle $\lambda (X):=\det H^*(X, \CO_X)_{L^2}$ does not depend on the choice of the metric on $\omega$.
In fact, if the genus $g(X)=0$, then $\lambda (X)$ is the trivial bundle. If $g(X)>0$, then $\lambda(X)$ is canonically isomorphic to $\det \Gamma (X, \omega )$
with the norm defined by the $L^2$-product
$$\pair{\alpha, \beta}=\frac i 2\int _{X}\alpha \wedge \bar \beta.$$

\subsection{Arakelov pairings}
 Let $X$ be a compact Riemann surface with a K\"ahler from $d\mu$ of total volume $1$. Then we  have a Green function $g:X\times X\setminus \Delta \to \RR$ of smooth type
which is symmetric and satisfies the equation of distribution:
$$\frac {\partial _x\bar\partial _x}{\pi i }g(x, y)=\delta _y(x)-d\mu (x).$$
This defines a hermitian line bundle  $\ol \CO_{X\times X}(\Delta)$ with a smooth metric given by the formula $-\log \|1\|(x, y)=g(x, y)$.
By linearity, this gives an Arakelov pairing 
$\pair {D, E}_A=g(D, E)$ for two divisors $D,E$ on $X$ with disjoint supports.

A hermitian line bundle $\ol L$ with a smooth metric on $X$ is called {\em admissible} if 
$c_1(\ol L)=(\deg L) \cdot d\mu$. Here is an example of admissible line bundles. For any $P\in X$, we have a metrized line bundle  $\ol \CO(P)$ by taking restriction of $\ol \CO(\Delta)$ to $\{P\}\times X$. 
We extend this to define $\ol \CO(D)$ for any divisor $D$ on $X$ by linearity.

Restricting to the diagonal, we also get a {\em dualizing sheaf}    $\bar \omega =(\Delta^*\ol \CO_{X\times X}(\Delta))^\vee$
for $d\mu$.
A main result of Arakelov \cite{Ara74} is that there is a unique K\"ahler form $d\mu$ on $X$ of total volume 1 so that $\bar\omega$ is admissible for $d\mu$. 
In fact, it can be expressed by
$$d\mu =\frac i{2g}\sum _i \alpha _i\wedge \bar \alpha_i,$$
where $\alpha _1, \cdots, \alpha _g$ is an orthonormal base for $\Gamma (X, \omega)_{L^2}$.

Now for any two divisors $D$ and $E$ with disjoint supports, 
 the Deligne pairing agrees with the Arakelov pairing in the sense that  
 $$\pair{\ol \CO(D), \ol \CO(E)}\simeq \BC (g(D, E))$$
 holds as an isometry of metrized complex lines.
Moreover, we have the  adjunction formula
$$\pair{\bar \omega(P), \ol \CO(P)} \iso \bar \omega(P)|_P\iso \BC,$$
where the second arrow is defined by taking residue at $P$.

For Arakelov metrics,  Faltings \cite{Fal1} constructed a metrized determinant for admissible line bundles without using Quillen metrics.
In fact, for any line bundle $L$ on $X$ and any point $P\in X$, 
the exact sequence
$$0\lra L\lra L(P)\lra L(P)|_P\lra 0$$
induces canonical  isomorphisms 
$$\det H^*(X, L(P))\iso \det H^*(X, L)\otimes L(P)|P\iso \det H^*(X, L)\otimes (L\otimes \omega^{-1})|_P.$$
Hence, for any admissible hermitian line bundle $\ol L$ on $X$, there is a unique way to attach a norm on the complex line $\det H^*(X, L)$ such that the resulting metrized line $\det H^*(\ol X, \ol L)_F$ satisfies the following two rules:
$$\det H^*(\ol X, \ol L(P))_F\iso  \det H^*(\ol X, \ol L)_F\otimes (\ol L\otimes \bar \omega^{-1})|_P.$$
$$\det H^*(\ol X, \ol \CO(c))_F\iso \det H^*(\ol X, \CO_X)_{L^2}((1-g)c), \qquad \forall\, c\in \BR.$$
By this construction, we have the same Riemann--Roch theorem given by
$$2\det H^*(\ol X, \ol L)_F\iso \pair{\ol L, \ol L-\bar \omega}+2\det H^*(\ol X, \CO_X)_{L^2}.$$
Thus, Faltings's determinants agree with  Deligne's modified determinants.

\subsection{Arithmetic surfaces}
Now we consider a projective and flat morphism  $\pi: \CX\to S=\Spec O_K$ with $\CX$ regular of dimension $2$, and $O_K$ the ring of integers of a number field $K$. We call $\CX$ an {\em arithmetic surface} over $O_K$. 
 Let $\wh\Picc (\CX)_\sm$ denote the category of hermitian line bundles with smooth metrics on $\CX$ in which morphisms are given by isometries, 
and let $\wh\Pic (\CX)_\sm$ be the group of isomorphism classes of $\wh\Picc (\CX)_\sm$.
Then we have the Deligne pairing
$$ \quad \wh\Picc (\CX)_\sm\times \wh \Picc (\CX)_\sm\lra \wh\Picc (S), \quad  (\ol \CL, \ol \CM)\longmapsto 
\pair{\ol \CL, \ol \CM}.$$
Taking arithmetic degrees, we obtain Deligne's intersection pairing
$$ \quad \wh\Pic (\CX)_\sm\times \wh \Pic (\CX)_\sm\lra \RR, \quad (\ol \CL, \ol \CM)\longmapsto \ol \CL\cdot \ol \CM:=\wh\deg
(\pair{\ol \CL, \ol \CM}).$$

Let $\omega=\omega_{\CX/O_K}$ be the relative dualizing sheaf. 
Endow $\omega$ with a smooth hermitian metric, so that we have a hermitian line bundle $\bar\omega$.  We write $\ol \CX$ for the pair $(\CX, \bar\omega)$. 
For each hermitian line bundle $\ol \CL$ with a smooth metric, we have a hermitian line bundle $\det H^*(\ol \CX, \ol \CL)$ on $S$ with the normalized Quillen metric in (\ref{eq-normQ}).
 We apply the Riemann--Roch formulae (\ref{eq-RR}), (\ref{eq-RRQ}) and  the Noether formulae (\ref{eq-noet}), (\ref{eq-noetQ}) to obtain the following arithmetic Riemann--Roch theorem and arithmetic Noether formulae:
$$
 \wh\deg (\det H^*(\ol \CX, \ol \CL))=\frac 12 \ol \CL\cdot (\ol \CL-\bar \omega)+\wh\deg
 (\det H^*(\CX, \CO_\CX)_{\CL^2}),
  $$
$$
 \wh\deg( \det H^*(\CX, \CO_\CX)_{\CL^2})=\frac 1{12}
 (\bar \omega^2+\deg \bar \Delta(\CX)),
$$
 where $\ol\Delta(\CX)$ is the  arithmetic discriminant divisor of $\ol \CX/ S$ on $S$ defined by 
$$
 \ol \Delta(\CX) =\sum _{v\nmid \infty}\ord _v\Delta(\CX) [v]-12 \sum_{v\mid\infty}\tau (\CX_v)[v] .
$$
 
 Now assume that $\bar \omega $ carries the Arakelov metrics, and that $\ol \CL$ is admissible. Then the above formulae are exactly the arithmetic Riemann--Roch theorem and the arithmetic Noether formula proved by Faltings \cite{Fal1}.

As in the case of arithmetic curves, the first application of the arithmetic Riemann--Roch theorem  is a criterion for the effectivity of hermitian line bundles.
An \emph{effective section} of a hermitian line bundle $\ol\CL$ on $\CX$ is a section $s\in H^0(\CX, \CL)$ such that 
$$\|s\|_{\sup}=\sup_{x\in \CX(\CC)} \|s\| (x)\leq1.$$

For a hermitian line bundle $\ol \CL$ with a smooth metric, we say that $\ol \CL$ is {\em relatively semipositive} if $\ol \CL$ has a semipositive curvature form at infinite places and non-negative degree with all curves in the special fibers.
Then Faltings \cite{Fal1} proved the following theorem.

\begin{thm}
Assume that $\ol \CL$ is relatively semipositive on $\CX$,  that $\CL_K$ is ample on $\CX_K$, and that the self-intersection number $\ol\CL^2>0$.
Then for sufficiently large integer $n$, the tensor power $n  \CL$ has a non-zero effective section.
\end{thm}

The second application of the arithmetic Riemann--Roch theorem is the arithmetic Hodge index theorem for the intersection pairing on $\wh \Pic (\CX)$.
We call a hermitian line bundle $\ol\CL$ on $\CX$ vertical if the generic fiber $\CL_K$ is isomorphic to the trivial line bundle $\CO_{\CX_K}$ on $\CX_K$.

\begin{thm}[arithmetic Hodge index, Faltings \cite{Fal1}, Hriljac \cite{Hri}]\label{app hodge index1}  
Let $\ol \CL$ be a hermitian line bundle with a smooth metric on $\CX$ such that $\deg \CL_K=0$. Then $\ol \CL^2\le 0$, and  $\ol \CL^2=0$ if and only if $\ol \CL=\pi^*\ol \CM$ for a hermitian line bundle $\ol \CM$ on $S$. Moreover, if $\ol \CL$ is perpendicular to all vertical hermitian line bundles, then
$$\ol \CL^2=-2[K:\QQ]\,\wh h(\CL_K),$$
where $\wh h(\CL_K)$ is the Neron--Tate height of $\CL_K\in \Jac (\CX_K)(K)$ with respect to the theta divisor. 
\end{thm}

This theorem was proved independently by Faltings \cite{Fal1} and Hriljac \cite{Hri}. 
As a beautiful application of this Hodge index theorem, Faltings \cite{Fal1} proved the semi-positivity (or equivalently nefness) of the relative dualizing sheaf with the Arakelov metric; i.e.
$h_{\bar \omega} (x)\ge 0$ for all $x\in \CX(\ol K)$, and $\bar \omega^2\ge 0$.

\subsection{Arithmetic ampleness}
As in classical algebraic geometry, we have a numerical criterion for arithmetic ampleness as follows.
Let $(\CX, \ol \CL)$ as above. Namely, 
 $\ol \CL$ is a hermitian line bundle with a smooth metric on an arithmetic surfaces $\CX$. 
We say that $\ol \CL$ is {\em arithmetically ample} if  for any hermitian line bundle $\ol \CM$, the $O_K$-module $H^0(\CX, n\CL+\CM)$
is generated by effective sections.

\begin{thm}[Nakai--Moishezon theorem, \cite{Zha92}]
Assume that $\ol \CL$ is relatively semipositive. Then the following conditions are equivalent:
\begin{enumerate}
\item $\ol \CL$ is arithmetically ample;
\item $\wh \deg (\ol \CL|_\CY)>0$ for any arithmetic curve $\CY$ in $\CX$.
\end{enumerate}
\end{thm}

Motivated by Szpiro's idea in \cite{Szp90}, $\ol \CL^2$ should be related to the essential minimum $e'(\ol \CL)$ of the height function 
$$h_{\ol \CL}(x)=\frac {\wh \deg (\ol\CL|_{\bar x})}{\deg x}, \qquad x\in \CX_K(\ol K),$$
where $\bar x$ is the Zariski closure of $x$ in $\CX$.
The following theorem of successive minimum gives a quantitative relation between these heights. 
More precisely, we define the successive minimum of $\ol \CL$ by the following formulae:
$$e(\ol \CL)=\inf _{x\in \CX(\ol K)}h_{\ol \CL}(x), \qquad e'(\ol \CL)=\liminf _{x\in \CX(\ol K)}h_{\ol \CL}(x),$$
where $\liminf$ denotes the smallest limit point. 
Then we have the following theorem.
\begin{thm}[successive minimum, Zhang \cite{Zha92}] Assume that $\deg (\CL_K)>0$ and that $\ol \CL$ is relatively semipositive.
 Then
$$e'(\ol \CL)\geq  \frac{\ol \CL^2}{2\deg \CL_K}\geq
\frac12 (e'(\ol \CL)+e(\ol \CL)).$$
\end{thm}
One consequence of the theorem is that small points are dense in $\CX_K$ if and only if $\ol \CL^2=0$.

\subsection{Non-archimedean admissible pairings} 
\label{app subsec admissible}
Let $K$ be a discrete valuation field and $X/K$ a smooth projective curve of positive genus. 
Using an idea of Chinburg--Rumely \cite{CR93}, the work of Zhang
\cite{Zha93} gave an extension of the Deligne pairing to $\ol K$-metrized line bundles over $X$ using metrized graphs. 
Replacing $K$ by a finite extension, we may assume that $X$ has a semistable minimal regular model $\CX$ over $O_K$.
Then for any finite extension $K'$ of $K$, $X_{K'}$ has a semistable minimal regular model $\CX'$ over $O_{K'}$.

Let $G=G_K$ be the dual graph of the special fiber $\CX_k$ of 
$\CX$ over the residue field $k$ of $O_K$.
Then $G_K$ has the vertex set $V_K$ parameterizing the irreducible components of $\CX_k$ and the edge set $E_K$  parameterizing the double points of $\CX_k$.
We make $G$ a metrized graph so that each edge has length $1$. The reduction process defines a map
$$r_K: X(K)\lra G_K.$$
Let $F(G)$ be the space of functions on $G$, which is piecewise smooth with finite directional derivatives at non-smooth points.
Then we have a Laplacian operator $\Delta (f)$ as a measure $F(G)$ defined by $-f''(x)dx$ on smooth sides and by $f'(x)\delta _x$ at singular points $x$.
Denote by $\Div _V (\CX)$ the group of vertical $\RR$-divisors of $\CX$, i.e. $\RR$-linear combinations of the irreducible components of the special fiber $\CX_k$. 
Then we have a map 
$$F_K: \Div _V (\CX)\lra F (G)$$
which sends a vertical divisor $D$  of $\CX$ to the function $F_D\in F(G)$ which is linear on $E_K$ and whose value on a vertex $v$ is just $\ord _v(D)$.

Let $K'/K$ be a finite extension of ramification index $e$. Then we have a canonical identification between the dual metrized graphs $G_K$ and $G_{K'}$, where $V_{K'}$ is equal to the union of the sets of  $n$-section points of edges in $E_K$. The maps $r_K$ and $F_K$ are compatible with 
$r_{K'}$ and $F_{K'}$. Let $\wh\Div _V (X_{\ol K})$ be the direct limit of $\Div _V (\CX_{K'})$ over all such $K'$. Then we have  reduction maps
$$r: X(\ol K)\lra G, \qquad F: \wh \Div _V (X_{\ol K})\lra F(G).$$
Through $r$, we may consider elements of $F(G)$ as functions on $X(\ol K)$.

Now we denote by $\wh \Picc _G (X)$ the category of $\ol K$-metrized line bundles $\ol L=(L, \|\cdot\|)$ on $X(\ol K)$ so that there is a finite extension 
$K'$ over $K$ and a line bundle $\CL $ over $\CX'$ satisfying $\ol L=\CL (f)$ in the sense that 
$\|\cdot\|=\|\cdot\|_{\CL}\exp (-f)$ for some function $f\in F(G)$.
For each $\ol L\in \wh \Picc _G(X)$, we have a curvature map $c_1(\ol L)$ on $F(G)$ extending the intersection pairing between vertical cycles and line bundles on integral models $\CX'$. 
Now  for two metrized line bundles $\ol L$ and $\ol M$ represented by $\CL (f)$ and $\CM (g)$, we define their Deligne pairing as
$$\pair{\ol L, \ol M}=\pair{\CL, \CM}\left(\int _G g c_1(\CL)+\int_G f c_1(\CM)-\int _G f \Delta g\right).$$

Let $d\mu$ be a measure on $G$ of volume $1$. Then a metrized bundle in $\wh \Picc _G(X)$ is called {\em admissible for $d\mu$} if $c_1(\ol L)=(\deg L)\cdot d\mu$. For such a measure $d\mu$,  we can construct a Green function $g(x, y)$ on $G\times G$. This will give a metrized line bundle
$\ol \CO(\Delta )$ as the usual line bundle on $\CO_{\CX\times \CX}$ twisted by $g (x, y)$. In this way, for any $P\in X(\ol K)$, we have an admissible line bundle $\ol \CO(P)$. By linearity, we get an admissible metrized line bundle $\ol \CO(D)$ for any divisor $D$ of $X$. 

Restricting to the diagonal, we also get dualizing sheaf $\bar \omega =(\Delta ^* \ol \CO_{X^2}(\Delta))^\vee$. The first main result in \cite{Zha93} is the existence of a unique semipositive measure $d\mu_a $ on $G$ so that the dualizing sheaf $\bar \omega_a $ is also admissible.  We call such a metric an {\em admissible metric}. The second main result is a comparison between  $\bar\omega_a$ 
and the relative dualizing sheaf $\omega _{\CX/O_K}$ of the semistable model $\CX$ with the Arakelov metric.
The final result gives 
$$\pair{ \omega _{\CX/O_K},  \omega_{\CX/O_K}}=\pair{\bar\omega_a, \bar \omega_a}(c).$$
We have $c\ge 0$, and $c=0$ if and only if $g(X)=1$ or $X$ has a potentially good reduction. 

Now  we return to the global situation in which $K$ is a number field, and $X$ is a smooth and projective curve over $K$ of genus $ \ge 2$.
Then we can certainly define the Deligne pairings for {\em adelic metrized line bundles} in a natural way. We also have the Riemann--Roch theorem,  Hodge index theorem, and the theorem of successive minima.
Most importantly, we have an admissible adelic line bundle $\bar \omega_a$ so that for any $P\in X(K)$,
$$(\bar \omega_a -(2g-2)\ol \CO(P))^2=-2[K:\QQ]\wh h (\omega -(2g-2)P).$$
These facts imply that $\bar\omega _{\mathrm{Ar}}^2\ge \bar \omega _a^2$, and the equality holds if and only if $X$ has good reduction everywhere. 
Here $\bar \omega _{\mathrm{Ar}}=(\omega_{\CX/O_K},\|\cdot\|_{\mathrm{Ar}})$ is the relative dualizing sheaf of the global minimal regular model $\CX$ over $O_K$ with the Arakelov metric at archimedean places, and $\CX$ is assumed to be semistable, which can be achieved by passing to a finite extension of $K$.
They also imply that  
$\omega_a^2\ge 0$, and that $\omega _a^2>0$ if and only if the Bogomolov conjecture holds for the embedding $X\to \Jac (X)$ using the base class $\omega/(2g-2)$. Recall that the Bogomolov conjecture in this case  asserts that
$$\liminf_{P\in X(\ol K)} \wh h (\omega -(2g-2)P)>0.$$
The Bogomolov conjecture, in this case, is eventually proved by Ullmo \cite{Ull} using an equidistribution theorem of Szpiro--Ullmo--Zhang \cite{SUZ}.
  
It is worth mentioning that we can re-construct the term $\bar \omega_a$ as an adelic line bundle on $X$ in terms of certain pull-back of invariant adelic line bundles on the Jacobian variety.  We refer to the appendix of Yuan \cite{Yua4} for more details on this interpretation.

\section{Arithmetic intersection theory}\label{app sec intersection}

In this section, we sketch intersection numbers of hermitian line bundles essentially due to Deligne \cite{Del}, intersection theory of arithmetic Chow groups of Gillet--Soul\'e \cite{GS1},
the arithmetic Riemann--Roch theorem due to 
Gillet--Soul\'e \cite{GS2}, and an arithmetic Hilbert--Samuel formula due to 
Gillet--Soul\'e \cite{GS2}. 
For updated details, see Yuan--Guo \cite[\S4, \S6, \S7]{YG25}.

\subsection{Definitions of arithmetic intersections in history}

Let us briefly recall the history of development of the definitions of arithmetic intersections. 
In 1974,  the pioneering work of Arakelov \cite{Ara74} introduced an arithmetic intersection theory of hermitian line bundles with admissible metrics (or equivalently arithmetic divisors with admissible Green functions) on arithmetic surfaces.
The work of Faltings \cite{Fal1} was also based on Arakelov's original setting.
However, the assumption of admissible metrics seems very restrictive in application and very hard to generalize to high codimensions, and it is imperative to have a canonical intersection number of any two hermitian line bundles. 
In the work of Deligne \cite[\S6]{Del}, from the Deligne pairing of hermitian line bundles, we have immediately obtained a definition of such an intersection number. 
Moreover, Deligne proved an arithmetic Riemann--Roch theorem using the Quillen metrics, which generalized Faltings's arithmetic Riemann--Roch theorem. 
This gives a satisfactory intersection theory on arithmetic surfaces. 

For general dimensions, Deligne's intersection formula on arithmetic surfaces can be easily generalized to get a canonical definition of arithmetic intersection numbers of $d$ hermitian line bundles on an arithmetic variety of dimension $d$.
This can also be seen from Elkik \cite{Elk1, Elk2}, which generalizes the Deligne pairing to general relative dimensions.
This intersection theory is sufficient for the purpose of this book. 

To obtain a full arithmetic intersection theory, we need to extend the above intersection theory to arithmetic cycles of arbitrary codimensions. 
This was done by the foundational work of Gillet--Soul\'e \cite{GS1}, where they introduced the concept of arithmetic Chow cycles, where Green functions for divisors are generalized to Green currents for Chow cycles, and defined intersections of two arbitrary arithmetic Chow cycles as a new arithmetic Chow cycle. 
Moreover, Gillet--Soul\'e \cite{GS2} proved an arithmetic Riemann--Roch theorem in the style of Grothendieck's Riemann--Roch theorem, which generalizes Deligne's Riemann--Roch theorem to morphisms of arithmetic varieties.

\subsection{Intersection theory of hermitian line bundles} \label{app subsec intersection}
 
Let $K$ be a number field, and let $\CX$ be a projective arithmetic variety over $O_K$. 
Recall that a hermitian line bundle $\ol \CL=(\CL, \|\cdot\|)$ on $\CX$ consists of a line bundle $\CL$ on $\CX$ and a continuous metric of $\CL(\BC)$ on $\CX(\BC)$ invariant under the action of the complex conjugate. 
Note that we do not assume that $\CX(\CC)$ is smooth. 
We say that the metric $\|\cdot\|$ is smooth if there is an open covering $U_i$ of  $\CX(\BC)$ by complex open subsets and analytic embedding $U_i\to V_i$ into complex manifolds $V_i$ so that the hermitian line bundle $\ol \CL|_{U_i}$ is the restriction of some hermitian line bundle $\ol \CM_i$ with a smooth metric on $V_i$. 
Then we can define the curvature $c_1(\ol \CL)$ as a smooth form on $\CX(\BC)$ locally as the restriction of $\displaystyle c_1(\ol \CM_i)=\pp \log \|m_i\|$ where $m_i$ is an invertible section of $\CM_i$ on $V_i$.
 
Let $\CZ$ be a closed integral subscheme of $\CX$ of dimension $d\ge 1$. 
We say that $\CZ$ is horizontal if $\CZ\to \Spec O_K$ is surjective; we
say that $\CZ$ is vertical if the image of $\CZ\to \Spec O_K$ is a closed point.  
Let $\ol \CL_1, \cdots, \ol \CL_d$ be $d$ hermitian line bundles with smooth metrics on $\CX$. 
Then we define the intersection number $\ol \CL_{1}\cdot \ol \CL_{2}\cdots \ol \CL_{d}\cdot[\CZ]$ by induction on $d$ as follows. 

If $d=1$, denote by $\wt\CZ\to \CZ$ the normalization, and by $\varphi: \wt \CZ\to \CX$ the induced finite morphism.
If $\CZ$ is vertical, then $\wt\CZ$ is a smooth projective curve over a finite field 
$\FF_p$, and we define 
$\ol \CL_1\cdot[\CZ]=\deg_{\FF_p}(\varphi^* \CL_1)\log p$. 
If $\CZ$ is horizontal,  then $\wt\CZ$ is isomorphic to an arithmetic curve $\Spec O_{K'}$ for some finite extension $K'$ of $K$, and we define $\ol \CL_1\cdot[\CZ]=\wh\deg(\varphi^* \ol\CL_1)$.

If $d>1$, we take a nonzero rational section $\ell_d$ of $\CL_d$ over $\CZ$. 
Assume that $\div (\ell_d)=\sum _i a_i \CZ_i$ with $a_i\in \BZ$ and $\CZ_i$ integral subschemes of $\CZ$ of dimension $d-1$.
We define 
$$
\ol \CL_{1}\cdots \ol \CL_{d}\cdot[\CZ]=\sum _i a_i\,  \ol \CL_{1}\cdots  \ol \CL_{d-1}\cdot[\CZ_i]
- \int_{\CZ(\BC)}\log \|\ell_d\| c_1(\ol \CL_1)\cdots c_1(\ol \CL_{d-1}).
$$
Here we take the convention that the integral on the right-hand side is zero if $\CZ(\CC)=\emptyset$, which happens if $\CZ$ is vertical.

With commutative algebra and the Stokes formula, we can prove that the intersection number is independent of the choice of $\ell_d$ and the ordering of $\ol \CL_1,\cdots, \ol \CL_d$.
When $\CZ=\CX$, we simply write $\ol \CL_{1}\cdots \ol \CL_{d}$ for $\ol \CL_{1}\cdots \ol \CL_{d}\cdot[\CX]$. 

One main application of the intersection theory is to introduce the arithmetic degree of a subvariety $Z$ of $X=\CX_K$ by
$$\wh\deg_{\ol \CL}(Z)=\ol \CL ^{\dim Z+1}\cdot [\ol Z],$$
 where $\ol Z$ is the Zariski closure of $Z$ in $\CX$. This was the height used by Faltings \cite{Fal90} in his proof of the Mordell--Lang conjecture for subvarieties of abelian varieties. However, when studying the positivity of the heights as in \cite{Zha1, Zha2} and particularly its relation to small points, 
 it is more convenient to normalize this height function by
 \begin{equation*}\label{eq-h(Z)-herm}
 h_{\ol \CL}(Z)=\frac{\wh\deg_{\ol \CL}(\ol Z)}{(\dim Z+1) \deg _{\CL_K}(Z)},
 \end{equation*}
where $\deg _{\CL_K}(Z)=\CL_K^{\dim Z}\cdot [Z]$ is the degree of $Z$ for $\CL_K$ on the generic fiber. The definition works only if $\deg _{\CL_K}(Z)\neq 0$, but this condition is easy to satisfy. In fact, $\CL_K$ is usually ample on $X$ in practice.  

\subsection{Arithmetic divisors}

There is a notion of arithmetic divisors equivalent to the notion of hermitian line bundles. 
For completeness, we recall the definition briefly. 

An arithmetic divisor $\ol \CD=(\CD, g_\CD)$ consists of a Cartier divisor $\CD$ on $\CX$ and a Green function $g_\CD$ of $\CD$ on $\CX$ invariant under the action of the complex conjugate. 
Here a function 
$g_\CD:\CX(\BC)\setminus |\CD|(\CC)\to \RR$ is called a Green function (resp. Green function of smooth type) of 
$\CD$ on $\CX$ if for any open subset $U$ of $\CX(\CC)$ on which $\CD(\CC)$ is defined by a single equation $f$, the function $g_\CD+\log |f|$ on $U\setminus |\CD|(\CC)$ can be extended to a continuous function (resp. smooth function) on $U$. 

Any global rational function $f$ on $\CX$ induces a principal arithmetic divisor 
$\wh\div(f):=(\div(f), -\log |f|)$. 
Two arithmetic divisors are linearly equivalent if their difference is a principal arithmetic divisor. 

Denote by $\wh\Pic(\CX)$ the group of isometry classes of hermitian line bundles, and denote by $\wh{\CaCl}(\CX)$ the group of linear equivalence classes of arithmetic divisors. 
There is a canonical isomorphism $\wh\Pic(\CX)\to \wh{\CaCl}(\CX)$, which sends a hermitian line bundle $\CLL$ on $\CX$ to the linear equivalence class of the arithmetic divisor 
$\wh\div(\ell):=(\div(\ell), -\log \|\ell\|)$.
Here $\ell$ is any nonzero rational section of $\CL$ on $\CX$.
In this isomorphism, hermitian line bundles with smooth metrics corresponds to arithmetic divisors with Green functions of smooth type.

\subsection{Intersection theory of arithmetic Chow cycles}

Now let us sketch the arithmetic intersection theory of arithmetic Chow cycles of Gillet--Soul\'e \cite{GS1}. 

Let $\CX$ be a regular scheme which is flat and quasi-projective over $\ZZ$. 
Note that we do not require $\CX$ to be projective for the moment, but we require it to be regular here. 
An {arithmetic Chow cycle} of codimension $p$ on $\CX$ is a pair $\overline{\CZ}=(\CZ,g_\CZ)$, where $\CZ$ is a Chow cycle of codimension $p$ on $\CX$ and $g_\CZ$ is a Green current of $\CZ(\CC)$ on $\CX(\CC)$. 
Namely, $g_\CZ$ is a $(p-1,p-1)$-current on $\CX(\CC)$ such that
$$
dd^c g_\CZ+\delta_{\CZ(\CC)}=[\omega]
$$
for some smooth $(p,p)$-form $\omega$ on $\CX(\CC)$, and we further require the complex conjugate transfers $g_\CZ$ to $(-1)^{p-1} g_\CZ$. 

For a pair $(\CY,f)$ where $\CY\subset\CX$ is a closed integral subscheme of codimension $p-1$ and $f$ is a rational function on $\CY$, define the corresponding {principal arithmetic Chow cycle} to be a Chow cycle of codimension $p$ on $\CX$ given by
$$
\wh\div_\CX(f)=\left(\div_\CX(f),[-\log|f|]_{\CX(\CC)}\right).
$$
Here both terms on the right are understood via push-forward by $\CY\to \CX$. 

The arithmetic Chow group of codimension $p$ on $\CX$ is defined as 
$$
\widehat{\mathrm{CH}}^p(\CX)=\widehat{\mathrm{Z}}^p(\CX)/\widehat{\mathrm{R}}^p(\CX).
$$
Here $\widehat{\mathrm{Z}}^p(\CX)$ is the group of 
arithmetic Chow cycles of codimension $p$ on $\CX$, and 
$\widehat{\mathrm{R}}^p(\CX)$
is the subgroup of $\widehat{\mathrm{Z}}^p(\CX)$ generated by all principal arithmetic Chow cycles of codimension $p$, and the images of the operators $\partial$ and $\ol\partial$ on currents on $\CX(\CC)$. 

Finally, Gillet-Soul\'e's intersection theory gives a bilinear pairing 
$$
\widehat{\mathrm{CH}}^p(\CX)_\QQ\times\widehat{\mathrm{CH}}^q(\CX)_\QQ\longrightarrow
\widehat{\mathrm{CH}}^{p+q}(\CX)_\QQ,
$$
which refines the intersection pairing of Chow groups of $\CX$.
If $\CX$ is further projective (or just proper) over $O_K$, there is a natural arithmetic degree map
$$
\wh\deg: \widehat{\mathrm{CH}}^{\dim \CX}(\CX)\longrightarrow\RR. 
$$
Composing with this degree map, we can get an intersection pairing
$$
\widehat{\mathrm{CH}}^p(\CX)\times\widehat{\mathrm{CH}}^{\dim \CX-p}(\CX)\longrightarrow
\RR.
$$

\subsection{Arithmetic Riemann--Roch theorem}

Let us introduce Gillet--Soul\'e's arithmetic Riemann--Roch theorem. 
Let $f:\CX\to \CY$ be a morphism, where $\CX$ and $\CY$ are regular schemes which are flat and quasi-projective over $\ZZ$. 
Assume that $f$ is smooth at all points of $\CX_\QQ$. 
Take a smooth hermitian metric on the relative tangent bundle $Tf=T_{\CX(\CC)/\CY(\CC)}$
such that its restriction to every fiber of $\CX(\CC)\to \CY(\CC)$ is K\"ahler.
Let $\ol\CE$ be a hermitian vector bundle on $\CX$ with smooth metrics. 
Similar to the case of relative dimension 1, there is a determinant line bundle 
$\det \mathrm{R}f_*(\ol\CE)_Q$ endowed with the Quillen metric. 
Then the arithmetic Riemann--Roch theorem of Gillet--Soul\'e \cite[Thm. 7]{GS2} is a formula
in $\widehat{\mathrm{CH}}^{1}(\CY)_\QQ$ of the form
$$
\det \mathrm{R}f_*(\ol\CE)_Q=f_*(\wh{\mathrm{ch}}(\ol\CE)\wh{\mathrm{Td}}(f))^{(1)}
+a( f_*(\mathrm{ch}(\CE_\CC) \mathrm{Td}(Tf_\CC) R(Tf_\CC)  )^{(0)} ).
$$
We refer to \cite[p. 160, Thm 1']{Sou} for more details on the terms on the right-hand side. 
We also refer to Faltings \cite{Fal90} for an extension of the formula to an equality in $\widehat{\mathrm{CH}}^*(\CY)_\QQ$, which is an exact analogue of Grothendieck's Riemann--Roch theorem in algebraic geometry. 

In the case that $\CY=\Spec\ZZ$, and $\CEE=\CLL$ is a hermitian line bundle, the theorem 
computes $\wh\deg \det H^*(\ol\CL)_Q$ in terms of arithmetic intersection numbers involving $\CLL$ and $Tf$. This is an arithmetic analogue of Hirzebruch's Riemann--Roch theorem in algebraic geometry. 

An importance consequence of the arithmetic Riemann--Roch theorem is
 an arithmetic Hilbert--Samuel formula, which we will introduce with a precise statement in the following subsection.

In 1991, Vojta \cite{Voj} gave the second proof of the Mordell conjecture, coming out as a deep analogue of the proof of the classical Thue--Siegel--Roth theorem in Diophantine approximation. 
Vojta's key idea is to find a global section of a small norm of a hermitian line bundle on a 3-dimensional arithmetic variety to bound heights of rational points, and his tool to prove the existence of such a section is to apply the arithmetic Riemann--Roch theorem.

\subsection{Arithmetic Hilbert--Samuel formula}

Now we introduce an arithmetic Hilbert--Samuel formula, which is the foundation of the arithmetic positivity to be introduced in the next section.

Let $K$ be a number field, $\CX$ be a (projective) arithmetic variety over $O_K$, and $\ol \CL$ be a hermitian line bundle on $\CX$ (with continuous metrics). 
Think $H^0(\CX, \ol \CL)$ as the vector bundle $H^0(\CX, \CL)$ on $\Spec O_K$ endowed with the supremum norm
$$\|\ell\|_{v,\sup}:=\sup _{x\in \CX_v (\BC)} \|\ell(x)\|_v$$
on $H^0(\CX, \CL_v)$ for any archimedean place $v$ of $K$. 
 This gives a supremum norm 
$$\|\ell\|_{\sup}:=\sup _{x\in \CX(\BC)} \|\ell(x)\|
=\sup_v \|\ell\|_{v,\sup}$$
on the Euclidean space 
$$H^0(\CX,  \CL)_\RR=H^0(\CX,  \CL)\otimes_{\ZZ} \BR\simeq \bigoplus _{v\mid \infty} H^0(\CX, \CL_v).$$

Then we have the set of effective sections given by
$$\wh H^0(\CX, \ol \CL):=\{\ell \in H^0(\CX,  \CL):  \|\ell\|_{\sup}\le 1  \}.$$
We want to count the ``dimension''
$$\wh h^0(\CX,\ol \CL):=\log \#\wh H^0(\CX, \ol \CL).$$
We say that $\CLL$ is \emph{effective} if $\wh h^0(\CX,\ol \CL)>0$. 
As in the case of arithmetic curves, this is again in the setting of Minkowski's geometry of numbers. 
Thus we define the Euler characteristic of $\ol \CL$ by
$$\chi (\ol \CL)=-\log \vol (H^0(\CX, \CL)_\RR/H^0(\CX, \CL))$$
where the Haar measure on $H^0(\CX, \CL)_\RR$ is normalized so that the volume of the unit ball of $\|\ell\|_{\sup}$ in $H^0(\CX, \CL)_\RR$ 
is $1$. 

For a hermitian line bundle $\ol \CL$ with smooth metrics, we say that $\ol \CL$ is {\em relatively semipositive} if $\ol \CL$ has semipositive curvature forms at infinite places and non-negative degrees with all closed curves in special fibers.

\begin{thm}[arithmetic Hilbert--Samuel formula, Gillet--Soul\'e]\label{thm-hs-herm}
Let $\CX$ be an arithmetic variety over $O_K$ of absolute dimension $d$, and let 
$\ol \CL$ be a hermitian line bundle on $\CX$ with smooth metrics. 
Assume that $\CL_K$ is ample and that $\ol \CL$ is relatively semipositive. Then for any hermitian line bundle $\ol \CM$ on $\CX$,  we have as $n\to \infty$, 
$$\chi (n\ol \CL+\ol \CM)=\frac {n^{d}}{d!}\ol \CL^d +o(n^{d}),$$
$$\wh h^0 (n\ol \CL+\ol \CM)\ge \frac {n^{d}}{d!}\ol \CL^d +o(n^{d}).$$
\end{thm}

When $\CX$ is regular, $\CL$ is ample in $\CX$, and $\ol \CL$ has positive curvatures, this was first proved by 
Gillet--Soul\'e as a consequence of their arithmetic Riemann--Roch theorem in \cite{GS2}, where they also used an estimate of analytic torsions of Bismut--Vasserot \cite{BV}.  A much shorter proof by induction on $\dim \CX$ was later found by Abbes--Bouche \cite{AB95}. 
The general case of the theorem is reduced to the smooth case by \cite{Zha1} using  resolution of singularities of $\CX_K$. 

One consequence of the formula is that $\wh h^0(n\ol \CL)\ne 0$ for sufficiently large $n$ if $\ol \CL^d>0$. To get more information about the growth of $\wh h^0(n\ol \CL)$, we need stronger positivity of $\ol \CL$, which will be discussed in the following section.

\section{Arithmetic positivity}\label{app sec positivity}

In this section, we introduce various positivity results of hermitian line bundles. The materials include arithmetic ampleness of Zhang \cite{Zha1}, 
the arithmetic bigness theorem of Yuan \cite{Yua1}, 
and various properties on volumes of hermitian line bundles by 
Chen \cite{Che1, Che2, Che3}, Moriwaki \cite{Mor5}, and Yuan \cite{Yua2}. 
We also mention the arithmetic Hodge index theorem for hermitian line bundles by Moriwaki \cite{Mor96}. 
For more details on these subjects, see Yuan \cite{Yua3} and Yuan--Guo \cite[\S7]{YG25}.

The theory of positivity of hermitian line bundles is significantly inspired by the theory of positivity of line bundles in algebraic geometry. 
For a complete exposition of the latter, we refer to Lazarsfeld \cite{Laz1}.

\subsection{Arithmetic ampleness} \label{app subsec ample}
As in classical algebraic geometry, we study the numerical criterion for arithmetic ampleness.

Let $K$ be a number field.
Let $\CX$ be an arithmetic variety over $O_K$ of absolute dimension $d$, and  
$\CLL$ be a hermitian line bundle on $\CX$ with smooth metrics. 
We say that $\ol \CL$ is {\em arithmetically ample} if for any hermitian line bundle $\ol \CM$, the $O_K$-module $H^0(\CX, n\CL+\CM)$
is generated by the subset $\wh H^0(\CX, n\CLL+\CMM)$ for every sufficiently large integer $n$. 
We say that $\ol \CL$ is {\em arithmetically positive} if $\CL_K$ is ample on $\CX_K$, $\ol \CL$ is relatively semipositive, and the height function 
has a positive lower bound, i.e. for some $c>0$, 
$h_{\ol \CL}(x)>c$  for all $x\in \CX(\ol K)$. 

\begin{thm}[arithmetic Nakai--Moishezon criterion, Zhang \cite{Zha1}] \label{app NM}
Assume that $\CL_K$ is ample on $\CX_K$, and that $\ol \CL$ is relatively semipositive. Then the following conditions are equivalent:
\begin{enumerate}
\item $\ol \CL$ is arithmetically positive;
\item $\ol \CL$ is arithmetically ample;
\item for any integral subvariety $Y$ of $\CX_K$, the height $h_{\ol \CL}(Y)>0$.
\end{enumerate}
Moreover, under these equivalent conditions, we have the following arithmetic Hilbert--Samuel formula for any hermitian line bundle $\CMM$ on $\CX$:
$$\wh h^0(n\ol \CL+\ol \CM)=\frac {n^{d}}{d!}\ol \CL^{d} +o(n^{d}), \qquad n\to\infty.$$
\end{thm}

This theorem is proved by Zhang \cite{Zha1} under a strong condition that the metric on $\ol \CL$ is semi-ample. Moriwaki \cite{Mor7} proved that this condition is equivalent to semipositivity of curvatures. 

Motivated by Szpiro's idea in \cite{Szp90}, $h_{\ol \CL}(\CX)$ should be related to the essential minimum of the height function $h_{\ol \CL}$ on $\CX(\ol K)$.
More precisely, we define the successive minima of $\ol \CL$ by the following formula:
$$e_i(\ol \CL)=\sup_{\cod Y=i}\inf _{x\in (\CX_K\setminus Y)(\ol K)}h_{\ol \CL}(x), \qquad i=1, \cdots {d},$$
where $Y$ runs through closed subvarieties of $\CX_K$ of codimension $i$ for each $i$, and the empty subvariety has codimension $d$.
The following theorem of successive minima gives a quantitative relation between these terms. 

\begin{thm}[successive minima, Zhang \cite{Zha1}] \label{app sm1}
Assume that $\CL_K$ is ample, and that $\ol \CL$ is relatively semipositive. Then
$$e_1(\ol \CL)\geq h_{\ol \CL}(\CX)\geq
\frac 1{d}\left(e_1(\ol \CL)+\cdots +e_{d}(\ol \CL)\right).$$
\end{thm}

We say that a hermitian line bundle $\CLL$ on an arithmetic variety $\CX$ is \emph{nef} if the following conditions hold:
\begin{itemize}
\item[(1)] the (continuous) metric of $\CLL$ is semipositive in the sense that the curvature current is positive,
\item[(2)] $\CLL$ has a non-negative (arithmetic) degree on any 1-dimensional closed integral  subscheme of $\CX$.
\end{itemize}
A hermitian line bundle (or its metric) is \emph{integrable}
 if it is isometric to the difference of two nef hermitian line bundles. 
The intersection theory of hermitian line bundles extends to nef hermitian line bundles by approximation, and thus extends to integrable hermitian line bundles by linearity.  

A quick consequence of the theorem is the following arithmetic analogue of Kleiman's theorem. 

\begin{cor}
Let $\CX$ be an arithmetic variety of dimension $d$. Let $\CLL_1, \cdots, \CLL_d$ be nef hermitian line bundles on $\CX$.
Then for any closed integral subscheme $\CY$ of $\CX$, we have 
$\CLL_1\cdots \CLL_{\dim \CY}\cdot [\CY]\geq 0$. 
\end{cor}

\begin{proof}
It suffices to assume $\CY=\CX$. 
We first prove the case that $\CLL_1, \cdots, \CLL_d$ are arithmetically positive. 
By Theorem \ref{app NM}, by replacing $\CLL_d$ by a positive multiple if necessary,
we can assume that there is a nonzero effective section $\ell_d\in \wh H^0(\CX,\CLL_d)$.   
Consider the intersection formula
$$
\ol \CL_{1}\cdots \ol \CL_{d}= \ol \CL_{1}\cdots  \ol \CL_{d-1}\cdot [\div(\ell_d)]
- \int_{\CX(\BC)}\log \|\ell_d\| c_1(\ol \CL_1)\cdots c_1(\ol \CL_{d-1}).
$$
It is non-negative by induction on $d$.

For general nef $\CLL_1, \cdots, \CLL_d$, fix an arithmetically positive hermitian line bundle $\CAA$ on $\CX$.
For any $i=1,\cdots, d$, and for any positive integer $m$, $m\CLL_i+\CAA$ is arithmetically positive. Then we have
$(m\CLL_1+\CAA)\cdots (m\CLL_d+\CAA)\geq 0$. 
This implies $\CLL_1\cdots \CLL_{d}\geq 0$. 
\end{proof}

\subsection{Arithmetic bigness}

In practice, we usually need an estimate of $\wh h^0$ similar to the arithmetic Hilbert--Samuel formula under weaker positivity conditions. 
The following arithmetic version of Siu's inequality gives an estimate, which is sufficient in application for many situations. 

\begin{thm}[arithmetic bigness, Yuan \cite{Yua1}] \label{app bigness}
Let $\CX$ be an arithmetic variety of absolute dimension $d$, and let $\ol \CL=\ol \CL_1-\ol \CL_2$ be the difference of two arithmetically positive hermitian bundles $\ol \CL_1$ and $\ol \CL_2$ on $\CX$. 
Then for any hermitian line bundle $\ol \CM$ on $\CX$, we have 
$$\chi(n\ol \CL+\ol \CM)\ge \frac {n^d}{d!} (\ol \CL_1^d-d\, \ol \CL_1^{d-1}\ol \CL_2)+o(n^d),$$
$$\wh h^0(n\ol \CL+\ol \CM)\ge \frac {n^d}{d!} (\ol \CL_1^d-d\, \ol \CL_1^{d-1}\ol \CL_2)+o(n^d).$$
\end{thm}

The proof of the theorem is inspired by the proof of the arithmetic Hilbert--Samuel formula by Abbes--Bouche \cite{AB95}. 

The theorem can be applied to the case that $\CLL$ is a small perturbation of an arithmetically positive hermitian line bundle, in which case the bound is very sharp. 
As we will see, an application of this theorem is an extension of the equidistribution theorem of Szpiro--Ullmo--Zhang from strictly positive metrics to semipositive metrics. 
Another interesting application is to simplify Vojta's proof of the Mordell conjecture in \cite{Voj}. 
As mentioned in the previous section, Vojta's proof applied Gillet--Soul\'e's arithmetic Riemann--Roch theorem to certain arithmetic variety of dimension 3, where he made great efforts to bound the contribution of analytic torsion and higher cohomology. 
The bigness theorem can be used to simplify this step drastically.
We refer to \cite[\S10]{YG25} for more details on this approach.

\subsection{Arithmetic volumes}

Let $\CX$ be an arithmetic variety of absolute dimension $d$, and let $\ol \CL$ be  a hermitian bundle on $\CX$. 
As in the geometric case, it is natural to define the \emph{arithmetic volume} of $\ol \CL$ by
$$\wh \vol (\ol \CL)=\limsup _{n\to\infty} \frac {d!}{n^d} \wh h^0(\CX,n\ol \CL).$$
We first have the following result. 

\begin{thm}[convergence, Chen \cite{Che1}, Yuan \cite{Yua2}] \label{app volume1}
The ``$\limsup$'' in the definition of $\wh \vol (\ol \CL)$ is a limit; i.e.
$$\wh \vol (\ol \CL)=\lim _{n\to\infty} \frac {d!}{n^d} \wh h^0(\CX, n\ol \CL).$$
\end{thm}

We say that the hermitian line bundle $\ol \CL$ is \emph{big} if $\vol (\ol \CL)>0$. 
Note that arithmetic bigness and arithmetic nefness are stable under pull-back via generically finite and surjective morphisms. 

It is always convenient to extend the volume function $\wh\vol:\wh\Pic(\CX)\to \RR$  to a volume function $\wh\vol:\wh\Pic(\CX)_\QQ\to \RR$ by setting $\wh \vol (a \ol \CL)=a^d\wh \vol (\ol \CL)$ for positive rational numbers $a$. 
An element of $\wh\Pic(\CX)_\QQ$ is called a hermitian $\QQ$-line bundles, and it is called arithmetically positive (resp. nef, big, effective) if some positive integer multiple of it is an arithmetically positive (resp. nef, big, effective) hermitian line bundle on $\CX$. 

It is easy to extend the arithmetic Hilbert--Samuel formula and the arithmetic bigness theorem to nef hermitian line bundles. 

\begin{cor} \label{app volume2}
Let $\CX$ be an arithmetic variety of absolute dimension $d$. Then the following are true. 
\begin{enumerate}
\item (arithmetic Hilbert--Samuel) 
Let $\ol \CL$ be a nef hermitian line bundle on $\CX$. Then 
$$\wh\vol (\ol \CL)=\ol \CL^d.$$
\item (arithmetic bigness) 
Let $\ol \CL_1$ and $\ol\CL_2$ be nef hermitian line bundles on $\CX$. Then 
$$ \wh\vol (\ol \CL_1-\ol \CL_2)\ge \ol \CL_1^d-d\, \ol \CL_1^{d-1}\ol \CL_2.$$
\end{enumerate}
\end{cor}
\begin{proof}
Note that we do not assume the metrics of $\CLL, \ol \CL_1, \ol\CL_2$ to be smooth, but we can convert them to the smooth case by the regularization theorem in \cite[Cor. C]{CGZ}, which asserts that any semipositive continuous metric on a complex line bundle is an increasing limit of semipositive smooth metrics. 

Take an arithmetically positive hermitian line bundle $\CAA$ on $\CX$.
Let $t$ be a small positive rational number.  
For (2), write 
$$
\ol \CL_1-\ol \CL_2=(\ol \CL_1+t\CAA)-(\ol \CL_2+t\CAA),
$$
apply Theorem \ref{app bigness} to the difference on the right-hand side, and take limit at $t\to 0$. 
Moriwaki \cite{Mor5} obtained part (1) by the continuity of the volume function in Theorem \ref{app volume3}(2). 
A simpler proof of part (1) is to take limits at $t\to 0$ of
$$\wh\vol (\ol \CL)\leq \wh\vol (\ol \CL+t\CAA)=(\ol \CL+t\CAA)^d$$
and
$$\wh\vol (\ol \CL)=\wh\vol ((\ol \CL+t\CAA)-t\CAA)\geq (\ol \CL+t\CAA)^d-d\, (\ol \CL+t\CAA)^{d-1}\cdot t\CAA.$$
\end{proof}

Besides the above quantitative results, 
 we also have the following properties of the volume function. 

\begin{thm} \label{app volume3}
Let $\CX$ be an arithmetic variety of absolute dimension $d$, and let $\ol \CL, \ol \CL_1, \ol\CL_2$ be hermitian bundles on $\CX$. Then the following are true. 
\begin{enumerate}

\item (birational invariance, Moriwaki \cite{Mor5}) 
If $\pi:\CX'\to \CX$ is a birational morphism of arithmetic varieties, then $\wh\vol(\CX',\pi^*\CLL)=\wh\vol(\CX,\CLL)$.

\item (continuity, Moriwaki \cite{Mor5}) The volume function $\wh\vol:\wh\Pic(\CX)_\QQ\to \RR$ is continuous in the sense that 
$$\lim_{t\to 0 }\wh\vol (\ol \CL_1+t\ol \CL_2)=\wh\vol (\ol \CL_1).$$

\item (differentiability, Chen \cite{Che3}) 
If $\ol \CL_1$ is big, then
$$\frac d{dt} \Big |_{t=0}\wh\vol (\ol \CL_1+t\ol \CL_2)=d \, \langle\ol \CL_1^{d-1} \rangle
\cdot \ol \CL_2,$$
where the right-hand side is the positive intersection number. 
In particular, if $\ol\CL_1$ is big and nef, then 
$$\frac d{dt} \Big |_{t=0}\wh\vol (\ol \CL_1+t\ol \CL_2)=d \, \ol \CL_1^{d-1} \cdot \ol \CL_2.$$

\item (log-concavity, Yuan \cite{Yua2})  If $\ol \CL_1$ and $\ol \CL_2$ are effective, then 
$$\wh \vol (\ol \CL_1+\ol \CL_2)^{ 1/d}\ge \wh \vol (\ol \CL_1)^{1/d}+\wh\vol (\ol \CL_2)^{1/d}.$$

\item (arithmetic Fujita approximation, Yuan \cite{Yua2}, Chen \cite{Che2}) 
If $\CLL$ is big, then for any $\epsilon>0$, there exist an arithmetic variety $\CX'$ with a birational morphism $\pi:\CX'\to \CX$, 
and an arithmetically positive hermitian $\QQ$-line bundle $\CAA$ on $\CX'$, such that 
$\pi^*\CLL-\CAA$ is effective on $\CX'$ and that 
$$ \wh\vol(\CAA)\geq\wh\vol(\CLL)-\epsilon. $$
\end{enumerate}
\end{thm}

The proof of Theorem \ref{app volume3}(1)(2) by Moriwaki \cite{Mor5} is inspired by the proof of Theorem \ref{app bigness} by Yuan \cite{Yua1}. 
The approach of Theorem \ref{app volume1} and Theorem \ref{app volume3}(4)(5) by Yuan \cite{Yua2} is based on an arithmetic version of the Okounkov body constructed by 
Okounkov \cite{Ok1, Ok2}, Lazarsfeld--Musta\c t\v a \cite{LM}, and Kaveh--Khovanskii \cite{KK1, KK2}. 
The proof of Theorem \ref{app volume3}(3) by Chen \cite{Che3} is given by a beautiful combination of Corollary \ref{app volume2}(1)(2) and Theorem \ref{app volume3}(4)(5).

\subsection{Hodge index theorem for hermitian line bundles}

We close this section by Moriwaki's generalization of the arithmetic Hodge index theorem of Faltings \cite{Fal1} and Hriljac \cite{Hri} (cf. Theorem \ref{app hodge index1})
 to high-dimensional arithmetic varieties.

\begin{thm}[Hodge index Theorem for hermitian line bundles, Moriwaki \cite{Mor96}] 
\label{app hodge index2}
Let $\pi: \CX\to \Spec O_K$ be an arithmetic variety such that $\CX_K$ is smooth and geometrically connected of dimension $d$. 
Let $\ol \CL_1, \cdots, \ol \CL_{d-1}, \ol \CM$ be integrable hermitian line bundles on $\CX$ such that $\ol \CL_1,\cdots, \ol \CL_{d-1}$ are arithmetically positive. 
Assume that $\CL_{1K}\cdots \CL_{d-1, K}\cdot \CM_K=0$ on the generic fiber. Then 
$$\ol \CL_1\cdots \ol \CL_{d-1}\cdot\ol \CM^2\le 0,$$
Moreover, the equality holds if and only if $\ol \CM=\pi^*\ol \CN$ for some hermitian line bundle $\ol\CN$ on $\Spec O_K$.
\end{thm}

\section{Adelic line bundles on projective varieties}
\label{app sec adelic}

In this section, we review the theory of adelic line bundles on projective varieties introduced by Zhang \cite{Zha2}, which is used in \S\ref{sec restricted} to compare the new notion of adelic line bundles on quasi-projective varieties.
Various positivity results for hermitian line bundles imply similar statements for adelic line bundles by approximation, but the Hodge index theorem for adelic line bundle of 
Yuan--Zhang \cite{YZ1} is significantly more difficult than its counterpart. 
 For updated details, see Yuan \cite{Yua3} and Yuan--Guo  \cite[\S8]{YG25}.

 \subsection{Metrized line bundles over local fields}
 \label{subsec-metrics}
Let $(K,|\cdot|)$ be a complete valuation field. The absolute value extends to a unique absolute value on the algebraic closure $\OK$ of $K$. 

Let $X$ be a projective variety over $K$, and $L$ be a line bundle on $X$. 
A \emph{$K$-metric} $\|\cdot\|$ of $L$ on $X$ is a collection of a $\overline K$-norm $\|\cdot\|$ on the fiber $L(x)=L_{\overline K}(x)$ of each algebraic point $x\in X(\overline K)$ which is \emph{Galois invariant} in the sense that, for any section $s$ of $L$ on an open subvariety $U$ of $X$, one has $\|s(x^\sigma)\|=\|s(x)\|$ for any $x\in U(\OK)$ and $\sigma\in \mathrm{Aut}(\overline K/K)$. 
 This is equivalent to a collection of norms of $L$ over the set $|X|_0$ of closed points of $X$.
 The pair $\ol L=(L, \|\cdot\|)$ is called a {\em metrized line bundle}. 
To introduce a continuity property, we separate the archimedean case and the non-archimedean case in the following. 

We first assume that $K$ is archimedean. By complex geometry, we have the usual notions of smooth metrics and continuous metrics of line bundles, and all continuous metrics are uniform limits of smooth 
metrics.  A continuous metric (or a metrized line bundle) is called {\em semipositive} (or {\em relatively nef}) if the metric is the uniform limit of a sequence of smooth metrics on $L$ with semipositive curvature forms. A metrized line bundle (or just its metric) is called {\em integrable} if it is isometric to the difference of two semipositive metrized line bundles.

  Assume that $K$ is non-archimedean in the following. Suppose $(\CX,\CM)$ is a \emph{projective model} of $(X, L^{\otimes e})$ over the valuation ring $O_K$ for some positive integer $e$; i.e.
 $\CX$ is a flat and projective integral scheme over $O_K$ with generic fiber $\CX_K$ isomorphic to $X$, and  $\CM$ is a line bundle on $\CX$ such that the generic fiber $(\CX_K,\CM_K)$ is isomorphic to $(X, L^{\otimes e})$.

Any point $x$ of $X(\overline K)$ extends to a point $\bar x$ of $\CX(O_{\overline K})$ by taking Zariski closure. 
Then the fiber $\CM(\bar x)$ is an $O_{\overline K}$-lattice of the one-dimensional $\overline K$-vector space $L^{\otimes e}(x)$.
It induces a $\overline K$-norm $\|\cdot\|'$ on $L^{\otimes e}(x)$ by the standard rule that 
$$
\CM(\bar x)=\{s\in L^{\otimes e}(x): \|s\|'\leq 1\}. 
$$
It thus gives a $\overline K$-norm $\|\cdot\|=\|\cdot\|'^{1/e}$ on $L(x)$.
Patching together, we obtain a  $K$-metric of $L$ on $X(\OK)$, and we call it \textit{the model metric induced by} $(\CX,\CM)$. 

The model metric (or the model metrized line bundle) is called \emph{semipositive} (or {\em relatively nef}) if $\CM$ is relatively nef in the sense that $\deg _{\CM}(\CZ)\ge 0$ for any projective curve $\CZ$ in the special fiber of $\CX$.

The induced model metrics are invariant under birational morphisms of projective models. Namely, if we further have a projective model $(\CX',\CM')$ of $(X, L^{\otimes e})$ dominating $(\CX,\CM)$, then they induce the same model metric on $L$. Here we say that $(\CX',\CM')$ dominates $(\CX,\CM)$ if there is a morphism $\pi: \CX'\to\CX$, extending the identity morphism on the generic fiber, such that $\pi^*\CM=\CM'$.

A $K$-metric $\|\cdot\|$ of a line bundle $L$ on $X$ is called  {\em m-continuous}  if it is the uniform limit of a sequence of model metrics $\|\cdot\|_n$ on $L$ in the sense that the function $\|\cdot\|_n/\|\cdot\|$ on $X(\OK)$ converges uniformly to 1. 
 An m-continuous $K$-metric $\|\cdot\|$ on a line bundle $L$ on $X$ (or the metrized line bundle $(L, \|\cdot\|)$) is called  {\em semipositive} (or {\em relatively nef}) if it is the uniform limit of a sequence of semipositive model metrics $\|\cdot\|_n$ on $L$. 
 A metrized line bundle (or just its metric) is called {\em integrable} if it is isometric to the difference of two  semipositive metrized line bundles. 
 
 In the name ``m-continuous'', the first letter comes from the word ``model''.
Note that our notion  ``m-continuous''  is strictly stronger than ``continuous on $X(\OK)$'' (in the non-archimedean case), where the latter means that for any section $s$ of $L$ on an open subvariety $U$ of $X$, the function $\|s\|$ is continuous on $U(\OK)$ under the  topology induced by the valuation of $\OK$. 
However, being m-continuous is equivalent to the condition that for all such $(U, s)$, the function $\|s\|$ extends to a continuous function on the Berkovich space $U^\an$. 
This relation is explored in the comparison in \S\ref{sec restricted}. 

To have a uniform terminology, in the archimedean case, smooth metrics are also called model metrics, and continuous metrics are also called m-continuous metrics.

 \subsection{Adelic line bundles} \label{app subsec adelic}

Let $K$ be a number field, and denote by $M_K$ the set of places of $K$. 
Let $X$ be a projective variety over $K$, and $L$ be a line bundle on $X$. 
For convenience, we write $X_{M_K}=\coprod_{v\in M_K} X_{K_v}$ for the disjoint union over all places $v$ of $K$.
By an $M_K$-\emph{metric} of $L$ on $X_{M_K}$, we mean a collection 
$(\|\cdot\|_v)_v$ of {\em m-continuous} $K_v$-metrics $\|\cdot\|_v$  of $L_{K_v}$ on $X_{K_v}$ over all places $v$ of $K$. 
The pair $\OL=(L, (\|\cdot\|_v)_v)$ is called a \emph{metrized line bundle on} $X$. 

We say that a sequence $\{(\|\cdot\|_{i,v})_v\}_{i\geq1}$ of $M_K$-metrics of $L$ on $X_{M_K}$ \emph{converges} to an $M_K$-metric $(\|\cdot\|_v)_v$ of $L$ on $X_{M_K}$ if 
there is a finite set $S$ of non-archimedean places of $K$, such that $\|\cdot\|_{i,v}=\|\cdot\|_{v}$ for any place $v\notin S$ and any $i$, and such that $\|\cdot\|_{i,v}/\|\cdot\|_{v}$ converges uniformly to 1 on $X(\OK_v)$ for every place $v$. 
In this case, we also say that the sequence $\{(L,(\|\cdot\|_{i,v})_v)\}_{i\geq1}$ 
\emph{converges} to $(L, (\|\cdot\|_v)_v)$ and that $(L, (\|\cdot\|_v)_v)$ is a limit of 
the sequence $\{(L,(\|\cdot\|_{i,v})_v)\}_{i\geq1}$. 

By an \emph{adelic metric} of $L$ on $X$, we mean an $M_K$-metric  
$(\|\cdot\|_v)_v$ of $L$ on $X_{M_K}$ satisfying the {\em coherence condition} 
 that there are a (non-empty) open subscheme $\CV$ of $\Spec O_K$, a projective and flat integral model $\CU\to \CV$ of $X$ over $\CV$, and a line bundle $\CL$ on $\CU$ extending $L$, so that for any closed point $v\in \CV$, the $K_v$-metric  $\|\cdot\|_v$ is induced by the integral model $(\CU\times_\CV O_{K_v}, \CL\times_\CV O_{K_v})$ of $(X,L)$. 
The pair $\OL=(L, (\|\cdot\|_v)_v)$ is called an \emph{adelic line bundle} on $X$.

Denote by $\wh\Picc (X)$ the groupoid of adelic line bundles on $X$, in which morphisms are defined to be isometries of adelic line bundles. 
Let $(\CV, \CU)$ be a pair as above; namely, 
$\CV$ is an open subscheme of $\Spec O_K$, and $\CU\to \CV$ is a projective and flat integral model of $X$ over $\CV$. 
Denote by $\wh\Picc (\CU)$ the groupoid of pairs $(\CL, (\|\cdot\|_v)_{v\notin \CV})$, where
 $\CL$ is a line bundle on $\CU$ extending $L$, and $\|\cdot\|_v$ is a m-continuous $K_v$-metric on $X(\OK_v)$ for every place $v\notin\CV$.
Then there is a fully faithful functor $\wh\Picc (\CU)\to \wh\Picc (X)$.  
In a suitable sense, 
 $\wh\Picc (X)$ is equivalent to the direct limit of $\wh\Picc (\CU)$ when varying $(\CU,\CV)$.

Let $(\CX, \CMM)$ be a {projective model of} $(X,L^{\otimes e})$ over $O_K$ for some positive integer $e$; i.e. $\CX$ is an arithmetic variety over $O_K$, and $\CMM=(\CM, \|\cdot\|)$ is a hermitian line bundle on $\CX$ with smooth metrics, such that the generic fiber $(\CX_K, \CM_K)=(X,L^{\otimes e})$. Then $(\CX, \CMM)$ induces an $M_K$-metric $(\|\cdot\|^{(\CX,\CMM)}_v)_v$ of $L$ on $X_{M_K}$. 
In fact, the metrics at archimedean places are just given by the $e$-th root of the hermitian metric of $\CMM$, and the metric at any non-archimedean place $v$ is just the model metric induced by the projective model $(\CX_{O_{K_v}},\CM_{O_{K_v}})$ of 
$(X_{K_v},L_{K_v}^{\otimes e})$. 
This induced $M_K$-metric $(\|\cdot\|^{(\CX,\CMM)}_v)_v$ is called a \emph{model adelic metric} of $L$ on $X_{M_K}$, and  $(L,(\|\cdot\|^{(\CX,\CMM)}_v)_v)$ is called a \emph{model adelic line bundle} on $X$. 

A model adelic line bundle on $X$ is said to be \emph{nef} if it can be induced by a projective model $(\CX, \CMM)$ over $O_K$ such that $\CMM$
is a nef hermitian line bundle on $\CX$. 
An  adelic line bundle on $X$ is said to be \emph{nef} if it is isometric to the limit of a sequence of nef model adelic line bundles on $X$.   
An  adelic line bundle on $X$ is said to be \emph{integrable} if it is isometric to the difference of two nef adelic line bundles on $X$.   
An adelic line bundle on $X$ is said to be \emph{relatively semipositive} (or \emph{relatively nef}) if its metric $\|\cdot\|_v$ is semipositive at every place $v$ of $K$.

Denote by 
$\wh\Picc(X)_\rmod$  (resp. $\wh\Picc(X)_\nef$,  $\wh\Picc(X)_\intb$)
the full subcategory of model (resp. nef, integrable) adelic line bundles in $\Picc(X)$.
Denote by 
$$
\wh\Pic(X), \quad \wh\Pic(X)_\rmod, \quad \wh\Pic(X)_\nef,\quad
\wh\Pic(X)_\intb
$$
the corresponding groups or semigroups of isomorphism classes of the categories.

The original motivation to introduce adelic line bundles comes from algebraic dynamics.
Let $(X,f,L)$ be a \emph{polarized dynamical system} over $K$; namely, $X$ is a  projective variety over $K$ with an endomorphism $f:X\to X$, and $L$ is an ample line bundle on $X$  
such that
$f^*L\simeq qL$ from some integer number $q>1$.
By Tate's limiting argument, Zhang \cite{Zha2} constructs an adelic line bundle $\overline L_f$ on $X$ with underlying line bundle $L$ such that $f^*\overline L_f\simeq q\overline L_f$. 
It turns out that $\OL_f$ is always nef with top self-intersection number 0.

 \subsection{Intersection numbers of integrable adelic line bundles}
 
Let us extend the intersection theory to adelic line bundles. So let $X$ be a projective variety over a number field $K$. 
By intersection theory of hermitian line bundles with smooth metrics,
for any closed subvariety $Z$ of $X$ of dimension $d$, the intersection pairings on integral models define a pairing 
$$[Z]:\wh\Pic (X)_\rmod ^{d+1}\lra \BR, \qquad (\ol L_1, \cdots, \ol L_{d+1})\longmapsto  \ol L_1\cdots  \ol L_{d+1}\cdot [Z].$$
More precisely, assume that $\ol L_i$ is induced by a projective model 
$(\CX_i, \CMM_i)$ over $O_K$ with $\CM_i|X=L_i^{e_i}$. 
By blowing-up and pull-back, we can assume that all $\CX_1=\cdots=\CX_{d+1}$.   
Then we can define 
$$ \ol L_1\cdots \ol L_{d+1}\cdot [Z]=\frac 1{e_1\cdots e_{d+1}} \CMM_1\cdots
\CMM_{d+1}\cdot [\ol Z],$$
where $\ol Z$ is the Zariski closure of $Z$ on $\CX_1$.

As each adelic line bundle is a limit of model adelic line bundles, the problem is about the continuity of the intersection pairing on the subspace $\wh\Pic (X)_\rmod$. However, this continuity does not hold even at archimedean places in general. Fortunately, this continuity holds for semipositive metrics. One main observation  in \cite{Zha2} is that the intersection pairing on model metrics induces an intersection pairing 
$$[Z]:\wh\Pic (X)_\intb ^{d+1}\lra \BR, \qquad (\ol L_1, \cdots, \ol L_{d+1})\longmapsto \ol L_1\cdots \ol L_{d+1} \cdot [Z].$$

Finally, for an integrable adelic line bundle $\ol L$ with $L$ ample on $X$, 
and for any closed subvariety $Z$ of $X$, 
we can define the arithmetic degree and the height by
$$\wh\deg_{\ol L}(Z)=\ol L ^{\dim Z+1}\cdot[Z], \qquad h_{\ol L}(Z)=\frac{\wh\deg_{\ol L}(Z)}{(\dim Z+1) \deg _{L}(Z)}.$$
If $Y$ is a closed subvariety of $X_\OK$, we define 
$$h_{\ol L}(Y)=h_{\ol L}(Y'),$$
where $Y'$ is the closed subvariety of $X$ corresponding to the Galois orbit of $Y$, i.e. the image of 
the composition $Y\to X_\OK\to X$. 

In particular, for any point $x\in X(\OK)$, we have 
$$
h_{\ol L}(x)=-\frac{1}{[K(x):K]}\sum_{v\in M_K} \sum_{\sigma:K(x)\to \OK_v} \log \|s(x)\|_v^{e_v},
$$
where $K(x)$ is the coefficient field of $x$ over $K$, $\sigma$ runs over all $K$-linear embedding $K(x)\to \OK_v$, 
and $s$ is any rational section of $L$ on $X$ without zero or pole at $x$. 
For the sake of the product formula, the factor $e_v=1$ for non-archimedean $v$ and for real $v$; $e_v=2$ for complex $v$.
This gives a height function 
$$
h_{\ol L}: X(\OK)\lra \RR,
$$
which is a Weil height for $L$ (up to a normalizing factor $[K:\QQ]$).

\subsection{Positivity of adelic line bundles}

Now we translate the previous positivity results from hermitian line bundles to adelic line bundles.
Let $X$ be a projective variety over a number field $K$ of dimension $d-1$, 
and  $\ol L$ be an adelic line bundle on $X$. 
We consider $H^0(X, \ol L)$ as the
vector space $H^0(X, L)$ over $K$ with a supremum norm on $H^0(X, L)_{K_v}$ at each place $v$ of $K$ given by
$$\|\ell\|_{v,\sup}:=\sup _{x\in X (\ol K_v)} \|\ell(x)\|_v.$$
Then  we have the set of effective sections by
$$\wh H^0(X, \ol L):=\{\ell \in H^0(X,  L):  \|\ell\|_{v,\sup}\le1, \forall v\},$$
and its {\em dimension} by
$$\wh h^0(X,\ol L):=\log \#\wh H^0(X, \ol L).$$
The {\em volume} of $\OL$ is defined by
$$\wh \vol (\ol L):=\limsup _{n\to \infty} \frac {d!}{n^{d}} \wh h^0(X,n\ol L).$$

We say that $\ol L$ is {\em effective} if $\wh h^0(X,\ol L)>0$. 
We say that $\ol L$ is {\em big} if $\wh \vol (\ol L)>0$. 
Recall that $\ol L$ is nef if it is a limit of model adelic line bundles induced by nef hermitian line bundles. 
We say that $\ol L$ is {\em arithmetically positive}  
if $L$ is ample, $\ol L$ is relatively semipositive, and $h_{\OL}(Y)>0$ for any closed subvariety $Y$ of $X$.
With some efforts, we can prove that  ``arithmetically positive'' implies ``nef''. 

By writing $\OL$ as the limit of model adelic line bundles induced by hermitian line bundles $\CMM_i$ with generic fiber $\CM_{i,K}\simeq e_iL$ on $X$, we can prove that Theorem \ref{app volume1} holds for adelic line bundles. 
 We actually have 
$$\wh \vol (\ol L)=\lim _{n\to \infty} \frac {d!}{n^{d}} \wh h^0(X,n\ol L)$$ 
and
$$
\wh \vol (\ol L)=\lim _{n\to \infty} \frac{1}{e_i^{d}} \wh \vol ( \CMM_i).
$$
Then we can also extend the function $\wh\vol:\wh\Pic(X)\to \RR$
to $\wh\vol:\wh\Pic(X)_\QQ\to \RR$ by homogeneity. 

With these notions, Corollary \ref{app volume2} 
and Theorem \ref{app volume3} hold for adelic line bundles. 
For convenience, we combine the statements in the following.

\begin{thm} \label{app volume4}
Let $X$ be a projective variety of dimension $d-1$ over a number field $K$, and let $\ol L, \ol L_1, \ol L_2$ be adelic line bundles on $X$. Then the following are true. 
\begin{enumerate}
\item (arithmetic Hilbert--Samuel) 
If $\ol L$ is nef, then 
$$\wh\vol (\ol L)=\ol L^d.$$
\item (arithmetic bigness) 
If $\ol L_1$ and $\ol L_2$ are nef, then 
$$ \wh\vol (\ol L_1-\ol L_2)\ge \ol L_1^d-d\, \ol L_1^{d-1}\ol L_2.$$
\item (log-concavity)  If $\ol L_1$ and $\ol L_2$ are effective, then 
$$\wh \vol (\ol L_1+\ol L_2)^{ 1/d}\ge \wh \vol (\ol L_1)^{1/d}+\wh\vol (\ol L_2)^{1/d}.$$
\item (birational invariance) 
If $\pi:X'\to X$ is a birational morphism of projective varieties over $K$, then $\wh\vol(X',\pi^*\OL)=\wh\vol(X,\OL)$.
\item (continuity) The volume function $\wh\vol:\wh\Pic(X)_\QQ\to \RR$ is continuous in the sense that 
$$\lim_{t\to 0 }\wh\vol (\OL_1+t\OL_2)=\wh\vol (\OL_1).$$
\item (differentiability) 
If $\OL_1$ is big, then
$$\frac d{dt} \Big |_{t=0}\wh\vol (\OL_1+t\OL_2)=d \, \langle\OL_1^{d-1} \rangle
\cdot \OL_2,$$
where the right-hand side is a generalized positive intersection number. 
In particular, if $\OL_1$ is big and nef, then 
$$\frac d{dt} \Big |_{t=0}\wh\vol (\OL_1+t\OL_2)=d \, \OL_1^{d-1} \cdot \OL_2.$$
\item (arithmetic Fujita approximation) 
If $\OL$ is big, then for any $\epsilon>0$, there exist a projective variety $X'$ over $K$ with a birational morphism $\pi:X'\to X$, 
and an arithmetically positive adelic $\QQ$-line bundle $\overline A$ on $X'$, such that 
$\pi^*\OL-\OA$ is effective on $X'$ and that 
$$ \wh\vol(\OA)\geq\wh\vol(\OL)-\epsilon. $$
Moreover, we can take $\OA$ to be a rational multiple of a model adelic line bundle on $X'$ induced by an arithmetically positive hermitian line bundle. 
\end{enumerate}
\end{thm}

\subsection{Theorem of successive minima}
Let $X$ be a projective variety of dimension $d-1$ over a number field $K$, and let $\ol L$ be an adelic line bundles on $X$. 
As in the hermitian case, we still have  the successive minimum
$$e_i(\ol L)=\sup_{\cod Y_K=i}\inf _{x\in (X\setminus Y)(\ol K)}h_{\ol L}(x), \qquad i=1, \cdots, {d},$$
where $Y$ runs through all closed subvarieties of $X$ of codimension $i$, and the empty subvariety has codimension $d-1$.
Then Theorem \ref{app sm1} implies the following adelic version. 

\begin{thm}[successive minimum, Zhang \cite{Zha2}] \label{app sm2} 
Assume that $L$ is ample, and that $\ol L$ is relatively semipositive. Then
$$e_1(\ol L) \geq h_{\ol L}(X)  \geq \frac 1{d}\left(e_1(\ol L)+\cdots +e_{d}(\ol L)\right).$$
\end{thm}

The most interesting case is when $e_1(\ol L)=e_2(\ol L)=\cdots=e_{d}(\ol L)$. In this case, we will describe the equidistribution theorem in the next section.

\subsection{Arithmetic Hodge index theorem}

We close this section by recalling the Hodge index theorem for adelic line bundles proved by Yuan--Zhang \cite{YZ1}, generalizing that of  
 Faltings \cite{Fal1} and Hriljac \cite{Hri} (cf. Theorem \ref{app hodge index1}) and 
that of Moriwaki \cite{Mor96} (cf. Theorem \ref{app hodge index2}). 
 
 For two integrable line bundles $\ol L$ and $\ol M$ on a projective variety $X$ over a number field $K$, we say that $\ol M$ is bounded by $\ol L$ if there is a positive integer $n$ such that both $n\ol L+\ol M$ and $n\ol L-\ol M$ are nef. 

\begin{thm}[Hodge index Theorem for adelic line bundles, Yuan--Zhang \cite{YZ1}]
\label{app hodge index3}
Let $\pi: X\to \Spec K$ be a geometrically connected and normal projective variety of dimension $d$. 
Let $\ol L_1, \cdots, \ol L_{d-1}, \ol M$ be integrable adelic  line bundles on $X$ such that $\ol L_1,\cdots, \ol L_{d-1}$ are arithmetically positive, and that
 $\ol M$ is bounded by each $\ol L_i$. 
 Assume that  $L_1\cdots L_{d-1}\cdot M=0$.  Then 
$$\ol L_1\cdots \ol L_{d-1}\cdot \ol M^2\le 0.$$
Moreover, the equality holds if and only if $\ol M=\pi^*\ol N$ for some adelic line bundle $\ol N$ on $\Spec K$.
\end{thm}

The difficult part of the theorem is the condition of the equality. We refer to \cite{YZ1}
for applications of this theorem to algebraic dynamics.

\section{Measures and Equidistribution}
\label{app sec equidistribution}

In this section, we introduce the Berkovich space of \cite{Ber1}, the
Chambert-Loir measure of \cite{CL}, 
the non-archimedean Calabi theorem of Yuan--Zhang \cite{YZ1}, 
the equidistribution theorem of Yuan \cite{Yua1}, and the proof of the Bogomolov conjecture by Ullmo \cite{Ull} and Zhang \cite{Zha4}.
For updated details, see Yuan \cite{Yua3}, and
Yuan--Guo 
\cite[\S8, \S9]{YG25}.

\subsection{Berkovich spaces}\label{app sec berkovich}

Let $(K,|\cdot|)$ be a complete valuation field, and $X$ be a scheme of finite type over $K$. 
Berkovich \cite{Ber1} introduced a canonical analytic space $X^\an$ associated to $X$ over $K$.  
Moreover, 
$X^\an$ is Hausdorff (resp. compact, path-connected) if and only if $X$ is separated
(resp. proper,  connected). 

If $K=\BC$, then $X^\an$ is just $X(\BC)$. 
If $K=\BR$, then $X^\an$ is the quotient of $X(\BC)$ by the action of the 
complex conjugate.

Now we recall Berkovich's definition of $X^\an$, which actually works for both non-archimedean $K$ and archimedean $K$. 
Let $U=\Spec(A)$ be an affine scheme of finite type over $K$.
Then $U^\an$ is defined to be the set of multiplicative semi-norms on $A$ extending the absolute value of $K$. 
Namely, $U^\an$ is the set of maps $\rho: A\to \RR_{\geq 0}$ satisfying:
\begin{enumerate}
\item (compatibility) $\rho|_{K}=|\cdot|$,
\item (triangle inequality) $\rho(a+b)\leq \rho(a)+\rho(b), \ \forall a,b\in A$,
\item (multiplicativity) $\rho(ab)= \rho(a)\rho(b), \ \forall a,b\in A$.
\end{enumerate}
Any $f\in A$ defines a map
$$|f|: U^\an \longrightarrow \RR, \quad \rho\longmapsto |f|_\rho:=\rho(f) $$
Endow $U^\an$ with the coarsest topology such that $|f|$
is continuous for all $f\in A$. 

For a general $K$-variety $X$, cover it by affine open schemes $U$.
Then $X^\an$ is obtained by glueing the corresponding $U^\an$ in the natural way. 
Each $U^\an$ is an open subspace of $X^\an$ by definition.

Denote by $|X|$ the underlying topological space of the scheme $X$, and by $|X|_0$ the subset of closed points of $|X|$.
There is a natural surjective map $X^\an\to |X|$. Every point of $|X|_0$ has a unique preimage, and thus there is a natural inclusion $|X|_0\hookrightarrow X^\an$. 

To describe the maps, it suffices to consider the affine case $U=\Spec(A)$.  
Then the map $U^\an\to |U|$ is just $\rho\mapsto \mathrm{ker}(\rho)$. 
The kernel of $\rho:A\to \RR_{\geq 0}$ is a prime ideal of $A$ by the multiplicativity. 
Let $x\in |U|_0$ be a closed point corresponding to a maximal idea $m_x$ of $A$. 
Then $A/m_x$ is a finite field extension of $K$, and thus has a unique valuation extending the valuation of $K$. 
This extension gives the unique preimage of $x$ in $U^\an$.

\subsection{Chambert-Loir measures} \label{app subsec measure}

Let $(K,|\cdot|)$ be a complete valuation field with a non-trivial absolute value, and $X$ be a projective variety of dimension $d$ over $K$. 
Let $\OL_1,\cdots, \OL_d$ be integrable metrized line bundles on $X$.
There is a signed Radon measure $c_1(\OL_1)\cdots c_1(\OL_d)$ on the analytic space $X^\an$ with total volume
$$
\int_{X^\an} c_1(\OL_1)\cdots c_1(\OL_d) =\deg_L(X). 
$$
The measure is positive if $\OL_1, \cdots, \OL_d$ are semipositive.

When $K=\CC$, then $X^\an$ is just $X(\CC)$. 
If both $X$ and $\|\cdot\|$ are smooth, then 
$c_1(\OL_1)\cdots c_1(\OL_d)$
is just the usual Monge--Amp\`ere measure on $X^\an$, defined as the wedge products of the Chern forms in complex analysis. 
Without the smoothness conditions, the measure is constructed by regularization by Bedford--Taylor \cite[Thm. 2.1]{BT}.  

When $K=\RR$, as $X^\an$ is a quotient of $X(\CC)$, the measure can be obtained by this quotient process. 

When $K$ is non-archimedean, $c_1(\OL_1)\cdots c_1(\OL_d)$ is the \emph{Chambert-Loir measure}, which is constructed by Chambert-Loir \cite{CL} when $K$ has a dense and countable subfield and extended to general $K$ by Gubler \cite{Gub1}. 

To illustrate the construction of $c_1(\OL_1)\cdots c_1(\OL_d)$, 
we further assume that the valuation of $K$ is discrete for simplicity.
By a limit process, 
we can assume that the metric of $\OL_i$ on $X$ is a model metric, i.e. it is induced by a projective model $(\CX,\CM_i)$ of $(X,L^{\otimes e_i})$ over $O_K$. 
Here we take a uniform $\CX$ for all $\OL_i$, and 
assume that $\CX$ is integrally closed in $X$ by passing to the normalization. 
Let $C_1,\cdots, C_r$ be the irreducible components of the special fiber of $\CX$. 
Then $\CX$ is regular at the generic points of $C_j$, so the multiplicity function 
$\ord_{C_j}: K(X)^\times \to \ZZ$ defines a norm on $K(X)$ by exponentiation. 
This gives a point $\xi_{C_j}\in X^\an$, which is called a Shilov point or a divisorial point. 
Then Chambert-Loir defines the measure by 
$$
c_1(\OL_1)\cdots c_1(\OL_d)= \frac{1}{e_1\cdots e_d}\sum_{j=1}^r  (\CM_1\cdots \CM_d\cdot [C_j])\ \delta_{\xi_{C_j}}
$$
Here $\delta_{\xi_{C_j}}$ is the Dirac measure supported at $\xi_{C_j}$. 

We can also interpret the measure in terms of functional analysis. 
Recall that ``m-continuous metrics'' and ``model metrics'' are introduced in \S \ref{subsec-metrics}, where m-continuous metrics are just uniform limits of model metrics. 
We say that a function $f:|X|_0\to \RR$ is {\em m-continuous} (resp. {\em a model function})
if $f=-\log \|1\|$ for some m-continuous $K$-metric (resp. model $K$-metric) of the trivial line bundle $\CO_X$.  
By Gubler \cite[Theorem 7.12]{Gub3}
(cf. \cite[Lem. 3.5]{Yua1}), the space of model functions on $|X|_0$ extends to $X^\an$ as a dense subset of the space $C(X^\an)$ of continuous functions on $X^\an$. 
Thus m-continuous functions on $|X|_0$ correspond to continuous functions on $X^\an$ exactly. 
To define the measure $c_1(\OL)^d$, it suffices to define the integral of continuous functions for this measure. Then we simply set 
$$\int_{X^\an} 
 f\, c_1( \ol L_1)\cdots  c_1(\ol L_{d})
 =\ol M\cdot \ol L_1\cdots  \ol L_{d},$$
where $\ol M=(\CO_X, \|\cdot\|)$ is an m-continuous metrized line bundle such that $f=-\log \|1\|$. 
The intersection can be defined for model metrized line bundles and extends to integrable ones by a limit process.

 \subsection{Non-archimedean Monge--Amp\`ere equation}
 
In this subsection, let $(K,|\cdot|)$ be a complete valuation field with a nontrivial absolute value, and $X$ be a geometrically integral projective variety of dimension $d$ over $K$. 
Denote by $\pi:X\to \Spec K$ the structure morphism.
 
Recall that Theorem \ref{app hodge index3} is a Hodge index theorem for adelic line bundles proved by Yuan--Zhang \cite{YZ1}. Its local version is the following theorem. 

\begin{thm}[local Hodge index theorem for metrized line bundles, \cite{YZ1}]
\label{app hodge index4}
Let $\ol L_1, \cdots, \ol L_{d-1}$ be semipositive metrized line bundles on $X$ such that the underlying line bundles $L_1,\cdots,  L_{d-1}$ are ample on $X$. 
Let  $\ol M$ be an integrable metrized line bundle on $X$ with underlying line bundle $M= \CO_X$.  
Assume that
 $\ol M$ is bounded by each $\ol L_i$. 
Then 
$$\ol L_1\cdots \ol L_{d-1}\cdot \ol M^2\le 0.$$
Moreover, the equality holds if and only if $\ol M\simeq \pi^*\ol N$ for some metrized line bundle $\ol N$ on $\Spec K$.
\end{thm}

Note that as $M$ is trivial, the intersection number is well-defined. 
The archimedean case of the theorem was previously proved by  
Kolodziej \cite{Kol03}. 
A quick consequence of the theorem is as follows. 
 
\begin{cor}[Calabi theorem for non-archimedean fields, \cite{YZ1}] \label{app calabi}
Let $L$ be an ample line bundle on $X$, and $\|\cdot\|_1$ and
$\|\cdot\|_2$ be two semipositive metrics on $L$. Then
$$
c_1(L, \|\cdot\|_1)^d=c_1(L, \|\cdot\|_2)^d
$$
if and only if $\ds\frac{\|\cdot\|_1}{\|\cdot\|_2}$ is a constant.
\end{cor}
\begin{proof}
Denote $\OL_i=(L, \|\cdot\|_i)$ as a metrized line bundle on $X$, and denote 
$\ds f=-\log \frac{\|\cdot\|_1}{\|\cdot\|_2}$ as a continuous function on $X^\an$. 
The equality of the measures gives
$$
\int_{X^\an} f c_1(\OL_1)^d=\int_{X^\an} f c_1(\OL_2)^d.
$$
This is just 
$$
(\OL_1-\OL_2)\cdot \OL_1^d=(\OL_1-\OL_2)\cdot \OL_2^d, 
$$
or equivalently 
$$
\sum_{i=0}^{d-1} (\OL_1-\OL_2)^2\cdot \OL_1^{i}\cdot \OL_2^{d-1-i}=0. 
$$
By Theorem \ref{app hodge index4}, every intersection number on the left-hand side is at most 0, and thus equal to 0. 
Then the theorem further implies that $f$ is a constant. 
\end{proof}

The history of Theorem \ref{app calabi} in the complex case is as follows.
In the 1950s, Calabi \cite{Ca1, Ca2} made the following famous conjecture.

\begin{quote}
\emph{Let $X$ be a compact complex manifold of dimension $d$ endowed with a K\"ahler form
$\omega$, 
and let $\Omega$ be a positive smooth volume form on $X$ such that 
$\int_X \Omega =\int_X \omega^d$.
Then there exists a smooth real-valued function $\phi$ on $X$ such that
$(\omega+i\partial \overline\partial \phi)^d=\Omega.$}
\end{quote}
Calabi proved that the function $\phi$ is unique up to constants (if it
exists).
The existence of the function is much deeper, and was finally solved by S. T.
Yau in the seminal paper \cite{Ya} in 1977.

Inspired by the above complex Monge--Amp\`ere equation, it is natural to consider the following non-archimedean analogue. 
\begin{quote}
\emph{Let $K$ be a non-archimedean field, $X$ be a geometrically integral projective variety of dimension $d$ over $K$, $L$ be an ample line bundle on $X$, and $d\mu$ be a positive Radon measure on $X^\an$ with total volume equal to $\deg_L(X)$. 
Find a necessary and sufficiently condition for $d\mu$ such that there exists a semipositive metric $\|\cdot\|$ on $L$ satisfying
$c_1(L, \|\cdot\|)^d=d\mu.$}
\end{quote}
From the work of Guedj--Zeriahi \cite{GZ07} in the complex case, it is natural to guess that a necessary and sufficiently condition in the non-archimedean case is still that $d\mu$ has measure zero at every pluripolar subset of $X^\an$ in a suitable sense. 
Many cases of this far-reaching problem are known, which include the work of 
Liu \cite{Liu11} on abelian varieties with totally degenerate reduction, and the work of Boucksom--Favre--Jonsson \cite{BFJ15} and its sequel \cite{BGJKM20, BGM22} for a breakthrough for general $X$.

\subsection{Equidistribution theorem}

Let $X$ be a projective variety of dimension $d$ over a number field $K$.
Let $\OL$ be a relatively semipositive adelic line bundle on $X$ such that $L$ is ample. 
For any place $v$ of $K$, the equilibrium measure  
$$\mu_{\overline L,v}:=\frac{1}{\deg_L(X)}c_1(\OL)^d$$
 on $X_{K_v}^\an$ is the normalized Monge--Amp\`ere measure. 

For any point $x\in X(\ol K)$, denote by $x'\in X$ the closed point corresponding to the Galois orbit of $x$. Then the base change $(x')_{K_v}$ is a finite subset of closed points of $X_{K_v}$, viewed as a finite subset of $X_{K_v}^\an$ via the natural injection $|X_{K_v}|_0\to X_{K_v}^\an$. This finite subset is just the image of functorial injection $(x')_{K_v}^\an\to X_{K_v}^\an$. 
We introduce a probability measure $\mu_{x,v}$ on $X_{K_v}^\an$ by 
$$\mu_{x,v}:=\frac{1}{\deg(x)}\sum_{z\in (x')_{K_v}^\an } \delta_z.$$

An infinite sequence $\{x_m\}_{m\geq 1}$ of $X(\overline K)$ is called \emph{generic} if any infinite subsequence is Zariski dense in $X$. 
An infinite sequence $\{x_m\}_{m\geq 1}$ of $X(\overline K)$ is called \emph{$h_{\overline L}$-small} if $h_{\overline L}(x_m)\to h_{\overline L}(X)$ as $m\rightarrow\infty$.

\begin{thm}[Equidistribution of small points, Yuan \cite{Yua1}] 
\label{app adelic equidistribution}
Let $X$ be a projective variety over a number field $K$ and let $\overline L$ be a relatively semipositive adelic line bundle on $X$ whose underlying line bundle $L$ is ample on $X$.
Let $\{x_m\}$ be a generic and $h_{\overline L}$-small sequence of $X(\overline K)$.
Then for any place $v$ of $K$, the probability measure $\mu_{x_m,v}$ converges weakly to the equilibrium measure 
$\mu_{\overline L,v}$ on $X_{K_v}^\an$.
\end{thm}

The equidistribution of small points originated in the landmark work of Szpiro--Ullmo--Zhang \cite{SUZ}, where they proved the equidistribution theorem assuming that $X$ is smooth, that $v$ is archimedean, and that the metric of $\OL$ at $v$ is smooth and strictly positive. 
Their proof is based on a variational principle on the arithmetic Hilbert--Samuel formula (cf. Theorem \ref{thm-hs-herm}).  
Since then, there were lots of works to generalize the theorem in dimension one, which include Bilu \cite{Bi97}, Autissier \cite{Au01}, 
Baker--Hsia \cite{BH05},  Baker--Rumely \cite{BR06}, Chambert-Loir \cite{CL}, and Favre--Rivera-Letelier \cite{FR06}. 
In particular, Chambert-Loir \cite{CL} introduced equilibrium measures on Berkovich spaces of any dimension, and his proof worked for general dimensions under a positivity assumption at $v$ similar to that of \cite{SUZ}.  
Finally, Yuan \cite{Yua1} proved the current version of the theorem, where the key is to apply 
 his  arithmetic bigness theorem (cf. Theorem \ref{app bigness}) in place of the arithmetic Hilbert--Samuel formula.

The function field analogue of the equidistribution theorem was obtained by Faber \cite{Fab} and Gubler \cite{Gub2} independently in slightly different settings. It is also worth mentioning that Moriwaki \cite{Mor3} prove an equidistribution over finitely generated fields over number fields based on the Moriwaki height.

While the original equidistribution theorem of \cite{SUZ} applies to abelian varieties over number fields, the current theorem applies to polarized dynamical systems over number fields. In fact, let $(X, f, L)$ be a {polarized dynamical system} over $K$.
As in \S\ref{app subsec adelic}, there is an adelic line bundle $\overline L_f$ on $X$ with underlying line bundle $L$ such that $f^*\overline L_f\simeq q\overline L_f$. 
It turns out that $\OL_f$ is relatively semipositive with top self-intersection number 0, so it satisfies the condition of the theorem. 
As a consequence, the equidistribution theorem is widely used in algebraic dynamics. 
We refer to \cite{Yua4} for more discussions on equidistribution in algebraic dynamics. 

Let $\OL$ be as in Theorem \ref{app adelic equidistribution}. 
 Note that Zhang's theorem of successive minima (cf. Theorem \ref{app sm2}) implies 
$e_1(\overline L) \geq h_{\overline L}(X)$. 
However, the existence of a generic and $h_{\overline L}$-small sequence is equivalent to the equality $e_1(\overline L) = h_{\overline L}(X)$. 
In general it is hard to construct $\OL$ with such a property, except that the dynamical case $\OL=\OL_f$ gives 
$$h_{\overline L}(X)=e_1(\ol L)=e_2(\ol L)=\cdots=e_{d}(\ol L)=0.$$

\subsection{Bogomolov conjecture}

One major application of the equidistribution theorem of Szpiro--Ullmo--Zhang \cite{SUZ} is the solution of the Bogomolov conjecture by Ullmo \cite{Ull} and Zhang \cite{Zha4}. 

\begin{thm}[Bogomolov conjecture,  Zhang \cite{Zha4}] 
\label{app bog}
Let $A$ be an abelian variety over a number field $K$, and let $L$ be a symmetric and ample line bundle on $A$. 
Let $X$ be a closed subvariety of $A$.
Assume that $X$ is not equal to the translation of an abelian subvariety by a torsion point. 
Then there is $\epsilon >0$ such that the set
$$\{ x\in X(\ol K): \wh h_L (x)<\epsilon\}$$
is not Zariski dense in $X$.
\end{thm}

The history of the Bogomolov conjecture is as follows. 
In \cite{Bog81}, Bogomolov asked some questions about finiteness of algebraic points 
of small N\'eron--Tate heights on a smooth projective curve $C$ of genus at least 2 over a number field. 
This is the case of curves in their Jacobian varieties in the form of Theorem \ref{app bog}.
 Zhang \cite{Zha1} solved an analogue of the conjecture for torus using his theory of arithmetic ampleness.
With the equidistribution theorem of Szpiro--Ullmo--Zhang \cite{SUZ}, 
Ullmo \cite{Ull} prove the original conjecture for curves in their Jacobian varieties, 
 and Zhang \cite{Zha4} prove the generalized conjecture for subvarieties of abelian varieties. 

We sketch the proof of \cite{Zha4}  as follows.  
Consider the morphism 
$$
\phi:X^m\lra A^{m-1}, \quad (t_1, \cdots, t_m)\longmapsto (t_1-t_2, \cdots, t_{m-1}-t_m). 
$$
For sufficiently large $m$, $\phi$ induces a birational morphism $X^m\to \phi(X^m)$. 
For the sake of contradiction, assume that the Bogomolov conjecture fails, so 
$X(\ol K)$ contains a generic and small sequence
$\{x_n\}_n$.
This sequence induces a generic and small sequence
$\{y_n\}_n$ of $X^m(\OK)$ with respect to the Neron--Tate height of $A^m$, where all $m$ components of $y_n$ 
in $X(\OK)$ are chosen from the sequence $\{x_n\}_n$ properly.
Then $\{\phi(y_n)\}_n$ is a generic and small sequence
 of $\phi(X^m)(\OK)$ with respect to the Neron--Tate height of $A^{m-1}$. 
 Apply the equidistribution theorem to $\{y_n\}_n$, and also apply it to $\{\phi(y_n)\}_n$. 
 We get a limit measure $\mu$ on $(X^m)_v(\CC)$
 and a limit measure $\mu'$ on $\phi(X^m)_v(\CC)$, where we fix an archimedean place $v$ of $K$.
 The compatibility of the sequences forces $\phi_*\mu=\mu'$. 
Note that the canonical metric of $L$ on $A_v(\CC)$ is smooth and strictly positive. 
Then $\mu$ and $\mu'$ actually come from strictly positive smooth top differential forms.
The equality $\phi_*\mu=\mu'$ of measures implies an equality $\mu=\phi^*\mu'$ of differential forms on $\phi^{-1}(U)$, where $U$ is any Zariski open subset of $(X^m)_v(\CC)$ such that $\phi^{-1}(U)\to U$ is an isomorphism.
By continuity, 
$\mu=\phi^*\mu'$ holds everywhere on $(X^m)_v(\CC)$.
As $\phi$ maps the diagonal $\Delta$ of $(X^m)_v(\CC)$ to a point, $\phi^*\mu'$ is not strictly positive at points of $\Delta$. 
This contradicts to the fact that $\mu$ is strictly positive. 

The geometric Bogomolov conjecture is an analogue of the Bogomolov conjecture over finitely generated function fields $K/k$.
Although the equidistribution theorem has an analogue over function fields as mentioned above, but there is no archimedean place to use, and thus 
 the proof of \cite{Ull, Zha4} does not work over function fields directly. 
 Another extra complication of the geometric case is given by the $(K/k)$-trace of $A$. 
However, by  careful analysis of the equilibrium measures at non-archimedean places, 
 Gubler \cite{Gub07b} proved the geometric Bogomolov conjecture for 
 abelian varieties with totally degenerate reduction at some place. 
Following the line, 
Yamaki \cite{Yamaki2016c, Yamaki2016a, Yamaki2017} reduced the geometric Bogomolov conjecture to the case of abelian varieties with good reduction everywhere, and proved the conjecture for $\dim(X)=1$ or $\mathrm{codim}(X)=1$.
By a completely different method using Betti maps of complex abelian schemes, Gao--Habegger \cite{GH} and Cantat--Gao--Habegger--Xie \cite{CGHX} proved the geometric Bogomolov conjecture for $\mathrm{char}(K) =0$. 
Finally, Xie--Yuan \cite{XY} proved the geometric Bogomolov conjecture for abelian varieties with good reduction everywhere, and thus complete the  geometric Bogomolov conjecture in all cases.

 \end{document}